\documentclass[11pt]{article}
\usepackage{amssymb}
\usepackage{amsmath,amsthm}
\usepackage{amsfonts}
\usepackage{leftidx}
\usepackage{mathrsfs}
\usepackage{authblk}
\usepackage{enumerate}  
\usepackage{abstract}
\usepackage{stmaryrd}
\usepackage{ulem}
\usepackage{slashed}
 \usepackage{float}
\usepackage{graphicx}
 \usepackage{hyperref}
  \usepackage{xcolor}
\usepackage{slashed}

\usepackage[top=1in, bottom=1.2in, left=1in, right=1in]{geometry}

\def\R{\mathbb{R}}

\def\Z{\mathbb{Z}}

\def\d{|\nabla|}

\def\p{\partial}

\def\be{\begin{equation}}
\def\ee{\end{equation}}

\newtheorem{theorem}{Theorem}[section]

\newtheorem{lemma}{Lemma}[section]
\newtheorem{proposition}{Proposition}[section]

\theoremstyle{definition}
\newtheorem{definition}{Definition}[section]
\newtheorem{convention}{Convention}[section]

\theoremstyle{remark}
\newtheorem{remark}{Remark}[section]

\numberwithin{equation}{section}
\title{Large data global solution of the 3D RVM system with cylindrical symmetry II: Pointwise estimates}
\author[$\star$]{Xuecheng Wang}
\affil[$\star$]{\small Tsinghua University \& BIMSA}

\date{  }

\begin{document}

 \maketitle 

\begin{abstract}
This is the second  part of a two-paper sequence  establishing the global existence of $3D$ relativistic Vlasov-Maxwell system (RVM) for \emph{arbitrarily large} smooth localized initial data with   {cylindrical symmetry}.

 In Part II, independently of the results in Part I \cite{PartI}, we establish two sets of pointwise estimates for the localized electromagnetic field. This is achieved by exploiting cylindrical symmetry, conservation laws, smoothing effects, and structural observations of the electromagnetic field, including its null structure, double null structure, and the magnetic field dichotomy. The first set provides a uniform upper bound for the localized electromagnetic field. The second set offers an upper bound that depends on the distance from the z-axis, providing better control when the electromagnetic field is located further away from the $z$-axis.

These two sets of estimates are fundamental to the iterative smoothing scheme (ISS), which is originated and inspired by   the work of Klainerman-Staffilani\cite{Klainerman3}.  These two sets were treated as black boxes in Part I.      A defining feature of these estimates is that they are assessed concerning the magnitude of the oscillation phase over characteristic time, effectively quantifying the smoothing effect.  It follows from these estimates that  cases without any smoothing effect are error terms. Addressing this constitutes the initial challenge in the bootstrap argument.
\end{abstract}

\setcounter{tocdepth}{2}

\tableofcontents
\section{Introduction}
In Part II,   we continue our study of the large data problem of  the following $3D$ relativistic Vlasov-Maxwell system, 
\be\label{mainequation}
(\textup{RVM})\qquad \left\{\begin{array}{l}
\p_t f + \hat{v} \cdot \nabla_x f + (E+ \hat{v}\times B)\cdot \nabla_v f =0,\\
\nabla \cdot E = 4 \pi \displaystyle{\int_{\R^3}  f(t, x, v) d v}, \qquad \nabla \cdot B =0, \\
\p_t E = \nabla \times B - 4\pi \displaystyle{\int_{\R^3} f(t, x, v) \hat{v} d v }, \\ 
 \p_t B  =- \nabla\times E,
\end{array}\right.
\ee
where $f:\R_t\times \R_x^3\times \R_v^3\longrightarrow [0, \infty)$ denotes the distribution function of particles, $ E, B  :\R_t\times \R_x^3 \longrightarrow \R^3 $ denote  the electromagnetic fields, and  $\hat{v}:=v/\sqrt{1+|v|^2}$.

   The main theorem of the  two-paper sequence is stated as follows.

\begin{theorem}\label{maintheorem}
  Given  any    cylindrical  symmetric initial data $E_0, B_0\in H^{s}(\R^3)$ ,  $  f_0(x,v)\in H^s(\R_x^3 \times \R_v^3)$, $s\in \mathbb{Z}_{+},s\geq 6$. Without loss of generality, we assume the  initial data are cylindrical symmetric with respect to the $z$-axis  in the following sense, 
\be\label{april25eqn1}
\begin{split}
   & \forall  \theta\in [0,2\pi], x,v\in \R^3,  \quad R=\begin{bmatrix}
\cos \theta & -\sin \theta & 0 \\ 
\sin \theta &  \cos\theta & 0\\ 
0 & 0 & 1\\  
\end{bmatrix},\\
 & f_0(Rx, Rv) = f_0( x,v), \quad  E_0(  Rx)= RE_0(x), \quad B_0( Rx)= RB_0(x). \\
 \end{split}
\ee
If the initial distribution function decays polynomially in the following sense,  
\be\label{assumptiononinitialdata}
\sum_{\alpha \in \mathbb{Z}_{+}^6,|\alpha|\leq s} \|(1+| x|+|v|)^{N_0}\nabla_{x,v}^\alpha f_0(x,v)\|_{L^2_{x,v}}< +\infty, \quad  N_0:= 10^{10},
\ee
then the relativistic Vlasov-Maxwell system \eqref{mainequation}  admits a global solution $ (f(t,\cdot, \cdot), E(t,\cdot), B(t,\cdot))\in  H^s(\R_x^3 \times \R_v^3)\times H^s(\R_x^3)\times H^s(\R_x^3)$. Moreover, the $L^\infty_x$-norm of the electromagnetic field grows at most polynomially over time. 

\end{theorem} 

 The theorem stated above was established in Part I \cite{PartI} based on 
 \textbf{fixed-time} pointwise estimates, which are derived in this paper (Section \ref{mainresultsthispaper}). The present work is dedicated to the rigorous derivation of these estimates. A comprehensive introduction, including a review of the current literature and open questions, can be found in Part I \cite{PartI}.

The rest of this section is organized as follows.
\begin{enumerate}
	\item[$\bullet$]   In subsection \ref{notationpre}, to ensure this paper's self-sufficiency, we reintroduce the notation used in Part I to ensure consistency between the two papers and to facilitate a better understanding of the main results. Additionally, we provide preliminary definitions, three important conservation laws,  and some key facts about the distribution function.

\item[$\bullet$] In subsection \ref{mainresultsthispaper}, we present   this paper's main results: two sets of pointwise estimates used in Part I \cite{PartI}.

\item[$\bullet$] In subsection \ref{mainingredientsfirstpart}, we  discuss the main ingredients and ideas used in obtaining these estimates. 
 
\item[$\bullet$] In subsection \ref{outlinefirstpaper}, we give the outline of this paper. 
\end{enumerate}

\subsection{Notation and Preliminaries}\label{notationpre}

\subsubsection{Notation}\label{notationsubsection}

\begin{definition}\label{relationdef}
 For any two numbers $A$ and $B$, we use  $A\lesssim B$, $A\approx B$,  and $A\ll B$ to denote  $A\leq C B$, $|A-B|\leq c A$, and $A\leq c B$ respectively, where $C$ is an absolute constant and $c$ is a sufficiently small absolute constant. We use $A\sim B$ to denote the case when $A\lesssim B$ and $B\lesssim A$.   For an integer $k\in\mathbb{Z}$, we use ``$k_{+}$'' to denote $\max\{k,0\}$ and  use ``$k_{-}$'' to denote $\min\{k,0\}$. For $x_0\in \R^3$, $r\in \R_{+}$, we use both the notation $B(x_0, r)$ and $B_r(x_0)$ to denote the set $\{x:|x-x_0|<r, x\in \R^3\}$. 
\end{definition}

\begin{convention}\label{conventionconst}
  We use the convention that all constants which only depend on the  \textbf{initial data}, e.g., the conserved quantities in   \eqref{conservationlaw}, will be treated as absolute constants. There are several absolute constants we use  throughout this two-paper sequence, which are $N_0:=10^{10},\epsilon:=10^{-7}, \iota:=10^{-4}. $
\end{convention}

We  fix an even smooth function $\tilde{\psi}:\R \rightarrow [0,1]$, which is supported in $[-3/2,3/2]$ and equals to one  in $[-5/4, 5/4]$. For any $k, k_1,k_2\in \mathbb{Z}$, we define $\psi_k, \psi_{\leq k},\psi_{\geq k}, \psi_{[k_1,k_2]} :\cup_{n\in \Z_+}\R^n\longrightarrow\R$, 
\[
\begin{split}
&\psi_{k}(x) := \tilde{\psi}(|x|/2^k) -\tilde{\psi}(|x|/2^{k-1}), \quad \psi_{\leq k}(x):= \sum_{l\leq k}\psi_{l}(x), \\
 &\psi_{\geq k}(x):= 1-\psi_{\leq k-1}(x), \quad \psi_{[k_1,k_2]}(x)=\sum_{k\in[k_1,k_2]\cap \Z}\psi_k(x). \\
\end{split}
\]
Moreover, we define the cutoff function $\psi_{l;\bar{l}}:\cup_{n\in \Z_+}\R^n\longrightarrow\R$ with the threshold $\bar{l}\in \Z$,  as follows,
  \be\label{cutoffwiththreshold}
\varphi_{l;\bar{l}}(x):=\left\{\begin{array}{ll}
\psi_{\leq \bar{l}}(x) & \textup{if\,\,} l=\bar{l}\\
\psi_l(x) & \textup{if\,\,} l>\bar{l}\\
\end{array}
\right., \quad \forall l_1,l_2\in \Z, \quad  \varphi_{[l_1, l_2];\bar{l}}(x):= \sum_{l\in [\max\{l_1,\bar{l}\}, l_2]\cap \Z } \varphi_{l;\bar{l}}(x), 
  \ee
In particular, if the threshold $\bar{l}=0$, we use the following notation, 
\be\label{cutoffwiththreshold100}
\forall k, k_1, k_2 \in \mathbb{Z}_{+},\quad  \varphi_k(x):=\varphi_{k;0}(x), \quad \varphi_{[k_1,k_2]}(x):=  \varphi_{[k_1,k_2];0}(x).
\ee

For a function $F(x)\in L^1$, we use $F^{+}:=F$ and $F^{-}:=\bar{F}$. Moreover, we use both $\widehat{F}(\xi)$ and $\mathcal{F}(F)(\xi)$ to denote the Fourier transform of $F$, which is defined as follows,
\[
\mathcal{F}(F)(\xi)= \int e^{-ix \cdot \xi} F(x) d x.
\]

We use $\mathcal{F}^{-1}(G)$ to denote the inverse Fourier transform of $G(\xi)$.
Moreover, for any $k\in \Z,$ we use  $P_{k}$, $P_{\leq k}$ and $P_{\geq k}$ to denote the projection operators  by the Fourier multipliers $\psi_{k}(\cdot),$ $\psi_{\leq k}(\cdot)$ and $\psi_{\geq k }(\cdot)$ respectively. For convenience in notation, we also use  $f_{k}(x)$ to abbreviate $P_{k} f(x)$.
 
We define the following class of symbol, 
\be\label{symbolclassdefinition}
\mathcal{S}^\infty:=\{m(\xi): m: \R^3\longrightarrow \R, \quad \|\mathcal{F}^{-1}[m](x)\|_{L^1_x}< +\infty\}.
\ee
Moreover, the $\mathcal{S}^\infty$-norm of symbols is defined as follows, 
\be
\| m(\cdot)\|_{\mathcal{S}^\infty} := \|\mathcal{F}^{-1}(m) (x)\|_{L^1_x}.
\ee

\begin{definition}\label{varioureldef}
For any vector $u=(u_1, u_2,u_3)\in \R^3, v\in \R^3/\{0\}$, we define 
\[
 \tilde{v}:=v/|v|,  \quad  \langle u \rangle :=(1+|u|^2)^{1/2}, \quad   \hat{u}:=u/\langle u \rangle, \quad     u_{\bot}:=(u_1, u_2) . 
\] 
In particular, we have $ {\hat{u}}_{\bot}=(u_1, u_2)/\langle u\rangle$.  Moreover, we define the following projection maps, 
\be\label{definitionprojection}
\mathbf{P}:\R^3\longrightarrow \R^2,  \mathbf{P}_i:\R^3\longrightarrow \R, \quad  \mathbf{P}(u)=  u_{\bot}, \quad \mathbf{P}_i(u)=u_i, i\in\{1,2,3\}.
\ee
\end{definition}

 We define the unit vectors of the Cartesian coordinate  system in $\R^3$ as follows, 
\be\label{2020feb18eqn1}
e_1:=(1,0,0), \quad e_2:=(0,1,0), \quad e_3:=(0,0,1).
\ee

Let $v\in \R^3/\{0\}$ be fixed. With the above notation, for any  $u \in \mathbb{R}^3,$ we have an important decomposition of  $u $ with respect to $v$,  which follows from a straightforward computation,
 \be\label{nov24eqn1}
 \begin{split}
&  u= (\tilde{v}\cdot u) \tilde{v}     + \sum_{i=1,2,3}  ((\tilde{v}\times e_i) \cdot u)  (\tilde{v}\times e_i),\\
\end{split}
\ee 
where ``$\cdot$'' denotes the standard inner product and  ``$\times$'' denotes the standard cross product in $\R^3$. 

The above decomposition is motivated from the following fact that the radial   derivative of $\hat{v}$ is much better than the rotational derivatives of $\hat{v}$ when $v$ is large. This fact follows from the following explicit computation, 
\be
\tilde{v}\cdot\nabla_v \hat{v}= \frac{\tilde{v}}{(1+|v|^2)^{3/2}},  \quad   (\tilde{v}\times e_i) \cdot\nabla_v \hat{v}= \frac{ (\tilde{v}\times e_i) }{(1+|v|^2)^{1/2}}, \quad i\in\{1,2,3\}.
\ee

For the rest of this paper, we will use the   decomposition in  \eqref{nov24eqn1}    constantly in the estimate of  derivatives of $\hat{v}$ and the relativistic velocity characteristics $\hat{V}(t)$without further explanation.

   \subsubsection{Main driving force: conservation laws}\label{conservationlawssubsection}
There are three conservation laws available to us, which are   the main driving force of studying the large data problem.   The first conservation law is  due to the transport nature of the Vlasov equation. More precisely, we have
 \be\label{2024nov9eqn111}
   \forall p \in[1,\infty], \qquad \| f(t,x,v)\|_{L^p_{x,v}}= \| f_0(x,v)\|_{L^p_{x,v}}. 
   \ee

   The second conservation law comes from the following  result of direct computations, see \eqref{mainequation},  
\be\label{dec3eqn1}
\frac{d}{dt}\big[\frac{1}{2}(|E|^2+|B|^2) + 4\pi \int_{\R^3} \sqrt{1+|v|^2}f(t,x,v)   d v  \big] = - \nabla_x\cdot\big(  (E\times B) + 4\pi \int_{\R^3} vf(t,x,v) d v  \big).
\ee
After integrating  the above  equality in the whole space $\R^3$, we have 
 \be\label{conservationlaw}
\mathcal{H}(t):= \int_{\R^3}|E(t,x)|^2 + |B(t,x)|^2 d x + 8\pi  \int_{\R^3} \int_{\R^3} \sqrt{1+|v|^2} f(t,x,v)  d x d v = \mathcal{H}(0).
\ee
Therefore,  $L^2$-norm of the electromagnetic fields $(E, B)$  and the   first momentum of the distribution function remain  bounded over time  within the lifespan of solution. 

Lastly, from   the work of Luk-Strain\cite{luk2}[Proposition 2.2],  the conservation law, as established in Lemma \ref{conservationlawlemma}, is satisfied by integrating \eqref{dec3eqn1} over the spacetime region bounded by the hyperplane $t=0$ and the backward light cone emanating from the point $(x,t)$.

 The improved properties of certain electric and magnetic field combinations (see \eqref{2024Dec7eqn21}) are reminiscent of the  ``electric-magnetic decomposition''  central to Christodoulou and Klainerman's seminal work \cite{ChristKlainerman}[Section 7.2] on the global stability of Minkowski spacetime.  
\begin{lemma}[\cite{luk2},Proposition 2.2]\label{conservationlawlemma}
For any $t\in [0,T^{ })$, we have
\be\label{march18eqn31}
\begin{split}
\sup_{x\in \R^3} &  \int_0^t \int_{\mathbb{S}^2} (t-s)^2 K_g^2(s, x+(t-s)\omega, \omega) d \omega ds \\
& + \int_0^t \int_{\R^3} \int_{\mathbb{S}^2} (t-s)^2 \langle v \rangle(1+\hat{v}\cdot \omega) f(s, x+(t-s)\omega, v) d \omega d v ds  \lesssim_{data} 1,\\
\end{split}
\ee
where
\be\label{2024Dec7eqn21}
K_g^2:= |E\cdot \omega|^2 + |B\cdot \omega |^2 +|E-\omega \times B|^2 +|B+\omega \times E|^2. 
\ee
\end{lemma}
 
\subsubsection{Backward characteristics}
 
Recall  \eqref{mainequation}.   If $s\leq t$ ($s\geq t$),    the following system of equations satisfied by the backward (forward)  characteristics, 
\be\label{backward} 
\left\{\begin{array}{rl}
\p_s X(x,v,s,t) &= \widehat{V}(x,v,s,t),\\ 
 \p_s V(x,v, s,t) &= E(s, X(x,v,s,t)) + \widehat{V}(x,v,s,t)\times B(s, X(x,v,s,t))\\
& := K(s,X(x,v,s,t), V(x,v,s,t)),\\ 
X(x,v,t,t)&=x, \quad V(x,v,t,t)=v.\\ 
\end{array}\right. 
\ee
 In subsequent discussions, we will  refer $K(s,X(x,v,s,t), V(x,v,s,t))$ as the acceleration force of characteristics.  Note that, due to the transport nature of the Vlasov equation,    for any $s, t\in [0, T^{  })$, we have
\be\label{conservation}
f(t,x,v) = f(s, X(x,v,s,t), V(x,v,s,t)),
\ee
which also gives us the conservation law for the $L^\infty_{x,v}$-norm of the distribution function $f$.

\subsubsection{Dyadically measuring the velocity characteristics}

To control the  electromagnetic field $(E, B)$ over time, see also Theorem \ref{maintheorem1part1},  we use the classic momentum method. More precisely, we define 
 \be\label{may2eqn1}
\mathfrak{M}_{}(t ):= \int_{\R^3} \int_{\R^3} (1+ |v|)^{N_0/10} f(t,x,v) d v d x,\quad    \overline{\mathfrak{M}}(t) := (1+t)^{(N_0)^3}+\sup_{s\in[0,t]} \mathfrak{M}(s), 
 \ee
 where $N_0=10^{10}$ is defined in \eqref{assumptiononinitialdata}. From the above definition, it's clear that $\overline{\mathfrak{M}}(t) $ is an increasing function with respect to time $t$ and the size of $t$ is much smaller than   $\overline{\mathfrak{M}}(t)$.

Since we will do dyadic decomposition for many variables, for convenience of comparing different objects,   we  define the following dyadic quantity to capture the size of $\overline{\mathfrak{M}}(t)$.
\begin{definition}\label{scaleofv}
For any $\forall t\in [0, T ),$ we define 
\be\label{dec2eqn1}
 M_t:= \inf\{k: k\in \Z, 2^k\geq (  \overline{\mathfrak{M}}(t) )^{1/  (N_0/10-1 )} \}.
\ee
\end{definition}
\begin{remark}\label{scaleofvphilosophical}
We provide a philosophical  interpretation for the important quantity $M_t$  defined above.  Given that we lack any a priori information about the distribution function $f$, we  assume that $f$ is supported in a annulus of size $R$, i.e., $supp(f)\subset\{v: |v|\in [R, 2R]\}$. Our goal is to extract information about $R$ from the moment $\mathfrak{M}_{}(t )$.  Note that, in view of the conservation laws \eqref{2024nov9eqn111} and \eqref{conservationlaw}, the first momentum has both the upper bound and the lower bound, which are independent of time.  Therefore, we have
\[
\mathfrak{M}_{}(t ) \sim R^{N_0/10-1}\int_{\R^3} \int_{\R^3} (1+ |v|)^{ } f(t,x,v) d v d x \sim R^{N_0/10-1}, \quad \Longrightarrow R\sim \big(\mathfrak{M}_{}(t ) \big)^{1/(N_0/10-1)}, 
\]
where $\sim$ is defined in \textbf{Definition} \ref{relationdef}.
 Recall the definition of $M_t$ in \eqref{dec2eqn1}. Roughly speaking, $M_t$ captures the information about the dyadic interval such that $v$ lives. 

\end{remark}
\begin{remark}
  Note that, for any $t_1, t_2\in [0, T^{ }),$ s.t., $t_1\leq t_2,$ we have $M_{t_1}\leq M_{t_2}$ because $ \overline{\mathfrak{M}}(t) $ is an increasing function with respect to $t$. 
\end{remark}

In view of the definition $M_t$ in \eqref{dec2eqn1}, to control the moment $\overline{\mathfrak{M}}(t) $, it suffices to control $M_t$ over time. Since we are working with the initial data with   polynomial decay, we would like to show   that 
the tail part doesn't play an essential role. To this end, we are motivated to define a time dependent majority set  as follows.

\begin{definition}[$t$-majority set]\label{tmajorityset}
For any $t\in [0, T^{  })$,   we define  \textbf{$t$-majority} set  of particles, which initially localize around zero, at time $s\in [0, T^{  })$ as follows,  
\be\label{may10majority}
R_t(s):= \{(X(x,v,s,0),V(x,v,s,0)): |V (x,v,0,0)|+| X(x,v,0,0) |\leq  2^{ M_t/2}  \}.
\ee
In particular, we have
\[
  R_t(0):=\{(x,v): |(x,v)|\leq 2^{ M_t/2} \}. 
\]
Moreover, to measure the $t$-majority set,  we define 
\be\label{may9en21}
\begin{split}
\alpha_t(s,x,v)&:=  \sup_{\tau \in[0,s] } \inf\{ k: k\in \R_+, |   V_{\bot}(x,v,\tau,0)| \leq 2^{k M_t  }   \}  ,\quad \alpha_t(s):=\sup_{(x,v)\in R_t(0)} \alpha_t(s,x,v), \\ 
\beta_t(s, x,v)&:= \sup_{\tau\in[0,s] } \inf\{ k: k\in \R_+, | V(x,v,\tau,0)|\leq 2^{k M_t}  \}, \quad \beta_t(s):=\sup_{(x,v)\in R_t(0)} \beta_t(s,x,v), \\ 
 \alpha_t&:=\alpha_t(t), \quad  \tilde{\alpha}_t=\min\{\alpha_t, (1+2\epsilon) \},\quad  \beta_t:=\beta_t(t), \quad  \tilde{\beta}_t=\min\{\beta_t, (1+2\epsilon) \},
\end{split}
\ee
 where $ {V}_{\bot}$ is defined in \textbf{Definition} \ref{varioureldef} and the velocity characteristics $V(x,v,s,t)$ is defined in \eqref{backward}.  In subsequent discussions, we will refer to the region of characteristics  that start from the $t$-majority set as the majority part and refer the region of characteristics  that start from the  complement  of the  $t$-majority set as the tail part. 
\end{definition}

 \begin{remark}\label{sizeofvphilosophical}

We offer a geometric interpretation of 
$\alpha_t$ and $\beta_t$, as much of this two-paper sequence aims to elucidate these quantities. Broadly speaking, if we assume that we are working with compactly supported initial data, we can identify a ball in $\R^3$
  that encompasses all forward velocity characteristics $V(x,v,s,0)$  originating from the compact support; in this case, 
 $\beta_t$
  represents the radius of this ball. Additionally, we can find another ball in 
$\R^2$
  that covers all projections of the forward velocity characteristics, denoted as  $\mathbf{P}\big(V(x,v,s,0)\big)$; here, 
 $\alpha_t$
  measures the radius of this ball in 
$\R^2$
 . To compare these radii in relation to 
$M_t$
 , we are motivated to adjust the scales in the definitions of 
 $\alpha_t$ and $\beta_t$. A priori, we don't know the size of $\alpha_t$ and $\beta_t$.  It turns out that $\tilde{\alpha}_t$ and $\tilde{\beta}_t$, which are bounded from the  above, are very useful  during the estimates. A posteriori, $\tilde{\alpha}_t$ ($\tilde{\beta}_t$) is same as $\alpha_t$ ($\beta_t$).

 \end{remark}
\begin{remark}

Since $M_t$ is an increasing function, $\forall t\in [0, T^{  }), s\in [0, t], $ from the above definition, we know that 
\be\label{2021dec18eqn1}
 R_s(0)\subset R_t(0), \quad  \alpha_s M_s \leq \alpha_t(s) M_t\leq \alpha_t M_t, \quad \beta_s M_s\leq \beta_t(s) M_t\leq \beta_t M_t. 
\ee
\end{remark}
\begin{remark}
Note that, $\forall s\in[0, t]\subset [0, T), x, v \in \R^3,$ we have 
\[
  X( X(x,v,0,s), V(x,v,0,s),s,0)= x, \quad V( X(x,v,0,s), V(x,v,0,s),s,0)= v. 
\]
Therefore, in view of the above definition,  for any $v\in \R^3$, s.t., either $|  v_{\bot}|\geq 2^{\alpha_t M_t+\epsilon M_t}$ or $| v|\geq 2^{\beta_t M_t+\epsilon M_t}$, we have
\be\label{nov24eqn27}
\forall s\in [0, t], x \in \R^3,  \quad (X(x,v,0,s), V(x,v,0,s))\notin R_t(0). 
\ee
 Thanks to the rapid polynomial decay rate of the initial data in   \eqref{assumptiononinitialdata}, from \eqref{nov24eqn27} and \eqref{may10majority}, the following estimate holds for any $x,v\in \R^3$, s.t., either $| v_{\bot}|\geq 2^{\alpha_t M_t+\epsilon M_t}$ or $| v|\geq 2^{\beta_t M_t+\epsilon M_t}$, 
 \be\label{nov24eqn41}
 \forall s\in [0, t],  \quad \big|f(s,x,v)\big| = \big|f_0( X(x,v,0,s),V(x,v,0,s) )\big|\lesssim 2^{- 2N_0 M_t/5}. 
 \ee
Thanks to the above rapid decay estimate,   the contribution from the case $|  v_{\bot}|\geq 2^{\alpha_t M_t+\epsilon M_t}$ and the case $| v|\geq 2^{\beta_t M_t+\epsilon M_t}$ becomes negilible\footnote{In the compactly supported initial data case, an analogue of the estimate \eqref{nov24eqn41} would simply be $f(s,x,v)=0.$}.  As a result,  it suffices  to consider the case $| v_{\bot}|\leq  2^{\alpha_t M_t+\epsilon M_t}$ and $| v|\leq 2^{\beta_t M_t+\epsilon M_t}$  when dealing with the distribution function  in subsequent discussions.
\end{remark}

\subsection{Main results}\label{mainresultsthispaper}

Our first result provides an $L^\infty_x$-estimate for the electromagnetic field. Although it is quite rough, it demonstrates  the effectiveness of the moment method.  
 \begin{theorem}\label{maintheorem1part1}
 Let $\epsilon=10^{-7}$.For  any  $t\in[0, T^{})$,   we have 
\be\label{maintheoremroughest}
 \| E(s,x)\|_{L^\infty_{[0,t]}L^\infty_x} +   \| B(s,x)\|_{L^\infty_{[0,t]}L^\infty_x}   \lesssim 2^{  (1+2\tilde{\alpha}_t +6\epsilon)M_t}.
\ee
Moreover, for any $x\in \R^3$, s.t., $|  x_{\bot}|\neq 0$, the following pointwise estimate holds, 
\be\label{pointwiseest}
 \sup_{s\in [0, t]}| E(s,x)| +   | B(s,x)|   \lesssim   | x_{\bot}|^{-1/2} 2^{ 5 \epsilon M_t}(2^{M_t}+ 2^{2 \tilde{\alpha}_t M_t }) + |  x_{\bot}|^{-1/4}  2^{5M_t/4+\tilde{\alpha}_t M_t/4+5\epsilon M_t}    +1.
\ee
 \end{theorem}
  \begin{proof}
 See section \ref{proofoftheoremrough}.
 \end{proof}

 Unfortunately, the above estimate is too rough to draw any immediate conclusions.  To better understand the structure of the electromagnetic field, we will examine the angularly localized acceleration force, which is the primary focus of this paper.
  Before presenting our main results for the localized acceleration force, we  first introduce the relevant notation.

Recall  \eqref{electromagnetic}.  As a result of direct computations, after doing dyadic decomposition for the size of $v$, 
we can reduce the Maxwell system into the standard wave equations as follows, 
\be\label{electromagnetic}
\forall U\in \{E,B\}, \quad (\p_t^2-\Delta) U =\sum_{j\in \Z_+} \mathcal{N}^j_U(t, x),
\ee
where
\be\label{sep5eqn1}
\begin{split}
\mathcal{N}^j_E(t, x)&:= - 4\pi \int_{\R^3}  \big(\hat{v} \p_t f(t, x, v)  + \nabla_x f(t, x , v)\big) \varphi_j(v) d v, \\
 \mathcal{N}^j_B(t, x)&:=-4\pi\int_{\R^3} \hat v\times \nabla_x f(t, x, v) \varphi_j(v)d v.\\
 \end{split}
\ee
After substituting the $\p_t f$ by using the the Vlasov equation in  \eqref{mainequation}   and doing integration by parts in $v$, for any $j\in \Z_+$, we have
\be\label{sep5eqn2}
\begin{split}
\mathcal{N}^j_E(t, x)&= 4\pi\big[  \int_{\R^3}\int_{\R^3}  \big(\hat{v} (\hat{v}\cdot \nabla_x) f(t, x, v) -\nabla_x f(t, x , v) \big)\varphi_j(v) \\
&\qquad - f(t,x,v)(E(t,x)+\hat{v}\times B(t,x))\cdot \nabla_v\big(\varphi_j(v) \hat{v}\big) d x d v\big]. \\
\end{split}
\ee
From the Duhamel's formula, $\forall  U\in\{E, B\},$ we have 
\be\label{march14eqn1}
\begin{split}
U(t)  
  & =  \cos(t\d)U_0 + \frac{\sin(t\d)}{\d}  U_1 + \sum_{  j\in\Z_{+}}    U_j(t), \\ 
  U_j(t)
&= \sum_{ \mu\in \{+,-\}}\int_0^t \frac{e^{i\mu(t-s)\d}}{ 2i \mu  \d}     \mathcal{N}^j_U(s) d s,
  \end{split}
\ee
where $U_0= U(t)\big|_{t=0} $ and $U_1:=(\p_t  U )\big|_{t=0}$.

Since the contribution from the initial data is determined and clear, we only need to consider the contribution from the nonlinear part  $U_j(t,x)$.  For the tail-part, $j\geq (1+2\epsilon)M_t,$ we have
\begin{theorem}\label{roughesttailpart}
 Let $\epsilon=10^{-7}$.For  any  $t\in[0, T^{})$,   we have 
 \be\label{2024nov13eqn21}
 \sum_{j\in [(1+2\epsilon)M_t,\infty)\cap \Z} \sum_{U\in \{E, B\}}\| U_j(t,x)\|_{L^\infty}\lesssim  2^{-5M_t}.  
 \ee
\end{theorem}
  \begin{proof}
 See section \ref{proofoftheoremrough}.
 \end{proof}

 Now, we focus on the main case $j\in [0,(1+2\epsilon)M_t]\cap\Z.$ To  better understand the role of the localized angle and the advantages of the \emph{null structure}, as motivated by the Duhamel's formula \eqref{march14eqn1}, we define the following angle-localized electromagnetic field with a general Fourier symbol.

\begin{definition}[Angular localization for the electromagnetic field]\label{angularlocalizationelect}

  For $ U\in \{E, B\} $, any fixed $ \mu\in \{+,-\}, k\in \Z_+,  j\in [0,(1+2\epsilon)M_t]\cap\Z, n\in [-M_t, 2]\cap \Z, \zeta\in \R^3/\{0\}$, and Fourier symbol $\mathfrak{m}(\xi, \zeta)\in \mathcal{S}^\infty$ (see \eqref{symbolclassdefinition}), we define the following localized electromagnetic field,   in which the   frequency   and its angle   relative to   a fixed direction  are localized,  
  \be\label{sep5eqn66}
  \begin{split}
 T_{k,j;n}^{\mu}(\mathfrak{m}, U)(t,x, \zeta) &:= \int_0^t \int_{\R^3} e^{i x\cdot \xi + i \mu(t-s)|\xi | } |\xi|^{-1}  \mathfrak{m}(\xi, \zeta)\widehat{\mathcal{N}^j_U}(s, \xi) \varphi_k(\xi)\varphi_{n;-M_t}\big( \tilde{\xi} + \mu \tilde{\zeta}\big) d \xi d s. \\
T_{k;n}^{\mu}(\mathfrak{m}, U)(t,x, \zeta)    &:=\sum_{j\in\Z_+}  T_{k,j;n}^{\mu}(\mathfrak{m}, U)(t,x, \zeta),\\
  \end{split}
  \ee
  where $\forall u\in \R^3/\{0\}, \tilde{u}:=u/|u|$.
\end{definition}
\begin{remark}
 ``$\zeta$'' in the above definition  will play the role of velocity characteristics $V(x,v,s,t)$ in the   iterative smoothing scheme in Part I \cite{PartI}. We localize the angle between $\xi$ and $\mu \zeta$ to capture the   information of the size of the oscillation phase, see Part I \cite{PartI}. 
\end{remark}

 \subsubsection{Refined decomposition of the localized electromagnetic field}

While the localization in \eqref{sep5eqn66} is sufficient for measuring the smoothing effect, the angle between 
$v$ and $\zeta$ is also crucial, as it quantifies the strength of the double null structure. 

 Roughly speaking, for any $U\in \{E, B\}$,  the symbol of the linear part of  $\mathcal{N}_U^j$ is of size $\tilde{v}\times \tilde{\xi}$, see \eqref{sep5eqn1}. In the case where $v$ and $\xi$  are parallel, which constitutes the  worst case scenario, the symbol vanishes.  This is commonly referred to as null structure. For a deeper understanding of null structure, see Klainerman's classic works \cite{Kla83,Kl85,Kla86} on the null structure of nonlinear wave equations.
    In particular,  the symbol of the linear part of  $\mathcal{N}_E^j+\hat{\zeta}\times \mathcal{N}_B^j$ is of size $(\hat{v}-\hat{\zeta})\times (\tilde{v}\times \tilde{\xi})$, which provides extra smallness when $v$ and $\zeta$ are very close. For this reason, we refer this structure of the acceleration force  as the \textbf{double null structure}. Hence, instead of estimating the electromagnetic field separately, we  estimate  the combination of  the electric field and the magnetic field as in  the acceleration force $K(s,X(x,v,s,t), V(x,v,s,t))$, see \eqref{backward}.

To  fully leverage the advantages of  the double null structure,   we decompose nonlinearity  $\mathcal{N}^j_U(t,x)$ for  $U\in \{E, B\}$ (see \eqref{sep5eqn1}) based on the   sizes of $v$ and  the angle between $v$ and $\zeta$  
  as follows, 
\be\label{sep5eqn10} 
\begin{split}
 \mathcal{N}^j_U(t,x)&=\sum_{i=0,1,2,3,4 } \mathcal{N}^{j,n}_{U;i}(t,x, \zeta)  , \\ 
\mathcal{N}^{j,n}_{E;i}(t,x,\zeta)&:= - 4\pi\big[ \int_{\R^3} \big( \hat{v} \p_t f(t, x, v)   +  \nabla_x f(t, x , v)\big)\varphi_{j,n}^i(v, \zeta) d v\big], \\ 
\mathcal{N}^{j,n}_{B;i}(t,x,\zeta) &:=-4\pi\int_{\R^3} \hat v\times \nabla_x f(t, x, v) \varphi_{j,n}^i(v, \zeta)d v,
\end{split}
\ee
where the cutoff functions $\varphi_{j,n}^i(v, \zeta),i\in\{0,1,2,3,4 \},$ are defined as follows, 
\be\label{sep4eqn6}
\begin{split}
\varphi_{j,n}^0(v, \zeta)& := \varphi_j(v)\psi_{\geq  \vartheta^{\star}_0}(\tilde{\zeta}- \tilde{v})\psi_{ \leq    (1/2+3\iota + 55\epsilon ) M_{t^{\star  }} +2}(v), \\
  \varphi_{j,n}^1(v, \zeta)&:= \varphi_j(v)\psi_{\geq  \vartheta^{\star}_0 }(\tilde{\zeta}- \tilde{v})\psi_{ >    (1/2+3\iota + 55\epsilon ) M_{t^{  \star}}+2 }(v),\\ 
\varphi_{j,n}^2(v, \zeta)&:= \varphi_j(v) \psi_{ [\vartheta^{\star}_1, \vartheta^{\star}_0) }( \tilde{\zeta}- \tilde{v} )  ,  \\
 \varphi_{j,n}^3(v, \zeta)&:= \varphi_j(v)  \psi_{<  \vartheta^{\star}_2   }( \tilde{\zeta} - \tilde{v} )
,\\
  \varphi_{j,n}^4(v, \zeta)&:= \varphi_j(v)  \psi_{ [\vartheta^{\star}_2, \vartheta^{\star}_1)   }( \tilde{\zeta} - \tilde{v} ), \\ 
  \vartheta^{\star}_0 &:=\max\{  n + \epsilon M_{t^\star}, -  10^{-3}  M_{t^\star} \},  \quad   \vartheta^{\star}_1:=  n + \epsilon   M_{t^\star}, \quad \vartheta^{\star}_2:=  n -\epsilon  M_{t^\star}.
\end{split}
\ee

We briefly discuss the motivation for decomposing into five components, as outlined in \eqref{sep4eqn6}. We compare   the angle between $\zeta$ and $v$, which is a free variable, with respect to $2^n$, which measures $\angle(\xi, -\mu\zeta)$ due to the cutoff function in \eqref{sep5eqn66}.  Roughly speaking, there are four cases: (i) either  $\angle(v, \zeta)$ is not too small or the case $n$ is not too small; (ii) $n$ is small, $\angle(v, \zeta)$ and $2^n$ are not comparable and  $\angle(v, \zeta)$ is bigger than  $2^n$; (iii) $n$ is small, $\angle(v, \zeta)$ and $2^n$ are not comparable and  $\angle(v, \zeta)$ is smaller than  $2^n$; (iv) $n$ is small, $\angle(v, \zeta)$ and $2^n$ are almost comparable.

 This intuition aligns closely with the partition of unity in \eqref{sep4eqn6}, aside from a technical issue related to the size of 
$v$, which motivates the definitions of 
 $\varphi_{j,n}^0(v, \zeta)$ and $\varphi_{j,n}^1(v, \zeta)$. We choose to split based on the size of the angle  $\angle(v, \zeta)$
  because this quantity is part of the double null structure. Furthermore, since 
$\zeta$ is fixed, it simplifies measurement with a constant reference point.

  Correspondingly, we decompose the localized electromagnetic field into five parts as follows, 
  \be\label{sep17eqn32}
 \begin{split}
T_{k,j;n}^{\mu}( \mathfrak{m}, U)(t,x, \zeta)&:=\sum_{ i\in\{0,1,2,3,4 \}} T_{k,j;n}^{\mu,i}(  \mathfrak{m}, U)(t,x, \zeta), \\ 
T_{k,j;n}^{\mu,i}( \mathfrak{m}, U)(t,x, \zeta)& :=  \int_0^t \int_{\R^3} e^{i x\cdot \xi + i \mu(t-s)|\xi | } |\xi|^{-1}  \mathfrak{m}(\xi, \zeta)\widehat{ \mathcal{N}^{j,n}_{U;i} }(s, \xi, \zeta) \\
&\qquad \times \varphi_k(\xi)\varphi_{n;-M_t}\big(  \tilde{\xi} + \mu \tilde{\zeta} \big) d \xi d s. 
 \end{split}
 \ee

 Since   the nonlinearities $\mathcal{N}^j_U(t,x)$, where $U \in {E, B}$, include linear terms dependent on the distribution function $f$ that result in the loss of one derivative, we can leverage the smoothing effect to obtain the following elliptic-hyperbolic type decomposition for the localized electromagnetic field, 
\be\label{oct7eqn1}
\begin{split}
  & T_{k,j;n}^{\mu,i}(  \mathfrak{m}, E)(t,x, \zeta) + \hat{\zeta}\times  T_{k,j;n}^{\mu,i}(  \mathfrak{m}, B)(t,x, \zeta)  \\
   & =  Ini_{k,j,n}^{\mu,i}( \mathfrak{m})(t, x, \zeta) + \mathfrak{H}_{k,j;n}^{\mu,i}( \mathfrak{m})(t, x, \zeta )     +  \mathfrak{E}^{\mu, i}_{k,j;n}( \mathfrak{m})(t,  x, \zeta) \mathbf{1}_{i\in\{0,1,2,3\} },
\end{split}
\ee
where $Ini_{k,j,n}^{\mu}( \mathfrak{m})( s, x, \zeta)$ depends only on initial data, which doesn't play a role in the argument. 

The elliptic part is absent for the case $i=4$ in \eqref{oct7eqn1}, because the resulting symbol of doing  normal form transformation  $(\mu |\xi|+ \hat{v}\cdot \xi  )^{-1} $ can be very singular if $v\in supp( \varphi_{j,n}^4(\cdot, \zeta))$, see \eqref{sep4eqn6}. It's not wise to do  normal form transformation for the case $i=4$ at the initial stage.

We  refer to the terms  $\mathfrak{E}^{\mu, i}_{k,j;n}( \mathfrak{m})( t, x, \zeta)$  as the elliptic parts, as   they resemble  the structure of the electric field of the Vlasov-Poisson system, where the electric field satisfies an elliptic equation. Correspondingly, we refer to the terms  $\mathfrak{H}_{k,j;n}^{\mu,i}( \mathfrak{m})( t, x, \zeta)$ as the hyperbolic parts because they solve nonlinear wave equations.

 The elliptic parts  $\mathfrak{E}^{\mu, i}_{k,j;n}( \mathfrak{m})( t, x, \zeta), i\in\{0,1,2,3\}, $ are defined as follows, 
\be\label{2022feb24eqn81} 
\begin{split}
   \mathfrak{E}^{\mu, i}_{k,j;n}( \mathfrak{m})(t, x, \zeta)& :=  \int_{\R^3}\int_{\R^3} e^{i x\cdot \xi} \hat{f}( t, \xi, v)   \mathfrak{m}(\xi, \zeta)  m_i(\xi, v, \zeta) \\
   &\quad\times \psi_k(\xi) \varphi_{j,n}^i(v,\zeta)  \varphi_{n;-M_t}( \tilde{\xi} + \mu \tilde{\zeta})    d\xi d v,\\
 m_i(\xi, v, \zeta)&:= \frac{    4\pi\big( {(\hat{v}-\hat{\zeta} )\times (\hat{v}\times \xi)+\xi(1-|\hat{v}|^2)} \big)   }{   |\xi|(\mu |\xi|+ \hat{v}\cdot \xi  )  },  
\end{split}
\ee 
and   the hyperbolic parts $\mathfrak{H}_{k,j;n}^{\mu,i}( \mathfrak{m})(t, x, \zeta), i \in\{0,1,2,3,4\},  $ are defined as follows,
\be\label{2024nov13hyperbolic}
    \mathfrak{H}_{k,j;n}^{\mu,i}( \mathfrak{m})(t, x, \zeta) := \int_{\R^3} e^{ix\cdot \xi + i \mu t|\xi|} \mathfrak{m}(\xi, \zeta) \mathcal{F}[  \mathfrak{H}_{k,j;n}^{\mu,i}] (t, \xi, \zeta) d \xi,
\ee
where
\be\label{2022feb22eqn81} 
\begin{split}
   \mathcal{F}[  \mathfrak{H}_{k,j;n}^{\mu, a }] (t, x, \zeta)& :=\int_0^t \int_{\R^3}  e^{  - i\mu \tau |\xi| } \mathcal{F}\big((E+\hat{v}\times B)f\big)(\tau, \xi, v)  \cdot \nabla_v \big[  m_{j,n}^a  (v, \zeta, \xi) \big]  \\
&\quad \times   \varphi_{n;-M_t} ( \tilde{\xi} + \mu \tilde{\zeta} )  \varphi_k(\xi) d v  d \tau,\quad a\in\{0,1,2,3\}, \\ 
  m_{j,n}^a  (v , \xi , \zeta)&:=   \big(
 \frac{   4\pi\big( {(\hat{v}-\hat{\zeta} )\times (\hat{v}\times \xi)+\xi(1-|\hat{v}|^2)} \big) }{    |\xi|(\mu |\xi|+ \hat{v}\cdot \xi  )  }+  \frac{\hat{v}}{|\xi|} \big) \varphi_{j,n}^a  (v, \zeta),\quad a\in\{0,1,2,3\}, \\ 
 \mathcal{F}[  \mathfrak{H}_{k,j;n}^{\mu,4  }] (t, x, \zeta)&:=  \int_0^t\int_{\R^3} e^{  - i \mu \tau |\xi | } \frac{\varphi_k(\xi)}{|\xi|} \varphi_{n;-M_t}( \tilde{\xi} + \mu \tilde{\zeta}  ) \\
&\quad \times \big(\widehat{ \mathcal{N}^{j,n}_{E;4} }(\tau, \xi, \zeta) + \hat{\zeta} \times  \widehat{ \mathcal{N}^{j,n}_{B;4} }(\tau, \xi, \zeta)\big) d \xi d \tau , \\ 
\end{split}
\ee
where   $   \mathfrak{H}_{k,j;n}^{\mu,a  }(s, \xi, \zeta)$ with $a\in \{0,1,2,3\}$   have two sources: one source is  from the normal form transformation when $\p_s $ hits   the profile  $g(s,x,v):=f(s,x+s\hat{v},v),$  
and the other source is from the nonlinearity of the electric field $\mathcal{N}^j_E(t, x)$. 
Meanwhile,  $  \mathfrak{H}_{k,j;n}^{\mu,4  }(s, \xi, \zeta)$ only comes from the nonlinearity of the electric field $\mathcal{N}^j_E(t, x)$ because we didn't do normal form transformation at the initial stage.

Our first main results concern the elliptic parts $ \mathfrak{E}^{\mu, i}_{k,j;n}( \mathfrak{m})(s,  x, \zeta)$. More precisely, we have 

\begin{theorem}\label{maintheoremellipitic}
Let  Fourier symbol $m(\xi, \zeta)\in  L^\infty_\zeta \mathcal{S}^\infty$,   $\mu\in\{+,-\},$  $ t^{\star}\in [0, T), \alpha^{\star}:= 2/3+\iota, \iota:=10^{-4},   \zeta \in\R^3/\{0\}, t \in [0,  t^{\ast}]$ be fixed     s.t., $\alpha_t M_{t_{ }}\leq  \alpha^{\star} M_{t^{\star}}, M_t\gg 1$. 
For any fixed   $k\in \Z_{+},  j\in [0,(1+2\epsilon)M_t]\cap\Z,   n\in [-  M_t, 2]\cap \Z$, $i\in\{0,1,2,3\},$  
\begin{enumerate}
   \item[(i)] The following $L^\infty_x$-type rough estimate holds, 
\be\label{2022feb25eqn1}
\begin{split}
    & \big\|\mathfrak{E}^{\mu, i}_{k,j;n}(\mathfrak{m})(t, \cdot, \zeta)\big\|_{L^\infty_x}\\
      &\lesssim \| \mathfrak{m}(\cdot, \zeta)\|_{\mathcal{S}^\infty}  \min\big\{ 2^{ 4{\alpha}^{\star}  M_{t^{\star}}/3 + 5\epsilon M_{t^{\star}} } \big(2^{{\alpha}^{\star}  M_{t^{\star}}/3 } + (2^{M_{t^{\star}} }\frac{| \zeta_{\bot}|}{|\zeta|})^{1/3}\big),    2^{(1-20\epsilon)M_{t^{\star}} } +2^{50\epsilon M_{t^{\star}}} \\
     &\quad \times  \mathbf{1}_{n\geq  (1-2 {\alpha}^{\star})M_{t^{\star}}-40\epsilon M_{t^{\star}} }    \min\{2^{(k+2n)/2+  {\alpha}^{\star}  M_{t^{\star}} },2^{(k+4n)/2+(3{\alpha}^{\star}  -1)   M_{t^{\star}} }\}     \big\}. \\
  \end{split}
\ee
  \item[(ii)]  The following pointwise estimates hold for any $x, \zeta\in \R^3$ s.t., $|  x_{\bot}|\in (0, 2^{M_{t^{\star}}/2}]$,  $|  \zeta_{\bot}|\sim   2^{\gamma_1 M_{t^{\star}}  } , |\zeta|\sim  2^{\gamma_2 M_{t^{\star}}   } ,$  where $\gamma_1\geq   \alpha^{\star} -4\epsilon    $, $\gamma_2 \leq (1-2\epsilon) $, 
\begin{itemize}
  \item[(ii-a)]If $ |  x_{\bot}|\leq  -k-n +2\epsilon M_{t}$, we have 
  \be\label{2024oct27eqn1}
   \big|   \mathfrak{E}^{\mu, i}_{k,j;n} (\mathfrak{m})(t, x,\zeta)\big| \lesssim   \| \mathfrak{m}(\cdot, \zeta)\|_{\mathcal{S}^\infty}  |  x_{\bot}|^{-1} 2^{(\gamma_1-\gamma_2)M_t} 2^{ 5 {\alpha}^{\star} M_{t^{\star}}/6} \| \mathfrak{m}(\cdot, \zeta)\|_{\mathcal{S}^\infty}.\\
  \ee
    \item[(ii-b)]If  
 $ | x_{\bot}|\geq  -k-n +3\epsilon M_{t}/2$, we have 
\be\label{2022feb24eqn1}
\begin{split}
  \big|   \mathfrak{E}^{\mu, i}_{k,j;n} (t, x,\zeta)\big| &\lesssim  \| \mathfrak{m}(\cdot, \zeta)\|_{\mathcal{S}^\infty}   |  x_{\bot}|^{-1/2} 2^{(  {\alpha}^{\star}+10\epsilon)M_{t^{\star}}}, \\
  \big|   \mathfrak{E}^{\mu, i}_{k,j;n} ( t, x,\zeta)\big| &\lesssim \| \mathfrak{m}(\cdot, \zeta)\|_{\mathcal{S}^\infty}   |  x_{\bot}|^{-1/2}2^{  (\gamma_1-\gamma_2)M_{t^{  }}/2}   \big[ 2^{(  {\alpha}^{\star}-10\epsilon)M_{t^{\star}}}   + 
  \mathbf{1}_{n\geq  -\alpha^\star M_{t^{  }}/2  }\\
  &\quad \times  \min\{2^{(k+2n)/2+ 2{\alpha}^{\star}M_{t^{  }}/3}, 2^{(k+4n)/2+{\alpha}^{\star}M_{t^{  }}-(\gamma_1-\gamma_2)M_{t^{  }}/3}\}\big].
  \end{split}
\ee
\end{itemize}
 
\end{enumerate}
\end{theorem}
 \begin{proof}
 See section \ref{proofpartofmaintheorems}.
 \end{proof}

For the hyperbolic parts $\mathfrak{H}_{k,j;n}^{\mu,i}( \mathfrak{m})(s, x, \zeta ), i\in\{0,1,2,3,4\}$, we have 
 \begin{theorem}\label{mainresultsfirstpart}
Let  Fourier symbol $m(\xi, \zeta)\in  L^\infty_\zeta \mathcal{S}^\infty$,   $\mu\in\{+,-\},$  $ t^{\star}\in [0, T), \alpha^{\star}:= 2/3+\iota, \iota:=10^{-4},   \zeta \in\R^3/\{0\}, t \in [0,  t^{\ast}]$ be fixed     s.t., $\alpha_t M_{t_{ }}\leq  \alpha^{\star} M_{t^{\star}}, M_t\gg 1$. 
For any fixed   $k\in \Z_{+},  j\in [0,(1+2\epsilon)M_t]\cap\Z, n\in [-  M_t, 2]\cap \Z$,  
 \begin{enumerate}
 \item[(i)] The following $L^\infty_x$-norm type estimates hold for any $\zeta\in \R^3/\{0\}$,
\be\label{2024oct8eqn1}
\begin{split}
 \sum_{i=0,2}\| \mathfrak{H}_{k,j;n}^{\mu,i}( \mathfrak{m})(t, x, \zeta )\|_{L^\infty_x} 
 &\lesssim \| \mathfrak{m}(\cdot, \zeta)\|_{\mathcal{S}^\infty}   \big[2^{(1-19\epsilon)M_{t^{\star}} } 
+   2^{  128\epsilon M_{t^{\star}} }    \mathbf{1}_{n\geq  \mathfrak{n}_1}\\
&\quad \times \min\{ 2^{(k+2n)/2+ (\alpha^{\star}+3\iota) M_{t^{\star}} }  , 2^{(k+4n)/2 +(1+6\iota)M_{t^{\star}}} \} \big], \\
 \| \mathfrak{H}_{k,j;n}^{\mu,1}( \mathfrak{m})(t, x, \zeta )\|_{L^\infty_x} 
 &\lesssim \| \mathfrak{m}(\cdot, \zeta)\|_{\mathcal{S}^\infty}   \big[2^{(1-19\epsilon)M_{t^{\star}} } 
+     2^{   121\epsilon M_{t^{\star}} }    \mathbf{1}_{n\geq  \mathfrak{n}_1 } \\
  &\quad \times  \min\{   2^{(k+2n)/2 + (2{\alpha}^{\star} -1)M_{t^{\star}}}  , 2^{(k+4n)/2 +(1+6\iota)M_{t^{\star}}} \} \big],\\
  \sum_{i=3,4}\| \mathfrak{H}_{k,j;n}^{\mu,i}( \mathfrak{m})(t, x, \zeta )\|_{L^\infty_x} 
 & \lesssim    \| \mathfrak{m}(\cdot, \zeta )\|_{\mathcal{S}^\infty}   \big[2^{(1-19\epsilon) M_{t^\star} } + 2^{128\epsilon M_{t^\star} } \mathbf{1}_{n\geq \mathfrak{n}_2 } \\
  &\quad \times \min\{2^{(k+2n)/2 +  {\alpha}^{\star} M_{t^\star}}, 2^{(k+4n)/2 +(7/6+5\iota/2)M_{t^\star}}\} \big],  \\
\end{split}
\ee
 where 
 \be
  \mathfrak{n}_1:= -(\alpha^{\star}+3\iota+60\epsilon) M_{t^{\star}}, \quad \mathfrak{n}_2:= - ( (1 +3\iota)/2+40\epsilon) M_{t^{\star}}. 
 \ee
 \item[(ii)] The following pointwise estimates hold for any $x, \zeta\in \R^3,$ s.t., $|x_{\bot}|\in (0, 2^{M_{t^{\star}}/2}]$,  $| \zeta_{\bot}|\sim   2^{\gamma_1 M_{t^{\star}}  } , |\zeta|\sim  2^{\gamma_2 M_{t^{\star}}   } ,$  where $\gamma_1\geq   \alpha^{\star} -4\epsilon    $, $\gamma_2 \leq (1-2\epsilon) $, 
 \begin{itemize}
  \item[ (ii-a) ] If $|  x_{\bot} |\leq 2^{-k-n+2\epsilon M_{t^{\star}}}$, we have 
\be\label{2024oct8eqn2}
\begin{split}
&\sum_{i=0,1,2,3,4} 2^{(\gamma_1-\gamma_2)M_{t^{\star}}}\big|\mathfrak{H}_{k,j;n}^{\mu,i}( \mathfrak{m})(t, x, \zeta ) \big|+\big|\mathbf{P}\big(\mathfrak{H}_{k,j;n}^{\mu,i}( \mathfrak{m})(t, x, \zeta ) \big)\big|\\  
&\lesssim   \| \mathfrak{m}(\cdot, \zeta )\|_{\mathcal{S}^\infty}  \big[\sum_{a\in \mathcal{T}+\mathcal{T}} | x_{\bot}|^{-a}2^{a(\gamma_1-\gamma_2)M_{t^{\star}} } 2^{( \alpha^{\star}   -10\epsilon)M_{t^{\star}}}  \big],\\
\end{split}
\ee
  \item[(ii-b)] If $|  x_{\bot} |\geq 2^{-k-n+\epsilon M_{t^{\star}}}$, we have 
  \be\label{2024oct8eqn5}
\begin{split}
 &\sum_{i=0,1,2 } 2^{(\gamma_1-\gamma_2)M_{t^{\star}}}\big|\mathfrak{H}_{k,j;n}^{\mu,i}( \mathfrak{m})(t, x, \zeta ) \big|+\big|\mathbf{P}\big(\mathfrak{H}_{k,j;n}^{\mu,i}( \mathfrak{m})(t, x, \zeta ) \big)\big|\\
&\lesssim  \| \mathfrak{m}(\cdot, \zeta )\|_{\mathcal{S}^\infty} \Big[   2^{    40\epsilon M_{t^{\star}}}     \big( \min\{ | x_{\bot}|^{-1/2} 2^{7 \alpha^{\star}  M_{t^{\star}} /8}, | x_{\bot}|^{-1} 2^{3 \alpha^{\star}  M_{t^{\star}} /8  } \} \big)  \\
& \quad  +  \sum_{a\in \mathcal{T}} | x_{\bot}|^{-a}2^{a(\gamma_1-\gamma_2)M_{t^{\star}} } \big[2^{( \alpha^{\star}   -10\epsilon)M_{t^{\star}}}+    \mathbf{1}_{n\geq -(\alpha^\star/2 + 30\epsilon)M_{t^\star}}   \\
&\quad   \times  \min\big\{2^{(k+2n)/2+ 2\alpha^{\star}M_{t^\star}/3- (\gamma_1-\gamma_2)M_{t^{\star} }/6},   2^{(k+2n)/2+  \alpha^{\star}M_{t^\star} - (\gamma_1-\gamma_2)M_{t^{\star} }/3 }\big\}  \big] \Big], \\ 
&\\
 &\sum_{i=3,4 } 2^{(\gamma_1-\gamma_2)M_{t^{\star}}}\big|\mathfrak{H}_{k,j;n}^{\mu,i}( \mathfrak{m})(t, x, \zeta ) \big|+\big|\mathbf{P}\big(\mathfrak{H}_{k,j;n}^{\mu,i}( \mathfrak{m})(t, x, \zeta ) \big)\big|\\
 &\lesssim  \| \mathfrak{m}(\cdot, \zeta )\|_{\mathcal{S}^\infty} \Big[   2^{    40\epsilon M_{t^{\star}}}     \big( \min\{ |  x_{\bot}|^{-1/2} 2^{7 \alpha^{\star}  M_{t^{\star}} /8}, |  x_{\bot}|^{-1} 2^{3 \alpha^{\star}  M_{t^{\star}} /8  } \}\big)  \\
& \quad  +  \sum_{a\in \mathcal{T}} | x_{\bot}|^{-a}2^{a(\gamma_1-\gamma_2)M_{t^{\star}} } \big[2^{( \alpha^{\star}   -10\epsilon)M_{t^{\star}}}+    \mathbf{1}_{n\in  \mathcal{N}_{t^{\star} }^1} \min\big\{2^{(k+2n)/2}   \\
&\qquad   \times  2^{ 2\alpha^{\star}M_{t^\star}/3- (\gamma_1-\gamma_2)M_{t^{\star} }/6},   2^{(k+2n)/2+  \alpha^{\star}M_{t^\star} - (\gamma_1-\gamma_2)M_{t^{\star} }/3 }\big\}  \\
& \quad +     \mathbf{1}_{n\in \mathcal{N}_{t^{\star}}^2  } \min\{   2^{(k+3n)/2 +3\alpha^{\star} M_{t^{\star} }/4  },   2^{(k+4n)/2 + \alpha^{\star} M_{t^{\star} }   - (\gamma_1-\gamma_2)M_{t^{\star} }/3}  \} \big]    \Big], \\
\end{split}
\ee
\end{itemize}
where $\mathbf{P}$ denotes the projection operator, see \eqref{definitionprojection}, and the sets $\mathcal{T}$, $\mathcal{N}_{t^{\star} }^1$, and $\mathcal{N}_{t^{\star} }^2$ are defined as follows, 
\be\label{2024oct8eqn8}
\begin{split}
&\mathcal{T}:=\{0,1/8, 1/6, 1/4,1/3, 3/8, 1/2\},\\
&\mathcal{N}_{t^{\star} }^1:=\{n: n\in [ (-  \alpha^{\star} +\gamma_1-\gamma_2) M_{t^{\star} }/2 -30\epsilon M_{t^{\star}}    , 2], \\ 
&\qquad n\notin  [    (\gamma_1-\gamma_2-2\epsilon)M_{t^{\star}},  (\gamma_1-\gamma_2+2\epsilon)M_{t^{\star}}] \}\cap \Z,  \\
&\mathcal{N}_{t^{\star} }^2:=\{n:  n\in  [    (\gamma_1-\gamma_2-2\epsilon)M_{t^{\star}},  (\gamma_1-\gamma_2+2\epsilon)M_{t^{\star}}] \}\cap \Z. 
\end{split}
\ee

 \end{enumerate}
\end{theorem}
 \begin{proof}
 See section \ref{proofpartofmaintheorems}.
 \end{proof}

A detailed discussion of the motivation underlying the above estimates is provided in Part I \cite{PartI}[Section 2.1].

\subsection{Main ingredients of proof of main results}\label{mainingredientsfirstpart}

Since the conservation laws in subsection \ref{conservationlawssubsection} are expressed in the physical space,  we need  representation formulas in the physical space to effectively utilize these conservation laws and  estimate the electromagnetic field. A fundamental tool for this purpose is the classic Kirchhoff's formulas as follows.  
\begin{lemma}[Kirchhoff's formulas]\label{Kirchhoff}
 For any $t\in \R,x\in \R^3$, the following equality holds, 
\be\label{dec20eqn1}
\d^{-1}   \sin( t\d ) h( x)  = \frac{1}{4\pi}   t  \int_{\mathbb{S}^2} h(x+t \theta) d \theta,
\ee
\be\label{march4eqn100}
\d^{-1}   \cos(t\d) h(x)= \frac{1}{4\pi}     \int_{\mathbb{S}^2}  \d^{-1}  h(x+t \theta) d \theta     + \frac{1}{4\pi}     \int_{\mathbb{S}^2} t \theta\cdot \frac{\nabla}{\d} h(x+t \theta) d \theta. 
\ee
\end{lemma}
 \begin{proof}
Note that 
\be\label{march4eqn41}
\begin{split}
 &\int_{\R^3} e^{-i   x\cdot \xi } \int_{\mathbb{S}^2} h(x+t\omega) d \omega d x = \int_{\mathbb{S}^2 } e^{i  t \xi \cdot \omega} \hat{h}(\xi) d\omega\\
 &=  2\pi \hat{h}(\xi)  \int_{0}^{\pi}e^{it|\xi|\cos(\phi)} \sin(\phi) d \phi =  \frac{4\pi \sin(  t |\xi|)}{t|\xi|}  \hat{h}(\xi)  .\\
 \end{split}
\ee
 Hence finishing the proof of the desired formula  \eqref{dec20eqn1}. Our desired equality  \eqref{march4eqn100}  holds after taking derivative with respect to ``$t$'' for  the equality  \eqref{dec20eqn1}.
\end{proof}

Unfortunately, directly applying the above Kirchhoff's formulas to the nonlinearities $\mathcal{N}_U^j(t,x), U \in \{E, B\}$(see \eqref{sep5eqn1}) results in a loss of one derivative. However, there is a smoothing effect available because, in the massive and relativistic setting, the speed of particles is less than that of the electromagnetic field.

The classic Glassey-Strauss decomposition \cite{glassey3}, also known as the S-T decomposition, was introduced to leverage this smoothing effect. In this framework, the T-part depends solely on $f$, while the S-part relies on both $(E, B)$ and $f$. However, this is not a free smoothing effect. Generally, while there is a gain in regularity, there is also a loss of singularity characterized by the factor $(1+\hat{v}\cdot \omega)^{-1}, \omega \in \mathbb{S}^2$. This loss can become significant when $v$ is large and the angle between $v$ and $-\omega$ is small.

To \emph{quantify} the smoothing effect and the advantages derived from the conservation laws, we adopt a similar approach to that in \cite{glassey3} to establish a \textbf{dyadically localized version of the Glassey-Strauss decomposition}. The precise formulation is not critical at this stage; we will elaborate on it further in the upcoming section. 

Furthermore, to leverage the conservation law in \eqref{march18eqn31}, which effectively controls the ``good part'' of the electromagnetic field, we further decompose the S-part into two components based on the following scheme:
\be\label{july1eqn11}
 E(s, x  )  + \hat{v}\times B(s, x  )= EB^1(t,s,x ,\omega, v)+  EB^2(t,s,x  ,\omega, v),
\ee
where
\be\label{july1eqn13}
\begin{split}
EB^1(t,s,x  ,\omega, v)&=  E(s, x  )  -\omega \times B(s, x )\\
&+\big((\omega\cdot  \mathbf{e}_3)(\hat{v}+\omega)_2    (B\cdot \omega )  ,  - (\hat{v}+\omega)_1(\omega\cdot \mathbf{e}_3)  (B\cdot \omega )   ,0 \big)(s,x ),\\
EB^2(t,s,x ,\omega, v)&= \big(  \sum_{i=1,2,3}  -  (\hat{v}+\omega)_2( \omega \times  \mathbf{e}_3)\cdot  \mathbf{e}_i (B\cdot(\omega \times  \mathbf{e}_i))-  (\hat{v}+\omega)_3 B_2,  \\
 & \qquad \sum_{i=1,2,3}  (\hat{v}+\omega)_1( \omega \times  \mathbf{e}_3)\cdot  \mathbf{e}_i (B\cdot(\omega \times  \mathbf{e}_i))+ (\hat{v}+\omega)_3 B_1 ,\\
 &\qquad (\hat{v}_1+ \omega_1) B_2 - (\hat{v}_2 + \omega_2)B_1 \big)(s, x +(t-s)\omega).
  \end{split}
\ee
\begin{remark}\label{decompositionremark}
The key points of the above decomposition are
\begin{enumerate}
\item[(i)] For $EB^1(t,s,x ,\omega, v)$, we can use the conservation law in \eqref{march18eqn31} directly. 
\item[(ii)]  For $EB^2(t,s,x ,\omega, v)$,   which depends solely on the magnetic field $B$, there is an additional gain in smallness related to $\angle(v, -\omega)$. This gain compensates for, but does not completely offset, the loss from the singularity  $(1+\hat{v}\cdot \omega)^{-1}$. 
\item[(ii)]  The projection  $\mathbf{P}\big(EB^2(t,s,x ,\omega, v)\big)$ is better than $ EB^2(t,s,x ,\omega, v)$  because it offers an additional gain in smallness, arising from  $\omega\times \mathbf{e}_3$ and $(\hat{v}+\omega)_3$,  particularly when  $\omega$ and $-v$
 are close to the north and south poles,  which happens when $|v|$ is very large.  The relevance of these poles is due to the fact that  $|  v_{\bot}|$  is generally bounded above by $2^{\alpha_t M_t}$ (with $ \alpha_t\leq 2/3+10^{-4}$), while $|v|$ can  be on the order of $2^{M_t}$. 
\end{enumerate}
\end{remark}

Building on the observations above, we can derive a rough estimate of the electromagnetic field in \eqref{maintheoremroughest} from Theorem \ref{maintheorem1part1} after optimizing various strategies. This includes employing different conservation laws and leveraging cylindrical symmetry by transforming to cylindrical coordinates.

At first glance, it appears that the only tool available to control $EB^2(t,s,x ,\omega, v)$ is the $L^2$-conservation law for the magnetic field. Unfortunately, this approach leads to a very rough estimate.  
As a by-product of proving the estimate \eqref{maintheoremroughest}, we establish a dichotomy for the localized magnetic field. This dichotomy is crucial for estimating the contribution from $EB^2(t,s,x ,\omega, v)$, which depends solely on the magnetic field and does not rely on the conservation law in \eqref{march18eqn31}.

Our key observation in the dichotomy is that, after decomposing the magnetic field into atomic components, \textbf{the 
$L^2_x$-norm and the 
$L^\infty_x$-norm of these atomic pieces cannot be large simultaneously.} Consequently, we can separate the magnetic field into two parts: one part achieves a better 
 $L^2_x$-estimate \footnote{To utilize this, we actually prove the estimate for the geometric mean of the 
$L^2_x$-norm and the 
$L^\infty_x$-norm.} than the conservation law, while the other part provides a better 
$L^\infty_x$-estimate than the rough estimate in Theorem \ref{maintheorem1part1}.

In addition to the localized Glassey-Strauss decomposition, another method for utilizing the smoothing effect is through the Fourier perspective. Generally, this can be derived from the detailed nonlinearity formula in \eqref{sep5eqn1} and Duhamel's formula in \eqref{march14eqn1}. In Fourier space, considering the magnetic field as an example and disregarding absolute constants, the nonlinear contribution is given by
\be\label{2024oct12eqn11}
\int_0^t \int_{\R^3}  \int_{\R^3} e^{ix\cdot \xi + i\mu(t-s)|\xi| - is \hat{v}\cdot \xi}  \frac{\hat{v}\times \xi}{|\xi|   }  \widehat{g}(s, \xi, v) \varphi_j(v) d \xi d v d s,
\ee
where $g(t,x,v)\footnote{  ``$g$''  is typically referred to as the \textit{profile} of the distribution function $f$.} :=f(t,x+t\hat{v}, v)$ and $\p_t g = \textit{nonlinearity}$.

We observe that the phase is oscillating in time $s$. More precisely, we have
\be\label{2024oct12eqn13}
e^{ix\cdot \xi + i\mu(t-s)|\xi| - i s \hat{v}\cdot \xi}  =\frac{\p_s\big( e^{ix\cdot \xi + i\mu(t-s)|\xi| - i s \hat{v}\cdot \xi}  \big)}{-i\mu |\xi | - i \hat{v}\cdot \xi }.
\ee

From the computation above, we observe that doing integration by parts in $s$ in \eqref{2024oct12eqn11} by using the identity \eqref{2024oct12eqn13} results in a gain in regularity, albeit with a loss of singularity of size 
 $(1+ \mu \hat{v}\cdot \tilde{\xi})^{-1},$ where $ \tilde{\xi}:=\xi/|\xi|$.
  This loss can be substantial if 
  $v$ is large and 
$\angle(v, -\mu \xi)$ is very small. Additionally, we observe that the symbol  $\hat{v}\times \xi$
 in \eqref{2024oct12eqn11} functions as a null structure.  

 This procedure is now commonly called the normal form transformation, see the classic works of Germain-Masmoudi-Shatah \cite{GerMasSha09,GerMasSha12}. Further details regarding this procedure and the general Fourier method, both of which have significantly contributed to recent advances in the study of global stability for nonlinear dispersive PDEs, may be found in the  pioneering works of   Ionescu-Pausader \cite{IP1,IP2,IP3,IP}.

 The two methods of exploiting the smoothing effect—the localized Glassey-Strauss decomposition and the normal form transformation—are fundamentally similar but differ technically in their loss functions. One loss function is 
 $(1+ \hat{v}\cdot \omega)^{-1}$ 
 , while the other is 
 $(1+ \mu \hat{v}\cdot \tilde{\xi})^{-1}$
 . Considering the oscillatory integral in \eqref{march4eqn41}, we find that $\omega$
  must be almost parallel to $\xi$
 ; otherwise, the integral is non-stationary. This situation leaves us with two cases: either $\omega$
 and $\xi$
  are aligned or they are oppositely directed. In the essential case where $\tilde{\xi}=-\mu \omega$ 
 , the loss functions are identical. In the other case, where $\tilde{\xi}= \mu \omega$,
  this leads to a minor effect, which we refer to as the error case.

 The strength of the conservation law is measured by the function $1+\hat{v}\cdot \omega$, while the strength of the double null structure,   is assessed on the Fourier side. If we rely on only one method to exploit the smoothing effect, we must contend with an additional error case. Although this case is philosophically simpler, it remains challenging to manage in practice. Therefore, to effectively leverage the double null structure, we will categorize different regions based on the possible sizes of $\angle(v, \tilde{\zeta})$
 and utilize both methods of exploiting the smoothing effect.

In summary, the key elements can be outlined as follows: 
\begin{enumerate}
\item[(i)] Two methods of leveraging the smoothing effect:
\begin{itemize}
 \item[(i-a)] Dyadically  localized Glassy-Strauss   decomposition.
  \item[(i-b)]  Normal form transformation.
\end{itemize}
\item[(ii)] Structure observations:
\begin{itemize}
 \item[(ii-a)] Null structure in both $E$ and $B$.
 \item[(ii-b)] The double null structure for the acceleration force $E+\hat{\zeta}\times B$.
 \item[(ii-c)]   A good decomposition of  the acceleration force in the S-part of the representation formula, see \eqref{july1eqn11}. 
\end{itemize}
\item[(iii)] 
A dichotomy for the localized magnetic field offers improved control compared to the 
$L^2_x$
  conservation law. This will be utilized to manage the contribution from the 
$EB^2$ 
  component of the S-part.

\item[(iv)] $\mathbf{P}_3(B)$ is better than $ B $
  due to the fact that $\mathbf{P}_{3}(\hat{v}\times \xi)$ is better than $\mathbf{P} (\hat{v}\times \xi)$ when  $\tilde{v}$ and $\tilde{\xi}$ are close to the north pole or the south pole.
\end{enumerate}

\subsection{Table of notations and the outline of this paper}\label{outlinefirstpaper}

The universally applicable notations and definitions central to this two-paper series are systematically delineated in Table   \ref{essential_global_notations_}.

\begin{table}[H]
\centering
\resizebox{\columnwidth}{!}{%
\begin{tabular}{ |c|c|c|c| } 
 \hline
 Notation & Definition & Appearance  & Remarks \\ 
 \hline
$N_0$ & $10^{10}$ & Theorem \ref{maintheorem} &  \\ 
  \hline
  $\epsilon$ & $10^{-7}=100/N_0$ & Convention \ref{conventionconst}  & The smallest absolute constant\\ 
  \hline
   $\iota$ & $10^{-4}=1000\epsilon  $ & Convention \ref{conventionconst}  & The next level of small constant\\ 
  \hline
 $\tilde{v}$ & $\tilde{v}:=v/|v|$ & Definition \ref{varioureldef} & The direction of $v=(v_1, v_2,v_3)$\\ 
 \hline
 $\hat{v}$ & $\hat{v}:=v/\sqrt{1+|v|^2} $ & Equation \ref{mainequation} & The relativistic velocity\\ 
 \hline
 $v_{\bot}$ & $v_{\bot}:=(v_1,v_2)$ & Definition \ref{varioureldef} & The first two projection of $v$\\ 
 \hline
  $\mathbf{P}(v)$  &$\mathbf{P}(v):=(v_1,v_2)$ & Definition \ref{varioureldef} &  $\mathbf{P}$ is mainly used for long notations\\ 
 \hline
 $\mathbf{P}_i(v)$  &$\mathbf{P}_i(v):=v_i$, $i\in\{1,2,3\}$ & Definition \ref{varioureldef} & $\mathbf{P}_i$ is mainly used for long notations \\ 
 \hline
 $  X(x,v,s,t ) $ & The space characteristics  & \eqref{backward} & Dependence of $x,v,t$ are often   \\
  $V(x,v,s,t )$& The velocity characteristics   & & dropped since they are fixed\\ 
 \hline
$K(s,X(s), V(s)) $ & The acceleration force  & \eqref{backward} &  \\
  \hline
$ M_t$ & Definition \ref{scaleofv} & \eqref{dec2eqn1} &  Remark \ref{scaleofvphilosophical} \\ 
 \hline
$\alpha_t, \beta_t  $ &  Definition \ref{tmajorityset} & \eqref{may9en21} & Remark \ref{sizeofvphilosophical}\\ 
 \hline
 $\tilde{\alpha}_t  $ & $\min\{\alpha_t, 1+2\epsilon\}$ & Definition \ref{tmajorityset}  & A posteriori, $\tilde{\alpha}_t$ ($\tilde{\beta}_t$) is same  \\ 
  $\tilde{\beta}_t  $ & $\min\{\beta_t, 1+2\epsilon\}$ &    & as $\alpha_t$ ($\beta_t$). It was not clear a priori.\\ 
 \hline
 $\alpha^{\star}$ & $2/3+\iota$ &  Theorem \eqref{maintheoremellipitic} & A posteriori,  it's the upper bound 
 \\
 & & & of  $ \alpha_t$, independent of time.\\ 
 \hline
\end{tabular}%
}
\caption{Essential global notations.}\label{essential_global_notations_}
\end{table}

Moreover, the following table of important notations used in this section, which is also consistent with the table in \cite{PartI}[Section 2],
\begin{table}[H]
\centering
\resizebox{\columnwidth}{!}{%
\begin{tabular}{ |c|c|c|c| } 
 \hline
 Notation & Definition & First appearance  & Remarks \\ 
 \hline
$T_{k,j;n}^{\mu}(\mathfrak{m}, U)(t,x, \zeta)$ & Definition \ref{angularlocalizationelect}  & \eqref{sep5eqn66} & $U\in\{E, B\}$; Angular localization \\
 & & & of the electromagnetic field\\
 \hline
$\varphi_{j,n}^i(v, \zeta)$ & \eqref{sep4eqn6}  & \eqref{sep5eqn10} & Partition of unity based on the size of $v$\\
 & & & and the angle between $v$ and $\zeta$\\
 \hline
 $ \mathfrak{E}^{\mu, i}_{k,j;n}( \mathfrak{m})(t, x, \zeta)$  & \eqref{2022feb24eqn81} & \eqref{oct7eqn1} &  The elliptic parts of the localized acceleration force\\
   \hline
    $\mathfrak{H}_{k,j;n}^{\mu,i}( \mathfrak{m})(t, x, \zeta )$  & \eqref{2022feb22eqn81} & \eqref{oct7eqn1} &  The hyperbolic parts of the localized acceleration force\\
   \hline
   $\gamma_1$ & Theorem \ref{maintheoremellipitic}[(ii)] & \eqref{2024oct27eqn1} & Measuring $|\zeta_{\bot}|$; local notation only \\
    $\gamma_2$ & & & Measuring $|\zeta |$; local notation only\\
       \hline
\end{tabular}%
}
\end{table}

To facilitate comprehensive understanding of the technical content,  we accompany each technical section with a notation table, which is systematically positioned at the conclusion of each section. Readers are advised to consult these tables during detailed analysis to ensure contextual consistency.

 The rest of this paper is organized as follows, 
 \begin{enumerate}
\item[$\bullet$] In section \ref{localizedGSdec}, we prove a   localized version   of the Glassey-Strauss decomposition. 
\item[$\bullet$] In section \ref{roughestelectromagnet},  we use the localized version   of the Glassey-Strauss decomposition obtained in section \ref{localizedGSdec} to prove rough pointwise estimates for the electromagnetic field. 
\item[$\bullet$] In section \ref{linfacceloc}, we prove   $L^\infty_{x}$-type estimates for the localized acceleration force. In particular, a refined structure of the third component of magnetic field $\mathbf{P}_3(B)$ is reveled as a by-product.

\item[$\bullet$] In section \ref{horizonestpotw}, we prove pointwise estimates for    the localized acceleration force, which provides better control if the position $x$ locates far away from the $z$-axis, i.e., $|  x_{\bot} |$ is big. 
\end{enumerate}

\noindent \textbf{Acknowledgment}\qquad 
The author thanks Yu Deng, Pin Yu, and Fan Zheng for helpful comments and suggestions on the early version of draft.  The author acknowledges support from  NSFC-12322110,12141102, 12326602,   and MOST-2020YFA0713003, 2024YFA1015000.

\section{A localized version   of the Glassey-Strauss decomposition}\label{localizedGSdec}

 Since the angle between $v$ and $-\omega$ and the size of $t-s$ will play   essential roles in later argument, to track down precisely the role of the   angle between $v$ and $-\omega$ and the size of $t-s$, we first localize the electromagnetic field and then obtain the corresponding Glassey-Strauss  decomposition for the localized electromagnetic field.

Note that, due to the fact that $t\leq 2^{\epsilon M_t-10}$, see  \eqref{may2eqn1},  we have 
\be\label{july11eqn33}
\mathbf{1}_{[0, t]}(x) = \sum_{m\in [-10M_t, \epsilon M_t]\cap \Z} \varphi_{m;-10M_t}(x) \mathbf{1}_{[0, t]}(x).
\ee
From the above equality  \eqref{july11eqn33}, after      localizing the angle between $v$ and $-\omega$ and the size of $t-s$, we have 
\be\label{july5eqn1}
\begin{split}
E_j(t,x)& = \sum_{  l\in [-j, 2]\cap \Z}    \sum_{m\in [-10M_t, \epsilon M_t]\cap \Z} E_{j,l}^{ m}(t, x)   ,\\
 B_j(t,x) & = \sum_{  l\in [-j,2]\cap \Z}  \sum_{m\in [-10M_t, \epsilon M_t]\cap \Z}  B_{j,l}^{ m}(t, x)   ,\\ 
E_{j,l}^{ m}(t, x)&:=     -   \int_0^t    \int_{\R^3} \int_{\mathbb{S}^2}  (t-s)   \big(  \hat{v} \p_t f(s, x+(t-s)  {\omega}, v)    + \nabla_x f(s, x+(t-s)  {\omega}, v)\big)\\
&\quad \times  \varphi^{ }_{j,  l}(v,    {\omega})  \varphi_{m;-10M_t}(t-s) d {\omega}   d v d s,\\ 
B_{j,l}^{ m}(t, x)&:=  \int_0^t   \int_{\R^3}\int_{\mathbb{S}^2} (t-s) \hat{v}\times \nabla_x f(s, x+(t-s)    {\omega} , v) \varphi_{j, l}^{ }(v,   {\omega})   \varphi_{m;-10M_t}(t-s)  d   {\omega}  d  v d s,
\end{split}
\ee
where the cutoff function $\varphi^{  }_{j,  l}(v,    {\omega}  ) $ is defined as follows, 
\be\label{cutoffwiththre}
\varphi_{j,  l  }^{ }(v, \omega)=\varphi_j(v)\varphi_{l;-j}(\tilde{v}+   {\omega}  )  . 
\ee
 
Building on the work of Glassey and Strauss \cite{glassey3}, we apply the same 
$S-T$ decomposition to the localized magnetic field. The key difference lies in the introduction of a new cutoff function, which necessitates adjustments to the formulas from the Glassey-Strauss decomposition. For clarity, we provide a detailed proof in the following proposition.

\begin{proposition}\label{GSdecomloc}
For any fixed $j\in \Z_+, l\in [-j,  2]\cap \Z, m\in [-10M_t, \epsilon M_t]\cap \Z, \forall U\in \{E, B\}, $ the following Glassey-Strauss decomposition holds for the localized electromagnetic field, 
\be \label{july5eqn31}
 U^m_{j,l}(t,x)= U^m_{S;j,l}(t ,x)  + U^m_{T;j,l}(t ,x)  + U_{0;j,l}^m(t,x) , 
\ee
 \be\label{nov26eqn41}
\begin{split}
U^m_{S;j,l}(t,x)&= \int_0^t  \int_{\R^3} \int_{\mathbb{S}^2} (t-s)  \big(E(s,x+(t-s)\omega) +\hat{v}\times  B(s,x+(t-s)\omega) \big)\\
&\qquad \cdot \nabla_v \big(\frac{  m_U(v, \omega) \varphi_{j, l}(v, \omega)}{1+\hat{v}\cdot \omega} \big)  f(s,x+(t-s)\omega, v)    \varphi_{m;-10M_t}(t-s  ) d \omega d v d s,
\end{split} 
\ee 
\be\label{july11eqn90}
\begin{split}
U^m_{T;j,l}(t,x)&=  \int_0^t \int_{\R^3} \int_{\mathbb{S}^2}  f(s,x+(t-s)\omega, v) \big[ \sum_{ q=1,2,3}  \omega^{m,U }_{j;q }(t-s,v,\omega) \p_{x_q}\varphi_{l;-j}( \tilde{v} +\omega)  \\
&\qquad  +\omega^{m,U}_{j;0 }(t-s,v,\omega) \varphi_{l;-j}( \tilde{v} +\omega)  \big]     d \omega d v ds ,
\end{split}
  \ee 
  \be\label{2022jan1eqn1}
 U^m_{0;j,l}(t,x)= \left\{\begin{array}{ll}
 \displaystyle{\int_{\R^3} \int_{\mathbb{S}^2} t f(0,x+t\omega, v)\frac{ m_U(v, \omega)}{1+\hat{v}\cdot \omega}\varphi_{j,-l}(v,\omega) d\omega d v} & \textup{if\,} t\in supp(\varphi_{m;-10M_t}(\cdot))\\ 
 &\\ 
 0 & \textup{if\,} t\notin supp(\varphi_{m;-10M_t}(\cdot)),
 \end{array}\right. 
\ee 
where the coefficients $m_U(v, \omega), U\in\{E, B\},$ and   the coefficients $\omega^{m, U}_{j ;a}(v,\omega), U\in\{E, B\}, a\in\{0,1,2,3\},$ are defined  as follows, 
\be\label{july9eqn11}
\begin{split}
m_E(v, \omega)&:= \hat{v}+\omega, \quad m_B(v, \omega)= \hat{v}\times \omega, \\
\omega^{m,E}_{j ;0}(t-s,v,\omega)&:= \frac{(1-|\hat{v}|^2) (\hat{v} +  \omega) }{(1+\hat{v}\cdot\omega)^2} \varphi_{j }(v ) \varphi_{m;-10M_t}(t-s  )-   \big(\omega -\frac{(\omega +\hat{v} )\hat{v}\cdot\omega }{1+\hat{v}\cdot \omega}\big)\\
&\quad \times (t-s) \p_x \varphi_{m;-10M_t}(t-s  )  \varphi_{j }(v ), \\ 
 \omega^{m,E}_{j;q }(t-s,v,\omega) & = -  \big(  (\delta_{1q}, \delta_{2q}, \delta_{3q})  + w_q \hat{v} -   \frac{(\omega +\hat{v} )(\omega_q + \hat{v}_q) }{1+\hat{v}\cdot \omega}\big)\varphi_{j }(v ) \varphi_{m;-10M_t}(t-s  ),\\
 \omega^{m,B}_{j ;0}(t-s,v,\omega)&:= \frac{(1-|\hat{v}|^2) (\hat{v} \times   \omega) }{(1+\hat{v}\cdot\omega)^2}\varphi_{j }(v ) \varphi_{m;-10M_t}(t-s  )- \frac{\hat{v}\times \omega}{1+\hat{v}\cdot \omega} \\
 &\qquad \times  (t-s) \p_x \varphi_{m;-10M_t}(t-s  )  \varphi_{j }(v ), \\
   \omega^{m,B}_{j;q }(t-s,v,\omega):&= - \big( \hat{v}\times(\delta_{1p}, \delta_{2p}, \delta_{3p})   - \frac{\hat{v}\times \omega}{1+\hat{v}\cdot \omega}\big(\hat{v}_p +\omega_p) \big)\varphi_{j }(v ) \varphi_{m;-10M_t}(t-s  ),
\end{split}
\ee
where $q\in \{1,2,3\}.$
\end{proposition}
\begin{proof}
We first prove the desired decomposition \eqref{july5eqn31} for the magnetic field. Recall \eqref{july5eqn1}, we  have 
\[
B^m_{j,l}(t,x)=  \int_{\R^3} \int_{|y-x|\leq t }  \varphi_{m;-10M_t}  (|y-x|) \frac{\varphi_{j, l}(v, \omega)}{|y-x|} \hat{v}\times \nabla_x f(t-|y-x|, y,  v)   dy  d  v,
\]
where
\[
\omega:= \frac{y-x}{|y-x|}. 
\] 
As in the \cite{glassey3}, we have 
\be\label{july9eqn14}
\forall i\in \{1,2,3\}, \quad \p_i = \omega_i(1+\hat{v}\cdot \omega)^{-1}S + \big(\delta_{ip}-\frac{\omega_i\hat{v}_p}{1+\hat{v}\cdot \omega}\big) T_p, 
\ee
where 
\[
S:=\p_t + \hat{v}\cdot \nabla_x, \quad\forall p\in \{1,2,3\}, \quad  T_p:= -\omega_p \p_t + \p_{x_p}. 
\]

By using the above equalities, as a typical example, the following decomposition holds for the first component of the vector $B^m_{j,l}(t,x)$, 
\be\label{july9eqn20}
\begin{split}
\big( B^m_{j,l}(t,x)\big)_1 &=  \int_{\R^3} \int_{|y-x|\leq t } \frac{    \varphi_{m;-10M_t}  (|y-x|)  }{|y-x|} (\hat{v}_2\p_{x_3}-\hat{v}_3 \p_{x_2}) \\
&\qquad \times f(t-|y-x|, y,  v) \varphi_{j, l}(v, \omega)  dy  d  v \\
&= I_{j,l}^{1;1}(t,x) + I_{j,l}^{1;2}(t,x),\\
\end{split}
\ee
where 
\be
\begin{split}
I_{j,l}^{1;1}(t,x)&=  \int_{\R^3} \int_{|y-x|\leq t } \frac{   \varphi_{m;-10M_t}  (|y-x|)}{|y-x|} \frac{\hat{v}_2\omega_3-\hat{v}_3\omega_2}{1+\hat{v}\cdot \omega}  Sf(t-|y-x|, y,  v) \varphi_{j, l}(v, \omega)  dy  d  v \\
&= \int_{\R^3} \int_{|y-x|\leq t } (E(t-|y-x|, y )+\hat{v}\times B(t-|y-x|, y ))\cdot \nabla_v\big(\frac{(\hat{v}_2\omega_3-\hat{v}_3\omega_2)\varphi_{j, l}(v, \omega)}{1+\hat{v}\cdot \omega} \big)  \\
&\qquad \times \frac{   \varphi_{m;-10M_t}  (|y-x|)}{|y-x|}  f(t-|y-x|, y,  v)dy  d  v,\\
I_{j,l}^{1;2}(t,x)&= \sum_{p=1,2,3} \int_{\R^3} \int_{|y-x|\leq t } \frac{  \varphi_{m;-10M_t}  (|y-x|)}{|y-x|} c_p T_p  f(t-|y-x|, y,  v) \varphi_{j, l}(v, \omega)  dy  d  v.
\end{split}
\ee
The coefficients $c_p, p\in\{1,2,3\},$ appeared above are given as follows,  
\[
c_p:= (\hat{v}_2\delta_{3p}-\hat{v}_3\delta_{2p})   - \frac{\hat{v}_2\omega_3-\hat{v}_3\omega_2}{1+\hat{v}\cdot \omega}\hat{v}_p.
\]
Now, we focus on the $I_{j,l}^{1;2}(t,x)$. By the divergence theorem, we have
\be\label{july9eqn22}
\begin{split}
I_{j,l}^{1;2}(t,x)&= - \sum_{p=1,2,3}\int_{\R^3} \int_{|y-x|\leq t } \p_{y_p} \big( \frac{ c_p \varphi_{j, l}(v, \omega)  \varphi_{m;-10M_t}  (|y-x|)}{|y-x|}  \big)   f(t-|y-x|, y,  v)  dy  d  v \\
&\quad + \int_{\R^3} \int_{|y-x|=t }  \frac{\omega_p c_p}{|y-x|}      f(t-|y-x|, y,  v) \varphi_{j, l}(v, \omega)    \varphi_{m;-10M_t}  (|y-x|) dy  d  v,
\end{split}
\ee
which  splits naturally into the $T$-part and $0$-part (i.e., initial data part).  As a result of direct computations, we have
\[
\begin{split}
\sum_{p=1,2,3} \p_{y_p} \big(  {|y-x|^{-1}} c_p \big) & =-|y-x|^{-2}(\hat{v}_2\omega_3 - \hat{v}_3\omega_2)\frac{1-|\hat{v}|^2}{1+(\hat{v}\cdot\omega)^2},\\
\sum_{p=1,2,3} c_p \p_{y_p}\big( \varphi_{l;-j}(\tilde{v}+\omega)\big) & = \sum_{p=1,2,3} \p_{x_q}\varphi_{l;-j}(\tilde{v}+\omega)\big(\frac{c_q -c_p\omega_p\omega_q}{|y-x|}   \big) =  \sum_{q=1,2,3}\frac{\p_{x_q}\varphi_{l;-j}( \tilde{v}+\omega)}{|y-x|} \\
&\qquad \times \big[(\hat{v}_2\delta_{3q}-\hat{v}_3\delta_{2q})   - \frac{\hat{v}_2\omega_3-\hat{v}_3\omega_2}{1+\hat{v}\cdot \omega}\big(\hat{v}_q +\omega_q)\big], \\
\sum_{p=1,2,3} c_p\p_{y_p}( \varphi_{m;-10M_t}(|y-x|))&= \sum_{p=1,2,3} \p_x  \varphi_{m;-10M_t}(|y-x|)c_p\omega_p\\
&= \p_x \varphi_{m;-10M_t}(|y-x|)\frac{\hat{v}_2 \omega_3-\hat{v}_3\omega_2}{1+\hat{v}\cdot \omega}.  
\end{split}
\]
From the results of the  above computations,  after changing coordinates, our desired decomposition in  \eqref{july5eqn31}  holds for the magnetic field.

We can do similar procedures for the electric field. Note that, from the   equalities in  \eqref{july9eqn14}, $\forall i\in\{1,2,3\}$, we have
\be\label{july9eqn21}
\hat{v}_i\p_t + \p_{x_i} = \hat{v}_i S +\p_{x_i}- \hat{v}_i \hat{v}_q\p_{x_q}=    \frac{\hat{v}_i + \omega_i  }{1+\hat{v}\cdot \omega }S +   d_p  T_p, \quad d_p:=  \big(\delta_{ip}-\frac{(\omega_i+\hat{v}_i)\hat{v}_p}{1+\hat{v}\cdot \omega}\big).
\ee
Note that 
\[
\p_{y_p}\big(|y-x|^{-1}d_p \big)= - |y-x|^{-2}\omega_p d_p + |y-x|^{-1}\p_{y_p}  d_p= - |y-x|^{-2} \frac{\hat{v}_i+\omega_i  }{(1+\hat{v}\cdot \omega)^2}\big(1 -|\hat{v}^2|   \big).
\]
Moreover, direct computations conclude that
\[
\begin{split}
\sum_{p=1,2,3} d_p \p_{y_p}\big( \varphi_{l;-j}( \tilde{v}+\omega)\big) & = \sum_{p,q=1,2,3}  \frac{ \p_{x_q}\varphi_{l;-j}(\tilde{v}+\omega) }{  |y-x|} \big( {d_q -d_p\omega_p\omega_q}  \big)\\
&=\sum_{q=1,2,3}\frac{ \p_{x_q}\varphi_{l;-j}(\tilde{v}+\omega) }{  |y-x|}  \big( \delta_{iq}  + w_q \hat{v}_i-   \frac{(\omega_i+\hat{v}_i)(\omega_q + \hat{v}_q) }{1+\hat{v}\cdot \omega}\big),\\
 \sum_{p=1,2,3} d_p\p_{y_p}( \varphi_{m;-10M_t} (|y-x|))&=\sum_{p=1,2,3}  \p_x  \varphi_{m;-10M_t}(|y-x|)d_p\omega_p\\
 &= \p_x  \varphi_{m;-10M_t}(|y-x|)\big(\omega_{i }-\frac{(\omega_i+\hat{v}_i)\hat{v}\cdot\omega }{1+\hat{v}\cdot \omega}\big).  
\end{split}
\]
From the above results of computations, after using the decomposition    \eqref{july9eqn21}  and redoing steps in  \eqref{july9eqn20}  and  \eqref{july9eqn22}, our desired decomposition in  \eqref{july5eqn31}   also holds for the localize electric field. Hence finishing the proof.

\end{proof}

To exploit the conservation law in \eqref{march18eqn31}, which provides a good control for the ``good part'' of electromagnetic field, we decompose the $S$ part of the electromagnetic field further into two parts by using  the decomposition of $E+\hat{v}\times B$ in  \eqref{july1eqn11}  as follows, 
\be\label{july9eqn41}
E^m_{S;j,l}(t,  x)=\sum_{i=1,2 }E^{m;i}_{S;j,l}(t, x), \quad B^m_{S;j,l}(t,  x)=\sum_{i=1,2 }B^{m;i}_{S;j,l}(t, x),
\ee
where $\forall\,U\in \{E, B\}$, 
\be\label{july9eqn31}
\begin{split}
U^{m;i}_{S;j,l}(t, x) &=  \int_0^t \int_{\R^3} \int_{\mathbb{S}^2}  EB^i(t,s,x ,\omega, v) \cdot \nabla_v \big(\frac{ m_{U}(v, \omega) \varphi_{j, l}(v, \omega)}{1+\hat{v}\cdot \omega} \big) \\
&\qquad \times (t-s)  f(s,x+(t-s)\omega, v)  \varphi_{m;-10M_t }(t-s) d \omega d v d s.\\
\end{split}
\ee

Moreover, recall  \eqref{july5eqn1},  we can define the  frequency  localized version of $U^m_{j,l}(t,x), U\in\{E, B\}$, and their corresponding G-S decomposition in Lemma \ref{GSdecomloc} and the decomposition of $S$ part in  \eqref{july9eqn41}  as follows,  
\be\label{july5eqn60}
\begin{split}
U^m_{k;j,l}(t, x)&:= \int_{\R^3} K_k(y)U^m_{j,l}(t, x-y) d y,\quad \forall \ast\in\{S,T,0\}, \\
 U^m_{k;\ast;j,l}(t,s,x)&:= \int_{\R^3} K_k(y) U^m_{\ast;j,l}(t,s,x-y) d y,\\   
 \forall i \in \{1,2 \}, \quad U^{m;i}_{k;\ast;j,l}(t,s,x)&=\int_{\R^3} K_k(y) U^{m;i}_{\ast;j,l}(t,s,x-y) d y,\\
\end{split}
\ee
where the kernel $ K_k(y)$ is defined as follows, 
\be
  K_k(y) :=\int_{\R^3} e^{iy\cdot\xi} \varphi_k(\xi) d \xi. 
\ee

The essential notations employed in this section are systematically detailed in Table \ref{tablesection2}.
\begin{table}[H]
\centering
\resizebox{\columnwidth}{!}{%
\begin{tabular}{ |c|c|c|c| } 
 \hline
 Notation & Definition    & Remarks \\ 
 \hline
 $U_{j,l}^{ m}(t, x)$, $U\in\{E, B\}$  & \eqref{july5eqn1} & The dyadically localized electromagnetic field  \\
\hline
$U^m_{k;j,l}(t, x)$, $U\in\{E, B\}$ & \eqref{july5eqn60} & Frequency localization of $U_{j,l}^{ m}(t, x)$ \\
\hline
$U^m_{S;j,l}(t ,x)$, $U\in \{E,B\}$ & Proposition \ref{GSdecomloc},  \eqref{nov26eqn41} & The $S$-part of G-S decomposition; depends on both $f$\\ 
& &  and the electromagnetic field \\  
\hline
$U^m_{T;j,l}(t ,x)$, $U\in \{E,B\}$ & Proposition \ref{GSdecomloc},  \eqref{july11eqn90} & The $T$-part of G-S decomposition; depends only on $f$ \\  
\hline
$U^{m;i}_{S;j,l}(t, x), i\in\{1,2\}$, & \eqref{july9eqn31} & Refined decomposition of the $S$-part; \\
$U\in \{E, B\}$ & & the decomposition \eqref{july1eqn11} is applied.\\
\hline
$ U^m_{k;\ast;j,l}(t,s,x)$, $\ast\in\{T, S\}$& \eqref{july5eqn60} & Frequency localization of $ U^m_{ \ast;j,l}(t,s,x)$  \\
\hline
$m_U(v, \omega)$, $U\in\{E,B\}$& \eqref{july9eqn11}  & coefficients appeared in the S-part and the initial data part  \\
\hline
$\omega^{m,U}_{j ;a}(t-s,v,\omega)$, $U\in\{E, B\},$  &\eqref{july9eqn11}  &coefficients appeared in the T-part \\
$ a\in\{0,1,2,3\}$& & \\
\hline
\end{tabular}%
}
\caption{Essential notations in section \ref{localizedGSdec}.}\label{tablesection2}
\end{table}

\section{ Rough pointwise estimates   for the electromagnetic field}\label{roughestelectromagnet}
The primary objective of this section is to establish Theorems \ref{maintheorem1part1} and \ref{roughesttailpart}. These theorems demonstrate that the high-momentum component of the distribution function governs the $L^\infty_x$ norm of the electromagnetic field, while the contribution from the tail component is insignificant.

\subsection{The tail part: the large  velocity case, i.e., $j\geq (1+2\epsilon)M_t$. }

We first consider the case $j\geq (1+2\epsilon)M_t,$ where $ \epsilon =10^{-7}=100/N_0$ was defined in section \ref{notationsubsection}. Note that,   the following $L^1_{x,v}$-estimate of the localized density holds from  the definition  \eqref{may2eqn1},  
\be\label{nov26eqn21}
\begin{split}
\sum_{j\geq (1+2\epsilon) M_t }\big| \int_{\R^3}\int_{\R^3} f(t,x,v) \psi_{j}(v) d x d v \big|& \lesssim \sum_{j\geq (1+2\epsilon)  M_t} 2^{- N_0 j/10 } \mathfrak{M}_{}(t )\\
& \lesssim 2^{-N_0(1+2\epsilon)M_t/10 + (N_0/10-1)M_t}\lesssim 2^{-100M_t},\\
\end{split}
\ee
which is very small and much better than the control provided by the conservation law of the first momentum.  

From the above discussion, intuitively speaking,  the contribution from the tail part, i.e., the case $j\geq (1+2\epsilon)M_t$, is negligible. To rigorously prove this point, in this subsection, we prove the following lemma.

\begin{lemma}\label{erroresti}
   For any   $t\in[0, T^{})$,   $j\in [ (1+2\epsilon)M_t, \infty)\cap \Z_+, m\in [-10 M_t, \epsilon M_t]\cap\Z,  l\in [-j,2]\cap \Z$, $U\in \{E, B\},$  we have 
\be\label{nov26eqn31}
|U^m_{j,l}(t,x)|\lesssim 2^{-9j}\big[ 1 + 2^m\big(\| E\|_{L^\infty_{[0,t]} L^\infty_x}+ \| B\|_{L^\infty_{[0,t]} L^\infty_x} \big)\big]. 
\ee
\end{lemma}
\begin{proof}
Recall the G-S decomposition in  \eqref{july5eqn31}.  Recall  \eqref{july11eqn90}  and the detailed formulas of $\omega^{m, U}_{j;a}(\cdot, \cdot, \cdot), U\in \{E, B\},a\in\{0,1,2,3\},$ in   \eqref{july9eqn11}. The following estimate holds for the T-part, 
\be\label{nov28eqn31}
|U^m_{T;j,l}(t,x)|\lesssim 2^{j} \int_{0}^t \int_{\mathbb{S}^2}\int_{\R^3} f(s, x+(t-s)\omega, v ) \varphi_j(v) \varphi_{m;-10M_t } (t-s) d v d \omega d s .
\ee
Similarly, recall \eqref{nov26eqn41}, after putting the electromagnetic field in $L^\infty_x$, we have 
\be\label{nov28eqn32}
\begin{split}
|U^m_{S;j,l}(t,x)|&\lesssim  2^{m+j}\big(\| E\|_{L^\infty_{[0,t]} L^\infty_x}+ \| B\|_{L^\infty_{[0,t]} L^\infty_x} \big)\\
&\quad \times \big[  \int_{0}^t \int_{\mathbb{S}^2}\int_{\R^3} f(s, x+(t-s)\omega, v ) \varphi_j(v) \varphi_{m;-10M_t }(t-s) d v d \omega d s\big]. \\
\end{split}
\ee

Now, we focus on the estimate of the common term in  \eqref{nov28eqn31}  and  \eqref{nov28eqn32}. 
Note that, after doing dyadic decomposition for the distribution function, we have
\be\label{nov28eqn33}
\int_0^t  \int_{\mathbb{S}^2} \int_{\R^3} f(s, x+(t-s)\omega, v ) \varphi_j(v) \varphi_{m;-10M_t }(t-s) d v d \omega d s \lesssim \sum_{k\in \Z_+} \sum_{m'\in [-k,\epsilon M_t]\cap \Z} J_{k,m'}(t, x),
\ee
where
\be 
 J_{k,m'}(t, x): = \int_0^t   \int_{\R^3} \int_{\mathbb{S}^2}    P_k(f)(s, x+(t-s)\omega, v ) \varphi_j(v)\varphi_{m';-k}(t-s) \varphi_{m;-10M_t }(t-s)d \omega   d v  d s .
\ee 

By doing integration by parts in $\xi$ many times, the following estimate holds for the kernel associated with  the Fourier multiplier operator $P_k$, 
\be\label{nov26eqn51}
\big| \int_{\R^3}e^{iy \cdot \xi}\varphi_{k}(\xi) d\xi  \big|\lesssim 2^{3k} (1+2^k|y|)^{-N_0^3}. 
\ee

From the above estimate of the kernel,  the definition of the high momentum of the distribution function in \eqref{may2eqn1}  and the volume of support of $v$,  we can rule out the case $k\leq 10j$ or  $ m'+k/2\leq \epsilon M_t$ as follows, 
\be\label{nov28eqn51}
\begin{split}
\sum_{\begin{subarray}{c}
k\in \Z_+, m\in [-k, \epsilon M_t]\cap \Z\\ 
k\leq 10j \textup{or\,} m'+k/2\leq \epsilon M_t
\end{subarray}}|  J_{k,m'}(t, x)|&\lesssim \sum_{\begin{subarray}{c}
k\in \Z_+, m\in [-k, \epsilon M_t]\cap \Z\\ 
k\leq 10j \textup{or\,} m'+k/2\leq \epsilon M_t
\end{subarray}} 2^{m'} \min\{2^{3k-nj+(n-1)M_t}, 2^{3j+3\epsilon M_t}\} \\
& \lesssim 2^{-10j}.\\
\end{split}
\ee

Now, we focus on the case $k\geq 10j, m'+k/2\geq \epsilon M_t.$ Note that, in terms of the profile $g(s,x,v):=f(s,x+s\hat{v}, v)$, the following equality holds on the Fourier side, 
\be\label{nov28eqn55}
\begin{split}
J_{k,m'}(t,x) &=\sum_{\star\in\{ess, err\}} J^{\star}_{k,m'}(t, x), \\
 J^{\star}_{k,m'}(t, x)  &:= \int_0^t \int_{\R^3} \int_{\R^3} \int_{\mathbb{S}^2} e^{ix\cdot \xi +i (t-s)\omega \cdot \xi - is \hat{v}\cdot \xi} \hat{g}(s,\xi, v) \varphi^{\star}(\omega, \xi) \varphi_{m';-k}(t-s)   \\
&\qquad  \times \varphi_{m;-10M_t }(t-s)  \varphi_j(v) \varphi_k(\xi) d \omega d v   d\xi  d s, \\
  \varphi^{ess}(\omega, \xi)&:= \psi_{\leq -(1-\epsilon)k/2-m'/2}(\omega\times \xi/|\xi|), \quad \varphi^{err}(\omega, \xi)= 1-  \varphi^{ess}(\omega, \xi). 
\end{split}
\ee
Note that,  each time we do integration by parts in $\omega$, we gain at least $2^{-\epsilon k}$. After doing integration by parts many times in $\omega$, the following estimate holds for the error type term, 
\be\label{nov28eqn56}
\begin{split}
 \sum_{\begin{subarray}{c}
 k\in [10j, \infty)\cap  \Z_+\\
   m'\in [-k, \epsilon M_t]\cap \Z\\
   m'+k/2\geq \epsilon M_t
 \end{subarray}} |J^{ err }_{k,m'}(t, x) |  & \lesssim  \sum_{ k\in [10j, \infty)\cap  \Z_+ } 2^{-100k} \int_0^t \int_{\R^3} \int_{\R^3} \int_{\mathbb{S}^2}    f(s, x-y+(t-s)\omega, v ) \\
&\qquad \times (1+2^{(1+\epsilon)k/2} |y|)^{-100} \varphi_j(v)  d \omega d v   dy d s\\
&\lesssim \sum_{ k\in [10j, \infty)\cap  \Z_+ } 2^{-100k -nj+(n-1)M_t}\lesssim 2^{-10j}. 
\end{split}
\ee
For the essential part, we observe that 
\[
|\omega \cdot \xi +\hat{v}\cdot \xi|\geq 2^{k-2j-10} - \min_{\mu\in \{+, -\}} |\omega\cdot \xi + \mu |\xi||\geq 2^{k-2j-10} - 2^{-m'-k+2\epsilon k +k}\gtrsim 2^{k-2j}.
\]

From the above estimate of the phase, after doing integration by parts in $s$ once,   using the equation satisfied by the Vlasov equation in  \eqref{mainequation}, and doing integration by parts in $v$, we have
\be\label{nov28eqn57}
J^{ess}_{k,m'}(t, x) = J^{ess;0}_{k,m'}(t, x) + J^{ess;1}_{k,m'}(t, x),
\ee 
where
\be\label{nov28eqn2}
\begin{split}
J^{ess;0}_{k,m'}(t, x)& :=   \int_{\R^3} \int_{\R^3} \int_{\mathbb{S}^2} e^{ix\cdot \xi +i (t-s)\omega \cdot \xi  } \frac{\varphi^{ess}(\omega, \xi) }{i(\omega+\hat{v})\cdot \xi} \hat{f}(s,\xi, v) \\
&\qquad \times \varphi_k(\xi) \varphi_j(v)\varphi_{m';-k}(t-s)  \varphi_{m;-10M_t }(t-s) d \omega d v    d\xi   \big|_{s=0}^t \\
&+   \int_0^t   \int_{\R^3} \int_{\R^3} \int_{\mathbb{S}^2} e^{ix\cdot \xi +i (t-s)\omega \cdot \xi  } \frac{\varphi^{ess}(\omega, \xi) }{i(\omega+\hat{v})\cdot \xi} \hat{f}(s,\xi, v) \varphi_k(\xi) \varphi_j(v)\\
&\qquad \times  \p_s\big(  \varphi_{m';-k}(t-s) \varphi_{m;-10M_t }(t-s) \big) d \omega  d v    d\xi   ds, \\ 
J^{ess;1}_{k,m'}(t, x)&:=   \int_0^t \int_{\R^3} \int_{\R^3} \int_{\mathbb{S}^2} e^{ix\cdot \xi +i (t-s)\omega \cdot \xi  }  \mathcal{F}[(E+\hat{v}\times B) f]\cdot \nabla_v \big[  \frac{\varphi^{ess}(\omega, \xi) }{i(\omega+\hat{v})\cdot \xi}  \varphi_k(\xi) \varphi_j(v)\big]\\
 & \qquad \times  \varphi_{m';-k}(t-s)  \varphi_{m;-10M_t }(t-s) d \omega d v   d\xi   d s.\\
\end{split}
\ee
 
In terms of kernel, we have
\be\label{nov28eqn3}
\begin{split}
J^{ess;0}_{k,m'}(t, x)&=  \int_{\R^3} \int_{\R^3} \int_{\mathbb{S}^2} \mathcal{K}_{k,j}^{ess }(y, \omega, v) f(s, x-y+(t-s)\omega, v) \varphi_{m';-k}(t-s)\\
&\qquad \times \varphi_{m;-10M_t }(t-s)   d \omega d y d v \big|_{s=0}^t  \\
& + \int_0^t  \int_{\R^3} \int_{\R^3} \int_{\mathbb{S}^2} \mathcal{K}_{k,j}^{ess }(y, \omega, v) f(s, x-y+(t-s)\omega, v)\\ 
&\qquad \times  \p_s\big( \varphi_{m';-k}(t-s)\varphi_{m;-10M_t }(t-s)  \big)  d \omega d y d v  d s,\\
J^{ess;1}_{k,m'}(t, x)&  =  \int_0^t  \int_{\R^3} \int_{\R^3} \int_{\mathbb{S}^2}     \nabla_v\mathcal{K}_{k,j}^{ess }(y, \omega, v)\cdot \big(E(s, x-y+(t-s)\omega) \\
&\quad +\hat{v}\times B(s, x-y+(t-s)\omega) \big)  f(s, x-y+(t-s)\omega, v) \\
&\qquad   \varphi_{m';-k}(t-s) \varphi_{m;-10M_t }(t-s)    d \omega d y d v  d s,\\
\end{split}
\ee
where
\be\label{nov28eqn7}
\mathcal{K}_{k,j}^{ess }(y, \omega, v):= \int_{\R^3} e^{i y\cdot \xi}  \frac{\varphi^{ess}(\omega, \xi) }{i(\omega+\hat{v})\cdot \xi} \varphi_k(\xi) \varphi_j(v) d \xi . 
\ee
By doing integration by parts in $\xi$ in $\omega$ direction and directions perpendicular to $\omega$, the following estimate holds for the kernel, 
\be\label{nov28eqn9}
\begin{split}
&|\mathcal{K}_{k}^{ess }(y, \omega, v)| + |\nabla_v\mathcal{K}_{k}^{ess }(y, \omega, v) | \\
& \lesssim 2^{-m'+(1+ \epsilon)k+2j } \varphi_{[j-1, j+1]}(v)( 1+ 2^{(1+\epsilon)k/2-m'/2}|y\times \omega |)^{-N_0^3}(1+2^{k-2j}|y\cdot \omega|)^{-N_0^3}. \\
\end{split}
\ee 
From the above estimate of kernel, we have 
\be\label{nov28eqn52}
\begin{split}
 \sum_{\begin{subarray}{c}
 k\in [10j, \infty)\cap  \Z_+\\
 m'\in [-k, \epsilon M_t]\cap \Z\\
 m'+k/2\geq \epsilon M_t\\
 \end{subarray}}\big|J^{ess;0}_{k,m'}(t, x)\big| 
 \lesssim  \sum_{ k\in [10j, \infty)\cap  \Z_+} 2^{-k+10j - n j +(n-1)M_t+\epsilon M_t}\lesssim 2^{-10j}.
 \end{split}
\ee

Now, we focus on the estimate of $ J^{ess;1}_{k,m'}(t, x)$.  In terms of the spherical coordinates, we have $\omega=(\sin \theta\cos\phi, \sin \theta\sin\phi, \cos \theta) $. After  using the Cauchy-Schwarz inequality, from the estimate of kernel in  \eqref{nov28eqn9},  we have
\be\label{nov28eqn49}
\begin{split}
\big|J^{ess;1}_{k,m'}(t, x)\big|&\lesssim\sum_{U\in \{E, B\}}2^{-m'+(1+\epsilon) k+2j} \Big( \int_0^t \int_{\R^3}\int_{\R^3}\int_0^{2\pi}\int_0^{ \pi}   |U(s, x-y+(t-s)\omega)|^2 \\
&\qquad  \times  ( 1+ 2^{(1+\epsilon)k/2-m'/2}|y  |)^{-N_0^3}  \varphi_{[j-1, j+1]}(v) \sin \theta d \theta d \phi  d y d v d s \Big)^{1/2}\\
&\qquad \times \Big( \int_0^t \int_{\R^3}\int_{\R^3}\int_{\mathbb{S}^2}  f(s, x-y+(t-s)\omega, v)   (1+2^{k-2j}|y\cdot \omega|)^{-N_0^3}\\
&\qquad \times ( 1+ 2^{(1+\epsilon)k/2-m'/2}|y\times \omega |)^{-N_0^3} \varphi_{[j-1, j+1]}(v)  d \omega d y d v d s \Big)^{1/2}. \\
\end{split}
\ee
From the above estimate,  the Jacobian of changing of coordinates $(y_1, y_2, \theta)\rightarrow x- y + (t-s)\omega$,    the volume of support of $v, \omega$,   the conservation law  \eqref{conservationlaw}, and the assumption that $j\geq (1+2\epsilon)M_t$,  we have
\be\label{nov28eqn58}
\begin{split}
 \sum_{\begin{subarray}{c}
 k\in [10j, \infty)\cap  \Z_+\\
 m'\in [-k, \epsilon M_t]\cap \Z\\
 m'+k/2\geq \epsilon M_t\\
 \end{subarray}} &\big|J^{ess;1}_{k,m'}(t, x)\big|\lesssim   2^{-10k} +   \sum_{\begin{subarray}{c}
 k\in [10j, \infty)\cap  \Z_+\\
 m'\in [-k, \epsilon M_t]\cap \Z\\
 \end{subarray}} 2^{-m'+(1+\epsilon) k+2j} \big( 2^{-(1+\epsilon)k /2 + m'/2  + 3j}  \big)^{1/2}\\
  &  \times \big(\min\{ 2^{m'+2p+5j} 2^{-(2+\epsilon)k+ m'}, 2^{-(1+\epsilon)k /2 + m'/2  }2^{m'-nj + (n-1)M_t}\} \big)^{1/2} \\
   &   \lesssim 2^{-10j}+  \sum_{\begin{subarray}{c}
 k\in [10j, \infty)\cap  \Z_+\\
 m'\in [-k, \epsilon M_t]\cap \Z\\
 \end{subarray}}  2^{-m'+(1+\epsilon) k+2j}2^{-m'+(1+\epsilon) k+2j} \\
  &\quad \times \big( 2^{-(1+\epsilon)k /2 + m'/2  + 3j}  \big)^{1/2} \big(   2^{2m'+2p+5j-(2+\epsilon)k} \big)^{3/8} \big( 2^{-(1+\epsilon)k /2 +3 m'/2-100j  }\big)^{1/8} \\
  & \lesssim 2^{-10j}. \\
 \end{split} 
\ee
To sum up, our desired estimate  \eqref{nov26eqn31}  holds after combining the obtained estimates  \eqref{nov28eqn31}  and  \eqref{nov28eqn32}, the decompositions  \eqref{nov28eqn33},   \eqref{nov28eqn55}  and  \eqref{nov28eqn57}, and  the  obtained estimates  \eqref{nov28eqn51},   \eqref{nov28eqn56},  \eqref{nov28eqn52}, and  \eqref{nov28eqn58}.

\end{proof}
\subsection{The main part: $j\leq (1+2\epsilon)M_t$. }
Now, we focus on the case $j\leq (1+2\epsilon)M_t$. Recall the G-S decomposition in  \eqref{july5eqn31} in Lemma \ref{GSdecomloc}, we first estimate the $T$-part. We have
\begin{lemma}\label{pointroughT}
 For any   $t\in[0, T^{})$,   $j\in [0, (1+2\epsilon)M_t]\cap \Z_+, m\in [-10 M_t , \epsilon M_t]\cap\Z,  l\in [-j,2]\cap \Z$, $U\in \{E, B\},$  we have 
\be\label{july10eqn90}
\begin{split}
    |U_{T;j,l}^{m}(t, x)   | \lesssim 1 + 2^{5\epsilon M_t} \min\big\{ & 2^{m + l + j+2\tilde{\alpha}_t M_t} ,  2^{4\tilde{\alpha}_t M_t/3+ M_t/3-l/3},\\
    & \qquad  2^{-m+2\tilde{\alpha}_t M_t/3-j/3-5l/3},  |  x_{\bot}|^{-1/2} 2^{(1+2\epsilon)M_t} \big\}.  \\
    \end{split}
\ee
\end{lemma}
\begin{proof}
Recall  \eqref{july11eqn90}. Due to the cylindrical symmetry, without loss of generality, we assumed that $x=(|  x_{\bot}|, 0, x_3)$. 

From the detailed formulas of coefficients $\omega^{m,U}_{j;a}(\omega, v), $ in  \eqref{july9eqn11},  the volume of support of $\omega, v$, the estimate  \eqref{nov24eqn41}  if $|  v_{\bot} |\geq 2^{\alpha_t M_t+\epsilon M_t}$,  and the estimate  \eqref{march18eqn31}  in Lemma \ref{conservationlawlemma}, we have  
\be\label{july9eqn55}
\begin{split}
| U_{T;j,l}^{m}(t,x)   |   &\lesssim   \int_{0}^{t } \int_{\R^3} \int_{\mathbb{S}^2} 2^{-l} f(s,x+(t-s)\omega, v) \varphi_{[j-1, j+1]}(v) \\
&\quad \times  \varphi_{[l-2,l+2];-j}( \tilde{v}+\omega)\varphi_{m;-10M_t }(t-s)  d \omega d v d s\\
&  \lesssim  2^{2\epsilon M_t}\min\big\{ 2^{m-l} 2^{\epsilon M_t + 2l + j+2\tilde{\alpha}_t M_t}  , 2^{-2m-j-3l}\big\}\mathbf{1}_{m\in (-10M_t, \epsilon M_t]\cap \Z} \\
&\qquad + 2^{m+4j} \mathbf{1}_{m=-10M_t} \\
&\lesssim 2^{2\epsilon M_t}\min\{ 2^{m + l + j+2\tilde{\alpha}_t M_t+\epsilon M_t} , \big(2^{m + l + j+2\tilde{\alpha}_t M_t+\epsilon M_t} \big)^{2/3}\big(2^{-2m-j-3l} \big)^{1/3}, \\
&\qquad \big(2^{m + l + j+2\tilde{\alpha}_t M_t+\epsilon M_t} \big)^{1/3}\big(2^{-2m-j-3l} \big)^{2/3}  \} +1 \\
&\lesssim  2^{3\epsilon M_t}\min\{ 2^{m + l + j+2\tilde{\alpha}_t M_t} ,  2^{4\tilde{\alpha}_t M_t/3+ M_t/3-l/3},  2^{-m+2\tilde{\alpha}_t M_t/3-j/3-5l/3} \} +1.\\
\end{split}
\ee

Alternatively,  if we use the spherical coordinate s.t., $\omega=(\sin \theta\cos\phi, \sin \theta \sin \phi, \cos\theta)$ and do  dyadic decomposition for the size of $\sin \theta$ and $\sin\phi$, we have
\be\label{2026may2eqn11}
\begin{split}
|  U_{T;j,l}^{m}(t,x)   |  & \lesssim \sum_{p, q\in [-10 M_t, 2]\cap \Z} U^{m;p,q}_{j,l}(t,x)\\
  U^{m;p,q}_{j,l}(t,x)&=  \int_{0}^{t } \int_{\R^3} \int_{\mathbb{S}^2} 2^{-l} f(s,x+(t-s)\omega, v) \varphi_{[j-1, j+1]}(v) \varphi_{[l-2,l+2];-j}( \tilde{v}+\omega) \\ 
&\quad \times \varphi_{m;-10M_t }(t-s) \varphi_{p;-10M_t}(\sin \theta) \varphi_{q;-10M_t}(\sin \phi) d \omega d v d s.
\end{split}
\ee
From the volume of support of $s, v, \omega$, the following rough estimate holds, 
\[
\big| U_{T;j,l}^{m}(t,x) \big|\lesssim 2^{-l} 2^{2l+3j+m + 2p+q}. 
\]

The above estimate is sufficient to rule out the case either $p=-10M_t,$ or $q=-10M_t$, or $m=-10M_t$.  Now it suffices to consider the case   $p, q\in(-10  M_t, 2]\cap \Z$ and $m\in (-10M_t, \epsilon M_t]\cap \Z$.

Let 
 \be\label{jan18eqn11}
 \begin{split}
w&:= x_3+(t-s)\cos \theta,\\
  z&:= \sqrt{| x_{\bot}|^2 + \big( (t-s) \sin \theta \big)^2 + 2 |  x_{\bot}| (t-s) \sin \theta  \cos\phi}.\\
\end{split}
\ee
As a result of direct computations, the size of the Jacobian of changing coordinates $(\theta, \phi) \longrightarrow (z,\omega)$ is given as follows,   
\be\label{march18eqn66}
|J_{(\theta, \phi)\longrightarrow( w, z) }|=|\frac{\partial z}{\partial \phi}||\frac{\partial w}{\partial \theta}| \sim  z^{-1}| x_{\bot}| (t-s)^2 \big( \sin \theta\big)^2 \sin \phi. 
\ee
Therefore, from the above estimate of the Jacobian, the conservation law  \eqref{conservationlaw}, and the volume of support of $v, \omega$, we have
\be\label{nov29eqn1}
\begin{split}
| U^{m;p,q}_{j,l}(t,x)| &\lesssim 2^{-l}\min\big\{|  x_{\bot}|^{-1} 2^{-m-p-q-j}, 2^{m+2p+q + 3j+2l}  \big\}\\
&\lesssim 2^{-l} \big(| x_{\bot}|^{-1} 2^{-m-p-q-j} \big)^{1/2}\big(2^{m+2p+q + 3j+2l}   \big)^{1/2}\\
&\lesssim |  x_{\bot}|^{-1/2} 2^{(1+2\epsilon)M_t}. \\
\end{split}
\ee
Hence   the desired estimate  \eqref{july10eqn90} holds after combining the obtained estimates  \eqref{july9eqn55}  and the above estimate  \eqref{nov29eqn1}. 
\end{proof}

\begin{lemma}\label{pointroughS1}
  For any   $t\in[0, T^{})$,   $j\in [0, (1+2\epsilon)M_t]\cap \Z_+, m\in [-10M_t , \epsilon M_t]\cap\Z,  l\in [-j,2]\cap \Z$, $U\in \{E, B\},$  we have   
\be\label{july9eqn51}
\begin{split}
 |   U_{S;j,l}^{m;1}(t, x)   |  &   \lesssim  2^{2\epsilon M_t}\min\{  2^{m/2+j+ \tilde{\alpha}_t M_t}  , 2^{-m-2l}  \} \\ 
 &\qquad +2^{-10M_t}\big(\| E(s,x)\|_{L^\infty_{[0,t]}L^\infty_x } + \| B(s,x)\|_{L^\infty_{[0,t]}L^\infty_x } \big). \\
 \end{split}
\ee
\end{lemma}
\begin{proof}
Recall \eqref{july9eqn31}, and the detailed formulas of $ E_{S;j,l}^{m,i}(t,s,x)$ and $ B_{S;j,l}^{m,i}(t,s,x)$ in  \eqref{july9eqn31}    and the detailed formula of $EB^1(t,s,x ,\omega, v)$ in    \eqref{july1eqn13}.  Due to the cylindrical symmetry, without loss of generality, we assumed that $x=(|  x_{\bot}|, 0, x_3)$. Moreover, in terms of spherical coordinates, we have $\omega=(\sin\theta \cos\phi, \sin \theta \sin \phi, \cos\theta)$.

From the Cauchy-Schwarz inequality, after localizing the size of $t-s$, the following estimate holds for any $U\in \{E, B\}$,
\be 
\begin{split}
 |   U_{S;j,l}^{m;1}(t,s,x)   | &  \lesssim 2^{m-j-2l}\Big( \int_{0}^{t } \int_{\R^3} \int_{\mathbb{S}^2} \varphi_{m;-10M_t } (t-s)  \varphi_{[j-1,j+1]}(v)\varphi_{[l-2,l+2];-j}( \tilde{v}+\omega)\\
&\qquad \times \big[| E(s, x +(t-s)\omega)  -\omega \times B(s, x +(t-s)\omega) |^2\\
&\qquad +|\omega \cdot B(s, x +(t-s)\omega) |^2  \big] d \omega d v ds \Big)^{1/2}\\
&   \times \Big( \int_{0}^{t }\int_{\R^3} \int_{\mathbb{S}^2} \varphi_{m;-10M_t }(t-s)  f(s, x +(t-s)\omega,v)  \varphi_{j,l}(v, \omega) d \omega d v ds  \Big)^{1/2}.\\
\end{split}
\ee 

From the   estimate  \eqref{march18eqn31}  in Lemma \ref{conservationlawlemma} and the  volume of support of $\omega, v$, the estimate  \eqref{nov24eqn41}  if $|  v_{\bot} |\geq 2^{\alpha_t M_t+\epsilon M_t}$, the following estimate holds if $m\in (-10M_t, \epsilon M_t]\cap \Z, $
\be\label{aug23eqn31}
\begin{split}
 |   U_{S;j,l}^{m;1}(t,s,x)   |  &  \lesssim 2^{-j-2l} 2^{(3j+2l)/2+2\epsilon M_t }\big( \min\big\{   2^{m+ 2l + j+2\tilde{\alpha}_t M_t}    , 2^{-2m-j-2l}\big\} \big)^{1/2} \\
 & \lesssim \min\{  2^{m/2+j+( \tilde{\alpha}_t+ \epsilon) M_t}    , 2^{-m-2l} \}.
 \end{split}
\ee
Moreover, by using the volume of support of $v, \omega$, the following estimate holds if $m=-10M_t$,
\[
\begin{split}
 |   U_{S;j,l}^{m;1}(t,s,x)   |  & \lesssim 2^{2m-j-2l} 2^{3j+2l}\big(\| E(s,x)\|_{L^\infty_{[0,t]}L^\infty_x } + \| B(s,x)\|_{L^\infty_{[0,t]}L^\infty_x } \big)\\
 & \lesssim 2^{-10 M_t}\big(\| E(s,x)\|_{L^\infty_{[0,t]}L^\infty_x } + \| B(s,x)\|_{L^\infty_{[0,t]}L^\infty_x } \big).
\end{split}
\]
 Hence the desired estimate  \eqref{july9eqn51}  holds from the above two estimates.

\end{proof}
 
\begin{lemma}\label{pointS2}
   For any   $t\in[0, T^{  })$,   $j\in [0, (1+2\epsilon)M_t]\cap \Z_+, m\in [-10 M_t , \epsilon M_t]\cap\Z,  l\in [-j,2]\cap \Z$, $U\in \{E, B\},$  we have    
\be\label{july10enq5}
  |   U_{S;j,l}^{m;2}(t, x)   | \lesssim 2^{ 2 (\tilde{\alpha}_t+5\epsilon)M_t} |  x_{\bot}|^{-1/2} +2^{-10 M_t}\|B(s,x)\|_{L^\infty_{[0,t]}L^\infty_x}.
\ee
Moreover, we have
\be\label{nov28eqn66}
   |   U_{S;j,l}^{m;2}(t, x)   | \lesssim   2^{(1+2\tilde{\alpha}_t +6\epsilon)M_t} + 2^{-8M_t}\big(\| E(s, \cdot)\|_{L^\infty_{[0, t]}L^\infty_x } + \| B (s, \cdot)\|_{L^\infty_{[0, t]}L^\infty_x  } \big).
 \ee
\end{lemma}
\begin{proof}
 Due to the cylindrical symmetry, without loss of generality, we assumed that $x=(|  x_{\bot}|, 0, x_3)$ and $\omega=(\sin\theta \cos\phi, \sin \theta \sin \phi, \cos\theta)$.  Recall  \eqref{july9eqn31}  and \eqref{july1eqn13}.  Note that,  the associated coefficient of $\mathbf{P}(EB^2(t,s,x ,\omega, v))$ is better than  the associated coefficient of $\mathbf{P}_3(EB^2(t,s,x ,\omega, v))$. Moreover, we observe that ``$\p_{v_3}$'' derivative of the  symbol is better than the $\nabla_{  v_{\bot}}$ derivative of the symbol in $  U_{S;j,l}^{m;2}(t, x)$. Thanks to this observation, after doing  dyadic decomposition for  the size of $ \sin \theta, $ and $\phi $, the following estimate holds as a result of  direct computations, 
\be\label{july10eqn20}
   |   U_{S;j,l}^{m;2}(t, x)   |   \lesssim \sum_{p,q\in \Z\cap [-10 M_t,2]\cap \Z  }    I_{j;l}^{m;p, q}(t,x),
\ee
where 
\be\label{july3eqn36}
\begin{split}
I_{j;l}^{m;p,q}(t,x)& :=   \int_{0}^{t } \int_{\R^3} \int_{0}^{2\pi}\int_0^\pi 2^{-j-l+\max\{l,p\}} (t-s) \varphi_{m;-10M_t }(t-s) \varphi_{[j-1,j+1]}(v)  \\
&\quad \times  |B(s,|  x_{\bot} | +(t-s)\sin\theta\cos\phi, (t-s)\sin\theta \sin \phi, x_3+(t-s)\cos\theta)|\\
&\quad\times  f(s,   |  x_{\bot} | +(t-s)\sin\theta\cos\phi, (t-s)\sin\theta \sin \phi, x_3+(t-s)\cos\theta, v )  \\
& \quad \times  \varphi_{[l-2,l+2];-j}(\tilde{v}+\omega) \varphi_{p;-10M_t}(\sin \theta ) \varphi_{q;-10M_t}(\sin \phi ) \sin \theta d \theta d\phi d v d s.
 \end{split}
\ee 

By using the volume of support, the following rough estimate holds, 
\be\label{july10eqn26}
\|I_{j;l}^{m;p,q}(t,x)\|_{L^\infty_x}\lesssim \|B(s,x)\|_{L^\infty_{[0,t]}L^\infty_x} 2^{2m+2p+q} 2^{-j-l} 2^{3j+2l}. 
\ee
From the above estimate, we can rule out the case either $m$ or $p$ or $q$ equals $-10M_t$.

 Now, it remains to consider the fixed $m, p, q$ s.t.,   $ (m, p, q)\in  \big((-10M_t, \epsilon M_t]\cap \Z\big)^{3} $.  After using the Cauchy-Schwarz inequality,  we have
\[
\begin{split}
& \big|I_{j;l}^{m;p, q}(t,x)\big| \lesssim  2^{m-j-l+\max\{l,p\}}  \Big( \int_{t-2^{m}}^{t-2^{m-1}} \int_{\R^3} \int_{0}^{2\pi}\int_0^\pi   \varphi_{j,l}(v, \omega)    \varphi_{p;-10M_t}(\sin \theta ) \varphi_{q;-10M_t}(\sin \phi ) \\
&  \quad \times |B(s,|  x_{\bot} | +(t-s)\sin\theta\cos\phi, (t-s)\sin\theta \sin \phi, x_3+(t-s)\cos\theta)|^2     \sin \theta  d \theta d\phi d v d s\Big)^{1/2} \\
& \quad \times  \Big(  \int_{t-2^{m}}^{t-2^{m-1}} \int_{\R^3} \int_{0}^{2\pi}\int_0^\pi      \sin \theta \varphi_{p;-10M_t}(\sin \theta ) \varphi_{q;-10M_t}(\sin \phi )  \varphi_{j,l}(v, \omega) \\
&  \quad \times  f(s,   | x_{\bot} | +(t-s)\sin\theta\cos\phi, (t-s)\sin\theta \sin \phi, x_3+(t-s)\cos\theta, v ) d \theta d\phi d v d s\Big)^{1/2}.
\end{split}
\]

From the   estimate of the Jacobian of changing coordinates $(\theta, \phi)\longrightarrow(z,w)$ in  \eqref{jan18eqn11},  the volume of support of $v,\omega$, the estimate  \eqref{nov24eqn41},    the conservation law in  \eqref{conservationlaw}, the fact that $|  v_{\bot}|/|v|\sim 2^{p}$  if $p\geq l +10$, and the estimate  \eqref{nov24eqn41}   if $|  v_{\bot}|\geq 2^{(\alpha_t+\epsilon)M_t}$,  the following estimate holds 
\be\label{july10eqn11}
\begin{split}
\big|I_{j;l}^{m;p, q}(t,x)\big| & \lesssim 2^{m-j-l+\max\{l,p\}}\big( \frac{2^{m } \min\{2^{3j+2l}, 2^{j+2(\tilde{\alpha}_t+ \epsilon)M_t}   \} }{2^{2m+p+q}|  x_{\bot} |} \big)^{1/2} \\
&\quad \times \big( 2^{m+2p+q} \min\{2^{3j+2l}, 2^{j+2(\tilde{\alpha}_t+ \epsilon)M_t} \} \big)^{1/2}\\
  & \lesssim 2^{m+\max\{l,p\}+p/2}|  x_{\bot}|^{-1/2} 2^{-j-l} \min\{2^{3j+2l}, 2^{j+2(\tilde{\alpha}_t+4\epsilon)M_t} \}\\
  &\lesssim 2^{m+p/2 +2 (\tilde{\alpha}_t+4\epsilon)M_t} | x_{\bot}|^{-1/2}  . 
  \end{split}
\ee
Hence finishing the proof of our  
 desired estimate  \eqref{july10enq5}.

Note that, if $| x_{\bot}| \geq 2^{m+p-10}$, then the following estimate holds  from the obtained estimate  \eqref{july10eqn11}, 
\[
\big|I_{j;l}^{m;p, q}(t,x)\big|\lesssim  2^{m+p/2 +2 (\tilde{\alpha}_t+4\epsilon)M_t} |  x_{\bot}|^{-1/2}\leq 2^{m/2+2 (\tilde{\alpha}_t+4\epsilon)M_t}.
\]

   For  the case $|  x_{\bot}| \leq 2^{m+p-10}$,  we have   $|  x_{\bot} +(t-s)  \omega_{\bot}| \sim 2^{m+p }$. 
Recall the G-S decomposition in Lemma \ref{GSdecomloc}. From the estimate  \eqref{march18eqn31}  in Lemma \ref{conservationlawlemma}, the estimate  \eqref{nov26eqn31}  in Lemma \ref{erroresti},   the pointwise estimate   \eqref{july10eqn90}  in Lemma \ref{pointroughT},  \eqref{july9eqn51}  in Lemma \ref{pointroughS1},   the estimate  \eqref{july10enq5}  in Lemma \ref{pointS2}, the estimate in  \eqref{nov24eqn41}  if $|  v_{\bot}|\geq 2^{\alpha_t M_t+\epsilon M_t}$,   we have
\be\label{nov28eqn60}
\begin{split}
\big|I_{j;l}^{m;p, q}(t,x)\big|&\lesssim  2^{-j-l+\max\{l,p\}}   \min\{2^{-m-j-2l},  2^{2m+2p}\min\{2^{j+2(\tilde{\alpha}_t +5\epsilon)M_t }, 2^{3j+2l}\}  \} \\
 &\quad \times\big[2^{(1+\tilde{\alpha}_t +5\epsilon)M_t} + \frac{2^{2(\tilde{\alpha}_t +5\epsilon)M_t}}{(2^{m+p})^{1/2}} + 2^{-9M_t}\big(\| E(s, \cdot)\|_{L^\infty_{[0, t]}L^\infty_x }\\
 &\qquad  + \| B (s, \cdot)\|_{L^\infty_{[0, t]}L^\infty_x  } \big) \big].\\
\end{split}
\ee
 If $p\leq l+10$, then from the above obtained estimate,  we have
\[
\begin{split}
\big|I_{j;l}^{m;p, q}(t,x)\big|& \lesssim 2^{-j }  \big(2^{-m-j-2l}\big)^{1/2}   \big(2^{2m+2p+3j+2l}\big)^{1/2}   \\
&\quad \times\big[2^{(1+\tilde{\alpha}_t +5\epsilon)M_t} + \frac{2^{2(\tilde{\alpha}_t +5\epsilon)M_t}}{(2^{m+p})^{1/2}} + 2^{-9M_t}\big(\| E(s, \cdot)\|_{L^\infty_{[0, t]}L^\infty_x } + \| B (s, \cdot)\|_{L^\infty_{[0, t]}L^\infty_x  } \big) \big]  \\
&  \lesssim 2^{(1+\tilde{\alpha}_t +5\epsilon)M_t} 2^{p/2}+   2^{-8M_t}\big(\| E(s, \cdot)\|_{L^\infty_{[0, t]}L^\infty_x } + \| B (s, \cdot)\|_{L^\infty_{[0, t]}L^\infty_x  } \big).\\
\end{split}
\]
 If $p\geq l+10$,   we have $|  v_{\bot} |\sim 2^{j+p} $. From the obtained estimate  \eqref{nov28eqn60}  and the estimate  \eqref{nov24eqn41}  if $|  v_{\bot} |\geq  2^{(\alpha_t +\epsilon)M_t}$,  we have 
 \[
\begin{split}
\big|I_{j;l}^{m;p, q}(t,x)\big|& \lesssim 2^{-j-l-m/2+p/2 } (2^{-m-j-2l})^{1/4} ( \min\{ 2^{2m+2p+3j+2l}, 2^{2m+2l+j+2(\tilde{\alpha}_t+\epsilon) M_t }\} )^{3/4}  \\
&\quad \times \big[ 2^{(1+\tilde{\alpha}_t +5\epsilon)M_t}+ 2^{-9M_t}\big(\| E(s, \cdot)\|_{L^\infty_{[0, t]}L^\infty_x } + \| B (s, \cdot)\|_{L^\infty_{[0, t]}L^\infty_x  } \big) \big] \\
& \lesssim 2^{-j-l-m/2+p/2 } (2^{-m-j-2l})^{1/4} \big(   2^{2m+2p+3j+2l}\big)^{1/4}  \big(2^{2m+2l+j+2(\tilde{\alpha}_t+\epsilon) M_t } \big)^{1/2}   \\
&\quad \times\big[ 2^{(1+\tilde{\alpha}_t +5\epsilon)M_t}    + 2^{-9M_t}\big(\| E(s, \cdot)\|_{L^\infty_{[0, t]}L^\infty_x } + \| B (s, \cdot)\|_{L^\infty_{[0, t]}L^\infty_x  } \big) \big] \\
 & \lesssim  2^{(1+2\tilde{\alpha}_t +6\epsilon)M_t} + 2^{-8M_t}\big(\| E(s, \cdot)\|_{L^\infty_{[0, t]}L^\infty_x } + \| B (s, \cdot)\|_{L^\infty_{[0, t]}L^\infty_x  } \big).\\
\end{split}
\]
To sum up, our desired estimate  \eqref{nov28eqn66}  holds from the above estimates and the estimates  \eqref{july10eqn20}  and  \eqref{july10eqn26}. 

\end{proof} 

To sum up, we have 
\begin{proposition}\label{Linfielec}
Let $\epsilon:=100/n=10^{-7}$. The following rough estimates hold for  any   $U\in \{E, B\},$ $t\in[0, T^{})$,  
\be\label{july10eqn89}
 \| U(s,x)\|_{L^\infty_{[0,t]}L^\infty_x}  + \| U_{j,l}^m(s,x)\|_{L^\infty_{[0,t]}L^\infty_x}  \lesssim 2^{  (1+2\tilde{\alpha}_t +6\epsilon)M_t}.
\ee
\end{proposition}
\begin{proof}
Recall   \eqref{july5eqn1}. From the estimate \eqref{nov26eqn31}  in Lemma \ref{erroresti},   the estimate  \eqref{july10eqn90}  in Lemma \ref{pointroughT},   the estimate  \eqref{july9eqn51}  in Lemma \ref{pointroughS1}, and the   estimate  \eqref{nov28eqn66}  in Lemma \ref{pointS2},  we have
\be
\begin{split}
&\sum_{U\in \{E, B\}} \| U(s,x)\|_{L^\infty_{[0,t]}L^\infty_x}  + \| U_{j,l}^m(t,x)\|_{L^\infty_{[0,t]}L^\infty_x}  \\
& \lesssim 2^{-6M_t}\big( \sum_{U\in \{E, B\}} \| U(s,x)\|_{L^\infty_{[0,t]}L^\infty_x} \big) +  2^{  (1+2\tilde{\alpha}_t +6\epsilon)M_t}, \\
 \Longrightarrow &  \sum_{U\in \{E, B\}} \| U(s,x)\|_{L^\infty_{[0,t]}L^\infty_x}  + \| U_{j,l}^m(t,x)\|_{L^\infty_{[0,t]}L^\infty_x}   \lesssim  2^{  (1+2\tilde{\alpha}_t +6\epsilon)M_t}. 
\end{split}
\ee
Hence finishing the proof. 
\end{proof}
\subsection{A dichotomy  for the localized magnetic  field}

  The main goal of this subsection is to improve our understanding about     the electromagnetic field. Our key observation is  that the $L^2_x$-norm and the $L^\infty_x$-norm of the localized magnetic field can not be very big at the same time.  For example, the main enemy of the $L^\infty$-estimate are  the localized pieces of the magnetic field  with $|v|\sim 2^j, \angle(v, -\omega)\sim 2^{-j}, j\approx (1+2\epsilon)M_t$. From the estimate  \eqref{july5eqn66}  in Lemma \ref{firstL2},  we observe that  the $L^2$-norm of those enemy are very small.

With this intuition,  the main goal of this subsection is to give a  dichotomy  for the localized magnetic  field. We show that either the $L^\infty$-estimate of the localized magnetic field is improved if comparing with the rough estimate obtained in Proposition \ref{Linfielec} or the geometric mean of the $L^\infty$-norm and the $L^2$-norm is not very big. 

To state our main result,    we need some notation. For $\forall  t\in [0, T^{}),$ we define
\be\label{indexsetL2}
\begin{split}
   \mathcal{S}(t)&:= \{ (m,k,j,l): m\in [-10 M_t ,\epsilon M_t]\cap \Z, k,j\in \Z_+,l\in [-j,2]\cap \Z \}, \\
\mathcal{T}_0(t)&:=\{ (m,k,j,l): k\leq 20M_t, j\geq (1+2\epsilon)M_t\}\cup \{ (m,k,j,l): k\leq 20M_t, m =-10 M_t \}\\  
\mathcal{T}_1(t)&:=\{ (m,k,j,l): m\in (-10 M_t ,\epsilon M_t]\cap \Z,    m+k\leq -2(1+2\epsilon)l +\epsilon M_t, \\
&\qquad \qquad m+l\leq -(1+100\epsilon)M_t, k\leq 20M_t\},\\ 
\mathcal{T}_2(t)&:=   \{ (m,k,j,l):   m\in (-10 M_t ,\epsilon M_t]\cap \Z, m+k\geq -2(1+2\epsilon)l +\epsilon M_t, \\
&\qquad \qquad m+j+2l\geq  0,  k\leq 20M_t \}, \\
    \mathcal{S}_1(t)&:= \mathcal{S}(t)\cap \big(\cup_{i=0,1,2}\mathcal{T}_i(t) \big),\qquad \mathcal{S}_2(t):= \mathcal{S}(t)\cap  \big(\cup_{i=0,1,2}\mathcal{T}_i(t) \big)^c. 
\end{split}
\ee

Recall \eqref{july5eqn1} and \eqref{july5eqn60}. With the above notation, for any $s\in [0,t]$, the following decomposition holds  for the magnetic field,  
\be
B(s,x)= \sum_{(m,k,j,l)\in \mathcal{S}_1(t)\cup \mathcal{S}_2(t)} B^{  {m}}_{  {k};  {j},    {l}}(s,x).
\ee

Now, we are ready to state the main result of this subsection. 
 \begin{proposition}\label{meanLinfest}
For any $t\in [0, T^{})$,  we have 
\be\label{aug4eqn10}
\begin{split}
\sum_{(m,k,j,l)\in \mathcal{S}_1(t)} \sup_{s\in [0,t]}\|B^{  {m}}_{  {k};  {j},    {l}}(s,\cdot)\|_{L^\infty}  & \lesssim   2^{11\epsilon M_t} \big( 2^{ 2 \tilde{\alpha}_t  M_t }   + 2^{7M_t/6+\tilde{\alpha}_t M_t/4} \big),  \\
\sum_{(m,k,j,l)\in \mathcal{S}_2(t)}\big(\sup_{s\in [0,t]}\|B^{  {m}}_{  {k};  {j},    {l}}(s,\cdot)\|_{L^\infty} &  \big)^{\frac{1}{2}  } \big( \sup_{s\in [0,t]}\|B^{ {m}}_{  {k};  {j},  {l}}(s,\cdot)\|_{L^2}\big)^{ \frac{1}{2}}\\
&  \lesssim  2^{10\epsilon M_t} \big( 2^{\tilde{\alpha}_t M_t}   + 2^{7M_t/12+ \tilde{\alpha}_t M_t/8} \big). 
 \end{split}
\ee
\end{proposition}
\begin{remark}
 We remark that, due to the decomposition of the electromagnetic field in  \eqref{july1eqn11}, we only use the dichotomy of the magnetic field for the estimate of $EB^2(t,s,x +(t-s)\omega ,\omega, v)$, see  \eqref{july1eqn13}. For conciseness of presentation, we only present the  dichotomy for the localized  magnetic field in this paper. Similar dichotomy for the localized electric field   is also available but not pursued here. 
\end{remark}
\begin{proof}
For better presentation, we postpone the proof of this proposition   to section \ref{dictomyproposition}. In this section, we mainly estimate the $L^2_x$-norm and the $L^\infty_x$-norm of  the localized magnetic field in different regions, which are the main ingredients in the proof of the Proposition \ref{meanLinfest}. 
\end{proof}

In the following Lemma,  we   show that the rough estimate of the electromagnetic field can be improved if $m+j+2l\geq 0.$
\begin{lemma}\label{bulkroughpoi}
 For any $j\in [0, (1+2\epsilon)M_t]\cap \Z_+, m\in [-10 M_t ,\epsilon M_t]\cap\Z,  l\in [-j,2]\cap \Z$, $U\in \{E, B\},$ the following rough estimates hold for $t\in[0, T^{})$, 
\be\label{sep17eqn11}
\big| U^m_{ j,l}(t,x)\big|\lesssim 2^{10\epsilon M_t}\big( 2^{ 2 \tilde{\alpha}_t  M_t }    + 2^{7M_t/6+\tilde{\alpha}_t M_t/4} \big)(1+2^{-m-j-2l}).
\ee
Moreover, for any $x\in \R^3, $ s.t., $  x_{\bot}\neq 0, $ we have
\be\label{sep21eqn31}
| U(t,x)|+| U_{j,l}^m(t,x)| \lesssim   |  x_{\bot}|^{-1/2} 2^{ 5 \epsilon M_t}(2^{M_t}+ 2^{2 \tilde{\alpha}_t M_t }) + |  x_{\bot}|^{-1/4}  2^{5M_t/4+\tilde{\alpha}_t M_t/4+5\epsilon M_t}    +1. 
\ee
 \end{lemma}

 \begin{proof}
Recall the G-S decomposition in  \eqref{july5eqn31}  in  Lemma \ref{GSdecomloc} and the decomposition  \eqref{july9eqn41}. Note that, after combining the   estimate  \eqref{july10eqn89}  in Proposition \ref{Linfielec},  the estimate  \eqref{july10eqn90}  in Lemma \ref{pointroughT},   the estimate  \eqref{july9eqn51}   in Lemma \ref{pointroughS1}, and the estimate \eqref{july10enq5}  in Lemma \ref{pointS2},  we have
\be\label{dec1eqn49}
\begin{split}
   | U_{0;j,l}^m(t,x)|+  |U_{T;j,l}^{m}(t, x)   | +     |U_{S;j,l}^{m;1}(t, x)   | & \lesssim 2^{5\epsilon M_t}(2^{2\tilde{\alpha}_t M_t /3+2j/3}+2^{j}) (1+2^{-m-j-2l}),\\ 
      | U_{0;j,l}^m(t,x)|+  |U_{T;j,l}^{m}(t, x)   | +     |U_{S;j,l}^{m;2}(t, x)   | &  \lesssim
 2^{ 10 \epsilon M_t}(2^{M_t}+ 2^{2 \tilde{\alpha}_t M_t }) |  x_{\bot}|^{-1/2}   +1.
 \end{split}
 \ee
Therefore, it suffices to estimate 
  $  U_{S;j,l}^{m;2}(t, x) $ for the proof of the desired estimate  \eqref{sep17eqn11}  and estimate  $  U_{S;j,l}^{m;1}(t, x) $ for the proof of the desired estimate  \eqref{sep21eqn31}. 

  Due to the cylindrical symmetry, without loss of generality, we assumed that $x=(|  x_{\bot}|, 0, x_3)$. Moreover, in terms of spherical coordinates, we have $\omega=(\sin\theta \cos\phi, \sin \theta \sin \phi, \cos\theta), \theta\in [0,\pi], \phi\in[0,2\pi]$. 

   Recall   \eqref{july10eqn20}. From the obtained estimate  \eqref{july10eqn26}, we know that it suffices to  consider the case   $m, p, q$ are fixed s.t.,    $ (m, p, q)\in  \big( (-10M_t, \epsilon M_t]\cap \Z \big)^{3} $.

 Moreover, note that, if  $| x_{\bot} |\geq 2^{m+p-10}$, then from the obtained estimate \eqref{july10eqn11},   we have
 \be\label{july10eqn51}
 \big|  I_{j;l}^{m;p,q}(t,x)\big|\lesssim 2^{m+p/2 +2 (\tilde{\alpha}_t+ \epsilon)M_t} 2^{-(m+p)/2} dy \lesssim  2^{ 2 (\tilde{\alpha}_t+2\epsilon)M_t} .
\ee

Therefore, it suffices to consider the case   $| x_{\bot} |\leq 2^{m+p-10}$.  {For this case, we have $  |  x_{\bot}  + (t-s) \omega_{\bot}|\sim (t-s)|  \omega_{\bot}|\sim 2^{m+p}$.} Recall  \eqref{july3eqn36}. After using the G-S decomposition in \eqref{july5eqn31}  in  Lemma \ref{GSdecomloc} and the decomposition \eqref{july9eqn41}  for the magnetic field,  we have
\be\label{dec1eqn41}
\big|  I_{ j;l}^{m;p,q}(t,x)\big| \lesssim  I_{ j;l;0}^{m;p,q}(t,x)  +   I_{ j;l;1}^{m;p,q}(t,x)  +  I_{ j;l;2}^{m;p,q}(t,x),  
\ee
where
\be\label{dec1eqn1}
\begin{split}
 I_{ j;l;0}^{m;p,q}(t,x) &:=   \int_{0}^{t } \int_{\R^3} \int_{0}^{2\pi}\int_0^\pi 2^{m-j-l+\max\{l,p\}}    \sum_{\begin{subarray}{c} 
  \tilde{j}\in [(1+2\epsilon)M_t, \infty)\cap \Z_+\\
   \tilde{l}\in [-\tilde{j}, 2]\cap \Z \\ 
   \tilde{m}\in [ -10 M_t ,\epsilon M_t]\cap \Z\\
\end{subarray}
   } |B_{ \tilde{j}, \tilde{l}}^{\tilde{m}}(s, x+(t-s)\omega) | \\
 &\quad  \times f(s,   x+(t-s)\omega, v )  \varphi_{m;-10M_t }(t-s)    \varphi_{[j-1,j+1]}(v)\\
 &\quad \times \varphi_{[l-2,l+2];-j}( \tilde{v}+\omega)  \varphi_{p;-10M_t}(\sin \theta ) \varphi_{q;-10M_t}(\sin \phi ) \sin \theta d \theta d\phi d v d s, 
  \end{split}
\ee
\be\label{2024nov11eqn31}
\begin{split}
   I_{ j;l;1}^{m;p,q}(t,x) &:=   \int_{0}^{t } \int_{\R^3} \int_{0}^{2\pi}\int_0^\pi 2^{m-j-l+\max\{l,p\}} \varphi_{m;-10M_t }(t-s)    \varphi_{[j-1,j+1]}(v)   \\
    & \quad \times \varphi_{[l-2,l+2];-j}( \tilde{v}+\omega)\Big[\sum_{\begin{subarray}{c} 
   \tilde{m}\in [ -10 M_t ,\epsilon M_t]\cap \Z\\
    \tilde{j}\in  [0,(1+2\epsilon)M_t)\cap \Z_+,  
  \tilde{l}\in [-\tilde{j}, 2]\cap \Z\\
\end{subarray}
   } |B_{0;\tilde{j}, \tilde{l}}^{\tilde{m}}(s, x+(t-s)\omega)\\
& \quad+ |B_{T;\tilde{j}, \tilde{l}}^{\tilde{m}}(s, x+(t-s)\omega) |   +|B_{S;\tilde{j}, \tilde{l}}^{\tilde{m};2}(s, x+(t-s)\omega) | \Big] \\
&\quad \times  f(s,   x+(t-s)\omega, v ) \varphi_{p;-10M_t}(\sin \theta ) \varphi_{q;-10M_t}(\sin \phi ) \sin \theta  d \theta d\phi d v d s, 
\end{split}
\ee
\be\label{2024nov11eqn32}
\begin{split}
 I_{ j;l;2}^{m;p,q}(t,x) &:= \int_{0}^{t } \int_{\R^3} \int_{0}^{2\pi}\int_0^\pi 2^{m-j-l+\max\{l,p\}} \varphi_{m;-10M_t }(t-s)    \varphi_{[j-1,j+1]}(v) \\
 &\quad \times  \varphi_{[l-2,l+2];-j}( \tilde{v} +\omega)  \Big[\sum_{\begin{subarray}{c}    \tilde{m}\in [-10 M_t ,\epsilon M_t]\cap \Z\\
     \tilde{j}\in   [0,(1+2\epsilon)M_t)\cap \Z_+, \tilde{l}\in [-\tilde{j}, 2]\cap \Z\\
\end{subarray}
   } |B_{S;\tilde{j}, \tilde{l}}^{\tilde{m};1}(s, x+(t-s)\omega) |\Big]\\
   &\quad  \times   f(s,   x+(t-s)\omega, v ) \varphi_{p;-10M_t}(\sin \theta ) \varphi_{q;-10M_t}(\sin \phi )   \sin \theta d \theta d\phi d v d s.
 \end{split}
\ee
 
We estimate $I_{ j;l;i}^{m;p,q}(t,x), i\in\{0,1,2\},$ one by one and proceed in steps as follows. 

\medskip

\noindent \textbf{Step 1.}\qquad The estimate of $I_{ j;l;0}^{m;p,q}(t,x). $

\medskip

Recall  \eqref{dec1eqn1}. From  the estimate  \eqref{nov26eqn31}  in Lemma \ref{erroresti},   the estimate     \eqref{july10eqn89}   in Proposition \ref{Linfielec},  and the volume of support of $\omega, v$,  we have
\be\label{dec1eqn42}
| I_{ j;l;0}^{m;p,q}(t,x)| \lesssim 2^{-j-l} 2^{2m+3j+2l} 2^{-9M_t + (3+10\epsilon)M_t} \lesssim 1. 
\ee

\medskip

\noindent \textbf{Step 2.}\qquad The estimate of $I_{ j;l;1}^{m;p,q}(t,x). $

\medskip

 Recall \eqref{2024nov11eqn31}. From  the estimate   \eqref{dec1eqn49},   the estimate  \eqref{march18eqn31}  in Lemma \ref{conservationlawlemma}, and the volume of support of $\omega, v$,  we have
\be\label{sep17eqn1}
\begin{split}
 I_{ j;l;1}^{m;p,q}(t,x) &\lesssim  2^{-j -l+\max\{l,p\}+5\epsilon M_t}   \big(2^{ M_t}+  2^{2 \tilde{\alpha}_t  M_t } \big) 2^{-m/2-p/2}  \min\{2^{-m-j-2l}, 2^{2m+2p }  2^{3j+2l}  \}\\
& \lesssim 2^{-j -l+\max\{l,p\}+5\epsilon M_t} \big(2^{ M_t}+  2^{2 \tilde{\alpha}_t  M_t } \big) 2^{-m/2-p/2} \\
&\quad\times  \min\{\big( 2^{-m-j-2l}\big)^{5/6}\big( 2^{2m+2p }  2^{3j+2l}   \big)^{1/6}, \big( 2^{-m-j-2l}\big)^{ 1/2}\big( 2^{2m+2p }  2^{3j+2l}   \big)^{1/2}\}\\
&\lesssim  2^{ 5\epsilon M_t}\big(2^{ M_t}+  2^{2 \tilde{\alpha}_t  M_t } \big)  \min\{  2^{-m} 2^{-4j/3-7l/3+\max\{l,p\}-p/6},  2^{  -l+\max\{l,p\}}   2^{ p/2}\}\\
 & \lesssim 2^{ 5\epsilon M_t} \big(2^{ M_t}+  2^{2 \tilde{\alpha}_t  M_t } \big)(1+2^{ -m-j-2l} ).\\
 \end{split}
\ee

\medskip

\noindent \textbf{Step 3.}\qquad The estimate of $I_{ j;l;2}^{m;p,q}(t,x). $

\medskip

 Recall \eqref{2024nov11eqn32}.  To estimate $I_{ j;l;2}^{m;p,q}(t,x)$, we decompose $ B_{S;  {j},   {l}}^{ {m};1}(t, x) $ further. More precisely, after doing  dyadic decomposition for the size of $\sin \theta$ and $\sin\phi$ and using the Cauchy-Schwartz inequality and the conservation law \eqref{march18eqn31}, for any $U\in \{E, B\}$,  we have
\be\label{2024oct13eqn1}
\begin{split}
 |  U_{S;  {j},   {l}}^{ {m};1}(t, x)  |  & \lesssim \sum_{ \tilde{p}, \tilde{q}\in   [-10M_t,2 ] \cap\Z }    U_{  {j};  {l}}^{  {m};\tilde{p}, \tilde{q}}(t,x),\\
   U_{  {j};  {l}}^{  {m};\tilde{p}, \tilde{q }}(t,x)& := 2^{  {m}-  {j}-2 {l}
 } (2^{-2   {m} +3 {j}+2  {l}})^{1/2}    \Big( \int_{0}^{t }\int_{\R^3} \int_{\mathbb{S}^2}  f(s, x +(t-s)\omega,v)  \\
 &\times \varphi_{  {m};-10M_t }(t-s)  \varphi_{\tilde{p};-10M_t}(\sin \theta ) \varphi_{\tilde{q};-10M_t}(\sin \phi ) \varphi_{  {j},  {l}}(v, \omega) d \omega d v d s \Big)^{1/2}.
\end{split}
\ee
From the above estimate, we have   
\be\label{aug23eqn20}
\begin{split}
  I_{ j;l;2}^{m;p,q}(t,x)    \lesssim  &\sum_{ \tilde{p}, \tilde{q}\in   [-10M_t,2 ] \cap\Z }  I_{ j;l;\tilde{p},\tilde{q}}^{m;p,q}(t,x),\\
  I_{ j;l;\tilde{p},\tilde{q}}^{m;p,q}(t,x)  : = &\int_{0}^{t } \int_{\R^3} \int_{0}^{2\pi}\int_0^\pi 2^{m-j-l+\max\{l,p\}}  \varphi_{  {m};-10M_t } (t-s) \varphi_{[j-1,j+1]}(v) \\
  & \times \varphi_{[l-2,l+2];-j}( \tilde{v} +\omega)    \Big[\sum_{\begin{subarray}{c} 
     \tilde{j}\in   [0,(1+2\epsilon)M_t)\cap \Z_+, \tilde{l}\in [-\tilde{j}, 2]\cap \Z\\
   \tilde{m}\in [-10 M_t ,\epsilon M_t]\cap \Z
\end{subarray}
   } |    B_{\tilde{j};\tilde{l}}^{\tilde{m};\tilde{p}, \tilde{q }} (s, x+(t-s)\omega) \Big] \\
   &\times   f(s,   x+(t-s)\omega, v )       \varphi_{p;-10M_t}(\sin \theta ) \varphi_{q;-10M_t}(\sin \phi )  \sin \theta d \theta d\phi d v d s.\\
\end{split}
\ee

 As before, the case either $\tilde{m}=-10M_t$, or $\tilde{p}=-10M_t$, or $\tilde{q}=-10M_t$ is trivial,  it suffices to  consider the case   $\tilde{m},\tilde{p}, \tilde{q}$ are fixed s.t.,    $ (\tilde{m},\tilde{p}, \tilde{q})\in  \big( (-10M_t, \epsilon M_t]\cap \Z \big)^{3} $.    Note that, $|  v_{\bot} |/|v|\sim 2^{\tilde{p}}$ if  $\tilde{p}\geq  {l} +10$.  From this fact,  the conservation law  \eqref{conservationlaw}, the estimate  \eqref{march18eqn31}  in Lemma \ref{conservationlawlemma},   the volume of support of $\omega, v$, and the estimate  \eqref{nov24eqn41}  if $|  v_{\bot} |\geq 2^{\alpha_t M_t+\epsilon M_t}$, we have 
\be\label{nov30eqn31}
\begin{split}
 |   U_{  {j};  {l}}^{  {m};\tilde{p}, \tilde{q }}(t,x) |   & \lesssim       2^{ {m}- {j}-2  {l}
 } (2^{-2   {m} +3 {j}+2  {l}})^{1/2} \big( \min\{2^{- {m}} 2^{-\tilde{p}-\tilde{q }} | x_{\bot}|^{-1} 2^{-  {j}} , 2^{  {m}+2\tilde{p}+ \tilde{q }+3  {j}+2  {l}} \} \big)^{1/2}  \\
 & \lesssim |  x_{\bot}|^{-1/4} 2^{  {j}-  {l}/2+\tilde{p}/4+2\epsilon M_t} \\
 &\lesssim |  x_{\bot}|^{-1/4} 2^{5\epsilon M_t}\min\big\{ 2^{- \tilde{p}/4+ M_t+\tilde{\alpha}_t M_t/2 }, 2^{ 5M_t/4+\tilde{\alpha}_t M_t/4 }\big\}. \\
 \end{split}
\ee

As a byproduct, together with the estimates \eqref{dec1eqn49} and  \eqref{2024oct13eqn1},  the above estimate concludes   the desired estimate \eqref{sep21eqn31}.

Now, we focus on the estimate of $I_{ j;l;2}^{m;p,q}(t,x).$  From the above obtained  estimate  \eqref{nov30eqn31},   the estimate  \eqref{march18eqn31}  in Lemma \ref{conservationlawlemma}, and the volume of support of $\omega$ and  $v$,  we have
 \be\label{nov30eqn68}
 \begin{split}
 |I_{ j;l;\tilde{p},\tilde{q}}^{m;p,q}(t,x)|\lesssim 2^{ -j-l+\max\{l,p\}}2^{-m/4-p/4} & 2^{-\tilde{p}/4+ M_t+\tilde{\alpha}_t M_t/2+6\epsilon M_t} \\
 & \times \min\{2^{-m-j-2l}, 2^{2m+2p }   2^{3j+2l}  \}.
 \end{split}
\ee

\medskip 

Based on the possible size of $\tilde{p}$, we split into two sub-steps as follows. 

\medskip

\noindent \textbf{Step 3A.}\qquad   If $\tilde{p}\geq p -\epsilon M_t $.
\medskip

From the obtained estimate  \eqref{nov30eqn68}, we have  \be\label{dec1eqn45}
  \begin{split}
  |I_{ j;l;\tilde{p},\tilde{q}}^{m;p,q}(t,x)| & \lesssim      2^{ -j-l+\max\{l,p\}}2^{-m/4-p/2+ M_t+\tilde{\alpha}_t M_t/2+6.25\epsilon M_t} \\
  &\quad \times  \big(2^{-m-j-2l}\big)^{3/4} \big( 2^{2m+2p }   2^{3j+2l}  \big)^{1/4}\\
&\lesssim 2^{M_t+\tilde{\alpha}_t M_t/2+7\epsilon M_t } (2^{-m-j-2l}).\\ 
\end{split}
\ee

\medskip

\noindent \textbf{Step 3B.}\qquad  If $\tilde{p}\leq p -\epsilon M_t$. 

\medskip 

Recall \eqref{aug23eqn20}. From   the obtained estimates    \eqref{2024oct13eqn1}, the Cauchy-Schwarz inequality,   the conservation law  \eqref{march18eqn31}  in Lemma \ref{conservationlawlemma},   and the volume of support of $v $, we have
\be\label{nov30eqn63}
\begin{split}
|I_{ j;l;\tilde{p},\tilde{q}}^{m;p,q}(t,x)|& \lesssim   
  2^{m-j-l+\max\{l,p\}} \big(2^{-2m-j-2l} \big)^{1/2}2^{\tilde{j}/2-\tilde{l} } 2^{(3  {j}+2  {l})/2+\epsilon M_t} \\
  &\quad \times \Big( \int_0^t \int_{\tau}^{t} \int_{\R^3} \int_{\mathbb{S}^2}  \int_{\mathbb{S}^2}    f(\tau, x+(t-s)\omega +(s-\tau)\tilde{\omega},u)  \\
  & \quad  \times     \varphi_{\tilde{j},\tilde{l}}(u,  \tilde{\omega} ) \varphi_{  {m};-10M_t } (t-s)  \varphi_{p;-10M_t}(\sin \theta ) \\
&\quad  \times     \varphi_{q;-10M_t}(\sin \phi )\varphi_{\tilde{p};-10M_t}(\sin \tilde{\theta} ) \varphi_{\tilde{q};-10M_t}(\sin \tilde{\phi} ) d \omega d \tilde{\omega}    d u ds  d \tau  \big)^{1/2}.
\end{split}
\ee
After doing dyadic decomposition for the size of $s$, we have 
\[
\begin{split}
\int_0^t \int_{\tau}^{t} \int_{\R^3} \int_{\mathbb{S}^2}  \int_{\mathbb{S}^2} &f(\tau, x+(t-s)\omega +(s-\tau)\tilde{\omega},u)   \varphi_{p;-10M_t}(\sin \theta ) \varphi_{ q;-10M_t}(\sin \phi )\varphi_{\tilde{p};-10M_t}(\sin \tilde{\theta} )   \\
& \times \varphi_{\tilde{q};-10M_t}(\sin \tilde{\phi} )  \varphi_{\tilde{j},\tilde{l}}(u,  \tilde{\omega} ) d \omega d \tilde{\omega}    d u ds  d \tau \lesssim \sum_{\kappa\in [-10M_t, \epsilon M_t]\cap \Z} H_{\kappa}(t,x),  
\end{split}
\]
where
\be\label{nov30eqn51} 
\begin{split}
  H_{\kappa}(t,x):= \int_0^t \int_{\tau}^{t} \int_{\R^3} \int_{\mathbb{S}^2}  \int_{\mathbb{S}^2} &f(\tau, x+(t-s)\omega +(s-\tau)\tilde{\omega},u) \varphi_{\tilde{j},\tilde{l}}(u,  \tilde{\omega} )\\
   &\times \varphi_{  {m};-10M_t } (t-s)  \varphi_{\kappa;-10M_t}(s-\tau)  \varphi_{p;-10M_t}(\sin \theta )  \\
   &\times \varphi_{q;-10M_t}(\sin \phi )\varphi_{\tilde{p};-10M_t}(\sin \tilde{\theta} ) \varphi_{\tilde{q};-10M_t}(\sin \tilde{\phi} )    d \omega d \tilde{\omega}    d u ds  d \tau. 
   \end{split}
\ee
From the volume of support of $H_{\kappa}(t,x)$, we have
\be\label{nov30eqn60}
|H_{\kappa}(t,x)|\lesssim 2^{ m+ \kappa + 2\tilde{p}+2p + 3\tilde{j}+2\tilde{l} +5\epsilon M_t}.
\ee 

From the above estimate, we can rule out the case $\kappa=-10M_t.$ Moreover, for any $\kappa\in (-10M_t, \epsilon M_t]\cap \Z, $ s.t., $\kappa\leq m $, the Jacobian of changing coordinates $(\phi, \theta, s)\longrightarrow  x+(t-s)\omega +(s-\tau)\tilde{\omega}$ can be estimated as follows, 
\[
\begin{split}
&\big| det\big(\begin{bmatrix}
\sin \theta \cos\phi - \sin  \tilde{\theta} \cos\tilde{\phi} & \sin \theta \sin\phi - \sin \tilde{\theta}\sin \tilde{\phi} & \cos \theta - \cos  \tilde{\theta}\\ 
-(t-s) \sin \theta \sin\phi & (t-s)  \sin \theta \cos\phi & 0\\ 
(t-s) \cos\theta\cos\phi & (t-s) \cos\theta \sin \phi & -(t-s) \sin\theta\\ 
\end{bmatrix}\big)\big|\\
&= (t-s)^2\sin \theta\big[ 1-\cos(\theta-  \tilde{\theta})+\sin \theta\sin  \tilde{\theta}(1-\cos(\phi- \tilde{\phi}) )\big]\\
&\gtrsim (t-s)^2\sin \theta | 1-\cos(\theta-  \tilde{\theta})|. \\
\end{split}
\]

Note that, as $\tilde{p}\leq p -\epsilon M_t$, we have  $|\theta-\tilde{\theta}|\sim 2^{ {p}}$.  From the above estimate of Jacobian, the volume of support of $\tilde{\omega}$, and the conservation law  \eqref{conservationlaw},  we have 
\be\label{nov30eqn61}
\big| H_{\kappa}(t,x) \big| \lesssim \int_0^t 2^{-2m - 2p +2\tilde{p}-\tilde{j}} \psi_{\leq \kappa+2}(s-\tau)d\tau  \lesssim 2^{-2m + \kappa -\tilde{j}}. 
\ee
If $\kappa\geq m $, then we change  coordinates $(\tilde{\phi}, \tilde{\theta}, s)\longrightarrow  x+(t-s)\omega +(s-\tau)\tilde{\omega}$. As a result, we have
\[
\big| H_{\kappa}(t,x) \big| \lesssim \int_0^t 2^{-2\kappa - 2p +2{p}-\tilde{j}} \psi_{\leq \kappa+2}(s-\tau)d\tau  \lesssim 2^{- \kappa -\tilde{j}}. 
\]
To sum up, in whichever case, we have
\[
\big| H_{\kappa}(t,x) \big| \lesssim \int_0^t 2^{-2\kappa - 2p +2{p}-\tilde{j}} \psi_{\leq \kappa+2}(s-\tau)d\tau  \lesssim 2^{- \max\{m,\kappa\} -\tilde{j}}. 
\]

After combining the obtained estimates  \eqref{nov30eqn63} , \eqref{nov30eqn51},   \eqref{nov30eqn60}, and   \eqref{nov30eqn61}, we have
\[
\begin{split}
|I_{ j;l;\tilde{p},\tilde{q}}^{m;p,q}(t,x)|&\lesssim  2^{ -l+\max\{l,p\}}  2^{\tilde{j}/2-\tilde{l} }  \Big( \sum_{\kappa\in [-10M_t, \epsilon M_t]\cap \Z} \min\{2^{2\kappa + 2\tilde{p}+2p + 3\tilde{j}+2\tilde{l} +5\epsilon M_t}, 2^{-\kappa -\tilde{j}}\} \\
&\qquad \times \mathbf{1}_{\kappa\in (-10M_t, \epsilon M_t]\cap \Z } + 2^{-10M_t}  \Big)^{1/2}\\
& \lesssim 2^{ -l+\max\{l,p\}}  2^{2\tilde{j}/3-2\tilde{l}/3} 2^{ p/3+ \tilde{p}/3 +3\epsilon M_t } +1 \\
&\lesssim  2^{ -l+\max\{l,p\}}    2^{4M_t/3+ p/3+ \tilde{p}/3 +5\epsilon M_t}+1.  \\
\end{split}
\]
From the above estimate and the obtained estimate  \eqref{nov30eqn68}, we have
 \be\label{dec1eqn46}
\begin{split}
|I_{ j;l;\tilde{p},\tilde{q}}^{m;p,q}(t,x)|&\lesssim \min\{ 2^{-m/4-p/4}    2^{-\tilde{p}/4+ M_t+\tilde{\alpha}_t M_t/2  }  2^{ 6\epsilon M_t} 2^{-m-2j-3l} ,       2^{ -l+4M_t/3+ p/3+ \tilde{p}/3 } 2^{ 5\epsilon M_t} \}\\
&\lesssim \big(2^{-m/4-p/4}    2^{-\tilde{p}/4+ M_t+\tilde{\alpha}_t M_t/2  }  2^{ 6\epsilon M_t} 2^{-m-2j-3l}\big)^{1/2}\big(  2^{ -l+4M_t/3+ p/3+ \tilde{p}/3 } 2^{ 5\epsilon M_t} \big)^{1/2}\\
&\lesssim 2^{7M_t/6+\tilde{\alpha}_t M_t/4+6\epsilon M_t}2^{-m-j-2l}.
\end{split}
\ee
Hence, after combining the obtained estimates \eqref{aug23eqn20},  \eqref{dec1eqn45}, and   \eqref{dec1eqn46},  we have
\be\label{dec1eqn52}
I_{ j;l;2}^{m;p,q}(t,x) \lesssim  2^{7\epsilon M_t }\big( 2^{M_t+\tilde{\alpha}_t M_t/2}  + 2^{7M_t/6+\tilde{\alpha}_t M_t/4}\big) 2^{-m-j-2l}.
\ee

To sum up, our desired estimate \eqref{sep17eqn11} holds after  combining the obtained first estimate in   \eqref{dec1eqn49},  \eqref{july10eqn20},   \eqref{july10eqn51},    \eqref{dec1eqn41},  \eqref{dec1eqn42},   \eqref{sep17eqn1}, and the above estimate  \eqref{dec1eqn52}.

 \end{proof}

Now, we consider the $L^2$-estimate of the localized magnetic field for the case $m\in (-10M_t, \epsilon M_t]\cap \Z$, in which we have $t-s\sim 2^m$. Recall  \eqref{july5eqn1}. From the stationary phase point of view, we know that the angle between $\omega$ and $\xi$  in $ B_{k;j,l}^m(t,x)$ are localized within the size of $2^{\max\{-m-k-l, -m/2-k/2\}+\epsilon M_t}$. More precisely, 
we have
\be\label{sep16eqn13} 
B_{k;j,l}^m(t,x)= B_{k;j,l}^{ess;m}(t,x)+ B_{k;j,l}^{err;m}(t,x),
\ee
where $\forall \star\in \{ess,err\},$ we have 
\be\label{sep16eqn1} 
\begin{split}
B_{k;j,l}^{\star;m}(t,x):=\int_0^t \int_{\R^3}   \int_{\R^3}\int_{\mathbb{S}^2}& e^{i(x+(t-s)\omega)\cdot \xi} (t-s) {i (\hat{v}\times \xi) }    \hat{f}(s, \xi , v) \\
&\times \varphi_{j, l}(v,  {\omega}) \varphi_k(\xi)  \varphi_{m,k,l}^{\star}(\omega , \xi) \varphi_{m;-10M_t }(t-s)  d  {\omega}  d\xi d  v d s,\\
\end{split}
\ee
where  the cutoff functions $\varphi_{m,k,l}^{\star}(\omega , \xi), \star\in \{ess,err\}$, are defined as follows, 
\be\label{sep16eqn2} 
 \varphi_{m,k,l}^{ ess }(\omega , \xi)  = \psi_{\leq   {c}(m,k,l)}( \omega \times \tilde{\xi}),\quad  \varphi_{m,k,l}^{err}(\omega , \xi) =  \psi_{ >  {c}(m,k,l)}( \omega \times \tilde{\xi}),  
 \ee 
where $ {c}(m,k,l):= \max\{-m-k-l, -m/2-k/2\}+\epsilon^2 (M_t +k_{+}).$ Recall  \eqref{march4eqn41}.  Note that, by doing integration by parts in $\omega$ once, we gain at least $2^{-\epsilon^2 (M_t +k_{+})}$ for the error type term  $B_{k;j,l}^{err;m}(t,x)$. After doing integration by parts in $\omega$ many times, e.g., $\epsilon^{-3}$ times, we have 
\be\label{sep16eqn30}
\begin{split}
 \| B_{k;j,l}^{err;m}(t,x)\|_{L^\infty\cap L^2}  & \lesssim 2^{-100(M_t+k_+)} 2^{4k } \min\{2^{-j}, 2^{-nj + (n-1)M_t}\} \\
 & \lesssim  2^{-100M_t-k}(1+2^{j-(1+2\epsilon)M_t})^{-10}. \\
 \end{split}
\ee

Now, we focus on the estimate of the essential part $B_{k;j,l}^{ess;m}(t,x)$. We first consider the case $k\in [ 0,20M_t]\cap \Z$ and   $m+k\geq -2 l +\epsilon M_t$. 
\begin{lemma}[$L^2$-estimate of the localize magnetic field]\label{firstL2}
 Let $t\in [0, T)$, $k\in [ 0,20M_t]\cap \Z$,  $j\in  [0, (1+2\epsilon)M_t]\cap  \Z_+, l\in [-j, 2]\cap \Z$, $m\in (-1 0M_t , \epsilon M_t]\cap \Z$, s.t., $m+k\geq -2 l +\epsilon M_t$.  Then  the following $L^2$-type  estimate holds, 
\be\label{july5eqn66}
\| B^m_{k;j,l}(t,x)\|_{L^2} \lesssim   2^{m+j+2l+ \epsilon M_t}+2^{-100 M_t}. 
\ee
\end{lemma}
\begin{proof}
 Note that, for the case we are considering, we have ${c}(m,k,l)=-m/2-k/2+\epsilon^2 (M_t+k_{+})\leq l - \epsilon M_t/3$.  Recall  \eqref{sep16eqn13}  and  \eqref{sep16eqn1}.   Similar to what we did in  \eqref{sep16eqn22}, from the volume of support of $v, \omega$, and the conservation law  \eqref{conservationlaw}, we have 
\be 
\begin{split}
\|B_{k;j,l}^{ess;m}(t,x)\|_{L^2}&\lesssim \sup_{g\in L^2, \|g\|_{L^2}=1} \big|\int_0^t \int_{\R^3}   \int_{\R^3}\int_{\mathbb{S}^2}   (t-s) \big|{  (\hat{v}\times \xi) } \big| \big|\overline{\hat{g}(\xi)}\big|  \big|\hat{f}(s, \xi , v)\big|  \\
&\qquad  \times  \varphi_{j, l}(v,  {\omega})    \varphi_{m,k,l}^{ ess }(\omega , \xi) \varphi_{m;-10M_t }(t-s)  \varphi_k(\xi)     d  {\omega} d\xi  d  v d s\big|  \\
&\lesssim \sup_{s\in [0,t]} \sup_{g\in L^2, \|g\|_{L^2}=1}  2^{2m+k+l} 2^{2 {c} (m,k,l)}    \int_{\R^3}   \int_{\R^3}\big|\hat{f}(s, \xi , v)\big|  \big|\overline{\hat{g}(\xi)}\big| \\
&\qquad  \times  \varphi_k(\xi)   \psi_{[l-2,l+2]}(\tilde{v}\times \tilde{\xi})  \varphi_j(v) d v d \xi \\
& \lesssim \sup_{s\in [0,t]} \sup_{g\in L^2, \|g\|_{L^2}=1} 2^{m+l+ \epsilon M_t} (2^{3j+2 l })^{1/2}\big(  \int_{\R^3}   \int_{\R^3}\big|\hat{f}(s, \xi , v)\big|^2 \varphi_j(v)   d v d \xi  \big)^{1/2}\\
&\lesssim    2^{m+j+2l+ \epsilon M_t}.\\
\end{split}
\ee
After combining the above estimate and the obtained estimate  \eqref{sep16eqn30}, our desired estimate  \eqref{july5eqn66}  holds.

\end{proof}

Now, we   consider the case $k\in [ 0,20M_t]\cap \Z$ and   $m+k\leq -2 l +\epsilon M_t$.
 
\begin{lemma}\label{smallk}
Let $t\in [0, T)$, $k\in [0,2 0M_t]\cap \Z$,  $ j\in  [0, (1+2\epsilon)M_t]\cap  \Z_+, l\in [-j, 2]\cap \Z$, $m\in ( -1 0M_t , \epsilon M_t]\cap \Z$, s.t., $m+k\leq -2 l +\epsilon M_t$. If $m+l \geq -M_t-100\epsilon M_t$, then  the following estimate holds, 
\be\label{sep16eqn41} 
\big(\|B_{k;j,l}^{ m}(t,x)\|_{L^\infty}  \big)^{1/2}\big(\|B_{k;j,l}^{ m}(t,x)\|_{L^2}   \big)^{1/2}\lesssim 2^{ M_t/2+ 10\epsilon M_t}. 
\ee 
If $m+l \leq -M_t-100\epsilon M_t$   then we have
\be\label{sep16eqn42} 
\|B_{k;j,l}^{ m}(t,x)\|_{L^\infty}\lesssim   2^{2\tilde{\alpha}_t M_t-50\epsilon M_t}.
\ee

\end{lemma}
\begin{proof}

 Recall  \eqref{sep16eqn1}. In terms of kernel, we have
\be\label{sep16eqn11}
\begin{split}
|B_{k;j,l}^{ess;m}(t,x)|\lesssim \int_0^t     \int_{\R^3} \int_{\R^3}\int_{\mathbb{S}^2} &(t-s) K_{m,k,l}(y, v, \omega ) f(s, x-y +(t-s)\omega) \\
&\times  \varphi_{j, l}(v,  {\omega})   \varphi_{m;-10M_t }(t-s)  d\omega d y d v d s, \\
\end{split}
\ee
where the kernel $  K_{m,k,l}(y, v)$ is defined as follows, 
\be\label{sep16eqn40} 
 K_{m,k,l}(y, \omega,v ):= \int_{\R^3} e^{i y \cdot \xi} i(\hat{v}\times \xi)   \varphi_{m,k,l}^{ess}(\omega , \xi) \varphi_k(\xi) d \xi.
\ee  
 By doing integration by parts in $\xi$ along $v$ direction and directions perpendicular to $v$ many times, we have
\be\label{sep16eqn12} 
 \big| K_{m,k,l}(y, \omega, v )\big|\lesssim 2^{4k+3{c}(m,k,l)+\epsilon M_t} (1+2^{k}|y\cdot \tilde{v}|)^{-100}  (1+2^{k+ {c}(m,k,l)  }|y\times \tilde{v}|)^{-100},
\ee
where $ {c}(m,k,l):= \max\{-m-k-l, -m/2-k/2\}+\epsilon^2 (M_t +k_{+}).$ 

Recall  \eqref{sep16eqn13}  and \eqref{sep16eqn11}. From the estimate of kernel in   \eqref{sep16eqn12}, the conservation law in  \eqref{conservationlaw}, and the volume of support of $\omega, v$, the estimate \eqref{nov24eqn41}  if $|  v_{\bot} |\geq 2^{(\alpha_t +\epsilon)M_t}$,  we have 
\be\label{sep16eqn21}
\begin{split}
|B_{k;j,l}^{ess;m}(t,x)|& \lesssim 2^{2m+ 5\epsilon M_t}\min\{2^{4k+3(-m-k-l) } 2^{2l-j}, 2^{k+(-m-k-l) } 2^{2l+j+2 \tilde{\alpha}_t  M_t} \}\\
&\lesssim 2^{ 5\epsilon M_t}\min\{2^{-m + k -l-j}  , 2^{m+l+j+2 \tilde{\alpha}_t M_t}\}\\
& \lesssim 2^{10\epsilon M_t}\min\{2^{-m + k -l-j}  , 2^{m+l+j+2 \tilde{\alpha}_t M_t}, \big(2^{-m + k -l-j}\big)^{1/3} \big(2^{m+l+j+2 \tilde{\alpha}_t M_t} \big)^{2/3} \} \\
& \lesssim 2^{10\epsilon M_t}\min\{2^{-m + k -l-j}  , 2^{m+l+j+2 \tilde{\alpha}_t M_t},2^{2M_t/3+4\tilde{\alpha}_t M_t/3+2\epsilon M_t}  \}.\\
\end{split}
\ee

Recall  \eqref{sep16eqn1}. From the volume of support of $v, \omega$, and the conservation law  \eqref{conservationlaw},  we also have the following $L^2$-type estimate, 
\be\label{sep16eqn22}
\begin{split}
\|B_{k;j,l}^{ess;m}(t,x)\|_{L^2}& \lesssim \sup_{g\in L^2, \|g\|_{L^2}=1} \big|\int_0^t \int_{\R^3}   \int_{\R^3}\int_{\mathbb{S}^2}   (t-s) \big|{  (\hat{v}\times \xi) } \big| \big|\overline{\hat{g}(\xi)}\big|  \big|\hat{f}(s, \xi , v)\big|    \\
&\quad  \times   \varphi_{m,k,l}^{ ess }(\omega , \xi)  \varphi_{m;-10M_t }(t-s)   \varphi_k(\xi) \varphi_{j, l}(v,  {\omega})  d  {\omega} d\xi d  v d s\big|  \\
 &\lesssim \sup_{s\in [0,t]} \sup_{g\in L^2, \|g\|_{L^2}=1} 2^{ 2m+k+2l+  {c} (m,k,l) }  \int_{\R^3}   \int_{\R^3}\big|\hat{f}(s, \xi , v)\big|  \big|\overline{\hat{g}(\xi)}\big| \\
 &\quad \times  \varphi_k(\xi)      \psi_{\leq {c}(m,k,l)+\epsilon M_t}(\tilde{v}\times \tilde{\xi})  \varphi_j(v) d v d \xi \\
& \lesssim \sup_{s\in [0,t]} 2^{ 2m+k+2l}  2^{   {c} (m,k,l) +2\epsilon M_t} (2^{3j+2   {c}(m,k,l) })^{1/2}\\
& \quad \times \big(  \int_{\R^3}   \int_{\R^3}\big|\hat{f}(s, \xi , v)\big|^2 \varphi_j(v)   d v d \xi  \big)^{1/2}\\
& \lesssim  2^{-k+j+20\epsilon M_t}. 
\end{split}
\ee

 Based on the possible size of $m+l$, we split into two cases as follows. 

\medskip

\noindent \textbf{Case 1.} \qquad If $m+l\geq -M_t-100\epsilon M_t$. 

\medskip

From the obtained estimates  \eqref{sep16eqn21}  and  \eqref{sep16eqn22}, we know that
\[
\begin{split}
&\big(\|B_{k;j,l}^{ess;m}(t,x)\|_{L^\infty}  \big)^{1/2}\big(\|B_{k;j,l}^{ess;m}(t,x)\|_{L^2}   \big)^{1/2}\\
&\lesssim ( 2^{20\epsilon M_t} 2^{-m + k -l-j})^{1/2}(2^{-k+j+20\epsilon M_t})^{1/2}  \lesssim 2^{100\epsilon M_t+M_t/2} .\\
\end{split}
\]
 After combining  the above estimate, the estimate  \eqref{sep16eqn30}, and the rough estimate in   \eqref{sep16eqn21}, our desired estimate  \eqref{sep16eqn41}  holds. 

 \medskip

\noindent \textbf{Case 2.} \qquad If $m+l\leq -M_t-100\epsilon M_t$. 

\medskip 

From the obtained estimate  \eqref{sep16eqn21}, we have
\be\label{sep16eqn44}
\|B_{k;j,l}^{ess;m}(t,x)\|_{L^\infty}\lesssim 2^{2\tilde{\alpha}_t M_t-50\epsilon M_t}.
\ee
 Therefore,   our desired estimate  \eqref{sep16eqn42}  holds after combining the above estimate, the estimate  \eqref{nov26eqn31}  in Lemma \ref{erroresti}, and the estimate  \eqref{july10eqn89}  in Proposition \ref{Linfielec}. 
\end{proof}

 Lastly, we estimate the very high frequency case, i.e., $k\in [ 20M_t,\infty)\cap \Z$. For this case, we have 
\begin{lemma}\label{largefrel2}
For any $t\in [0, T)$, $k\in   [20M_t, \infty)\cap \Z$,  $j\in     \Z_+, l\in [-j, 2]\cap \Z$, $m\in [ -10 M_t  ,  \epsilon M_t]\cap \Z$,  we have 
 \be\label{july5eqn111}
  \|B^m_{k;j,l}(t,x)\big\|_{L^2}\lesssim  2^{-(1-3\epsilon)k+4j +20\epsilon M_t}.
 \ee
\end{lemma}
\begin{proof}

Recall  \eqref{july5eqn60}  and  the decomposition of $B^m_{j,l}(t,x)$ in \eqref{july5eqn31}. After localizing the size of $t-s$ further with threshold $-k$,   we have 
\be\label{july5eqn110}
\begin{split}
\|B^m_{k;j,l}(t,x)\big\|_{L^2}\lesssim  &\sup_{g\in L^2, \|g\|_{L^2}=1}  \big|J^{0 }_{m,k;j,l}(g)(t)\big| \\
& + \big[ \sum_{m'\in [-k,\epsilon M_t]\cap \Z } \big|J^{T;m'}_{m,k;j,l}(g)(t)\big| + \big|J^{S;m'}_{m,k;j,l}(g)(t)\big|\big], \\
\end{split}
\ee
where 
\[
\begin{split}
J^0_{k;j,l}(g)(t)& :=  \int_{\R^3}\overline{g(x)}B_{0;j,l}^m(t,x) d x ,\\
 J^{S;m'}_{m,k;j,l}(g)(t)&:= \int_{0}^{t }  \int_{\R^3} \int_{\R^3}  \int_{\mathbb{S}^2} (t-s)   e^{   i(t-s)\omega \cdot \xi} \overline{\hat{g}(\xi)} \varphi_k(\xi)   \varphi_{m;-10M_t }(t-s) \varphi_{m';-k}(t-s)  \\ 
  &\qquad \times \mathcal{F}( f(E+\hat{v}\times B)  )(s, \xi, v)\cdot \nabla_v \big(\frac{ \hat{v}\times \omega \varphi_{j,l}(v, \omega)}{1+\hat{v}\cdot \omega} \big)    d\omega d\xi d v d s, \\
J^{T;m'}_{m,k;j,l}(g)(t)&:= \int_{0}^{t } \int_{\R^3} \int_{\R^3}  \int_{\mathbb{S}^2}  e^{   i(t-s)\omega \cdot \xi} \overline{\hat{g}(\xi)}  \hat{f}(s, \xi, v)   \varphi_k(\xi)   \\
&\qquad \times   \omega^{m,b}_{j,l}(t-s,v,\omega) \varphi_{m';-k}(t-s)   d\omega d\xi d v d s. 
\end{split}
\]

Recall  \eqref{2022jan1eqn1}. We  first estimate $J^0_{k;j,l}(g)(t)$.   Similar to what we did in  \eqref{sep16eqn13}, from the stationary phase point of view, we know that the angle between $\xi$ and $\omega$ is localized around the size of $ c(m,k,l):= \max\{-m-k-l, -m/2-k/2\}+\epsilon^2 (M_t +k_{+})$.  More precisely,  for the error type term, we do integration by parts in $\omega$ many times, which gives us a similar estimate as in  \eqref{sep16eqn30}.  For the main part, we use the  volume of support of $\omega$ and the Cauchy-Schwarz inequality. As a result, we have 
\be\label{july5eqn100}
\begin{split}
 |J^0_{k;j,l}(g)(t)| & \lesssim  2^{-k+2\epsilon k + 2\epsilon M_t}2^{-2l} \| \int_{\R^3} f(0,x,v)\varphi_j(v)d  v \|_{L^2_x } + 2^{-10M_t-10k} \\
 & \lesssim 2^{-(1-2\epsilon)k  + 10\epsilon M_t}.\\
 \end{split}
\ee

Now, we estimate the $T$ part and the $S$ part. Note that,  from the explicit formula of the coefficient $\omega^{m,b}_{j;a}(t-s,v,\omega), a\in \{0,1,2,3\}, $ in  \eqref{july9eqn11}, we have
\[
|\omega^{m,b}_{j;a}(t-s,v,\omega)|\varphi_{l;-j}( \tilde{v} +\omega)\lesssim 2^{-l}. 
\]
 
Following the same stationary phase  type argument, from the volume of support of $\omega$ of the main part, the Cauchy-Schwarz inequality,   the conservation law  \eqref{conservationlaw}, and the estimate  \eqref{nov24eqn41}  if $| v_{\bot}|\geq 2^{(\alpha_t+\epsilon)M_t}$,  we have
\be\label{july5eqn101}
\begin{split}
|J^{S;m'}_{m,k;j,l}(g)(t)| + |J^{T;m'}_{m,k;j,l}(g)(t)| & \lesssim \sup_{s\in [0, t]} \sum_{U\in\{E, B\} }  2^{m'-k+2\epsilon k + 2\epsilon M_t}2^{-j-2l} \\
&\qquad \times \| \int_{\R^3} U(s,x) f(s,x,v) \varphi_j(v) d v  \|_{L^2_x}\\
 & + 2^{m'-l} 2^{-m'-k+2\epsilon k + 2\epsilon M_t}  \| \int_{\R^3} f(s,x,v)\varphi_j(v)d  v \|_{L^2_x }\\
 &\lesssim 2^{-(1-2\epsilon)k+4j +20\epsilon M_t}. \\
\end{split}
\ee
 Recall  \eqref{july5eqn110}, our desired estimate \eqref{july5eqn111} holds from the above estimates 
 \eqref{july5eqn100}  and   \eqref{july5eqn101}.
\end{proof}
\subsection{Proof of Position \ref{meanLinfest}}\label{dictomyproposition}

 Recall the decomposition in  \eqref{july5eqn60}  and   \eqref{july5eqn31}. We first prove the desired estimate.  Recall the index set  \eqref{indexsetL2} and  the G-S decomposition  \eqref{GSdecomloc}. From the volume of support of $v, \omega$, the estimate  \eqref{nov26eqn31}  in Lemma \ref{erroresti}, and  the estimate \eqref{july10eqn89}  in Proposition \ref{Linfielec}, the following estimate holds if $(m,k,j,l)\in\mathcal{S}(t)\cap   \mathcal{T}_0(t) $, 
\[
\begin{split}
\sum_{ (m,k,j,l)\in  \mathcal{S}(t)\cap   \mathcal{T}_0(t) }  \sup_{s\in [0,t]}\|B^{  {m}}_{  {k};  {j},    {l}}(s,\cdot)\|_{L^\infty}  &\lesssim \big( \sum_{j\geq (1+2\epsilon)M_t} 2^{-5j} \big) + 1 + 2^{20\epsilon M_t} \mathbf{1}_{m=-10 M_t,  j\leq (1+2\epsilon)M_t}  \\
 & \quad \times \big(2^{2m+(1+2\tilde{\alpha}_t )M_t  + 3j}+  2^{m+ 3j}\big)  \lesssim 1. \\
 \end{split}
\]  
Therefore,  the desired first estimate in   \eqref{aug4eqn10}  holds from the above obtained estimate, 
  the estimate \eqref{sep16eqn42}  in Lemma \ref{smallk} for the case $(m,k,j,l)\in\mathcal{S}(t)\cap   \mathcal{T}_1(t)   $ and the estimate  \eqref{sep17eqn11}  in Proposition \ref{bulkroughpoi} for the case $(m,k,j,l)\in\mathcal{S}(t)\cap   \mathcal{T}_2(t) .    $  
 
 Now we focus on the proof of  the desired second estimate in   \eqref{aug4eqn10}. Based on the possible size of $k$, we split into two cases as follows.

 \medskip

\noindent \textbf{Case 1.} \qquad  If $k\leq 20M_t. $

\medskip

From the obtained estimate \eqref{july5eqn66}  in Lemma \ref{firstL2} and the estimate  \eqref{sep17eqn11}  in Proposition \ref{bulkroughpoi},  the following estimate holds if $ m+k\geq -2(1+2\epsilon)l +\epsilon M_t, m+j+2l <  0$, 
\[
\begin{split}
&\big(\sup_{s\in [0,t]}\|B^{  {m}}_{  {k};  {j},    {l}}(s,\cdot)\|_{L^\infty}  \big)^{\frac{1}{2}  } \big( \sup_{s\in [0,t]}\|B^{ {m}}_{  {k};  {j},  {l}}(s,\cdot)\|_{L^2}\big)^{ \frac{1}{2}} \\
&\lesssim    \big( 2^{m+j+2l+ \epsilon M_t}+2^{-100 M_t}\big)^{1/2} \big( 2^{10\epsilon M_t}\big( 2^{ 2 \tilde{\alpha}_t  M_t }  + 2^{7M_t/6+\tilde{\alpha}_t M_t/4} \big)(1+2^{-m-j-2l})\big)^{1/2}. 
\end{split}
\]
After combining the above estimate and the  obtained estimate  \eqref{sep16eqn41}  in Lemma \ref{smallk}, we have
\be\label{sep19eqn2}
\begin{split}
\sum_{\begin{subarray}{c} 
 k\leq 20M_t, 
(m,k,j,l)\in \mathcal{S}_2(t)\\ 
 \end{subarray}}& \big(\sup_{s\in [0,t]}\|B^{  {m}}_{  {k};  {j},    {l}}(s,\cdot)\|_{L^\infty}  \big)^{\frac{1}{2}  } \big( \sup_{s\in [0,t]}\|B^{ {m}}_{  {k};  {j},  {l}}(s,\cdot)\|_{L^2}\big)^{ \frac{1}{2}}\\ 
 &\lesssim  2^{8\epsilon M_t}\big( 2^{\tilde{\alpha}_t M_t}   + 2^{7M_t/12+ \tilde{\alpha}_t M_t/8} \big). 
 \end{split}
\ee

 \medskip

\noindent \textbf{Case 2.} \qquad  If $k\geq 20M_t. $

\medskip

From the estimate  \eqref{july5eqn111}  in Lemma \ref{largefrel2},   the estimate  \eqref{july10eqn89}  in Proposition \ref{Linfielec}, and the estimate  \eqref{nov26eqn31}  in Lemma \ref{erroresti},  we have
\[
\begin{split}
\sum_{ (m,k,j,l)\in \mathcal{S}_2(t), k\geq 20M_t} & \big(\sup_{s\in [0,t]}\|B^{  {m}}_{  {k};  {j},    {l}}(s,\cdot)\|_{L^\infty}  \big)^{\frac{1}{2}  } \big( \sup_{s\in [0,t]}\|B^{ {m}}_{  {k};  {j},  {l}}(s,\cdot)\|_{L^2}\big)^{ \frac{1}{2}}\\
&\lesssim \sum_{k\geq20 M_t}2^{-k/3+5M_t}\lesssim 1. \\
\end{split}
\]
Hence finishing the desired estimate  \eqref{aug4eqn10}.

 \subsection{Proof of   Theorem \ref{maintheorem1part1}  and Theorem \ref{roughesttailpart}}\label{proofoftheoremrough}

 Our desired estimate  \eqref{maintheoremroughest} in   Theorem \ref{maintheorem1part1} follows directly from the obtained estimate \eqref{july10eqn89} in Proposition \ref{Linfielec}. Meanwhile, the desired estimate \eqref{pointwiseest} holds from the obtained estimate  \eqref{sep21eqn31}.

   Moreover, from  the decomposition of the electromagnetic field in  \eqref{july5eqn1}  and the estimate \eqref{nov26eqn31}  in Lemma \ref{erroresti}, we have 
\be\label{dec2eqn31}
\sum_{j\in \Z, j\geq (1+2\epsilon)M_t}\sum_{U\in \{E, B\}}\| U_j(t,x  )\|_{L^\infty_x} \lesssim \sum_{j\in \Z, j\geq (1+2\epsilon)M_t} 2^{-9j + (3+10\epsilon)M_t} \lesssim 2^{-5M_t}. 
\ee
Hence finishing the proof of Theorem \ref{roughesttailpart}. 

The essential notations employed in this section are systematically detailed in Table \ref{tablesection3}.
\begin{table}[H]
\centering
\resizebox{\columnwidth}{!}{%
\begin{tabular}{ |c|c|c|c| } 
 \hline
 Notation & Definition    & Remarks \\ 
 \hline
 $ \varphi^{\star}(\omega, \xi)$, $\star\in\{ess,err\}$ & \eqref{nov28eqn55} & Localizing the angle between $\omega$ and $\xi$ by \\
 $ \varphi_{m,k,l}^{ \star }(\omega , \xi)$& \eqref{sep16eqn2} & using the stationary phase analysis\\
 \hline
 $(w, z)$  & \eqref{jan18eqn11} & Cylindrical coordinates\\
 \hline
 $J_{(\theta, \phi)\longrightarrow( w, z) } $ &\eqref{march18eqn66} & The Jacobian of changing coordinates $(\theta, \phi) \longrightarrow (z,\omega)$\\
 \hline
  $ U_{j;l}^{m;p, q}(t,x)$ & \eqref{2026may2eqn11} & Further dyadic decomposition for the $T$-part based on the \\
  & & size of $(\theta, \phi)$ in the Kirchhoff's formulas\\
\hline
 $ I_{j;l}^{m;p, q}(t,x)$ & \eqref{july3eqn36} & Further dyadic decomposition for the $S_2$-part based on the \\
  & & size of $(\theta, \phi)$ in the Kirchhoff's formulas\\
\hline
$ \mathcal{S}_i(t), i\in\{1,2\}$  &  \eqref{indexsetL2}  & Sets of parameters for the dichotomy of the magnetic field\\
\hline
$ B_{k;j,l}^{\star;m}(t,x), \star\in\{ess,err\}$ & \eqref{sep16eqn1} & Further localization of $ B_{k;j,l}^{m}(t,x)$ based on\\
& &  the   angle between $\omega$ and $\xi$\\
\hline
\end{tabular}%
}
\caption{Essential notations in section \ref{roughestelectromagnet}.}\label{tablesection3}
\end{table}

    \section{ $L^\infty_x$-type estimates of the   localized acceleration force}\label{linfacceloc}

  The main goal of this section is to prove the part (i) in Theorem \ref{maintheoremellipitic} for the elliptic parts and   the part (i) in Theorem \ref{mainresultsfirstpart} for the hyperbolic parts. 

The elliptic components can be easily estimated, as detailed in Lemma \ref{ellplinf} in section \ref{Linityelliptic}.   To estimate the hyperbolic parts (given the decomposition in \eqref{oct7eqn1}), it suffices to estimate the localized acceleration force $T_{k,j;n}^{\mu,i}(\mathfrak{m}, E)(t,x, \zeta) + \hat{\zeta} \times T_{k,j;n}^{\mu,i}(\mathfrak{m}, B)(t,x, \zeta)$. This is done using two refined decompositions, detailed in Section \ref{2024Decrefineddec}.

\subsection{Refined decomposition of the localized electromagnetic field}\label{2024Decrefineddec}

As discussed in section \ref{mainingredientsfirstpart}, we use two methods to leverage the smoothing effect, rather than just one, to effectively utilize the double null structure. The primary objective of this subsection is to achieve a good decomposition of $T_{k,j;n}^{\mu}(\mathfrak{m}, U)(t,x, \zeta)$
  by taking advantage of both the smoothing effect and the double null structure, as summarized in Lemma \ref{locdeclemm}.

Recall \eqref{sep17eqn32} and the definition of cutoff functions $\varphi_{j,n}^i(v, \zeta)$ for $i\in\{0,1,2,3,4\},$ in \eqref{sep4eqn6}. For the rest of this paper, we refer the region with velocity far away from  the fixed direction determined by $\zeta$, i.e., determined by the $\varphi_{j,n}^i(v, \zeta), i\in\{0,1\}$, as \textbf{the large angle region}
and call the other region, i.e., determined by the $\varphi_{j,n}^i(v, \zeta), i\in \{ 2,3,4 \}$, as \textbf{the small angle region}.  

 We do   normal form transformation in the large angle region  and do the Glassey-Strauss type decomposition  in the region  in the small angle region. More precisely, the main result of using two ways  of exploiting the smoothing effect   is summarized as follows, 
\begin{lemma}\label{locdeclemm}
Let $T_{k,j;n}^{\mu,i}( \mathfrak{m}, U)(t,x, \zeta) $ be defined as in \eqref{sep17eqn32}. 
\begin{enumerate}
\item[(i)] For the large angle region, i.e., $i\in\{0,1\}$, we have
\be\label{2024oct14eqn41}
 T_{k,j;n}^{\mu,i}( \mathfrak{m}, U)(t,x, \zeta)= \sum_{\star\in\{ini, ell,bil\}} \widetilde{T}_{k,j;n}^{\star;\mu,i }( \mathfrak{m}, U)(t,x, \zeta) ,
 \ee
where $\widetilde{T}_{k,j;n}^{\star;\mu,i }( \mathfrak{m}, U)(t,x, \zeta), \star\in \{ini, ell,bil\}$ are defined in \eqref{oct1eqn1}. 
\item[(ii)] For the small angle region, i.e., $i\in\{2,3,4\},$ we have 
\be\label{sep18eqn50}
\begin{split}
T_{k,j;n }^{\mu,i}(\mathfrak{m}, U)(t,x, \zeta)&= \sum_{\begin{subarray}{c}
(l,r)\in \mathcal{B}_i\\
m\in [-10 M_t ,\epsilon M_t]\cap \Z\\
\end{subarray}}   T_{k,j;n,l,r}^{\mu,m, i  }(\mathfrak{m}, U)(t,x, \zeta),   \\ 
\forall ( l,r)\in \mathcal{B}_i, \quad T_{k,j;n;l,r}^{\mu,m,i}(\mathfrak{m}, U)(t,x, \zeta)&= \sum_{\star\in \{0, T,S\}} \widetilde{T}_{k,j;n;l,r}^{\star; \mu, m,i}(\mathfrak{m}, U)(t,x, \zeta),
\end{split}
\ee
where $T_{k,j;n,l,r}^{\mu,m, i  }(\mathfrak{m}, U)(t,x, \zeta)$ are defined in \eqref{sep18eqn31} and $ \widetilde{T}_{k,j;n;l,r}^{\star; \mu, m,i}(\mathfrak{m}, U)(t,x, \zeta), \star\in \{0, T,S\}$, are defined in \eqref{sep18eqn44} and  \eqref{sep19eqn61} and the index sets $\mathcal{B}_i, i\in\{2,3,4\},$ are defined as follows, 
\be\label{2024oct14eqn1}
\begin{split}
\mathcal{B}_2&:= \{(l,r): l\in[  \bar{l}_2 ,2]\cap \Z, r\in [  \bar{l}_2 , l ]\cap \Z\},\\
\mathcal{B}_3:=\mathcal{B}_4&:=\{(l,r): l\in[\bar{l}_3,2]\cap \Z, r=\bar{l}_3 \},\\
\bar{l}_2&:=\lfloor n + \epsilon M_t/2\rfloor, \bar{l}_3:=\bar{l}_4:=-j.\\
\end{split}
\ee
 
\end{enumerate}
\end{lemma}
\begin{remark}
Roughly speaking, (i) The ``$ini$''-part in \eqref{sep4eqn30} plays the same role of the ``$0$''-part in the G-S decomposition. They all only depends on the initial data. Therefore, they can safely be discarded; (ii) The ``$ell$''-part in \eqref{sep4eqn30} plays the same role of the ``$T$''-part in the G-S decomposition. They all only depends on $f$; (iii) The ``$bil$''-part in \eqref{sep4eqn30} plays the same role of the ``$S$''-part in the G-S decomposition. They are all bilinear with respect to $f$ and the electromagnetic field $(E, B).$ 
\end{remark}
\begin{remark}
There are two extra parameters, $l$ and $r$, in the decomposition \eqref{sep18eqn50}, which corresponds to an additional inhomogeneous dyadic decomposition for the size of $\angle(v, -\omega)$. More precisely, we use the partition of unity $1=\sum_{l\in[r,2]\cap \Z} \varphi_{l;r}(\tilde{v}+\omega) $, see \eqref{cutoffwiththreshold} for the definition of cutoff function  $\varphi_{l;r}(\cdot)$. Roughly speaking, $l$ measures the  size of $\angle(v, -\omega)$ and $r$ measures the threshold of the inhomogeneous dyadic decomposition.

\end{remark}
\begin{remark}

In the small angle region, we do not consistently apply the G-S decomposition (the second decomposition in \eqref{sep18eqn50}). We utilize it only when the frequency is large, as the smoothing effect outweighs the singular loss   $(1+\hat{v}\cdot\omega)^{-1}$. When the frequency is smaller, we revert to the original form of the nonlinearity, which incorporates only cutoff functions, as seen in the first decomposition in \eqref{sep18eqn50}.
 
\end{remark}
\subsection{Proof of Lemma \ref{locdeclemm}}

We proceed in steps as follows. 

\medskip

\noindent \textbf{Step 1.}\qquad The large angle case, i.e., $i\in\{0,1\}.$

\medskip

 Recall   \eqref{sep5eqn10}.  Note that, after using the Vlasov equation in  \eqref{mainequation}  to substitute $\p_t f$ and doing integration by parts in $v$ once,  we have
\[
\begin{split}
\mathcal{N}^{j,n}_{E;i}(t,x,\zeta) := - 4\pi\big[ \int_{\R^3}  &\big(- \hat{v} \hat{v}\cdot \nabla_x  f(t, x, v)   +  \nabla_x f(t, x , v)\big)\varphi_{j,n}^i(v, \zeta) \\
&  + f(t,x,v) (E+\hat{v}\times B)\cdot \nabla_v\big( \hat{v}\varphi_{j,n}^i(v, \zeta)\big) d v   \big].
\end{split}
\]
 
 Let $g(t,x,v):=f(t,x+t\hat{v},v)$ be the profile of  $f$. In terms of  $g$, the following equality holds on the Fourier side from the above equality,   \eqref{sep5eqn10}, and \eqref{sep5eqn66},   
\be\label{sep5eqn70}
\begin{split}
T_{k,j;n}^{\mu,i}( \mathfrak{m}, U)(t,x, \zeta):= \int_0^t & \int_{\R^3} \int_{\R^3} e^{i x\cdot \xi + i \mu(t-s)|\xi | - is\hat{v}\cdot \xi } |\xi|^{-1}   \hat{g}(s,\xi, v) \\
&\times \mathfrak{m}(\xi, \zeta) m^{0}_U(\xi, v)  \varphi_k(\xi)\varphi_{n;-M_t}\big(  \tilde{\xi} + \mu \tilde{\zeta}  \big)\varphi_{j,n}^i(v, \zeta) \\
& +   e^{i x\cdot \xi + i \mu(t-s)|\xi |  } |\xi|^{-1}   \mathcal{F}[(E+\hat{v}\times B) f](s,\xi, v)\cdot\nabla_v \big(\hat{v}\varphi_{j,n}^i(v, \zeta)\big) \\
 &\times \mathfrak{m}(\xi, \zeta) m^{1}_U(\xi, v)  \varphi_k(\xi)\varphi_{n;-M_t}\big(   \tilde{\xi} + \mu \tilde{\zeta}   \big)d \xi  d v d s,
\end{split} 
\ee
where
\be\label{sep5eqn96}
\begin{split}
m_{E}^{0}(\xi, v) & =-i4\pi (\hat{v}\hat{v}\cdot \xi -\xi ),\quad m_{B}^{0}(\xi, v) = i 4\pi\hat{v}\times \xi,\\
 m_{E}^{1}(\xi, v) & =4\pi \hat{v}, \quad m_{B}^{1}(\xi, v) =(0,0,0).\\
\end{split}
\ee
Due to the cutoff function $\varphi_{j,n}^i(v, \zeta), i\in\{0,1\}$, see \eqref{sep4eqn6},  we have $n\leq - \epsilon M_t+10$ and $|\tilde{\zeta}-\tilde{v}|\gtrsim  2^{  \vartheta^{\star}_0 }\gg 2^{n+\epsilon M_t/4}$, which implies further that  
 \be\label{2022feb8eqn51}
|\tilde{\xi} + \mu \tilde{v}|\gtrsim 2^{   \vartheta^{\star}_0 }, \quad |\xi| + \mu \hat{v}\cdot \xi \gtrsim 2^{k+ 2 \vartheta^{\star}_0}. 
 \ee

 To exploit high oscillation in ``$s$'' when the frequency is very large, we  do integration by parts in ``$s$'' once for the   first integral in  \eqref{sep5eqn70}. As a result, for any $ i\in\{0,1\}$, we have
 \be\label{sep4eqn22}
 \begin{split}
  &\int_0^t \int_{\R^3} \int_{\R^3} e^{i x\cdot \xi + i \mu(t-s)|\xi | - is\hat{v}\cdot \xi } |\xi|^{-1}  \hat{g}(s,\xi, v)  \mathfrak{m}(\xi, \zeta) m^0_U(\xi, v)  \varphi_k(\xi)\\
  &\qquad \times \varphi_{n;-M_t}\big(  \tilde{\xi} + \mu \tilde{\zeta}  \big)\varphi_{j,n}^i(v, \zeta) d \xi d v d s \\
  &  =  \int_{\R^3}\int_{\R^3} e^{i x\cdot \xi} e^{i\mu(t-s)|\xi|} 
\frac{   \mathfrak{m}(\xi, \zeta) m^0_U(\xi, v)  }{  i  |\xi|( \mu|\xi|+ \hat{v}\cdot \xi  ) }
 \hat{f}(s, \xi, v)  \varphi_{j,n}^i(v, \omega)\varphi_k(\xi)  \varphi_{n;-M_t}\big(  \tilde{\xi} + \mu \tilde{\zeta}   \big) d\xi d v  \Big|_{s=0}^t\\
 &+ \int_0^t \int_{\R^3}\int_{\R^3} e^{i x\cdot \xi+i\mu (t-s)|\xi|}\varphi_k(\xi)  \varphi_{n;-M_t}\big(   \tilde{\xi} + \mu \tilde{\zeta}   \big) \\
 &\qquad\times   \mathcal{F}\big((E+\hat{v}\times B)f\big)(s, \xi, v)  \cdot \nabla_v \big[
 \frac{  \mathfrak{m}(\xi, \zeta) m^0_U(\xi, v)  }{   i |\xi|(\mu |\xi|+ \hat{v}\cdot \xi  ) }
 \varphi_{j,n}^i (v, \omega)\big]  d\xi d v  d s.\\
 \end{split}
 \ee
 
 After combining the decomposition  \eqref{sep5eqn70}  and  \eqref{sep4eqn22}, we have 
 \be\label{sep4eqn30}
 T_{k,j;n}^{\mu,i}( \mathfrak{m}, U)(t,x, \zeta)= \sum_{\star\in\{ini, ell,bil\}} \widetilde{T}_{k,j;n}^{\star;\mu,i }( \mathfrak{m}, U)(t,x, \zeta) ,
 \ee
 where
  \be\label{oct1eqn1}
  \begin{split}
  \widetilde{T}_{k,j;n}^{ini;\mu ,i}( \mathfrak{m}, U)(t,x, \zeta) &= - \int_{\R^3}\int_{\R^3} e^{i x\cdot \xi} e^{i\mu t |\xi| } 
\frac{  \mathfrak{m}(\xi, \zeta) m^0_U(\xi, v)  }{  i |\xi|(\mu |\xi|+ \hat{v}\cdot \xi  )  }
 \hat{f}(0, \xi, v)   \\
 &\qquad \times  \varphi_{j,n}^i(v, \zeta)\varphi_k(\xi)   \varphi_{n;-M_t}\big(  \tilde{\xi} + \mu \tilde{\zeta}  \big) d\xi d v,  \\
  \widetilde{T}_{k,j;n}^{ell;\mu ,i}( \mathfrak{m}, U)(t,x, \zeta) &=  \int_{\R^3}\int_{\R^3} e^{i x\cdot \xi}  \frac{  \mathfrak{m}(\xi, \zeta) m^0_U(\xi, v)  }{   i |\xi|(\mu |\xi|+ \hat{v}\cdot \xi  ) }
 \hat{f}(t, \xi, v) \\
 & \qquad \times  \varphi_{j,n}^i(v, \zeta)\varphi_k(\xi)   \varphi_{n;-M_t}\big(   \tilde{\xi} + \mu \tilde{\zeta}   \big) d\xi d v ,  \\
  \widetilde{T}_{k,j;n}^{bil;\mu,i }( \mathfrak{m}, U)(t,x, \zeta)&= \int_0^t \int_{\R^3}\int_{\R^3} e^{i x\cdot \xi + i\mu (t-s) |\xi| }   |\xi|^{-1} \varphi_{n;-M_t}\big(   \tilde{\xi} + \mu \tilde{\zeta}   \big)  \varphi_k(\xi)\\
  &  \qquad  \times  \mathcal{F}\big((E+\hat{v}\times B)f\big)(s, \xi, v)  \cdot \nabla_v \big[ \big(
 \frac{ m^0_U(\xi, v)  }{ i ( \mu |\xi|+ \hat{v}\cdot \xi  ) }\\
&  \qquad  +   m^{1}_U(\xi,v) \big) \mathfrak{m}(\xi, \zeta)
 \varphi_{j,n}^i (v, \zeta)\big]   d\xi d v  d s.
 \end{split}
\ee

Moreover,  we can represent the elliptic part  $\widetilde{T}_{k,j;n}^{ell;\mu,i }( \mathfrak{m}, U)(t,x, \zeta)$  in terms of kernel  as follows, 
\be\label{sep5eqn86}
\widetilde{T}_{k,j;n}^{ell;\mu,i }(\mathfrak{m}, U)(t,x, \zeta)=\int_{\R^3} \int_{\R^3} \mathcal{K}^{ell;U,i }_{k,j,n}(\mathfrak{m})(y, v,  \zeta) f(t, x-y, v) d y d v,
\ee 
where the kernel   $  \mathcal{K}^{ell;U,i}_{k,j,n}(y, v,   \zeta)$ is defined as follows,  
\be\label{sep5eqn88}
\mathcal{K}^{ell;U ,i}_{k,j,n}(\mathfrak{m})(y, v,   \zeta):= \int_{\R^3} e^{i y\cdot \xi}   
  \frac{ \mathfrak{m}(\xi, \zeta)  m^0_U(\xi, v)  }{  i |\xi| ( \mu|\xi|+ \hat{v}\cdot \xi  ) } 
 \varphi_{j,n}^i (v, \zeta) \varphi_{n;-M_t}\big(  \tilde{\xi} + \mu \tilde{\zeta}  \big)  \varphi_k(\xi)  d \xi.
\ee
From the Kirchhoff's formulas in Lemma \ref{Kirchhoff}, we can represent the bilinear part $\widetilde{T}_{k,j;n}^{bil;\mu ,i}( \mathfrak{m}, U)(t,x, \zeta)$ in terms of kernels  as follows, 
\be\label{sep5eqn87}
\begin{split}
\widetilde{T}_{k,j;n}^{bil;\mu,i }(  \mathfrak{m}, U)(t,x, \zeta)= \int_{0}^{t}  \int_{\R^3}& \int_{\R^3} \int_{\mathbb{S}^2}  \big(E(s, x-y+(t-s)\omega)+ \hat{v}\times B(s, x-y+(t-s)\omega) \big)  \\
&\cdot \nabla_v \big( (t-s)\mathcal{K}^{U ,i}_{k,j,n}(\mathfrak{m})(y, v, \omega, \zeta) + \mathcal{K}^{err;U ,i}_{k,j,n}(\mathfrak{m})(y, v,  \zeta)\big)\\
& \times  f(s,x-y+(t-s)\omega, v) d \omega d y d v ds,\\
\end{split}
\ee
where the kernels    $\mathcal{K}^{U ,i}_{k,j,n}(y, v, \omega, \zeta)$ and $  \mathcal{K}^{err;U,i }_{k,j,n}(y, v,   \zeta)$ are defined as follows, 
\be\label{sep5eqn85}
\begin{split}
\mathcal{K}^{U,i}_{k,j,n}(\mathfrak{m})(y, v, \omega, \zeta):= \int_{\R^3}& e^{i y\cdot \xi} i(    {\omega \cdot \tilde{\xi} }  +   \mu  )  \varphi_{n;-M_t}\big(  \tilde{\xi} + \mu \tilde{\zeta}  \big)  \varphi_k(\xi)  \\
&\times  \big[ \big(
  \frac{m^0_U(\xi, v)  }{   i ( \mu|\xi|+ \hat{v}\cdot \xi  ) } +      m^{1}_U(\xi,v) \big) \mathfrak{m}(\xi, \zeta) 
 \varphi_{j,n}^i (v, \zeta)\big] d \xi,\\
 \mathcal{K}^{err;U,i}_{k,j,n}(\mathfrak{m})(y, v,   \zeta)  := \int_{\R^3} & e^{i y\cdot \xi}   |\xi|^{-1}   \varphi_{n;-M_t}\big(   \tilde{\xi} + \mu \tilde{\zeta} \big)  \varphi_k(\xi) \\
 &\times \big[ \big(
  \frac{m^0_U(\xi, v)  }{   i ( \mu|\xi|+ \hat{v}\cdot \xi  ) } +      m^{1}_U(\xi,v) \big)  \mathfrak{m}(\xi, \zeta) 
 \varphi_{j,n}^i (v, \zeta)\big]  d \xi. \\
 \end{split}
\ee

\medskip

\noindent \textbf{Step 2.}\qquad The small angle case, i.e., $i\in\{2,3,4\}.$

\medskip

We first start with the Kirchhoff's formula  in Lemma \ref{Kirchhoff}. Note that, after first doing dyadic decompositions for the angle between $\omega$ and $-v$  and then using the partition of the characteristic function $\mathbf{1}_{[0,t]}$ in  \eqref{july11eqn33},  for any $U\in\{E, B\}, i\in\{  2,3,4  \}, $ we have
 \be\label{sep17eqn61}
\int_0^t \d^{-1} \sin((t-s)\d)  \mathcal{N}^{j,n}_{U;i}(s,x, \zeta)  d s =    \sum_{\begin{subarray}{c}
(l,r)\in \mathcal{B}_i\\
m\in [-10 M_t ,\epsilon M_t]\cap \Z\\  
\end{subarray}}  - 4\pi\big[  H^{U;m,i }_{j,n,l,r}(t, x, \zeta )\big],
\ee  
 where
  \be\label{sep17eqn63}
  \begin{split}
  H^{E;m,i}_{j,n,l,r}(t, x, \zeta )&:= \int_0^t    \int_{\R^3} \int_{\mathbb{S}^2} (t-s)  \big( \hat{v} \p_s    +  \nabla_x \big) f(s, x+(t-s)\omega , v)\\
  &\quad \times \varphi_{m;-10M_t }(t-s) \varphi_{j,n}^{i; r}(v, \zeta)  \varphi_{l; r}(\tilde{v}+\omega )  d\omega d v ds , \\
H^{B;m,i}_{j,n,l,r}(t, x, \zeta )&:=    \int_0^t    \int_{\R^3} \int_{\mathbb{S}^2} (t-s)\big( \hat{v} \times   \nabla_x f(s, x+(t-s)\omega  , v)\big)\varphi_{j,n}^{i; r}(v, \zeta)   \\
&\quad \times   \varphi_{l; r}(\tilde{v}+\omega )  \varphi_{m;-10M_t }(t-s)   d\omega d v d s,\\
\varphi_{j,n}^{2; r}(v, \zeta)&:=  \varphi^2_{j,n}(v, \zeta)\varphi_{r;\bar{l}_2}(\tilde{v}-\tilde{\zeta}), \quad \varphi_{j,n}^{3; r}(v, \zeta):=  \varphi^3_{j,n}(v, \zeta), \quad \varphi_{j,n}^{4; r}(v, \zeta):=  \varphi^4_{j,n}(v, \zeta). \\
\end{split}
 \ee

   After taking ``$\p_t$'' derivative on both hand sides of \eqref{sep17eqn61} and then using the partition of the characteristic function $\mathbf{1}_{[0,t]}$ in  \eqref{july11eqn33}, we have
\be\label{sep17eqn71}
\int_0^t   \cos((t-s)\d)  \mathcal{N}^{j,n}_{U;i}(s,x, \zeta)  d s =    \sum_{\begin{subarray}{c}
(l,r)\in \mathcal{B}_i\\
m\in [-10 M_t ,\epsilon M_t]\cap \Z\\  
\end{subarray}} - 4\pi\big[ \widetilde{H}^{U;m,i}_{j,n,l,r}(t, x, \zeta )\big],
\ee
where  
\be\label{sep19eqn49}
\begin{split}
 \widetilde{H}^{E;m,i}_{j,n,l,r}(t, x, \zeta )=   \int_0^t    \int_{\R^3} \int_{\mathbb{S}^2}&  \big( 1+(t-s)\omega \cdot \nabla \big)  \big( \hat{v} \p_s f(s, x+(t-s)\omega , v)   \\
 &+  \nabla_x f(s, x+(t-s)\omega  , v)\big)\varphi_{j,n}^{i; r}(v, \zeta)\\
& \times    \varphi_{l; r}(\tilde{v}+\omega ) \varphi_{m;-10M_t }(t-s)  d\omega d v ds  ,\\
\widetilde{H}^{B;m,i}_{j,n,l}(t, x, \zeta )=  \int_0^t    \int_{\R^3} \int_{\mathbb{S}^2} & \big( 1+(t-s)\omega \cdot \nabla \big)  \big( \hat{v} \times \nabla_x f(s, x+(t-s)\omega , v)   \big) \\
&\times  \varphi_{j,n}^{i; r}(v, \zeta)  \varphi_{l; r}(\tilde{v}+\omega )   \varphi_{m;-10M_t }(t-s)  d\omega d v ds. 
\end{split}
\ee
 After combining  \eqref{sep17eqn61}  and  \eqref{sep17eqn71}, $ \forall U\in \{E, B\},$ from   \eqref{sep17eqn32}, we have 
\be\label{nov5eqn10}
\begin{split}
T_{k,j;n}^{\mu,i}( \mathfrak{m}, U)(t,x, \zeta) & =   \int_0^t  T^{\mu}_{k,n}\big(\d^{-1}  \cos((t-s)\d)  \mathcal{N}^{j,n}_{U;i}\big)(s,x, \zeta) \\
&\qquad  + i\mu  T^{\mu}_{k,n}\big( \d^{-1}  \sin((t-s)\d)  \mathcal{N}^{j,n}_{U;i}\big)(s,x, \zeta)  d s\\
&  =  \sum_{\begin{subarray}{c}
(l,r)\in \mathcal{B}_i\\
m\in [-10 M_t ,\epsilon M_t]\cap \Z\\  
\end{subarray}}       - 4\pi\big[  T_{k,j;n,l,r}^{\mu,m,i}(\mathfrak{m}, U)(t,x, \zeta)\big], \\
\end{split}
\ee
where  the symbol of the Fourier multiplier operator $T^{\mu}_{k,n}$ is $ \mathfrak{m}(\xi, \zeta)  \varphi_k(\xi)\varphi_{n;-M_t}\big( \tilde{\xi} + \mu \tilde{\zeta} \big)$ and
\be\label{sep19eqn48}
 T_{k,j;n,l,r}^{\mu,m,i}(\mathfrak{m}, U)(t,x, \zeta) = T^{\mu}_{k,n} H^{U;m,i}_{j,n,l,r}(t, x, \zeta )+ \d^{-1}  T^{\mu}_{k,n} \widetilde{H}^{U;m,i}_{j,n,l,r}(t, x, \zeta ).
\ee 
 In terms of kernels, from    \eqref{sep17eqn63}  and  \eqref{sep19eqn49},  we have 
 \be\label{sep18eqn31}
\begin{split}
T_{k,j;n,l,r}^{\mu,m,i}(\mathfrak{m},B)(t,x, \zeta)&= \int_{0}^t \int_{\R^3} \int_{\R^3} \int_{\mathbb{S}^2} \big( (t-s) \mathfrak{K}^{\mu, B}_{k;n}(y,\omega, v, \zeta) +\mathfrak{K}^{err;\mu,B}_{k;n}(y, v, \zeta)\big)\\
&\qquad \times f(s, x-y+(t-s)\omega, v) \varphi_{j,n}^{i; r}(v, \zeta)  \varphi_{l; r}(\tilde{v}+\omega )   \\
 &\qquad \times   \varphi_{m;-10M_t }(t-s) d\omega dy d v ds, \\
 T_{k,j;n,l,r}^{\mu,m,i}(\mathfrak{m},E)(t,x, \zeta)&= \int_{0}^t \int_{\R^3} \int_{\R^3} \int_{\mathbb{S}^2} \big( (t-s) \mathfrak{K}^{\mu,E}_{k;n}(y,\omega, v, \zeta)+\mathfrak{K}^{err;\mu,E}_{k;n}(y, v, \zeta)\big) \\
 &\qquad\times f(s, x-y+(t-s)\omega, v)  \varphi_{j,n}^{i; r}(v, \zeta)  \varphi_{l; r}(\tilde{v}+\omega ) \varphi_{m;-10M_t }(t-s) \\
 &  \qquad       +    (E(s,  x-y+(t-s)\omega)+\hat{v}\times B(s,  x-y+(t-s)\omega))\\
 &\qquad \cdot  \nabla_v \big( (t-s)   \mathfrak{H}^{\mu,E,i}_{k,j;n,l,r}(y,\omega, v, \zeta)  +\mathfrak{H}^{err;\mu,E,i}_{k,j;n,l,r}(y, v, \zeta)\big)\\
 &\qquad\times  f(s, x-y+(t-s)\omega, v)\varphi_{m;-10M_t }(t-s)  d \omega dy d v ds,
 \end{split}
\ee
where the kernels appeared above are defined as follows, 
\be\label{sep19eqn85}
\begin{split}
 \mathfrak{K}^{\mu, U}_{k;n}( \mathfrak{m})(y,\omega, v, \zeta)&:=\int_{\R^3 } e^{i y\cdot \xi }  \mathfrak{m}(\xi, \zeta) \varphi_k(\xi)\varphi_{n;-M_t}\big(   \tilde{\xi} + \mu \tilde{\zeta}   \big)   m_{U}^{0}(\xi, v)(i \omega \cdot \tilde{\xi} +i \mu)d \xi, \\
\mathfrak{K}^{err;\mu, U}_{k;n}( \mathfrak{m})(y,\omega, v, \zeta)&:=i \int_{\R^3 } e^{i y\cdot \xi }  \mathfrak{m}(\xi, \zeta) \varphi_k(\xi)\varphi_{n;-M_t}\big(    \tilde{\xi} + \mu \tilde{\zeta}  \big)  m_{U}^{0}(\xi, v)|\xi|^{-1} d \xi,\\ 
  \mathfrak{H}^{\mu, U ,i}_{k,j;n,l,r}( \mathfrak{m})(y,\omega, v, \zeta)&:=\int_{\R^3 } e^{i y\cdot \xi }  \mathfrak{m}(\xi, \zeta) \varphi_k(\xi)\varphi_{n;-M_t}\big(    \tilde{\xi} + \mu \tilde{\zeta}   \big)  m_{U}^{1}(\xi, v)   \\
  &\qquad \times (i \omega \cdot   \tilde{\xi} -\mu)    \varphi_{j,n}^{i; r}(v, \zeta)  \varphi_{l; r}(\tilde{v}+\omega )  d \xi, \\ 
\mathfrak{H}^{err;\mu, U,i}_{k,j;n,l,r}( \mathfrak{m})(y,\omega, v, \zeta)&:=i \int_{\R^3 } e^{i y\cdot \xi }  \mathfrak{m}(\xi, \zeta) \varphi_k(\xi)\varphi_{n;-M_t}\big(   \tilde{\xi} + \mu \tilde{\zeta} \big)  m_{U}^{1}(\xi, v)  \\
&\qquad \times  |\xi|^{-1}  \varphi_{j,n}^{i; r}(v, \zeta)  \varphi_{l; r}(\tilde{v}+\omega )  d \xi. 
\end{split}
\ee

The obtained formulas in  \eqref{sep18eqn31}   will be used in the low frequency regime, i.e., when $k$ is relatively small.  For the high frequency regime, as in Lemma \ref{GSdecomloc},  we can also do G-S type decomposition for $ T_{k,j;n,l}^{\mu,m,i}(\mathfrak{m}, U)(t,x, \zeta), U\in \{E, B\}.$  To this end,   from   \eqref{sep19eqn48}, it suffices to do   G-S type decomposition for $H^{U;m,i}_{j,n,l,r}(t, x, \zeta )$ and  $ \widetilde{H}^{U;m,i}_{j,n,l,r}(t, x, \zeta )$.

We first do the  G-S type decomposition for $H^{U;m,i}_{j,n,l}(t, x, \zeta ). $  Note that,    the only difference between $H^{U;m,i}_{j,n,l}(t, x, \zeta ) $ and the localized  G-S   decomposition in Lemma \ref{GSdecomloc} comes from  by the new introduced cutoff functions $  \varphi_{j,n}^{i; r}(v, \zeta), i\in\{ 2,3,4 \}  $, which doesn't depend on  the angular variable $\omega$. Therefore, with minor modifications in the proof of Lemma \ref{GSdecomloc}, 
  for any fixed $j\in \Z_+, l\in [ \bar{l}_i ,  2]\cap \Z, $  the following Glassey-Strauss decomposition holds for  any $H^{U;m,i}_{j,n,l,r}(t, x, \zeta ), U\in \{E, B\},  ( l,r)\in \mathcal{B}_i,$
\be \label{sep19eqn11}
  H^{U;m,i}_{j,n,l,r}(t, x, \zeta )=  \sum_{ \star\in \{0, T,S\} } H^{\star;U;m,i}_{j,n,l,r}(t, x, \zeta ), 
\ee
where
  \be\label{sep17eqn64}
  \begin{split}
  H^{S;U;m,i}_{j,n,l,r}(t, x, \zeta )= \int_0^t  \int_{\R^3} \int_{\mathbb{S}^2} &  (t-s)  \big(E(s,x+(t-s)\omega) +\hat{v}\times  B(s,x+(t-s)\omega) \big)\\
  &\cdot \nabla_v \big(\frac{  m_{U}(v, \omega )    \varphi_{j,n}^{i; r}(v, \zeta)  \varphi_{l; r}(\tilde{v}+\omega )   }{1+\hat{v}\cdot \omega}    \big)\\
& \times  f(s,x+(t-s)\omega, v)  \varphi_{m;-10M_t }(t-s) d \omega d v d s, \\
\end{split}
\ee 
  \be\label{sep17eqn65}
    \begin{split}
 H^{T;U;m,i}_{j,n,l,r}(t, x, \zeta )= \sum_{  j'\in[j-2,j+2]\cap \Z } \int_0^t \int_{\R^3} \int_{\mathbb{S}^2} &    \omega^{m;U}_{j,l,r}(t-s,v,\omega)  \varphi_{ j' ,n}^{i; r}(v, \zeta)\\
 & \times    f(s,x+(t-s)\omega, v)   d \omega d v ds ,
 \end{split}
  \ee 
  and $   H^{0;U;m,i}_{j,n,l,r}(t, x, \zeta )=0$ if $t\notin supp(\varphi_{m;-10M_t}(\cdot)),$ while $t\in supp(\varphi_{m;-10M_t}(\cdot)),$ we have 
\be\label{sep17eqn66}
   \begin{split}
   H^{0;U;m,i}_{j,n,l,r}(t, x, \zeta )= \int_{\R^3} \int_{\mathbb{S}^2} &t f(0,x+t\omega, v)\frac{ m_{U}(v, \omega)}{1+\hat{v}\cdot \omega}     \varphi_{l; r}(\tilde{v}+\omega )   \varphi_{j,n}^{i; r}(v, \zeta)  d\omega d v.
   \end{split}
   \ee
The symbols $m_{U}(v, \omega), U\in \{E, B\}$ appeared above are defined in  \eqref{july9eqn11} and $\omega^{m;U}_{j,l,r}(t-s,v,\omega) $ are defined as follows, 
 \be\label{2022feb16eqn4}
\omega^{m;U}_{j,l,r}(t-s,v,\omega)  = \omega^{m,U}_{j;0 }(t-s,v,\omega)   \varphi_{l; r}(\tilde{v}+\omega )    +\sum_{ q=1,2,3}  \omega^{m,U }_{j;q }(t-s,v,\omega) \p_{x_q} \varphi_{l; r}(\tilde{v}+\omega ),
\ee
where $\omega^{m,U}_{j ;q}(t-s,v,\omega), q\in \{0,1,2,3\},$ are defined in   \eqref{july9eqn11}.

Lastly,  we   do the  G-S type decomposition for   $\widetilde{H}^{U;m,i}_{j,n,l}(t, x, \zeta ) $ in  \eqref{sep19eqn48}.    The main difference for   this case comes from the new introduced coefficient ``$\omega$'' in  \eqref{sep19eqn49}. With minor modification of the argument used in    \eqref{july5eqn31}  in Lemma \ref{GSdecomloc},  we 
have 
\be\label{sep19eqn33}
\forall U\in \{E, B\}, ( l,r)\in \mathcal{B}_i, \quad  \widetilde{H}^{U;m,i}_{j,n,l,r}(t, x, \zeta ) = \sum_{\star\in\{0, T, S\}} \widetilde{H}^{\star;U;m,i}_{j,n,l,r}(t, x, \zeta ),
\ee
where 
\[
\begin{split}
\widetilde{H}^{S; U;m,i}_{j,n,l,r}(t, x, \zeta )  & =  \int_0^t  \int_{\R^3} \int_{\mathbb{S}^2}   \big( 1+(t-s)\omega \cdot \nabla \big)\Big[\big(E(s,x+(t-s)\omega) \\
&\qquad+\hat{v}\times  B(s,x+(t-s)\omega) \big)  \cdot \nabla_v \big(\frac{  m_{U}(v, \omega )     \varphi_{j,n}^{i; r}(v, \zeta)  \varphi_{l; r}(\tilde{v}+\omega )   }{1+\hat{v}\cdot \omega}  \big)  \\
&\qquad   \times  f(s,x+(t-s)\omega, v) \Big]  \varphi_{m;-10M_t }(t-s) d \omega d v d s, \\
\widetilde{H}^{T; U;m,i}_{j,n,l,r}(t, x, \zeta ) & =\sum_{   j' \in [j-2, j+2]}  \int_0^t  \int_{\R^3} \int_{\mathbb{S}^2} (t-s)^{-1} \big[ (t-s) c^{m,U }_{j,l,r;q}(t-s,v,\omega)  \\
 &  \qquad \times \p_{x_q}f(s, x+(t-s)\omega, v)  +   c^{m,U }_{j,l,r;err}(t-s,v,\omega)  \\
 &\qquad \times     f(s, x+(t-s)\omega, v)\big]  \varphi_{   j' ,n}^{i; r}(v, \zeta)    d \omega d v d s, 
\end{split}
 \]
   and $ \widetilde{H}^{0; U;m,i}_{j,n,l,r}(t, x, \zeta )=0$ if $t\notin supp(\varphi_{m;-10M_t}(\cdot)),$ while $t\in supp(\varphi_{m;-10M_t}(\cdot)),$ we have 
\be\label{2021jan2eqn3}
\widetilde{H}^{0; U;m,i}_{j,n,l,r}(t, x, \zeta )  =  
     {\int_{\R^3} \int_{\mathbb{S}^2}    \frac{  m_{U}(v, \omega )}{1+\hat{v}\cdot \omega} (1+ t\omega \cdot \nabla )f(0 ,x+t\omega, v )     \varphi_{j,n}^{i; r}(v, \zeta)     \varphi_{l; r}(\tilde{v}+\omega )  d\omega d v    }.   
\ee
The coefficients $ m_{U}(v, \omega ), U\in \{E, B\}$,   are defined in  \eqref{july9eqn11} and  the coefficients $c^{m,U }_{j,l,r;q}(t-s,v,\omega)$  and $  c^{m,U }_{j,l,r;err}(t-s,v,\omega), U\in \{E, B\}$, appeared above are given as follows, 
\be\label{sep20eqn42}
\begin{split}
c^{m,E  }_{j,l,r;q}(t-s,v,\omega) :=   \omega_q &\big[ \omega^{m,E}_{j;0 }(t-s,v,\omega)   \varphi_{l; r}(\tilde{v}+\omega )    \\
&\quad  +\sum_{ q=1,2,3}  \omega^{m,E }_{j;q }(t-s,v,\omega) \p_{x_q} \varphi_{l; r}(\tilde{v}+\omega )   \big] \\
 &\quad-  \big[(\delta_{1q},\delta_{2q},\delta_{3q} )-\frac{(\omega +\hat{v} )\hat{v}_q}{1+\hat{v}\cdot \omega}  -\omega_q\big( \omega     - \frac{(\omega +\hat{v} )\hat{v}\cdot \omega}{1+\hat{v}\cdot \omega}\big)    \big] \\
 &\quad \times    \varphi_{l; r}(\tilde{v}+\omega )   \varphi_j(v)  \varphi_{m;-10 M_t }(t-s),
 \end{split}
\ee
\be\label{sep20eqn44}
\begin{split}
c^{m,B  }_{j,l,r;q}(t-s,v,\omega):=     \omega_q  & \big[ \omega^{m,B}_{j;0 }(t-s,v,\omega)  \varphi_{l; r}(\tilde{v}+\omega )   \\
&\quad  +\sum_{ q=1,2,3}  \omega^{m,B }_{j;q }(t-s,v,\omega) \p_{x_q}  \varphi_{l; r}(\tilde{v}+\omega )   \big] \\
 &\quad+     \big(\hat{v}\times (\delta_{1q},\delta_{2q},\delta_{3q} )-\frac{(\hat{v}\times \omega)  (\hat{v}_q+\omega_q)}{1+\hat{v}\cdot \omega}\big) \\
 &\quad\times  \varphi_{l; r}(\tilde{v}+\omega )   \varphi_j(v)  \varphi_{m;-10 M_t }(t-s),\\
 c^{m, E }_{j,l,r;err}(t-s,v,\omega):=  &\big(\frac{\hat{v} +\omega  }{(1+\hat{v}\cdot \omega)^2}\big(1 -|\hat{v}^2|   \big) + \omega   - \frac{(\omega +\hat{v} )\hat{v}\cdot \omega}{1+\hat{v}\cdot \omega}\big) \\
&\quad \times    \varphi_{l; r}(\tilde{v}+\omega )   \varphi_j(v)   \varphi_{m;-10 M_t }( t-s )\\
&\quad - \sum_{q=1,2,3}       \big( (\delta_{1q}, \delta_{2q}, \delta_{3q})  + w_q \hat{v} -   \frac{(\omega +\hat{v} )(\omega_q + \hat{v}_q) }{1+\hat{v}\cdot \omega}\big) \\
&\quad  \times \varphi_j(v)\p_{x_q}    \varphi_{l; r}(\tilde{v}+\omega ) \varphi_{m;-10 M_t }( t-s ) - (t-s)    \varphi_{l; r}(\tilde{v}+\omega )   \\
&\quad   \times   \varphi_j(v)\p_x \varphi_{m;-10 M_t }(t-s)\big(\omega -\frac{(\omega +\hat{v} )\hat{v}\cdot\omega }{1+\hat{v}\cdot \omega}\big),\\
c^{m,B }_{j,l,r;err}(t-s,v,\omega):= &\big[   \big(\frac{\hat{v} \times \omega  }{(1+\hat{v}\cdot \omega)^2}\big(1 -|\hat{v}^2|   \big) +  \frac{\hat{v}\times \omega}{1+\hat{v}\cdot \omega}\big)      \varphi_{l; r}(\tilde{v}+\omega )  \varphi_j(v) \\
&\quad    - \sum_{q=1,2,3}    \big( \hat{v}\times (\delta_{1q}, \delta_{2q}, \delta_{3q})   -   \frac{ \hat{v}\times \omega (\omega_q + \hat{v}_q) }{1+\hat{v}\cdot \omega}\big)   \\
  &\quad \times \varphi_j(v)  \p_{x_q} \varphi_{l; r}(\tilde{v}+\omega )\big]\varphi_{m;-10 M_t }(t-s) \\
  &\quad  - (t-s)   \varphi_{l; r}(\tilde{v}+\omega )    \varphi_j(v)  \p_x \varphi_{m;-10 M_t }(t-s) \frac{\hat{v}\times \omega}{1+\hat{v}\cdot \omega},\\
 \end{split} 
\ee
 where the coefficients $\omega^{m,U}_{j,l;i}(t-s,v,\omega), i\in\{0,1,2,3\}, U\in\{E,B\},$ are defined in \eqref{july9eqn11}. 

To sum up,  after combining  the decomposition in \eqref{sep19eqn11}  and  \eqref{sep19eqn33},  from  \eqref{sep19eqn48},  we have
\be\label{sep19eqn60}
T_{k,j;n,l,r}^{\mu,m, i}(\mathfrak{m}, U)(t,x, \zeta)= \sum_{\star\in \{0, T,S\}} \widetilde{T}_{k,j;n,l,r}^{\star;  \mu,m, i }(\mathfrak{m}, U)(t,x, \zeta),
\ee
where
\be\label{sep18eqn44}
\begin{split}
\widetilde{T}_{k,j;n,l,r}^{S; \mu ,m, i}(\mathfrak{m}, U)(t,x, \zeta)&=  \int_0^t  \int_{\R^3} \int_{\R^3} \int_{\mathbb{S}^2} \big[ (t-s)  K_{k;n}^{ \mu}(\mathfrak{m})(y,  \zeta, \omega)+ \widetilde{K}_{k;n}^{ \mu}(\mathfrak{m})(y,   \zeta )\big] \\
&\quad \times    \big(E(s,x-y+(t-s)\omega)+\hat{v}\times  B(s,x-y+(t-s)\omega) \big) \\
 &\quad \cdot \nabla_v \big(\frac{  m_{U}(v, \omega )      \varphi_{j,n}^{i; r}(v, \zeta)     \varphi_{l; r}(\tilde{v}+\omega )      }{1+\hat{v}\cdot \omega}   \big)\\
 &\quad \times f(s,x-y+(t-s)\omega, v) \varphi_{m;-10M_t}(t-s)  d \omega dy  d v d s, \\
  \widetilde{T}_{k,j;n,l,r}^{T; \mu,m , i}(\mathfrak{m}, U)(t,x, \zeta)=&\sum_{\tilde{j}\in [j-2, j+2]}   \int_0^t \int_{\R^3} \int_{\R^3} \int_{\mathbb{S}^2} \big[ i\mu \omega^{m;U}_{j,l,r}(t-s,v,\omega)  K_{k;n}^{ }(\mathfrak{m})(y,  \zeta ) \\
  &\quad + c^{m,U}_{j,l,r;q}(t-s,v,\omega) K_{k;n}^{ q }(\mathfrak{m})(y, \zeta  )\\
  & \quad + (t-s)^{-1}   c^{m,U}_{j,l,r;err}(t-s,v,\omega)   \tilde{K}_{k;n}^{ \mu}(\mathfrak{m})(y,   \zeta ) \big] \\
  &\quad \times  f(s,x-y+(t-s)\omega, v)   \varphi_{j,n}^{i; r}(v, \zeta)   d \omega dy  d v ds, \\
 \end{split}
\ee
and $ \widetilde{T}_{k,j;n,l,r}^{0; \mu,m, i}(\mathfrak{m}, U)(t,x, \zeta)=0$ if   $t\notin supp(\varphi_{m;-10M_t}(\cdot)),$ while $t\in supp(\varphi_{m;-10M_t}(\cdot)),$ we have 
\be\label{sep19eqn61}
\begin{split}
 \widetilde{T}_{k,j;n,l,r}^{0; \mu,m, i}(\mathfrak{m}, U)(t,x, \zeta) =   \int_{\R^3} \int_{\R^3} \int_{\mathbb{S}^2} & \big(t K_{k ;n}^{ \mu}(y,   \zeta, \omega) +  \tilde{K}_{k ;n}^{ \mu}(y,  \zeta )  \big) \frac{  m_{U}(v, \omega ) }{1+\hat{v}\cdot \omega} \\ 
 & \times f(0,x-y+t\omega, v)      \varphi_{j,n}^{i; r}(v, \zeta)     \varphi_{l; r}(\tilde{v}+\omega )   d\omega d v. \\ 
 \end{split}
\ee
where the kernels appeared above are defined as follows, 
\be\label{sep5eqn48}
\begin{split}
K_{k;n}^{ \mu}(\mathfrak{m})(y,  \zeta, \omega)&:= \int_{\R^3} e^{iy \cdot \xi }    \mathfrak{m}_{k,n}(\xi, \zeta)( i   {\omega \cdot \tilde{\xi} }  + i \mu  ) d \xi,\\ 
 \widetilde{K}_{k;n}^{ \mu}(\mathfrak{m})(y,   \zeta )&:= \int_{\R^3} e^{iy \cdot \xi }  |\xi|^{-1}    \mathfrak{m}_{k,n}(\xi, \zeta) d \xi,\\ 
 K_{k;n}^{  }(\mathfrak{m})(y,  \zeta )&:= \int_{\R^3} e^{iy \cdot \xi }     \mathfrak{m}_{k,n}(\xi, \zeta)  d \xi, \\
 K_{k;n}^{ q }(\mathfrak{m})(y,  \zeta )&:= \int_{\R^3} e^{iy \cdot \xi }     \mathfrak{m}_{k,n}(\xi, \zeta) \frac{i\xi_q}{|\xi|} d \xi, \\ 
\end{split}
 \ee
 where the symbol is defined as follows, 
 \be\label{2024oct14eqn31}
  \mathfrak{m}_{k,n}(\xi, \zeta) := \mathfrak{m}(\xi, \zeta)\varphi_k(\xi)\varphi_{n;-M_t}\big(  \tilde{\xi} + \mu \tilde{\zeta}\big). 
 \ee
This concludes the proof of Lemma \ref{locdeclemm}.

\subsection{$L^\infty_x$-type estimates  of the elliptic parts}\label{Linityelliptic}
In this section, we mainly estimate the $L^\infty_x$-norm of the elliptic parts. This includes  the elliptic part $\widetilde{T}_{k,j;n}^{ell;\mu ,i}( \mathfrak{m},U)(t,x,  \zeta )$ of $T_{k,j;n}^{\mu,i}(m, U)(t,x,  \zeta) $ in \eqref{2024oct14eqn41}  for the case $i\in\{0,1\}$, and the  elliptic parts $\mathfrak{E}^{\mu, i'}_{k,j;n}(\mathfrak{m})(s, \cdot, \zeta)$ in \eqref{oct7eqn1} for the case $i'\in \{0,1,2,3\}.$ More precisely, we have
\begin{lemma}\label{ellplinf}
Under the assumption of Theorem \ref{maintheoremellipitic},    for any $i\in\{0,1\}, i'\in\{0,1,2,3\},$ we have
 \be\label{sep5eqn100}
\begin{split}
&\big\|\widetilde{T}_{k,j;n}^{ell;\mu ,i}( \mathfrak{m}, E)(t,x,  \zeta )+ \hat{\zeta}\times  \widetilde{T}_{k,j;n}^{ell;\mu ,i}( \mathfrak{m}, B )(t,x,\zeta)\big\|_{L^\infty_x}+\|\mathfrak{E}^{\mu, i'}_{k,j;n}(\mathfrak{m})(t, x, \zeta)\|_{L^\infty_x}\\
&\lesssim \| \mathfrak{m}(\cdot, \zeta)\|_{\mathcal{S}^\infty} \big[2^{(1-20\epsilon)M_{t^{\star}} } +2^{50\epsilon M_{t^{\star}}}   \mathbf{1}_{n\geq  (1-2 {\alpha}^{\star})M_{t^{\star}}-40\epsilon M_{t^{\star}} } \\
&\qquad \times \min\{2^{(k+2n)/2+  {\alpha}^{\star}  M_{t^{\star}} },2^{(k+4n)/2+(3{\alpha}^{\star}  -1)   M_{t^{\star}} }\} \big]. 
\end{split}
\ee
Moreover, we have the following rough estimate for the elliptic part, 
\be\label{sep8eqn2}
\begin{split}
&\big\|\widetilde{T}_{k,j;n}^{ell;\mu ,i }( \mathfrak{m}, E)(t,x, \zeta)+ \hat{\zeta}\times  \widetilde{T}_{k,j;n}^{ell;\mu,i }( \mathfrak{m}, B )(t,x, \zeta)\big\|_{L^\infty_x}+\|\mathfrak{E}^{\mu, i'}_{k,j;n}(\mathfrak{m})(t, x, \zeta)\|_{L^\infty_x}\\
& \lesssim  2^{ 4{\alpha}^{\star}  M_{t^{\star}}/3 + 5\epsilon M_{t^{\star}} } \big(2^{{\alpha}^{\star}  M_{t^{\star}}/3 } + (2^{M_{t^{\star}} }\frac{|  \zeta_{\bot}|}{|\zeta|})^{1/3}\big). 
\end{split}
\ee
\end{lemma}

\begin{proof}
Recall \eqref{2022feb24eqn81}. In terms of kernel, we have
\be\label{2024nov14eqn11}
\begin{split}
\mathfrak{E}^{\mu, i'}_{k,j;n}(\mathfrak{m})(t, x, \zeta)&= \int_{\R^3} \int_{\R^3} \mathfrak{K}^{i'}_{k,j,n}(\mathfrak{m})(y, v,  \zeta) f(t, x-y, v) d y d v,\\
\mathfrak{K}^{i'}_{k,j,n}(\mathfrak{m})(y, v,  \zeta) &=  \int_{\R^3} e^{iy \cdot \xi}  \mathfrak{m}(\xi, \zeta)  m_{i'}(\xi, v, \zeta)   \psi_k(\xi) \varphi_{j,n}^a(v,\zeta)  \varphi_{n;-M_t}( \tilde{\xi} + \mu \tilde{\zeta}) d \xi .
\end{split}
\ee

As a result of direct computation, from \eqref{sep5eqn96}, we have
\be\label{sep5eqn101}
m_{E}^0(\xi,v) + \hat{\zeta}\times m_{B}^0(\xi,v)= -i 4\pi\big( {(\hat{v}-\hat{\zeta} )\times (\hat{v}\times \xi)+\xi(1-|\hat{v}|^2)} \big). 
\ee

Recall \eqref{2022feb24eqn81}, the definition of cutoff functions $  \varphi_{j,n}^{i}(v, \zeta)$ in \eqref{sep4eqn6},   \eqref{sep5eqn86}  and the definition of the kernel in  \eqref{sep5eqn88}.  From the above equality and the estimate  \eqref{2022feb8eqn51},  after  doing integration by part in $\zeta$ direction and  directions perpendicular to $\zeta$, the following estimate holds for the kernel,
\be\label{2021dec21eqn80}
\begin{split}
&\big| \mathcal{K}^{ell;E,i}_{k,j,n}(\mathfrak{m})(y, v, \zeta)+\hat{\zeta}\times  \mathcal{K}^{ell;B,i}_{k,j,n}(\mathfrak{m})(y, v,   \zeta)\big| +\big|\mathfrak{K}^{i'}_{k,j,n}(\mathfrak{m})(y, v,  \zeta)\big|\\
& \lesssim  2^{2k+2n+\epsilon M_t}\| m(\cdot, \zeta)\|_{\mathcal{S}^\infty}(1+2^k|y\cdot \tilde{\zeta}|)^{-N_0}(1+2^{k+n}|y\times \tilde{\zeta}|)^{-N_0}. 
\end{split}
\ee
From the above estimate of kernel, the conservation law in  \eqref{conservationlaw}, and the estimate  \eqref{nov24eqn41}  if $|  v_{\bot}|\geq 2^{(\alpha_t+\epsilon)M_t}$,  we have
\[
\begin{split}
&\big\|\widetilde{T}_{k,j;n}^{ell;\mu ,i}(m, E)(t,x, \zeta)+ \hat{\zeta}\times  \widetilde{T}_{k,j;n}^{ell;\mu,i }(m, B )(t,x, \zeta)\big\|_{L^\infty_x}+\|\mathfrak{E}^{\mu, i'}_{k,j;n}(\mathfrak{m})(t, x, \zeta)\|_{L^\infty_x}\\
&\lesssim \| \mathfrak{m}(\cdot, \zeta)\|_{\mathcal{S}^\infty} 2^{\epsilon M_t }\min\{2^{-k+j+2(\tilde{\alpha}_t+\epsilon )M_t}\big(\mathbf{1}_{2^n\leq 2^{\epsilon M_t}|  \zeta_{\bot}|/|\zeta|}+ 2^{(\tilde{\alpha}_t M_t - j-n)_{-}} \mathbf{1}_{2^n\geq 2^{\epsilon M_t}|  \zeta_{\bot}|/|\zeta|} \big), 2^{2(k+n)-j}\} \\
&\lesssim \| \mathfrak{m}(\cdot, \zeta)\|_{\mathcal{S}^\infty} 2^{\epsilon M_t }  \min\big\{\big(2^{-k+j+2(\tilde{\alpha}_t+\epsilon )M_t}\big)^{1/2}\big( 2^{2(k+n)-j}\big)^{1/2},\\
&\qquad \big(2^{-k+j+2(\tilde{\alpha}_t+\epsilon )M_t} \big(\mathbf{1}_{2^n\leq 2^{\epsilon M_t}| \zeta_{\bot}|/|\zeta|}+ 2^{(\tilde{\alpha}_t M_t - j-n)_{-}} \mathbf{1}_{2^n\geq 2^{\epsilon M_t}|  \zeta_{\bot}|/|\zeta|} \big) \big)^{2/3}\big( 2^{2(k+n)-j}\big)^{1/3}  \big\} \\
&\lesssim \| \mathfrak{m}(\cdot, \zeta)\|_{\mathcal{S}^\infty} 2^{3\epsilon M_t} \min\{2^{(k+2n)/2+ \tilde{\alpha}_t  M_t},   2^{ 4{\alpha}^{\star}  M_{t^{\star}}/3 + 5\epsilon M_{t^{\star}} }\big(2^{{\alpha}^{\star}  M_{t^{\star}}/3 } + (2^{M_{t^{\star}} }\frac{|\zeta_{\bot}|}{|\zeta|})^{1/3}\big) \}.\\
\end{split}
\]
Hence finishing the  proof of the desired estimates \eqref{sep5eqn100} and \eqref{sep8eqn2}.  
\end{proof}

\subsection{$L^\infty_x$-type estimates of  the bilinear parts   in the large angle region}\label{Linfbigre}

Recall \eqref{2024oct14eqn41}. In this  section, we mainly estimate the $L^\infty_x$-norm of  the bilinear part of the localized acceleration force  in the large angle region. 

  Recall  \eqref{sep5eqn87}. After  using  the decomposition of the acceleration force in  \eqref{july1eqn11}  and    the partition of the characteristic function $\mathbf{1}_{[0,t]}$ in  \eqref{july11eqn33}  to localize the size of $t-s$,  doing dyadic decomposition for the size of $\omega +\tilde{v}$ and the size of $\tilde{v}-\tilde{\zeta}$, and  then  localizing $\omega$ around $\tilde{\zeta}$, we have
 \be\label{sep5eqn99}
 \begin{split}
 &\widetilde{T}_{k,j;n}^{bil;\mu ,i }( \mathfrak{m}, E)(t,x, \zeta )+ \hat{\zeta}\times  \widetilde{T}_{k,j;n}^{bil;\mu,i }(\mathfrak{m}, B )(t,x, \zeta) \\
 & =  \sum_{\begin{subarray}{c}
m\in[-10 M_t,\epsilon M_t]\cap \Z\\
 l\in [-j,2]\cap \Z, a \in \{1,2\}\\ 
\star\in\{main,err\}, r\in [\vartheta^\star_0, 2]\cap \Z \\ 
\end{subarray}}    \widetilde{K}^{\star;m,i;a}_{k,j ; n,l,r }(t, x, \zeta )   ,
\end{split}
 \ee
 where
\be\label{2022feb8eqn53}
\begin{split}
\widetilde{K}^{\star;m,i;a}_{k,j ;n,l,r }(t,   x,\zeta) =  \int_0^t \int_{\R^3} \int_{\R^3}& \int_{\mathbb{S}^2}   EB^a(t,s,x-y  ,\omega, v) \cdot \nabla_v \big[(t-s)\\
& \times \big( \mathcal{K}^{E ,i}_{k,j,n}(\mathfrak{m})(y, v, \omega,\zeta)+ \hat{\zeta}\times \mathcal{K}^{B,i }_{k,j,n}(\mathfrak{m})(y, v, \omega, \zeta)\big)\\
&  +    \mathcal{K}^{err;E,i}_{k,j,n}(\mathfrak{m})(y, v,   \zeta)+ \hat{\zeta}\times \mathcal{K}^{err;B,i}_{k,j,n}(\mathfrak{m})(y, v,  \zeta)  \big] \\
&\times   f(s, x-y+(t-s)\omega , v )  \varphi_{r; \vartheta^\star_0 }(\tilde{v}-\tilde{\zeta}) \\
&\times
 \varphi_{m;-10M_t }(t-s) \varphi_{l;-j}(\tilde{v}+\omega)  \varphi_{ m; k,n}^{\star }(\omega, \zeta)  d\omega d y d v d s, 
 \end{split}
\ee 
 where the kernels $\mathcal{K}^{U,i}_{k,j,n}(\mathfrak{m})(y, v, \omega, \zeta)$ and $\mathcal{K}^{err;U,i}_{k,j,n}(\mathfrak{m})(y, v, \omega, \zeta), U\in \{E, B\},$  are defined in  \eqref{sep5eqn85}     and the cutoff functions $  \varphi_{ m; k,n}^{\star}(\omega, \zeta), \star\in\{main,err\},$ are defined as follows,  
\be\label{sep8eqn21}
\begin{split}
  \varphi_{ m; k,n}^{main}(\omega, \zeta)& := \psi_{\leq \max\{n, -m/2-k/2\}+\epsilon M_t/2 }(\omega  \times \tilde{\zeta}), \\
     \varphi_{ m; k,n}^{err}(\omega,\zeta)&:= \psi_{> \max\{n, -m/2-k/2\}+\epsilon M_t/2 }(\omega \times   \tilde{\zeta}).  \\
\end{split}
\ee

Recall  \eqref{sep5eqn85}.  From  \eqref{sep5eqn101}  and  the estimate  \eqref{2022feb8eqn51},   after doing integration by parts in $\xi$ in $ \zeta $ direction and directions perpendicular to $ \zeta $ or $e_i,i\in\{1,2,3\}$ directions many times, the following estimate holds for any $U\in \{E, B\}, i\in \{0,1\},$
\be\label{july23eqn54}
\begin{split}
 \varphi_{r; n +\epsilon M_t/2 }&(\tilde{v}-\tilde{\zeta}) \Big[|\nabla_v \big( {\mathcal{K}}^{E,i }_{k,j ;n}(y, v,  \zeta )+ \hat{\zeta}\times {\mathcal{K}}^{B,i }_{k,j ;n}(y, v,  \zeta ) \big) | \\
 & + 2^k| \nabla_v \big(  {\mathcal{K}}^{err;E,i }_{k,j ;n}(y, v,  \zeta  ) + \hat{\zeta}\times  {\mathcal{K}}^{err;B,i }_{k,j ;n}(y, v,  \zeta  ) \big) | \\
&+ 2^{-\max\{r, (\gamma_1-\gamma_2)M_{t^\star}\}}\big[ |\nabla_v \mathbf{P}\big(\big( {\mathcal{K}}^{E,i }_{k,j ;n}(y, v,  \zeta ) + \hat{\zeta}\times {\mathcal{K}}^{B,i }_{k,j ;n}(y, v,  \zeta ) \big) \big)  |\\
& + 2^k| \nabla_v \mathbf{P}\big( \big(  {\mathcal{K}}^{err;E,i }_{k,j ;n}(y, v,  \zeta  ) + \hat{\zeta}\times  {\mathcal{K}}^{err;B,i }_{k,j ;n}(y, v,  \zeta  ) \big)  \big) | \big] \Big] \\
& \lesssim 2^{3k+2n} 2^{-j- r }\| \mathfrak{m}(\cdot,   \zeta )\|_{\mathcal{S}^\infty}  \min\big\{(1+2^k|y\cdot \tilde{ \zeta } |)^{-N_0^3} (1+2^{k+n}(|y\times\tilde{\zeta} |))^{-N_0^3},\\
& \qquad  (1+2^{k+n} |  y_{\bot} | )^{-N_0^3}  (1+2^{k+n}(2^{n}+|\tilde{\zeta}\times e_3|)^{-1} |y_3| )^{-N_0^3}\big\}.\\
\end{split}
\ee

 Similar to formula  \eqref{march4eqn41}, for any fixed $\xi, v$, by doing integration by parts in $\omega$ once, we gain at least $2^{-\epsilon M_t/40}$ once for the error type term. After doing integration by parts in $\omega$ many times, from the rough estimate of the electromagnetic field in  \eqref{july10eqn89}  in Proposition \ref{Linfielec} and the estimate of kernel in  \eqref{july23eqn54}, $\forall i\in\{0,1\}, a \in\{1,2\},$ we have
\be\label{sep8eqn51} 
  \big| \widetilde{K}^{err;m,i;a }_{k,j ;l,n}(t, x, \zeta )  \big|\lesssim \| \mathfrak{m}(\cdot, \zeta )\|_{\mathcal{S}^\infty} . 
\ee

Moreover, from the estimate of kernels in  \eqref{july23eqn54}, the estimate  \eqref{july10eqn89}  in Proposition \ref{Linfielec}, and the volume of support of $v, \omega$,  $\forall i\in\{0,1\}, a \in\{1,2\}, $ the following estimate holds if $m=-10M_t,$ 
\be\label{2022feb13eqn1}
  \big|   \widetilde{K}^{main;m, i;a}_{k,j ;l,  n}(t, x, \zeta )\big|\lesssim 2^{m+3j-\vartheta^{\star}_0 +3M_t+20\epsilon M_t} \| \mathfrak{m}(\cdot, \zeta )\|_{\mathcal{S}^\infty} \lesssim  \| \mathfrak{m}(\cdot, \zeta )\|_{\mathcal{S}^\infty}. 
\ee

From now on, we focus on the estimate of the main parts  $  \widetilde{K}^{main;m,i;a }_{k,j ;n}(t, x, \zeta )   $  for the  case $m\in (-10M_t, \epsilon M_t]\cap \Z$, in which we have $t-s\sim 2^{m}. $ 

\begin{lemma}\label{largeregime1}
Let $i\in \{0,1\},$ $  l\in     [-j,2]\cap \Z , m\in (-10M_t, \epsilon M_t]\cap \Z $. Under the assumption of Theorem \ref{mainresultsfirstpart},        the following estimate holds,  
\be\label{2022feb8eqn58}
|\widetilde{K}^{main;m,i;1}_{k,j ;n,l,r} (t,   x, \zeta )|\lesssim  \| \mathfrak{m}(\cdot,   \zeta)\|_{\mathcal{S}^\infty}   2^{  (1-20\epsilon) M_{t^{\star}} }. 
\ee
\end{lemma}
\begin{proof}
Recall  \eqref{2022feb8eqn53}. From the estimate of kernels in  \eqref{july23eqn54}, we have
\be\label{2024oct14eqn61}
\begin{split}
|\widetilde{K}^{main;m,i;1}_{k,j ;n,l,r} (t,   x, \zeta )|&\lesssim  \int_0^t \int_{\R^3} \int_{\R^3} \int_{\mathbb{S}^2}   2^{3k+2n} 2^{-j- \vartheta^{\star}_0+2\epsilon M_t}\| \mathfrak{m}(\cdot,   \zeta )\|_{\mathcal{S}^\infty}(2^{m}+2^{-k})\\
&\quad \times  \big[ |E(s, x-y+(t-s)\omega)-\omega\times B(s, x-y+(t-s)\omega)| \\
&\quad + |\omega\cdot    B(s, x-y+(t-s)\omega)|\big]f(s,x-y+(t-s)\omega,v )\\
&\quad \times  (1+2^k|y\cdot \tilde{ \zeta } |)^{-N_0^3} (1+2^{k+n}(|y\times\tilde{\zeta} |))^{-N_0^3} \\
&\quad \times  \varphi_{j,n}^0 (v, \zeta)\varphi_{m;-10M_t }(t-s) \varphi_{l;-j}(\tilde{v}+\omega)  d\omega d y d v ds. \\
\end{split}
\ee
  From the Cauchy-Schwarz inequality,  the estimate  \eqref{march18eqn31}  in Lemma \ref{conservationlawlemma}, the estimate  \eqref{nov24eqn41}  if $|  v_{\bot}|\geq 2^{(\alpha_t + \epsilon)M_t}$, and the volume of support of $\omega$, we have
 \be
 \begin{split}
\eqref{2024oct14eqn61}& \lesssim \| \mathfrak{m}(\cdot,   \zeta )\|_{\mathcal{S}^\infty}   2^{-j- \vartheta^{\star}_0+2\epsilon M_t}(2^{m }+2^{-k})\\
&\quad  \times  \big(\min\{ 2^{m+3k+2n}, 2^{-2m} \} \min\{2^{3j+2l}, 2^{j+2(\alpha_t + \epsilon)M_t}\}\big)^{1/2}  \\ 
&\quad \times \big(  \min\{2^{m+3j+2l}, 2^{m+j+2(\alpha_t + \epsilon)M_t}, 2^{-2m-j-2l}, 2^{m+3k+2n-j}\}\big)^{1/2}\\
&\lesssim \| \mathfrak{m}(\cdot,   \zeta )\|_{\mathcal{S}^\infty}   2^{-j- \vartheta^{\star}_0+2\epsilon M_t}\big[2^{m} \min\big\{\big(2^{-2m+3j+2l}\big)^{1/2}\big(2^{-2m-j-2l}\big)^{1/2}, \\
&\quad  \big(2^{-2m+j+2(\alpha_t + \epsilon)M_t}\big)^{1/2}\big( 2^{m+j+2(\alpha_t + \epsilon)M_t}\big)^{1/2} \big\} \\
&\quad +2^{-k}\min\big\{ \big((2^{m+3k+2n})^{2/3}(2^{-2m})^{1/3}2^{3j+2l} \big)^{1/2} \big(2^{-2m-j-2l}\big)^{1/2},\\
&\quad \big((2^{m+3k+2n})^{2/3}(2^{-2m})^{1/3}2^{j+2(\alpha_t + \epsilon)M_t} \big)^{1/2} \big( 2^{m+j+2(\alpha_t + \epsilon)M_t}\big)^{1/2}  \big\}  \big]\\
&\lesssim \| \mathfrak{m}(\cdot,   \zeta )\|_{\mathcal{S}^\infty} 2^{- \vartheta^{\star}_0+ 10\epsilon  M_{t^{\star}}} \min\{2^{-m}, 2^{m/2+2\alpha^{\star}M_{t^{\star}}} \}\\
&\lesssim \| \mathfrak{m}(\cdot,   \zeta )\|_{\mathcal{S}^\infty} 2^{4\alpha^{\star}M_{t^{\star}}/3- \vartheta^{\star}_0+10\epsilon  M_{t^{\star}}}\lesssim  \| \mathfrak{m}(\cdot,   \zeta)\|_{\mathcal{S}^\infty}   2^{  (1-20\epsilon) M_{t^{\star}} }. \\
\end{split}
 \ee
Hence finishing the proof of our desired estimate  \eqref{2022feb8eqn58}. 
\end{proof}

\begin{lemma}\label{largeregime2}
 Let  $i\in \{0,1\}$, $  l\in     [-j,2]\cap \Z , m\in (-10M_t, \epsilon M_t]\cap \Z $. Under the assumption of Theorem \ref{mainresultsfirstpart},        the following estimate holds   if $n\geq (1-2\alpha^{\star}-35\epsilon)M_{t^{\star}}  $, $    j\leq   (1/2+3\iota + 55\epsilon) M_{t^{\star}}  $,  
\be\label{aug5eqn33}
\begin{split}
|\widetilde{K}^{main;m,i;2}_{k,j ; n,l,r  } (t,   x, \zeta )| &\lesssim  \| \mathfrak{m}(\cdot,   \zeta)\|_{\mathcal{S}^\infty} \big[  2^{  (1-19.5\epsilon) M_{t^{\star}} } +2^{ 120.5\epsilon M_{t^{\star}} }  \mathbf{1}_{n\geq   (1-2\alpha^{\star}-35\epsilon)M_{t^{\star}}} \\
  & \qquad \times \min\{2^{(k+2n)/2+ (\alpha^{\star}+3\iota) M_{t^{\star}} },2^{(k+4n)/2+(1+6\iota )M_{t^{\star}} } \}  \big]. \\
\end{split}
\ee
Moreover, if  $k+2j\leq    (2-50\epsilon)  M_{t^{\star}}$, then we have 
\be\label{2022feb11eqn151}
 |\widetilde{K}^{main;m,i;2}_{k,j ; n,l,r  } (t,   x, \zeta )| \lesssim  \| \mathfrak{m}(\cdot,   \zeta)\|_{\mathcal{S}^\infty}    2^{  (1-19.5\epsilon) M_{t^{\star}} }.
\ee
\end{lemma}
\begin{proof}
Recall \eqref{2022feb8eqn53}, \eqref{july1eqn13},   and the estimate of kernels in  \eqref{july23eqn54}. Let $\omega =(\sin \theta \cos\phi, \sin \theta \sin \phi, \cos\theta).$  After 
 localizing the sizes of  $\theta, \phi $, the following estimate holds, 
\be\label{aug5eqn1}
\begin{split}
|\widetilde{K}^{main;m,i;2}_{k,j ; n,l,r } (t,   x, \zeta )| &\lesssim \sum_{p, q\in [-10M_t,2]\cap \Z } \big|{H}^{m, i ;p,q  }_{k,j ; n,l,r  } (t,   x,  \zeta)\big|,  \\
{H}^{m, i;p,q }_{k,j ; n,l,r  }(t,   x,  \zeta)&:=     \int_{0}^{t }\int_{\R^3}\int_{\R^3} \int_0^{2\pi} \int_{0}^{\pi}  \| \mathfrak{m}(\cdot,   \zeta )\|_{\mathcal{S}^\infty}  (2^m+2^{-k})  2^{3k+2n} 2^{-j- r  +2\epsilon M_t}\\
&\quad \times |\hat{v}+\omega|   |B(s, x  -y+(t-s)\omega )|f(s, x-y+(t-s)\omega, v)\\
&\quad  \times \min\big\{(1+2^k|y\cdot \tilde{ \zeta } |)^{-N_0^3} (1+2^{k+n}(|y\times\tilde{\zeta} |))^{-N_0^3}, \\
& \qquad  (1+2^{k+n} | y_{\bot} | )^{-N_0^3} (1+2^{k+n}(2^{n}+|\tilde{\zeta}\times e_3|)^{-1} |y_3| )^{-N_0^3}\big\}\\
&\quad \times  \varphi_{r; \vartheta^\star_0 }(\tilde{v}-\tilde{\zeta})  \varphi_{l;-j}(\tilde{v}+\omega)   \varphi_{m;-10M_t }(t-s) \varphi_{p;-10M_t}(\sin \theta )\\
&\quad  \times  \varphi_{q;-10M_t}(\sin\phi )   \varphi_{ m; k,n}^{main }(\omega, \zeta)  \sin \theta d\theta d \phi  d v d y  d s. 
\end{split}
\ee

From the estimate  \eqref{july10eqn89}  in Proposition \ref{Linfielec}, and the volume of support of $v, \omega$, we can rule out the case $p=-10M_t$  and the case $q=-10M_t$. Hence, it suffices to consider the case $p, q\in (-10M_t, 2]\cap \Z,$ in which we have $\sin \theta \sim 2^p, \sin \phi \sim 2^q. $

Based on  different sizes of parameters, e.g., $m+k,l,n,m+p, k+n, | x_{\bot}|$,  we separate into four cases as follows.

\medskip 

 \noindent \textbf{Case 1.} \qquad If $   m+k\leq -2n+\epsilon M_t $. 
 
\medskip 

Recall  \eqref{sep8eqn21}. For this case,  the volume of support of $\omega$ for  the main part is bounded by $2^{-m-k+\epsilon M_t}$. Therefore, by using the estimate \eqref{march18eqn31}  in Lemma \ref{conservationlawlemma},  the volume of support of $\omega$ and $ v$, and the estimate if $| v_{\bot} |\geq 2^{(\alpha_t +\epsilon)M_t},$ we have
\be\label{aug5eqn13}
\begin{split}
&|{H}^{m, i;p,q }_{k,j ; n,l,r  }(t, x, \zeta)|\\
&\lesssim  \sup_{s\in [0,t]}\|B (s,\cdot)\|_{L^2} \| \mathfrak{m}(\cdot, \zeta)\|_{\mathcal{S}^\infty}(2^m+2^{-k})   2^{-j- r  +l+2\epsilon M_t} \\ 
&\quad\times  \big(2^{m+3k+2n}\min\{1, 2^{-m-k}\}  \min\{ 2^{3j+2\min\{l,r\}}, 2^{j+2 (\alpha_t +\epsilon)M_t}\}  \big)^{1/2} \\
&  \quad \times   \big(\min\{ 2^{m}\min\{1, 2^{-m-k}\} \min\{ 2^{3j+2\min\{l,r\}}, 2^{j+2 (\alpha_t +\epsilon)M_t}\} , \\
& \,\,\quad 2^{3k+2n-j  } 2^{m}\min\{1, 2^{-m-k}\}, 2^{-2m -j-2l}   \} \big)^{1/2}  \\
 &\lesssim   \| \mathfrak{m}(\cdot, \zeta)\|_{\mathcal{S}^\infty}  2^{-j-r+l+2\epsilon M_t}\big[ 2^{m}\min\big\{\big(2^{2k+2n+3j+2\min\{l,r\}}\big)^{1/2}\\
 &\quad \times \big(  2^{-2m-j-2l} \big)^{1/2},  \big(2^{k-m+2\epsilon M_t}\min\{ 2^{3j+2\min\{l,r\}}, 2^{j+2 (\alpha_t +\epsilon)M_t}\} \big)^{1/2}\\
 &\quad \times \big(2^{-k}\min\{ 2^{3j+2\min\{l,r\}}, 2^{j+2 (\alpha_t +\epsilon)M_t}\} \big)^{1/2}\big\}+2^{-k}\min\big\{\big(  2^{2k+2n-j} \big)^{1/2} \\
 &\quad \times \big(2^{2k+2n+3j+2\min\{l,r\}}\big)^{1/2} ,\big(2^{m}\min\{ 2^{3j+2\min\{l,r\}}, 2^{j+2 (\alpha_t +\epsilon)M_t}\} \big)^{1/2} \\
 &\quad \times   \big(2^{2k+2n+2\epsilon M_t}\min\{ 2^{3j+2\min\{l,r\}}, 2^{j+2 (\alpha_t +\epsilon)M_t}\} \big)^{1/2}  \big\}\big]\\
  &\lesssim  \| \mathfrak{m}(\cdot, \zeta)\|_{\mathcal{S}^\infty} 2^{  5\epsilon M_t} 2^{-r+l}\min\big\{ 2^{k+n+\min\{l,r\}-l}, 2^{m/2 }  \min\{ 2^{ 2j+2\min\{l,r\}}, 2^{ 2 (\alpha_t +\epsilon)M_t}\} \big\} .
\end{split}
\ee

\medskip 

 \noindent \textbf{Case 2.} \qquad  If $m+k\geq -2n+\epsilon M_t$ and  $m+p\leq -k-n+2\epsilon M_t$. 
 
\medskip

  For this case,  the volume of support of $\omega$ for  the main part is bounded by $2^{2n+\epsilon M_t}$. From the Cauchy-Schwarz inequality,  the estimate \eqref{march18eqn31}  in Lemma \ref{conservationlawlemma}, and   the volume of support of $(\omega, v)$,    we have
\be\label{sep8eqn48}
\begin{split}
&|{H}^{m, i;p, q  }_{k,j ; n,l,r  }(t, x, \zeta )|\\
&\lesssim  \| \mathfrak{m}(\cdot, \zeta)\|_{\mathcal{S}^\infty}   2^{m -j- r +l+2\epsilon M_t}  \big(\min\{2^{m+2\min\{p,n\} }  2^{ 3j+2l} ,  2^{-2m  -j-2l}   \} \big)^{1/2}\\
&\qquad \times \big(2^{m+3k+2n+2\min\{p,n\}}   2^{ 3j+2\min\{l,r\}}  \big)^{1/2}  \big( \sup_{s\in [0,t]}\|B (s,\cdot)\|_{L^2} \big)\\
&\lesssim   \| \mathfrak{m}(\cdot, \zeta)\|_{\mathcal{S}^\infty}  2^{ 4\epsilon M_t}\min\big\{  2^{2m+3k/2+n + 2\min\{p,n\}}  2^{2j+2l} ,   2^{m/2+3k/2+n+\min\{p,n\}} \big\}\\
&  \lesssim  \| \mathfrak{m}(\cdot, \zeta)\|_{\mathcal{S}^\infty} 2^{     4\epsilon M_t} \min\big\{    2^{2j+2l}   2^{-k/2-n},   2^{m/2+2j+2l} ,  2^{k+n/2+\min\{p,n\}/2} \big\} \\
&\lesssim \| \mathfrak{m}(\cdot, \zeta)\|_{\mathcal{S}^\infty}   2^{     4\epsilon M_t} \min\big\{ \big(2^{m/2+2j+2l} \big)^{1/2} \big( 2^{k+n}\big)^{1/2}, \big(  2^{2j+2l}   2^{-k/2-n}\big)^{2/3}\big( 2^{k+n}\big)^{1/3}\big\} \\
&\lesssim  \| \mathfrak{m}(\cdot, \zeta)\|_{\mathcal{S}^\infty} 2^{ 5\epsilon M_t}\min\{2^{(k+n)/2+j},  2^{4(j+l)/3  -n/3}\}  .
\end{split}
\ee

 \medskip 

 \noindent \textbf{Case 3.} \qquad  If $m+k\geq -2n+\epsilon M_t$, $m+p\geq -k-n+2\epsilon M_t$  and $|  x_{\bot}|\geq 2^{m+p-10}$.  
 
\medskip

Recall  \eqref{aug5eqn1}. Note that   we have $|\tilde{\zeta} \times \mathbf{e}_3|\lesssim 2^{\max\{n,p\}+\epsilon M_t}$ for the case we are considering.  After   using the Cauchy-Schwarz inequality,  the volume of support of $v,\omega$, the estimate  \eqref{nov24eqn41}  if $| v_{\bot}|\geq 2^{(\alpha_t+\epsilon)M_t}$,    the estimate of Jacboian $(y_1,y_2, \theta)\longrightarrow x-y+(t-s)\omega$,    or the estimate of Jacobian in   \eqref{march18eqn66} for changing coordinates $(\theta, \phi) \longrightarrow (z(y),w(y))$, 
  where
  \[
  \begin{split}
   z(y)&:=\sqrt{\big(|  x_{\bot}- y_{\bot}+(t-s)\sin\theta\cos\phi\big)^2 + (t-s)^2(\sin \theta \sin \phi)^2 }, \\
     w(y)&:= x_3 -y_3+(t-s)\cos\theta, \\
   \end{split}
  \]
   we have
\be\label{aug5eqn19}
\begin{split}
\big|{H}^{m , i ;p,q  }_{k,j ; n,l,r }(t, x, \zeta )\big| & \lesssim \sup_{s\in [0,t]}\|B (s,\cdot)\|_{L^2} 2^{m-j-r    +l+2\epsilon M_t} \| \mathfrak{m}(\cdot, \zeta)\|_{\mathcal{S}^\infty} \\
&\qquad\times  \big( \min\{ \frac{2^{-p-q} }{2^{m}| x_{\bot}|}, 2^{2k+n +\max\{n,p\}+q  } \} \min\{2^{3j+2 \min\{l,r\} }, 2^{j+2\tilde{\alpha}_t M_t}\} \}\big)^{1/2}  \\
&\qquad \times  \big( \min\{ 2^{m+2\min\{n,p\}   } 2^{3j+2 \min\{l,r\} }, 2^{-2m-j-2l} \} \big)^{1/2} \\
&\lesssim  \| \mathfrak{m}(\cdot, \zeta)\|_{\mathcal{S}^\infty}  2^{m-j- r/2+l/2 +2\epsilon M_t} \big(  2^{-m/2+j +\min\{n,p\}} \big)^{1/2}\\
&\qquad \times   (2^{-m+k+\max\{n,p\}-p  } \min\{2^{3j+2 \min\{l,r\} }, 2^{j+2\tilde{\alpha}_t M_t}\})^{1/2}\\
&\lesssim  2^{(k+n)/2    +2\epsilon M_t}\min\{2^{j+l}, 2^{\tilde{\alpha}_t M_t- r/2+l/2}\} .
\end{split}
 \ee

  \medskip 

 \noindent \textbf{Case 4.} \qquad   If $m+k\geq -2n+\epsilon M_t$, $m+p\geq -k-n+2\epsilon M_t$  and $|  x_{\bot}|\leq 2^{m+p-10}$. 
 
\medskip

 Note that, for this case, we have $| x_{\bot}-  y_{\bot} + (t-s)  \omega_{\bot}|\sim 2^{m+p}$ for any $y\in B(0, 2^{-k-n+\epsilon M_t/2}).$    After using the Cauchy-Schwarz inequality,  the cylindrical symmetry of solution, and changing coordinates $\theta\longrightarrow x_3 -y_3+(t-s)\cos\theta$, the estimate  \eqref{nov24eqn41}  if $| v_{\bot}|\geq 2^{\alpha_t M_t + \epsilon M_t}, $ from the estimate  \eqref{aug5eqn1} and the estimates of kernels in    \eqref{july23eqn54},  we have 
\be\label{aug5eqn44}
\begin{split}
|{H}^{m,i;p,q }_{k,j ; n,l,r }(t, x, \zeta )| &
\lesssim  \big(\sup_{s\in [0,t]}\|B (s,\cdot)\|_{L^2}\big)  2^{m -j- r  +l+2\epsilon M_t}  \| \mathfrak{m}(\cdot,   \zeta)\|_{\mathcal{S}^\infty} \\
&\quad \times \big(  \min\{   2^{m+2\min\{p,n\} }  \min\{2^{3j+2  \min\{l,r\} }, 2^{j+2 (\tilde{\alpha}_t+\epsilon) M_t}\} ,   2^{-2m-j-2l} \} \big)^{1/2}\\
&\quad  \times \big( \frac{2^{k+\max\{n,p\}+\epsilon M_t}}{2^{m+p}}  \min\{2^{3j+2  \min\{l,r\} }, 2^{j+2 (\tilde{\alpha}_t+\epsilon) M_t}\} \big)^{1/2} \\
&  \lesssim   2^{- r/2+l/2+ m/4+k/2 +(\max\{n,p\}+\min\{n,p\}-p)/2+5\epsilon M_t}  \\
&\quad \times \min\{2^{j+\min\{l,r\}},2^{\tilde{\alpha}_t M_t}\} \| \mathfrak{m}(\cdot, \zeta)\|_{\mathcal{S}^\infty} \\
& \lesssim 2^{(k+n)/2  +5\epsilon M_t}\min\{2^{j+l},2^{\tilde{\alpha}_t M_t- r/2+l/2}\} \| \mathfrak{m}(\cdot, \zeta)\|_{\mathcal{S}^\infty}.\\
\end{split}
\ee

Therefore, if $n\geq (1-2\alpha^{\star}-35\epsilon)M_{t^{\star}}  $, $    j\leq   (1/2+3\iota + 55\epsilon) M_{t^{\star}}  $, then the following estimate holds after combining  the obtained estimates  \eqref{aug5eqn13},  \eqref{sep8eqn48},  \eqref{aug5eqn19}, and \eqref{aug5eqn44},  
\be 
\begin{split}
 |{H}^{m,i;p,q }_{k,j ; n,l,r }(t, x, \zeta )| 
& \lesssim \| \mathfrak{m}(\cdot,   \zeta)\|_{\mathcal{S}^\infty} \big[  2^{  (1-20\epsilon) M_{t^{\star}} }+2^{  120 \epsilon M_{t^{\star}} } \mathbf{1}_{n\geq   (1-2\alpha^{\star}-35\epsilon)M_{t^{\star}}} \\
 & \quad \times   \min\{2^{(k+2n)/2+ (\alpha^{\star}+3\iota) M_{t^{\star}} },2^{(k+4n)/2+(1+6\iota )M_{t^{\star}} } \} \big].\\
  \end{split}
 \ee
 Hence finishing the proof of our desired estimate  \eqref{aug5eqn33}. 

 Moreover, if   $k+2j\leq   (2-50\epsilon)  M_{t^{\star}}$, then the desired estimate  \eqref{2022feb11eqn151}  holds from the obtained estimates\eqref{aug5eqn13},  \eqref{sep8eqn48},  \eqref{aug5eqn19}, and \eqref{aug5eqn44}. 

\end{proof}

\begin{lemma}\label{mainlemmasecpart}
 Let      $  m\in (-10M_t, \epsilon M_t]\cap \Z $, $i\in \{0,1\},$   $l\in [-j,2]\cap \Z, j\in [0, (1+2\epsilon)M_t]\cap \Z, k\in \Z_+$, s.t.,  $k+2j> (2-50\epsilon)  M_{t^{\star}}$.  Under the assumption of Theorem \ref{mainresultsfirstpart},        the following estimate holds if $n<  (1-2\alpha^{\star}-35\epsilon)M_{t^{\star}}  $ or  $    j\geq   (1/2+3\iota + 55\epsilon) M_{t^{\star}}   $, 
\be\label{aug3eqn14}
\begin{split}
|\widetilde{K}^{main;m, i ;2}_{k,j ; n,l,r}(t, x,\zeta )| & \lesssim \| \mathfrak{m}(\cdot, \zeta)\|_{\mathcal{S}^\infty}  \big[ 2^{(1-19.5\epsilon)M_{t^{\star}}} +  2^{  100.5\epsilon M_{t^{\star}} } \mathbf{1}_{n\geq  -(\alpha^{\star}+3\iota+60\epsilon) M_{t^{\star}} } \\
  &\qquad \times \min\{2^{(k+2n)/2  +(2{\alpha}^{\star} -1)M_{t^{\star}}}  , 2^{(k+4n)/2 +(1+6\iota)M_{t^{\star}}} \}  \big]. \\
\end{split}
\ee
\end{lemma}
\begin{proof}
Recall  \eqref{july1eqn13},   \eqref{2022feb8eqn53}, and the decomposition of $\widetilde{K}^{main;m,i;2}_{k,j ;l,n} (t,   x, \zeta )$ in \eqref{aug5eqn1}.   We localize ${H}^{m, i ;p,q  }_{k,j ;l,n} (t,   x,  \zeta)$ further by  using the decomposition of the magnetic  field  in  \eqref{july5eqn1} (see also \eqref{july5eqn60}) and the   index sets defined in  \eqref{indexsetL2}. More precisely, from the estimate  of kernels in  \eqref{july23eqn54},  we have
\be\label{aug10eqn31}
\big|{H}^{m, i ;p,q  }_{k,j ;  n,l,r } (t,   x,  \zeta)\big|\lesssim \sum_{ (\tilde{m}, \tilde{k}, \tilde{j}, \tilde{l}) \in \mathcal{S}_1(t)\cup  \mathcal{S}_2(t) }   {H}^{m,i;p,q;\tilde{m}, \tilde{k},\tilde{j},\tilde{l} }_{k,j ; n,l,r }(t,x, \zeta ),
\ee
where 
\be\label{aug3eqn32}
\begin{split}
\big| {H}^{m,i;p,q;\tilde{m}, \tilde{k},\tilde{j},\tilde{l} }_{k,j ; n,l,r }(t,x, \zeta )\big|&\lesssim (2^{m}+2^{-k})2^{3k+2n-j - r +l +2\epsilon M_t} \| \mathfrak{m}(\cdot, \zeta)\|_{\mathcal{S}^\infty} \\
&\quad \times  \int_{0}^{t }\int_{\R^3}\int_{\R^3} \int_0^{2\pi}\int_{0}^{\pi} \big|B^{\tilde{m}}_{\tilde{k};\tilde{j},  \tilde{l}} (s, x-y+(t-s)\omega)\big|\\
&\quad \times \min\big\{(1+2^k|y\cdot \tilde{ \zeta } |)^{-N_0^3} (1+2^{k+n}(|y\times\tilde{\zeta} |))^{-N_0^3},\\
&\qquad (1+2^{k+n} | y_{\bot} | )^{-N_0^3} (1+2^{k+n}(2^{n}+|\tilde{\zeta}\times e_3|)^{-1} |y_3| )^{-N_0^3}\big\}\\
&\quad\times f(s,x-y+(t-s)\omega, v)     \varphi_{r; \vartheta^\star_0 }(\tilde{v}-\tilde{\zeta})   \varphi_{l;-j}(\tilde{v}+\omega) \\
& \quad \times    \varphi_{m;-10M_t}(t-s)   \varphi_{p;-10M_t}(\sin \theta ) \varphi_{q;-10M_t}(\sin\phi ) \sin \theta d\theta d \phi  d v d y  d s.
\end{split} 
\ee

Based on  different sizes of parameters, e.g., $m+k,l,n,m+p, k+n, |  x_{\bot}|$,  we separate into four cases as follows.

\medskip

\noindent \textbf{Case 1.}\qquad If $m+k\leq -2n+ \epsilon M_t  $. 

\medskip

Similar to the obtained estimate in  \eqref{aug5eqn13},  from the Cauchy-Schwarz inequality,  the volume of support of $v, \omega, $ the estimate  \eqref{march18eqn31}  in Lemma \ref{conservationlawlemma}, the conservation law  \eqref{conservationlaw}, and the estimate  \eqref{nov24eqn27}  if $| v_{\bot} |\geq 2^{(\alpha_t+\epsilon)M_t}$,  the following estimate holds after putting $B^{\tilde{m}}_{\tilde{k};\tilde{j}, -\tilde{l}}(s,\cdot)$ in $L^2_x$,  
\be\label{aug4eqn41}
\begin{split}
|   {H}^{m,i;p,q;\tilde{m}, \tilde{k},\tilde{j},\tilde{l} }_{k,j ; n,l,r  }(t,x, \zeta ) | & \lesssim  \sup_{s\in [0,t]}\|B^{\tilde{m}}_{\tilde{k};\tilde{j},  \tilde{l}}(s,\cdot)\|_{L^2} 2^{ -r+l + 4\epsilon M_t} \| \mathfrak{m}(\cdot, \zeta)\|_{\mathcal{S}^\infty}  \\
& \quad \times \min\big\{ 2^{k+n+\min\{l,r\}-l}, 2^{m/2 } \min\{2^{2j+2\min\{r,l\}}, 2^{2\tilde{\alpha}_t M_t}\} \big\}\\
&   \lesssim\sup_{s\in [0,t]}\|B^{\tilde{m}}_{\tilde{k};\tilde{j},  \tilde{l}}(s,\cdot)\|_{L^2} 2^{ -r+l + 4\epsilon M_t} \big( 2^{k+n+\min\{l,r\}-l}\big)^{1/2} \\
&\quad \times  \big( 2^{-(k+2n)/2+\epsilon M_t/2 }\min\{2^{2j+2\min\{r,l\}}, 2^{2\tilde{\alpha}_t M_t}\} \big)^{1/2}\\
  &     \lesssim 2^{-r/2+l/2  + k/4 +5\epsilon M_t} \sup_{s\in [0,t]}\|B^{\tilde{m}}_{\tilde{k};\tilde{j}, \tilde{l}}(s,\cdot)\|_{L^2} \\
  &\qquad \times \| \mathfrak{m}(\cdot, \zeta)\|_{\mathcal{S}^\infty} \min\{2^{j+\min\{r,l\}}, 2^{\tilde{\alpha}_t M_t}\} . \\
\end{split}
\ee

Alternatively, if we put the localized magnetic field in $L^\infty$, then from the estimate \eqref{aug3eqn32},  the volume of support of $v,\omega$, and the estimate  \eqref{march18eqn31}  in Lemma \ref{conservationlawlemma},  we have
\be
\begin{split}
|   {H}^{m,i;p,q;\tilde{m}, \tilde{k},\tilde{j},\tilde{l} }_{k,j ; n,l,r }(t,x, \zeta ) |&\lesssim  \sup_{s\in [0,t]}\|B^{\tilde{m}}_{\tilde{k};\tilde{j}, \tilde{l}}(s,\cdot)\|_{L^\infty}    (2^{m}+2^{-k}) 2^{-j- r  +l+3\epsilon M_t} \\
&\quad\times    \min\{2^{-k} 2^{3j+2\min\{l,r\} } , 2^{2k+2n -j  }, 2^{-2m  -j-2l}   \}     \| \mathfrak{m}(\cdot, \zeta)\|_{\mathcal{S}^\infty}\\
&\lesssim \sup_{s\in [0,t]}\|B^{\tilde{m}}_{\tilde{k};\tilde{j}, \tilde{l}}(s,\cdot)\|_{L^\infty}   2^{-j- r  +l+3\epsilon M_t}\big[ 2^m \min\big\{( 2^{-k+3j+2\min\{l,r\} } )^{1/2} \\
&\quad \times (2^{-2m  -j-2l})^{1/2} ,   2^{2k+2n -j  } \big\}+2^{-k} \min\big\{( 2^{-k+3j+2\min\{l,r\} } )^{1/2}\\
 &\quad \times (2^{2k+2n -j  })^{1/2} ,   2^{2k+2n -j  } \big\} \big] \| \mathfrak{m}(\cdot, \zeta)\|_{\mathcal{S}^\infty}.\\
\end{split}
\ee
From the above estimate, we concludes that
\be\label{aug4eqn42}
\begin{split}
  |   {H}^{m,i;p,q;\tilde{m}, \tilde{k},\tilde{j},\tilde{l} }_{k,j ; n,l,r }(t,x, \zeta ) | & \lesssim \sup_{s\in [0,t]}\|B^{\tilde{m}}_{\tilde{k};\tilde{j}, \tilde{l}}(s,\cdot)\|_{L^\infty }\| \mathfrak{m}(\cdot, \zeta)\|_{\mathcal{S}^\infty} \\
  &\quad \times 2^{-r+l  + 4\epsilon M_t }\min\{2^{-k/2+r-l}, 2^{k-2j}    \} .\\
  \end{split}
\ee

If $(\tilde{m},\tilde{k},\tilde{j},\tilde{l})\in \mathcal{S}_1(t)$, then from the first estimate in   \eqref{aug4eqn10}  in Proposition \ref{meanLinfest} and the estimate  \eqref{aug4eqn42}, we have 
\be\label{aug4eqn51}
\begin{split}
 \sum_{(\tilde{m},\tilde{k},\tilde{j},\tilde{l})\in \mathcal{S}_1(t)} & | {H}^{m,i;p,q;\tilde{m}, \tilde{k},\tilde{j},\tilde{l} }_{k,j ; n,l,r }(t, x, \zeta )|\\
 &\lesssim \min\{2^{-k/2},  2^{-2j/3}\}  2^{- \vartheta^{\star}_0  + 15 \epsilon M_t} \big( 2^{ 2 \tilde{\alpha}_t  M_t }  + 2^{7M_t/6+\tilde{\alpha}_t M_t/4} \big) \| \mathfrak{m}(\cdot, \zeta)\|_{\mathcal{S}^\infty}\\
 &\lesssim \big(2^{-k/2}\big)^{2/5} \big(2^{-2j/3}\big)^{3/5} 2^{- \vartheta^{\star}_0  + 15 \epsilon M_t} \big( 2^{ 2 \tilde{\alpha}_t  M_t }  + 2^{7M_t/6+\tilde{\alpha}_t M_t/4} \big) \| \mathfrak{m}(\cdot, \zeta)\|_{\mathcal{S}^\infty}\\
&  \lesssim 2^{-(k+2j)/5- \vartheta^{\star}_0  + 15\epsilon M_t} \big( 2^{ 2 \tilde{\alpha}_t  M_t }   + 2^{7M_t/6+\tilde{\alpha}_t M_t/4} \big) \| \mathfrak{m}(\cdot, \zeta)\|_{\mathcal{S}^\infty}\\
& \lesssim  2^{(1-20\epsilon)M_{t^{\star}}}\| \mathfrak{m}(\cdot, \zeta)\|_{\mathcal{S}^\infty}.\\
\end{split}
\ee

If  $(\tilde{m},\tilde{k},\tilde{j},\tilde{l})\in \mathcal{S}_2(t)$, after combining the estimates  \eqref{aug4eqn41}  and \eqref{aug4eqn42}  and using the second estimate in \eqref{aug4eqn10}  in Proposition \ref{meanLinfest}, we have
\be\label{aug4eqn50}
\begin{split}
\sum_{(\tilde{m},\tilde{k},\tilde{j},\tilde{l})\in \mathcal{S}_2(t)} & | {H}^{m,i;p,q;\tilde{m}, \tilde{k},\tilde{j},\tilde{l} }_{k,j ; n,l,r }(t, x, \zeta )|\\
&\lesssim  \sum_{(\tilde{m},\tilde{k},\tilde{j},\tilde{l})\in \mathcal{S}_2(t)}  2^{ -r    + 5\epsilon M_t}  \big( (2^{-k/2 })^{5/6} (2^{k-2j})^{1/6} \sup_{s\in [0,t]}\|B^{\tilde{m}}_{\tilde{k};\tilde{j}, -\tilde{l}}(s,\cdot)\|_{L^\infty }\big)^{1/2} \\ 
 &\quad \times  \big( 2^{k/4  } \min\{2^{j+\min\{l,r\}}, 2^{\tilde{\alpha}_t M_t}\}  \sup_{s\in [0,t]}\|B^{\tilde{m}}_{\tilde{k};\tilde{j}, \tilde{l}}(s,\cdot)\|_{L^2 }\big)^{1/2}  \| \mathfrak{m}(\cdot, \zeta)\|_{\mathcal{S}^\infty} \\
 &\lesssim  2^{- \vartheta^{\star}_0  + \tilde{\alpha}_t M_t/3+40\epsilon M_t} \big( 2^{\tilde{\alpha}_t M_t}  + 2^{7M_t/12+ \tilde{\alpha}_t M_t/8} \big)\| \mathfrak{m}(\cdot, \zeta)\|_{\mathcal{S}^\infty}\\
 &\lesssim 2^{(1-20\epsilon)M_{t^{\star}}}\| \mathfrak{m}(\cdot, \zeta)\|_{\mathcal{S}^\infty}.
 \end{split}
\ee
To sum up,   we have
\be\label{aug4eqn52}
 \sum_{(\tilde{m},\tilde{k},\tilde{j},\tilde{l})\in \mathcal{S}_1(t) \cup \mathcal{S}_2(t) }  | {H}^{m,i;p,q;\tilde{m}, \tilde{k},\tilde{j},\tilde{l} }_{k,j ; n,l,r }(t, x, \zeta)|\lesssim   2^{(1-20\epsilon)M_{t^{\star}}}\| \mathfrak{m}(\cdot, \zeta)\|_{\mathcal{S}^\infty}. 
\ee

\medskip

\noindent \textbf{Case 2.}\qquad If   $m+k\geq -2n+ \epsilon M_t  $  and  $m+p\leq -k-n+2\epsilon M_t$.

\medskip

Similar to the obtained estimate in  \eqref{sep8eqn48}, from the estimate  \eqref{nov24eqn41}  if $|  v_{\bot}|\geq 2^{(\alpha_t+\epsilon )M_t}$, we have 
\be\label{july23eqn61}
\begin{split}
 | {H}^{m,i;p,q;\tilde{m}, \tilde{k},\tilde{j},\tilde{l} }_{k,j ; n,l,r }(t, x, \zeta )| &\lesssim  2^{-n/3- r  + 5\epsilon M_t} \| \mathfrak{m}(\cdot, \zeta)\|_{\mathcal{S}^\infty} \sup_{s\in [0,t]}\|B^{\tilde{m}}_{\tilde{k};\tilde{j}, \tilde{l}}(s,\cdot)\|_{L^2}\\
 &\quad \times \min\{2^{4(j+l)/3}, 2^{4\tilde{\alpha}_t M_t/3}\} . \\
 \end{split}
\ee
Alternatively, we have
\be\label{july23eqn62} 
\begin{split}
  | {H}^{m,i;p,q;\tilde{m}, \tilde{k},\tilde{j},\tilde{l} }_{k,j ;n,l,r}(t, x, \zeta )|  
&\lesssim  2^{m -j- r   +l+2\epsilon M_t}  \| \mathfrak{m}(\cdot, \zeta)\|_{\mathcal{S}^\infty} \sup_{s\in [0,t]}\|B^{\tilde{m}}_{\tilde{k};\tilde{j}, \tilde{l}}(s,\cdot)\|_{L^\infty}  \\
&\quad \times \min\{2^{m+3j+2\min\{l,r\}+2\min\{p,n\}} , 2^{-2m -j-2l}   \} \\
 &\lesssim \min\big\{\big( 2^{m+3j+2\min\{l,r\}+2\min\{p,n\}}\big)^{1/2}\big( 2^{-2m -j-2l}\big)^{1/2} , \\
 &\quad  \big( 2^{m+3j+2\min\{l,r\}+2\min\{p,n\}}\big)^{1/3}\big( 2^{-2m -j-2l}\big)^{2/3} \big\} \\
 & \quad \times   2^{m -j- r   +l+2\epsilon M_t}   \| \mathfrak{m}(\cdot, \zeta)\|_{\mathcal{S}^\infty} \sup_{s\in [0,t]}\|B^{\tilde{m}}_{\tilde{k};\tilde{j}, \tilde{l}}(s,\cdot)\|_{L^\infty}\\
&\lesssim 2^{   2\epsilon M_t} \| \mathfrak{m}(\cdot, \zeta)\|_{\mathcal{S}^\infty} \sup_{s\in [0,t]}\|B^{\tilde{m}}_{\tilde{k};\tilde{j}, \tilde{l}}(s,\cdot)\|_{L^\infty} \\
&\quad \times \min\{2^{m/2  +\min\{p,n\}},  2^{-2j/3 +2n/3-2r/3} \}. \\
\end{split}
\ee 

If $(\tilde{m},\tilde{k},\tilde{j},\tilde{l})\in \mathcal{S}_1(t)$, then from the first estimate in  \eqref{aug4eqn10}  in Proposition \ref{meanLinfest} and the estimate  \eqref{july23eqn62}, we have
\be\label{2022feb11eqn97}
\begin{split}
 \sum_{(\tilde{m},\tilde{k},\tilde{j},\tilde{l})\in \mathcal{S}_1(t)} | {H}^{m,i;p,q;\tilde{m}, \tilde{k},\tilde{j},\tilde{l} }_{k,j ; n,l,r }(t, x , \zeta )|&\lesssim \min\{2^{m/2  + n},  2^{-2j/3} \}  2^{  13\epsilon M_t} \| \mathfrak{m}(\cdot, \zeta)\|_{\mathcal{S}^\infty} \\
 &\quad \times  \big( 2^{ 2 \tilde{\alpha}_t  M_t }   + 2^{7M_t/6+\tilde{\alpha}_t M_t/4} \big)\\
  &\lesssim  2^{(1-20\epsilon)M_{t^{\star}}}\| \mathfrak{m}(\cdot, \zeta)\|_{\mathcal{S}^\infty},\\
 \end{split}
\ee
where we used the fact that $r\geq \vartheta^\star_0\geq n+\epsilon M_t$ and  the assumption that  either  $n<  (1-2\alpha^{\star}-35\epsilon)M_{t^{\star}}     $ or  $    j\geq   (1/2+3\iota + 55\epsilon) M_{t^{\star}} $. 

If  $(\tilde{m},\tilde{k},\tilde{j},\tilde{l})\in \mathcal{S}_2(t)$, then after combining the estimates \eqref{july23eqn61}  and  \eqref{july23eqn62},  and the estimate  \eqref{aug4eqn10}  in Proposition \ref{meanLinfest}, we have
\be\label{2022feb11eqn96}
\begin{split}
\sum_{(\tilde{m},\tilde{k},\tilde{j},\tilde{l})\in \mathcal{S}_2(t)}  &| {H}^{m,i;p,q;\tilde{m}, \tilde{k},\tilde{j},\tilde{l} }_{k,j ; n,l,r }(t, x , \zeta )|\\
&\lesssim \sum_{(\tilde{m},\tilde{k},\tilde{j},\tilde{l})\in  \mathcal{S}_2(t) }  2^{- \vartheta^{\star}_0  + 5\epsilon M_t}   \big(  2^{-2j/3+2n/3} \sup_{s\in [0,t]}\|B^{\tilde{m}}_{\tilde{k};\tilde{j}, \tilde{l}}(s,\cdot)\|_{L^\infty} \big)^{1/2} \\
&\quad \times  \big( 2^{ -n/3} 2^{2(j+l)/3 +2\tilde{\alpha}_t M_t/3} \sup_{s\in [0,t]}\|B^{\tilde{m}}_{\tilde{k};\tilde{j}, \tilde{l}}(s,\cdot)\|_{L^2} \big)^{1/2}\| \mathfrak{m}(\cdot, \zeta)\|_{\mathcal{S}^\infty} \\
  & \lesssim2^{ - \vartheta^{\star}_0  +  20\epsilon M_t +4 {\alpha}^{\star}  M_{t^{\star}}/3 }    \| \mathfrak{m}(\cdot, \zeta)\|_{\mathcal{S}^\infty} \\
 &\lesssim  2^{(1-20\epsilon)M_{t^{\star}}}\| \mathfrak{m}(\cdot, \zeta)\|_{\mathcal{S}^\infty}.\\
 \end{split}
\ee

Therefore,  after combining the above two estimates, we have 
\be\label{aug5eqn80}
 \sum_{(\tilde{m},\tilde{k},\tilde{j},\tilde{l})\in \mathcal{S}_1(t) \cup \mathcal{S}_2(t) }   | {H}^{m,i;p,q;\tilde{m}, \tilde{k},\tilde{j},\tilde{l} }_{k,j ;n,l,r}(t, x, \zeta)|\lesssim   2^{(1-20\epsilon)M_{t^{\star}}}\| \mathfrak{m}(\cdot, \zeta)\|_{\mathcal{S}^\infty}.
\ee

 \medskip

\noindent \textbf{Case 3.}\qquad If  $m+k\geq -2n+ \epsilon M_t  $,   $m+p\geq -k-n+2\epsilon M_t$,   and $|  x_{\bot}|\geq 2^{m+p-10}$.

\medskip

 Similar to the obtained estimate \eqref{aug5eqn19}, we have
\be\label{aug4eqn60}
\begin{split}
&\big| {H}^{m,i;p,q;\tilde{m}, \tilde{k},\tilde{j},\tilde{l} }_{k,j ; n,l,r }(t,x, \zeta )\big| \\
& \lesssim  2^{m-j- \vartheta^{\star}_0  +l+2\epsilon M_t} \big(\frac{2^{-p-q} }{2^{m}|x_{\bot}|} \min\{2^{3j+2l}, 2^{j+2\tilde{\alpha}_t M_t}\} \big)^{1/2}\| \mathfrak{m}(\cdot, \zeta)\|_{\mathcal{S}^\infty}\\
&\quad \times  \big( 2^{m+2p+q  }  \min\{2^{3j+2l}, 2^{j+2\tilde{\alpha}_t M_t}\} \big)^{1/2}  \big( \sup_{s\in [0,t]}\|B^{\tilde{m}}_{\tilde{k};\tilde{j}, \tilde{l}}(s,\cdot)\|_{L^2}\big) \\
& \lesssim 2^{ - \vartheta^{\star}_0  + m/2+ l+2\epsilon M_t}\min\{2^{2j+2l}, 2^{2\tilde{\alpha}_t M_t}\}\| \mathfrak{m}(\cdot, \zeta)\|_{\mathcal{S}^\infty}\sup_{s\in [0,t]}\|B^{\tilde{m}}_{\tilde{k};\tilde{j}, \tilde{l}}(s,\cdot)\|_{L^2}. 
\end{split}
\ee
 Alternatively, if we put the localized electric field in $L^\infty_x$, then from the volume of support of $v,\omega$ and the estimate  \eqref{march18eqn31} in Lemma \ref{conservationlawlemma},   we have
 \be\label{aug4eqn64}
\begin{split}
|{H}^{m,i;p,q;\tilde{m}, \tilde{k},\tilde{j},\tilde{l} }_{k,j ; n,l,r }(t, x, \zeta )|&\lesssim \sup_{s\in [0,t]}\|B^{\tilde{m}}_{\tilde{k};\tilde{j}, \tilde{l}}(s,\cdot)\|_{L^\infty}  \| \mathfrak{m}(\cdot, \zeta)\|_{\mathcal{S}^\infty}
 2^{m-j- r  +l+2\epsilon M_t}\\ 
 &\quad \times \min\{    2^{m +2\min\{p,n\} } 2^{3j+2\min\{l,r\}} , 2^{-2m-j-2l}  \}\\
  &  \lesssim  \sup_{s\in [0,t]}\|B^{\tilde{m}}_{\tilde{k};\tilde{j}, \tilde{l}}(s,\cdot)\|_{L^\infty} 2^{ 2\epsilon M_t}
 \| \mathfrak{m}(\cdot, \zeta)\|_{\mathcal{S}^\infty}\\
 &\quad \times \min\{2^{m/2+\min\{p,n\} },    2^{-2j/3 }, 2^{-m/2+n/3-4j/3-l/3-r}\} . 
\end{split}
\ee

If $(\tilde{m},\tilde{k},\tilde{j},\tilde{l})\in \mathcal{S}_1(t)$, then from the obtained estimate  \eqref{aug4eqn64}  and the first estimate in \eqref{aug4eqn10} in Proposition \ref{meanLinfest}, we have
\be\label{aug4eqn74}
\begin{split}
 \sum_{(\tilde{m},\tilde{k},\tilde{j},\tilde{l})\in  \mathcal{S}_1(t) } & | {H}^{m,i;p,q;\tilde{m}, \tilde{k},\tilde{j},\tilde{l} }_{k,j ; n,l,r }(t, x, \zeta ) |\\
 &\lesssim   \min\{2^{m/2  + n},  2^{-2j/3} \} 2^{  13\epsilon M_t} \big( 2^{ 2 \tilde{\alpha}_t  M_t }   + 2^{7M_t/6+\tilde{\alpha}_t M_t/4} \big) \| \mathfrak{m}(\cdot, \zeta)\|_{\mathcal{S}^\infty}\\
&\lesssim   2^{(1-20\epsilon)M_{t^{\star}}}\| \mathfrak{m}(\cdot, \zeta)\|_{\mathcal{S}^\infty},
\end{split}
\ee
where, again,  we used the assumption that either  $n<  (1-2\alpha^{\star}-35\epsilon)M_{t^{\star}}    $ or  $    j\geq   (1/2+3\iota + 55\epsilon) M_{t^{\star}}  $.

If $(\tilde{m},\tilde{k},\tilde{j},\tilde{l})\in \mathcal{S}_2(t)$, then after combining the estimates \eqref{aug4eqn60}  and \eqref{aug4eqn64},  and the second estimate \eqref{aug4eqn10}  in Proposition \ref{meanLinfest}, we have 
\be\label{2022feb11eqn105}
\begin{split}
 \sum_{(\tilde{m},\tilde{k},\tilde{j},\tilde{l})\in  \mathcal{S}_2(t) } &\big| {H}^{m,i;p,q;\tilde{m}, \tilde{k},\tilde{j},\tilde{l} }_{k,j ; n,l,r }(t, x, \zeta)\big|\\ & \lesssim \sum_{(\tilde{m},\tilde{k},\tilde{j},\tilde{l})\in  \mathcal{S}_2(t) } \| \mathfrak{m}(\cdot, \zeta)\|_{\mathcal{S}^\infty}  2^{- \vartheta^{\star}_0  + n/6+ \tilde{\alpha}_t M_t/3+5\epsilon M_t} \\
 &\quad \times  \big(\sup_{s\in [0,t]}\|B^{\tilde{m}}_{\tilde{k};\tilde{j}, \tilde{l}}(s,\cdot)\|_{L^2}\big)^{1/2} \big(\sup_{s\in [0,t]}\|B^{\tilde{m}}_{\tilde{k};\tilde{j}, \tilde{l}}(s,\cdot)\|_{L^\infty}\big)^{1/2} \\
&   \lesssim 2^{ - \vartheta^{\star}_0  + \tilde{\alpha}_t M_t/3+35\epsilon M_t} \big( 2^{\tilde{\alpha}_t M_t}   + 2^{7M_t/12+ \tilde{\alpha}_t M_t/8} \big) \| \mathfrak{m}(\cdot, \zeta)\|_{\mathcal{S}^\infty}\\ &\lesssim   2^{(1-20\epsilon)M_{t^{\star}}}\| \mathfrak{m}(\cdot, \zeta)\|_{\mathcal{S}^\infty}.\\
\end{split}
\ee
To sum up,   we have
\be\label{aug5eqn81}
 \sum_{(\tilde{m},\tilde{k},\tilde{j},\tilde{l})\in \mathcal{S}_1(t) \cup \mathcal{S}_2(t) }   | {H}^{m,i;p,q;\tilde{m}, \tilde{k},\tilde{j},\tilde{l} }_{k,j ; n,l,r  } (t, x, \zeta )|\lesssim   2^{(1-20\epsilon)M_{t^{\star}}}\| \mathfrak{m}(\cdot, \zeta)\|_{\mathcal{S}^\infty}. 
\ee

   \medskip

\noindent \textbf{Case 4.}\qquad  If $m+k\geq -2n+ \epsilon M_t  $,   $m+p\geq -k-n+2\epsilon M_t$,   and $|  x_{\bot}|\leq  2^{m+p-10}$.

\medskip 

 Note that, for this case, we have $|   x_{\bot}- y_{\bot} + (t-s) \omega_{\bot}|\sim 2^{m+p}$ for any $y\in B(0, 2^{-k-n+\epsilon M_t/2})\subset \R^3$,  and the volume of support of $\omega$ is bounded by $2^{2\min\{p,n\}}$. Based on the possible size of $n$, we split into two sub-cases as follows.

 \medskip

  \textbf{Sub-case 4A.}\qquad  If $n\geq  -(\alpha^{\star}+3\iota+60\epsilon) M_t $.

\medskip

Similar to the obtained estimate  \eqref{aug5eqn44},  we have 
\be\label{aug4eqn63}
\begin{split}
| {H}^{m,i;p,q;\tilde{m}, \tilde{k},\tilde{j},\tilde{l} }_{k,j ;n,l,r}(t, x, \zeta )| &
\lesssim   2^{m -j -r+l + 2\epsilon  M_t}\| \mathfrak{m}(\cdot, \zeta)\|_{\mathcal{S}^\infty}\sup_{s\in [0,t]}\|B^{\tilde{m}}_{\tilde{k};\tilde{j}, \tilde{l}}(s,\cdot)\|_{L^2}\\
&\quad \times \big( \frac{2^{k+\max\{n,p\}}2^{3j+2\min\{l,r\}}}{2^{m+p}}  \big)^{1/2}  \big(    2^{-2m-j-2l}   \big)^{1/2}\\
&\lesssim     2^{-m/2+k/2  +\max\{n,p\}/2-p/2+ 2\epsilon  M_t}\| \mathfrak{m}(\cdot, \zeta)\|_{\mathcal{S}^\infty} \sup_{s\in [0,t]}\|B^{\tilde{m}}_{\tilde{k};\tilde{j}, \tilde{l}}(s,\cdot)\|_{L^2} .
\end{split}
\ee

Moreover, note that the obtained estimate  \eqref{aug4eqn64}  is also valid for the case we are considering. 
From the obtained estimate  \eqref{aug4eqn64}  and the  first estimate in \eqref{aug4eqn10}  in Proposition \ref{meanLinfest}, the following estimate holds  if $(\tilde{m},\tilde{k},\tilde{j},\tilde{l})\in \mathcal{S}_1(t)$, we have 
 \be\label{aug4eqn77}
 \begin{split}
 \sum_{(\tilde{m},\tilde{k},\tilde{j},\tilde{l})\in \mathcal{S}_1(t)}  \big| {H}^{ m,i;p,q; \tilde{m}, \tilde{k},\tilde{j},\tilde{l} }_{k,j ; n,l,r  }(t, x, \zeta)\big| &
\lesssim \min\{2^{m+n}, 2^{-2j/3}\}  2^{ 13\epsilon M_t}  \| \mathfrak{m}(\cdot, \zeta)\|_{\mathcal{S}^\infty}\\
&\quad \times \big( 2^{ 2 \tilde{\alpha}_t  M_t }   + 2^{7M_t/6+\tilde{\alpha}_t M_t/4} \big)\\
& \lesssim   2^{(1-20\epsilon)M_{t^{\star}}}\| \mathfrak{m}(\cdot, \zeta)\|_{\mathcal{S}^\infty}.\\
\end{split}
\ee

 After combining the estimates  \eqref{aug4eqn63}  and  \eqref{aug4eqn64},  and the estimate  \eqref{aug4eqn10}  in Proposition \ref{meanLinfest}, the following estimate holds if   $(\tilde{m},\tilde{k},\tilde{j},\tilde{l})\in \mathcal{S}_2(t)$, we have 
 \be\label{aug4eqn76}
 \begin{split}
 \sum_{(\tilde{m},\tilde{k},\tilde{j},\tilde{l})\in \mathcal{S}_2(t)} &\big| {H}^{m,i;p,q;\tilde{m}, \tilde{k},\tilde{j},\tilde{l} }_{k,j ;n,l,r}(t, x, \zeta )\big| \\
&\lesssim \sum_{(\tilde{m},\tilde{k},\tilde{j},\tilde{l})\in \mathcal{S}_2(t)} \| \mathfrak{m}(\cdot, \zeta)\|_{\mathcal{S}^\infty}  2^{ 2\epsilon  M_t}\big( \sup_{s\in [0,t]}\|B^{\tilde{m}}_{\tilde{k};\tilde{j}, \tilde{l}}(s,\cdot)\|_{L^\infty} 
2^{m/2+\min\{n,p\}} \big)^{1/2} \\
 &\quad  \times  \big(2^{ -m/2+  k/2+\max\{l,p\}+ \max\{n,p\}/2-p/2 } \sup_{s\in [0,t]}\|B^{\tilde{m}}_{\tilde{k};\tilde{j}, \tilde{l}}(s,\cdot)\|_{L^2}\big)^{1/2} \\
 &  \lesssim 2^{ (k+2n)/4}   2^{ {\alpha}^{\star}  M_{t^{\star}} + 13\epsilon M_t}  \| \mathfrak{m}(\cdot, \zeta)\|_{\mathcal{S}^\infty}.\\
 \end{split}
 \ee

  \medskip

  \textbf{Sub-case 4B.}\qquad  If $ n\leq   -(\alpha^{\star}+3\iota +60\epsilon) M_t $.

\medskip

Recall     \eqref{aug5eqn1}. From the estimate  \eqref{sep21eqn31}  in Lemma \ref{bulkroughpoi},     and the volume of support of $\omega$ and $v$, we have 
\be\label{aug4eqn75}
\begin{split}
  \big|{H}^{m, i ;p,q  }_{k,j ;n,l,r} (t,   x,  \zeta)\big|&\lesssim   \| \mathfrak{m}(\cdot, \zeta)\|_{\mathcal{S}^\infty}  2^{m-j- r +l+ 9\epsilon  M_t}\min\{2^{m+2\min\{n,p\}+3j+2\min\{l,r\}}, 2^{-2m-j-2l} \}  \\
  &\quad \times ( 2^{2\tilde{\alpha}_t M_t } 2^{-(m+p)/2} + 2^{-(m+p)/4}  2^{5M_t/4+\tilde{\alpha}_t M_t/4 } ) \\
  & \lesssim \| \mathfrak{m}(\cdot, \zeta)\|_{\mathcal{S}^\infty} 2^{m-j- r +l+ 9\epsilon  M_t}\big( 2^{m+2\min\{n,p\}+3j+2\min\{l,r\}}\big)^{1/2}\\
&\quad \times  \big( 2^{-2m-j-2l}\big)^{1/2} ( 2^{2\tilde{\alpha}_t M_t } 2^{-(m+p)/2} + 2^{-(m+p)/4}  2^{5M_t/4+\tilde{\alpha}_t M_t/4 } )\\
&  \lesssim \| \mathfrak{m}(\cdot, \zeta)\|_{\mathcal{S}^\infty}  2^{   9\epsilon  M_t}( 2^{2\tilde{\alpha}_t M_t  } 2^{\min\{n,p\}/2}  +     2^{5M_t/4+\tilde{\alpha}_t M_t/4} 2^{3\min\{n,p\}/4}    ) \\
  & \lesssim 2^{(1-20\epsilon)M_{t^{\star}}}\| \mathfrak{m}(\cdot, \zeta)\|_{\mathcal{S}^\infty}. 
  \end{split}
\ee

Recall  \eqref{aug10eqn31}. To sum up,   from   \eqref{aug4eqn75},  \eqref{aug4eqn76}, and   \eqref{aug4eqn77},   we have
\be\label{aug5eqn82}
\begin{split}
   \big|{H}^{m, i ;p,q  }_{k,j ;n,l,r} (t,   x,  \zeta)\big|&\lesssim \| \mathfrak{m}(\cdot, \zeta)\|_{\mathcal{S}^\infty} \big[ 2^{(1-20\epsilon)M_{t^{\star}}} 
  +  2^{  100\epsilon M_{t^{\star}} } \mathbf{1}_{n\geq  -(\alpha^{\star}+3\iota+60\epsilon) M_{t^{\star}} } \\
  &\qquad \times \min\{2^{(k+2n)/2 }2^{(2{\alpha}^{\star} -1)M_{t^{\star}}}  , 2^{(k+4n)/2 +(1+6\iota)M_{t^{\star}}} \} \big] . 
\end{split}
\ee

Recall  \eqref{aug5eqn1} and   \eqref{aug10eqn31}. The desired estimate  \eqref{aug3eqn14}  holds from the   estimates   \eqref{aug4eqn52},  \eqref{aug5eqn80},   \eqref{aug5eqn81}, and  \eqref{aug5eqn82}.

\end{proof}
To sum up, the main results of this subsection can be summarized in the following Proposition.

\begin{proposition}\label{finalestfirst}
 Under   the assumption of Theorem \ref{mainresultsfirstpart},      we have
  \be\label{sep27eqn1}
  \begin{split}
&   \big\|\widetilde{T}_{k,j;n}^{bil;\mu,0}( \mathfrak{m}, E)(t,x, \zeta )  + \hat{\zeta}\times  \widetilde{T}_{k,j;n}^{bil;\mu,0 }(\mathfrak{m}, B )(t,x, \zeta)\big\|_{L^\infty_x} + \big\|  \mathfrak{H}_{k,j;n}^{\mu,0}( \mathfrak{m})(t, x, \zeta)\big\|_{L^\infty_x} \\
&  \lesssim  \| \mathfrak{m}(\cdot, \zeta)\|_{\mathcal{S}^\infty}  \big[   2^{(1-19\epsilon)M_{t^{\star}} }   +  2^{  121\epsilon M_{t^{\star}} }  \mathbf{1}_{n\geq  -(\alpha^{\star}+3\iota+60\epsilon) M_{t^{\star}} } 
 \\
&\quad  \times  \min\{ 2^{(k+2n)/2+ (\alpha^{\star}+3\iota) M_{t^{\star}} }  , 2^{(k+4n)/2 +(1+6\iota)M_{t^{\star}}} \} \big],\\
& \big\|\widetilde{T}_{k,j;n}^{bil;\mu,1}( \mathfrak{m}, E)(t,x, \zeta )  + \hat{\zeta}\times  \widetilde{T}_{k,j;n}^{bil;\mu,1 }(\mathfrak{m}, B )(t,x, \zeta)\big\|_{L^\infty_x}  + \big\|  \mathfrak{H}_{k,j;n}^{\mu,1}( \mathfrak{m})(t, x, \zeta)\big\|_{L^\infty_x} \\
    &\lesssim  \| \mathfrak{m}(\cdot, \zeta)\|_{\mathcal{S}^\infty}  \big[   2^{(1-19\epsilon)M_{t^{\star}} }  +  2^{   121\epsilon M_{t^{\star}} }    \mathbf{1}_{n\geq  -(\alpha^{\star}+3\iota+60\epsilon) M_{t^{\star}} } \\
  &\quad  \times \min\{   2^{(k+2n)/2 + (2{\alpha}^{\star} -1)M_{t^{\star}}}  , 2^{(k+4n)/2 +(1+6\iota)M_{t^{\star}}} \} \big].
\end{split}
 \ee
 
\end{proposition}
 \begin{proof}

 Recall \eqref{2024nov13hyperbolic} and \eqref{oct1eqn1}. Since we did the same normal form transformation for the same object for the cases where  $i\in\{0,1\}$, we have
 \be\label{2024nov14eqn32}
 \widetilde{T}_{k,j;n}^{bil;\mu ,i }( \mathfrak{m}, E)(t,x, \zeta )+ \hat{\zeta}\times  \widetilde{T}_{k,j;n}^{bil;\mu,i }(\mathfrak{m}, B )(t,x, \zeta)= \mathfrak{H}_{k,j;n}^{\mu,i}( \mathfrak{m})(t, x, \zeta).
 \ee

  Our desired first estimate in   \eqref{sep27eqn1}  holds after combining  the estimate  \eqref{sep5eqn99},  the estimate \eqref{sep8eqn51},    the estimate  \eqref{2022feb8eqn58}  in Lemma  \ref{largeregime1}, the estimates  \eqref{aug5eqn33}  and  \eqref{2022feb11eqn151}  in Lemma \ref{largeregime2},   and the estimate 
  \eqref{aug3eqn14}  in Lemma \ref{mainlemmasecpart}.

Moreover, recall the definition of the cutoff function $\varphi_{j,n}^1 (v, \zeta)$ in  \eqref{sep4eqn6}, we know that   the the size of $j$ is greater than $  (1/2+3\iota + 55\epsilon ) M_{t^{\star}}  $. Hence, our desired second estimate  in  \eqref{sep27eqn1}  holds from  the estimate  \eqref{sep8eqn51},    the estimate  \eqref{2022feb8eqn58}  in Lemma  \ref{largeregime1}, the estimate   \eqref{2022feb11eqn151}  in Lemma \ref{largeregime2},   and the estimate 
 \eqref{aug3eqn14}  in Lemma \ref{mainlemmasecpart}.

\end{proof}

 \subsection{$L^\infty$-type estimates of   the acceleration force      in the small angle region}\label{linfsmallre}

 In this  section, we estimate  the small angle region case. More precisely, we estimate the $L^\infty_x$-norm of $T_{k,j;n}^{\mu,i}(\mathfrak{m}, E)(t,x, \zeta ) + \hat{\zeta}\times T_{k,j;n}^{\mu,i}(\mathfrak{m}, B)(t,x, \zeta ), i\in \{2,3, 4\} $, see  \eqref{sep17eqn32}. 

Recall the decomposition of $T_{k,j;n}^{\mu,i}(\mathfrak{m}, E)(t,x, \zeta ) + \hat{\zeta}\times T_{k,j;n}^{\mu,i}(\mathfrak{m}, B)(t,x, \zeta ) $   as presented  in  \eqref{sep18eqn50}  in Lemma \ref{locdeclemm}. We first consider the case   $m+k$ is relatively small, for which we use the formulas in \eqref{sep18eqn31}. 
\begin{lemma}\label{smallfrepartII}
Let $i\in \{2,3,4\},  ( l,r)\in \mathcal{B}_i $, see  \eqref{sep18eqn50}.   Under the assumption of Theorem \ref{mainresultsfirstpart}, the following estimates hold  if  $m+k\leq -2 l +4\epsilon M_t$, 
\be\label{sep19eqn71}
\begin{split}
&\big\|  T_{k,j;n,l,r}^{\mu,m,2}(\mathfrak{m},E)(t,x, \zeta) + \hat{\zeta}\times T_{k,j;n,l,r}^{\mu,m, 2}(\mathfrak{m},B)(t,x, \zeta)  \big\|_{L^\infty_x}\\
&\lesssim    \| \mathfrak{m}(\cdot, \zeta)\|_{\mathcal{S}^\infty}\big[2^{(1-19.5\epsilon)M_{t^{\star}}} + 2^{127.5\epsilon  M_{t^{\star}} } \mathbf{1}_{n\geq (1-2  {\alpha}^{\star}-90\epsilon  )M_{t^{\star}}  } \\
&\qquad \times \min\{2^{(k+2n)/2+ {\alpha}^{\star} M_{t^{\star}} }, 2^{(k+4n)/2+(1+3\iota)M_{t^{\star}} }\}  \big], 
\end{split}
\ee
\be\label{2022feb12eqn11}
\begin{split}
\sum_{i=3,4}&\big\|  T_{k,j;n,l,r}^{\mu,m,i}(\mathfrak{m},E)(t,x, \zeta) + \hat{\zeta}\times T_{k,j;n,l,r}^{\mu,m, i}(\mathfrak{m},B)(t,x, \zeta)  \big\|_{L^\infty_x}\\
&\lesssim  \| \mathfrak{m}(\cdot, \zeta)\|_{\mathcal{S}^\infty} \big[  2^{(1- 19.5\epsilon)M_{t^{\star}}}  
 +   2^{127.5\epsilon  M_{t^{\star}} }   \mathbf{1}_{n\geq -(1/2+3\iota/2+40\epsilon)M_{t^{\star}} } \\
 &\qquad \times \min\{ 2^{(k+2n)/2+\alpha^{\star} M_{t^{\star}}}, 2^{(k+4n)/2+ (7/6+5\iota/2) M_{t^{\star}}} \}  \big].
 \end{split}
\ee
\end{lemma}
\begin{proof}
Recall  \eqref{sep18eqn31}. From the decomposition of the acceleration force in  \eqref{july1eqn11}, for $i\in\{2,3,4\},$ we have  
\be\label{sep19eqn81}
\begin{split}
&T_{k,j;n,l,r}^{\mu,m,i}(\mathfrak{m},E)(t,x, \zeta) + \hat{\zeta}\times T_{k,j;n,l,r}^{\mu,m, i}(\mathfrak{m},B)(t,x, \zeta)\\
& = H_{k,j;n,l,r}^{\mu,m,i;lin}(t, x,\zeta) +\sum_{a=1,2} H_{k,j;n,l,r}^{\mu,m,i;non,a}(t, x, \zeta),
\end{split}
\ee 
where
\be\label{sep19eqn82}
\begin{split}
 H_{k,j;n,l,r}^{\mu,m,i;lin}(t, x,\zeta)  &:= \int_{0}^t \int_{\R^3} \int_{\R^3} \int_{\mathbb{S}^2} \varphi_{m;-10M_t }(t-s) \big[(t-s) \big( \mathfrak{K}^{\mu, E}_{k;n}(y,\omega, v,  \zeta)\\
 &\quad   + \hat{\zeta}\times   \mathfrak{K}^{\mu, B}_{k;n}(y,\omega, v, \zeta)\big) +\mathfrak{K}^{err;\mu,E}_{k;n}(y, v, \zeta)  + \hat{\zeta}\times \mathfrak{K}^{err;\mu,B}_{k;n}(y, v, \zeta)\big] \\
  &\quad \times f(s, x-y+(t-s)\omega, v)   \varphi_{j,n}^{i; r}(v, \zeta)  \varphi_{l; r}(\tilde{v}+\omega )   d\omega dy d v ds, \\
  H_{k,j;n,l,r}^{\mu,m,i;non,a}(t, x, \zeta)&=  
\int_{0}^t \int_{\R^3} \int_{\R^3} \int_{\mathbb{S}^2}  EB^a(t,s,x-y +(t-s)\omega,\omega, v) \\
 &\quad\cdot  \nabla_v \big( (t-s) \mathfrak{H}^{\mu,E,i}_{k,j;n,l,r}(y,\omega, v, \zeta)+\mathfrak{H}^{err;\mu,E,i}_{k,j;n,l,r}(y, v, \zeta)\big) \\
 &\quad \times  f(s, x-y+(t-s)\omega, v)  \varphi_{m;-10M_t }(t-s)   d \omega dy d v ds.
\end{split}
\ee
 
  Similar to what we did in  \eqref{sep16eqn13}, we localize frequency around $\omega$     and do dyadic decompositions for $ \omega_{\bot}$ for the nonlinear part as follows, 
\be\label{sep22eqn1}
\begin{split}
H_{k,j;n,l,r}^{\mu,m,i;lin}(t, x,\zeta)&= \sum_{\star\in\{ess,err\}}H_{k,j;n,l,r;\star}^{\mu,m,i;lin}(t, x,\zeta),\\
 H_{k,j;n,l,r}^{\mu,m,i;non,a}(t, x,\zeta)&=\sum_{\begin{subarray}{c}
\star\in \{ess, err\}\\
p\in[-10M_t,2]\cap \Z
\end{subarray}}  H_{k,j;n,l,r;p;\star}^{\mu,m,i;non,a}(t, x,\zeta),
\end{split}
\ee
where
\be\label{2022feb8eqn1}
\begin{split}
H_{k,j;n,l,r;\star}^{\mu,m,i;lin}(t, x,\zeta)&= \int_{0}^t \int_{\R^3} \int_{\R^3} \int_{\mathbb{S}^2}\big[ (t-s) \mathcal{K}_{k,n,l,r}^{\star;\mu,m}(y, v,\omega,  \zeta) +\mathcal{K}_{k,n,l,r}^{\star;err,\mu,m}(y, v,\omega,\zeta) \big] \\
&\quad \times \varphi_{j,n}^{i; r}(v, \zeta)    \varphi_{l; r}(\tilde{v}+\omega )   \varphi_{m;-10M_t }(t-s)    \\
&\quad \times  f(s, x-y+(t-s)\omega, v) d\omega dy d v ds, \\
H_{k,j;n,l,r;p;\star}^{\mu,m,i;non,a}(t, x,\zeta)&= \int_{0}^t \int_{\R^3} \int_{\R^3} \int_{\mathbb{S}^2}  \varphi_{m;-10M_t }(t-s)    EB^a(t,s,x-y +(t-s)\omega,\omega, v)\\
&\quad \cdot  \nabla_v \big( (t-s) \mathcal{H}^{\star;\mu,E,i}_{k,j; n,l,r }(y,\omega, v, \zeta)  +\mathcal{H}^{\star;err;\mu,E,i}_{k,j;n,l,r }(y, v,\omega,  \zeta)\big) \\
 &\quad \times f(s, x-y+(t-s)\omega, v) \varphi_{p;-10M_t}(  \omega_{\bot})  d \omega dy d v ds,\\
\end{split}
\ee
and  the kernels appeared above are defined as follows, 
\be\label{2022feb8eqn11}
\begin{split}
\mathcal{K}_{k,n,l,r}^{\star;\mu,m}(y, v, \omega, \zeta) & := \int_{\R^3} e^{iy \cdot \xi} \mathfrak{m}(\xi, \zeta ) \varphi_k(\xi)\varphi_{n;-M_t}\big( \tilde{\xi} + \mu \tilde{\zeta} \big)  (i \omega \cdot \tilde{\xi}-\mu) \\
&\qquad \times \varphi_{c(m,k,l)}^{\star}(\omega , \xi)  \big( m_{E}^{0}(\xi, v) + \hat{\zeta}\times m_{B}^{0}(\xi, v) \big) d\xi, \\ 
\mathcal{K}_{k,n,l,r}^{\star;err;\mu,m}(y, v, \omega, \zeta) &:= \int_{\R^3} e^{iy \cdot \xi} \mathfrak{m}(\xi, \zeta ) \varphi_k(\xi)\varphi_{n;-M_t}\big(  \tilde{\xi} + \mu \tilde{\zeta} \big) |\xi|^{-1} \\
&\qquad \times \varphi_{c(m,k,l)}^{ \star}(\omega , \xi) \big( m_{E}^{0}(\xi, v) + \hat{\zeta}\times m_{B}^{0}(\xi, v) \big)d\xi, \\
  \mathcal{H}^{\star; \mu, E,i }_{k,j;n,l,r }( \mathfrak{m})(y,\omega, v, \zeta )&:=\int_{\R^3 } e^{i y\cdot \xi }  \mathfrak{m}(\xi,  \zeta) \varphi_k(\xi)\varphi_{n;-M_t}\big(\tilde{\xi} +\mu  \tilde{\zeta}\big)  \hat{v}   (   {i \omega \cdot \tilde{\xi} } -\mu)  \\
  &\qquad \times   \varphi_{j,n}^{i; r}(v, \zeta)  \varphi_{l; r}(\tilde{v}+\omega )  \varphi_{c(m,k,l)}^{\star}(\omega , \xi)  d \xi, \\ 
\mathcal{H}^{\star;err;\mu, E,i}_{k,j;n,l,r}( \mathfrak{m})(y,\omega, v, \zeta)& :=i \int_{\R^3 } e^{i y\cdot \xi }  \mathfrak{m}(\xi, \zeta) \varphi_k(\xi)\varphi_{n;-M_t}\big(\tilde{\xi} +\mu  \tilde{\zeta}\big) \hat{v}  |\xi|^{-1} \\
&\qquad \times    \varphi_{j,n}^{i; r}(v, \zeta)  \varphi_{l; r}(\tilde{v}+\omega ) \varphi_{c(m,k,l)}^{\star}(\omega , \xi) d \xi,\\ 
\end{split}
\ee
 \[
 \begin{split}
 c(m,k,l)&:=\min\{\max\{-m/2-k/2, -m-k-l\}+\epsilon M_t/20, 2\},\\  
\varphi_{c(m,k,l)}^{ess}(\omega , \xi)&:=\psi_{\leq c(m,k,l)}(\omega \times \tilde{\xi}),\qquad  \varphi_{c(m,k,l)}^{err}(\omega , \xi):=\psi_{> c(m,k,l)}(\omega \times \tilde{\xi}).
\end{split}
 \]

 Since $m+k\leq - 2l+4\epsilon M_t$, we have $c(m,k,l)\leq -m-k-l+3\epsilon M_t. $   Similar to the obtained estimate for the error type estimate in  \eqref{sep16eqn30}, by doing integration by parts in ``$\omega$'' many times, we have 
\be\label{sep20eqn21}
\big|H_{k,j;n,l,r;err}^{\mu,m,i;lin}(t, x,\zeta)\big|\lesssim \| \mathfrak{m}(\cdot, \zeta)\|_{\mathcal{S}^\infty}. 
\ee
Moreover, recall the rough estimate of the electromagnetic field   \eqref{july10eqn89}  in Proposition \ref{Linfielec}. By doing integration by parts in ``$\omega$'' many times,  similar to the obtained estimate for the error type estimate in  \eqref{sep16eqn30}, the following estimate holds, 
\be\label{sep20eqn27}
\sum_{a=1,2}\big| \sum_{p\in[-10M_t,2]\cap \Z} H_{k,j;n,l,r;p;err}^{\mu,m,i;non,a}(t, x,\zeta)\big|\lesssim \| \mathfrak{m}(\cdot, \zeta)\|_{\mathcal{S}^\infty}. 
\ee

Now, we focus on the estimate of essential parts $H_{k,j;n,l;ess}^{\mu,m,i;lin}(t, x,\zeta)$ and $H_{k,j;n,l;p;ess}^{\mu,m,i;non,a}(t, x,\zeta)$. Additionally, again from  the rough estimate of the electromagnetic field   \eqref{july10eqn89}  in Proposition \ref{Linfielec} and the volume of support of $v$ and $\omega$, we can rule out the case $p=-10M_t$.  It suffices to consider the case $p\in (-10M_t, 2]\cap \Z.$

To analyze the above terms individually,  we proceed in steps as follows. 

\medskip

\noindent \textbf{Step 1.}\qquad The estimate of $H_{k,j;n,l,r ;ess}^{\mu,m,i;lin}(t, x,\zeta)  $, $i=2 .$

\medskip

 Due to the  cutoff functions $\varphi^i_{j,n}(v, \zeta)$, see  \eqref{sep4eqn6}, we only have to consider the case $n\leq -  \epsilon M_t +10. $ Moreover,  we have $   2^{n+\epsilon M_t}\lesssim |\tilde{v}-\tilde{\zeta}| \lesssim 2^{-\epsilon M_t/2}$, which implies that $\max\{l, c(m,k,l)\}\geq n +3\epsilon M_t/4$. It implies further that $m+k+l\leq -n+3\epsilon M_t.$

Recall  \eqref{2022feb8eqn11}  and the   result of computation  of symbols  in  \eqref{sep5eqn101}.  After doing integration by parts in $\xi$ along $\zeta$ direction and directions perpendicular to $\zeta$ many times, the following estimate holds, 
 \be 
 \begin{split}
 &\big|  \mathcal{K}_{k,n,l,r}^{ess;\mu,m}(y, v, \omega, \zeta)  \big| + 2^k \big| \mathcal{K}_{k,n,l,r}^{ess;err;\mu,m}(y, v, \omega, \zeta)\big|\\
  &\lesssim 2^{4k+ 2n   }(2^{2l} + 2^{2c(m,k,l)}) \| \mathfrak{m}(\cdot, \zeta)\|_{\mathcal{S}^\infty}  (1+2^k|y\cdot \tilde{\zeta}|)^{-N_0^3} (1+2^{k+n }|y\times \tilde{\zeta}|)^{-N_0^3} .
  \end{split}
\ee

From the above estimate of kernels, the conservation law  \eqref{conservationlaw},   the volume of support of $\omega, v$, and the estimate  \eqref{nov24eqn41}  if $|  v_{\bot}|\geq 2^{(\alpha_t+\epsilon)M_t}$, we have
\be\label{2022feb12eqn15}
\begin{split}
\big|H_{k,j;n,l,r;ess}^{\mu,m,i;lin}(t, x, \zeta)\big|&\lesssim 2^{m }(2^{2l} + 2^{2c(m,k,l)}) (2^{m+k}+1) \| \mathfrak{m}(\cdot, \zeta)\|_{\mathcal{S}^\infty} \\
&\qquad\times  \min\{2^{3k+2n} 2^{2l-j},  2^{ 2l} 2^{j+2(\tilde{\alpha}_t+\epsilon)M_t}  \}\\
 &\lesssim   \| \mathfrak{m}(\cdot, \zeta)\|_{\mathcal{S}^\infty} 2^{ 6\epsilon M_t}\min\{2^{2k+2n-j},   2^{-k+j+2(\tilde{\alpha}_t+\epsilon)M_t}  \} \\
 & \lesssim  2^{8 \epsilon M_t} \min\{2^{4 \tilde{\alpha}_t M_t/3 +j/3 + 2n/3 }, 2^{(k+2n)/2+  \tilde{\alpha}_t M_t}\}\| \mathfrak{m}(\cdot, \zeta)\|_{\mathcal{S}^\infty}\\
& \lesssim  \| \mathfrak{m}(\cdot, \zeta)\|_{\mathcal{S}^\infty}\big[2^{(1-20\epsilon)M_{t^{\star}}} + 2^{60\epsilon  M_{t^{\star}} } \mathbf{1}_{n\geq (1-2  {\alpha}^{\star}-50\epsilon  )M_{t^{\star}}  } \\
&\quad \times \min\{2^{(k+2n)/2+ {\alpha}^{\star} M_{t^{\star}} }, 2^{(k+4n)/2+(1+3\iota)M_{t^{\star}} }\} \big]. \\
\end{split}
\ee

\medskip

\noindent \textbf{Step 2.}\qquad
 The estimate of $H_{k,j;n,l}^{\mu,m,i;lin}(t, x,\zeta)  $, $i=3,4.$
\medskip

 Due to the  cutoff functions $\varphi^i_{j,n}(v, \zeta), i\in\{3,4\}$, see  \eqref{sep4eqn6}, we have $|\tilde{v}- \tilde{\zeta}|\lesssim 2^{n+\epsilon M_t}. $  Based on the possible size of $m+k$, we split into two sub-cases as follows.

 \medskip

\textbf{Step 2A.}\qquad If $m+k\leq \min\{-2l, - l-n\} +4\epsilon M_t$. 

 \medskip

Note that, for this case, we have $c(m,k,l)\geq n -4\epsilon M_t.$ Recall     \eqref{sep5eqn101}.  After doing integration by parts in $\xi$ along $\zeta$ direction and directions perpendicular to $\zeta$ many times, we have  
\be 
\begin{split}
 &\big|  \mathcal{K}_{k,n,l,r}^{ess;\mu,m}(y, v, \omega, \zeta)  \big| + 2^k \big| \mathcal{K}_{k,n,l,r}^{ess;err;\mu,m}(y, v, \omega, \zeta)\big| \\
 & \lesssim 2^{4k+ 4n + 8\epsilon M_t } \| \mathfrak{m}(\cdot, \zeta)\|_{\mathcal{S}^\infty}  (1+2^k|y\cdot \tilde{\zeta}|)^{-N_0^3} (1+2^{k+n }|y\times \tilde{\zeta}|)^{-N_0^3} .
 \end{split}
\ee

From the above estimate of kernels, the conservation law  \eqref{conservationlaw},   the volume of support of $\omega, v$, and the estimate  \eqref{nov24eqn41}  if $|  v_{\bot}|\geq 2^{(\alpha_t+\epsilon)M_t}$, we have
\be\label{sep19eqn100}
\begin{split}
\big|H_{k,j;n,l,r;ess}^{\mu,m,i;lin}&(t, x, \zeta)\big|\\
&\lesssim 2^{m+2n +8\epsilon M_t}(2^{m+k}+1) \| \mathfrak{m}(\cdot, \zeta)\|_{\mathcal{S}^\infty} \min\{2^{3k+2n} 2^{2l-j},  2^{ 2l} 2^{j+2(\tilde{\alpha}_t+\epsilon)M_t}  \}\\
&\lesssim  \min\{ 2^{-2m-2l -j + 24\epsilon M_t}, 2^{2k+2n-j+16\epsilon M_t},    2^{ m+l+n} 2^{j+2 \tilde{\alpha}_t M_t+14\epsilon M_t}\}\| \mathfrak{m}(\cdot, \zeta)\|_{\mathcal{S}^\infty}\\
&\lesssim \| \mathfrak{m}(\cdot, \zeta)\|_{\mathcal{S}^\infty} \min\big\{ \big(2^{-2m-2l -j + 24\epsilon M_t}\big)^{1/3}\big(  2^{ m+l+n +j+2 \tilde{\alpha}_t M_t+14\epsilon M_t} \big)^{2/3},\\ 
&\qquad  \big(2^{-2m-2l -j + 24\epsilon M_t}\big)^{1/4}\big(2^{2k+2n-j+16\epsilon M_t}\big)^{1/4} \big(  2^{ m+l+n +j+2 \tilde{\alpha}_t M_t+14\epsilon M_t} \big)^{1/2}  \big\}\\
&\lesssim 2^{20\epsilon M_t} \min\{2^{4\tilde{\alpha}_t M_t/3 +j/3+2n/3} ,  2^{(k+2n)/2+\tilde{\alpha}_t M_t}  \}\| \mathfrak{m}(\cdot, \zeta)\|_{\mathcal{S}^\infty}\\
 &\lesssim  \| \mathfrak{m}(\cdot, \zeta)\|_{\mathcal{S}^\infty}\big[2^{(1-20\epsilon)M_{t^{\star}}} + 2^{80\epsilon  M_{t^{\star}} } \mathbf{1}_{n\geq (1-2  {\alpha}^{\star}-60\epsilon  )M_{t^{\star}}  } \\
 &\qquad \times \min\{2^{(k+2n)/2+ {\alpha}^{\star} M_{t^{\star}} }, 2^{(k+4n)/2+(1+3\iota)M_{t^{\star}} }\} \big]. \\
 \end{split}
\ee

 \medskip

\textbf{Step 2B.}\qquad  If $  \min\{-2l, - l-n\} +4\epsilon M_t  <  m+k\leq -2l+4\epsilon M_t$. 

 \medskip

Note that, for this sub-case, we have $l\leq n$ and  $c(m,k,l)\leq n $.  
 Recall     \eqref{sep5eqn101}  and the definition of kernels in  \eqref{2022feb8eqn11}.  After doing integration by parts in $\xi$ along $v$ direction and directions perpendicular to $v$ many times, we have 
\be
\begin{split}
&\big(\big|  \mathcal{K}_{k,n,l,r}^{ess;\mu,m}(y, v, \omega, \zeta)  \big| + 2^k \big| \mathcal{K}_{k,n,l,r}^{ess;err;\mu,m}(y, v, \omega, \zeta)\big| \big)\varphi_{l;-j}(\omega  +  \tilde{v})   \\
&\lesssim 2^{4k+ n+3c(m,k,l) +\epsilon M_t } \| \mathfrak{m}(\cdot, \zeta)\|_{\mathcal{S}^\infty}  (1+2^k|y\cdot \tilde{v}|)^{-N_0^3 } (1+2^{k+c(m,k,l) }|y\times \tilde{v}|)^{-N_0^3 } .\\
\end{split}
\ee
From the above estimate of kernels, the conservation law  \eqref{conservationlaw},   the volume of support of $\omega, v$, and the estimate  \eqref{nov24eqn41}  if $|  v_{\bot}|\geq 2^{(\alpha_t+\epsilon)M_t}$, we have
\be\label{sep19eqn101}
\begin{split}
&\big|H_{k,j;n,l,r;ess}^{\mu,m;lin}(t, x,\zeta)\big|\\
&\lesssim 2^{m+n+c(m,k,l)+\epsilon M_t }(2^{m+k}+1)  \min\{2^{3k+2 c(m,k,l)} 2^{2l-j},  2^{ 2l} 2^{j+2(\tilde{\alpha}_t+\epsilon)M_t}  \}\| \mathfrak{m}(\cdot, \zeta)\|_{\mathcal{S}^\infty}\\
&\lesssim 2^{n}   \min\{  2^{-m+k-l-j+14\epsilon M_t}, 2^{-2m-j-3l+18\epsilon M_t},  2^{m+l+7\epsilon M_t} 2^{j+2\tilde{\alpha}_t M_t}\}\| \mathfrak{m}(\cdot, \zeta)\|_{\mathcal{S}^\infty}\\
&\lesssim  2^{n} \| \mathfrak{m}(\cdot, \zeta)\|_{\mathcal{S}^\infty}  \min\big\{ \big( 2^{-2m-j-3l+18\epsilon M_t}\big)^{1/4}\big(2^{m+l+7\epsilon M_t} 2^{j+2\tilde{\alpha}_t M_t} \big)^{3/4},\\
&\qquad  \big(  2^{-m+k-l-j+14\epsilon M_t}\big)^{1/2} \big(2^{m+l+7\epsilon M_t} 2^{j+2\tilde{\alpha}_t M_t} \big)^{1/2}\big\}\\
&  \lesssim 2^{n+15\epsilon M_t}\min\{ 2^{3\tilde{\alpha}_t M_t/2 +j/2 + \epsilon M_t}, 2^{k/2 + \tilde{\alpha}_t M_t} \} \| \mathfrak{m}(\cdot, \zeta)\|_{\mathcal{S}^\infty}\\
 & \lesssim \| \mathfrak{m}(\cdot, \zeta)\|_{\mathcal{S}^\infty}\big[ 2^{(1-20\epsilon)M_{t^{\star}}} + 2^{15\epsilon M_{t^{\star}}} 2^{(k+2n)/2+\alpha^{\star} M_{t^{\star}}} \mathbf{1}_{n\geq -(1+3\iota+80\epsilon)M_{t^{\star}}/2 }\big]. \\
 \end{split}
\ee
To sum up, after combining the obtained estimates  \eqref{sep20eqn21},  \eqref{sep19eqn100}, and  \eqref{sep19eqn101}, we have
\be\label{sep20eqn24}
\begin{split}
\|H_{k,j;n,l,r;ess}^{\mu,m,i;lin}(t, x,\zeta)\|_{L^\infty_x} &\lesssim  \| \mathfrak{m}(\cdot, \zeta)\|_{\mathcal{S}^\infty}   \big[  2^{(1-20\epsilon)M_{t^{\star}}}  
+   2^{60\epsilon  M_{t^{\star}} }   \mathbf{1}_{n\geq -(1+3\iota+120\epsilon)M_{t^{\star}}/2 } \\
&\quad \times  \min\{ 2^{(k+2n)/2+\alpha^{\star} M_{t^{\star}}}, 2^{(k+4n)/2+ (7/6+5\iota/2) M_{t^{\star}}} \} \big]. \\
\end{split}
\ee

\medskip

\noindent \textbf{Step 3.}\qquad The estimate of $H_{k,j;n,l,r;p; ess}^{\mu,m,i;non,1}(t, x,\zeta)$, $i=2,3,4.$

\medskip

Recall \eqref{2022feb8eqn1},  \eqref{2022feb8eqn11},   and the definition of   cutoff functions $ \varphi_{j,n}^{i; r}(v, \zeta)$ in  \eqref{sep17eqn63} and   \eqref{sep4eqn6}. As a result of direct computations, after doing integration by parts in $\xi$ along $\zeta$ direction and directions perpendicular to $\zeta$ many times, we have 
 \be\label{sep20eqn1}
 \begin{split}
 &\big| \nabla_v \mathfrak{H}^{ess; \mu,E,i}_{k,j;n,l,r}(y,\omega, v, \zeta)\big| + 2^k|\nabla_v  \mathfrak{H}^{ess;  err;\mu,E,i}_{k,j;n,l,r}(y,\omega, v, \zeta)|\\
 &\lesssim  \| \mathfrak{m}(\cdot, \zeta)\|_{\mathcal{S}^\infty} 2^{3k+2n-j + \epsilon M_t} (2^{-l} +    \sum_{b\in \{ \vartheta^\star_0, \vartheta^\star_1,\vartheta^\star_2,r \}} 2^{- b}\mathbf{1}_{|\tilde{v}- \tilde{\zeta}|\lesssim 2^{b}}\big) \\
 &\qquad \times   (1+2^k|y\cdot \tilde{\zeta}|)^{-N_0^3} (1+2^{k+n }|y\times \tilde{\zeta}|)^{-N_0^3} . \\
 \end{split}
\ee

  From the above estimate of kernels,  the estimate  \eqref{march18eqn31}  in Lemma \ref{conservationlawlemma},  the volume of support of $v, \omega$,  and the estimate  \eqref{nov24eqn41}  if $|  v_{\bot}|\geq 2^{(\alpha_t+\epsilon)M_t}$, we have 
\be\label{sep19eqn103}
\begin{split}
&\big| H_{k,j;n,l,r;p; ess}^{\mu,m;non,1}(t, x,\zeta)\big|\\
&\lesssim \sum_{b\in \{l, \vartheta^\star_0, \vartheta^\star_1,\vartheta^\star_2,r   \}} \| \mathfrak{m}(\cdot, \zeta)\|_{\mathcal{S}^\infty}  (2^m+2^{-k} ) 2^{-j- b +2\epsilon M_t }   \big( \min\{2^{-2m}, 2^{m+3k+2n}\} 2^{ 3j+2b} \big)^{1/2}\\
&\quad \times  \big(\min\{ 2^{m+3k+2n+2l-j},    2^{m+2l+j+2(\tilde{\alpha}_t+\epsilon)M_t} \} \big)^{1/2}\\
&\lesssim   \| \mathfrak{m}(\cdot, \zeta)\|_{\mathcal{S}^\infty}  2^{j/2+2\epsilon M_t}\big[\min\big\{\big(2^{2k+2n-j} \big)^{1/3} \big(2^{-k+j+2(\tilde{\alpha}_t+\epsilon)M_t}  \big)^{2/3}, 2^{2k+2n-j}\big\}\big]^{1/2}\\
&  \lesssim 2^{6\epsilon M_t}\min\{  2^{2\tilde{\alpha}_t M_t/3+  2j/3+n/3 }, 2^{k+n}\}\| \mathfrak{m}(\cdot, \zeta)\|_{\mathcal{S}^\infty}\\
 &\lesssim  \| \mathfrak{m}(\cdot, \zeta)\|_{\mathcal{S}^\infty}\big[2^{(1-20\epsilon)M_{t^{\star}}} + 2^{127\epsilon M_{t^{\star}}} \mathbf{1}_{n \geq  (1-2\alpha^{\star}-90\epsilon)M_{t^\star}  } \\
 &\quad \times \min\{ 2^{(k+2n)/2+\alpha^{\star}M_{t^{\star}} },2^{(k+4n)/2+ (1+3\iota) M_{t^{\star}} } \}  \big]. \\
 \end{split}
\ee

\medskip

\noindent \textbf{Step 4.}\qquad  The estimate of $H_{k,j;n,l,r}^{\mu,m,i;non,2}(t, x,\zeta)$, $i=2 .$

\medskip

Due to the  cutoff functions $\varphi^i_{j,n}(v, \zeta)$ (see  \eqref{sep4eqn6}), we  have $m+k+l\leq -n+2\epsilon M_t$ and  $\max\{l, c(m,k,l)\}\geq n +3\epsilon M_t/4$.   Note that, the volume of support of $\omega$ is bounded by $2^{2 \max\{l,c(m,k,l)\} }$.  From the Cauchy-Schwarz inequality,  the estimate of kernels in  \eqref{sep20eqn1}, the volume of support of $v, \omega$,  we have
\be\label{sep20eqn25}
\begin{split}
&|  H_{k,j;n,l,r;p; ess}^{\mu,m;non,2}(t, x,\zeta)| \\
& \lesssim  \sum_{b\in \{l, \vartheta^\star_0, \vartheta^\star_1,r \}}   \| \mathfrak{m}(\cdot, \zeta)\|_{\mathcal{S}^\infty}  2^{ -j- b+ l +4\epsilon M_t} (2^{m}+2^{-k})\\
&\qquad \times \big(2^{m+3k+2n} 2^{2 l }  ( 2^{3j+2b}\mathbf{1}_{b\neq l} + 2^{3j+2c(m,k,l)}\mathbf{1}_{b = l}  )\big)^{1/2}\\
&\qquad\times  \big(\min\{2^{m+3k+2n+2  l  -j}, 2^{m+2l+3j +2 \max\{l,c(m,k,l)\}} \} \big)^{1/2}\\
&\lesssim 2^{  16\epsilon M_t}  \min\{2^{k+n},   2^{-k/2+2j+n} \}\| \mathfrak{m}(\cdot, \zeta)\|_{\mathcal{S}^\infty}\\
&\lesssim 2^{  16\epsilon M_t} \min\big\{(2^{k+n})^{1/3}( 2^{-k/2+2j+n} )^{2/3}, (2^{k+n})^{2/3}( 2^{-k/2+2j+n} )^{1/3}\big\}\| \mathfrak{m}(\cdot, \zeta)\|_{\mathcal{S}^\infty}\\
 &\lesssim 2^{  16\epsilon M_t}  \min\{ 2^{4j/3+n}, 2^{k/2+n+2j/3}\}\| \mathfrak{m}(\cdot, \zeta)\|_{\mathcal{S}^\infty}\\
&\lesssim  \| \mathfrak{m}(\cdot, \zeta)\|_{\mathcal{S}^\infty}\big[ 2^{(1-20\epsilon)M_{t^{\star}}}  +  2^{60\epsilon M_{t^{\star}}}   \mathbf{1}_{n\geq   -M_{t^{\star}}/3-40\epsilon M_{t^{\star}}}\\
&\qquad \times  \min\{ 2^{(k+2n)/2+2M_{t^{\star}}/3 },2^{(k+4n)/2+  M_{t^{\star}} } \} \big].\\
\end{split}
\ee

\medskip

\noindent \textbf{Step 5.}\qquad The estimate of $H_{k,j;n,l,r}^{\mu,m,i;non,2}(t, x,\zeta)$, $i=3,4.$

\medskip

 Due to the  cutoff functions $\varphi^i_{j,n}(v, \zeta)$ for the cases $i \in\{3,4\}$ (see  \eqref{sep4eqn6}), we have $|\tilde{v}- \tilde{\zeta}|\lesssim 2^{n+\epsilon M_t}. $   As a result of direct computations, after doing integration by parts in $\xi$ along $\zeta$ direction and directions perpendicular to $\zeta$ many times, we have 
\be\label{2022feb9eqn11}
\begin{split}
&\big| \nabla_v \mathfrak{H}^{ess; \mu,E,i}_{k,j;n,l,r}(y,\omega, v, \zeta)\big| + 2^k|\nabla_v  \mathfrak{H}^{ess;  err;\mu,E,i}_{k,j;n,l,r}(y,\omega, v, \zeta)|\\
& \lesssim  \| \mathfrak{m}(\cdot, \zeta)\|_{\mathcal{S}^\infty} 2^{3k+2n-j + \epsilon M_t}  2^{-\min\{l,n\}}  \varphi_{[l-1,l+1];-j}(\omega + \tilde{v}) \\
 &\qquad \times  \psi_{< n+\epsilon M_t+2}( \tilde{\zeta} - \tilde{v} )\varphi_j(v)  (1+2^k|y\cdot \tilde{\zeta}|)^{-N_0^3}  (1+2^{k+n }|y\times \tilde{\zeta}|)^{-N_0^3}  . \\
 \end{split}
\ee

  Based on the possible size of $m+k$, we   proceed in two sub-steps as follows.

\medskip
\textbf{Step 5A.}\qquad If $m+k\leq \min\{-2l, - l-n\} +4\epsilon M_t$. 
\medskip

 From the estimate of kernels in  \eqref{2022feb9eqn11}, the volume of support of $v, \omega$,  we have
\be\label{2022feb9eqn12}
\begin{split}
&|  H_{k,j;n,l,r;p; ess}^{\mu,m,i;non,2}(t, x,\zeta)|\\
& \lesssim   \| \mathfrak{m}(\cdot, \zeta)\|_{\mathcal{S}^\infty}  2^{\max\{m,-k\}}  2^{ -j- \min\{l,n\}+ l +4\epsilon M_t}\big(2^{m+3k+2n + 2l+3j+2n}  \big)^{1/2} \\
&\qquad \times \big(\min\{2^{m+3k+2n+2  l  -j}, 2^{m+2l+3j +2 n \}} \} \big)^{1/2}  \\
&\lesssim 2^{  12\epsilon M_t}  \min\{2^{k+n},   2^{-k/2+2j+n} \}\| \mathfrak{m}(\cdot, \zeta)\|_{\mathcal{S}^\infty}\\
&\lesssim 2^{  12\epsilon M_t}\min\big\{(2^{k+n})^{1/3}(2^{-k/2+2j+n})^{2/3},(2^{k+n})^{2/3}(2^{-k/2+2j+n})^{1/3}   \big\}  \| \mathfrak{m}(\cdot, \zeta)\|_{\mathcal{S}^\infty} \\
&\lesssim 2^{  12\epsilon M_t}  \min\{ 2^{4j/3+n}, 2^{k/2+n+2j/3}\}\| \mathfrak{m}(\cdot, \zeta)\|_{\mathcal{S}^\infty} \\
&\lesssim  \| \mathfrak{m}(\cdot, \zeta)\|_{\mathcal{S}^\infty}\big[ 2^{(1-20\epsilon)M_{t^{\star}}}  +  2^{50\epsilon M_{t^{\star}}}   \mathbf{1}_{n\geq   -M_{t^{\star}}/3-60\epsilon M_{t^{\star}}} \\
&\quad \times \min\{ 2^{(k+2n)/2+2M_{t^{\star}}/3 },2^{(k+4n)/2+  M_{t^{\star}} } \}  \big].
\end{split}
\ee

\medskip
\textbf{Step 5B.}\qquad If $  \min\{-2l, - l-n\} +4\epsilon M_t< m + k \leq -2l + 4\epsilon M_t$. 
\medskip

For this sub-case, we have $l\leq n$ and $c(m,k,l)\leq n -\epsilon M_t. $ Recall   \eqref{2022feb8eqn11}.  After doing integration by parts in $\xi$ along $v$ direction and directions perpendicular to $v$ many times, for $i\in\{3,4\},$ we have 
 \be\label{2022feb6eqn1}
\begin{split}
&| \nabla_v  \mathcal{H}^{ess;\mu, E, i }_{k,j;n,l,r}( \mathfrak{m})(y,\omega, v, \zeta)| + 2^k| \nabla_v  \mathcal{H}^{ess;err;\mu, E,i }_{k,j;n,l,r}( \mathfrak{m})(y,\omega, v, \zeta)|\\
& \lesssim   2^{3k+2 c(m,k,l) -j-l } \| \mathfrak{m}(\cdot, \zeta)\|_{\mathcal{S}^\infty}  \varphi_{[l-1,l+1];-j}(\omega + \tilde{v}) \psi_{< n+\epsilon M_t+2}( \tilde{\zeta} - \tilde{v} )\varphi_j(v)   \\
&\quad \times (1+2^{k-4\epsilon M_t}|y\cdot \tilde{v}|)^{-N_0^3 }(1+2^{k+ c(m,k,l) -4\epsilon M_t}|y\times \tilde{v}|)^{- N_0^3 } .\\
\end{split}
\ee

Moreover, an important observation here is that  $\p_{v_3}(symbol)$ is better than $\nabla_{  v_{\bot} }(symbol)$.    More precisely,  the following improved estimate holds  for the third component of the kernels, 
 \be\label{2022feb9eqn16}
\begin{split}
&| \p_{v_3} \mathcal{H}^{ess;\mu, E,i }_{k,j; n,l,r}( \mathfrak{m})(y,\omega, v, \zeta)| + 2^k |\p_{v_3}  \mathcal{H}^{ess;err;\mu, E, i }_{k,j;n,l,r}( \mathfrak{m})(y,\omega, v, \zeta)|\\
& \lesssim   2^{3k+2 c(m,k,l)   +2\epsilon M_t-j } (2^{-l +\max\{l, p\}} +2^{-n+\max\{n,p\}})\\
&\quad \times  \| \mathfrak{m}(\cdot, \zeta)\|_{\mathcal{S}^\infty}\varphi_{[l-1,l+1];-j}(\omega + \tilde{v})  \psi_{< n+\epsilon M_t+2}( \tilde{\zeta} - \tilde{v} )\varphi_j(v)   \\
&\quad \times (1+2^{k-4\epsilon M_t}|y\cdot \tilde{v}|)^{-N_0^3 }(1+2^{k+ c(m,k,l) -4\epsilon M_t}|y\times \tilde{v}|)^{- N_0^3 }    .\\
\end{split}
\ee
Based on the possible size of $p$,   we  proceed in two sub-steps as follows.

\medskip
\quad \textbf{Step 5B-a.} \qquad If $p\leq l+\epsilon M_t$. \qquad 

\medskip

Recall that  the associated coefficient of $\mathbf{P}(EB^2(t,s,x ,\omega, v))$ is better than  the associated coefficient of $\mathbf{P}_3(EB^2(t,s,x ,\omega, v))$ (see \eqref{july1eqn13}). 
From the Cauchy-Schwarz inequality,  the     estimate of kernels in  \eqref{2022feb6eqn1} and  \eqref{2022feb9eqn16}, the volume of support of $v, \omega$,  we have
\be\label{sep20eqn28}
\begin{split}
 &|  H_{k,j;n,l,r;ess}^{\mu,m;p;non}(t, x,\zeta)|\\ &\lesssim \| \mathfrak{m}(\cdot, \zeta)\|_{\mathcal{S}^\infty}  2^{m-j +\max\{l,p\}  +2\epsilon M_t } \big(2^{m+3k+2 c(m,k,l) } 2^{2n} 2^{3j+2l} \big)^{1/2}  \\
 &\qquad \times  \big(\min\{2^{m+3k+2 c(m,k,l)+2l-j}, 2^{m+2l+3j+2n+12\epsilon M_t} \} \big)^{1/2}\\
&   \lesssim 2^{l+16\epsilon M_t}  \min\{2^{k+n},   2^{m/2+2j+ 2n },   2^{-(k+2l)/2+2j+ 2n } \} \| \mathfrak{m}(\cdot, \zeta)\|_{\mathcal{S}^\infty}\\
&\lesssim \| \mathfrak{m}(\cdot, \zeta)\|_{\mathcal{S}^\infty}\big[2^{(1-20\epsilon)M_{t^{\star}} } +  2^{(k+2n)/2+2M_{t^{\star}}/3+20\epsilon M_{t^{\star}} }  \mathbf{1}_{n\geq - M_{t^{\star}}/3-40\epsilon M_{t^{\star}}}\big].\\
\end{split}
\ee

\medskip
\quad \textbf{Step 5B-b.} \qquad If $p\geq  l+\epsilon M_t$.

\medskip

\qquad Due to the cutoff function $\varphi_{l;-j}(\tilde{v}+\tilde{\omega})\varphi_{p;-10M_t}( \omega_{\bot})$, for this sub-case, we have $| v_{\bot}|/|v|\sim 2^p$. Moreover, from the estimate  \eqref{nov24eqn41}, we can rule out the case $| v_{\bot}|\geq 2^{(\alpha_t+\epsilon) M_t}$. Hence, it suffices to consider the case $p+j\leq (\alpha_t+\epsilon)M_t$. From the Cauchy-Schwarz inequality,  the     estimate of kernels in  \eqref{2022feb6eqn1}  and  \eqref{2022feb9eqn16}, the volume of support of $v, \omega$,  we have
\be\label{2022feb6eqn3}
\begin{split}
|  H_{k,j;n,l,r;ess}^{\mu,m,i;p;non}(t, x,\zeta)|&\lesssim  \| \mathfrak{m}(\cdot, \zeta)\|_{\mathcal{S}^\infty} 2^{m-j +\max\{l,p\} }\big(2^{m+3k+2 c(m,k,l) } 2^{2n} 2^{3j+2l} \big)^{1/2}  \\
&\quad \times    \big(\min\{2^{m+3k+2 c(m,k,l)+2l-j}, 2^{m+2l+3j+2n+12\epsilon M_t} \} \big)^{1/2} \\
& \lesssim \| \mathfrak{m}(\cdot, \zeta)\|_{\mathcal{S}^\infty}  2^{p +16\epsilon M_t}  \min\{2^{k+n },   2^{m/2+2j+ 2n} \}\\
& \lesssim  \| \mathfrak{m}(\cdot, \zeta)\|_{\mathcal{S}^\infty}\big[ 2^{(1-20\epsilon) M_{t^{\star}} } + 2^{30\epsilon M_{t^{\star}}}\mathbf{1}_{n\geq -(\alpha^{\star}+40\epsilon)M_{t^{\star}}/2 }\\
&\quad \times \min\{2^{(k+2n)/2+ \alpha^{\star} M_{t^{\star}}}, 2^{(k+4n)/2+(1+2\iota)M_{t^{\star}}} \}  \big]. \\
\end{split}
\ee

Recall the decomposition in  \eqref{sep19eqn81}.  To sum up,  our desired estimate  \eqref{sep19eqn71}    holds from the obtained estimates  \eqref{2022feb12eqn15},  \eqref{sep19eqn103}, and \eqref{sep20eqn25}. The desired estimate \eqref{2022feb12eqn11}  holds from the obtained estimates \eqref{sep20eqn24},  \eqref{sep19eqn103},   \eqref{2022feb9eqn12}, \eqref{sep20eqn28},  and    \eqref{2022feb6eqn3}.
\end{proof}

Next, we turn our attention to the relatively high frequency case, specifically when $m+k\geq -2l+4\epsilon M_t$. In this scenario, we will utilize the second decomposition from \eqref{sep18eqn50} in Lemma \ref{locdeclemm}.

 In the following Lemma, we first estimate the ``$T$-part''. 
\begin{lemma}\label{secondTpartlarge}
For any $i\in \{2,3,4\}, ( l,r)\in \mathcal{B}_i $, see  \eqref{sep18eqn50}.    Under the assumption of Theorem \ref{mainresultsfirstpart}, the following estimate holds if  $m+k\geq -2 l +4\epsilon M_t$, 
\be\label{sep6eqn49}
\begin{split}
&\big\|   \widetilde{T}_{k,j;n,l,r }^{T; \mu,m, i }(\mathfrak{m}, E)(t,x,  \zeta) + \hat{\zeta}\times   \widetilde{T}_{k,j;n,l,r}^{T; \mu ,m, i}(\mathfrak{m},B)(t,x, \zeta )\big\|_{L^\infty_x}\\
&\lesssim \| \mathfrak{m}(\cdot, \zeta)\|_{\mathcal{S}^\infty}\big[2^{(1-20\epsilon)M_{t^{\star}} } + 2^{50\epsilon  M_{t^{\star}} } \mathbf{1}_{n\geq -(3/7+30\epsilon)M_{t^{\star}}  } \\
&\quad \times \min\{2^{(k+2n)/2+  \alpha^{\star} M_{t^{\star}} } , 2^{(k+4n)/2+  (1+3\iota) M_{t^{\star}} }\} \big]. 
\end{split}
\ee
\end{lemma}

\begin{proof}
Recall   \eqref{sep18eqn44}. Let $\omega=(\sin \theta\cos\phi, \sin \theta \sin \phi, \cos\theta)$. After doing dyadic decomposition for  the sizes of $\sin \theta$, $\sin \phi$, localizing the frequency around the fixed direction of $\tilde{\zeta}$, and using the volume of support of $v, \omega$ and estimate  \eqref{july10eqn89}  in Proposition \ref{Linfielec} for the case $p=-10M_t$ or $q=-10M_t$, we have 
\be\label{sep7eqn51}
\begin{split}
&\big|   \widetilde{T}_{k,j;n,l,r}^{T; \mu ,m, i}(\mathfrak{m}, E)(t,x,  \zeta ) + \hat{\zeta}\times   \widetilde{T}_{k,j;n,l,r}^{T; \mu,m, i }(\mathfrak{m},B)(t,x, \zeta )\big| \\
& \lesssim \sum_{\begin{subarray}{c}
\star\in \{ess,err\}  \\
p, q\in ( -10M_t, 2]\cap \Z\\
\end{subarray}} \big| G_{k,j; n,l,r }^{\star; m, i;p,q}(t,x,\zeta)\big|  +\| \mathfrak{m}(\cdot, \zeta)\|_{\mathcal{S}^\infty},
\end{split}
\ee
where 
\be\label{sep9eqn21}
\begin{split}
 G_{k,j; n,l,r }^{\star; m, i;p,q}(t,x,\zeta) &:=  \int_0^t  \int_{\R^3}  \int_{\R^3}  \int_{ 0}^{2\pi}    \int_0^\pi \big[ i\mu \big( \omega^{m;E}_{j,l,r}(t-s,v,\omega) +\hat{\zeta}\times \omega^{m;B}_{j,l,r}(t-s,v,\omega)\big)   K_{k;n}^{   }(\mathfrak{m})(y, \zeta   )\\
 &\quad   + \big(c^{a;m,E}_{j,l,r}(t-s,v,\omega) +\hat{\zeta}\times c^{a;m,B}_{j,l,r}(t-s,v,\omega) \big)   K_{k;n}^{  a }(\mathfrak{m})(y, \zeta    ) + (t-s)^{-1} \\
 & \quad \times \big(  c^{err;m,E}_{j,l,r}(t-s,v,\omega) + \hat{\zeta}\times  c^{err;m,B}_{j,l,r}(t-s,v,\omega)\big)    \tilde{K}_{k;n}^{  \mu}(\mathfrak{m})(y,    \zeta   ) \big]\\
 & \quad \times  f(s,x-y+(t-s)\omega, v)   \varphi_{j,n}^{i; r}(v, \zeta) \varphi_{m;-10M_t}(t-s)\\
 & \quad \times    \varphi_{p;-10M_t}(\sin \theta)  \varphi_{q;-10M_t}(\sin \phi) \psi_{\star}(\omega, \zeta) \sin \theta d\theta d\phi    d y d v d s. \\
 \end{split}
\ee

The above kernels $K_{k;n}^{   }(\mathfrak{m})(y,   \zeta ),   K_{k;n}^{  a } (\mathfrak{m})(y,   \zeta  ),  $ and $  \tilde{K}_{k;n}^{ \mu}(\mathfrak{m})(y,    \zeta ) $ are defined  in  \eqref{sep5eqn48}  and the cutoff functions $\psi_{\star}(\omega, \zeta),$ $ \star\in \{ess,err\}$, are defined as follows, 
\be\label{feb9eqn71}
\begin{split}
\psi_{ess}(  \omega,\zeta)& :=\psi_{\leq \max\{n, -(m+k)/2\} +\epsilon M_t/2 }(\omega \times \tilde{\zeta}),\\
    \psi_{err}(   \omega,\xi)&:=\psi_{ >  \max\{n, -(m+k)/2\} +\epsilon M_t/2  } (\omega \times \tilde{\zeta}).\\
\end{split}
\ee

  Recall  \eqref{sep5eqn48}. After doing integration by parts in $\xi$ in the direction of  $  {\zeta}$   and  directions perpendicular to $\zeta$ or doing integration by parts in $\xi$ in $e_1,e_2,e_3$ directions, we have 
\be\label{sep6eqn31}
\begin{split}
&\big|   K_{k;n}^{  \mu}(\mathfrak{m})(y,  \zeta, \omega) \big|+ \big|   K_{k;n}^{   }(\mathfrak{m})(y,     \zeta) \big|+ \big|   K_{k;n}^{  a}(\mathfrak{m})(y , \zeta) \big| + 2^k \big| \widetilde{K}_{k;n}^{  \mu}(\mathfrak{m})(y,  \zeta)\big|\\
&\lesssim 2^{3k+2n}\| \mathfrak{m}(\cdot, \zeta)\|_{\mathcal{S}^\infty} \min\big\{ (1+2^k|y\cdot\tilde{\zeta}|)^{- N_0^{ 3}} (1+2^{k+n}|y\times \tilde{\zeta}| )^{- N_0^{ 3} }, \\
&\qquad  (1+2^{k+n}|  y_{\bot}| )^{- N_0^{ 3} } (1+2^{k+n } (2^n + |\tilde{\zeta}\times e_3|)^{-1} |y_3|)^{-N_0^3} \big\}.\\
\end{split}
\ee

We first rule out the error type terms. Due to the cutoff function $ \psi_{err}(   \omega,\xi)$, see  \eqref{feb9eqn71},  after representing $ G_{k,j;l,n}^{err; m, i;p,q}(t,x,\zeta)$ on the Fourier side and  doing integration by parts in $\omega$ once, c.f.,  \eqref{march4eqn41},  we gain at least $2^{-\epsilon M_t/40}$ for the error type term. Hence, after doing integration by parts in $\omega$ many times, we have
\be\label{2022feb9eqn22}
\big| G_{k,j;  n,l,r }^{err; m, i;p,q}(t,x,\zeta)\big|\lesssim  \| \mathfrak{m}(\cdot, \zeta)\|_{\mathcal{S}^\infty}.
\ee

Now, we focus on the estimate of the essential part $ G_{k,j;l,n}^{ess; m, i;p,q}(t,x,\zeta)$. Based on the possible size of $i$, we   proceed in two  main steps  as follows. 

\medskip
\noindent \textbf{Step 1.} \qquad   If $i=2.$
\medskip

Due to the cutoff functions $\varphi_{j,n}^i(v, \zeta)$, see  \eqref{sep4eqn6}, we have $l\geq n+3\epsilon M_t/4$ and $|\tilde{v}-\tilde{\zeta}|\sim 2^l$. Since $u\times(u \times v)= uu\cdot v - v u\cdot u $,   we have 
\be\label{2022feb9eqn71}
\begin{split}
 ( \delta_{1r}, \delta_{2r}, \delta_{3r}) + \omega_r \hat{v}&=    - \hat{v}\times (\hat{v} \times ( \delta_{1r}, \delta_{2r}, \delta_{3r}))\\
 &\qquad  + ( \delta_{1r}, \delta_{2r}, \delta_{3r})\cdot( \hat{v}+\omega) \hat{v} + (1-\hat{v}\cdot \hat{v}) ( \delta_{1r}, \delta_{2r}, \delta_{3r})), \\ 
\hat{v} + \omega  &= -\hat{v}\hat{v}\cdot \omega + \omega \hat{v}\cdot\hat{v} +\hat{v}(1+\hat{v}\cdot \omega)+\omega(1-\hat{v}\cdot \hat{v})\\
& = - \hat{v}\times(\hat{v}\times \omega)  +\hat{v}(1+\hat{v}\cdot \omega)+\omega(1-\hat{v}\cdot \hat{v}).
\end{split}
\ee
 From the above equalities,  the detailed formulas of coefficients $\omega^{m;U}_{j,l}(t-s,v,\omega),$ $c^{q;m,U}_{j,l}(t-s,v,\omega),$ $  c^{err;m,U}_{j,l}(t-s,v,\omega),  U\in \{E, B\} $, in  \eqref{2022feb16eqn4}, \eqref{sep20eqn42}, and  \eqref{sep20eqn44},  we have 
\be\label{2022feb6eqn11}
\begin{split}
&\big|c^{ a;m,E}_{j,l,r}(t-s,v,\omega)  
+\hat{\zeta}\times c^{a;m,B}_{j,l,r}(t-s,v,\omega)\big| \\
&+\big|c^{err;m,E}_{j,l,r}(t-s,v,\omega)  
+\hat{\zeta}\times c^{err;m,B}_{j,l,r}(t-s,v,\omega)\big|\\
&+\big|\omega^{m;E}_{j,l,r}(t-s,v,\omega) +\hat{\zeta}\times \omega^{m;B}_{j,l,r}(t-s,v,\omega)\big|\\
& \lesssim 2^{ \epsilon M_t} \varphi_{[l-1,l+1];-j}(v/|v|+\omega)\varphi_j(v)\varphi_{[m-1,m+1];-10M_t }(t-s). \\
\end{split}
\ee

Based on the possible size of $m+k$, we    proceed in two  sub-steps   as follows. 

\medskip

 \textbf{Step 1A.}\qquad If $m+k\leq -2n+\epsilon M_t$.

 \medskip 

From  the   estimate of coefficients in \eqref{2022feb6eqn11}, the  estimate of kernels in  \eqref{sep6eqn31},   the volume of support of $(v, \omega)$, and the estimate  \eqref{nov24eqn41} if $| v_{\bot} |\geq 2^{(\alpha_t + \epsilon)M_t}$,   we have
\be\label{aug8eqn23}
\begin{split}
&\big| G_{k,j;   n,l,r }^{ess; m, i;p,q}(t,x,\zeta)\big|\\
&\lesssim \| \mathfrak{m}(\cdot,\zeta)\|_{\mathcal{S}^\infty} 2^{ m + \epsilon M_t}    (2^{2n}+2^{-m-k+\epsilon M_t})  \min\{ 2^{3k+2n-j},   2^{j+2(\tilde{\alpha}_t+\epsilon) M_t}\}\\
&\lesssim 2^{2\epsilon M_t}  \| \mathfrak{m}(\cdot,\zeta)\|_{\mathcal{S}^\infty}   \min\big\{(2^{2k+2n-j})^{1/3}(2^{-k+j+2(\tilde{\alpha}_t+\epsilon) M_t})^{2/3},\\
&\qquad  (2^{2k+2n-j})^{1/2}(2^{-k+j+2(\tilde{\alpha}_t+\epsilon) M_t})^{1/2} \big\}\\
&\lesssim  \| \mathfrak{m}(\cdot,\zeta)\|_{\mathcal{S}^\infty} 2^{4\epsilon M_t} \min\{2^{2n/3+ (4\tilde{\alpha}_t M_t+j)/3  }, 2^{(k+2n)/2 + \tilde{\alpha}_t M_t}\}\\
 & \lesssim  \| \mathfrak{m}(\cdot,\zeta)\|_{\mathcal{S}^\infty}\big[ 2^{(1-20\epsilon)M_{t^{\star}}}+   2^{(k+2n)/2+\alpha^{\star}M_{t^\star} + 4\epsilon M_{t^\star}}  \mathbf{1}_{n\geq (1-2\alpha^{\star}-40\epsilon) M_{t^{\star}}  }\big].
 \end{split}
\ee
 
\medskip

 \textbf{Step 1B.}\qquad If $m+k\geq -2n+\epsilon M_t$. 

 \medskip

From  the   estimate of coefficients in  \eqref{2022feb6eqn11}, the  estimate of kernels in  \eqref{sep6eqn31},   the estimate \eqref{march18eqn31}  in Lemma \ref{conservationlawlemma}, the volume of support of $(v, \omega)$, and the estimate  \eqref{nov24eqn41}  if $| v_{\bot} |\geq 2^{(\alpha_t + \epsilon)M_t}$,   we have
\be\label{2022feb9eqn70}
\begin{split}
&\big| G_{k,j; n,l,r}^{ess; m, i;p,q}(t,x,\zeta)\big|\\
&\lesssim \| \mathfrak{m}(\cdot,\zeta)\|_{\mathcal{S}^\infty} 2^{  \epsilon M_t}    \min\{   2^{m +3j+2l + 2n}, 2^{m  + 2n+j+2(\tilde{\alpha}_t+\epsilon) M_t}, 2^{-2m-j-2l}\}\\
& \lesssim  \| \mathfrak{m}(\cdot, \zeta)\|_{\mathcal{S}^\infty}2^{\epsilon M_t} \min\big\{\big(2^{m  + 2n+j+2(\tilde{\alpha}_t+\epsilon) M_t} \big)^{1/2} \big( 2^{-m-j+k}\big)^{1/2} , \\
&\qquad  \big( 2^{m +3j+2l + 2n}\big)^{1/6}\big(2^{m  + 2n+j+2(\tilde{\alpha}_t+\epsilon) M_t} \big)^{1/2} \big( 2^{-2m-j-2l}\big)^{1/3}\big\}\\
& \lesssim  \| \mathfrak{m}(\cdot, \zeta)\|_{\mathcal{S}^\infty}2^{4\epsilon M_t} \min\{2^{2j/3+\tilde{\alpha}_t M_t + n  } ,   2^{k/2 + n  }  2^{\tilde{\alpha}_t M_t }  \} \\
& \lesssim \| \mathfrak{m}(\cdot,\zeta)\|_{\mathcal{S}^\infty}\big[ 2^{(1-20\epsilon)M_{t^{\star}}}+   2^{(k+2n)/2+\alpha^{\star}M_{t^\star} + 4\epsilon M_{t^\star}}  \mathbf{1}_{n\geq (1-2\alpha^{\star} -40\epsilon ) M_{t^{\star}}  }\big].\\
 \end{split}
\ee
We remark that, in the above estimate, we used the fact that $-m-2l\leq  k .$

\medskip
\noindent \textbf{Step 2.} \qquad   If $i=3,4.$
\medskip

Due to the cutoff functions $\varphi_{j,n}^i(v, \zeta), i\in\{3,4\}$, see  \eqref{sep4eqn6}, the support of the $v, \omega, \xi$,    implies that $|\tilde{v}\times \omega|\lesssim  2^{n+\epsilon M_t} +2^{l-\epsilon M_t}, |\tilde{v}\times\omega |\sim 2^l$.  If $l> -j$, then it follows that $l\leq n +\epsilon M_t$, since in this case,   $ |\omega \times\tilde{\xi}|\sim 2^l$. Therefore, it's always true that  $l\in [-j, \max\{-j, n+\epsilon M_t\}]\cap \Z. $

From the equalities in  \eqref{2022feb9eqn71}, the following estimate holds for the case we are considering, 
\be\label{2022feb9eqn61}
\begin{split}
&\big|c^{ a;m,E}_{j,l,r}(t-s,v,\omega)  
+\hat{\zeta}\times c^{a;m,B}_{j,l,r}(t-s,v,\omega)\big| \\
& +\big|c^{err;m,E}_{j,l,r}(t-s,v,\omega)  
+\hat{\zeta}\times c^{err;m,B}_{j,l,r}(t-s,v,\omega)\big|\\
&+\big|\omega^{m;E}_{j,l,r}(t-s,v,\omega) +\hat{\zeta}\times \omega^{m;B}_{j,l,r}(t-s,v,\omega)\big| \\
& \lesssim 2^{ -l + n + \epsilon M_t} \varphi_{[l-1,l+1];-j}(v/|v|+\omega)\varphi_j(v)\varphi_{[m-1,m+1];-10M_t }(t-s). 
\end{split}
\ee

From  the above estimate of coefficients, the  estimate of kernels in  \eqref{sep6eqn31},   the volume of support of $(v, \omega)$, and the estimate \eqref{nov24eqn41}  if $| v_{\bot} |\geq 2^{(\alpha_t + \epsilon)M_t}$,   we have
\be\label{2022feb9eqn62}
\begin{split}
\big| G_{k,j;  n,l,r }^{ess; m, i;p,q}(t,x,\zeta)\big| & \lesssim \| \mathfrak{m}(\cdot,\zeta)\|_{\mathcal{S}^\infty} 2^{ - l+n+\epsilon M_t}    2^{2l} \\
&\qquad \times  \min\{2^{m+ 3j+2p+q},2^{m+ 3j+2\max\{l,n+\epsilon M_t\}},  2^{ m+ j+2(\tilde{\alpha}_t+\epsilon) M_t}\} .\\
\end{split}
\ee
 Moreover, by using  the estimate  \eqref{march18eqn31}  in Lemma \ref{conservationlawlemma} for any fixed $y$, and the estimate of kernels in  \eqref{sep6eqn31}, we have
\be\label{2022feb9eqn63}
\big| G_{k,j;  n,l,r}^{ess; m, i;p,q}(t,x,\zeta)\big| \lesssim  \| \mathfrak{m}(\cdot,\zeta)\|_{\mathcal{S}^\infty} 2^{-2m- j-3l+n+\epsilon M_t}. 
\ee

Based on the possible size  of $p$, we      proceed in two  sub-steps  as follows.

\medskip

 \textbf{Step 2A.}\qquad  If $p\leq \max\{n  + \epsilon M_t, -j\}+\epsilon M_t$. 

\medskip

 After combining   the estimates  \eqref{2022feb9eqn62}   and    \eqref{2022feb9eqn63},  we have 
\be\label{sep3eqn6}
\begin{split}
\big| G_{k,j; n,l,r}^{ess; m, i;p,q}(t,x,\zeta)\big|&\lesssim 2^{3\epsilon M_t} \| \mathfrak{m}(\cdot, \zeta)\|_{\mathcal{S}^\infty}\min\big\{(  2^{m+ 3j+l+2p+ n})^{2/3}(2^{-2m- j-3 l+n })^{1/3},  \\
&\qquad\times  ( \min\{   2^{m+n+l} 2^{j+2\tilde{\alpha}_t M_t} ,  2^{m+n+l} 2^{3j+2p} \}  )^{1/2}  (2^{-2m- j-3 l+n })^{1/2 }\big\}   \\
& \lesssim  \| \mathfrak{m}(\cdot, \zeta)\|_{\mathcal{S}^\infty}2^{7\epsilon M_t}\min\{2^{5j/3+n+4p/3-l/3} ,   2^{k/2 + n  } \min\{2^{\tilde{\alpha}_t M_t }, 2^{j+p}\}  \}\\
  &\lesssim  \| \mathfrak{m}(\cdot, \zeta)\|_{\mathcal{S}^\infty}\big[2^{(1-20\epsilon)M_{t^{\star}} } + 2^{10\epsilon  M_{t^{\star}} } \mathbf{1}_{n\geq -3M_{t^{\star}}/7-30\epsilon M_{t^{\star}} } \\
  &\qquad \times \min\{2^{(k+2n)/2+  \alpha^{\star} M_{t^{\star}} } , 2^{(k+4n)/2+   M_{t^{\star}} }\} \big].\\
\end{split}
\ee
We remark that, in the above estimate, we used the fact that $-m-2l\leq  k .$

\medskip

 \textbf{Step 2B.}\qquad If $p\geq  \max\{n  + \epsilon M_t, -j\}+\epsilon M_t$. 
\medskip

Note that, for this case, we have $|  v_{\bot}|\sim 2^{p+j}.$  Moreover, we have  $m+p\geq m+l+\epsilon M_t\geq -k-l+2\epsilon M_t$. From the estimate in  \eqref{nov24eqn41}, we can rule out the case $|  v_{\bot}|\geq 2^{(\alpha_t+\epsilon)M_t}$, i.e.,  $p+j\geq \alpha_t M_t+\epsilon M_t$. Now, it would be
sufficient to consider the case $p+j\leq  ({\alpha}_t  + \epsilon) M_t. $

Based on the possible size of $|  x_{\bot}|$, we   proceed in two  sub-steps as follows. 

\medskip

\quad  \textbf{Step 2B-a.}\quad If $ |  x_{\bot}|\geq 2^{m+p-\epsilon M_t}$. 

 \medskip

Note that, for the sub-case we are considering,  we have   
\[
\forall y\in B(0, 2^{-k-l+\epsilon M_t/2}), \quad |  x_{\bot}- y_{\bot}| \sim | x_{\bot} |.
\]

From the estimate of the Jacobian of changing coordinates $(\theta, \phi)\longrightarrow (z,w)$
 in \eqref{march18eqn66}, the estimate of kernels in  \eqref{sep6eqn31}, we have
 \[
 \begin{split}
\big| G_{k,j; n,l,r}^{ess; m, i;p,q}(t,x,\zeta)\big|& \lesssim  \| \mathfrak{m}(\cdot, \zeta  )\|_{\mathcal{S}^\infty} 2^{-l+n+\epsilon M_t}2^{-m-p-q} |  x_{\bot}|^{-1} 2^{-j}\\
 & \lesssim  \| \mathfrak{m}(\cdot, \zeta)\|_{\mathcal{S}^\infty} 2^{-2m-2p-q-j-l+n+2\epsilon M_t}.\\
 \end{split}
 \] 
 After combing the above estimate and the estimate 
  \eqref{2022feb9eqn62}, we have
\be\label{sep3eqn10}
\begin{split}
\big| G_{k,j;  n,l,r}^{ess; m, i;p,q}(t,x,\zeta)\big|&\lesssim \| \mathfrak{m}(\cdot, \zeta)\|_{\mathcal{S}^\infty} 2^{2\epsilon M_t} \min\{ 2^{-2m-2p-q-j-l+n },2^{m+3j+l+2p+q+n  } \}\\
&\lesssim 2^{-m/2+j+n+2\epsilon M_t}\| \mathfrak{m}(\cdot,  \zeta)\|_{\mathcal{S}^\infty}. 
\end{split}
\ee
Moreover, after combining the above estimate, and the obtained estimate  \eqref{2022feb9eqn62}, we have
\be\label{sep9eqn31}
\begin{split}
 &\big| G_{k,j; n,l,r}^{ess; m, i;p,q}(t,x,\zeta)\big| \\
 & \lesssim   \| \mathfrak{m}(\cdot, \zeta)\|_{\mathcal{S}^\infty} 2^{n+ \epsilon M_t} \min\{ 2^{-m/2+j+ \epsilon M_t }, 2^{m } 2^{3j+2\max\{l,n+\epsilon M_t\}+l }   \}\\
&\lesssim \| \mathfrak{m}(\cdot, \zeta)\|_{\mathcal{S}^\infty} 2^{n+ \epsilon M_t} \min\big\{ 2^{(k+2l)/2+j+ \epsilon M_t}, \\
&\qquad  ( 2^{-m/2+j+ \epsilon M_t })^{2/3}(2^{m+ 3j+2\max\{l,n+\epsilon M_t\} +l })^{1/3} \big\}\\
& \lesssim  \| \mathfrak{m}(\cdot, \zeta)\|_{\mathcal{S}^\infty}\big[ 2^{(1-20\epsilon) M_{t^{\star}}} + 2^{10\epsilon  M_{t^{\star}} }\mathbf{1}_{n\geq -  M_{t^{\star}}/3-30\epsilon M_{t^{\star}}} \\
&\qquad \times \min\{2^{(k+2n)/2  +   {\alpha}^{\star}  M_{t^{\star}}  }    , 2^{(k+4n)/2+    M_{t^{\star}} }\}  \big].\\
\end{split}
\ee
In the above estimate, we used the fact that $-m-2l\leq  k .$

\medskip

\quad   \textbf{Step 2B-b.}\quad  If  $| x_{\bot}|\leq  2^{m+p-\epsilon M_t}$. 

 \medskip

Note that, for the sub-case we are considering, $\forall y\in B(0, 2^{-k-l+\epsilon M_t/2}),  $ we have
\[
 |  x_{\bot}- y_{\bot}+(t-s) \omega_{\bot}| \sim | (t-s)  \omega_{\bot} |\sim 2^{m+p}.
\]
From the cylindrical symmetry of the distribution function, the estimate of kernels in  \eqref{sep6eqn31}, the estimate of coefficients in \eqref{2022feb9eqn61},  and the volume of support of $\omega$ for fixed $v$,  we have 
\be\label{sep3eqn11}
\begin{split}
 \big| G_{k,j; n,l,r}^{ess; m, i;p,q}(t,x,\zeta)\big| &\lesssim  \| \mathfrak{m}(\cdot, \zeta)\|_{\mathcal{S}^\infty} 2^{m+3k+3n-l+\epsilon M_t} 2^{2l} \frac{2^{-k-n+\epsilon M_t}}{2^{m+p}} 2^{-j} \\
 & \lesssim  \| \mathfrak{m}(\cdot,\zeta)\|_{\mathcal{S}^\infty} 2^{2k+  2n+l-p-j+2\epsilon M_t}.\\
 \end{split} 
\ee
From the above estimate and  the obtained    estimates  \eqref{2022feb9eqn62}   and  \eqref{2022feb9eqn63},  we have
\be\label{sep3eqn16}
\begin{split}
 &\big| G_{k,j; n,l,r}^{ess; m, i;p,q}(t,x,\zeta)\big| \\
 &\lesssim  \| \mathfrak{m}(\cdot, \zeta)\|_{\mathcal{S}^\infty} 2^{\epsilon M_t} \min\{ 2^{-2m- j-3l+n },  2^{m+l + n+ 3j+2\min\{ \max\{l,n+\epsilon M_t\} ,p\}}, 2^{2k+  2n+l-p-j+ \epsilon M_t}\}\\
&\lesssim  \| \mathfrak{m}(\cdot, \zeta)\|_{\mathcal{S}^\infty}2^{\epsilon M_t} \min\big\{ \big(2^{-2m- j-3l+n }\big)^{1/4}  \big(  2^{m+l + n+ 3j+2\min\{ \max\{l,n+\epsilon M_t\} ,p\}}\big)^{3/4}, \\
&\qquad  \big( 2^{m+l + n+ 3j+  2\min\{ \max\{l,n+\epsilon M_t\} ,p\}}\big)^{1/2} \big( 2^{2k+  2n+l-p-j+ \epsilon M_t }\big)^{1/4} \big( 2^{-2m- j-3l+n }\big)^{1/4}\big\}\\
&\lesssim \| \mathfrak{m}(\cdot, \zeta)\|_{\mathcal{S}^\infty} 2^{2\epsilon M_t}\min\{  2^{2j+n+3\max\{l,n+\epsilon M_t\}/2 }, 2^{k/2+5n/4+3 \min\{ \max\{l,n+\epsilon M_t\} ,p\}/4+j }\}
  \\
 &\lesssim  \| \mathfrak{m}(\cdot, \zeta)\|_{\mathcal{S}^\infty}\big[ 2^{(1-20\epsilon)M_{t^{\star}}} + 2^{30\epsilon M_{t^{\star}}} \mathbf{1}_{n\geq -(2/5+30\epsilon)M_{t^{\star}} } \\
 &\qquad \times \min\{2^{(k+2n)/2+ \alpha^{\star}  M_{t^{\star}} } , 2^{(k+4n)/2+   M_{t^{\star}} }\} \big]. \\
 \end{split}
\ee

 To sum up, our desired estimate  \eqref{sep6eqn49}  holds from the obtained estimates  \eqref{aug8eqn23},  \eqref{2022feb9eqn70},    \eqref{sep3eqn6},  \eqref{sep9eqn31}, and  \eqref{sep3eqn16}.  

\end{proof}

It remains to      estimate the $S$-part of the   decomposition   in  \eqref{sep18eqn50}  in Lemma \ref{locdeclemm}.  Similar to the   decomposition we did in   \eqref{sep7eqn51},  
 after using the decomposition of the acceleration force in  \eqref{july1eqn11},  
doing dyadic decomposition   for $\sin \theta, \sin \phi$,   using the volume of support of $v, \omega$ and the rough estimate of the electromagnetic field  \eqref{july10eqn89}  in Proposition \ref{Linfielec} for the case $p=-10M_t$ or $q=-10M_t$, and ruling out the error type case, i.e., the case when $|\omega\times \tilde{\zeta}|\gtrsim 2^{\max\{n, -(m+k)/2\}+\epsilon M_t/2}$,  we have 
\be\label{sep7eqn61}
\begin{split}
 & \big| \widetilde{T}_{k,j;  n,l,r }^{S; \mu ,m, i }(\mathfrak{m}, E)(t,x, \zeta ) +\hat{ \zeta}\times  \widetilde{T}_{k,j;  n,l,r }^{S; \mu ,m, i }(\mathfrak{m}, B)(t,x,  \zeta)\big| \\
 & \lesssim \sum_{\begin{subarray}{c}
  a\in \{1,2\} \\
 p, q\in (-10 M_t,2]\cap \Z\\
 \end{subarray}}\widetilde{K}^{ m;p,q,i;a}_{k,j ;  n,l,r }(t, x,   \zeta) +\| \mathfrak{m}(\cdot, \zeta)\|_{\mathcal{S}^\infty}, 
 \end{split}
\ee
where
\be\label{aug8eqn31}
\begin{split}
\widetilde{K}^{m;p,q,i;a}_{k,j ;  n,l,r }(t, x, \zeta )&= \int_0^t  \int_{\R^3} \int_{\R^3}\int_0^{2\pi} \int_0^{\pi} \big[ (t-s)  K_{k;n}^{ \mu}(\mathfrak{m})(y,  \zeta, \omega)+ \widetilde{K}_{k;n}^{ \mu}(\mathfrak{m})(y,   \zeta )\big] \\
&\quad\times   \varphi_{m;-10M_t }(t-s) f(s,x-y+(t-s)\omega, v) \\
&\quad \times \psi_{ess}(  \omega,\zeta) EB^a(t,s,x  -y +(t-s)\omega ,\omega, v)   \\
&\quad  \cdot \nabla_v \big(\frac{  m_{E}(v, \omega )  + \hat{ \zeta}\times  m_{B}(v, \omega )}{1+\hat{v}\cdot \omega} \varphi_{j,n}^{i; r}(v, \zeta)  \varphi_{l; r}(\tilde{v}+\omega )  \big)   \\ 
  &\quad \times    \varphi_{p;-10 M_t}(\sin \theta) \varphi_{q;-1 0M_t}(\sin\phi) \sin \theta  d\theta d\phi  dy  d v d s.\\
\end{split}
\ee
The cutoff function $  \psi_{ess}(  \omega,\zeta)$ appeared above is defined in  \eqref{feb9eqn71}. 

In the following Lemma, we first estimate $\widetilde{K}^{m;p,q,i;1}_{k,j ;l,n}(t, x, \zeta )$ for the case $i=2.$  
\begin{lemma}\label{smallangleS1part1}
    Let $  (l,r)\in \mathcal{B}_2 $, see  \eqref{sep18eqn50}.     Under the assumption of Theorem \ref{mainresultsfirstpart}, the following estimate holds if  $m+k\geq -2 l +4\epsilon M_t$, 
\be\label{2022feb10eqn35}
\begin{split}
\big| \widetilde{K}^{m;p,q,2;1}_{k,j ; n,l,r}(t, x, \zeta )\big| & \lesssim \| \mathfrak{m}(\cdot,  \zeta)\|_{\mathcal{S}^\infty} \big[ 2^{(1-20\epsilon)M_{t^{\star}} } +   2^{ 100\epsilon M_{t^{\star}} }   \mathbf{1}_{n\geq   (1-2 {\alpha}^{\star}-75\epsilon)M_{t^{\star}} } \\
&\qquad \times  \min\{2^{(k+2n)/2+{\alpha}^{\star} M_{t^{\star}} }, 2^{(k+4n)/2 + (1+3\iota)  M_{t^{\star}}} \}    \big] .\\
\end{split}
\ee
\end{lemma}

\begin{proof}
Recall  \eqref{aug8eqn31}. 
Due to the cutoff functions $\varphi_{j,n}^{2;r}(v, \zeta)$, see \eqref{sep17eqn63} and   \eqref{sep4eqn6}, and the assumption that $m+k\geq -2 l +4\epsilon M_t$, we have $l\geq n+3\epsilon M_t/4>\bar{l}_2$,     and $|\tilde{v}-\tilde{\zeta}|\sim 2^l$. This   implies further that  $2^l \sim 2^r$. Recall   \eqref{2022feb9eqn71}.  As a result of direct computations, we have 
\be\label{2022feb10eqn1}
\big|   \nabla_v \big(\frac{  m_{E}(v, \omega )  + \hat{ \zeta}\times  m_{B}(v, \omega )}{1+\hat{v}\cdot \omega}       \varphi_{j,n}^{i; r}(v, \zeta)  \varphi_{l; r}(\tilde{v}+\omega )    \big)   \big| \lesssim 2^{-j  +\epsilon M_t} 2^{-l}    . 
\ee
  Based on the possible size of $m, k,l,n$, we proceed in two steps as follows. 

\medskip

\noindent \textbf{Step 1.}  \qquad If $m+k\leq -2n+2\epsilon M_t$.  

\medskip

From the estimate  \eqref{2022feb10eqn1}, the estimate of kernel in  \eqref{sep6eqn31},    the estimate  \eqref{march18eqn31}  in Lemma \ref{conservationlawlemma}, the volume of support of $\omega, v$, 
    the Cauchy-Schwartz inequality, and the estimate  \eqref{nov24eqn41}  if $|  v_{\bot}|\geq 2^{(\alpha_t+\epsilon)M_t}$,  we have
\be\label{2022feb10eqn31}
\begin{split}
&\big| \widetilde{K}^{m;p,q,2;1}_{k,j ;  n,l,r }(t, x, \zeta )\big| \\
& \lesssim  (2^m + 2^{-k})  2^{-j-l  +\epsilon M_t}  \| \mathfrak{m}(\cdot,  \zeta )\|_{\mathcal{S}^\infty}     \big(\min\{2^{-2m}, 2^{m+3k+2n}\} 2^{ 3j+2l} \big)^{1/2} \\
&\qquad \times \big(\min\big\{ 2^{m+3k+2n-j}  (2^{2n} + 2^{-m-k})2^{2 \epsilon M_t}  ,   \\
&\qquad    2^{m}\min\{2^{  3j+2l}, 2^{j+2(\tilde{\alpha}_t+\epsilon)M_t}\}(2^{2n} + 2^{-m-k})\big\} \big)^{1/2} \\
&\lesssim 2^{j/2+4\epsilon M_t}   \big(\min\{2^{2k+2n-j}, 2^{-k  }\min\{2^{  3j+2l}, 2^{j+2\tilde{\alpha}_t M_t}\}\} \big)^{1/2}\| \mathfrak{m}(\cdot, \zeta)\|_{\mathcal{S}^\infty} \\
& \lesssim  2^{j/2+4\epsilon M_t}\| \mathfrak{m}(\cdot, \zeta)\|_{\mathcal{S}^\infty}  \min\big\{ \big(2^{2k+2n-j}\big)^{1/4} \big(2^{-k  }\min\{2^{  3j+2l}, 2^{j+2\tilde{\alpha}_t M_t}\}  \big)^{1/4}, \\
  &\qquad  \big(2^{2k+2n-j}\big)^{1/6} \big(2^{-k  }\min\{2^{  3j+2l}, 2^{j+2\tilde{\alpha}_t M_t}\}  \big)^{1/3} \big\}\\
&\lesssim \| \mathfrak{m}(\cdot, \zeta)\|_{\mathcal{S}^\infty} 2^{ 2\epsilon M_t} \min\{2^{k/4+n/2+j/2}\min\big\{2^{(j+l)/2}, 2^{\tilde{\alpha}_t M_t/2}\},  \\
&\qquad   2^{2j/3 +n/3}\min\{2^{2\tilde{\alpha}_t M_t/3},2^{2j/3+2l/3}\} \big\} 
\end{split}
\ee 
From the above estimate, we conclude that 
\be\label{2022feb10eqn32}
\begin{split}
\big| \widetilde{K}^{m;p,q,2;1}_{k,j ;  n,l,r }(t, x, \zeta )\big| & 
\lesssim \| \mathfrak{m}(\cdot,  \zeta)\|_{\mathcal{S}^\infty} \big[ 2^{(1-20\epsilon)M_{t^{\star}} } +   2^{ 100\epsilon M_{t^{\star}} }     \mathbf{1}_{n\geq   (1-2 {\alpha}^{\star}-75\epsilon)M_{t^{\star}} } \\
&\qquad \times \min\{2^{(k+2n)/2+{\alpha}^{\star} M_{t^{\star}} }, 2^{(k+4n)/2 + (1+3\iota)  M_{t^{\star}}} \}  \big] .
\end{split}
\ee
 
\medskip

\noindent \textbf{Step 2.} \qquad If $m+k\geq -2n+\epsilon M_t$, i.e., $-m/2-k/2\leq n -\epsilon M_t/2$. 

\medskip

From the estimate of kernel in  \eqref{sep6eqn31},  the estimate  \eqref{march18eqn31}  in Lemma \ref{conservationlawlemma},   and the estimate  \eqref{nov24eqn41}  if $| v_{\bot}|\geq 2^{(\alpha_t+\epsilon)M_t}$,   we have   
\be\label{2022feb10eqn33}
\begin{split}
 &\big| \widetilde{K}^{m;p,q,2;1}_{k,j ;  n,l,r }(t, x, \zeta )\big| \\
 & \lesssim  \| \mathfrak{m}(\cdot, \zeta)\|_{\mathcal{S}^\infty} 2^{m-j-l+4\epsilon M_t}    \big(2^{-2m }\min\{2^{3j+2l},2^{j+2\tilde{\alpha}_t M_t}\} \big)^{1/2}     \\
&\qquad \times  \big( \min\{   2^{m+2n}\min\{2^{3j+2l},2^{j+2 \tilde{\alpha}_t  M_t}\}, 2^{-2m-j-2l}\}\big)^{1/2}     \\
& \lesssim 2^{4\epsilon M_t}\min\{2^{j/2}, 2^{\tilde{\alpha}_t M_t-j/2-l}\}  \| \mathfrak{m}(\cdot, \zeta)\|_{\mathcal{S}^\infty}   \min\big\{ \big( 2^{m+2n}\min\{2^{3j+2l},2^{j+2 \tilde{\alpha}_t  M_t}\} \big)^{1/3}\\
&\qquad \times  \big( 2^{-2m-j-2l}\big)^{1/6},   \big( 2^{m+2n} 2^{3j+2l}  \big)^{1/4} \big( 2^{-m+k+2n-2l-j}\big)^{1/4} \big\}\\
&\lesssim  \| \mathfrak{m}(\cdot, \zeta)\|_{\mathcal{S}^\infty} 2^{4\epsilon M_t} \min\{ 2^{ (k+4n)/4  +  j/2 }\min\{2^{j/2},2^{\tilde{\alpha}_t M_t/2-l/2}\} ,2^{j +2n/3+\tilde{\alpha}_t M_t/3}\}\\
& \lesssim  \| \mathfrak{m}(\cdot,  \zeta)\|_{\mathcal{S}^\infty} \big[ 2^{(1-20\epsilon)M_{t^{\star}} } +   2^{ 50\epsilon M_{t^{\star}} }  \mathbf{1}_{n\geq   - {\alpha}^{\star} M_{t^{\star}}/2-75\epsilon M_{t^{\star}} } \\
&\qquad \times  \min\{2^{(k+2n)/2+{\alpha}^{\star} M_{t^{\star}} }, 2^{(k+4n)/2 +   M_{t^{\star}}} \}\big].
\end{split}
\ee

In the above estimate, we used the fact that $-m\leq k + 2n.$  To sum up, our desired estimate   \eqref{2022feb10eqn35}  holds from the obtained estimates  \eqref{2022feb10eqn32} and  \eqref{2022feb10eqn33}.

\end{proof}

\begin{lemma}\label{smallangleS1part5}
 Let $ i\in\{3,4\}, (l,r)\in \mathcal{B}_i $, see  \eqref{sep18eqn50}.    Under the assumption of Theorem \ref{mainresultsfirstpart}, the following estimate holds if  $m+k\geq 
-2 l +4\epsilon M_t$, 
\be\label{2022feb10eqn66}
\begin{split}
\big| \widetilde{K}^{m;p,q,i;1}_{k,j ; n,l,r }(t, x, \zeta )\big|& \lesssim   \| \mathfrak{m}(\cdot, \zeta)\|_{\mathcal{S}^\infty}\big[2^{(1-20 \epsilon) M_{t^{\star}} } +  2^{40\epsilon M_{t^{\star}} } \mathbf{1}_{n\geq - (1/2+20\epsilon) M_{t^{\star}} } \\
&\qquad \times \min\{2^{(k+2n)/2+2M_{t^{\star}}/3},2^{(k+4n)/2+ 7 M_{t^{\star}}/6}\} \big].
\end{split}
\ee
\end{lemma}
\begin{proof}
Due to the cutoff functions $\varphi_{j,n}^i(v, \zeta), i\in\{3,4\}$, see  \eqref{sep4eqn6}, the support of the $v, \omega, \xi$,    implies that $|\tilde{v}\times \omega|\lesssim  2^{n+\epsilon M_t} +2^{l-\epsilon M_t}, |\tilde{v}\times\omega |\sim 2^l$.

If $l> -j$,   then it follows that $l\leq n +\epsilon M_t$, since  in this case,   $ |\omega \times\tilde{\xi}|\sim 2^l$.   This further implies that  $l\in [-j, \max\{-j, n+\epsilon M_t\}]\cap \Z. $ Moreover, recall  \eqref{2022feb9eqn71},   as a result of direct computations, we have 
\be\label{2022feb10eqn61}
\big|   \nabla_v \big(\frac{  m_{E}(v, \omega )  + \hat{ \zeta}\times  m_{B}(v, \omega )}{1+\hat{v}\cdot \omega}       \varphi_{j,n}^{i; r}(v, \zeta)  \varphi_{l; r}(\tilde{v}+\omega )     \big| \lesssim 2^{-j-l-\min\{l,n\}+ n +\epsilon M_t} . 
\ee

  Based on the possible sizes of $p $, we   proceed in two steps as follows. 

\medskip

\noindent \textbf{Step 1.} \quad  If $ p\leq l+ \epsilon M_t. $

 \medskip
 
From the estimate  \eqref{2022feb10eqn61}, the estimate of kernels $ K_{k;n}^{ \mu}(\mathfrak{m})(y,   \zeta, \omega)$ and $ \widetilde{K}_{k;n}^{ \mu}(\mathfrak{m})(y, \zeta)$ in  \eqref{sep6eqn31}, the estimate  \eqref{march18eqn31}  in Lemma \ref{conservationlawlemma}, and the volume of support of $v,\omega$, we have
\be\label{sep21eqn1}
\begin{split}
& \big|\widetilde{K}^{ m;p,q,3;1}_{k,j ;  n,l,r }(t,x,  \zeta ) \big|\\
 &\lesssim  2^{m-j-l-\min\{l,n\}+n+\epsilon M_t} \big( 2^{-2m+3j+2l} \big)^{1/2}\big(\min\{   2^{m+2l+3j+2 p}, 2^{-2m-j-2l} \} \big)^{1/2} \| \mathfrak{m}(\cdot,  \zeta)\|_{\mathcal{S}^\infty}\\
 &\lesssim 2^{n+j/2-\min\{l,n\}+2\epsilon M_t} \| \mathfrak{m}(\cdot,  \zeta)\|_{\mathcal{S}^\infty} \min\big\{(  2^{m+2l+3j+2 p})^{1/2},\\ 
 &\qquad  (2^{m+4l+3j})^{1/3} (2^{-2m-j-2l})^{1/6},   (2^{-2m-j-2l})^{1/2}  \big\}\\
&  \lesssim 2^{  4\epsilon M_t}\| \mathfrak{m}(\cdot,  \zeta)\|_{\mathcal{S}^\infty}  \big(   \min\{2^{n+4j/3},  2^{k+ n }\} \mathbf{1}_{l> -j} + \min\{ 1+ 2^{n+j}, 2^{k/2+n}\}\mathbf{1}_{l=-j} \big)\\
&\lesssim \| \mathfrak{m}(\cdot, \zeta)\|_{\mathcal{S}^\infty}\big[2^{(1-20\epsilon)M_{t^{\star}} } + 2^{35\epsilon M_{t^{\star}}}\mathbf{1}_{n\geq -M_{t^{\star}}/3-30\epsilon M_{t^{\star}} } \\
&\qquad \times \min\{2^{(k+2n)/2+ 2 M_{t^{\star}}/3 }, 2^{(k+4n)/2+ M_{t^{\star}} } \} \big]. 
\end{split}
\ee

\medskip

\noindent \textbf{Step 2.} \quad    If $ p\geq l+ \epsilon M_t .  $

\medskip

 Note that, due to the  cutoff function $\varphi_{l;-j}(\omega+\tilde{v})\varphi_{p;-10M_t}(\sin\theta)$, we have $|  v_{\bot}|\sim 2^{p+j}.   $
  From the estimate in  \eqref{nov24eqn41}, we can rule out the case $| v_{\bot}|\geq 2^{(\alpha_t+\epsilon)M_t}$, i.e.,  $p+j\geq \alpha_t M_t+\epsilon M_t$. Now, it would be
sufficient to consider the case $p+j\leq   {\alpha}_t M_t+\epsilon M_t $. 

Based on the possible size of $|  x_{\bot}|$, we  proceed in two sub-steps as follows.

\medskip

  \textbf{Step 2A.} \quad  If $| x_{\bot}|\geq 2^{m+p-10}$. 

\medskip

Similar to what we did in \eqref{march18eqn66}, after using the Cauchy-Schwarz inequality, the estimate  \eqref{march18eqn31}  in Lemma \ref{conservationlawlemma}, and the estimate of kernels in \eqref{sep6eqn31},  we have 
\be\label{sep21eqn2}
\begin{split}
&|\widetilde{K}^{  m;p,q,3;1}_{k,j ;  n,l,r}(t, x, \zeta  )| \\
& \lesssim   \| \mathfrak{m}(\cdot,  \zeta)\|_{\mathcal{S}^\infty} 2^{m-j -l-\min\{l,n\}+n+\epsilon M_t }  \big(2^{-2m} 2^{3j+2l}  \big)^{1/2} \\
&\qquad \times    \big(\min\{  | {x}_{\bot} |^{-1}  {2^{-m-p-q-j}} , 2^{-2m-j-2l},  2^{m+3j+2l}\min\{2^{2p+q}, 2^{2 n+2\epsilon M_t} \} \big)^{1/2}  
 \\
  &   \lesssim   \| \mathfrak{m}(\cdot, \zeta)\|_{\mathcal{S}^\infty} 2^{j/2-\min\{l,n\}+n+2\epsilon M_t}  \big(\min\{   2^{-m/2+j+l}, 2^{-2m-j-2l},  2^{m+3j+2l+2n} \}\big)^{1/2}  \\
&\lesssim  \| \mathfrak{m}(\cdot, \zeta )\|_{\mathcal{S}^\infty}  2^{   4\epsilon M_t} \big(1  + \min\{2^{k/4+j+n }, 2^{k+n },  2^{m/2+2j+2n} \}  \mathbf{1}_{l\leq n+\epsilon M_t} \big)\\
& \lesssim \| \mathfrak{m}(\cdot, \zeta)\|_{\mathcal{S}^\infty}\big[2^{(1-20\epsilon)M_{t^{\star}}} + 2^{30\epsilon M_{t^{\star}} }\mathbf{1}_{n\geq - (1+40\epsilon )M_{t^{\star}}/2 } \\
&\qquad \times \min\{2^{(k+2n)/2 + 2 M_{t^{\star}}/3 }  , 2^{(k+4n)/2+  M_{t^{\star}} } \} \big].\\
\end{split}
\ee 

\medskip

  \textbf{Step 2B.} \quad  If $| x_{\bot}|\leq 2^{m+p-10}$.

 \medskip

Let$\{\tilde{\zeta}, \tilde{\zeta}_1, \tilde{\zeta}_2\}$ be a fixed   orthonormal frame.  From   the Cauchy-Schwarz inequality,   the estimate of kernels in  \eqref{sep6eqn31}, and the estimate  \eqref{march18eqn31}  in Lemma \ref{conservationlawlemma},  we have
\be\label{sep21eqn52}
 |\widetilde{K}^{  m;p,q,3;1}_{k,j ;  n,l,r }(t, x, \zeta  )|\lesssim   2^{m-j-l-\min\{l,n\}+n+\epsilon M_t}  \big(2^{-2m} 2^{3j+2l}  \big)^{1/2} \big( J^{m;p, q}_{k,j ;  n,l  }(t, x,  \zeta ) \big)^{1/2}\| \mathfrak{m}(\cdot, \zeta)\|_{\mathcal{S}^\infty},
\ee
where
\be\label{sep26eqn21}
\begin{split}
  J^{m;p, q}_{k,j ;  n,l  }(t, x,\zeta)&:=\int_0^{t} \int_{\R^3}\int_{\R^3} \int_0^{2\pi} \int_0^{\pi} 2^{3k+2n}  f(s, x-y+(t-s)\omega,  \zeta )\\
  &\qquad  \times (1+2^k|y\cdot\tilde{\zeta}|)^{-N_0^3} (1+2^{k+n}(|y\cdot \tilde{\zeta}_1| + |y\cdot \tilde{\zeta}_2|) )^{-N_0^3}\\
&\qquad \times  \varphi_j(v)  \psi_{<  n+\epsilon M_t   }( \tilde{\zeta} - \tilde{v} )   \varphi_{l; -j}(\tilde{v}+\omega )    
  \varphi_{m;-10M_t }(t-s)\\
  &\qquad \times  \varphi_{p;-10 M_t}(\sin \theta) \varphi_{q;-10 M_t}(\sin \phi)  \sin \theta d \theta d \phi d y d v d s. \\
  \end{split}
\ee
Note that, on the Fourier side, after localizing the frequency around the fixed direction of $v$, we have
\be\label{sep21eqn71}
\begin{split}
 J^{m;p, q }_{k,j ;  n,l  }(t, x, \zeta )&=\sum_{\star\in\{ess,err\} }  J^{m;p,q;\star }_{k,j ;  n,l  }(t, x, \zeta ),\\ 
  J^{m;p, q;\star }_{k,j ;  n,l   }(t, x, \zeta )&:=\int_0^{t} \int_{\R^3}\int_{\R^3} \int_0^{2\pi} \int_0^{\pi}  e^{i(x+(t-s)\omega)\cdot \xi} \hat{f}(s, \xi, v)\widehat{\mathfrak{R}}_{k,n}(\xi)  \\
  &\qquad \times    \varphi_{l;-j}(\tilde{v}+\omega)    \varphi_j(v)  \psi_{<  n+\epsilon M_t   }( \tilde{\zeta} - \tilde{v} )   \varphi_{m;-10M_t }(t-s)
 \\ 
 & \qquad \times \varphi_{p;-10 M_t}(\sin \theta)   \varphi_{q;-10 M_t}(\sin \phi )  \psi_{\star} (v,\xi) \sin \theta d \theta d \phi d \xi d v d s, 
\end{split}
\ee
where $ \psi_{ess}(v,\xi):=\psi_{\leq l +5}(\tilde{v}\times \tilde{\xi}), \psi_{err}(v,\xi):=\psi_{ > l+5 } (\tilde{v}\times \tilde{\xi})$ and $\widehat{\mathfrak{R}}_{k,n}(\xi) $ is defined as follows, 
\[
\begin{split}
\widehat{\mathfrak{R}}_{k,n}(\xi) &:= \int_{\R^3}e^{i y\cdot \xi}  2^{3k+2n}  (1+2^k|y\cdot \tilde{\zeta}|)^{-N_0^3} (1+2^{k+n}( |y\cdot \tilde{\zeta}_1| + |y\cdot \tilde{\zeta}_2|) )^{-N_0^3} d y \\  &=\widehat{\mathfrak{R}}( 2^{-k}(\xi\cdot \tilde{\zeta}), 2^{-k-n} (\xi\cdot \tilde{\zeta}_1), 2^{-k-n}(\xi\cdot \tilde{\zeta}_2) ),
\end{split}
\]
where  
\be\label{sep5eqn89}
\widehat{\mathfrak{R}}(\xi):= \int_{\R^3} e^{iy \cdot \xi} (1+|y\cdot \tilde{\zeta} |)^{-N_0^3}(1+(|y\cdot \tilde{\zeta}_1| +|y\cdot \tilde{\zeta}_2| ) )^{-N_0^3}  d y. 
\ee

Note that $m+k\geq -2l+ 4\epsilon M_t  $ and $p\geq l +\epsilon M_t$. By doing integration by parts in $\omega$ many times for the error type term, we have
\be\label{sep21eqn51}
|  \sum_{q\in(-10M_t,2]\cap \Z} J^{m;p, q;err }_{k,j ;  n,l  }(t, x, \zeta)|\lesssim 2^{-100M_t }. 
\ee
For the essential part $ J^{m;p;ess }_{k,j ;l,n}(t, x, V )$, in terms of kernel, we have
\be\label{sep21eqn88}
\begin{split}
\big| J^{m;p, q;ess}_{k,j ;  n,l  }(t, x, \zeta)\big| & \lesssim \int_0^{t} \int_{\R^3}\int_{\R^3} \int_0^{2\pi} \int_0^{\pi} K_{c(m,k,l),n}(y, \zeta, v) f(s, x-y+(t-s)\omega, v) \\
&\qquad  \times   \varphi_j(v)  \psi_{<  n+\epsilon M_t   }( \tilde{\zeta} - \tilde{v} )   \varphi_{l; -j}(\tilde{v}+\omega ) \varphi_{m;-10M_t}(t-s) \\
&\qquad \times \varphi_{p;-1 0M_t}(\sin \theta) \varphi_{q;-1 0M_t}(\sin \phi)   \sin \theta d \theta d \phi d y d v d s,\\
\end{split}
\ee
where
\be\label{sep21eqn89}
K_{c(m,k,l),n}(y, \zeta, v):=\int_{\R^3} e^{iy\cdot \xi}  \psi_{\leq l }(\tilde{v}\times \tilde{\xi} )    \widehat{\mathfrak{R}}( 2^{-k}(\xi\cdot \tilde{\zeta}), 2^{-k-n} (\xi\cdot \tilde{\zeta}_1), 2^{-k-n}(\xi\cdot \tilde{\zeta}_2) )  d\xi. 
\ee

Note that $ | \tilde{v} \times \tilde{\zeta}|\lesssim 2^{n+\epsilon M_t}$ and $|  v_{\bot}|/|v|\sim 2^p$. By doing integration by parts in $\xi$ along $v$ direction and directions perpendicular to $v$ or alternatively, along $e_i, i\in\{1,2,3\},$ directions many times, the following estimate holds for the kernel $ K_{c(m,k,l),n}(y, \zeta, \omega)$, 
\be\label{sep26eqn23}
\begin{split}
|K_{c(m,k,l),n}(y, \zeta, v)|& \lesssim 2^{3k +\epsilon^2 (k+M_t) +2\min\{l,n\} }\min\big\{(1+2^{k-\epsilon M_t}|y  \cdot \tilde{v}|)^{-N_0^3}\\
  &\quad  \times  (1+2^{k+\min\{l,n\} }|y\times \tilde{v}|)^{-N_0^3}, (1+2^{k+\min\{l,n\}}|y_{\bot}|)^{-N_0^3}\\
 &\quad \times  (1+ (2^{-k-l+p} +2^{-k-n+\max\{n+\epsilon M_t,p\}})^{-1} |y_3| )^{-N_0^3} \big\}. \\
 \end{split}
\ee

Note that, for the case we are considering, we have $m+p\geq -k-l+\epsilon M_t$ and $| x_{\bot} -  y_{\bot} +(t-s)  \omega_{\bot}|\sim 2^{m+p}$ if $|  y_{\bot}|\leq 2^{-k-l+\epsilon M_t/2}$. Moreover,  the following bounds hold for the volumes of the supports of $\omega$($v$ fixed) and $\phi$ ($v$ and $\theta$ fixed):
\[
\begin{split}
|Vol_{ess,v}(supp(\omega))|&\lesssim \min\{2^{2p+q}, 2^{2l}, 2^{p+q+l}\},\\
 |Vol_{ess,v,\theta}(\phi)|&\lesssim  \min\{ 2^{l-p},1\}.
\end{split}
\]

From the above estimate of kernel, the Jacobian of changing coordinates $(y_1,y_2,\theta)\rightarrow x-y+(t-s)\omega$, the volume of support of $\omega$, and the cylindrical symmetry of solution,  we have
\be 
\begin{split}
\big| J^{m;p, q;ess}_{k,j ;  n,l  }(t, x, \zeta )\big| & \lesssim     2^{3k+2\min\{l,n\}   +\epsilon (k+2M_t)  }2^{-k-\min\{l,n\} } \\
&\quad \times  (2^{-k-l+p} +2^{-k-n+\max\{n+\epsilon M_t,p\}})  2^{  l-p } 2^{-j}   2^{-(m+p)  }\\
&\lesssim 2^{-m+k+l-p-j+\epsilon^2  k+ 3\epsilon M_t }. 
\end{split}
\ee

After combining the above estimate, the estimate  \eqref{sep21eqn51}, and the estimate  \eqref{sep21eqn52}  and using the volume of support of $v, \omega$,    we have
\be\label{2024Dec8eqn5}
\begin{split}
  &|\widetilde{K}^{  m;p,q,3;1}_{k,j ;  n,l,r }(t, x, \zeta  )| \\
  & \lesssim     \| \mathfrak{m}(\cdot, \zeta)\|_{\mathcal{S}^\infty}  2^{m-j-l-\min\{l,n\}+n+\epsilon M_t}  \big(2^{-2m} 2^{3j+2l}  \big)^{1/2}  \\
& \qquad \times \big(\min\{2^{-m+ k+l-p-j+ \epsilon^2  k+ 3\epsilon M_t} +2^{-100M_t }, 2^{m }  2^{3j+2l+2\min\{\max\{l,n+\epsilon M_t\} ,p\}} \}  \big)^{1/2}.  \\
\end{split}
\ee
From the above estimate, we have
\be\label{aug25eqn31}
\begin{split}
\eqref{2024Dec8eqn5}& \lesssim 2^{j/2+\max\{l,n\}-l+2\epsilon M_t} \| \mathfrak{m}(\cdot, \zeta)\|_{\mathcal{S}^\infty} \big[1+\big(\min\big\{2^{m+3j+2l+2\max\{l,n\}},\\
& \qquad  (2^{-m+ k+l-p-j+ \epsilon^2  k+ 3\epsilon M_t})^{(2-2\epsilon^2)/3} (2^{m+3j+2l+2p})^{(1+2\epsilon^2)/3}  \}\big)^{1/2}\big]\\
&\lesssim \| \mathfrak{m}(\cdot, \zeta)\|_{\mathcal{S}^\infty} \big( 1+ 2^{10\epsilon M_t}\min\{  2^{m/2+ 2j+2\max\{l,n \}}, 2^{-m/6+ k/3-l/3+\max\{l,n \}+2j/3}\} \big)\\
& \lesssim  \| \mathfrak{m}(\cdot, \zeta)\|_{\mathcal{S}^\infty} \big( 1+ 2^{10\epsilon M_t}\min\{  2^{m/2+ 2j+2\max\{l,n \}},  2^{    k/2 + \max\{l,n \}+2j/3}\}  \big)  \| \mathfrak{m}(\cdot, \zeta)\|_{\mathcal{S}^\infty}\\
 &\lesssim  \| \mathfrak{m}(\cdot, \zeta)\|_{\mathcal{S}^\infty}\big[2^{(1-20 \epsilon) M_{t^{\star}} } +  2^{40\epsilon M_{t^{\star}} } \mathbf{1}_{n\geq - (1/2+20\epsilon) M_{t^{\star}} }\\
 &\qquad \times \min\{2^{(k+2n)/2+2M_{t^{\star}}/3},2^{(k+4n)/2+ 7 M_{t^{\star}}/6}\} \big].
 \end{split}
\ee
Recall    \eqref{sep7eqn61}. To sum up,  after combining  the obtained estimates    \eqref{sep21eqn1}, \eqref{sep21eqn2}, and  \eqref{aug25eqn31}, our desired estimate  \eqref{2022feb10eqn66}  holds.
\end{proof}

\begin{lemma}\label{smallangleS2part1}
  Let $  (l,r)\in \mathcal{B}_2 $, see  \eqref{sep18eqn50}.     Under the assumption of Theorem \ref{mainresultsfirstpart} and the assumption  $m+k\geq -2 l +4\epsilon M_t$, the following estimate holds  if   $n \geq   (1-2\alpha^{\star}-35\epsilon)M_{t^{\star}} $,   $    j\leq    (1/2+3\iota + 55\epsilon) M_{t^{\star}} $, or  $k+2j\leq (2-50\epsilon)  M_{t^{\star}}$,
\be\label{2022feb11eqn21}
\begin{split}
  \big| \widetilde{K}^{m;p,q,2;2}_{k,j ; n,l,r  }(t, x, \zeta )\big|\lesssim \| \mathfrak{m}(\cdot, \zeta)\|_{\mathcal{S}^\infty}\big[ 2^{  (1-20\epsilon)   M_{t^{\star}}} + 2^{120\epsilon  M_{t^{\star}} }\mathbf{1}_{n \geq  (1-2\alpha^{\star}-35\epsilon)M_{t^{\star}} } \\
 \qquad \times \min\{ 2^{(k+2n)/2+ (\alpha^{\star}+3\iota) M_t },2^{(k+4n)/2+ (1+6\iota) M_{t^{\star}}  } \}\big].\\
 \end{split}
\ee
\end{lemma}
\begin{proof}
Due to the cutoff function  $\varphi_{j,n}^2(v, \zeta) $, see  \eqref{sep4eqn6}, and the assumption that $m+k\geq -2 l +4\epsilon M_t$, we have $l\geq n+3\epsilon M_t/4> \bar{l}_2$,      and $|\tilde{v}-\tilde{\zeta}|\sim 2^l$. This implies further that $2^l \sim 2^r$. Note that, besides the rough estimate of the gradient of coefficients in  \eqref{2022feb10eqn1}, the following improved estimate   holds for the $\p_{v_3}$ derivative of coefficients, 
\be\label{2022feb16eqn71}
\big|   \p_{v_3} \big(\frac{  m_{E}(v, \omega )  + \hat{ \zeta}\times  m_{B}(v, \omega )}{1+\hat{v}\cdot \omega}       \varphi_{j,n}^{i; r}(v, \zeta)  \varphi_{l; r}(\tilde{v}+\omega )     \big)   \big| \lesssim 2^{-j  -l+\max\{l,p\}+  \epsilon M_t}  .  
\ee

Recall  \eqref{aug8eqn31}  and  \eqref{july1eqn13}.    From the estimate of kernels in  \eqref{sep6eqn31},  the above estimate and the estimate of coefficients in  \eqref{2022feb10eqn1}, we have 
\be\label{2022feb11eqn86}
\begin{split}
&\big| \widetilde{K}^{m;p,q,2;2}_{k,j ; n,l,r  }(t, x, \zeta ) \big| \\
&\lesssim \int_{0}^{t }\int_{\R^3}\int_{\R^3} \int_0^{2\pi} \int_{0}^{\pi} 2^{m+3k+2n-j  + \max\{p, l  \}  +2\epsilon M_t } \varphi_{m;-10M_t  }(t-s)     \psi_{ess}(  \omega,\zeta)   \\
&\quad \times   f(s, x-y+(t-s)\omega, v) \big|B  (s, x-y+(t-s)\omega)\big|\\
&\quad  \times \| \mathfrak{m}(\cdot, \zeta)\|_{\mathcal{S}^\infty} \min\big\{ (1+2^k|y\cdot\tilde{\zeta}|)^{- N_0^{ 3}} (1+2^{k+n}|y\times \tilde{\zeta}| )^{- N_0^{ 3} }, \\
&\qquad  (1+2^{k+n}|  y_{\bot}| )^{- N_0^{ 3} } (1+2^{k+n } (2^n + |\tilde{\zeta}\times e_3|)^{-1} |y_3|)^{-N_0^3} \big\}\\
& \quad  \times \varphi_{p;-10M_t}(\sin \theta ) \varphi_{q;-10M_t}(\sin\phi )   \varphi_{[l-1,l+1]; \bar{l}_2 }(\tilde{v}+\omega)    \sin \theta d\theta d \phi  d v d y  d s. 
\end{split}
\ee

Based on  different sizes of parameters, e.g., $m+k,l,n,m+p, k+n, |  x_{\bot}|$,  we  proceed in four steps as follows.

 \medskip

\noindent \textbf{Step 1.} \quad   If $   m+k\leq -2n+\epsilon M_t $. 

 \medskip

Recall  \eqref{feb9eqn71}. For this case,  the volume of support of $\omega$ for  the main part is bounded by $2^{-m-k+\epsilon M_t}$. Therefore, after using the estimate  \eqref{march18eqn31}  in Lemma \ref{conservationlawlemma}, the volume of support of $\omega$ and $ v$, and the estimate \eqref{nov24eqn41}  if $| v_{\bot}|\geq 2^{(\alpha_t+\epsilon)M_t}$, we have
 \be\label{2022feb10eqn80}
\begin{split}
&\big| \widetilde{K}^{m;p,q,2;2}_{k,j ; n,l,r  }(t, x, \zeta ) \big|\\
&\lesssim   \sup_{s\in [0,t]}\|B (s,\cdot)\|_{L^2} \| \mathfrak{m}(\cdot, \zeta)\|_{\mathcal{S}^\infty}    2^{m-j + \max\{p, l  \}  +3 \epsilon M_t}  \big(2^{2k+2n}  \min\{2^{j+2(\alpha_t+\epsilon) M_t}, 2^{3j+2l }\}  \big)^{1/2}\\
&\qquad    \times \big(\min\{2^{  -k  }  \min\{2^{j+2(\alpha_t+\epsilon) M_t}, 2^{3j+2l }\} , 2^{2k+2n-j  }, 2^{-2m -j-2l}   \} \big)^{1/2} \\
&  \lesssim  \| \mathfrak{m}(\cdot, \zeta)\|_{\mathcal{S}^\infty}  2^{m-j    +4 \epsilon M_t} \min\big\{ (2^{ 2k+2n+3j+2l})^{1/2}(2^{-2m-j-2l})^{1/2}, \\
& \qquad   \big(2^{  -k  }  \min\{2^{j+2\alpha_t M_t}, 2^{3j+2l }\} \big)^{1/2}  \big(2^{-m+k}  \min\{2^{j+2\alpha_t M_t}, 2^{3j+2l }\}  \big)^{1/2}\big\}\\
 &   \lesssim   2^{5\epsilon M_t} \min\big\{ 2^{k+n}, 2^{m/2 } \min\{2^{2j+2l}, 2^{2\alpha_t M_t} \}  \big\} \| \mathfrak{m}(\cdot, \zeta)\|_{\mathcal{S}^\infty} 
\end{split}
\ee
 From the above estimate, we conclude that  
\be\label{2022feb10eqn81}
\begin{split}
\big| \widetilde{K}^{m;p,q,2;2}_{k,j ; n,l,r  }(t, x, \zeta ) \big| & \lesssim \| \mathfrak{m}(\cdot, \zeta)\|_{\mathcal{S}^\infty} \big[ 2^{(1-20\epsilon)M_{t^{\star}}} +   2^{120\epsilon  M_{t^{\star}} }\\
&\qquad \times \min\{ 2^{(k+2n)/2+ (\alpha^{\star}+3\iota) M_{t^{\star}}  },2^{(k+4n)/2+ (1+6\iota) M_{t^{\star}}  } \} \big] .\\
\end{split}
\ee

 \medskip

\noindent \textbf{Step 2.} \quad If $m+k\geq -2n+\epsilon M_t$ and  $m+p\leq -k-n+2\epsilon M_t$. 

 \medskip

 For this case,  the volume of support of $\omega$ for  the main part is bounded by $2^{2n+2\epsilon M_t}$. From the Cauchy-Schwarz inequality,  the estimate  \eqref{march18eqn31}  in Lemma \ref{conservationlawlemma},    the volume of support of  $\omega$ and $ v$, and the estimate  \eqref{nov24eqn41}  if $|  v_{\bot}|\geq 2^{(\alpha_t+\epsilon)M_t}$,    we have
\be\label{2022feb10eqn82}
\begin{split}
 \big| \widetilde{K}^{m;p,q,2;2}_{k,j ; n,l,r  }(t, x, \zeta ) \big|  
& \lesssim    \big( \sup_{s\in [0,t]}\|B (s,\cdot)\|_{L^2} \big) \| \mathfrak{m}(\cdot, \zeta)\|_{\mathcal{S}^\infty}   2^{m-j + \max\{p, l  \}  +3 \epsilon M_t}\\
&\qquad \times  \big(2^{m+3k+2n+2\min\{p,n\}}    \min\{2^{j+2(\alpha_t+\epsilon) M_t}, 2^{3j+2l }\}   \big)^{1/2}\\
& \qquad  \times    \big(\min\{2^{m+2\min\{p,n\} }   \min\{2^{j+2(\alpha_t+\epsilon) M_t}, 2^{3j+2l }\}  ,  2^{-2m  -j-2l}   \} \big)^{1/2}\\
&   \lesssim   \| \mathfrak{m}(\cdot, \zeta)\|_{\mathcal{S}^\infty} 2^{4\epsilon M_t} \min\big\{      2^{-k/2-n} \min\{2^{2j+2l}, 2^{2\alpha_t M_t} \} , \\
&\qquad   2^{m/2 } \min\{2^{2j+2l}, 2^{2\alpha_t M_t} \} ,  2^{k+n/2+\min\{p,n\}/2} \big\}.
 \end{split}
 \ee
From the above estimate, we conclude that
\be\label{2022feb10eqn83}
\begin{split}
 &\big| \widetilde{K}^{m;p,q,2;2}_{k,j ; n,l,r  }(t, x, \zeta ) \big|  \\
 &\lesssim  \| \mathfrak{m}(\cdot, \zeta)\|_{\mathcal{S}^\infty} 2^{4\epsilon M_t}\min\big\{(2^{k+n})^{1/2} ( 2^{m/2+2j+2l})^{1/2},  (2^{-k/2-n+2j+2l})^{2/3}(2^{k+n})^{1/3} \big\}\\
&\lesssim \| \mathfrak{m}(\cdot, \zeta)\|_{\mathcal{S}^\infty}2^{4.5\epsilon M_t} \min\{2^{(k+n)/2+j}, 2^{4(j+l)/3}   2^{ -n/3+4\epsilon M_t} \}\\
& \lesssim 2^{  (1-20\epsilon)   M_{t^{\star}}}\| \mathfrak{m}(\cdot, \zeta)\|_{\mathcal{S}^\infty}. 
\end{split}
\ee

 \medskip

\noindent \textbf{Step 3.} \quad If $m+k\geq -2n+\epsilon M_t$, $m+p\geq -k-n+2\epsilon M_t$  and $| x_{\bot}|\geq 2^{m+p-10}$. 

 \medskip

 Recall \eqref{2022feb11eqn86}. Similar to the obtained estimate \eqref{aug5eqn44}, by using the same strategy,  we have   
  \be\label{2022feb10eqn86} 
  \begin{split}
 &\big| \widetilde{K}^{m;p,q,2;2}_{k,j ; n,l,r  }(t, x, \zeta ) \big| \\
  & \lesssim       \big( \sup_{s\in [0,t]}\|B (s,\cdot)\|_{L^2} \big)  \| \mathfrak{m}(\cdot, \zeta)\|_{\mathcal{S}^\infty}     2^{m-j  + \max\{p, l  \}  +4 \epsilon M_t} \\
 &\qquad  \times    \big( \min\{ 2^{m+2\min\{n,p\}  }\min\{2^{3j+2 l },  2^{j+2 \tilde{\alpha}_t  M_t} \}, 2^{-2m-j-2l} \} \big)^{1/2}\\
 & \qquad  \times \big( \min\{ \frac{2^{-p-q} }{2^{m}|  x_{\bot}|}, 2^{2k+n +\max\{n,p\}+q  } \} \min\{2^{3j+2l}, 2^{j+2\tilde{\alpha}_t M_t}\} \}\big)^{1/2} \\
 & \lesssim   \| \mathfrak{m}(\cdot, \zeta)\|_{\mathcal{S}^\infty} 2^{m-j  +4\epsilon M_t}       \big(  2^{-m/2 +\min\{n,p\}} \min\{2^{j}, 2^{\tilde{\alpha}_t M_t-l }\} \big)^{1/2} \\
 &\qquad \times  \big(2^{-m+k+\max\{n,p\}-p  } \min\{2^{3j+  2 l }, 2^{j+2\tilde{\alpha}_t M_t}\}\big)^{1/2}  \\
 &  \lesssim \| \mathfrak{m}(\cdot, \zeta)\|_{\mathcal{S}^\infty} 2^{(k+n)/2 +2\epsilon M_t}\min\{2^{j+l}, 2^{\tilde{\alpha}_t M_t}\}  \\
 & \lesssim  \| \mathfrak{m}(\cdot, \zeta)\|_{\mathcal{S}^\infty} \big[ 2^{(1-20\epsilon)M_{t^{\star}}} +   2^{120\epsilon  M_{t^{\star}} }\min\{ 2^{(k+2n)/2+ (\alpha^{\star}+3\iota) M_{t^{\star}}  },2^{(k+4n)/2+ (1+6\iota) M_{t^{\star}}  } \}\big].
 \end{split}
 \ee

 \medskip

\noindent \textbf{Step 4.} \quad  If $m+k\geq -2n+\epsilon M_t$, $m+p\geq -k-n+2\epsilon M_t$  and $|  x_{\bot}|\leq 2^{m+p-10}$. 

\medskip

 Recall  \eqref{2022feb10eqn81}. Similar to the obtained estimate \eqref{aug5eqn44}, by using the same strategy,  we have   
\be\label{2022feb11eqn24}
\begin{split}
  &\big| \widetilde{K}^{m;p,q,2;2}_{k,j ; n,l,r  }(t, x, \zeta ) \big|\\
&\lesssim     \big(\sup_{s\in [0,t]}\|B (s,\cdot)\|_{L^2}\big)   2^{m-j  + \max\{p, l  \}  + 5 \epsilon M_t}   \big( \frac{2^{k+\max\{n,p\} }}{2^{m+p}}   \min\{2^{3j+2  l }, 2^{j+2  \tilde{\alpha}_t M_t}\} \big)^{1/2}\\
&\qquad \times \big(  \min\{   2^{m+2\min\{p,n\} } \min\{2^{3j+2l }, 2^{j+2 \tilde{\alpha}_t  M_t}\},   2^{-2m-j-2l} \} \big)^{1/2} \| \mathfrak{m}(\cdot,   \zeta)\|_{\mathcal{S}^\infty} \\
& \lesssim  \| \mathfrak{m}(\cdot,   \zeta)\|_{\mathcal{S}^\infty}  2^{m-j  + \max\{p, l  \}  + 5 \epsilon M_t}  \big(  {2^{-m-p+k+\max\{n,p\} }}    \min\{2^{3j+2  l }, 2^{j+2  \tilde{\alpha}_t M_t}\} \big)^{1/2} \\
&\qquad  \times \big(  (2^{m+2\min\{p,n\} +3j+2l })^{1/2}(2^{-2m-j-2l})^{1/2} \big)^{1/2}\\
& \lesssim \| \mathfrak{m}(\cdot, \zeta)\|_{\mathcal{S}^\infty}  2^{(k+n)/2 +2\epsilon M_t}\min\{2^{j+l}, 2^{\tilde{\alpha}_t M_t}\} \\
&\lesssim \| \mathfrak{m}(\cdot, \zeta)\|_{\mathcal{S}^\infty} \big[ 2^{(1-20\epsilon)M_{t^{\star}}} +   2^{120\epsilon  M_{t^{\star}} }\min\{ 2^{(k+2n)/2+ (\alpha^{\star}+3\iota) M_{t^{\star}}  },2^{(k+4n)/2+ (1+6\iota) M_{t^{\star}}  } \}\big].
\end{split}
\ee

To sum up, our desired estimate  \eqref{2022feb11eqn21}  holds from the obtained estimates \eqref{2022feb10eqn81},   \eqref{2022feb10eqn83}, \eqref{2022feb10eqn86}  and  \eqref{2022feb11eqn24}.

\end{proof}
 \begin{lemma}\label{smallangleS2part3}
   Let $  (l,r)\in \mathcal{B}_2 $, see  \eqref{sep18eqn50}.    Under the assumption  of Theorem \ref{mainresultsfirstpart} and the assumptions that  $m+k\geq -2 l +4\epsilon M_t$ and $k+2j\geq (2-50\epsilon)  M_{t^{\star}}$, the following estimate holds   if   $n \leq    (1-2\alpha^{\star}-35\epsilon)M_{t^{\star}} $, or    $    j\geq    (1/2+3\iota + 55\epsilon) M_{t^{\star}} $, 
\be\label{2022feb11eqn41}
\begin{split}
 \big| \widetilde{K}^{m;p,q,2;2}_{k,j ; n,l,r  }(t, x, \zeta ) \big| & \lesssim \| \mathfrak{m}(\cdot, \zeta)\|_{\mathcal{S}^\infty}\big[ 2^{  (1-20\epsilon)   M_{t^{\star}}} +  2^{100\epsilon M_{t^{\star}} } \mathbf{1}_{n\geq  -(\alpha^{\star}+3\iota+60\epsilon) M_{t^{\star}}} \\
 &\qquad \times  \min\{2^{(k+2n)/2 }2^{(2{\alpha}^{\star} -1)M_{t^{\star}}}  , 2^{(k+4n)/2 +(1+6\iota)M_{t^{\star}} } \}\big].\\
 \end{split}
\ee
\end{lemma}
\begin{proof}
Recall  \eqref{aug8eqn31}. Again,   due to the cutoff function  $\varphi_{j,n}^2(v, \zeta)$, see \eqref{sep17eqn63} and   \eqref{sep4eqn6}, and the assumption that $m+k\geq -2 l +4\epsilon M_t$, we have  we have $l\geq n+3\epsilon M_t/4> \bar{l}_2$,        $|\tilde{v}-\tilde{\zeta}|\sim 2^l$, and $2^l\sim 2^r.$

 Similar to what we did in \eqref{aug10eqn31}, after   using the decomposition in \eqref{july5eqn1}  for the magnetic field, the following estimate holds  from the estimate of kernels in  \eqref{sep6eqn31}  and the estimates of coefficients in   \eqref{2022feb10eqn1}  and \eqref{2022feb16eqn71}, 
\be\label{2022feb11eqn130}
\big|\widetilde{K}^{m;p,q,2;2}_{k,j ; n,l,r }(t, x, \zeta )\big|\lesssim  \sum_{ (\tilde{m}, \tilde{k}, \tilde{j}, \tilde{l}) \in \mathcal{S}_1(t)\cup  \mathcal{S}_2(t) }  \big|\widetilde{K}^{m;p,q; \tilde{m}}_{k,j ; n,l  ; \tilde{k};\tilde{j}, \tilde{l}}(t, x, \zeta )\big|,
 \ee
 where 
\be\label{2022feb11eqn1}
\begin{split}
\widetilde{K}^{m;p,q; \tilde{m}}_{k,j ; n,l  ; \tilde{k};\tilde{j}, \tilde{l}}(t, x, \zeta ) & := \int_0^t  \int_{\R^3} \int_{\R^3}\int_0^{2\pi} \int_0^{\pi}   \| \mathfrak{m}(\cdot, \zeta)\|_{\mathcal{S}^\infty}  2^{m+3k+2n-j  + \max\{p, l  \}  +3 \epsilon M_t}        \\
&\quad \times   f(s,x-y+(t-s)\omega, v)   \big|B^{\tilde{m}}_{\tilde{k};\tilde{j},  \tilde{l}} (s, x-y+(t-s)\omega)\big| \\
& \quad   \times  \min\big\{ (1+2^k|y\cdot\tilde{\zeta}|)^{- N_0^{ 3}} (1+2^{k+n}|y\times \tilde{\zeta}| )^{- N_0^{ 3} },\\
&\qquad  (1+2^{k+n}|  y_{\bot}| )^{- N_0^{ 3} } (1+2^{k+n } (2^n + |\tilde{\zeta}\times e_3|)^{-1} |y_3|)^{-N_0^3} \big\} \\
    &\quad   \times   \varphi_{m;-10M_t  }(t-s)  \varphi_{[l-1,l+1];\bar{l}_2 }(\tilde{v}+\omega) \psi_{ess}(  \omega,\zeta) \\
    &\quad \times    \varphi_{p;-10 M_t}(\sin \theta) \varphi_{q;-1 0M_t}(\sin\phi)  \sin \theta  d\theta d\phi  dy  d v d s. \\
\end{split}
\ee

Based on  different sizes of parameters, e.g., $m+k,l,n,m+p, k+n, | x_{\bot}|$,  we proceed in four steps as follows.

 \medskip

\noindent \textbf{Step 1.}\quad If $m+k\leq -2n+ \epsilon M_t  $.

 \medskip

Recall  \eqref{2022feb11eqn1}. Similar to the obtained estimate  \eqref{2022feb10eqn80}, we have 
\be\label{2022feb11eqn90}
\begin{split}
\big|\widetilde{K}^{m;p,q; \tilde{m}}_{k,j ; n,l   ; \tilde{k};\tilde{j}, \tilde{l}}(t, x, \zeta )\big|&\lesssim  2^{5\epsilon M_t} \| \mathfrak{m}(\cdot, \zeta)\|_{\mathcal{S}^\infty}  \sup_{s\in [0,t]}\|B^{\tilde{m}}_{\tilde{k};\tilde{j},  \tilde{l}}(s,\cdot)\|_{L^2}\\
&\qquad \times \min\big\{ 2^{k+n}, 2^{m/2 } \min\{2^{2j+2l}, 2^{2\alpha_t M_t} \}   \big\} \\
\end{split}
\ee
Alternatively, similar to the obtained estimate  \eqref{aug4eqn42}, if we put the localized magnetic field in $L^\infty$, we have
 \be\label{2022feb11eqn91}
\begin{split}
 \big|\widetilde{K}^{m;p,q; \tilde{m}}_{k,j ; n,l  ; \tilde{k};\tilde{j}, \tilde{l}}(t, x, \zeta )\big| & \lesssim    \sup_{s\in [0,t]}\|B^{\tilde{m}}_{\tilde{k};\tilde{j},  \tilde{l}}(s,\cdot)\|_{L^\infty} \| \mathfrak{m}(\cdot, \zeta)\|_{\mathcal{S}^\infty}  2^{m-j  + 4\epsilon M_t}\\
 &\quad \times    \min\{2^{-k} 2^{3j+2 l } , 2^{2k+2n -j  }, 2^{-2m  -j-2l}      \}\\
& \lesssim  \sup_{s\in [0,t]}\|B^{\tilde{m}}_{\tilde{k};\tilde{j}, \tilde{l}}(s,\cdot)\|_{L^\infty }2^{3\epsilon M_t }\| \mathfrak{m}(\cdot, \zeta)\|_{\mathcal{S}^\infty} \\
&\quad \times \min\{2^{-k/2}, 2^{k-2j}, 2^{-2j/3}, 2^{4l/3+2n/3+2j/3}\}. 
\end{split}
\ee

Similar to the obtained estimates  \eqref{aug4eqn50}  and  \eqref{aug4eqn51}, after combining the above two estimates \eqref{2022feb11eqn90}   and  \eqref{2022feb11eqn91}  and the using   the estimates  in \eqref{aug4eqn10}  in Proposition \ref{meanLinfest}, we have
\be\label{2022feb11eqn92}
\sum_{  (\tilde{m},\tilde{k},\tilde{j},\tilde{l})\in \mathcal{S}_1(t)   \cup  \mathcal{S}_2(t)}  \big|\widetilde{K}^{m;p,q; \tilde{m}}_{k,j ; n,l  ; \tilde{k};\tilde{j}, \tilde{l}}(t, x, \zeta )\big|  \lesssim  2^{(1-20\epsilon)M_{t^{\star}}}\| \mathfrak{m}(\cdot, \zeta)\|_{\mathcal{S}^\infty}. 
\ee

 \medskip

\noindent \textbf{Step 2.}\quad   If $m+k\geq -2n+\epsilon M_t$ and  $m+p\leq -k-n+2\epsilon M_t$. 

 \medskip

 Recall  \eqref{2022feb11eqn1}. Similar to the obtained estimate  \eqref{2022feb10eqn82}, we have 
\be\label{2022feb11eqn93}
\begin{split}
 \big|\widetilde{K}^{m;p,q; \tilde{m}}_{k,j ; n,l  ; \tilde{k};\tilde{j}, \tilde{l}}(t, x, \zeta )\big| &\lesssim    2^{4\epsilon M_t} \| \mathfrak{m}(\cdot, \zeta)\|_{\mathcal{S}^\infty}  \sup_{s\in [0,t]}\|B^{\tilde{m}}_{\tilde{k};\tilde{j},  \tilde{l}}(s,\cdot)\|_{L^2} \\
 &\qquad \times \min\big\{   2^{-k/2-n} \min\{2^{2j+2l}, 2^{2\alpha_t M_t} \}     ,    2^{k+n } \big\}  . \\
 \end{split}
\ee

Alternatively, similar to the obtained estimate \eqref{july23eqn62}, after putting the localized magnetic field in $L^\infty$, we have
\be\label{2022feb11eqn94}
\begin{split}
 &\big|\widetilde{K}^{m;p,q; \tilde{m}}_{k,j ; n,l  ; \tilde{k};\tilde{j}, \tilde{l}}(t, x, \zeta )\big| \\
 & \lesssim    \sup_{s\in [0,t]}\|B^{\tilde{m}}_{\tilde{k};\tilde{j},  \tilde{l}}(s,\cdot)\|_{L^\infty} \| \mathfrak{m}(\cdot, \zeta)\|_{\mathcal{S}^\infty} 2^{m-j  + 4\epsilon M_t}   \min\{2^{m+3j+2l+2\min\{p,n\}} , 2^{-2m -j-2l}   \}\\
  &  \lesssim \sup_{s\in [0,t]}\|B^{\tilde{m}}_{\tilde{k};\tilde{j}, \tilde{l}}(s,\cdot)\|_{L^\infty}2^{2\epsilon M_t} \min\{2^{m/2  +n},  2^{-2j/3} \}\| \mathfrak{m}(\cdot, \zeta)\|_{\mathcal{S}^\infty}.\\
\end{split}
\ee

Similar to the obtained estimates  \eqref{2022feb11eqn96}  and  \eqref{2022feb11eqn97}, after combining the above two estimates  \eqref{2022feb11eqn93}   and  \eqref{2022feb11eqn94}  and the using   the estimates  in  \eqref{aug4eqn10}  in Proposition \ref{meanLinfest}, we have
\be\label{2022feb11eqn95}
\sum_{  (\tilde{m},\tilde{k},\tilde{j},\tilde{l})\in \mathcal{S}_1(t)   \cup  \mathcal{S}_2(t)}  \big|\widetilde{K}^{m;p,q; \tilde{m}}_{k,j ; n,l  ; \tilde{k};\tilde{j}, \tilde{l}}(t, x, \zeta )\big| \lesssim  2^{(1-20\epsilon)M_{t^{\star}}}\| \mathfrak{m}(\cdot, \zeta)\|_{\mathcal{S}^\infty}. 
\ee

 \medskip

\noindent \textbf{Step 3.}\quad   If $m+k\geq -2n+\epsilon M_t$, $m+p\geq -k-n+2\epsilon M_t$  and $|  x_{\bot}|\geq 2^{m+p-10}$.

 \medskip

 Recall  \eqref{2022feb11eqn1}. Similar to the obtained estimate    \eqref{2022feb10eqn86}, we have 
 \be\label{2022feb11eqn110}
 \begin{split}
 &\big|\widetilde{K}^{m;p,q; \tilde{m}}_{k,j ; n,l  ; \tilde{k};\tilde{j}, \tilde{l}}(t, x, \zeta )\big| \\
 & \lesssim \| \mathfrak{m}(\cdot, \zeta)\|_{\mathcal{S}^\infty}         \big(   \sup_{s\in [0,t]}\|B^{\tilde{m}}_{\tilde{k};\tilde{j},  \tilde{l}}(s,\cdot)\|_{L^2} \big)  2^{m-j   +4 \epsilon M_t}  \big(   \frac{2^{-p-q} }{2^{m}|  x_{\bot}|}  \min\{2^{3j+2 l}, 2^{j+2\tilde{\alpha}_t M_t}\} \}\big)^{1/2} \\
&\qquad \times  \big( \min\{ 2^{m+2 p+q  }\min\{2^{3j+2 l },  2^{j+2 \tilde{\alpha}_t  M_t} \}  \} \big)^{1/2}  \\
  & \lesssim    2^{m/2 +5\epsilon M_t}\min\{2^{2j+2l}, 2^{2\tilde{\alpha}_t M_t}\}   \| \mathfrak{m}(\cdot, \zeta)\|_{\mathcal{S}^\infty}  \sup_{s\in [0,t]}\|B^{\tilde{m}}_{\tilde{k};\tilde{j},  \tilde{l}}(s,\cdot)\|_{L^2}. \\
  \end{split}
\ee

Alternatively, similar to the obtained estimate  \eqref{aug4eqn64}, if we put the localized magnetic field in $L^\infty$, we have
  \be\label{2022feb11eqn111}
\begin{split}
&\big|\widetilde{K}^{m;p,q; \tilde{m}}_{k,j ; n,l  ; \tilde{k};\tilde{j}, \tilde{l}}(t, x, \zeta )\big| \\
& \lesssim    \sup_{s\in [0,t]}\|B^{\tilde{m}}_{\tilde{k};\tilde{j},  \tilde{l}}(s,\cdot)\|_{L^\infty}  2^{m-j  + 5\epsilon M_t}  
  \| \mathfrak{m}(\cdot, \zeta)\|_{\mathcal{S}^\infty} \\
  &\qquad \times  \min\big\{ 2^{m+2\min\{n,p\}  }\min\{2^{3j+2 l },  2^{j+2 \tilde{\alpha}_t  M_t} \} , 2^{-2m -j-2l}   \big\} \\
& \lesssim  \sup_{s\in [0,t]}\|B^{\tilde{m}}_{\tilde{k};\tilde{j}, -\tilde{l}}(s,\cdot)\|_{L^\infty}2^{6\epsilon M_t   } \min\{2^{m/2+\min\{n,p\} }, 2^{-2j/3 }, 2^{-m/2+n/3-4( j+ l)/3}\}\| \mathfrak{m}(\cdot, \zeta)\|_{\mathcal{S}^\infty}. \end{split}
\ee

Similar to the obtained estimates  \eqref{2022feb11eqn105}  and  \eqref{aug4eqn74}, after combining the above two estimates \eqref{2022feb11eqn110}   and \eqref{2022feb11eqn111}  and     using   the estimates  in   \eqref{aug4eqn10}  in Proposition \ref{meanLinfest}, we have
\be\label{2022feb11eqn103}
\sum_{  (\tilde{m},\tilde{k},\tilde{j},\tilde{l})\in \mathcal{S}_1(t)   \cup  \mathcal{S}_2(t)}  \big|\widetilde{K}^{m;p,q; \tilde{m}}_{k,j ; n,l  ; \tilde{k};\tilde{j}, \tilde{l}}(t, x, \zeta )\big|  \lesssim  2^{(1-20\epsilon)M_{t^{\star}}}\| \mathfrak{m}(\cdot, \zeta)\|_{\mathcal{S}^\infty}. 
\ee
 
 \medskip

\noindent \textbf{Step 4.}\quad  If $m+k\geq -2n+\epsilon M_t$, $m+p\geq -k-n+2\epsilon M_t$  and $| x_{\bot}|\leq 2^{m+p-10}$.

 \medskip
 
 Note that, for this case, we have $|   x_{\bot}-  y_{\bot} + (t-s) \omega_{\bot}|\sim 2^{m+p}$ for any $y\in B(0, 2^{-k-n+\epsilon M_t/2})\subset \R^3$,  and the volume of support of $\omega$ is bounded by $2^{2\min\{p,n\}}$.

  Based on the possible size of $n$, we proceed in two sub-steps as follows. 

 \medskip

 \textbf{Step 4A.}\quad   If $\min\{p,n\}\leq   -(\alpha^{\star}+3\iota +60\epsilon) M_t$.

 \medskip

  Recall  the estimate   \eqref{2022feb11eqn1}.  Similar to the obtained estimate  \eqref{aug4eqn75}, we have 
\be\label{2022feb11eqn141}
\begin{split}
\big|\widetilde{K}^{m;p,q,2;2}_{k,j ; n,l,r }(t, x, \zeta )\big|  &\lesssim   2^{m-j  + 9\epsilon  M_t}  ( 2^{2\tilde{\alpha}_t M_t } 2^{-(m+p)/2} + 2^{-(m+p)/4}  2^{5M_t/4+\tilde{\alpha}_t M_t/4 } )  \\
& \qquad \times \min\{2^{m+2\min\{n,p\}+3j+2l}, 2^{-2m-j-2l} \}  \| \mathfrak{m}(\cdot, \zeta)\|_{\mathcal{S}^\infty} \\
& \lesssim \| \mathfrak{m}(\cdot, \zeta)\|_{\mathcal{S}^\infty}  2^{ 8\epsilon  M_t}( 2^{2\tilde{\alpha}_t M_t  } 2^{\min\{n,p\}/2}  +     2^{5M_t/4+\tilde{\alpha}_t M_t/4} 2^{3\min\{n,p\}/4}    )\\
& \lesssim 2^{(1-20\epsilon)M_{t^{\star}}}\| \mathfrak{m}(\cdot, \zeta)\|_{\mathcal{S}^\infty}. 
\end{split}
\ee

 \medskip

 \textbf{Step 4B.}\quad   If $\min\{p,n\}\geq   -(\alpha^{\star}+3\iota +60\epsilon) M_t$.

 \medskip

Similar to the obtained estimate    \eqref{2022feb11eqn24}, we have 
\be\label{2022feb11eqn121}
\begin{split}
 \big|\widetilde{K}^{m;p,q; \tilde{m}}_{k,j ; n,l  ; \tilde{k};\tilde{j}, \tilde{l}}(t, x, \zeta )\big| & \lesssim      \big( \sup_{s\in [0,t]}\|B^{\tilde{m}}_{\tilde{k};\tilde{j},  \tilde{l}}(s,\cdot)\|_{L^2}\big)  2^{m-j  + \max\{p, l  \}  + 5 \epsilon M_t}   \| \mathfrak{m}(\cdot,   \zeta)\|_{\mathcal{S}^\infty} \\
&\quad \times   \big(     2^{-2m-j-2l} \big)^{1/2}  \big( \frac{2^{k+\max\{n,p\} }}{2^{m+p}}    2^{3j+2 l } \big)^{1/2} \\
  &  \lesssim 2^{-m/2+k/2 + \max\{n,p\}/2-p/2}\| \mathfrak{m}(\cdot, \zeta)\|_{\mathcal{S}^\infty}  \big( \sup_{s\in [0,t]}\|B^{\tilde{m}}_{\tilde{k};\tilde{j},  \tilde{l}}(s,\cdot)\|_{L^2}\big).\\
\end{split}
\ee

Additionally,  note that the obtained estimate  \eqref{2022feb11eqn111}  is also valid if we put the localized magnetic field in $L^\infty$.  Similar to the obtained estimates  \eqref{aug4eqn76}  and  \eqref{aug4eqn77}, after combining the above two estimates  \eqref{2022feb11eqn121}   and  \eqref{2022feb11eqn111}  and     using   the estimates in \eqref{aug4eqn10}  in Proposition \ref{meanLinfest}, we have
\be\label{2022feb11eqn123}
\begin{split}
&\sum_{(\tilde{m},\tilde{k},\tilde{j},\tilde{l})\in \mathcal{S}_1(t) \cup \mathcal{S}_2(t) }    \big|\widetilde{K}^{m;p,q; \tilde{m}}_{k,j ; n,l  ; \tilde{k};\tilde{j}, \tilde{l}}(t, x, \zeta )\big| \\& \lesssim  \| \mathfrak{m}(\cdot, \zeta)\|_{\mathcal{S}^\infty} \big[ 2^{(1-20\epsilon)M_{t^{\star}} } +  2^{100\epsilon M_{t^{\star}} } \mathbf{1}_{n\geq  -(\alpha^{\star}+3\iota+60\epsilon) M_t} \\
 &\qquad \times  \min\{2^{(k+2n)/2 }2^{(2{\alpha}^{\star} -1)M_{t^{\star}}}  , 2^{(k+4n)/2 +(1+6\iota)M_t} \}\big] . 
\end{split}
\ee
Recall  \eqref{2022feb11eqn130}. The desired estimate  \eqref{2022feb11eqn41}  holds from the above estimate, and the obtained estimates  \eqref{2022feb11eqn92}, \eqref{2022feb11eqn95}, \eqref{2022feb11eqn103}, and \eqref{2022feb11eqn141}. 
\end{proof}

\begin{lemma}\label{smallangleS2part6}
   Let $ i\in\{3,4\}, (l,r)\in \mathcal{B}_i $, see  \eqref{sep18eqn50}.     Under the assumption of Theorem \ref{mainresultsfirstpart}, the following estimate holds if  $m+k\geq -2 l +4\epsilon M_t$, 
\be\label{2022feb10eqn71}
\begin{split}
 \big| \widetilde{K}^{m;p,q,i;2}_{k,j ; n,l,r }(t, x, \zeta )\big| 
 \lesssim   \| \mathfrak{m}(\cdot, \zeta )\|_{\mathcal{S}^\infty}   &\big[2^{(1-20\epsilon) M_{t^\star} } + 2^{50\epsilon M_{t^\star} } \mathbf{1}_{n\geq   -((1+3\iota)/2  + 30\epsilon)M_{t^{\star}} }\\
 &\qquad \times  \min\{2^{(k+2n)/2 +  {\alpha}^{\star} M_{t^\star}}, 2^{(k+4n)/2 +(1+3\iota)M_{t^\star}}\}  \big].
 \end{split}
\ee
\end{lemma}
\begin{proof}

Due to the cutoff functions $\varphi_{j,n}^i(v, \zeta), i \in\{3,4\}$, see  \eqref{sep4eqn6}, the support of the $v, \omega, \xi$,    implies that $|\tilde{v}\times \omega|\lesssim  2^{n+\epsilon M_t} +2^{l-\epsilon M_t}, |\tilde{v}\times\omega |\sim 2^l$.  If $l> -j$,   then it follows that $l\leq n +\epsilon M_t$, since  in this case,   $ |\omega \times\tilde{\xi}|\sim 2^l$.   This further implies that  $l\in [-j, \max\{-j, n+\epsilon M_t\}]\cap \Z. $ 

 Moreover,  note that, besides the rough estimate of the gradient of coefficients obtained in \eqref{2022feb10eqn61}, the following improved estimate   holds for the $\p_{v_3}$ derivative of coefficients, 
\be\label{sep21eqn8}
\begin{split}
& \big|\p_{v_3} \big( \frac{  m_{E}(v, \omega )  + \hat{ \zeta}\times  m_{B}(v, \omega )}{1+\hat{v}\cdot \omega}    \varphi_{j,n}^{i; r}(v, \zeta)  \varphi_{l; r}(\tilde{v}+\omega )     \big) \big| \varphi_{p;-1 0M_t}(  \omega_{\bot}) \\
&\lesssim 2^{\epsilon M_t}\big( 2^{-j-2l+n+\max\{l,p\}} + 2^{-j-l+ \max\{p,n+\epsilon M_t\}}\big)\\
 &\lesssim 2^{-j-l-\min\{l,n\}+n + \max\{l,p\}+2\epsilon M_t}.\\
\end{split}
\ee

Recall  \eqref{aug8eqn31}  and   \eqref{july1eqn13}. 
From the   estimate of kernels in  \eqref{sep6eqn31}, the estimate of coefficients in  \eqref{2022feb10eqn61}  and   \eqref{sep21eqn8}, and the Cauchy-Schwarz inequality, we have
\be\label{sep21eqn91}
\begin{split}
  \big| \widetilde{K}^{m;p,q,i;2}_{k,j ; n,l,r }(t, x, \zeta )\big| &\lesssim     2^{m-j -\min\{l,n\}+n + \max\{l,p\}+2\epsilon M_t}  \| \mathfrak{m}(\cdot, \zeta )\|_{\mathcal{S}^\infty} \\
  &\quad \times  \big( H^{m;p,i;B }_{k,j ;l,n}(t, x,  \zeta  )  \big)^{1/2}\big( H^{m;p,i;f }_{k,j ;l,n}(t, x, \zeta )  \big)^{1/2},
  \end{split}
\ee
where
 \be\label{2022feb6eqn32}
 \begin{split}
H^{m;p,i;B }_{k,j ;l,n}(t, x, \zeta  )&:=\int_0^t \int_{\R^3} \int_{\R^3} \int_{\mathbb{S}^2}   2^{3k+2n}     |B(s, x-y+(t-s)\omega)|^2 \\
& \quad  \times   \min\big\{  (1+2^{k+n}|  y_{\bot}| )^{- N_0^{ 3} } (1+2^{k+n- \max\{-j, p, n+\epsilon M_t\} } |y_3|)^{-N_0^3} , \\
&\qquad (1+2^k|y\cdot\tilde{\zeta}|)^{- N_0^{ 3}} (1+2^{k+n}|y\times \tilde{\zeta}| )^{- N_0^{ 3} } \big\}\varphi_{m;-10M_t }(t-s) \\
&\quad \times    \varphi_{l;-j}(\tilde{v}+\omega)    \varphi_{j,n}^i(v, \zeta )  \varphi_{p;-10M_t}(\sin \theta)    d\omega d y d v d s,\\
H^{m;p,i;f }_{k,j ;l,n}(t, x,  \zeta  )&:=\int_0^t \int_{\R^3} \int_{\R^3} \int_{\mathbb{S}^2}   2^{3k+2n} \varphi_{l;-j}(\tilde{v}+\omega)  f(s, x-y+(t-s)\omega, v)  \\
&\quad \times  (1+2^k|y\cdot\tilde{\zeta}|)^{-N_0^3} (1+2^{k+n}(|y\cdot \tilde{\zeta }_1|+|y\cdot \tilde{\zeta }_2|) )^{-N_0^3} \\
& \quad  \times        \varphi_{j,n}^i (v,\zeta)     \varphi_{m;-10M_t }(t-s)\varphi_{p;-10M_t}(\sin \theta)    d\omega d y d v d s,
\end{split}
\ee
 where $\{\tilde{\zeta}, \tilde{\zeta}_1, \tilde{\zeta}_2\}$ is a fixed orthonormal frame. 

By using a similar argument as in the estimate of  $ J^{m;p }_{k,j ;l,n}(t, x, \zeta)$ in   \eqref{sep21eqn71},    after localizing the frequency around $v$, using the stationary phase type argument for the error type term, and using similar estimates as in  \eqref{sep21eqn88}--\eqref{sep26eqn23}  for the essential type term, we have
\be\label{sep21eqn90}
\begin{split}
&\big| H^{m;p,i;B }_{k,j ;l,n}(t, x, \zeta )\big| \\
 &\lesssim 
\int_0^t \int_{\R^3} \int_{\R^3} \int_{\mathbb{S}^2}    2^{3k+\epsilon^2 (k+ M_t)+2\min\{l,n\}+3\epsilon M_t} |B(s, x-y+(t-s)\omega)|^2 \\
&\qquad \times    (1+2^{k-\epsilon M_t}|y\cdot \tilde{v}|)^{-N_0^3} (1+2^{k+\min\{l,n\}-\epsilon M_t}(|y\times \tilde{v}|  )^{-N_0^3} \\
&\qquad \times   \varphi_{l;-j}(\tilde{v}+\omega)    \varphi_{j,n}^i(v, \zeta)   \varphi_{m;-10M_t }(t-s)\varphi_{p;-10M_t}(\sin \theta)      d\omega d y d v d s+2^{-100M_t }. \\
\end{split}
\ee
 
Based on the possible sizes of $m+p, k+l$, and $|  x_{\bot}|$ and the relative sizes of $p$ and $l$,  we proceed in three steps as follows.

 \medskip

\noindent \textbf{Step 1.}\quad   If $m+p\leq -k- \min\{l,n\}  + \epsilon M_t. $

\medskip

Note that, for this case, we have $p+\min\{l,n\}\leq 2l.$  From the obtained estimates \eqref{sep21eqn91}  and  \eqref{sep21eqn90},  the estimate  \eqref{march18eqn31} in Lemma \ref{conservationlawlemma}, the conservation law  \eqref{conservationlaw}, and the volume of support of $ \omega$ and $v $, we have
\be\label{aug19eqn31}
\begin{split}
 \big| \widetilde{K}^{m;p,q,i;2}_{k,j ; n,l,r }(t, x, \zeta )\big|  
&  \lesssim \sup_{s\in [0,t]}\|B  (s,\cdot)\|_{L^2}  2^{m -j- \min\{l,n\} +n+\max\{p, l\}  + 8\epsilon M_t} \\
&\qquad \times  \big(\min\{2^{m +3j+2\min\{l,n\} +2 p   },  2^{-2m  -j-2l}   \} \big)^{1/2 }   \\
&\qquad \times \big(  2^{m+ 3k +\epsilon (k+ M_t)+2\min\{l,n\}} 2^{2p} 2^{3j+ 2\min\{l,n\}} \big)^{1/2}   \| \mathfrak{m}(\cdot, \zeta ) \|_{\mathcal{S}^\infty} \\
& \lesssim     2^{n+\max\{p,l\}+10\epsilon M_t+\epsilon^2 k} \min\big\{  2^{-k/2+2j   },  2^{k } \big\}  \| \mathfrak{m}(\cdot, \zeta)\|_{\mathcal{S}^\infty}\\
&  \lesssim   2^{n+\max\{p,l\}+12\epsilon M_t}  \min\{ 2^{4j/3  }, 2^{k/2 +2j/3} \}\| \mathfrak{m}(\cdot, \zeta )\|_{\mathcal{S}^\infty} \\
&\lesssim \| \mathfrak{m}(\cdot, \zeta)\|_{\mathcal{S}^\infty} \big[2^{(1-20\epsilon)M_{t^{\star}}} + 2^{15\epsilon M_t} 2^{(k+4n)/2+2 M_{t^{\star}}/3} \mathbf{1}_{n\geq- (1/3+40\epsilon ) M_{t^{\star}}  }\big].\\
\end{split}
\ee

  \medskip

\noindent \textbf{Step 2.}\quad 
If $m+p\geq -k-\min\{l,n\}  + \epsilon M_t, $ and $|  x_{\bot}|\geq 2^{m+p-10}$. 

\medskip

Recall  \eqref{aug8eqn31}.  Similar to what we did in  \eqref{march18eqn66},  after  changing coordinates $(\theta, \phi) \longrightarrow (z(y),w(y))$,
  where
  \be
  \begin{split}
   z(y)&:=\sqrt{\big(| x_{\bot}-  y_{\bot}|+(t-s)\sin\theta\cos\phi\big)^2 + (t-s)^2(\sin \theta \sin \phi)^2 } \\
     w(y)&:= x_3 -y_3+(t-s)\cos\theta, 
     \end{split}
  \ee
 and   using  the estimate   \eqref{sep21eqn91},    the   estimate of kernels in  \eqref{sep6eqn31},   the estimate of coefficients in \eqref{2022feb10eqn61} and   \eqref{sep21eqn8}, the volume of support of $v, \omega$, and the estimate  \eqref{nov24eqn41}  if $|  v_{\bot} |\geq 2^{(\alpha_t+\epsilon)M_t}$,   
we have
\be\label{2021dec9eqn41}
\begin{split}
 &\big| \widetilde{K}^{m;p,q,i;2}_{k,j ; n,l,r }(t, x, \zeta )\big| \\
 & \lesssim \sup_{s\in [0,t]}\|B  (s,\cdot)\|_{L^2} \| \mathfrak{m}(\cdot, \zeta )\|_{\mathcal{S}^\infty}   2^{m-j-  \min\{l,n\}+n+\max\{l,p\} + 8\epsilon M_t}  \\
& \qquad  \times ( 2^{m+2p+q} \min\{ 2^{3j+2l}, 2^{j+2\tilde{\alpha}_t M_t }\})^{1/2}( 2^{-2m-2p-q} \min\{ 2^{3j+2l}, 2^{j+2\tilde{\alpha}_t M_t }\})^{1/2} \\
& \lesssim  2^{m/2+\max\{l,n\}+\max\{l,p\} + 9\epsilon M_t} \min\{2^{2j+l}, 2^{2\tilde{\alpha}_t M_t-l}\}\| \mathfrak{m}(\cdot, \zeta )\|_{\mathcal{S}^\infty} \\
&\lesssim \min\{2^{  n +2\tilde{\alpha}_t M_t + 10\epsilon M_t}, 2^{ l+\max\{l,n\}+2j + 5\epsilon M_t}\}\| \mathfrak{m}(\cdot, \zeta )\|_{\mathcal{S}^\infty}.
\end{split}
\ee
In the above estimate, we used the fact that $|  v_{\bot}| \sim 2^{p+j}$ if $p\geq l + \epsilon M_t$.

Additionally, from the estimate  \eqref{sep21eqn91}, the estimate  \eqref{sep21eqn90},  the estimate   \eqref{march18eqn31}  in Lemma \ref{conservationlawlemma}, we   have 
\[ 
\begin{split}
 \big| \widetilde{K}^{m;p,q,i;2}_{k,j ; n,l,r }(t, x, \zeta )\big|  
  & \lesssim 2^{m-j-  \min\{l,n\}+n +\max\{l,p\} + 8\epsilon M_t}      (2^{-2m-j-2l})^{1/2}\| \mathfrak{m}(\cdot, \zeta )\|_{\mathcal{S}^\infty} \\
  &\quad \times  \big( \min\{2^{m+3k+\epsilon^2 k+2\min\{l,n\}+2p+q+3j+2l}, 2^{-2m-2p-q + 3j+2l}\}\big)^{1/2} \\
 & \lesssim  2^{-3j/2-  \min\{l,n\}/2  +n+\max\{l,p\} +8\epsilon M_t} (2^{-m/2 +(3+\epsilon^2 )k/2 +   3j})^{1/2} \| \mathfrak{m}(\cdot, \zeta )\|_{\mathcal{S}^\infty}\\
 & \lesssim 2^{-m/4-l/2+ (3+\epsilon )k/4+n/2+\max\{l,n\}/2  +8\epsilon M_t} \| \mathfrak{m}(\cdot, \zeta )\|_{\mathcal{S}^\infty}  \\
&\lesssim 2^{(4+\epsilon )k/4+ n/2+\max\{l,n\}/2   +8\epsilon M_t}\| \mathfrak{m}(\cdot, \zeta )\|_{\mathcal{S}^\infty},
\end{split}
\]
where  we used the fact that $m+2 l  \geq - k. $ 

After combining the above estimate and the obtained estimates \eqref{2021dec9eqn41}, we have
\be\label{2022feb6eqn40}
\begin{split}
  \big| \widetilde{K}^{m;p,q,i;2}_{k,j ; n,l,r }(t, x, \zeta )\big| &\lesssim\| \mathfrak{m}(\cdot, \zeta )\|_{\mathcal{S}^\infty}\big[2^{(1-20\epsilon) M_{t^\star} }+ 2^{50\epsilon M_{t^\star} } \mathbf{1}_{n\geq   -(1/3+2\iota +30\epsilon )M_{t^{\star}} }  \\ 
&  \quad \times  \min\{2^{(k+2n)/2 +  {\alpha}^{\star} M_{t^\star}}, 2^{(k+4n)/2 +(1+3\iota)M_{t^\star}}\}\big]. 
\end{split}
\ee

  \medskip

\noindent \textbf{Step 3.}\quad   If $m+p\geq -k-\min\{l,n\}  + \epsilon M_t, $  and $| x_{\bot}|\leq 2^{m+p-10}$.

  \medskip

For this case, we have $|  x_{\bot}-  y_{\bot} +(t-s)  \omega_{\bot}|\sim 2^{m+p}$  if $|   y_{\bot}|\leq 2^{-k-n+\epsilon M_t/2}$.  Based on the possible size of $n$, we divide further into two sub-steps as follows.

  \medskip

 \textbf{Step 3A.}\quad   If $n\leq  - (1/2+3\iota/2 + 30\epsilon)M_{t^{\star}}.$ 

\medskip

 Similar to the obtained estimate  \eqref{aug4eqn75},  from  the estimate of coefficients in  \eqref{2022feb10eqn61}  and   \eqref{sep21eqn8},  the pointwise  estimate of magnetic field  \eqref{sep21eqn31} in Lemma \ref{bulkroughpoi},  the estimate  \eqref{march18eqn31}  in Lemma \ref{conservationlawlemma},  and the volume of support of $ \omega$ and $ v$, we have 
\be\label{aug25eqn51}
\begin{split}
   \big| \widetilde{K}^{m;p,q,i;2}_{k,j ; n,l,r }(t, x, \zeta )\big| & \lesssim   2^{m-j-\min\{l,n\}+n+\max\{l,p\}+8\epsilon M_t}  \| \mathfrak{m}(\cdot, \zeta)\|_{\mathcal{S}^\infty} \\
   &\quad \times   \min\{2^{ m+ 2\min\{ n ,p\}+3j+2l}, 2^{-2m-j-2l} \}\\
&\quad \times \big(     2^{-(m+p)/4}  2^{5M_t/4+\tilde{\alpha}_t M_t/4 } +2^{2 \tilde{\alpha}_t M_t}2^{-(m+p)/2} \big) \\
&\lesssim 2^{ \max\{l,n\}+\max\{l,p\}+ 8 \epsilon M_t}   \| \mathfrak{m}(\cdot, \zeta)\|_{\mathcal{S}^\infty} \big(2^{2m+     2\min\{ n ,p\}+2j+ l}\big)^{3/4}\\
&\quad   \times    \big( 2^{-m-2j-3l}\big)^{1/4}   (  2^{-(m+p)/4}  2^{5M_t/4+\tilde{\alpha}_t M_t/4 } +2^{2\tilde{\alpha}_t  M_t}2^{-(m+p)/2} \big)\\
&\lesssim  2^{  \max\{l,n\}+   \max\{l,p\}+ j+3 \min\{ n ,p\}/2+5m/4+8\epsilon M_t} \\
 &\quad \times (   2^{-(m+p)/4}  2^{5M_t/4+ \tilde{\alpha}_t M_t/4 } +2^{2 \tilde{\alpha}_t  M_t}2^{-(m+p)/2} \big)   \| \mathfrak{m}(\cdot, \zeta)\|_{\mathcal{S}^\infty} \\
 & \lesssim   2^{(1-20\epsilon)M_{t^{\star}}} \| \mathfrak{m}(\cdot, \zeta)\|_{\mathcal{S}^\infty}. 
 \end{split}
\ee
In the above estimate, we used the fact that $|  v_{\bot}| \sim 2^{p+j}$ if $p\geq l + \epsilon M_t$ and the estimate  \eqref{nov24eqn41}  if $| v_{\bot}|\geq 2^{(\alpha_t +\epsilon)M_t}. $

  \medskip

 \textbf{Step 3B.}\quad   If  $n\geq  - (1/2+3\iota/2 + 30\epsilon)M_{t^{\star}}.$  

  \medskip

From the obtained estimates   \eqref{sep21eqn91}  and  \eqref{2022feb6eqn32},   the cylindrical symmetry of the electromagnetic field, the estimate  \eqref{march18eqn31}  in Lemma \ref{conservationlawlemma}, the volume of support of $\omega,v$, the Jacobian of changing coordinates $(y_1, y_2, \theta)\longrightarrow x- y + (t-s)\omega$,  and the estimate  \eqref{nov24eqn41}  if $|  v_{\bot}|\geq 2^{(\alpha_t+\epsilon)M_t}$,  we have 
\be\label{aug25eqn52}
\begin{split}
 &\big| \widetilde{K}^{m;p,q,i;2}_{k,j ; n,l,r }(t, x, \zeta )\big| \\
&\lesssim \| \mathfrak{m}(\cdot, \zeta)\|_{\mathcal{S}^\infty} 2^{m-j-\min\{l,n\}+n+\max\{l,p\}+8\epsilon M_t} \\
&\quad \times   \big( \min\{   2^{m+2\min\{p,n\}+3j+2l},2^{m+2\min\{l,p\}+ j+2\tilde{\alpha}_t M_t} ,   2^{-2m-j-2l} \} \big)^{1/2}\\
& \quad \times \big(  2^{3k+ 2n} 2^{-k-n-k-n+\max\{-j, p, n+\epsilon M_t\} }  {2^{-(m+p)}} \min\{2^{j+2\tilde{\alpha}_t M_t},2^{3j+2\min\{l,n\}} \} \big)^{1/2} \\
&   \lesssim    \| \mathfrak{m}(\cdot, \zeta)\|_{\mathcal{S}^\infty}   2^{m/4+  k/2 + n+\max\{-j, p, n+\epsilon M_t\}/2-p/2 + \max\{l,p\} + 8\epsilon M_t }   \\
& \quad \times  \big(2^{\min\{n,p\}+  j}\mathbf{1}_{l=-j, l> n+\epsilon M_t} +2^{\min\{n,p\}/2} \min\{2^{j}, 2^{\tilde{\alpha}_t M_t -l }\} \mathbf{1}_{l\leq n+\epsilon M_t} \big)\\
& \lesssim 2^{(k+3n)/2+   {\alpha}^{\star} M_{t^{\star}} + 10\epsilon   M_{t^{\star}} } \| \mathfrak{m}(\cdot, \zeta)\|_{\mathcal{S}^\infty}\\
& \lesssim  2^{30\epsilon  M_{t^{\star}} } \min\{ 2^{(k+3n)/2+   {\alpha}^{\star} M_{t^{\star}}  }  , 2^{(k+4n)/2+  M_{t^{\star}}} \}  \| \mathfrak{m}(\cdot, \zeta)\|_{\mathcal{S}^\infty}.\\
\end{split}
\ee 
In the above estimate,  we used the fact that $|   v_{\bot}| \sim 2^{p+j}$ if $p\geq l + \epsilon M_t$ and the estimate  \eqref{nov24eqn41}  if $|   v_{\bot}|\geq 2^{(\alpha_t +\epsilon)M_t}. $ 

To sum up, our desired estimate  \eqref{2022feb10eqn71}  holds after combining the obtained estimates  \eqref{aug19eqn31},  \eqref{2022feb6eqn40},  \eqref{aug25eqn51}, and  \eqref{aug25eqn52}. 

\end{proof}

To sum up, the main results of this  section can be summarized in the following Proposition.

\begin{proposition}\label{finalestsmallprop}
   Under the assumption of Theorem \ref{mainresultsfirstpart}, we have 
  \be\label{2022feb12eqn21}
  \begin{split}
  & \big| T_{k,j;n}^{\mu,2}( \mathfrak{m}, E)(t,x, \zeta) +  \hat{\zeta}\times T_{k,j;n}^{\mu,2}( \mathfrak{m}, B)(t,x, \zeta)\big| \\
 & \lesssim  \| \mathfrak{m}(\cdot, \zeta)\|_{\mathcal{S}^\infty}  \big[   2^{(1-19\epsilon)M_{t^{\star}} } +  2^{128\epsilon M_{t^{\star}} } \mathbf{1}_{n\geq  -(\alpha^{\star}+3\iota+50\epsilon) M_{t^{\star}}}\\
 &\quad \times \min\{ 2^{(k+2n)/2+ (\alpha^{\star}+3\iota) M_{t^{\star}} }  , 2^{(k+4n)/2 +(1+6\iota)M_{t^{\star}}} \}   \big],\\
\end{split}
 \ee
  \be\label{2022feb12eqn22}
  \begin{split}
  &\sum_{i=3,4}\big| T_{k,j;n}^{\mu,i}( \mathfrak{m}, E)(t,x, \zeta) +  \hat{\zeta}\times T_{k,j;n}^{\mu,i}( \mathfrak{m}, B)(t,x, \zeta)\big|   \\
  & \lesssim    \| \mathfrak{m}(\cdot, \zeta )\|_{\mathcal{S}^\infty}   \big[2^{(1-19\epsilon) M_{t^\star} } + 2^{128\epsilon M_{t^\star} } \mathbf{1}_{n\geq   -(1/2+3\iota/2 + 40\epsilon)M_{t^{\star}} } \\
  &\quad \times \min\{2^{(k+2n)/2 +  {\alpha}^{\star} M_{t^\star}}, 2^{(k+4n)/2 +(7/6+5\iota/2)M_{t^\star}}\} \big]. 
  \end{split}
\ee
\end{proposition}
\begin{proof}
Recall the G-S type  decomposition in   \eqref{sep18eqn50}  and the decomposition of the S-part in  \eqref{sep7eqn61}. 

The desired estimate  \eqref{2022feb12eqn21}  holds from the estimate \eqref{sep19eqn71}  in Lemma \ref{smallfrepartII}, the estimate  \eqref{sep6eqn49}  in Lemma \ref{secondTpartlarge}, the estimate  \eqref{2022feb10eqn35}  in Lemma  \ref{smallangleS1part1},   the estimate  \eqref{2022feb11eqn21}  in Lemma \ref{smallangleS2part1}, and the estimate  \eqref{2022feb11eqn41}  in Lemma \ref{smallangleS2part3}. 

The desired estimate  \eqref{2022feb12eqn22}  holds from the estimate  \eqref{2022feb12eqn11}  in Lemma \ref{smallfrepartII}, the estimate \eqref{sep6eqn49}  in Lemma \ref{secondTpartlarge}, the estimate  \eqref{2022feb10eqn66}  in Lemma  \ref{smallangleS1part5},    and the estimate  \eqref{2022feb10eqn71}  in Lemma \ref{smallangleS2part6}.
\end{proof}

 \subsection{Pointwise estimates of   the projection  of the magnetic field $\mathbf{P}_3(B)$ }

As a byproduct of the previous sections, we demonstrate in this section that the projection of the magnetic field onto the 
$z$-axis, 
 $\mathbf{P}_3(B(t,x))$, is better than  $B(t,x)$. Furthermore, a portion of this projection exhibits similar estimates to the localized acceleration force, which possesses the double null structure. This is due to the fact that the associated symbol   $\mathbf{P}_{3}(\hat{v}\times \xi)$
 is better than  the symbol $  (\hat{v}\times \xi)$, which is associated with $B(t,x)$, when  $\tilde{v}$ and $\tilde{\xi}$ are close to the north or south pole.  These estimates are crucial to analyzing the elliptic part of the iterative smoothing scheme in Part I \cite{PartI}[Section 6].

 To clarify the above observation, we need to decompose $\mathbf{P}_3(B)$
  into components that make it evident that the symbol $\mathbf{P}_{3}(\hat{v}\times \xi)$
  serves as the double null structure for certain pieces. It is important to note that this symbol does not function as the double null structure for all components, as there is no significant improvement when 
 $\tilde{v}$ or $\tilde{\xi}$ are not close to the north or south pole.

Recall  \eqref{sep5eqn10}.  Let  $\zeta\in \R^n/\{0\}, k\in \Z_+, n, \tilde{n}\in [-M_t, 2]\cap \Z$ be fixed  such that $|  \zeta_{\bot}|/|\zeta|\leq 2^{\tilde{n}-10}$, we define the  localized version of $\mathbf{P}_3(B)$ as follows, 
\be\label{nov6eqn41}
\begin{split}
{}_{}^zT_{k,n}^{\mu}(B)(t,x,\zeta) & : =\int_0^t \int_{\R^3}\int_{\R^3} e^{i x\cdot \xi +i \mu(t-s)|\xi|} \mathbf{P}_3(\hat{v}\times \xi )|\xi|^{-1} \\
&\quad \times  \widehat{f}(s, \xi, v) \varphi_k(\xi)  \varphi_{n;-M_t}( \tilde{\xi}+ \mu \tilde{\zeta})  d \xi d v d s.\\
\end{split}
\ee

Based on the possible size of the angle between $v$ and $\pm \zeta$, we divide ${}_{}^zT_{k,n}^{\mu}(B)(t,x,\zeta)$ into two parts as follows,  
\be\label{nov6eqn21}
\begin{split}
 {}_{}^zT_{k,n}^{\mu}(B)(t,x,\zeta)&: = {}_{}^zT_{k,n}^{\mu;1}(B)(t,x,\zeta) + {}_{}^zT_{k,n}^{\mu;2}(B)(t,x,\zeta), \\ 
  {}_{}^zT_{k,n}^{\mu;i}(B)(t,x,\zeta)&: =\int_0^t \int_{\R^3} \int_{\R^3}e^{i x\cdot \xi  + i \mu(t-s)|\xi|} \mathbf{P}_3(\hat{v}\times \xi )|\xi|^{-1} \\
  &\quad \times  \widehat{f}(s, \xi, v) \varphi_k(\xi)   \varphi_{n;-M_t}( \tilde{\xi}+ \mu \tilde{\zeta}) \varphi_{\tilde{n}}^i(v, \zeta)  d \xi d v d s, \\ 
  \varphi_{\tilde{n}}^1(v, \zeta)&: =\psi_{\geq \tilde{n}+5}(\tilde{v}- \tilde{\zeta} ), \quad \varphi_{\tilde{n}}^2(v, \zeta): =\psi_{ < \tilde{n}+5}( \tilde{v}- \tilde{\zeta}). 
 \end{split}
 \ee

For the second part $ {}_{}^zT_{k,n}^{\mu;2}(B)(t,x,\zeta)$,  similar  to the obtained decomposition   in \eqref{nov5eqn10}, we have
\be\label{nov6eqn43}
 {}_{}^zT_{k,n}^{\mu;2}(B)(t,x,\zeta) = \sum_{j\in \Z_+, l\in[-j,2]\cap \Z} \sum_{m\in[-10 M_t  , \epsilon M_t]\cap \Z}  {}_{}^zT_{k,j;n,l}^{\mu;m,2}(B)(t,x,\zeta),
\ee
where, as in  \eqref{sep18eqn31}, we have 
 \be\label{nov5eqn23}
\begin{split}
&{}_{}^zT_{k,j;n,l}^{\mu;m,2}(B)(t,x,\zeta)\\
&= \int_{0}^t \int_{\R^3} \int_{\R^3} \int_{\mathbb{S}^2}  \big( (t-s) \mathbf{P}_3\big( \mathfrak{K}^{\mu, B}_{k;n}(y,\omega, v, \zeta)\big) +\mathbf{P}_3\big( \mathfrak{K}^{err;\mu,B}_{k;n}(y, v, \zeta) \big) \big)    \\
&\quad \times  f(s, x-y+(t-s)\omega, v)  \varphi_{\tilde{n}}^2(v, \zeta)  \varphi_j(v)  \varphi_{l;-j}(\tilde{v}+\omega)  \varphi_{m;-10M_t  }(t-s)    d\omega dy d v ds,\\
\end{split}
\ee
where the kernels $\mathfrak{K}^{\mu, B}_{k;n}(y,\omega, v, \zeta)$ and $\mathfrak{K}^{err;\mu, B}_{k;n}(y,\omega, v, \zeta)$ are defined in \eqref{sep19eqn85}.

To sum up, from  \eqref{nov6eqn41},  \eqref{nov6eqn21}, and  \eqref{nov6eqn43},   for any $\zeta\in \R^3$, after doing dyadic decomposition for the size of frequency and  the angle between $\xi$ and $\pm \zeta$, we have the following decomposition,
\be\label{nov6eqn47} 
\begin{split}
\mathbf{P}_3(B(t,x))&= \sum_{\begin{subarray}{c}
k\in \Z_+, n\in [-M_t, 2]\cap \Z \\ 
  \mu\in \{+,-\}, n\geq \tilde{n}\\ 
\end{subarray}}   {}_{}^zT_{k,n}^{\mu}(B)(t,x,\zeta) \\
&\quad +   \sum_{\begin{subarray}{c}
k\in \Z_+, n\in [-M_t, 2]\cap \Z \\ 
  \mu\in \{+,-\}, n\leq \tilde{n}\\ 
\end{subarray}} {}_{}^zT_{k,n}^{\mu;1}(B)(t,x,\zeta) + {}_{}^zT_{k,n}^{\mu;2}(B)(t,x,\zeta), \\ 
{}_{}^zT_{k,n}^{\mu;2}(B)(t,x,\zeta)&= \sum_{(m,k,j,l)\in \mathcal{S}_1(t)\cup \mathcal{S}_2(t)} {}_{}^zT_{k,j;n,l}^{\mu;m,2}(B)(t,x,\zeta), 
\end{split}
\ee
where the index sets $\mathcal{S}_i(t), i\in \{ 1,2\}$, are defined in  \eqref{indexsetL2}.

 In the following Lemma, we prove that  the ${}_{}^zT_{k,n}^{\mu;1}(B)(t,x,\zeta)$  part of $ \mathbf{P}_3(B(t,x))$ has an improved estimate because  the symbol  $\mathbf{P}_{3}(\hat{v}\times \xi)$   plays the role of the double null structure in the localized region. Moreover, since the difference between the null structure and the double null structure can be compensated by the smallness factor $2^n$. For any localized piece $ {}_{}^zT_{k,n}^{\mu}(B)(t,x,\zeta)$ of $B$, in which we have the null structure,  we have a slightly worse estimate (by a factor of $2^{-n}$) than the localized acceleration force. 
 
\begin{lemma}\label{goodpartprojmagn}
Let     $\mu\in\{+,-\},  k, j\in \Z_+, n, \tilde{n}\in [-M_t, 2]\cap \Z, $  $ t^{\star}\in [0, T), \alpha^{\star}:= 2/3+\iota, \iota:=10^{-4},   x\in \R^3, \zeta \in\R^3/\{0\}, t \in [0,  t^{\ast}]$ be fixed     s.t., $\alpha_t M_{t_{ }}\leq  \alpha^{\star} M_{t^{\star}}, M_t\gg 1$, and $| \zeta_{\bot}|/ |\zeta|\leq 2^{\tilde{n}-10}$.   If  $n\leq \tilde{n}$,  then  the following estimates holds, 
\be\label{nov5eqn1}
\begin{split}
\big\|  {}_{}^zT_{k,n}^{\mu;1}(B)(t,x,\zeta) \big\|_{L^\infty_x}  &\lesssim   2^{(1-19\epsilon)M_{t^{\star}} }  + 2^{  130\epsilon M_{t^{\star}} }  \mathbf{1}_{n\geq  -(\alpha^{\star}+3\iota+60\epsilon) M_t    }    \\
&\qquad \times  \min\{ 2^{(k+2n)/2+ (\alpha^{\star}+3\iota) M_t }  , 2^{(k+4n)/2 +(1+6\iota)M_t} \}, \\
\big\|  {}_{}^zT_{k,n}^{\mu }(B)(t,x,\zeta) \big\|_{L^\infty_x}    &  \lesssim   2^{-n}     \big[   2^{(1-19\epsilon)M_{t^{\star}} } +   2^{  130\epsilon M_{t^{\star}} }  \mathbf{1}_{n\geq  -(\alpha^{\star}+3\iota+60\epsilon) M_t} 
\\
 &\qquad\times      \min\{ 2^{(k+2n)/2+ (\alpha^{\star}+3\iota) M_t }  , 2^{(k+4n)/2 +(1+6\iota)M_t} \}    \big]. \\
\end{split}
\ee 
 
\end{lemma}
\begin{proof}
Recall that $|  \zeta_{\bot}|/ |\zeta|\leq 2^{\tilde{n}-10}$. For any $(v,\xi)\in supp \big( \psi_n\big( \tilde{\zeta}\times \tilde{\xi} \big) \psi_{l}( \tilde{v}\times \tilde{\zeta}) \big)$, where $l\geq \tilde{n}+5,$ we have
\[
\big| \mathbf{P}_3(\hat{v}\times \tilde{\xi} )   \big|\lesssim 2^{2l}.
\]
Thanks to the above   estimate,  the upper bound of  the symbol $\mathbf{P}_3(\hat{v}\times \tilde{\xi} )  $  is at least   as good as the upper bound $2^{l+\max\{l,n\}+\epsilon M_t}$ we used for the symbol   $(\hat{v}-\zeta)\times (\hat{v}\times \xi/|\xi|)$ of $T_{k,j;n}^{\mu,2}( \mathfrak{m}, E)(t,x, \zeta) +  \hat{\zeta}\times T_{k,j;n}^{\mu,2}( \mathfrak{m}, B)(t,x, \zeta), i\in \{0,1,2\},$   after rerunning the arguments used in subsection \ref{Linfbigre} and subsection \ref{linfsmallre},    our desired first estimate in  \eqref{nov5eqn1}  holds after combining the   estimate  \eqref{sep27eqn1}  in  Proposition \ref{finalestfirst} and the estimate  \eqref{2022feb12eqn21}  in Proposition \ref{finalestsmallprop}. 

Similarly,   the symbol $2^n\big(\hat{v}\times \xi/|\xi|\big)$ is at least as good as the upper bound $2^{l+\max\{l,n\}+\epsilon M_t}$ we used for the  symbol  $(\hat{v}-\zeta)\times (\hat{v}\times \xi/|\xi|)$ of  $T_{k,j;n}^{\mu,2}( \mathfrak{m}, E)(t,x, \zeta) +  \hat{\zeta}\times T_{k,j;n}^{\mu,2}( \mathfrak{m}, B)(t,x, \zeta), i\in \{0,1,2,3\},$  in subsection \ref{Linfbigre} and subsection \ref{linfsmallre}. Recall the decomposition  \eqref{nov6eqn21}.  After rerunning the arguments used  in subsection \ref{Linfbigre} and subsection \ref{linfsmallre},  the   desired second estimate  in  \eqref{nov5eqn1}  holds after combining the estimate  \eqref{nov5eqn1}  and the estimate \eqref{2022feb12eqn22}  in Proposition \ref{finalestsmallprop}. 
\end{proof}

For the other component of $ \mathbf{P}_3(B(t,x))$
 , specifically ${}_{}^zT_{k,n}^{\mu;2}(B)(t,x,\zeta)$
 (see \eqref{nov6eqn21}), although we do not have a similar double null structure, there remains an advantage in that the symbol 
  $\mathbf{P}_{3}(\hat{v}\times \xi)$ is better than $  (\hat{v}\times \xi)$.

 To leverage this benefit, comparing with the dichotomy  obtained in Proposition \ref{meanLinfest} for the magnetic field $B$, we have the following improved dichotomy (by a factor of $2^{\tilde{n}}$) for the ${}_{}^zT_{k,n}^{\mu;2}(B)(t,x,\zeta)$ part of $ \mathbf{P}_3(B(t,x))$.

 \begin{lemma}\label{set1goodPartP3B}
For any fixed $\zeta\in \R^3/\{0\},  t\in [0, T)$,  $ \mu\in \{+,-\}, k, j\in \Z_+, n, \tilde{n}\in [-M_t, 2]\cap \Z, m\in [-10M_t, \epsilon M_t]\cap \Z$,  s.t.,  $n\leq \tilde{n}$,  $|  \zeta_{\bot}|/ |\zeta|\leq 2^{\tilde{n}-10}$, we have 
\be\label{nov6eqn4}
\begin{split}
&\sum_{ (m,k,j,l)\in \mathcal{S}_1(t) }  \big\|     {}_{}^zT_{k,j;n,l}^{\mu;m,2}(B)(t,x,\zeta) \big\|_{L^\infty_x} \\
&\lesssim  2^{ \tilde{n} + 10\epsilon M_t }  \big( 2^{ 2 \tilde{\alpha}_t  M_t }   + 2^{7M_t/6+\tilde{\alpha}_t M_t/4} \big),\\ 
&\sum_{ (m,k,j,l)\in \mathcal{S}_2(t) }\big(\big\|     {}_{}^zT_{k,j;n,l}^{\mu;m,2}(B)(t,x,\zeta) \big\|_{L^\infty_x}\big)^{1/2} \big(\big\|     {}_{}^zT_{k,j;n,l}^{\mu;m,2}(B)(t,x,\zeta) \big\|_{L^2_x}\big)^{1/2} \\
& \lesssim   2^{ \tilde{n} + 10\epsilon M_t }  \big( 2^{  \tilde{\alpha}_t  M_t }   + 2^{7M_t/12+\tilde{\alpha}_t M_t/8} \big). 
\end{split}
\ee
\end{lemma}

\begin{proof}
Recall  \eqref{nov6eqn21}.  Based on the possible size of $m+k$, we proceed in two steps  as follows. 

\medskip

\noindent \textbf{Step 1.}\qquad If $m+k\leq -2l+4\epsilon M_t. $

\medskip

Recall  \eqref{nov5eqn23}. Due to the fact that $|\mathbf{P}_3(\hat{v}\times \xi)|$ is of size $2^{\tilde{n}} |(\hat{v}\times \xi)|$, with minor modifications in obtaining the estimates  in Lemma \ref{smallk},   the following estimate holds if   $m+l \geq -M_t-100\epsilon M_t$  and $j\in [0, (1+2\epsilon)M_t]\cap \Z ,$  
\be\label{2021dec10eqn11}
\big(\|  {}_{}^zT_{k,j;n,l}^{\mu;m,2}(B)(t,x,\zeta) \|_{L^\infty}  \big)^{1/2}\big(\|  {}_{}^zT_{k,j;n,l}^{\mu;m,2}(B)(t,x,\zeta) \|_{L^2}   \big)^{1/2}\lesssim 2^{\tilde{n} +  M_t/2+ 10\epsilon M_t}.
 \ee
If $m+l\leq -M_t-100\epsilon M_t, $   or  $j\in [  (1+2\epsilon)M_t, \infty)\cap \Z ,$  we have 
\be\label{2021dec10eqn12}
\|  {}_{}^zT_{k,j;n,l}^{\mu;m,2}(B)(t,x,\zeta) \|_{L^\infty} \lesssim  2^{\tilde{n} + 2\tilde{\alpha}_t M_t-50\epsilon M_t}  (1+2^{j-(1+2\epsilon)M_t})^{-3}. 
\ee

\medskip

\noindent \textbf{Step 2.}\qquad  If $m+k\geq -2l+4\epsilon M_t. $

\medskip

Again,  due to the fact that $|\mathbf{P}_3(\hat{v}\times \xi)|$ is of size $2^{\tilde{n}} |(\hat{v}\times \xi)|$, with minor modifications in obtaining the estimate  \eqref{july5eqn66}  in Lemma  \ref{firstL2} and the estimate  \eqref{july5eqn111}  in Lemma \ref{largefrel2},  the following estimate holds 
\be\label{2021dec10eqn8}
\|  {}_{}^zT_{k,j;n,l}^{\mu;m,2}(B)(t,x,\zeta) \|_{L^2} \lesssim \min\{  2^{\tilde{n} + m+j+2l+ \epsilon M_t}+2^{-1000M_t}, 2^{-(1-3\epsilon)k+4j +20\epsilon M_t}\}. 
\ee 

As in   \eqref{sep19eqn60},    after using the decomposition of the acceleration force in  \eqref{july1eqn11},  we have 
\be\label{nov5eqn59}
\begin{split}
\big|{}_{}^zT_{k,j;n,l}^{\mu;m,2}(B)(t,x,\zeta) \big|&\lesssim  1+ \big| {}_{}^zT_{k,j;n,l}^{T ;\mu;m,2}(B)(t,x,\zeta)\big| \\
&\quad + \big|{}_{}^zT_{k,j;n,l}^{S,1 ;\mu;m,2}(B)(t,x,\zeta)\big| + \big|{}_{}^zT_{k,j;n,l}^{S,2 ;\mu;m,2}(B)(t,x,\zeta)\big|,
\end{split}
\ee
where
\be\label{nov5eqn25}
\begin{split}
   &{}_{}^zT_{k,j;n,l}^{T;\mu;m,2}(B)(t,x,\zeta) \\
  &= \int_0^t \int_{\R^3} \int_{\R^3} \int_{\mathbb{S}^2} \big[ i\mu  \mathbf{P}_3\big(  \omega^{m;B}_{j,l}(t-s,v,\omega)\big)  K_{k;n}^{  }(1)(y,  \zeta ) \\
  &\quad  + \mathbf{P}_3\big( c^{q;m,B}_{j,l}(t-s,v,\omega)\big) K_{k;n}^{   q }(1)(y, \zeta  ) \\
  &\quad   + (t-s)^{-1}     \mathbf{P}_3\big(  c^{err;m,B}_{j,l}(t-s,v,\omega)  \big) \tilde{K}_{k;n}^{ \mu}(1)(y,   \zeta ) \big] \\
  &\quad \times f(s,x-y+(t-s)\omega, v)   \varphi_{\tilde{n}}^2(v, \zeta) \varphi_j(v)  \psi_{\leq \tilde{n}+\epsilon M_t}(  { \tilde{\zeta} + \omega } )  d \omega dy  d v ds,\\
 & {}_{}^zT_{k,j;n,l}^{S,i;\mu;m,2}(B)(t,x,\zeta) \\
 & =  \int_0^t  \int_{\R^3} \int_{\R^3} \int_{\mathbb{S}^2}\phi_{m;-10M_t,c_t}(t-s)f(s,x-y+(t-s)\omega, v) \\
   &\quad\times EB^i(t,s,x-y ,\omega, v)\cdot \nabla_v \big(\frac{ \mathbf{P}_3\big(m_{B}(v, \omega )\big)  \varphi_{j,-l}(v, \omega)}{1+\hat{v}\cdot \omega}  \varphi_{\tilde{n}}^2(v, \zeta) \varphi_j(v) \big)\\
  &\quad \times   \big[ (t-s)  K_{k;n}^{ \mu}(1)(y,  \zeta, \omega)+ \widetilde{K}_{k;n}^{ \mu}(1)(y,   \zeta )\big]      \psi_{\leq \tilde{n}+\epsilon M_t}(   { \tilde{\zeta} +  \omega } ) d \omega dy  d v d s, 
 \end{split}
\ee
where $i\in\{1,2\}, $ the coefficients  $ \omega^{m;B}_{j,l}(t-s,v,\omega), c^{q;m,B}_{j,l}(t-s,v,\omega),$ and $ c^{err;m,B}_{j,l}(t-s,v,\omega) $ are defined in    \eqref{sep20eqn44}, and  the kernels $ K_{k;n}^{ \mu}(1)(y,  \zeta, \omega)$,  $ \widetilde{K}_{k;n}^{ \mu}(1)(y,   \zeta )$, $K_{k;n}^{ }(1)(y,  \zeta ), $ and $ K_{k;n}^{q }(1)(y,  \zeta )$ are defined  in \eqref{sep5eqn48}. 

After doing integration by parts in $\xi$ along $\zeta$ direction and directions perpendicular to $\zeta$ many times, we have
\be\label{nov5eqn65}
\begin{split}
& |K_{k;n}^{ \mu}(1)(y,  \zeta, \omega)| +  \big|K_{k;n}^{  }(1)(y, \zeta  )\big| + \big| K_{k;n}^{   q }(1)(y, \zeta  )\big| +2^k | \widetilde{K}_{k;n}^{ \mu}(1)(y,   \zeta )| \\
& \lesssim 2^{3k+2n}(1+2^{k }|y\cdot \tilde{\zeta}|)^{-N_0^3 })(1+2^{k+n}|y\times \tilde{\zeta}| )^{- N_0^3 }.\\
\end{split}
\ee
 
Recall the fact (see \eqref{2022feb8eqn11}) that    the  angle between $\omega$ and $\pm \xi$ is less than $c(m,k,l):=\min\{ \max\{-m-k-l, -m/2-k/2\}+5\epsilon M_t, 10\}$ for the essential part.   Due to the cutoff function  $ \varphi_{\tilde{n}}^2(v, \zeta)$ (see \eqref{nov6eqn21}), the essential part vanishes if $l> \max\{-j, \tilde{n}+ \epsilon M_t\}$. Hence, it suffices to consider the case $l\leq  \max\{-j, \tilde{n}+ \epsilon M_t\}$. 

 Based on the estimate  \eqref{nov5eqn59}, we proceed in   three sub-steps as follows. 

\medskip
 
 \textbf{Step 2A}. \qquad The estimate of $T$ part.

\medskip

Recall  \eqref{nov5eqn25}. From the estimate of kernels in  \eqref{nov5eqn65},  the estimate in Lemma \ref{conservationlawlemma},  the volume of support of $v, \omega$, the estimate  \eqref{nov24eqn41}  if $|  v_{\bot} |\geq 2^{(\alpha_t+\epsilon)M_t}$, we have
\be\label{2021dec10eqn1}
\begin{split}
\big|   {}_{}^zT_{k,j;n,l}^{T;\mu;m,2}(B)(t,x,\zeta) \big| & \lesssim 2^{\tilde{n}- l+5\epsilon M_t} \min\{2^{-2m-j-2l}, 2^{m+2l+j+2\tilde{\alpha}_t M_t } \} \\
 & \lesssim  2^{-m+\tilde{n}+2\tilde{\alpha}_t M_t/3-j/3-5l/3+20\epsilon M_t}. \\
 \end{split}
\ee

\medskip

  \textbf{Step 2B}. \qquad The estimate of $S_1$ part.

\medskip

Recall    the definition of coefficients  $m_{B}(v, \omega )$ in  \eqref{july9eqn11}.  From the estimate of kernels in \eqref{nov5eqn65},  the estimate in Lemma \ref{conservationlawlemma} and the volume of support, we have
\be\label{2021dec10eqn2}
\begin{split}
\big|   {}_{}^zT_{k,j;n,l}^{S,1;\mu;m,2}(B)(t,x,\zeta) \big| & \lesssim 2^{m+\tilde{n}- j-l-\min\{l, \tilde{n}\}+5\epsilon M_t} \big(2^{-2m+3j+2 \min\{l, \tilde{n}\}}\big)^{1/2} \big(  2^{-2m-j-2l} \big)^{1/2} \\
& \lesssim 2^{-m-2l+\tilde{n}+5\epsilon M_t}.
\end{split}
\ee

\medskip

  \textbf{Step 2C}. \qquad  The estimate of $S_2$ part.

\medskip

  For this case, with minor modification in obtaining the estimate  \eqref{sep17eqn11}  in Proposition \ref{bulkroughpoi}, we have
\be\label{2021dec10eqn3}
\big|   {}_{}^zT_{k,j;n,l}^{S, 2;\mu;m,2}(B)(t,x,\zeta) \big| \lesssim 2^{ \tilde{n}  + 10\epsilon M_t } \big( 2^{ 2 \tilde{\alpha}_t  M_t }   + 2^{7M_t/6+\tilde{\alpha}_t M_t/4} \big)(1+2^{-m-j-2l}) .
\ee
Recall  \eqref{nov5eqn59}. After combining the obtained estimates  \eqref{2021dec10eqn1}--\eqref{2021dec10eqn3}, we have  
\be\label{nov5eqn61}
 \big\|  {}_{}^zT_{k,j;n,l}^{\mu;m,2}(B)(t,x,\zeta) \big\|_{L^\infty_x} \lesssim   2^{ \tilde{n} + 10\epsilon M_t }  \big( 2^{ 2 \tilde{\alpha}_t  M_t }   + 2^{7M_t/6+\tilde{\alpha}_t M_t/4} \big) (1+2^{-m-j-2l}).  
\ee

Hence, the desired estimates  in \eqref{nov6eqn4}  hold  after combining the obtained estimates  \eqref{2021dec10eqn11}--\eqref{2021dec10eqn8}, the estimate  \eqref{nov5eqn61}, the estimate  \eqref{nov26eqn31}  in Lemma \ref{erroresti}, and the estimate  \eqref{july10eqn89} in  Proposition \ref{Linfielec}.  
\end{proof}

The essential notations employed in this section are systematically detailed in Table \ref{tablesection4}.
\begin{table}[H]
\centering
\resizebox{\columnwidth}{!}{%
\begin{tabular}{ |c|c|c|c| } 
 \hline
 Notation & Definition    & Remarks \\ 
\hline
$ \widetilde{T}_{k,j;n}^{\star;\mu,i }( \mathfrak{m}, U)(t,x, \zeta), \star\in\{bil,ell\}$ & \eqref{oct1eqn1} & Using the normal form transformation; only appear\\
& &  in the large anglee region, i.e., $i\in\{0,1\}$. \\
\hline
 $\widetilde{T}_{k,j;n;l,r}^{\star; \mu, m,i}(\mathfrak{m}, U)(t,x, \zeta), \star\in\{T, S\}$&  \eqref{sep18eqn44} & Using the  G-S decomposition; only appear\\
& &  in the small angle  region, i.e., $i\in\{2,3,4\}$. \\
\hline
$\mathcal{B}_i, i\in\{2,3,4\}$ & \eqref{2024oct14eqn1} & Set of parameters $l$ and $r$, which are used for the  \\
& & inhomogeneous localization for the size of $\tilde{v}+w$ and $\tilde{v}-\tilde{\zeta}$  \\
\hline
$\mathcal{K}^{ell;U ,i}_{k,j,n}(\mathfrak{m})(y, v,   \zeta), U\in\{E, B\}$ & \eqref{sep5eqn88} & The corresponding kernel of $\widetilde{T}_{k,j;n}^{ell;\mu,i }( \mathfrak{m}, U)(t,x, \zeta)$, see \eqref{sep5eqn86} \\
\hline
$\mathcal{K}^{U,i}_{k,j,n}(\mathfrak{m})(y, v, \omega, \zeta)$, $ U\in\{E, B\}$  & \eqref{sep5eqn85} & The essential kernel and the error  kernel of $\widetilde{T}_{k,j;n}^{bil;\mu,i }(  \mathfrak{m}, U)(t,x, \zeta)$ \\
$\mathcal{K}^{err;U,i}_{k,j,n}(\mathfrak{m})(y, v, \omega, \zeta)$& & \\
\hline
$\varphi_{ m; k,n}^{\star}(\omega, \zeta), \star\in\{ess,err\}$ & \eqref{sep8eqn21} & Localizing the angle between $\omega$ and $\zeta$ by \\
 & & using the stationary phase analysis\\ 
\hline
$\widetilde{K}^{\star;m,i;a}_{k,j ;n,l ,r }(t,   x,\zeta), \star\in\{main,err\}$ & \eqref{2022feb8eqn53} & Refined decomposition of $\widetilde{T}_{k,j;n}^{bil;\mu ,i }( \mathfrak{m}, E)(t,x, \zeta ) $ \\
& &  $+ \hat{\zeta}\times  \widetilde{T}_{k,j;n}^{bil;\mu,i }(\mathfrak{m}, B )(t,x, \zeta)$ based on the size of $\tilde{v}+w$ and $\tilde{v}-\tilde{\zeta}$  \\
\hline
$ H_{k,j;n,l,r}^{\mu,m,i;lin}(t, x,\zeta) $ & \eqref{sep19eqn82} &   Only used in the small angle region ($i=2,3,4$)\\
$ H_{k,j;n,l,r}^{\mu,m,i;non,a}(t, x,\zeta), a\in\{1,2\} $ & & and the non-resonance case, i.e., $m+k\leq -2l +4\epsilon M_t$\\
\hline
$ G_{k,j; n,l,r }^{\star; m, i;p,q}(t,x,\zeta)$, $\star\in\{ess,err\}$ & \eqref{sep9eqn21} & Refined   decomposition of $   \widetilde{T}_{k,j;n,l,r}^{T; \mu ,m, i}(\mathfrak{m}, E)(t,x,  \zeta ) $ \\
& &  $ + \hat{\zeta}\times   \widetilde{T}_{k,j;n,l,r}^{T; \mu,m, i }(\mathfrak{m},B)(t,x, \zeta )$ based on the size of $\sin\theta $,  $\sin\phi$ and $\omega\pm \tilde{\zeta}$.   \\
\hline
$\widetilde{K}^{m;p,q,i;a}_{k,j ;  n,l,r }(t, x, \zeta )$ & \eqref{aug8eqn31} &  Refined   decomposition of  $\widetilde{T}_{k,j;  n,l,r }^{S; \mu ,m, i }(\mathfrak{m}, E)(t,x, \zeta ) $ $+\hat{ \zeta}\times  \widetilde{T}_{k,j;  n,l,r }^{S; \mu ,m, i }(\mathfrak{m}, B)(t,x,  \zeta)$\\
& &  based on the size of $\sin\theta $,  $\sin\phi$, and the decomposition in \eqref{july1eqn11} \\
\hline
  $ {}_{}^zT_{k,n}^{\mu}(B)(t,x,\zeta)$ & \eqref{nov6eqn41}   & Angular localization of $P_3(B)$\\
    \hline
      $ {}_{}^zT_{k,n}^{\mu;i}(B)(t,x,\zeta)$ & \eqref{nov6eqn21}    & Further decomposition of  $ {}_{}^zT_{k,n}^{\mu}(B)(t,x,\zeta)$\\
      &   & based on the angle between $v$ and $\zeta$\\
    \hline
    $  {}_{}^zT_{k,j;n,l}^{\mu;m,2}(B)(t,x,\zeta)$ & \eqref{nov5eqn23} & Atomic decomposition for $ {}_{}^zT_{k,n}^{\mu;i}(B)(t,x,\zeta)$;  \\
    & & an improved   dichotomy holds for $ {}_{}^zT_{k,n}^{\mu;i}(B)(t,x,\zeta)$, see Lemma \ref{set1goodPartP3B}\\
    \hline
\end{tabular}%
}
\caption{Essential notations in section \ref{linfacceloc}.}\label{tablesection4}
\end{table}

 \section{Pointwise estimates of the localized acceleration force}\label{horizonestpotw}

  The main goal of this section is to prove the part (ii) in Theorem \ref{maintheoremellipitic} for the elliptic parts and   the part (ii) in Theorem \ref{mainresultsfirstpart} for the hyperbolic parts. For the sake of simplicity, we will continue to use the notation and terminology established in section \ref{linfacceloc}. 

The elliptic components can be easily estimated, as detailed in Lemma \ref{ellpointestpartI} in section \ref{ellipticpartspointwiseest}.  In view of  the decomposition in \eqref{oct7eqn1}, to estimate the hyperbolic parts, it suffices to estimate  the localized acceleration force $T_{k,j;n}^{\mu,i}(  \mathfrak{m}, E)(t,x, \zeta) + \hat{\zeta}\times  T_{k,j;n}^{\mu,i}(  \mathfrak{m}, B)(t,x, \zeta)$. As in section \ref{linfacceloc}, we divide into the large angle region and the small angle region. For these two regions, as in section \ref{linfacceloc}, we use two different refined decompositions  obtained in Lemma  \ref{locdeclemm}.

\subsection{pointwise estimates of the elliptic parts}\label{ellipticpartspointwiseest}
\begin{lemma}\label{ellpointestpartI}
Under the assumption of  the part (ii) in Theorem \ref{maintheoremellipitic},    for any   $j\in [0, (1+2\epsilon)M_{t_{ }}]\cap \Z,$   $i\in\{0,1\}, i'\in\{0,1,2,3\},$ we have
\be\label{sep6eqn1}
\begin{split}
 &\big| \widetilde{T}_{k,j;n}^{ell;\mu,i }( \mathfrak{m}, E)(t,x, \zeta)+ \hat{\zeta}\times  \widetilde{T}_{k,j;n}^{ell;\mu,i }( \mathfrak{m}, B )(t,x, \zeta)  \big|+ \big|\mathfrak{E}^{\mu, i'}_{k,j;n}(\mathfrak{m})(t, x, \zeta) \big| \\
&  \lesssim  \min_{b\in \{1,2\}} \| \mathfrak{m}(\cdot, \zeta)\|_{\mathcal{S}^\infty}\big[|  x_{\bot}|^{-1/2} 2^{(\gamma_1-\gamma_2)M_{t^{\star}}/2} 2^{(  {\alpha}^{\star}-10\epsilon)M_{t^{\star}}}+    |  x_{\bot}|^{-1} 2^{(\gamma_1-\gamma_2)M_t} 2^{ 5 {\alpha}^{\star} M_{t^{\star}}/6} \\
&\quad\times  \mathbf{1}_{| x_{\bot}|\leq 2^{-k-n+b\epsilon M_{t^{\star}}}}   +  |  x_{\bot}|^{-1/2}2^{(\gamma_1-\gamma_2)M_{t^{\star}}/2} \mathbf{1}_{n\geq  (\gamma_1-\gamma_2-24\epsilon)M_{t^{\star}}   }\mathbf{1}_{| x_{\bot}|\geq 2^{-k-n+   b \epsilon M_{t^{\star}}}} \\
 &\quad \times  \min\{2^{(k+2n)/2+ 2{\alpha}^{\star}M_{t^{\star}}/3}, 2^{(k+4n)/2+{\alpha}^{\star}M_{t^{\star}}-(\gamma_1-\gamma_2)M_{t^{\star}}/3}\}\big].
 \end{split}
\ee
Moreover, the following rough estimate holds, 
\be\label{2024nov14eqn41}
\begin{split}
 &\big| \widetilde{T}_{k,j;n}^{ell;\mu,i }( \mathfrak{m}, E)(t,x, \zeta)+ \hat{\zeta}\times  \widetilde{T}_{k,j;n}^{ell;\mu,i }( \mathfrak{m}, B )(t,x, \zeta)  \big|+ \big|\mathfrak{E}^{\mu, i'}_{k,j;n}(\mathfrak{m})(t, x, \zeta) \big| \\
 &\lesssim \mathfrak{m}(\cdot, \zeta)\|_{\mathcal{S}^\infty} |  x_{\bot}|^{-1/2}  2^{(  {\alpha}^{\star}+2\epsilon)M_{t^{\star}}}.
 \end{split}
\ee
  
\end{lemma}
\begin{proof}
Recall  \eqref{sep5eqn86},  \eqref{sep5eqn88}, and  \eqref{2024nov14eqn11}. Based on the possible size of $|   x_{\bot} |$, we split it into two steps as follows.      

 \medskip

\noindent \textbf{Step 1.} \qquad If $| x_{\bot} |\leq 2^{-k-n+2\epsilon M_{t^{\star}}}$. 

\medskip

From the estimate of the kernel in  \eqref{2021dec21eqn80}, the conservation law  \eqref{conservationlaw}, and the volume of support of $v$,    the estimate  \eqref{nov24eqn41}  in if $|  v_{\bot} |\geq 2^{(\alpha_t + \epsilon)M_{t^{\star}}}$,
\be\label{2021dec13eqn41}
\begin{split}
&\big| \widetilde{T}_{k,j;n}^{ell;\mu,i }( \mathfrak{m}, E)(t,x, \zeta)+ \hat{\zeta}\times  \widetilde{T}_{k,j;n}^{ell;\mu,i }( \mathfrak{m}, B )(t,x, \zeta)  \big|+ \big|\mathfrak{E}^{\mu, i'}_{k,j;n}(\mathfrak{m})(t, x, \zeta) \big|\\
& \lesssim 2^{\epsilon M_t}\min\{2^{-k+2(  {\alpha}_t +\epsilon)M_t + j},2^{-k+3j} ,  2^{2k+2n-j}\} \| \mathfrak{m}(\cdot, \zeta)\|_{\mathcal{S}^\infty}\\
&\lesssim 2^{\epsilon M_t} \big( 2^{-k+2(  {\alpha}_t +\epsilon)M_t + j}\big)^{1/6} \big(2^{-k+3j} \big)^{1/6} \big( 2^{2k+2n-j}\big)^{2/3} \| \mathfrak{m}(\cdot, \zeta)\|_{\mathcal{S}^\infty}\\
&\lesssim 2^{3\epsilon M_{t^{\star}}}    2^{k+4n/3+ {\alpha}^{\star} M_{t^{\star}}/3   }      \| \mathfrak{m}(\cdot, \zeta)\|_{\mathcal{S}^\infty} \\
& \lesssim 2^{5\epsilon M_{t^{\star}}}  |  x|^{-1} 2^{ {\alpha}^{\star} M_{t^{\star}}/3}    \| \mathfrak{m}(\cdot, \zeta)\|_{\mathcal{S}^\infty}\\
& \lesssim   |  x_{\bot}|^{-1} 2^{(\gamma_1-\gamma_2)M_t} 2^{ 5 {\alpha}^{\star} M_{t^{\star}}/6} \| \mathfrak{m}(\cdot, \zeta)\|_{\mathcal{S}^\infty}.\\
\end{split}
\ee

 \medskip

\noindent \textbf{Step 2.} \qquad   If $|  x_{\bot} |\geq 2^{-k-n+\epsilon M_{t^{\star}}}$. 

 \medskip

From the estimate of the kernel in  \eqref{2021dec21eqn80},  the conservation law  \eqref{conservationlaw},   the cylindrical symmetry of the distribution function, and the estimate  \eqref{nov24eqn41}  in if $|  v_{\bot} |\geq 2^{(\alpha_t + \epsilon)M_t}$,  we have
\be\label{2024oct17eqn31}
\begin{split}
&\big| \widetilde{T}_{k,j;n}^{ell;\mu,i }( \mathfrak{m}, E)(t,x, \zeta)+ \hat{\zeta}\times  \widetilde{T}_{k,j;n}^{ell;\mu,i }( \mathfrak{m}, B )(t,x, \zeta)  \big|+ \big|\mathfrak{E}^{\mu, i'}_{k,j;n}(\mathfrak{m})(t, x, \zeta) \big|\\
&\lesssim \| \mathfrak{m}(\cdot, \zeta)\|_{\mathcal{S}^\infty} 2^{\epsilon M_t}\min\{2^{-k+2(   {\alpha}_t  +\epsilon)M_t + j},2^{-k+3j} ,  |  x_{\bot}|^{-1}2^{k+n-j}, 2^{2k+2n-j}\}  \\
& \lesssim \| \mathfrak{m}(\cdot, \zeta)\|_{\mathcal{S}^\infty} 2^{\epsilon M_t} \min\big\{ \big(  | x_{\bot}|^{-1}2^{k+n-j}\big)^{1/2}\big(2^{-k+2(   {\alpha}_t  +\epsilon)M_t + j}\big)^{1/2}, \\
&\qquad |  x_{\bot}|^{-1}2^{k+n-j}, \big(2^{-k+2(   {\alpha}_t  +\epsilon)M_t + j} \big)^{1/3}\big( 2^{-k+3j}\big)^{1/3}\big( 2^{2k+2n-j}\big)^{1/3}  \big\}\\
\end{split}
\ee
 
From the above estimate, we conclude the desired rough estimate \eqref{2024nov14eqn41} and the following estimate, 
\be\label{2021dec13eqn42}
\begin{split}
\eqref{2024oct17eqn31}&\lesssim \big[|  x_{\bot}|^{-1/2} 2^{(\gamma_1-\gamma_2)M_{t^{\star}}/2} 2^{( {\alpha}^{\star}-10\epsilon)M_{t^{\star}}} + 2^{5\epsilon M_{t^{\star}}} \mathbf{1}_{n\geq (\gamma_1-\gamma_2 - 24\epsilon ) M_{t^{\star}}} \\
&\qquad \times \min\big\{| x_{\bot}|^{-1}2^{k+n-j}, 2^{2n/3+j + 2 {\alpha}^{\star}  M_{t^{\star}}/3} \big\} \big]\| \mathfrak{m}(\cdot, \zeta)\|_{\mathcal{S}^\infty}\\
&\lesssim \| \mathfrak{m}(\cdot, \zeta)\|_{\mathcal{S}^\infty}\big[|  x_{\bot}|^{-1/2} 2^{(\gamma_1-\gamma_2)M_{t^{\star}}/2} 2^{(  {\alpha}^{\star}-10\epsilon)M_{t^{\star}}}+    | x_{\bot}|^{-1/2}2^{(\gamma_1-\gamma_2)M_{t^{\star}}/2} \\
 &\qquad \times \mathbf{1}_{n\geq  (\gamma_1-\gamma_2)M_{t^{\star}} - 24\epsilon M_{t^{\star}}}\min\{2^{(k+2n)/2+ 2{\alpha}^{\star}M_{t^{\star}}/3}, 2^{(k+4n)/2+{\alpha}^{\star}M_{t^{\star}}-(\gamma_1-\gamma_2)M_{t^{\star}}/3}\}\big].
 \end{split}
\ee
Hence   our desired estimate  \eqref{sep6eqn1}  holds after combining the above two estimates  \eqref{2021dec13eqn41}  and  \eqref{2021dec13eqn42}. 
\end{proof}

\subsection{pointwise estimates    of  the bilinear parts   in the large angle region}\label{pointwibigre}

Recall the decomposition of $T_{k,j;n}^{\mu,i}( \mathfrak{m}, U), i\in\{0,1\}, U\in \{E, B\}, $ in  \eqref{sep4eqn30}. In this section, we focus on the estimate of the bilinear part $ \widetilde{T}_{k,j;n}^{bil;\mu,i }( \mathfrak{m}, U)(t,x, \zeta) $. 

Recall the obtained estimates in \eqref{sep5eqn99}, \eqref{sep8eqn51}, and  \eqref{2022feb13eqn1}. It suffices to estimate  the main parts  $  \widetilde{K}^{main;m,i;a }_{k,j ;n, l, r}(t, x, \zeta )   $  for the  case $m\in (-10M_t, \epsilon M_t]\cap \Z$, in which we have $t-s\sim 2^{m}. $ 
\begin{lemma}\label{projlargeregime1}
Let $i\in \{0,1\},$ $  l\in     [-j,2]\cap \Z , m\in (-10M_t, \epsilon M_t]\cap \Z $.   Under the assumption of  the part (ii) in Theorem \ref{mainresultsfirstpart},         the following estimate holds,  
\be\label{2022feb13eqn12}
\begin{split}
|\mathbf{P}\big( \widetilde{K}^{main;m,i;1}_{k,j ; n, l, r } (t,   x, \zeta )\big)| &+ 2^{(\gamma_1-\gamma_2)M_{t^\star}}  | \widetilde{K}^{main;m,i;1}_{k,j ; n, l, r } (t,   x, \zeta )  | \\
 &  \lesssim  \min_{b\in\{1,2\}}    \big[ \sum_{a\in \{0, 1/4,1/2\}} | x_{\bot}|^{-a} 2^{  a(\gamma_1-\gamma_2)M_{t^{\star}} + (\alpha^{\star}-10\epsilon) M_{t^{\star}} } \\
 &\qquad +    |  x_{\bot}|^{-1} 2^{(\gamma_1-\gamma_2) M_{t^{\star}} } 2^{ 2 {\alpha}^{\star} M_{t^{\star}}/3}  \mathbf{1}_{|  x_{\bot}|\leq 2^{-k-n+b\epsilon M_{t^{\star}}}}   \\
 &\qquad    + \mathbf{1}_{n\geq  -  (\alpha^{\star}/2+30 \epsilon)M_{t^{\star}}    }  \mathbf{1}_{|  x_{\bot} |\geq 2^{-k-n+2\epsilon M_{t^{\star}}}}  | x_{\bot}|^{-a} 2^{  a(\gamma_1-\gamma_2)M_{t^{\star}}  }\\
&\qquad \times  \min\{2^{(k+2n)/2+2\alpha^{\star} M_{t_{\star}} /3},2^{(k+4n)/2+ \alpha^{\star} M_{t^{\star}}  } \}  \big]  \| \mathfrak{m}(\cdot, \zeta)\|_{\mathcal{S}^\infty}. \\
\end{split}
\ee
\end{lemma}
\begin{proof}
Recall  \eqref{2022feb8eqn53}. Due to the cutoff functions $\varphi_{r; \vartheta^\star_0 }(\tilde{v}-\tilde{\zeta})$ and $\varphi^{i}_{j,n}(v, \zeta)$ in \eqref{sep4eqn6}, we know that $|\tilde{v}-\tilde{\zeta}|\sim 2^r\gtrsim 2^{\vartheta^\star_0}$ and $|  v_{\bot} |\sim 2^{j+ \max\{r, (\gamma_1-\gamma_2)M_{t^\star}\}}$ if $r\notin [ (\gamma_1-\gamma_2-\epsilon)M_{t^\star},  (\gamma_1-\gamma_2+\epsilon)M_{t^\star}]$.  Based on the estimate in \eqref{nov24eqn41}, it suffices to consider the case where $j+ \max\{r, (\gamma_1-\gamma_2)M_{t^\star}\}\leq (\alpha_t+\epsilon) M_t+10$,  provided that  $r\notin [ (\gamma_1-\gamma_2-\epsilon)M_{t^\star},  (\gamma_1-\gamma_2+\epsilon)M_{t^\star}]$.

  Based on the size of $m,k,l$,  we proceed in two steps as follows. 

 \medskip 

\noindent \textbf{Step 1.} \qquad If $m+k\leq -2n+4\epsilon M_{t^{\star}}$.
 
  \medskip

Recall  \eqref{2022feb8eqn53}  and   \eqref{sep8eqn21}. Based on the possible size of $| x_{\bot}|$, we proceed in two sub-steps as follows. 

  \medskip 

  \textbf{Step 1A.} \qquad  If either $|  x_{\bot}|\leq 2^{-k-n+2\epsilon M_{t^{\star}}}$, or $|  x_{\bot}|\leq 2^{m+2\epsilon M_{t^{\star}}}$.

    \medskip  

From the estimate of kernels in  \eqref{july23eqn54}, the estimate  \eqref{march18eqn31}  in Lemma \ref{conservationlawlemma},  and the volume of support of  $\omega, v$,   we have 
\be\label{aug21eqn87}
\begin{split}
  |  {\widetilde{K}}^{main;m,i;1}_{k,j ; n, l, r }(t, x, \zeta )| & \lesssim  \| \mathfrak{m}(\cdot, \zeta)\|_{\mathcal{S}^\infty} 2^{m-j-r + \epsilon M_{t^{\star}}     } \big(2^{-2m} 2^{3j+2\min\{l, r\}}  \big)^{1/2} \\
  &\quad \times  \big(  \min\{ 2^{-m-k+4\epsilon M_{t^{\star}}}   2^{m+3k+2n-j}, 2^{-2m-j-2l} \} \big)^{1/2}\\
   &\lesssim  2^{ 4\epsilon M_{t^{\star}}} \min\{2^{k+n }, 2^{-m-\max\{l,r\}}\} \| \mathfrak{m}(\cdot, \zeta)\|_{\mathcal{S}^\infty} \\
   & \lesssim |  x_{\bot}|^{-1 }2^{(\gamma_1-\gamma_2)M_{t^{\star}} +   2\alpha^{\star}  M_{t^{\star}}/3 }   \| \mathfrak{m}(\cdot, \zeta)\|_{\mathcal{S}^\infty}.
\end{split}
\ee

  \medskip 

  \textbf{Step 1B.} \qquad  If $|  x_{\bot}|\geq 2^{ \epsilon M_{t^{\star}}} \max\{2^{m}, 2^{-k-n}\}$. 

  \medskip 

For this case, we have
 $
 |  x_{\bot} -  y_{\bot}+(t-s)  \omega_{\bot}|\sim |  x_{\bot}|,$ for any $  y\in B(0, 2^{-k-n+\epsilon M_{t^{\star}}/2 }). $  From   the estimate of kernels in  \eqref{july23eqn54},   the cylindrical symmetry of solution,  the estimate  \eqref{march18eqn31}  in Lemma \ref{conservationlawlemma},  the volume of support of  $\omega, v$, and the estimate  \eqref{nov24eqn41}  if $|  v_{\bot}|\geq 2^{(\alpha_t+\epsilon)M_{t }}$,  we have 
\be\label{2024oct17eqn61}
\begin{split}
 |\mathbf{P}\big( \widetilde{K}^{main;m,i;1}_{k,j ; n, l, r } (t, &  x, \zeta )\big)| + 2^{(\gamma_1-\gamma_2)M_{t^\star}}  | \widetilde{K}^{main;m,i;1}_{k,j ; n, l, r } (t,   x, \zeta )  | \\
 & \lesssim  \| \mathfrak{m}(\cdot, \zeta)\|_{\mathcal{S}^\infty}   2^{m-j-r  + \epsilon M_{t^{\star}}     }  2^{ \max\{r, (\gamma_1-\gamma_2)M_{t^\star}\}}\\
 &\quad \times    \big(\min\big\{|  x_{\bot}|^{-1} 2^{k+n-j},  2^{-k } \min\{  2^{3j+2\min\{l,r\}}, 2^{j+2(\alpha_t+\epsilon) M_t}\} \big\} \big)^{1/2}  \\
 &\quad \times  \big(2^{-2m} \min\{  2^{3j+2\min\{l,r\}}, 2^{j+2(\alpha_t+\epsilon) M_t}\} \big)^{1/2} \\
 &\lesssim \| \mathfrak{m}(\cdot, \zeta)\|_{\mathcal{S}^\infty} 2^{j/2-r+\max\{r, (\gamma_1-\gamma_2)M_{t^\star}\}+3\epsilon M_{t^{\star}} }\min\{2^{\min\{l,r\}}, 2^{ \alpha^{\star} M_{t^{\star}}-j}\}\\
  &\quad \times \big(\min\big\{ \big(|  x_{\bot}|^{-1} 2^{k+n-j} \big)^{1/2} \big(2^{-k }   2^{ j+2\alpha_t  M_t } \big)^{1/2}, | x_{\bot}|^{-1} 2^{k+n-j}\big\} \big)^{1/2}
 \end{split}
 \ee
From the above estimate, we conclude that 
\be\label{aug21eqn90}
\begin{split}
\eqref{2024oct17eqn61}& \lesssim   2^{-r + \max\{r, (\gamma_1-\gamma_2)M_{t^\star}\}  +5\epsilon M_{t^{\star}}} \min\{ 2^{ \alpha^{\star} M_{t^{\star}}-j}, 2^{r} \}    \| \mathfrak{m}(\cdot, \zeta)\|_{\mathcal{S}^\infty} \\
&\quad \times \min\big\{|  x_{\bot}|^{-1/4} 2^{   \alpha^{\star}   M_{t^{\star}}/2 +j/2 +n/4}, |  x_{\bot}|^{-1/2}  2^{(k+n)/2 }\big\}\\
&\lesssim  \| \mathfrak{m}(\cdot, \zeta)\|_{\mathcal{S}^\infty}\big[ |  x_{\bot}|^{-1/4} 2^{(\gamma_1 -\gamma_2)M_{t^{\star}}/4  + ( \alpha^{\star}- 10\epsilon)  M_{t^{\star}}}  + |  x_{\bot}|^{-1/2} 2^{ (\gamma_1 -\gamma_2) M_{t^{\star}}/2   }   \\
&\quad \times  \mathbf{1}_{|  x_{\bot}|\geq 2^{-k-n+2\epsilon M_{t^{\star}}} } \mathbf{1}_{n\geq  (\gamma_1 -\gamma_2-60\epsilon) M_{t^{\star} }  } \min\{2^{(k+2n)/2+ \alpha^{\star} M_{t^{\star}}/2  }, 2^{(k+4n)/2+ \alpha^{\star} M_{t^{\star}}    } \}  \big]. 
\end{split}
\ee
 
 \medskip 

\noindent \textbf{Step 2.} \qquad  If $m+k\geq -2n+4\epsilon M_{t^{\star}}$.

 \medskip 

 From the estimate of kernels in  \eqref{july23eqn54},  the volume of support of the essential part of $\omega$, the estimate  \eqref{nov24eqn41}  if $|  v_{\bot}|\geq 2^{(\alpha_t+\epsilon)M_{t^{ }}}$,  and the estimate  \eqref{march18eqn31}  in Lemma \ref{conservationlawlemma}, we have  
\be\label{2022feb13eqn21}
\begin{split}
 |\mathbf{P}\big( \widetilde{K}^{main;m,i;1}_{k,j ; n, l, r } (t,   & x, \zeta )\big)|  + 2^{(\gamma_1-\gamma_2)M_{t^\star}}  | \widetilde{K}^{main;m,i;1}_{k,j ; n, l, r } (t,   x, \zeta )  | \\
  & \lesssim \| \mathfrak{m}(\cdot, \zeta)\|_{\mathcal{S}^\infty}      2^{m-j-r +5\epsilon M_{t^{\star}}}2^{ \max\{r, (\gamma_1-\gamma_2)M_{t^\star}\}}\\
& \quad \times \big(   \min\{2^{-2m-j-2l}, 2^{m+2n } \min\{2^{3j+ 2\min\{l,r\}  }, 2^{j+2  {\alpha}_t M_{t^{  }}}\}  \big)^{1/2} \\
&\quad \times  \big(2^{-2m}\min\{2^{3j+2\min\{l,r\}  }, 2^{j+2  {\alpha}_t M_{t^{  }}}\} \big)^{1/2}\\
 &\lesssim     2^{ -r +6\epsilon M_{t^{\star}} +  \max\{r, (\gamma_1-\gamma_2)M_{t^\star}\}} \min\{2^{-m }, 2^{ m/2+ n    }  \| \mathfrak{m}(\cdot, \zeta)\|_{\mathcal{S}^\infty} \\
 &\quad \times \min\{2^{2j+ 2\min\{l,r\} }, 2^{2{\alpha}_t M_{t^{  }} } \} \}.
 \end{split}
\ee

  Based on the possible size of $|  x_{\bot}|$, we  proceed in two sub-steps  as follows. 

  \medskip 

  \textbf{Step 2A.} \qquad If $|  x_{\bot}|\leq 2^{ m + \epsilon M_{t^{\star}}   }$.

    \medskip 

From the obtained estimate  \eqref{2022feb13eqn21}, we have 
\be\label{aug21eqn85}
  |  {\widetilde{K}}^{main;m,i;1}_{k,j ;n,l,r}(t, x, \zeta )|\lesssim   | x_{\bot}|^{-1/2}2^{(\gamma_1-\gamma_2)M_{t^{\star}}/2 +  (\alpha^{\star}-10\epsilon)  M_{t^{\star}}}   \| \mathfrak{m}(\cdot, \zeta)\|_{\mathcal{S}^\infty}. 
\ee
 
 \medskip 

  \textbf{Step 2B.} \qquad  If $|  x_{\bot}|\geq 2^{  m+\epsilon M_{t^{\star}}  }$.

   \medskip 

 From the estimate of kernels in  \eqref{july23eqn54}, the cylindrical symmetry of solution, the volume of support of $\omega,v$,  the estimate  \eqref{nov24eqn41}  if $|  v_{\bot}|\geq 2^{(\alpha_t+\epsilon)M_{t^{ }}}$, and the estimate  \eqref{march18eqn31} in Lemma \ref{conservationlawlemma}, we have
\[
\begin{split}
  |\mathbf{P}\big( \widetilde{K}^{main;m,i;1}_{k,j ; n, l, r } (t,  & x, \zeta )\big)| + 2^{(\gamma_1-\gamma_2)M_{t^\star}}  | \widetilde{K}^{main;m,i;1}_{k,j ; n, l, r } (t,   x, \zeta )  |\\
&\lesssim    \| \mathfrak{m}(\cdot, \zeta)\|_{\mathcal{S}^\infty} 2^{ \max\{r, (\gamma_1-\gamma_2)M_{t^\star}\}} 2^{m-j-r +5\epsilon M_{t^{\star}}}\\
&\quad \times   \big(2^{-2m} 2^{3j+    2\min\{l,r\} }  \big)^{1/2}\big(    2^{m + 2k+3n-j} |  x_{\bot}|^{-1}   \big)^{1/2}  \\
&\lesssim |  x_{\bot}|^{-1/2}  2^{m/2+ k+3n/2 +\max\{r, (\gamma_1-\gamma_2)M_{t^\star}\}  +5\epsilon M_{t^{\star}} } \| \mathfrak{m}(\cdot, \zeta)\|_{\mathcal{S}^\infty}. \\
\end{split}
\]
After combining the above estimate and the obtained estimate  \eqref{2022feb13eqn21}, we have 
\be\label{aug21eqn86}
\begin{split}
   |\mathbf{P}\big( \widetilde{K}^{main;m,i;1}_{k,j ; n, l, r } (t, &  x, \zeta )\big)|  + 2^{(\gamma_1-\gamma_2)M_{t^\star}}  | \widetilde{K}^{main;m,i;1}_{k,j ; n, l, r } (t,   x, \zeta )  | \\
&     \lesssim  \| \mathfrak{m}(\cdot, \zeta)\|_{\mathcal{S}^\infty}    2^{  6\epsilon M_{t^{\star}}}  \min\big\{2^{2n/3+4\alpha^\star M_{t^\star}/3}, 2^{-m/2 +n/3 +2\alpha^\star M_{t^\star}/3 -r}, \\
&\qquad  |  x_{\bot}|^{-1/2}  2^{m/2+ k+3n/2 +\max\{r, (\gamma_1-\gamma_2)M_{t^\star}\}} \big\}\\
&\lesssim  \| \mathfrak{m}(\cdot, \zeta)\|_{\mathcal{S}^\infty}    2^{  6\epsilon M_{t^{\star}}}  \min\big\{2^{2n/3+4\alpha^\star M_{t^\star}/3}, \big(2^{-m/2 +n/3 +2\alpha^\star M_{t^\star}/3 -r}\big)^{1/2} \\
&\qquad \times  \big( |  x_{\bot}|^{-1/2}  2^{m/2+ k+3n/2 +\max\{r, (\gamma_1-\gamma_2)M_{t^\star}\}} \big)^{1/2}\big\}\\
&\lesssim \| \mathfrak{m}(\cdot, \zeta)\|_{\mathcal{S}^\infty}   \big[ 2^{(  \alpha^{\star}-10\epsilon)M_{t^{\star}}}+| x_{\bot}|^{-1/4} 2^{(\gamma_1-\gamma_2)M_{t^{\star}}/4}     \mathbf{1}_{ n\geq -  (\alpha^{\star}/2+30 \epsilon)M_{t^{\star}}   }   \\
&\qquad \times    \min\big\{2^{(k+2n)/2+ \alpha^{\star} M_{t^{\star}}/2  -3  \vartheta^{\star}_0}, 2^{(k+4n)/2+ \alpha^{\star} M_{t^{\star}}      } \big\} \big].
\end{split}
\ee
To sum up, our desired estimate  \eqref{2022feb13eqn12}  holds after combining the obtained estimates  \eqref{aug21eqn87},  \eqref{aug21eqn90},  \eqref{aug21eqn85}, and  \eqref{aug21eqn86}.

\end{proof}
\begin{lemma}\label{projlargeregime2}
Let $i\in \{0,1\},$ $  l\in     [-j,2]\cap \Z , m\in (-10M_t, \epsilon M_t]\cap \Z $.   Under the assumption of  the part (ii) in Theorem \ref{mainresultsfirstpart}, the following estimate holds if $  m+2j+l+r\leq  \alpha^{\star}  M_{t^{\star}}   +30\epsilon M_{t^{\star}}$  \textup{or} $ n > -( \alpha^{\star}    /2+30\epsilon) M_{t^{\star}}  $, 
\be\label{2022feb16eqn74}
\begin{split}
|\mathbf{P}\big( \widetilde{K}^{main;m,i;2}_{k,j ; n, l, r } (t, &   x, \zeta )\big)| + 2^{(\gamma_1-\gamma_2)M_{t^\star}}  | \widetilde{K}^{main;m,i;2}_{k,j ; n, l, r } (t,   x, \zeta )  | \\
&  \lesssim  \min_{b\in\{1,2\}}  \| \mathfrak{m}(\cdot, \zeta)\|_{\mathcal{S}^\infty} \big[   |  x_{\bot}|^{-1/2}2^{ (\gamma_1-\gamma_2) M_{t^{\star}} /2   } 2^{   ( \alpha^{\star} -10\epsilon) M_{t^{\star}}  }\\
& \quad + |  x_{\bot}|^{-1} 2^{(\gamma_1-\gamma_2) M_{t^{\star}} + 2\alpha^{\star}M_{t^{\star}}/3}  \mathbf{1}_{|  x_{\bot}|\leq  2^{-k-n+b\epsilon M_t}}  \\
&\quad +  \mathbf{1}_{|  x_{\bot}|\geq 2^{-k-n+b\epsilon M_t}} \big( 2^{ 40\epsilon M_{t^{\star}}} \min\{ |  x_{\bot}|^{-1/2} 2^{7 \alpha^{\star}  M_{t^{\star}} /8}, |  x_{\bot}|^{-1} 2^{3 \alpha^{\star}  M_{t^{\star}} /8} \}\\
  &\quad  +  |  x_{\bot}|^{-1/2}2^{ (\gamma_1-\gamma_2) M_{t^{\star}} /2   }      \mathbf{1}_{n \geq   -({\alpha}^{\star}/2 + 30\epsilon ) M_{t^{\star}}   } \\
  &\quad \times \min\{  2^{(k+2n)/2 + 2 {\alpha}^{\star} M_{t^{\star}}/3},2^{(k+4n)/2 + {\alpha}^{\star} M_{t^{\star}}-(\gamma_1-\gamma_2) M_{t^{\star}} /4  } \} \big) \big]. 
  \end{split}
\ee  

\end{lemma}
\begin{proof}
Due to the cutoff functions $\varphi_{r; \vartheta^\star_0 }(\tilde{v}-\tilde{\zeta})$ and $\varphi^{i}_{j,n}(v, \zeta)$ in  \eqref{sep4eqn6}, we know that $|\tilde{v}-\tilde{\zeta}|\sim 2^r\gtrsim 2^{\vartheta^\star_0}$ and $|  v_{\bot} |\sim 2^{j+ \max\{r, (\gamma_1-\gamma_2)M_{t^\star}\}}$ if $r\notin [ (\gamma_1-\gamma_2-\epsilon)M_{t^\star},  (\gamma_1-\gamma_2+\epsilon)M_{t^\star}]$. From the estimate  \eqref{nov24eqn41},  it suffices to consider the case where   $j+ \max\{r, (\gamma_1-\gamma_2)M_{t^\star}\}\leq (\alpha_t+\epsilon) M_t+10$ if $r\notin [ (\gamma_1-\gamma_2-\epsilon)M_{t^\star},  (\gamma_1-\gamma_2+\epsilon)M_{t^\star}]$.

Similar to the obtained estimate  \eqref{aug5eqn1}  for $  |\ {\widetilde{K}}^{main;m,i;2}_{k,j ;l,n}(t, x, \zeta  )| $, from  the estimate of kernels in  \eqref{july23eqn54},  we have
\be\label{2022feb23eqn31}
\begin{split}
& |\mathbf{P}\big( \widetilde{K}^{main;m,i;2}_{k,j ; n, l, r } (t,   x, \zeta )\big)| + 2^{(\gamma_1-\gamma_2)M_{t^\star}}  | \widetilde{K}^{main;m,i;2}_{k,j ; n, l, r } (t,   x, \zeta )  | \\
 & \lesssim  \sum_{p, q\in [-10M_t,2]\cap \Z } 2^{\max\{r, (\gamma_1-\gamma_2) M_{t^\star}\}} \big|{H}^{m, i ;p,q  }_{k,j ; n,l,r  } (t,   x,  \zeta)\big|,\\
 \end{split}
\ee
where ${H}^{m, i ;p,q  }_{k,j ; n,l,r  } (t,   x,  \zeta)$ is defined in  \eqref{aug5eqn1}.

Based on the possible size of $|  x_{\bot}|$, we proceed in three steps   as follows.

\medskip
  
\noindent \textbf{Step 1.}\qquad If $ |  x_{\bot}|\leq 2^{-k-n+2\epsilon M_t}  $.

\medskip

Based on the possible size of $m+k$, we  proceed in three sub-steps    as follows. 

\medskip
  
  \textbf{Step 1A.} \qquad If $m+k\leq -2n + \epsilon M_t. $

\medskip
  
Recall the definition of the cutoff function $  \varphi_{ m; k,n}^{main}(\omega, \zeta)$ in   \eqref{sep8eqn21}. From the  volume of support of $\omega, $  the estimate of kernels in \eqref{july23eqn54},  the Cauchy-Schwarz inequality, the conservation law  \eqref{conservationlaw},  the estimate   \eqref{march18eqn31} in Lemma \ref{conservationlawlemma},  and the volume of support of $(\omega, v)$, we have
\be\label{oct8eqn31}
\begin{split}
\big|{H}^{m, i;p,q }_{k,j ; n, l,r }(t,   x,  \zeta)\big|&\lesssim  2^{m-j- r +l +4\epsilon M_t} \| \mathfrak{m}(\cdot, \zeta)\|_{\mathcal{S}^\infty} \sup_{s\in[0,t]}\|B(s,\cdot)\|_{L^2_x} \\
&\quad \times  \big( 2^{m+3k+2n} 2^{-m-k}   2^{3j+2 \min\{l,r\} } \} \big)^{1/2}
\big(  2^{-2m-j-2l}   \big)^{1/2} \\
&\lesssim   \| \mathfrak{m}(\cdot, \zeta)\|_{\mathcal{S}^\infty}  2^{    6\epsilon M_{t^{\star}}  }  2^{k+n}\lesssim |  x_{\bot}|^{-1} 2^{   8\epsilon M_{t^{\star}} } \| \mathfrak{m}(\cdot, \zeta)\|_{\mathcal{S}^\infty} \\
& \lesssim | x_{\bot}|^{-1} 2^{(\gamma_1-\gamma_2) M_{t^{\star}}  + 2\alpha^{\star}M_{t^{\star}}/3} \| \mathfrak{m}(\cdot, \zeta)\|_{\mathcal{S}^\infty} . \\
\end{split}
\ee

\medskip
  
  \textbf{Step 1B.} \qquad  If $m+k\geq -2n + \epsilon M_t$, $m+p \leq - k -n + 4\epsilon M_t.$

\medskip

Recall again the definition of the cutoff function $  \varphi_{ m; k,n}^{main}(\omega, \zeta)$ in  \eqref{sep8eqn21}. Note that, the volume of support of $\omega$ is bounded by $2^{2\min\{n, p\}}$. Hence, from the estimate of kernels in  \eqref{july23eqn54},  the Cauchy-Schwarz inequality, the conservation law  \eqref{conservationlaw},  the estimate \eqref{march18eqn31}  in Lemma \ref{conservationlawlemma},  and the volume of support of $(\omega, v)$, we have
\be\label{2022feb14eqn1}
\begin{split}
\big|{H}^{m, i;p,q }_{k,j ; n, l,r }(t,   x,  \zeta)\big|&\lesssim  2^{m-j-r +l +4\epsilon M_t}  \sup_{s\in[0,t]}\|B(s,\cdot)\|_{L^2_x} \| \mathfrak{m}(\cdot, \zeta)\|_{\mathcal{S}^\infty}\\
&\quad \times  \big( 2^{m+3k+2n}  2^{2\min\{n, p\}}   2^{3j+2  \min\{l,r\} } \} \big)^{1/2}
\big(  2^{-2m-j-2l}   \big)^{1/2} \\
&\lesssim  \| \mathfrak{m}(\cdot, \zeta)\|_{\mathcal{S}^\infty}  2^{    6\epsilon M_{t^{\star}}  }  2^{k+n}\lesssim |  x_{\bot}|^{-1} 2^{    8\epsilon M_{t^{\star}} }  \| \mathfrak{m}(\cdot, \zeta)\|_{\mathcal{S}^\infty}   \\
&\lesssim |  x_{\bot}|^{-1} 2^{(\gamma_1-\gamma_2) M_{t^{\star}} + 2\alpha^{\star}M_{t^{\star}}/3} \| \mathfrak{m}(\cdot, \zeta)\|_{\mathcal{S}^\infty}  . 
\end{split}
\ee

\medskip
  
  \textbf{Step 1C.} \qquad  If $m+k\geq -2n + \epsilon M_t$, $m+p \geq - k -n +4 \epsilon M_t.$

\medskip

For this  case, we have $| x_{\bot} -   y_{\bot} + (t-s) \omega_{\bot}|\sim 2^{m+p}$ for any $y\in B(0, 2^{-k-n+\epsilon M_t/2})$. From  the estimate of kernels in  \eqref{july23eqn54}, the cylindrical symmetry of solution,   the Cauchy-Schwarz inequality, the conservation law  \eqref{conservationlaw},  the estimate  \eqref{march18eqn31}  in Lemma \ref{conservationlawlemma},  and the volume of support of $(\omega, v)$, we have  
\be\label{2022feb14eqn3}
\begin{split}
\big|{H}^{m, i;p,q }_{k,j ; n, l,r }(t,   x,  \zeta)\big|&\lesssim  2^{m-j-r  +l +4\epsilon M_t}  \| \mathfrak{m}(\cdot, \zeta)\|_{\mathcal{S}^\infty}  \sup_{s\in[0,t]}\|B(s,\cdot)\|_{L^2_x}  \\
&\quad \times \big( 2^{m+3k+2n}  \frac{2^{-k-n+\epsilon M_t} 2^{2\min\{n, p\}} }{2^{m+p}}   2^{3j+2 \min\{l,r\} } \} \big)^{1/2}
\big(  2^{-2m-j-2l}   \big)^{1/2}\\
& \lesssim  \| \mathfrak{m}(\cdot, \zeta)\|_{\mathcal{S}^\infty}  2^{    6\epsilon M_{t^{\star}}  }  2^{k+n} \\
& \lesssim |  x_{\bot}|^{-1} 2^{(\gamma_1-\gamma_2) M_{t^{\star}} + 2\alpha^{\star}M_{t^{\star}}/3} \| \mathfrak{m}(\cdot, \zeta)\|_{\mathcal{S}^\infty} .
\end{split} 
\ee

\medskip
  
\noindent \textbf{Step 2.}\qquad  \qquad If $ 2^{-k-n+ \epsilon M_t}\leq |  x_{\bot}| < 2^{\epsilon M_t}\max\{ 2^{m +p },2^{-k-n  }\}    $.

\medskip

Due to the restriction of $ 2^{-k-n+ \epsilon M_t}\leq |  x_{\bot}| < 2^{\epsilon M_t}\max\{ 2^{m +p },2^{-k-n  }\}    $, we only have to consider the case $m+p + k+n\geq 0$. Additionally, for any $y\in B(0, 2^{-k-n+\epsilon M_t/2}),$ we have$|  x_{\bot} -  y_{\bot}|\sim |  x_{\bot}|$. 

 Based on the possible size of $m+k$, we proceed in three sub-steps    as follows. 

\medskip
  
  \textbf{Step 2A.}\qquad If $m+k\leq -2n+\epsilon M_t, m+k\leq \alpha^{\star} M_{t^\star} +62\epsilon M_{t^\star}  . $

  \medskip

Recall  \eqref{aug5eqn1}. From   the Cauchy-Schwarz inequality, the conservation law  \eqref{conservationlaw},  the estimate  \eqref{march18eqn31} in Lemma \ref{conservationlawlemma}, the Jacobian of  changing coordinates $(\theta, \phi) \longrightarrow (z(y),w(y))$,  the estimate  \eqref{nov24eqn41}  if $|  v_{\bot}|\geq 2^{(\alpha_t+\epsilon)M_{t^{   }}}$, and the volume of support of $\omega$ of the essential part, we have
\be\label{2024Dec8eqn2}
\begin{split}
&2^{\max\{r, (\gamma_1-\gamma_2) M_{t^\star}\}} \big|{H}^{m, i;p,q }_{k,j ; n, l,r }(t,   x,  \zeta)\big| \\
& \lesssim  \sup_{s\in[0,t]}\|B(s,\cdot)\|_{L^2_x} \| \mathfrak{m}(\cdot, \zeta)\|_{\mathcal{S}^\infty}   2^{m-j- r +l +4\epsilon M_t}    2^{\max\{r, (\gamma_1-\gamma_2) M_{t^\star}\}}  \\
&\quad \times \big( \min\big\{2^{-m}2^{-p-q} | x_{\bot}|^{-1} , 2^{m+3k+2n} 2^{p+q}2^{-m/2-k/2+\epsilon M_{t^{\star}}/2} \} \\
& \quad  \times   \min\{2^{3j+2\min\{r,l\} },2^{j+2  \alpha_{t}  M_{t^{ }}  } \} \big)^{1/2}  \big(  \min\{  2^{-2m-j-2l},  2^{m-m-k+\epsilon M_{t^{\star}} + j+2  {\alpha}_t M_{t^{  }}}\} \big)^{1/2}.\\
\end{split}
\ee
From the above estimate, we have
\be\label{2024oct17eqn81}
\begin{split}
\eqref{2024Dec8eqn2}& \lesssim \| \mathfrak{m}(\cdot, \zeta)\|_{\mathcal{S}^\infty}   2^{m-j- r +l +4\epsilon M_t}    2^{\max\{r, (\gamma_1-\gamma_2) M_{t^\star}\}}\min\big\{ \big(2^{-2m-j-2l} \big)^{1/2}\\
& \quad  \times   \big(   (2^{-m}2^{-p-q} |  x_{\bot}|^{-1}  )^{1/2}   (2^{2k+p+q}2^{-m/2-k/2+\epsilon M_{t^{\star}}/2}   )^{1/2} 2^{3j+2r }  \big)^{1/2} , \big(2^{-k+\epsilon M_{t^{\star}} + j+2  {\alpha}_t M_{t^{  }}} \big)^{1/2} \\
&\quad \times  \big(   (2^{-m}2^{-p-q} |  x_{\bot}|^{-1})^{1/2}   (2^{2k+p+q}2^{-m/2-k/2+\epsilon M_{t^{\star}}/2} )^{1/2}\min\{ 2^{3j+2r }, 2^{j+2  {\alpha}_t M_{t^{  }}}\} \big)^{1/2}\big\}\\
 &\lesssim  |  x_{\bot}|^{-1/4}  2^{ \max\{r, (\gamma_1-\gamma_2) M_{t^\star}\}+   7\epsilon M_{t^{\star}}}  \| \mathfrak{m}(\cdot, \zeta)\|_{\mathcal{S}^\infty}  \\
 &\quad \times \min\big\{2^{-3m/4+3(m+k)/8}, 2^{3m/4-(m+k)/8} \min\{2^{2  \alpha^{\star}  M_{t^{\star}} - r }, 2^{j+ \alpha^{\star}  M_{t^{\star}}}\} \big\}\\
 &\lesssim  2^{ 40\epsilon M_{t^{\star}}}\min\{ |  x_{\bot}|^{-1/2} 2^{7 \alpha^{\star}  M_{t^{\star}} /8}, | x_{\bot}|^{-1} 2^{3 \alpha^{\star}  M_{t^{\star}} /8  } \}  \| \mathfrak{m}(\cdot, \zeta)\|_{\mathcal{S}^\infty}. 
\end{split}
\ee
 
 \medskip
  
  \textbf{Step 2B.}\qquad If $m+k\leq -2n+\epsilon M_t, m+k\geq \alpha^{\star} M_{t^\star} +62\epsilon M_{t^\star}   . $
 \medskip

Note that, for this case, we have $n\leq   -( \alpha^{\star}    /2+ 31\epsilon) M_{t^{\star}}  $. Hence, by the assumption in  the Lemma, we have $m+2j+r+l\leq   \alpha^{\star}  M_{t^{\star}}   +30\epsilon M_{t^{\star}}.$
From the estimate of kernels in  \eqref{july23eqn54},  the Cauchy-Schwarz inequality, the conservation law  \eqref{conservationlaw},  the estimate  \eqref{march18eqn31}  in Lemma \ref{conservationlawlemma}, the Jacobian of  changing coordinates $(\theta, \phi) \longrightarrow (z(y),w(y))$,    and the volume of support of $\omega$ of the essential part, we have
 \be\label{2022feb14eqn10}
 \begin{split}
  \big|{H}^{m, i;p,q }_{k,j ; n, l,r  }(t,   x,  \zeta)\big| &  \lesssim  2^{m-j-r  +l +4\epsilon M_t}    \| \mathfrak{m}(\cdot, \zeta)\|_{\mathcal{S}^\infty}  \big(    2^{m +p + q-(m+k)/2+ 3j+2 \min\{l, r\}}    \big)^{1/2}
\\
 &\quad \times   \big(  2^{-m}2^{-p-q} |  x_{\bot}|^{-1}      2^{3j+ 2 \min\{l, r\} } \} \big)^{1/2}     \big(\sup_{s\in[0,t]}\|B(s,\cdot)\|_{L^2_x}\big) \\ 
 & \lesssim |  x_{\bot}|^{-1/2}  2^{-(m+k)/4+ \alpha^{\star}  M_{t^{\star}}  +34\epsilon M_{t^{\star}}}  \| \mathfrak{m}(\cdot, \zeta)\|_{\mathcal{S}^\infty} \\
&\lesssim |  x_{\bot} |^{-1/2} 2^{3\alpha^{\star}  M_{t^{\star}}/4+34\epsilon  M_{t^{\star}}}  \| \mathfrak{m}(\cdot, \zeta)\|_{\mathcal{S}^\infty}\\ &   \lesssim |  x_{\bot}|^{-1/2} 2^{(\gamma_1-\gamma_2)M_{t^{\star}}/2}  2^{   (\alpha^{\star}-10\epsilon)  M_{t^{\star}}   }  \| \mathfrak{m}(\cdot, \zeta)\|_{\mathcal{S}^\infty}. 
\end{split}
\ee

 \medskip
  
  \textbf{Step 2C.}\qquad  If $m+k\geq -2n+\epsilon M_t. $
 \medskip

Recall again  \eqref{sep8eqn21}. Note that, the volume of support of $\omega$ is bounded by $2^{\min\{2n,   2p + q, p+q + n\}}$.  
From the estimate of kernels in  \eqref{july23eqn54},   the Jacobian of  changing coordinates $ (y_1, y_2, \theta) \longrightarrow x-y+(t-s)\omega$ and the volume of support of $ \phi$ and the volume of support of the essential part of $y_3$, we have
\be\label{aug23eqn80}
\begin{split}
 \big|{H}^{m, i;p,q }_{k,j ; n, l,r  }(t,   x,  \zeta)\big| 
 &  \lesssim 2^{m-j-r  +l +4\epsilon M_t}  \| \mathfrak{m}(\cdot, \zeta)\|_{\mathcal{S}^\infty}  \sup_{s\in[0,t]}\|B(s,\cdot)\|_{L^2_x} \\
 &\quad \times  \big( \min\{ 2^{m+p+q+n+3j+2 \min\{l,r\} },2^{-2m-j-2l} \}   \big)^{1/2} \big(  2^{3j+2\min\{l,r\}}  \\ 
&\quad \times   \min\{  2^{-m-p-q}|   x_{\bot} |^{-1}, 2^{ 3k+2n +q}2^{-k-n+\max\{(\gamma_1-\gamma_2)M_{t^\star}, n\}+\epsilon M_{t^{\star}}}   \}  \big)^{1/2} \\ 
& \lesssim  2^{   5\epsilon M_{t^{\star}}}  \sup_{s\in[0,t]}\|B(s,\cdot)\|_{L^2_x}   \| \mathfrak{m}(\cdot, \zeta)\|_{\mathcal{S}^\infty} \\
&\quad\times \min\big\{ |  x_{\bot}|^{-1/2}  2^{m+2j+r+ l+n/2},   |   x|^{-1/2} 2^{ (2k+n+  \max\{(\gamma_1-\gamma_2)M_{t^\star}, n\}  )/4} \big\}\\
\end{split}
 \ee

 From the above estimate, we conclude that 
 \be\label{2022feb14eqn11}
 \begin{split}
    \big|{H}^{m, i;p,q }_{k,j ; n, l,r  }(t,   x,  \zeta)\big|  &\lesssim \| \mathfrak{m}(\cdot, \zeta)\|_{\mathcal{S}^\infty}  |  x_{\bot}|^{-1/2}   2^{ (\gamma_1-\gamma_2) M_{t^{\star}}/2}\big[ 2^{  ({\alpha}^{\star}-10\epsilon) M_{t^{\star}}   } +  \mathbf{1}_{n \geq -({\alpha}^{\star}/2 + 30\epsilon) M_{t^{\star}}   } \\
  &\quad \times   \min\{  2^{(k+2n)/2  +  5 {\alpha}^{\star} M_{t^{\star}}/8 },2^{(k+4n)/2 + {\alpha}^{\star} M_{t^{\star}}-(\gamma_1-\gamma_2) M_{t^{\star}} /4 } \}\big].
  \end{split}
 \ee

 \medskip
  
 \noindent  \textbf{Step 3.}\qquad  If $  |  x_{\bot}| \geq  2^{\epsilon M_t}\max\{ 2^{m +p },2^{-k-n  }\}    $.

 \medskip

Note that,  for any $y\in B(0, 2^{-k-n+\epsilon M_t/2}),$ we have $| x_{\bot} - y_{\bot}+(t-s) \omega_{\bot} |\sim |  x_{\bot}|$. Based on the possible size of $m+k$, we proceed in two sub-steps    as follows.

 \medskip
  
   \textbf{Step 3A.} \qquad If $m+k\leq -2n+\epsilon M_t. $

 \medskip

 From the cylindrical symmetry of the distribution function, the estimate of kernels in (\eqref{july23eqn54}),  the Cauchy-Schwarz inequality, the conservation law   \eqref{conservationlaw},  the estimate  \eqref{march18eqn31}  in Lemma \ref{conservationlawlemma}, and the volume of support of $\omega$, we have
\be\label{aug22eqn10}
\begin{split}
  \big|{H}^{m, i;p,q }_{k,j ;  n, l,r  }(t,   x,  \zeta)\big|  &  \lesssim    \| \mathfrak{m}(\cdot, \zeta)\|_{\mathcal{S}^\infty} \sup_{s\in[0,t]}\|B(s,\cdot)\|_{L^2_x}  2^{m-j- r +l +5\epsilon M_t}  \\
  &\quad \times    \big(    2^{m+3k+2n} 2^{-m-k } 2^{-k-n}|  x_{\bot} |^{-1}    2^{3j+  2 \min\{l,r\}  } \big)^{1/2} \\
& \quad \times \big(   \min\{ 2^{m}2^{-m-k}2^{3j+   2 \min\{l,r\}   },  2^{-2m-j-2l} \} \big)^{1/2}\\
& \lesssim  |  x_{\bot}|^{-1/2} 2^{  7\epsilon M_{t^{\star}}} \min\{  2^{(k+n)/2}, 2^{m+2j+l+r+n/2}\}  \| \mathfrak{m}(\cdot, \zeta)\|_{\mathcal{S}^\infty}\\
& \lesssim   | x_{\bot} |^{-1/2} 2^{ 30\epsilon M_{t^{\star}}}  \| \mathfrak{m}(\cdot, \zeta)\|_{\mathcal{S}^\infty}\big[ 2^{ 3 \alpha^{\star}  M_{t^{\star}}  /4  }   +    2^{(k+n)/2} \mathbf{1}_{n \geq -({\alpha}^{\star}/2 + 30\epsilon ) M_{t^{\star}}   } \big]. 
\end{split}
\ee

 \medskip
  
   \textbf{Step 3B.} \qquad If $m+k\geq -2n+\epsilon M_t. $
  \medskip
  
From the estimate of kernels in  \eqref{july23eqn54},  the Cauchy-Schwarz inequality, the conservation law  \eqref{conservationlaw},  the Jacobian of  changing coordinates $(\theta, \phi) \longrightarrow (z(y),w(y))$ (see e.g.,  \eqref{march18eqn66}),   and the volume of support of $(\omega, v)$, we have
\be\label{aug22eqn55}
\begin{split}
  \big|{H}^{m, i;p,q }_{k,j ;  n, l,r  }(t,   x,  \zeta)\big| &    \lesssim   \sup_{s\in[0,t]}\|B(s,\cdot)\|_{L^2_x}   \| \mathfrak{m}(\cdot, \zeta)\|_{\mathcal{S}^\infty}  2^{m-j- r  +l +5\epsilon M_t}\\
  &\quad \times   \big( \min\{ 2^{m+p+q+n+3j+2 \min\{l,r\} },2^{-2m-j-2l} \}   \big)^{1/2}\\
& \quad \times \big( |  x_{\bot}|^{-1} \min\{  2^{-m-p-q}, 2^{m+3k+2n+p+q+n}2^{-k-n }  \} 2^{3j+2 \min\{l,r\} } \big)^{1/2} \\
&\lesssim   |  x_{\bot}|^{-1/2} 2^{  7\epsilon M_{t^{\star}}} \min\big\{   2^{m+2j+ r+ l+n/2},     2^{ (k+n)/2} \big\} \| \mathfrak{m}(\cdot, \zeta)\|_{\mathcal{S}^\infty}\\
  &\lesssim   |  x_{\bot} |^{-1/2} 2^{ 30\epsilon M_{t^{\star}}}  \| \mathfrak{m}(\cdot, \zeta)\|_{\mathcal{S}^\infty}\big[ 2^{ 3 \alpha^{\star}  M_{t^{\star}}  /4  }   +     2^{(k+n)/2} \mathbf{1}_{n \geq -({\alpha}^{\star}/2 + 30\epsilon ) M_{t^{\star}}   } \big].  \\
\end{split}
 \ee 
Recall  \eqref{2022feb23eqn31}. To sum up, the desired estimate  \eqref{2022feb16eqn74}  holds after combining the obtained estimates  \eqref{oct8eqn31}--\eqref{2022feb14eqn10} ,    and  \eqref{2022feb14eqn11}--\eqref{aug22eqn55}. 
 \end{proof}

\begin{lemma}\label{projlargeregime4}
Let $i\in \{0,1\},$ $  l\in     [-j,2]\cap \Z , m\in (-10M_t, \epsilon M_t]\cap \Z $.   Under the assumption of  the part (ii) in Theorem \ref{mainresultsfirstpart},   the following estimate holds if $  m+2j+l+r\geq  \alpha^{\star}  M_{t^{\star}}  +30\epsilon M_{t^{\star}}$  \textup{and} $ n \leq  -( \alpha^{\star}    /2+30\epsilon) M_{t^{\star}}  $, 
\be\label{2022feb14eqn31}
  |\ {\widetilde{K}}^{main;m,i;2}_{k,j ;n,l,r}(t, x, \zeta  )|  \lesssim  \sum_{a\in\{0,1/4,1/2\}}  |  x_{\bot}|^{-a}2^{ a(\gamma_1-\gamma_2) M_{t^{\star}}     } 2^{   ( \alpha^{\star} -10\epsilon) M_{t^{\star}}  }.
  \ee
\end{lemma}
\begin{proof}
For this case, we localize the magnetic  field further by  using the decomposition in \eqref{july5eqn1}  and use the decomposition obtained in  \eqref{aug10eqn31}. 

Recall  \eqref{aug3eqn32}.   From the estimate of kernels in  \eqref{july23eqn54}  
and the estimate  \eqref{march18eqn31}  in Lemma \ref{conservationlawlemma}, the following estimate holds,
\be\label{aug23eqn61}
 \big| {H}^{m,i;p,q;\tilde{m}, \tilde{k},\tilde{j},\tilde{l} }_{k,j ; n,l,r }(t,x, \zeta )\big| \lesssim  \| \mathfrak{m}(\cdot, \zeta)\|_{\mathcal{S}^\infty} \big( \sup_{s\in[0,t]}\|  B^{\tilde{m}}_{\tilde{k};\tilde{j},  \tilde{l}}(s,\cdot)\|_{L^\infty_x } \big)     2^{m-j- r+l +5\epsilon M_t}       2^{-2m-j-2l}  . 
\ee

From the above estimate and the estimate in  \eqref{aug4eqn10} in Proposition \ref{meanLinfest},    the following estimate holds  if $(\tilde{m}, \tilde{k}, \tilde{j}, \tilde{l}) \in \mathcal{S}_1(t)$,
\be\label{2022feb14eqn36}
\begin{split}
 \sum_{(\tilde{m}, \tilde{k}, \tilde{j}, \tilde{l}) \in \mathcal{S}_1(t) } \big| {H}^{m,i;p,q;\tilde{m}, \tilde{k},\tilde{j},\tilde{l} }_{k,j ;l,n}(t,x, \zeta )\big|  &\lesssim 2^{  15\epsilon M_t } 2^{-(m+2j+r+l)} 2^{ 2  {\alpha}^{\star}  M_t }  \| \mathfrak{m}(\cdot, \zeta)\|_{\mathcal{S}^\infty} \\
 & \lesssim \| \mathfrak{m}(\cdot, \zeta)\|_{\mathcal{S}^\infty} 2^{(\alpha^{\star}-10\epsilon)M_t}. 
 \end{split}
\ee

Therefore, from now on, we restrict ourselves to the case $(\tilde{m}, \tilde{k}, \tilde{j}, \tilde{l}) \in \mathcal{S}_2(t)$. Based on  the possible size of $|  x_{\bot}|$, we proceed in three steps   as follows.

\medskip

\noindent \textbf{Step 1.} \qquad If $ | x_{\bot}|\leq 2^{-k-n+2\epsilon M_t}  $.

\medskip

Similar to the obtained estimates  \eqref{oct8eqn31}--\eqref{2022feb14eqn3}, regardless the size of $m+k$, we have
\[
\big| {H}^{m,i;p,q;\tilde{m}, \tilde{k},\tilde{j},\tilde{l} }_{k,j ;l,n}(t,x, \zeta )\big| \lesssim \| \mathfrak{m}(\cdot, \zeta)\|_{\mathcal{S}^\infty} \big(\sup_{s\in[0,t]}\|  B^{\tilde{m}}_{\tilde{k};\tilde{j},  \tilde{l}}(s,\cdot)\|_{L^2_x }  \big)   |  x_{\bot}|^{-1} 2^{    8\epsilon M_{t^{\star}} }. 
\]
After combining the above estimate, the estimate  \eqref{aug23eqn61}), and the second estimate  in   \eqref{aug4eqn10}  in Proposition \ref{meanLinfest}, we have
\be\label{2022feb18eqn31}
\begin{split}
& \sum_{(\tilde{m}, \tilde{k}, \tilde{j}, \tilde{l}) \in \mathcal{S}_2(t) } \big| {H}^{m,i;p,q;\tilde{m}, \tilde{k},\tilde{j},\tilde{l} }_{k,j ;l,n}(t,x, \zeta )\big|\\
& \lesssim |  x_{\bot}|^{-1/2} 2^{ 8\epsilon M_t-(m+2j+r+l)/2}  2^{11\epsilon M_t} \big( 2^{   \tilde{\alpha}_t  M_t }   + 2^{7M_t/12+\tilde{\alpha}_t M_t/8} \big)\| \mathfrak{m}(\cdot, \zeta)\|_{\mathcal{S}^\infty} \\
&\lesssim  |  x_{\bot}|^{-1/2} 2^{(\gamma_1-\gamma_2)M_{t^{\star}}/2  +   ( \alpha^{\star} -10\epsilon) M_{t^{\star}}  }\| \mathfrak{m}(\cdot, \zeta)\|_{\mathcal{S}^\infty}. \\
\end{split}
\ee

\medskip

\noindent \textbf{Step 2.} \qquad If $ 2^{-k-n+ \epsilon M_t}\leq |  x_{\bot}| < 2^{\epsilon M_t}\max\{ 2^{m +p },2^{-k-n  }\}    $.

\medskip

 Based on the possible size of $m+k$, we proceed in three sub-steps   as follows.

\medskip

\textbf{Step 2A.}\qquad If $m+k\leq -2n+\epsilon M_t, m+k\leq  \alpha^{\star} M_{t^\star}  . $

\medskip

 Similar to the obtained estimate  \eqref{2024oct17eqn81}, we have
 \[
 \begin{split}
 \big| {H}^{m,i;p,q;\tilde{m}, \tilde{k},\tilde{j},\tilde{l} }_{k,j ;l,n}(t,x, \zeta )\big| 
  & \lesssim 2^{40\epsilon M_{t^{\star}}} \| \mathfrak{m}(\cdot, \zeta)\|_{\mathcal{S}^\infty}  \big(\sup_{s\in[0,t]}\|  B^{\tilde{m}}_{\tilde{k};\tilde{j},  \tilde{l}}(s,\cdot)\|_{L^2_x }  \big)\\
  &\quad \times \min\big\{ |  x_{\bot}|^{-1/2} 2^{7 \alpha^{\star}  M_{t^{\star}} /8}, |  x_{\bot} |^{-1} 2^{3 \alpha^{\star}  M_{t^{\star}} /8  } \big\} .
  \end{split}
 \]
After combining the above estimate, the estimate  \eqref{aug23eqn61}, and the second estimate in   \eqref{aug4eqn10}  in Proposition \ref{meanLinfest}, we have
\[
 \sum_{(\tilde{m}, \tilde{k}, \tilde{j}, \tilde{l}) \in \mathcal{S}_2(t) } \big| {H}^{m,i;p,q;\tilde{m}, \tilde{k},\tilde{j},\tilde{l} }_{k,j ;l,n}(t,x, \zeta )\big| \lesssim 2^{30\epsilon M_{t^{\star}}+\alpha^{\star} M_{t^\star} }  2^{-(m+2j+2l)/2}  | x_{\bot}|^{-1/2} 2^{3 \alpha^{\star}  M_{t^{\star}} /16  } \| \mathfrak{m}(\cdot, \zeta)\|_{\mathcal{S}^\infty} 
\]
\be
 \lesssim | x_{\bot}|^{-1/2}  
2^{(\gamma_1-\gamma_2)M_{t^{\star}}/2+(\alpha^{\star}-10\epsilon) M_{t^{\star}} }
\| \mathfrak{m}(\cdot, \zeta)\|_{\mathcal{S}^\infty}.
\ee

\medskip

\textbf{Step 2B.}\qquad  If $m+k\leq -2n+\epsilon M_t, m+k\geq \alpha^{\star} M_{t^\star}   . $

\medskip

Similar to the obtained estimate  \eqref{2022feb14eqn10}, we have
\[
 \big|{H}^{m, i;p,q }_{k,j ;l,n}(t,   x,  \zeta)\big|   \lesssim |  x_{\bot}|^{-1/2}  2^{-(m+k)/4  +34\epsilon M_{t^{\star}} +(m+2j+r+ l)}   \| \mathfrak{m}(\cdot, \zeta)\|_{\mathcal{S}^\infty} \big(\sup_{s\in[0,t]}\|  B^{\tilde{m}}_{\tilde{k};\tilde{j},  \tilde{l}}(s,\cdot)\|_{L^2_x }  \big).
\]
After combining the above estimate, the estimate  \eqref{aug23eqn61}, and the second estimate in   \eqref{aug4eqn10}  in Proposition \ref{meanLinfest}, we have
\be 
\begin{split}
\sum_{(\tilde{m}, \tilde{k}, \tilde{j}, \tilde{l}) \in \mathcal{S}_2(t) } \big| {H}^{m,i;p,q;\tilde{m}, \tilde{k},\tilde{j},\tilde{l} }_{k,j ;l,n}(t,x, \zeta )\big| & \lesssim |  x_{\bot}|^{-1/4} 2^{ \alpha^{\star} M_{t^\star}  +20\epsilon M_{t^\star}  } 2^{-(m+k)/8}  \| \mathfrak{m}(\cdot, \zeta)\|_{\mathcal{S}^\infty} \\
&\lesssim  |  x_{\bot}|^{-1/4} 2^{(\gamma_1-\gamma_2)M_{t^{\star}}/4  +   ( \alpha^{\star} -10\epsilon) M_{t^{\star}}  }\| \mathfrak{m}(\cdot, \zeta)\|_{\mathcal{S}^\infty}. 
\end{split}
\ee

 \medskip

\textbf{Step 2C.}\qquad  \qquad If $m+k\geq -2n+\epsilon M_t. $

 \medskip
 
Similar to the obtained estimate \eqref{aug23eqn80}, we have
\[
 \big|{H}^{m, i;p,q }_{k,j ;l,n}(t,   x,  \zeta)\big|   \lesssim 2^{   5\epsilon M_{t^{\star}}}  |  x_{\bot}|^{-1/2}  2^{m+2j+r+l+n/2}   \| \mathfrak{m}(\cdot, \zeta)\|_{\mathcal{S}^\infty} \big(\sup_{s\in[0,t]}\|  B^{\tilde{m}}_{\tilde{k};\tilde{j},  \tilde{l}}(s,\cdot)\|_{L^2_x }  \big).
\]
After combining the above estimate, the estimate  \eqref{aug23eqn61}, and the second estimate  \eqref{aug4eqn10}  in Proposition \ref{meanLinfest}, we have
\be
\begin{split}
 \sum_{(\tilde{m}, \tilde{k}, \tilde{j}, \tilde{l}) \in \mathcal{S}_2(t) } \big| {H}^{m,i;p,q;\tilde{m}, \tilde{k},\tilde{j},\tilde{l} }_{k,j ;l,n}(t,x, \zeta )\big| &\lesssim |  x_{\bot}|^{-1/4} 2^{ \alpha^{\star} M_{t^\star}  +20\epsilon M_{t^\star}  } 2^{n/4}  \| \mathfrak{m}(\cdot, \zeta)\|_{\mathcal{S}^\infty}\\
&\lesssim  |  x_{\bot}|^{-1/4} 2^{(\gamma_1-\gamma_2)M_{t^{\star}}/4  +   ( \alpha^{\star} -10\epsilon) M_{t^{\star}}  }\| \mathfrak{m}(\cdot, \zeta)\|_{\mathcal{S}^\infty}. \\
\end{split}
\ee

\medskip

\noindent \textbf{Step 3.} \qquad   If $  |  x_{\bot}| \geq  2^{\epsilon M_t}\max\{ 2^{m +p },2^{-k-n  }\}    $.

\medskip

Similar to the obtained estimates  \eqref{aug22eqn10}  and  \eqref{aug22eqn55}, regardless the size of $m+k$, we have
\[
\big| {H}^{m,i;p,q;\tilde{m}, \tilde{k},\tilde{j},\tilde{l} }_{k,j ;l,n}(t,x, \zeta )\big| \lesssim  |  x_{\bot}|^{-1/2} 2^{  7\epsilon M_{t^{\star}}}    2^{m+2j+r+l+n/2}  \| \mathfrak{m}(\cdot, \zeta)\|_{\mathcal{S}^\infty} \big(\sup_{s\in[0,t]}\|  B^{\tilde{m}}_{\tilde{k};\tilde{j},  \tilde{l}}(s,\cdot)\|_{L^2_x }  \big).
\]
After combining the above estimate, the estimate \eqref{aug23eqn61}, and the second estimate  \eqref{aug4eqn10}  in Proposition \ref{meanLinfest}, we have
\be\label{2022feb14eqn37}
\begin{split}
 \sum_{(\tilde{m}, \tilde{k}, \tilde{j}, \tilde{l}) \in \mathcal{S}_2(t) } \big| {H}^{m,i;p,q;\tilde{m}, \tilde{k},\tilde{j},\tilde{l} }_{k,j ;l,n}(t,x, \zeta )\big| &\lesssim | x_{\bot}|^{-1/4} 2^{ \alpha^{\star} M_{t^\star}  +20\epsilon M_{t^\star}  } 2^{n/4}  \| \mathfrak{m}(\cdot, \zeta)\|_{\mathcal{S}^\infty}\\
&\lesssim  |  x_{\bot}|^{-1/4} 2^{(\gamma_1-\gamma_2)M_{t^{\star}}/4  +   ( \alpha^{\star} -10\epsilon) M_{t^{\star}}  }\| \mathfrak{m}(\cdot, \zeta)\|_{\mathcal{S}^\infty}. 
\end{split}
\ee
To sum up, our desired estimate  \eqref{2022feb14eqn31} holds from the obtained estimates  \eqref{2022feb14eqn36}--\eqref{2022feb14eqn37}.

\end{proof}

To sum up, the main result of this  section can be summarized as follows. 
\begin{proposition}\label{pointestfinalpartI}
       Under the assumption of  the part (ii) in Theorem \ref{mainresultsfirstpart},  for any    $ j\in [0, (1+2\epsilon)M_{t_{ }}]\cap\Z,     n\in [-M_{t_{ }}, 2]\cap \Z,$  the following estimate holds for  $i\in\{0,1\}, $
\be\label{sep27eqn21} 
\begin{split}
&\big| \mathbf{P}\big(  \mathfrak{H}_{k,j;n}^{\mu,i}( \mathfrak{m})(t, x, \zeta)  \big) \big|   + \big|\mathbf{P}\big( \widetilde{T}_{k,j;n}^{bil;\mu,i}( \mathfrak{m}, E)(t,x, \zeta )+ \hat{\zeta}\times  \widetilde{T}_{k,j;n}^{bil;\mu, i}(\mathfrak{m}, B )(t,x, \zeta)\big)  \big| \\
&\quad+ 2^{(\gamma_1-\gamma_2)M_{t^\star}}\Big(\big|  \mathfrak{H}_{k,j;n}^{\mu,i}( \mathfrak{m})(t, x, \zeta) \big|   + \big| \widetilde{T}_{k,j;n}^{bil;\mu,i}( \mathfrak{m}, E)(t,x, \zeta )+ \hat{\zeta}\times  \widetilde{T}_{k,j;n}^{bil;\mu, i}(\mathfrak{m}, B )(t,x, \zeta)  \big| \Big) \\
&\lesssim  \min_{b\in\{1,2\}} \| \mathfrak{m}(\cdot, \zeta)\|_{\mathcal{S}^\infty}  \Big[ \sum_{a\in \{0, 1/4,1/2\}} |  x_{\bot}|^{-a} 2^{  a(\gamma_1-\gamma_2)M_{t^{\star}} + (\alpha^{\star}-10\epsilon) M_{t^{\star}} }  \\
&\quad  +  |  x_{\bot}|^{-1} 2^{  (\gamma_1-\gamma_2)M_{t^{\star}}} 2^{ 5{\alpha}^{\star} M_{t^{\star}}/6}   \mathbf{1}_{|  x_{\bot}|\leq 2^{-k-n+b\epsilon M_{t^{\star}}}} +2^{40\epsilon M_{t^{\star}}}  \mathbf{1}_{| x_{\bot} |\geq 2^{-k-n+b \epsilon M_{t^{\star}}}}\\
&\quad \times \big( \min\{ |  x_{\bot}|^{-1/2} 2^{7 \alpha^{\star}  M_{t^{\star}} /8}, |  x_{\bot}|^{-1} 2^{3 \alpha^{\star}  M_{t^{\star}} /8} \} +  |  x_{\bot}|^{-a} 2^{  a(\gamma_1-\gamma_2)M_{t^{\star}} }  \mathbf{1}_{n \geq   -({\alpha}^{\star}/2 + 30\epsilon ) M_{t^{\star}}   } \\
& \quad \times  \min\{2^{(k+2n)/2+2\alpha^{\star}M_{t^{\star}}/3},2^{(k+4n)/2+ \alpha^{\star}M_{t^{\star}}-(\gamma_1-\gamma_2)M_{t^{\star}} /3 } \}    \big)      \Big]. 
 \end{split}
\ee

\end{proposition}
\begin{proof}
 Recall the equality \eqref{2024nov14eqn32},  the decomposition in  \eqref{sep4eqn30}  and the estimate  \eqref{sep5eqn99}. Our desired estimate  \eqref{sep27eqn21}  holds after   combining the estimate  \eqref{sep6eqn1}  in Lemma \ref{ellpointestpartI},  the obtained estimates   \eqref{sep8eqn51}  and    \eqref{2022feb13eqn1},   the estimate  \eqref{2022feb13eqn12}  in Lemma \ref{projlargeregime1}, the estimate  \eqref{2022feb16eqn74}  in Lemma \ref{projlargeregime2},  and the estimate  \eqref{2022feb14eqn31}  in Lemma \ref{projlargeregime4}. 

\end{proof}

\subsection{pointwise estimates    in the small angle region}

 Because   the  estimate of the localized acceleration force with a loss $2^{-(\gamma_1-\gamma_2)M_{t^{\star}}}$ can be obtained as a byproduct, in this  section, we mainly estimate     $\mathbf{P}\big( T_{k,j;n,l,r}^{\mu,m, i}(\mathfrak{m},E)(t,x, \zeta ) $ $+ \hat{\zeta}\times T_{k,j;n,l,r}^{\mu,m, i}(\mathfrak{m},B)(t,x, \zeta)  \big)$, in the small angle region,  i.e., $ i\in\{2,3,4\}.$

We first consider the case   $m+k$ is relatively small, i.e., the case $m+k\leq -2l+\epsilon M_{t^{\star}}$. For this case, we mainly use the first formulation in \eqref{sep18eqn50} in Lemma \ref{locdeclemm}.  Recall the decomposition we did in  \eqref{sep19eqn81}  and  \eqref{sep22eqn1}. Let $\omega= (\sin \theta\cos\phi, \sin \theta \sin \phi, \cos\theta)$. After  doing dyadic decomposition for the sizes of $   \sin \theta,\sin\phi$ for the linear part and the nonlinear parts,   we have  
\be\label{sep24eqn46}
\begin{split}
& \mathbf{P}\big( T_{k,j;n,l,r}^{\mu,m, i}(\mathfrak{m},E)(t,x, \zeta ) + \hat{\zeta}\times T_{k,j;n,l,r}^{\mu,m, i}(\mathfrak{m},B)(t,x, \zeta) \big)\\
&= \sum_{ \begin{subarray}{c}
 \star\in\{ess , err \}   \end{subarray}} \big[ \sum_{ p,q \in [-10M_t, 2]\cap \Z}    \mathbf{P}\big(  H_{k,j;n,l,r ;  p,q;\star}^{\mu,m,i;lin}(t, x,\zeta) \big) \\
 &\quad  +\sum_{a=1,2} \sum_{ p,q \in [-10M_t, 2]\cap \Z}   \mathbf{P}\big(H_{k,j;n,l,r; p,q;\star}^{\mu,m,i;non,a}(t, x, \zeta)\big)\big],
 \end{split}
\ee
where 
\be\label{2022feb15eqn1}
\begin{split}
 H_{k,j;n,l,r; p,q;\star}^{\mu,m,i;lin}(t, x,\zeta)& :=   \int_{0}^t \int_{\R^3} \int_{\R^3} \int_{\mathbb{S}^2} \big[ (t-s) \mathcal{K}_{k,n,l,r}^{\star;\mu,m}(y, v,\omega,  \zeta) +\mathcal{K}_{k,n,l,r}^{\star;err,\mu,m}(y, v,\omega,\zeta) \big] \\
 &\quad \times  f(s, x-y+(t-s)\omega, v)  \varphi_{l; r}(\tilde{v}+\omega )  \varphi_{m;-10M_t }(t-s) \\
& \quad\times     \varphi_{p;-10M_t}(\sin \theta) \varphi_{q;-10M_t}(\sin \phi)  \varphi_{j,n}^{i; r}(v, \zeta) d\omega dy d v ds,
\end{split}
\ee
\be\label{2022feb15eqn2}
\begin{split}
H_{k,j;n,l,r; p,q;\star}^{\mu,m,i;non,a}(t, x,\zeta)&= \int_{0}^t \int_{\R^3} \int_{\R^3} \int_{\mathbb{S}^2}   \varphi_{m;-10M_t }(t-s)   EB^a(t,s,x-y +(t-s)\omega,\omega, v) \\
& \quad \cdot  \nabla_v \big( (t-s) \mathcal{H}^{\star;\mu,E,i}_{k,j;n,l,r}(y,\omega, v, \zeta)+\mathcal{H}^{\star;err;\mu,E,i}_{k,j;n,l,r }(y, v,\omega,  \zeta)\big) \\
& \quad \times f(s, x-y+(t-s)\omega, v)  \varphi_{p;-10M_t}(\sin \theta) \varphi_{q;-10M_t}(\sin \phi)      d \omega dy d v ds,
 \end{split}
\ee
where the kernels appeared above are defined in  \eqref{2022feb8eqn11}.

Similar to  the obtained estimates  \eqref{sep20eqn21}  and  \eqref{sep20eqn27}, it suffices to  focus on the estimate of essential parts $H_{k,j;n,l,r;p,q;ess}^{\mu,m,i;lin}(t, x,\zeta)$ and $H_{k,j;n,l,r;p,q;ess}^{\mu,m,i;non,a}(t, x,\zeta)$ for the case $p,q,r\in (-10M_t, 2]\cap \Z.$ 

\begin{lemma}\label{smallkfull1}
Let  $i=2,  ( l,r)\in \mathcal{B}_2 $, see \eqref{sep18eqn50}.    Under the assumption of  the part (ii) in Theorem \ref{mainresultsfirstpart},  the following estimate holds if  $m+k \leq -2l +4\epsilon M_t, $
\be\label{2022feb16eqn53}
\begin{split}
 &\big|\mathbf{P}\big( T_{k,j;n,l,r}^{\mu,m, i}(\mathfrak{m},E)(t,x, \zeta ) + \hat{\zeta}\times T_{k,j;n,l,r}^{\mu,m, i}(\mathfrak{m},B)(t,x, \zeta) \big) \big| \\
 & \lesssim  \min_{b\in\{1,2\}} \big[\sum_{a\in\{1/8,1/4, 3/8,1/2\}} |  x_{\bot}|^{-a }  2^{a(\gamma_1-\gamma_2)M_{t^{\star} }+(\alpha^{\star} -10\epsilon )M_{t^{\star}}  } +\mathbf{1}_{|  x_{\bot}|\leq 2^{-k-n+b\epsilon M_t}} \\
&\qquad \times  |  x_{\bot}|^{-1}2^{(\gamma_1-\gamma_2)M_{t^\star}}2^{5\alpha^{\star}M_{t^{\star}}/6} +  \mathbf{1}_{n\geq  ( \gamma_1-\gamma_2- 100\epsilon )M_{t^{\star} } } \mathbf{1}_{|  x_{\bot}| >  2^{-k-n+b\epsilon M_t}}   | x_{\bot}|^{-a }  2^{a(\gamma_1-\gamma_2)M_{t^{\star} } }\\
& \qquad \times   \min\{2^{(k+2n)/2+ 2\alpha^{\star}M_{t^\star}/3- (\gamma_1-\gamma_2)M_{t^{\star} }/6}, 2^{(k+2n)/2+  \alpha^{\star}M_{t^\star} - (\gamma_1-\gamma_2)M_{t^{\star} }/3 }\}\big]\| \mathfrak{m}(\cdot, \zeta)\|_{\mathcal{S}^\infty}.
\end{split}
\ee
\end{lemma}
\begin{proof}

Recall \eqref{2022feb15eqn1}. Due to the  cutoff functions $ \varphi_{j,n}^{i; r}(v, \zeta)  \varphi_{l; r}(\tilde{v}+\omega ) $, see \eqref{sep17eqn63}  and   \eqref{sep4eqn6}, we have $n +3\epsilon M_t/4\leq r\leq l, |\tilde{v}-\tilde{\zeta}|\sim 2^r.$  and   $n\leq -   M_t/100. $  

Moreover, if $r\notin [(\gamma_1-\gamma_2-\epsilon )M_{t^\star}, (\gamma_1-\gamma_2+\epsilon)M_{t^\star}]$, then we have $|  v_{\bot} |/|v|\sim 2^{\max\{r, (\gamma_1-\gamma_2  )M_{t^\star} \}+j}. $ From the estimate  \eqref{nov24eqn41}, we can rule out the case $\max\{r, (\gamma_1-\gamma_2  )M_{t^\star} \}+j\geq  (\alpha_t+\epsilon) M_t$ if $r\notin [(\gamma_1-\gamma_2-\epsilon )M_{t^\star}, (\gamma_1-\gamma_2+\epsilon)M_{t^\star}]$. Hence, it  suffices to consider the case  $\max\{r, (\gamma_1-\gamma_2  )M_{t^\star} \}+j\leq  (\alpha_t+\epsilon) M_t$ if $r\notin [(\gamma_1-\gamma_2-\epsilon )M_{t^\star}, (\gamma_1-\gamma_2+\epsilon)M_{t^\star}]$.

Recall  \eqref{2022feb8eqn11},  the   result of computation  of symbols  in  \eqref{sep5eqn101}  and  \eqref{sep5eqn96}. For the case we are considering, we have
\be\label{sep22eqn44}
\begin{split}
|\mathbf{P} \big( m_{E}^{0}(\xi, v) + \hat{\zeta}\times m_{B}^{0}(\xi, v) \big)|  & \lesssim  2^{  \max\{ r,   (\gamma_1-\gamma_2)M_{t^{\star}}\}} \min\{2^{2r}, 2^{2\max\{l, c(m,k,l)\}}\}, \\ 
 |\mathbf{P} \big( m_{E}^{1}(\xi, v) + \hat{\zeta}\times m_{B}^{1}(\xi, v) \big)| &  \lesssim  2^{ \max\{   r, (\gamma_1-\gamma_2)M_{t^{\star}}\}}.  
\end{split}
\ee

We proceed in two main steps to estimate the     linear part and the nonlinear part as follows. 

\medskip

\noindent \textbf{Step 1.}\quad  The estimate of $\mathbf{P}\big(   H_{k,j;n,l,r; p,q;ess }^{\mu,m,i;lin}(t, x,\zeta)\big)  $.

\medskip

Recall  \eqref{2022feb8eqn1}. From the   estimate of the symbol in \eqref{sep22eqn44}. The following estimate holds after doing integration by parts in $\xi$ along $\zeta$ direction and directions perpendicular to $\zeta$ many times,
\be\label{2022feb14eqn61}
\begin{split}
 &\big| \mathbf{P} \big( \mathcal{K}_{k,n,l,r}^{ess;\mu,m}(y, v, \omega, \zeta)\big)  \big| + 2^k \big|\mathbf{P} \big( \mathcal{K}_{k,n,l,r}^{ess;err;\mu,m}(y, v, \omega, \zeta)\big) \big| \\
  & \lesssim 2^{4k+ 2n   } 2^{  \max\{ r,   (\gamma_1-\gamma_2)M_{t^{\star}}\}}\min\{2^{2r}, 2^{2\max\{l, c(m,k,l)\}}\}\\
    & \qquad  \times    (1+2^k|y\cdot \tilde{\zeta}|)^{-N_0^3} (1+2^{k+n }|y\times \tilde{\zeta}|)^{-N_0^3} \| \mathfrak{m}(\cdot, \zeta)\|_{\mathcal{S}^\infty}  .\\
  \end{split}
\ee

Based on the size of $|  x_{\bot}|$, we proceed in three sub-steps as follows. 

\medskip

\textbf{Step 1A.}\quad   If $|  x_{\bot}|\leq 2^{-k-n+2\epsilon M_t}.$

\medskip

Recall  \eqref{2022feb8eqn1}. From the   estimate of kernels in  \eqref{2022feb14eqn61}  and the volume of support of $v,\omega$, we have 
\be\label{2022jan8eqn51}
\begin{split}
\big| \mathbf{P}\big(   H_{k,j;n,l,r; p,q;ess }^{\mu,m,i;lin}(t, x,\zeta)\big)\big|&\lesssim  \| \mathfrak{m}(\cdot, \zeta)\|_{\mathcal{S}^\infty}  2^{ 2  \max\{l, c(m,k,l)\} + 4\epsilon M_t} (2^{m+k}+1) \\
&\quad \times   \min\{2^{m+3k+2n+2l-j},2^{m+2l}\min\{2^{j+2\alpha^{\star} M_{t^{\star}}}, 2^{3j+2r}\}\} \\
&  \lesssim  2^{ 4\epsilon M_t} \min\{2^{2k+2n-j},2^{-k  } 2^{2j+\alpha^{\star}M_{t^{\star}}} \}\| \mathfrak{m}(\cdot, \zeta)\|_{\mathcal{S}^\infty} \\
& \lesssim  2^{ 4\epsilon M_t} \big(2^{2k+2n-j}\big)^{1/3}\big(2^{-k  } 2^{2j+\alpha^{\star}M_{t^{\star}}}\big)^{2/3}\| \mathfrak{m}(\cdot, \zeta)\|_{\mathcal{S}^\infty}\\
 &  \lesssim  |  x_{\bot}|^{-1}2^{(\gamma_1-\gamma_2)M_{t^\star}}2^{(5\alpha^{\star} /6-\epsilon) M_{t^\star}} \| \mathfrak{m}(\cdot, \zeta)\|_{\mathcal{S}^\infty} .
    \end{split}
\ee

\medskip

\textbf{Step 1B.}\quad   If   $ 2^{-k-n+\epsilon M_t} \leq  |  x_{\bot}|< 2^{\epsilon M_t}\max\{2^{m}, 2^{-k-n}\} $. 

\medskip

  From the   estimate of kernels in  \eqref{2022feb14eqn61}, the estimate  \eqref{march18eqn31}  in Lemma \ref{conservationlawlemma},    the volume of support of $v,\omega$, and the estimate   \eqref{nov24eqn41}  if $| v_{\bot}|\geq 2^{(\alpha_t +\epsilon)M_t}$,   we have 
\be
\begin{split}
&\big|\mathbf{P}\big( H_{k,j;n,l; ess}^{\mu,m,i;lin}(t, x,\zeta) \big)\big|\\
&\lesssim   \| \mathfrak{m}(\cdot, \zeta)\|_{\mathcal{S}^\infty}   2^{ 2  \max\{l, c(m,k,l)\} + 4\epsilon M_t}   2^{ \max\{  r,   (\gamma_1-\gamma_2)M_{t^{\star}}\} }(2^{m+k}+1)  
\\
&\qquad \times \min\{2^{m+3k+2n+2l-j}, 2^{m+2l+j+2\alpha_t M_t}, 2^{m+2l+2r +3j}    \} \\
&\lesssim   \| \mathfrak{m}(\cdot, \zeta)\|_{\mathcal{S}^\infty} 2^{ \max\{r,    (\gamma_1-\gamma_2)M_{t^{\star}}\}+ 6\epsilon M_t } \min\big\{2^{2k+2n-j},2^{-m-2k  }\min\{2^{3j+2r-2l }, 2^{ j+ 2\alpha_t M_t -2l}   \}   \big\}\\
&\lesssim  \| \mathfrak{m}(\cdot, \zeta)\|_{\mathcal{S}^\infty}  2^{ \max\{r,    (\gamma_1-\gamma_2)M_{t^{\star}}\}+ 6\epsilon M_t }  \min\{ |  x_{\bot}|^{-1/2} 2^{n  } \min\{2^{j }, 2^{ \alpha_t M_t-r  }\}, 2^{2k+2n-j} \}.
\end{split}
\ee
From the above estimate, we conclude that 
\be
\begin{split}
\big|\mathbf{P}\big( H_{k,j;n,l; ess}^{\mu,m,i;lin}(t, x,\zeta) \big)\big| &\lesssim  \| \mathfrak{m}(\cdot, \zeta)\|_{\mathcal{S}^\infty}\big[|  x_{\bot}|^{-1/2}  2^{(\gamma_1-\gamma_2)M_{t^{\star} }/2+(\alpha^{\star} -11\epsilon )M_{t^{\star}}  } \\
&\quad +  | x_{\bot}|^{-3/8}  2^{3(\gamma_1-\gamma_2)M_{t^{\star} }/8  -\epsilon  M_{t^{\star}} } \mathbf{1}_{n\geq  ( \gamma_1-\gamma_2- 30\epsilon )M_{t^{\star} }/2 } \\ 
&\quad\times  \min\{2^{(k+2n)/2+ 2\alpha^{\star}M_{t^\star}/3- (\gamma_1-\gamma_2)M_{t^{\star} }/6}, 2^{(k+4n)/2+  \alpha^{\star}M_{t^\star}}\}  \big]. 
\end{split}
\ee

\medskip

\textbf{Step 1C.}\quad If   $   |  x_{\bot}|\geq  2^{\epsilon M_t}\max\{2^{m}, 2^{-k-n}\} $.

\medskip

Note that, for any $y\in B(0, 2^{-k-n+\epsilon M_t/2}), $ we have $|  x_{\bot} -   y_{\bot} +(t-s)  \omega_{\bot} |\sim |     x_{\bot} |.$  From the   estimate of kernels in  \eqref{2022feb14eqn61}, the cylindrical symmetry of the distribution function,  the estimate  \eqref{march18eqn31}  in Lemma \ref{conservationlawlemma},    the volume of support of $v,\omega$, and the estimate   \eqref{nov24eqn41}  if $| v_{\bot}|\geq 2^{(\alpha_t +\epsilon)M_t}$,   we have 
\[
\begin{split}
\big|\mathbf{P}\big( H_{k,j;n,l; ess}^{\mu,m,i;lin}(t, x,\zeta) \big)\big|&\lesssim   \| \mathfrak{m}(\cdot, \zeta)\|_{\mathcal{S}^\infty}   2^{ 2  \max\{l, c(m,k,l)\} + 4\epsilon M_t}   2^{ \max\{  r,   (\gamma_1-\gamma_2)M_{t^{\star}}\}}(2^{m+k}+1)  \\ 
 &\quad \times  \min\big\{|  x_{\bot}|^{-1} 2^{m+2k+n+2l-j}, 2^{m+2r} \min\{2^{j+2\alpha_t M_t}, 2^{3j+2l}\}   \big\} \\
 &  
\lesssim 2^{ \max\{r,    (\gamma_1-\gamma_2)M_{t^{\star}}\}+ 6\epsilon M_t } \min\big\{|  x_{\bot}|^{-1}2^{ k+ n-j},2^{- k  }\min\{2^{3j+2r  }, 2^{ j+ 2\alpha_t M_t }   \}   \big\}\\
&\lesssim 2^{ \max\{r,    (\gamma_1-\gamma_2)M_{t^{\star}}\}+ 6\epsilon M_t } \min\big\{ \big( |  x_{\bot}|^{-1}2^{ k+ n-j}\big)^{1/2}\big(  2^{ -k+j+ 2\alpha_t M_t }\big)^{1/2},\\
&\quad  \big( |  x_{\bot}|^{-1}2^{ k+ n-j}\big)^{2/3}\big(  2^{ -k+2j+ \alpha_t M_t }\big)^{1/3}\big\}.\\
\end{split}
\]
From the above estimate, we conclude that
\be
\begin{split}
&\big|\mathbf{P}\big( H_{k,j;n,l; ess}^{\mu,m,i;lin}(t, x,\zeta) \big)\big|\\
&  \lesssim  \| \mathfrak{m}(\cdot, \zeta)\|_{\mathcal{S}^\infty}  2^{ \max\{r,    (\gamma_1-\gamma_2)M_{t^{\star}}\}+ 8\epsilon M_t } | x_{\bot}|^{-1/2} \\
&\quad \times  \min\big\{  2^{n/2  }  2^{ \alpha_t M_t  } , 2^{(k+n)/2+ n/3+ \alpha_t M_t/3 } \big\} \\
& \lesssim \| \mathfrak{m}(\cdot, \zeta)\|_{\mathcal{S}^\infty}  |  x_{\bot}|^{-1/2}  2^{(\gamma_1-\gamma_2)M_{t^{\star} }/2  -\epsilon  M_{t^{\star}}}   \big[ 2^{ (\alpha^{\star} -10\epsilon )M_{t^{\star}} }+      \mathbf{1}_{n\geq  ( \gamma_1-\gamma_2-  40\epsilon )M_{t^{\star} } }\\
&\quad \times \min\{2^{(k+2n)/2+ 2\alpha^{\star}M_{t^\star}/3}, 2^{(k+4n)/2+  \alpha^{\star}M_{t^\star} - (\gamma_1-\gamma_2)M_{t^{\star} }/3 }\} \big]. 
\end{split}
\ee

 \medskip

\noindent \textbf{Step 2.}  \qquad  The estimate of $\mathbf{P}\big(  H_{k,j;n,l,r; p,q;ess}^{\mu,m,i;non,a}(t, x,\zeta) \big), a \in \{1,2\}$. 

 \medskip

Recall the estimate of the  symbol in  \eqref{sep22eqn44}. After doing integration by parts in $\xi$ along $\zeta$ direction and directions perpendicular to $\zeta$ many times, we have 
 \be\label{2022feb16eqn41}
 \begin{split}
 & \big| \nabla_v \mathbf{P}\big( \mathfrak{H}^{ess; \mu,E,i}_{k,j;n,l,r}(y,\omega, v, \zeta) \big) \big| + 2^k|\nabla_v  \mathbf{P}\big(\mathfrak{H}^{ess;  err;\mu,E,i}_{k,j;n,l,r}(y,\omega, v, \zeta)\big) |\\
& \lesssim  \| \mathfrak{m}(\cdot, \zeta)\|_{\mathcal{S}^\infty} 2^{3k+2n-j + \max\{r,    (\gamma_1-\gamma_2)M_{t^{\star}}\}+ \epsilon M_t}    2^{- r}  (1+2^k|y\cdot \tilde{\zeta}|)^{-N_0^3} (1+2^{k+n }|y\times \tilde{\zeta}|)^{-N_0^3}. 
\end{split}
\ee

Based on the size of $|  x_{\bot}|$, we proceed in three sub-steps as follows. 

\medskip

\textbf{Step 2A.}\quad    If $|  x_{\bot}|\leq 2^{-k-n+2\epsilon M_t}.$

\medskip

Recall  \eqref{2022feb15eqn2} and \eqref{july1eqn13}. From the   estimate of kernels in  \eqref{2022feb16eqn41}, the estimate  \eqref{march18eqn31}  in Lemma \ref{conservationlawlemma}, and the volume of support of $v,\omega$, we have 
\be
\begin{split}
\big|\mathbf{P}\big(  H_{k,j;n,l,r; p,q;ess}^{\mu,m,i;non,1}(t, x,\zeta) \big) \big| & \lesssim\| \mathfrak{m}(\cdot, \zeta)\|_{\mathcal{S}^\infty}     (2^{m }+2^{-k}) 2^{-j-r+2\epsilon M_t} \\
&\quad\times \big( \min\{2^{-2m}, 2^{m+3k+2n}\} 2^{3j+2r } \big)^{1/2} \big( 2^{m+3k+2n+2l-j}    \big)^{1/2} 
\\
    &\lesssim   2^{k+n+10\epsilon M_t } \| \mathfrak{m}(\cdot, \zeta)\|_{\mathcal{S}^\infty} \\
  & \lesssim    |  x_{\bot}|^{-1}2^{(\gamma_1-\gamma_2)M_{t^\star}}2^{ \alpha^{\star}M_{t^{\star}}/2 -\epsilon  M_{t^{\star}}} \| \mathfrak{m}(\cdot, \zeta)\|_{\mathcal{S}^\infty},\\
  \big|\mathbf{P}\big(  H_{k,j;n,l,r; p,q;ess}^{\mu,m,i;non,2}(t, x,\zeta) \big) \big|& \lesssim \| \mathfrak{m}(\cdot, \zeta)\|_{\mathcal{S}^\infty}    (2^{m }+2^{-k}) 2^{-j-r+l+2\epsilon M_t} \\
&\quad\times \big(  2^{m+3k+2n+2l}  2^{3j+2r } \big)^{1/2} \big( 2^{m+3k+2n+2l-j}    \big)^{1/2} \\
   &\lesssim   2^{k+n+10\epsilon M_t} \| \mathfrak{m}(\cdot, \zeta)\|_{\mathcal{S}^\infty} \\
   & \lesssim    |  x_{\bot}|^{-1}2^{(\gamma_1-\gamma_2)M_{t^\star}}2^{ \alpha^{\star}M_{t^{\star}}/2 -\epsilon  M_{t^{\star}}} \| \mathfrak{m}(\cdot, \zeta)\|_{\mathcal{S}^\infty} .
\end{split}
\ee
 
\medskip

\textbf{Step 2B.}\quad  If   $ 2^{-k-n+\epsilon M_t} \leq  |  x_{\bot}|< 2^{\epsilon M_t}\max\{2^{m}, 2^{-k-n}\} $.
\medskip

  From the   estimate of kernels in  \eqref{2022feb16eqn41}, the estimate  \eqref{march18eqn31}in Lemma \ref{conservationlawlemma},    the volume of support of $v,\omega$, and the estimate   \eqref{nov24eqn41}  if $|  v_{\bot}|\geq 2^{(\alpha_t +\epsilon)M_t}$,   we have 
\be
\begin{split}
\big|\mathbf{P}\big(  H_{k,j;n,l,r; p,q;ess}^{\mu,m,i;non,1}(t, x,\zeta) \big) \big| & \lesssim \| \mathfrak{m}(\cdot, \zeta)\|_{\mathcal{S}^\infty}    (2^{m }+2^{-k}) 2^{-j-r+ \max\{  r,   (\gamma_1-\gamma_2)M_{t^{\star}}\}+2\epsilon M_t} \\
&\quad\times  \big( \min\{2^{-2m}, 2^{m+3k+2n}\} \min\{2^{3j+2r}, 2^{j+2\alpha_t M_t }\} \big)^{1/2}\\
&\quad\times \big( \min\{2^{m+3k+2n+2l-j}, 2^{m+2r}\min\{2^{3j+2l}, 2^{j+2\alpha_t M_t }\}    \big)^{1/2}   \\
& \lesssim  \| \mathfrak{m}(\cdot, \zeta)\|_{\mathcal{S}^\infty}2^{\max\{  r,   (\gamma_1-\gamma_2)M_{t^{\star}}\}+10\epsilon M_t} \\
&\quad \times  \min\big\{  2^{k+n}, |  x_{\bot} |^{-1/2}2^{-k}  \min\{2^{2j}, 2^{ 2\alpha_t M_t -2 r} \}\big\}\\
& \lesssim   \| \mathfrak{m}(\cdot, \zeta)\|_{\mathcal{S}^\infty}\big[ |  x_{\bot}|^{-1/4 }  2^{(\gamma_1-\gamma_2)M_{t^{\star} }/4+(\alpha^{\star} -11\epsilon )M_{t^{\star}}  }\\
&\quad +  |  x_{\bot}|^{-1/8}  2^{ (\gamma_1-\gamma_2)M_{t^{\star} }/8   -\epsilon  M_{t^{\star}} } \mathbf{1}_{n\geq  ( \gamma_1-\gamma_2- 200\epsilon )M_{t^{\star} }/2 }\\
&\quad \times \min\{2^{(k+2n)/2+ 5\alpha^{\star}M_{t^\star}/8 }, 2^{(k+4n)/2+  \alpha^{\star}M_{t^\star}}\}     \big],\\
\end{split}
\ee
By using the same strategy, we have 
\be
\begin{split}
\big|\mathbf{P}\big(  H_{k,j;n,l,r; p,q;ess}^{\mu,m,i;non,2}(t, x,\zeta) \big) \big|& \lesssim     \| \mathfrak{m}(\cdot, \zeta)\|_{\mathcal{S}^\infty}  (2^{m }+2^{-k}) 2^{-j-r+l+ \max\{  r,   (\gamma_1-\gamma_2)M_{t^{\star}}\}+2\epsilon M_t} \\
&\quad \times \big( 2^{m+3k+2n+2l} \min\{2^{3j+2r }, 2^{j+2\alpha_t M_t }\} \big)^{1/2}\\
&\quad \times \big( \min\{2^{m+3k+2n+2l-j}, 2^{m+2r}\min\{2^{3j+2l}, 2^{j+2\alpha_t M_t }\}    \big)^{1/2}  \\
&  \lesssim   \| \mathfrak{m}(\cdot, \zeta)\|_{\mathcal{S}^\infty}\big[ |  x_{\bot}|^{-1/4 }  2^{(\gamma_1-\gamma_2)M_{t^{\star} }/4+(\alpha^{\star} -11\epsilon )M_{t^{\star}}  }\\
&\quad+  |  x_{\bot}|^{-1/8}  2^{ (\gamma_1-\gamma_2)M_{t^{\star} }/8   -\epsilon  M_{t^{\star}} }   \mathbf{1}_{n\geq  ( \gamma_1-\gamma_2- 200\epsilon )M_{t^{\star} }/2 } \\
&\quad \times  \min\{2^{(k+2n)/2+ 5\alpha^{\star}M_{t^\star}/8 }, 2^{(k+4n)/2+  \alpha^{\star}M_{t^\star}}\} \big]. 
\end{split}
\ee

\medskip

\textbf{Step 2C.}\quad  If   $ | x_{\bot}|\geq 2^{\epsilon M_t}\max\{2^{m}, 2^{-k-n}\} $. 
\medskip

 Note that, for any $y\in B(0, 2^{-k-n+\epsilon M_t/2}), $ we have $|  x_{\bot} -   y_{\bot} +(t-s)  \omega_{\bot} |\sim |    x_{\bot} |.$  From the   estimate of kernels in \eqref{2022feb16eqn41}, the cylindrical symmetry of the distribution function,  the estimate  \eqref{march18eqn31}  in Lemma \ref{conservationlawlemma},    the volume of support of $v,\omega$, and the estimate   \eqref{nov24eqn41}  if $|  v_{\bot}|\geq 2^{(\alpha_t +\epsilon)M_t}$,   we have 
\be
\begin{split}
&\big|\mathbf{P}\big(  H_{k,j;n,l,r; p,q;ess}^{\mu,m,i;non,1}(t, x,\zeta) \big) \big| \\
 & \lesssim  \| \mathfrak{m}(\cdot, \zeta)\|_{\mathcal{S}^\infty}     (2^{m }+2^{-k}) 2^{-j-r+ \max\{  r,   (\gamma_1-\gamma_2)M_{t^{\star}}\}+2\epsilon M_t}\\
&\quad \times  \big( \min\{2^{-2m}, 2^{m+3k+2n}\} \min\{2^{3j+2r }, 2^{j+2\alpha_t M_t }\} \big)^{1/2}\\
&\quad \times \big( \min\{|   x_{\bot}|^{-1} 2^{m+2k+ n+2l-j}, 2^{m+2r}\min\{2^{3j+2l}, 2^{j+2\alpha_t M_t }\}    \big)^{1/2} \\
& \lesssim \| \mathfrak{m}(\cdot, \zeta)\|_{\mathcal{S}^\infty}2^{-j-r+\max\{  r,   (\gamma_1-\gamma_2)M_{t^{\star}}\}+8\epsilon M_t}\\
&\quad \times \min\big\{\big(2^{-2m+3j+2r} \big)^{1/2}\big( |  x_{\bot}|^{-1} 2^{k+n-j} \big)^{1/2}, \big( |  x_{\bot}|^{-1} 2^{k+n -j} \big)^{1/4}\\
&\quad \times  \big(  2^{-k+ 2r+3j} \big)^{1/4}  \big(2^{-2m }\min\{2^{3j+2r }, 2^{j+2\alpha_t M_t }\}  \big)^{1/2} \big\}\\
&\lesssim  \| \mathfrak{m}(\cdot, \zeta)\|_{\mathcal{S}^\infty}2^{\max\{  r,   (\gamma_1-\gamma_2)M_{t^{\star}}\}+8\epsilon M_t}\\
&\quad \times \min\big\{  | x_{\bot} |^{-1/2} 2^{(k+n)/2}, |  x_{\bot} |^{-1/4}2^{n/4}  \min\{2^{ j +r/2}, 2^{  \alpha_t M_t -  r/2} \}\big\}\\
&  \lesssim   \| \mathfrak{m}(\cdot, \zeta)\|_{\mathcal{S}^\infty}\big[ |   x_{\bot}|^{-1/4 }  2^{(\gamma_1-\gamma_2)M_{t^{\star} }/4+(\alpha^{\star} -11\epsilon )M_{t^{\star}}   }+  |  x_{\bot}|^{-1/2}  2^{ (\gamma_1-\gamma_2)M_{t^{\star} }/2   -\epsilon  M_{t^{\star}} }   \\
& \quad \times  \mathbf{1}_{n\geq  ( \gamma_1-\gamma_2- 100\epsilon )M_{t^{\star} }   }   \min\{2^{(k+2n)/2+ \alpha^{\star}M_{t^\star}/2 }, 2^{(k+2n)/2+  \alpha^{\star}M_{t^\star}}\} \big],
\end{split}
\ee

By using the same strategy, we have 
\be\label{2022feb16eqn54}
\begin{split}
&\big|\mathbf{P}\big(  H_{k,j;n,l,r; p,q;ess}^{\mu,m,i;non,2}(t, x,\zeta) \big) \big| \\
& \lesssim    \| \mathfrak{m}(\cdot, \zeta)\|_{\mathcal{S}^\infty}   (2^{m }+2^{-k}) 2^{-j-r+l+ \max\{  r,   (\gamma_1-\gamma_2)M_{t^{\star}}\}+2\epsilon M_t}\\
&\quad \times  \big( 2^{m+3k+2n+2l}\min\{2^{3j+2r }, 2^{j+2\alpha_t M_t }\} \big)^{1/2}\\
&\quad  \times \big( \min\{|  x_{\bot}|^{-1} 2^{m+2k+ n+2l-j}, 2^{m+2r}\min\{2^{3j+2l}, 2^{j+2\alpha_t M_t }\}    \big)^{1/2} \\
& \lesssim \| \mathfrak{m}(\cdot, \zeta)\|_{\mathcal{S}^\infty}2^{-j-r+\max\{  r,   (\gamma_1-\gamma_2)M_{t^{\star}}\}+8\epsilon M_t}
\\
& \quad \times \min\big\{\big(2^{-2m+3j+2r} \big)^{1/2}\big( |  x_{\bot}|^{-1} 2^{k+n-j} \big)^{1/2},  \big( |  x_{\bot}|^{-1} 2^{k+n -j} \big)^{1/4}\\
&\quad  \times \big(  2^{-k+ 2r+3j} \big)^{1/4} \big(2^{-2m }\min\{2^{3j+2r }, 2^{j+2\alpha_t M_t }\}  \big)^{1/2} \big\}
\\
&  \lesssim \| \mathfrak{m}(\cdot, \zeta)\|_{\mathcal{S}^\infty} \big[ |  x_{\bot}|^{-1/4 }  2^{(\gamma_1-\gamma_2)M_{t^{\star} }/4+(\alpha^{\star} -11\epsilon )M_{t^{\star}}  }+  | x_{\bot}|^{-1/2}  2^{ (\gamma_1-\gamma_2)M_{t^{\star} }/2 -\epsilon  M_{t^{\star}}   } \\
& \quad \times  \min\{2^{(k+2n)/2+ \alpha^{\star}M_{t^\star}/2 }, 2^{(k+2n)/2+  \alpha^{\star}M_{t^\star}}\}   \mathbf{1}_{n\geq  ( \gamma_1-\gamma_2- 100\epsilon )M_{t^{\star} }   }\big]  .
\end{split} 
\ee
 
 To sum up, from the decomposition in  \eqref{sep24eqn46}, our desired estimate  \eqref{2022feb16eqn53}  holds after combining the obtained estimates \eqref{2022jan8eqn51}--\eqref{2022feb16eqn54}.
\end{proof}

\begin{lemma}\label{smallkfull2}
Let  $ i\in\{3,4\}, ( l,r)\in \mathcal{B}_i $, see  \eqref{sep18eqn50}.    Under the assumption of   the part (ii) in Theorem \ref{mainresultsfirstpart},   the following estimate holds if  $m+k\leq \min\{- l -n,-2l\} + \epsilon M_{t^{\star}}/4, $
\be\label{2022feb16eqn58}
\begin{split}
  &\big|\mathbf{P}\big( T_{k,j;n,l,r}^{\mu,m, i}(\mathfrak{m},E)(t,x, \zeta ) + \hat{\zeta}\times T_{k,j;n,l,r}^{\mu,m, i}(\mathfrak{m},B)(t,x, \zeta) \big) \big|\\
  & \lesssim 
\min_{b\in\{1,2\}}  \big[\sum_{a\in \{ 1/4,3/8, 1/2\}}|  x_{\bot}|^{-a}2^{a(\gamma_1-\gamma_2)M_{t^{\star}}} 2^{( {\alpha}^{\star }-10\epsilon)M_{t^{\star}}} +  \mathbf{1}_{|  x_{\bot}|\leq 2^{-k-n+b\epsilon M_{t^{\star}}}}   \\
 &  \quad \times  |  x_{\bot}|^{-1}2^{(\gamma_1-\gamma_2)M_{t^{\star}}}  2^{ 5\alpha^{\star} M_{t^{\star}}/6}    +  |  x_{\bot}|^{-a} 2^{a(\gamma_1-\gamma_2)M_{t^{\star}}  }\mathbf{1}_{n\geq    - \alpha^{\star} M_{t^{\star}}/2  } \mathbf{1}_{|  x_{\bot} |\geq 2^{-k-n+b\epsilon M_{t^{\star}}}}\\
 &  \quad  \times    \min\{2^{(k+2n)/2 + 2\alpha^{\star} M_{t^{\star}}/ 3 - (\gamma_1-\gamma_2)M_t/24}, 2^{(k+4n)/2 + \alpha^{\star} M_{t^{\star}}}\} \big]\| \mathfrak{m}(\cdot, \zeta)\|_{\mathcal{S}^\infty}. 
 \end{split}
 \ee 
\end{lemma}
\begin{proof}
Recall \eqref{2022feb15eqn1}. Due to the  cutoff functions $ \varphi_{j,n}^{i; r}(v, \zeta)  \varphi_{l; r}(\tilde{v}+\omega ), i\in\{3,4\} $, see  \eqref{sep17eqn63}  and   \eqref{sep4eqn6}, we have $ r=-j, |\tilde{v}-\tilde{\zeta}|\lesssim  2^{n+\epsilon M_t}.$  
Recall  \eqref{2022feb8eqn11}. From   the   result of computation  of symbols  in  \eqref{sep5eqn101}  and  \eqref{sep5eqn96},  we have
\be\label{2022feb19eqn51}
\begin{split}
|\mathbf{P} \big( m_{E}^{0}(\xi, v) + \hat{\zeta}\times m_{B}^{0}(\xi, v) \big)| &  \lesssim  2^{  \max\{n , (\gamma_1-\gamma_2)M_{t^{\star}}\}+ n  +  c(m,k,l)+2\epsilon M_t}, \\ 
 |\mathbf{P} \big( m_{E}^{1}(\xi, v) + \hat{\zeta}\times m_{B}^{1}(\xi, v) \big)| &  \lesssim  2^{ \max\{  n, (\gamma_1-\gamma_2)M_{t^{\star}}\} +2\epsilon M_t}.  
\end{split}
\ee
Recall  \eqref{2022feb8eqn11}.  From the above estimate,  after doing integration by parts in $\xi$ along  $\zeta$ direction  and directions perpendicular to   $\zeta$ many times, we have 
\be\label{2022jan6eqn44}
\begin{split}
   \big| \mathbf{P} \big( \mathcal{K}_{k,n,l}^{ess;\mu,m}(y, v, \omega, \zeta)\big)  \big| &+ 2^k \big|\mathbf{P} \big( \mathcal{K}_{k,n,l}^{ess;err;\mu,m}(y, v, \omega, \zeta)\big) \big|  \\
 & \lesssim\| \mathfrak{m}(\cdot,  \zeta)\|_{\mathcal{S}^\infty}  2^{4k+ 3n+ c(m,k,l)+2\epsilon M_t}   2^{\max\{n, (\gamma_1-\gamma_2)M_{t^{\star}}\}}   \\
 & \quad \times (1+2^k|y\cdot \tilde{\zeta}|)^{ -N_0^3}(1+2^{k+n}|y\times \tilde{\zeta}|)^{- N_0^3}, \\
   \big| \nabla_v \mathbf{P} \big({\mathcal{K}}_{k,n,l}^{ non; \mu,m}\big)(y, v, \omega, \zeta)\big|   &   + 2^{ k} \big| \nabla_v  \mathbf{P} \big({\mathcal{K}}_{k,n,l}^{ non;err;\mu,m}\big)(y, v, \omega, \zeta)\big|\\ 
  &\lesssim  \| \mathfrak{m}(\cdot,  \zeta)\|_{\mathcal{S}^\infty} 2^{3k+ 2n +\max\{n, (\gamma_1-\gamma_2)M_{t^{\star}}\}+2\epsilon M_t} 2^{-j-\min\{l,n\}} \\
  &\quad \times     (1+2^k|y\cdot \tilde{\zeta}|)^{ -N_0^3}(1+2^{k+n}|y\times \tilde{\zeta}|)^{- N_0^3}. 
\end{split}
\ee

 Moreover, if $n\leq (\gamma_1-\gamma_2-2\epsilon )M_{t^\star} $, then we have $|  v_{\bot} |/|v|\sim 2^{ (\gamma_1-\gamma_2  )M_{t^\star} +j}. $ From the estimate  \eqref{nov24eqn41}, we can rule out the case $  (\gamma_1-\gamma_2  )M_{t^\star}  +j\geq  (\alpha_t+\epsilon) M_t$ if $r\leq(\gamma_1-\gamma_2-2\epsilon )M_{t^\star}$. Hence, it suffices to consider the case  $  (\gamma_1-\gamma_2  )M_{t^\star} +j\leq  (\alpha_t+\epsilon) M_t$ if $r\leq  (\gamma_1-\gamma_2-2\epsilon )M_{t^\star} $. 

 Based on the possible size of $| x_{\bot}|$, we proceed in three steps as follows. 

\medskip

\noindent \textbf{Step 1.}\quad   If $|  x_{\bot} |\leq 2^{-k-n+2\epsilon M_{t^{\star}}}$. 

\medskip

Recall  \eqref{2022feb15eqn1}. From the estimate of kernels in  \eqref{2022jan6eqn44}, and the volume of support of $\omega, v$,   we have 
\be\label{2022jan16eqn79}
\begin{split}
&\big|\mathbf{P}\big(  H_{k,j;n,l,r; p,q;ess}^{\mu,m,i;lin}(t, x,\zeta) \big)\big| \\ 
& \lesssim \| \mathfrak{m}(\cdot, \zeta)\|_{\mathcal{S}^\infty} 2^{ n+c(m,k,l)+4\epsilon M_{t^{\star}}} 2^{\max\{n, (\gamma_1-\gamma_2)M_{t^{\star}}\}}   (2^{m+k}+1)\\
&\quad \times   \min\{2^{m+3k+2n+2l-j},2^{m+2l}\min\{2^{j+2\alpha^{\star} M_{t^{\star}}}, 2^{3j+2n}\}\}\\ 
& \lesssim 2^{n+\max\{n, (\gamma_1-\gamma_2)M_{t^{\star}}\} +8\epsilon M_{t^{\star}}} \min\{2^{2k+n-j},2^{-k-n }\min\{2^{j+2\alpha^{\star}M_t}, 2^{3j+2n}\}\}\| \mathfrak{m}(\cdot, \zeta)\|_{\mathcal{S}^\infty} \\
&\lesssim 2^{n+\max\{n, (\gamma_1-\gamma_2)M_{t^{\star}}\} +8\epsilon M_{t^{\star}}} \big(2^{2k+n-j}\big)^{2/3}\big(2^{-k-n } 2^{2j+\alpha^{\star}M_t+n} \big)^{1/3}\| \mathfrak{m}(\cdot, \zeta)\|_{\mathcal{S}^\infty} \\
 &\lesssim 2^{10\epsilon M_{t^{\star}}} |  x_{\bot}|^{-1}2^{\alpha^{\star}M_{t^{\star}}/3} \| \mathfrak{m}(\cdot, \zeta)\|_{\mathcal{S}^\infty} \\
 & \lesssim  |  x_{\bot}|^{-1}2^{(\gamma_1-\gamma_2)M_{t^{\star}}} 2^{(5\alpha^{\star}/6-\epsilon) M_{t^{\star}}} \| \mathfrak{m}(\cdot, \zeta)\|_{\mathcal{S}^\infty} .
 \end{split}
\ee

Now, we estimate the nonlinear part. Recall  \eqref{2022feb15eqn2}. From the Cauchy-Schwarz inequality, the estimate of kernels in  \eqref{2022jan6eqn44},  the estimate  \eqref{march18eqn31}  in Lemma \ref{conservationlawlemma},   and the volume of support of $v, \omega$,  we have
\be\label{2022jan8eqn22}
\begin{split}
&\big| \mathbf{P}\big(H_{k,j;n,l,r; p,q;ess}^{\mu,m,i;non,1}(t, x,\zeta)\big)\big| \\ & \lesssim  \| \mathfrak{m}(\cdot, \zeta)\|_{\mathcal{S}^\infty} 2^{\max\{n, (\gamma_1-\gamma_2)M_{t^{\star}}\}+4\epsilon M_t}  (2^{m}+2^{-k}) 2^{-j- \min\{l, n \}  } \\
 &\quad \times  \big(\min\{ 2^{-2m}, 2^{m+3k+2n}\} 2^{3j+2\min\{l,n\}} \big)^{1/2}  \big( 2^{m+3k+2n+2l -j} \big)^{1/2} \\
& \lesssim 2^{k+n+10\epsilon M_{t^{\star}}} \| \mathfrak{m}(\cdot, \zeta)\|_{\mathcal{S}^\infty} \\
&\lesssim |  x_{\bot}|^{-1}2^{(\gamma_1-\gamma_2)M_{t^{\star}}} 2^{(  \alpha^{\star}/2-\epsilon) M_{t^{\star}}} \| \mathfrak{m}(\cdot, \zeta)\|_{\mathcal{S}^\infty}.\\
\end{split}
 \ee

Now, we estimate the second part  of the  nonlinearity.   From the Cauchy-Schwarz inequality, the estimate of kernels in  \eqref{2022jan6eqn44},   the conservation law  \eqref{conservationlaw}  for the electromagnetic field,    and the volume of support of $v, \omega$, we have
\be
\begin{split}
  &\big| \mathbf{P}\big(H_{k,j;n,l,r; p,q;ess}^{\mu,m,i;non,2}(t, x,\zeta)\big)\big|\\ &  \lesssim \| \mathfrak{m}(\cdot, \zeta)\|_{\mathcal{S}^\infty} 2^{\max\{n, (\gamma_1-\gamma_2)M_{t^{\star}}\}+4\epsilon M_{t^{\star}}}  (2^{m}+2^{-k})   2^{l-j- \min\{l, n \}  }
\\
&\quad \times \big(2^{ m+3k+ 2n   }  2^{2l}  2^{3j+2n} \big)^{1/2} \big( 2^{ m+3k+ 2n   }  2^{2l}  2^{-j}  \big)^{1/2}\\
& \lesssim  (2^{m}+2^{-k})   2^{ l+n-\min\{l,n\}+10\epsilon M_{t^{\star}} }2^{2k+n+\min\{l,n\}}\| \mathfrak{m}(\cdot, \zeta)\|_{\mathcal{S}^\infty}\\
& \lesssim 2^{k+n+14\epsilon M_{t^{\star}}} \| \mathfrak{m}(\cdot, \zeta)\|_{\mathcal{S}^\infty} \\
& \lesssim |  x_{\bot}|^{-1}2^{(\gamma_1-\gamma_2)M_{t^{\star}}} 2^{(  \alpha^{\star}/2-\epsilon) M_{t^{\star}}}   \| \mathfrak{m}(\cdot, \zeta)\|_{\mathcal{S}^\infty}.
\end{split}
\ee

\medskip

\noindent \textbf{Step 2.}\quad    If $ |  x_{\bot}| \geq 2^{-k-n+ \epsilon M_{t^{\star}}}  $   and $| x_{\bot} |\geq 2^{m+p+\epsilon M_{t^{\star}} /2}$. 

\medskip

For this case, we have $\forall y \in B(0, 2^{-k-n+\epsilon  M_{t^{\star}}/2 }), |  x_{\bot} -   y_{\bot} + (t-s) \omega_{\bot}|\sim |  x_{\bot}|$.  From the estimate of kernels in  \eqref{2022jan6eqn44},   the estimate  \eqref{march18eqn31}  in Lemma \ref{conservationlawlemma}, the cylindrical symmetry of the solution, and the volume of support of $\omega, v$,    we have
\be
\begin{split}
&\big|\mathbf{P}\big(  H_{k,j;n,l,r; p,q;ess }^{\mu,m,i;lin}(t, x,\zeta) \big)\big| \\
& \lesssim \| \mathfrak{m}(\cdot, \zeta)\|_{\mathcal{S}^\infty}  2^{ n+c(m,k,l)+4\epsilon M_{t^{\star}}} 2^{\max\{n, (\gamma_1-\gamma_2)M_{t^{\star}}\}}   (2^{m+k}+1) \\
&  \quad \times   \min\{2^{m+2k+n+2l-j}| x_{\bot}|^{-1},2^{m+2l+j+2\alpha^{\star} M_{t^{\star}} },  2^{m+2n+3j+2l}, 2^{m+2l} 2^{3k+2n-j} \}\\
& \lesssim 2^{n+\max\{n, (\gamma_1-\gamma_2)M_{t^{\star}}\} +8\epsilon M_{t^{\star}}} \| \mathfrak{m}(\cdot, \zeta)\|_{\mathcal{S}^\infty}\\
&\quad \times \min\big\{| x_{\bot}|^{-1} 2^{k-j},2^{-k-n } 2^{j+2\alpha^{\star} M_{t^{\star}} },   2^{j/3+4\alpha^{\star} M_{t^{\star}} /3-n/3}, 2^{5j/3+ n  }\big\}\\
& \lesssim \| \mathfrak{m}(\cdot, \zeta)\|_{\mathcal{S}^\infty}2^{n+\max\{n, (\gamma_1-\gamma_2)M_{t^{\star}}\} +8\epsilon M_{t^{\star}}}  \min\big\{\big( |  x_{\bot}|^{-1} 2^{k-j}\big)^{1/2}\big( 2^{-k-n } 2^{j+2\alpha^{\star} M_{t^{\star}} }\big)^{1/2}, \\
& \quad  \big( |  x_{\bot}|^{-1} 2^{k-j}\big)^{1/2}\big( 2^{j/3+4\alpha^{\star} M_{t^{\star}} /3-n/3}\big)^{1/4}\big(2^{5j/3+ n  }\big)^{1/4}\big\}\\
 & \lesssim  |  x_{\bot} |^{-1/2} 2^{(\gamma_1-\gamma_2)M_{t^{\star}}/2-\epsilon M_{t^{\star}} } \big[2^{(\alpha^{\star}-10\epsilon)M_t}  +  2^{(k+2n)/2 + 7\alpha^{\star}M_t/12} \mathbf{1}_{n\geq - \alpha^{\star}M_t/4} \big]\| \mathfrak{m}(\cdot, \zeta)\|_{\mathcal{S}^\infty}.
  \end{split}
\ee

Now, we estimate the first part of the  nonlinearity. From the estimate of kernels in  \eqref{2022jan6eqn44},   the cylindrical symmetry of the distribution function,    and  the volume of support of $v, \omega$, we have
\be 
\begin{split}
 &\big| \mathbf{P}\big(H_{k,j;n,l,r; p,q;ess}^{\mu,m,i;non,1}(t, x,\zeta)\big)\big| \\
 & \lesssim \| \mathfrak{m}(\cdot, \zeta)\|_{\mathcal{S}^\infty}  2^{\max\{n, (\gamma_1-\gamma_2)M_{t^{\star}}\}+5\epsilon M_t}  (2^{m}+2^{-k}) 2^{-j- \min\{l, n \} } \\
 &\quad \times \big( \min\{ 2^{-2m}, 2^{m+3k+2n}\}  2^{  3j+2\min\{l,n\}} \big)^{1/2}\\
 & \quad \times \big(\min\{ |  x_{\bot}|^{-1}2^{m+2k+n+2l -j},   2^{m+2l}  2^{3j+2n} \} \big)^{1/2} \\
&\lesssim \| \mathfrak{m}(\cdot, \zeta)\|_{\mathcal{S}^\infty} 2^{\max\{n, (\gamma_1-\gamma_2)M_{t^{\star}}\}+  6\epsilon M_{t^{\star}}  } \min\{ |  x_{\bot}|^{-1/2} 2^{(k+n)/2}, | x_{\bot}|^{-1/4} 2^{j+n/4} \}\\
&\lesssim \| \mathfrak{m}(\cdot, \zeta)\|_{\mathcal{S}^\infty} \big[|  x_{\bot}|^{-1/4} 2^{(\gamma_1-\gamma_2)M_{t^{\star}}/4} 2^{(\alpha^{\star}-11\epsilon)M_{t^{\star}}}  + \mathbf{1}_{n\geq (\gamma_1-\gamma_2-80\epsilon)M_{t^{\star}}}\\
&\quad  \times  | x_{\bot}|^{-1/2} 2^{(\gamma_1-\gamma_2)M_{t^{\star}}/2-\epsilon M_{t^{\star}}}  2^{(k+4n)/2+\alpha^{\star} M_{t^{\star}}/2}\big].\\
  \end{split}
\ee

Now, we estimate the second part of the  nonlinearity $\mathbf{P}\big(H_{k,j;n,l,r; p,q;ess}^{\mu,m,i;non,2}(t, x,\zeta)\big)$. From the cylindrical symmetry of the distribution function,  the Jacobian of changing coordinates $(\theta, \phi)\longrightarrow (z,w)$  in  \eqref{march18eqn66}, the estimate of kernels in \eqref{2022jan6eqn44}  and  the volume of support of $v, \omega$, we have 
 \be
   \begin{split}
&\big| \mathbf{P}\big(H_{k,j;n,l,r; p,q;ess}^{\mu,m,i;non,2}(t, x,\zeta)\big)\big| \\
& \lesssim \| \mathfrak{m}(\cdot, \zeta)\|_{\mathcal{S}^\infty}  2^{\max\{n, (\gamma_1-\gamma_2)M_{t^{\star}}\}+6\epsilon M_t}  (2^{m}+2^{-k}) 2^l 2^{-j- \min\{l, n \}  }\\
&\quad\times \big(\min\{ |  x_{\bot} |^{-1} {2^{-m-p-q  }} ,|  x_{\bot}|^{-1} 2^{m+2k+n+2p+q} \} 2^{3j+2\min\{l,n\}} \big)^{1/2} \\
&\quad \times  \big(\min\{     2^{m+2l}  2^{3j+2 n}, 2^{m+2l} 2^{3k+2n-j}, 2^{-2m-j-2l}\} \big)^{1/2}\\
&\lesssim \| \mathfrak{m}(\cdot, \zeta)\|_{\mathcal{S}^\infty} 2^{\max\{n, (\gamma_1-\gamma_2)M_{t^{\star}}\}+6\epsilon M_t}      |  x_{\bot}|^{-1/2} 2^{(2k+n)/4+l+j/2}\\
&\quad \times \big(\min\{ 2^{-j-2l}, ((2^{m}+2^{-k}) 2^{m+2l+3j+2n})^{1/2}(2^{-j-2l})^{1/2}     \} \big)^{1/2}  \\
&   \lesssim 2^{\max\{n, (\gamma_1-\gamma_2)M_{t^{\star}}\}+  10\epsilon M_{t^{\star}}    }  |  x_{\bot}|^{-1/2} \min\{2^{k/2+ n/4}, 2^{j+n/2}  \}\| \mathfrak{m}(\cdot, \zeta)\|_{\mathcal{S}^\infty}\\
&\lesssim \| \mathfrak{m}(\cdot, \zeta)\|_{\mathcal{S}^\infty} |  x_{\bot}|^{-1/2} 2^{(\gamma_1-\gamma_2)M_{t^{\star}}/2-\epsilon M_{t^{\star}}} \\
&\quad \times \big[ 2^{(\alpha^{\star}-10\epsilon)M_{t^{\star}}}  +  2^{(k+4n)/2+5\alpha^{\star} M_{t^{\star}}/8}\mathbf{1}_{n\geq (\gamma_1-\gamma_2-80\epsilon)M_{t^{\star}}}\big].
  \end{split}
\ee

\medskip

\noindent \textbf{Step 3.}\quad   If $ | x_{\bot}| \geq 2^{-k-n+ \epsilon M_{t^{\star}}}  $   and $|  x_{\bot} |\leq  2^{m+p+\epsilon M_{t^{\star}}/2}$. 

\medskip

For this sub-case, we have $p\geq l +\epsilon  M_{t^{\star}}/4$ and $|  v_{\bot} |/|v|\sim 2^p$. From the estimate \eqref{nov24eqn41}, we know that  it suffices to consider the case $j+p\leq (\alpha_t+\epsilon)M_t+10.$  Moreover,  $\forall y\in B(0, 2^{-k-n+\epsilon M_{t^{\star}}/2}), |  x_{\bot} -   y_{\bot} |\sim |  x_{\bot}|$.  From   the estimate of kernels in  \eqref{2022jan6eqn44},  the volume of support of $\omega, v$, and the Jacobian of changing coordinates $(\theta, \phi)\longrightarrow (z,w)$  in  \eqref{march18eqn66}, the following estimate holds, 
 \be
   \begin{split}
&\big|\mathbf{P}\big(  H_{k,j;n,l,r; p,q;ess}^{\mu,m,i;lin}(t, x,\zeta) \big)\big|\\
&\lesssim   \| \mathfrak{m}(\cdot, \zeta)\|_{\mathcal{S}^\infty} 2^{\max\{n, (\gamma_1-\gamma_2)M_{t^{\star}}\}}  2^{ n+c(m,k,l)} (2^{m+k}+1) \\
 &\quad \times \min\{\frac{2^{-m-p-q  }}{|  x_{\bot} |} 2^{-j} ,   2^{m+3j+2l +2p+q}, 2^{m+3k+2n+2l-j} \}
\\ 
& \lesssim \| \mathfrak{m}(\cdot, \zeta)\|_{\mathcal{S}^\infty}   2^{n+\max\{n, (\gamma_1-\gamma_2)M_{t^{\star}}\}+5\epsilon M_{t^{\star}}}
\min\{  |  x_{\bot}|^{-1/2} 2^{j/2 +\alpha_t M_t/2   }, 2^{2k+n-j}   \}\\
&\lesssim  \| \mathfrak{m}(\cdot, \zeta)\|_{\mathcal{S}^\infty} \big[| x_{\bot}|^{-1/2} 2^{(\gamma_1-\gamma_2)M_{t^{\star}}/2} 2^{ ( \alpha^{\star} -11 \epsilon ) M_{t^{\star}}  }  \\
&\quad + |  x_{\bot}|^{-3/8}2^{3(\gamma_1-\gamma_2) M_{t^{\star} }/8-\epsilon M_{t^{\star}} }\mathbf{1}_{    n\geq (\gamma_1-\gamma_2-80\epsilon)M_{t^{\star}}}\\
&\quad \times \min\{2^{(k+2n)/2 + 2\alpha^{\star} M_{t^{\star}}/ 3 - (\gamma_1-\gamma_2)M_t/24}, 2^{(k+4n)/2 + \alpha^{\star} M_{t^{\star}}}\}\big].
   \end{split}
\ee  

  Recall  \eqref{2022feb15eqn2}. From     the estimate of kernels in  \eqref{2022jan6eqn44}, the estimate  \eqref{march18eqn31}  in Lemma \ref{conservationlawlemma},  the volume of support of $\omega, v$,   the Jacobian of changing coordinates $(\theta, \phi)\longrightarrow (z,w)$  in  \eqref{march18eqn66}, and the Jacobian of changing coordinates $(y_1, y_2, \theta)\longrightarrow x-y+(t-s)\omega$,    we have 
 \be
   \begin{split}
  &\big| \mathbf{P}\big(H_{k,j;n,l,r; p,q;ess}^{\mu,m,i;non,1}(t, x,\zeta)\big)\big| \\
  & \lesssim \| \mathfrak{m}(\cdot, \zeta)\|_{\mathcal{S}^\infty}   2^{\max\{n, (\gamma_1-\gamma_2)M_{t^{\star}}\}+2\epsilon M_t}  (2^{m}+2^{-k}) 2^{-j- \min\{l, n \} } \\
  &\quad \times \big( \min\{ |  x_{\bot} |^{-1}{2^{-m-p-q  }} ,  2^{ 2k+n +q} \} 2^{-j}, 2^{m+2p+q+3j+2n} \big)^{1/2}\\
  &\quad \times  \big( \min\{2^{-2m}, 2^{m+3k+2n}\} 2^{ 3j+2\min\{l,n\}} \big)^{1/2}\\
& \lesssim \| \mathfrak{m}(\cdot, \zeta)\|_{\mathcal{S}^\infty} 2^{\max\{n, (\gamma_1-\gamma_2)M_{t^{\star}}\}+2\epsilon M_t+j/2} \big(\min\big\{ (|  x_{\bot}|^{-2}2^{-q+\epsilon M_{t^{\star}}-j} )^{1/2} (2^{2k+n+q-j})^{1/2}, \\
 & \quad   (|  x_{\bot}|^{-1}2^{-m-p-q+\epsilon M_{t^{\star}}-j} )^{1/2} (2^{m+2p+q+3j+2n})^{1/2}  \big\}\big)^{1/2}\\
&\lesssim 2^{\max\{n, (\gamma_1-\gamma_2)M_{t^{\star}}\}+  6\epsilon M_{t^{\star}}  } \min\{ |  x_{\bot}|^{-1/2}  2^{k/2+ n/4}, |  x_{\bot}|^{-1/4} 2^{j+n/2} \}\| \mathfrak{m}(\cdot, \zeta)\|_{\mathcal{S}^\infty} \\
&\lesssim \| \mathfrak{m}(\cdot, \zeta)\|_{\mathcal{S}^\infty} |  x_{\bot}|^{-1/2} 2^{(\gamma_1-\gamma_2)M_{t^{\star}}/2-\epsilon M_{t^{\star}}}\\
&\quad \times  \big[ 2^{(\alpha^{\star}-10\epsilon)M_{t^{\star}}}  +  2^{(k+4n)/2+5\alpha^{\star} M_{t^{\star}}/8}\mathbf{1}_{n\geq (\gamma_1-\gamma_2-80\epsilon)M_{t^{\star}}}\big].
  \end{split}
\ee

  By using similar strategy, for the second part of the nonlinearity $\mathbf{P}({H}_{k,j;n,l;ess}^{\mu,m;non,2;p,q})(t, x,\zeta)$,  we have
 \be\label{2022jan8eqn52}
   \begin{split}
   & \big| \mathbf{P}\big(H_{k,j;n,l,r; p,q;ess}^{\mu,m,i;non,2}(t, x,\zeta)\big)\big| \\
    & \lesssim \| \mathfrak{m}(\cdot, \zeta)\|_{\mathcal{S}^\infty}  2^{\max\{n, (\gamma_1-\gamma_2)M_{t^{\star}}\}+2\epsilon M_t}  (2^{m}+2^{-k}) 2^l 2^{-j- \min\{l, n \} }\\
    &\quad \times \big(\min\{     2^{m+2l}  2^{3j+2 n}, 2^{m+2l} 2^{3k+2n-j},  2^{-2m-j-2l} \} \big)^{1/2}\\
&\quad\times \big(\min\{ | x_{\bot} |^{-1}{2^{-m-p-q  }} ,  2^{ 2k+n +q} \} 2^{3j+2\min\{l,n\}} \big)^{1/2}\\
&\lesssim  \| \mathfrak{m}(\cdot, \zeta)\|_{\mathcal{S}^\infty}  2^{\max\{n, (\gamma_1-\gamma_2)M_{t^{\star}}\}+  10\epsilon M_{t^{\star}}    }  |  x_{\bot}|^{-1/2} \min\{2^{k/2+ n/4}, 2^{j+n/2}  \}\\
&\lesssim \| \mathfrak{m}(\cdot, \zeta)\|_{\mathcal{S}^\infty} | x_{\bot}|^{-1/2} 2^{(\gamma_1-\gamma_2)M_{t^{\star}}/2-\epsilon M_{t^{\star}}} \\
&\quad \times  \big[ 2^{(\alpha^{\star}-10\epsilon)M_{t^{\star}}}  +  2^{(k+4n)/2+5\alpha^{\star} M_{t^{\star}}/8}\mathbf{1}_{n\geq (\gamma_1-\gamma_2-80\epsilon)M_{t^{\star}}}\big].
  \end{split}
\ee
Hence, our desired estimate  \eqref{2022feb16eqn58}  holds after combining  the obtained estimates  \eqref{2022jan16eqn79}--\eqref{2022jan8eqn52}.
\end{proof}

\begin{lemma}\label{smallprojec1}
Let  $ i\in\{3,4\},  ( l,r)\in \mathcal{B}_i $, see  \eqref{sep18eqn50}.    Under the assumption of  the part (ii) in Theorem \ref{mainresultsfirstpart},    the following estimate holds if  $  \min\{- l -n,-2l\} + \epsilon M_{t^{\star}}/4 < m+k \leq -2l+ \epsilon M_{t^{\star}}/5, $
\be\label{2022feb16eqn59}
\begin{split}
  &\big|\mathbf{P}\big( T_{k,j;n,l,r}^{\mu,m, i }(\mathfrak{m},E)(t,x, \zeta ) + \hat{\zeta}\times T_{k,j;n,l,r}^{\mu,m, i}(\mathfrak{m},B)(t,x, \zeta) \big) \big|\\
&\lesssim \min_{b\in\{1,2\}} \| \mathfrak{m}(\cdot, \zeta)\|_{\mathcal{S}^\infty}  \big[\sum_{a\in \{0,1/8,1/4,3/8,1/2\}}|  x_{\bot}|^{-a}2^{a(\gamma_1-\gamma_2)M_{t^{\star}}}  2^{( {\alpha}^{\star }-10\epsilon)M_{t^{\star}}} \\
 & \quad  +  \mathbf{1}_{|  x_{\bot}|\leq 2^{-k-n+b\epsilon M_{t^{\star}}}}     |  x_{\bot}|^{-1}2^{(\gamma_1-\gamma_2)M_{t^{\star}}}  2^{ 5\alpha^{\star} M_{t^{\star}}/6}  \\
 &\quad  +  |  x_{\bot}|^{-a} 2^{a(\gamma_1-\gamma_2)M_{t^{\star}} }\mathbf{1}_{n\geq    - \alpha^{\star} M_{t^{\star}}/2  }   \mathbf{1}_{|  x_{\bot}|\geq 2^{-k-n+b \epsilon M_{t^{\star}}}} \\
 &\quad  \times  \min\{2^{(k+2n)/2 + 2\alpha^{\star} M_{t^{\star}}/ 3 - (\gamma_1-\gamma_2)M_t/24}, 2^{(k+4n)/2 + \alpha^{\star} M_{t^{\star}}}\} \big].
 \end{split}
 \ee  
\end{lemma}
 \begin{proof}

Note that, due to the constraint of  $ \min\{- l -n,-2l\}+ \epsilon M_{t^{\star}}/4 < m+k\leq -2l+ \epsilon M_{t^{\star}}/5, $ we only need to consider the case $l\leq n- M_{t^{\star}}/20. $  Moreover, we have $  c(m,k,l)\leq \max\{-m/2-k/2, -m-k-l\}+\epsilon M_t/20 \leq n-\epsilon  M_{t^{\star}}/20$. Hence,  $|\tilde{\zeta}\times \tilde{v}|\sim |\tilde{\zeta}\times \omega |\sim |\tilde{\zeta}\times \tilde{\xi} |\sim 2^n$ if $n> -M_t.$ Moreover, due to the cutoff function $\varphi^i_{j,n}(v,\zeta), i \in\{3,4\}$, see  \eqref{sep5eqn10},   we also have $|\tilde{v}-\tilde{\zeta}|\lesssim 2^{n+\epsilon M_t}$. 

Recall  \eqref{2022feb8eqn11}. After doing integration by parts in $\xi$ along $v$ (alternatively along $\zeta$) direction  and directions perpendicular to $v$ (alternatively to $\zeta$) many times, from the  estimate of   symbol in   \eqref{2022feb19eqn51}, we have 
\be\label{sep22eqn61} 
\begin{split}
&2^{-k-n-c(m,k,l)}\big(  \big| \mathbf{P} \big( \mathcal{K}_{k,n,l,r}^{ess;\mu,m}(y, v, \omega, \zeta)\big)  \big| + 2^k \big|\mathbf{P} \big( \mathcal{K}_{k,n,l,r}^{ess;err;\mu,m}(y, v, \omega, \zeta)\big) \big| \big) \\
 &+ 2^{j+l }\big(  \big| \nabla_v \mathbf{P}\big( \mathfrak{H}^{ess; \mu,E,i}_{k,j;n,l,r}(y,\omega, v, \zeta) \big) \big|+ 2^k|\nabla_v  \mathbf{P}\big(\mathfrak{H}^{ess;  err;\mu,E,i}_{k,j;n,l,r}(y,\omega, v, \zeta)\big) | \big) \\
 &  \lesssim 2^{3k+ 2c(m,k,l)+  \max\{n, (\gamma_1-\gamma_2)M_{t^{\star}}\}+2\epsilon M_t}  \| \mathfrak{m}(\cdot,  \zeta)\|_{\mathcal{S}^\infty}\\
 &\quad  \times  \min\big\{  (1+2^{k+c(m,k,l)-n-\epsilon M_t}|y\cdot \tilde{\zeta}|)^{-100} (1+2^{k+ c(m,k,l) }|y\times \tilde{\zeta}|)^{-100} , \\ 
  & \qquad (1+2^k|y\cdot \tilde{v}|)^{-100}   (1+2^{k+ c(m,k,l) }|y\times \tilde{v}|)^{-100}  \big\}.
    \end{split}
\ee
Moreover, the following improved estimate holds the $\p_{v_3}$ derivative of kernels, 
\be\label{sep22eqn93} 
\begin{split}
   & \big| \p_{v_3} \mathbf{P}\big( \mathfrak{H}^{ess; \mu,E,i}_{k,j;n,l,r}(y,\omega, v, \zeta) \big) \big| + 2^k| \p_{v_3}  \mathbf{P}\big(\mathfrak{H}^{ess;  err;\mu,E,i}_{k,j;n,l,r}(y,\omega, v, \zeta)\big) |\\
   &\lesssim  \| \mathfrak{m}(\cdot, \zeta)\|_{\mathcal{S}^\infty} 2^{3k+  2c(m,k,l,n)} 2^{-j-l+\max\{l,p\}+3\epsilon M_t}   2^{\max\{n, (\gamma_1-\gamma_2)M_{t^{\star}}\}} \\
   &\quad \times  \min\big\{(1+2^{k+c(m,k,l )-n-\epsilon M_t }|y\cdot \tilde{\zeta}|)^{-100} (1+2^{k+ c(m,k,l ) }|y\times \tilde{\zeta}|)^{-100} ,\\
   & \qquad     (1+2^k|y\cdot \tilde{v}|)^{-100}  (1+2^{k+ c(m,k,l ) }|y\times \tilde{v}|)^{-100}  \big\}. 
        \end{split}
\ee
 
  Note that,  $  |  v_{\bot}|/|v|\sim 2^{ (\gamma_1-\gamma_2)M_{t^{\star}}  }   $ if $n \notin [(\gamma_1-\gamma_2-2\epsilon)M_{t^\star}, (\gamma_1-\gamma_2+2\epsilon)M_{t^\star}]$. From the estimate  \eqref{nov24eqn41},  we can rule out the case $(\gamma_1-\gamma_2)M_{t^{\star}} + j\geq (\alpha_t + \epsilon)M_t $. Therefore, it would  be sufficient to consider the case $(\gamma_1-\gamma_2)M_{t^{\star}} + j\leq (\alpha_t + \epsilon)M_t $ if $n    \notin [(\gamma_1-\gamma_2-2\epsilon)M_{t^\star}, (\gamma_1-\gamma_2+2\epsilon)M_{t^\star}]$.

Based on the possible size of $| x_{\bot}|$ and $p$, we proceed in four steps    as follows. 

\medskip

\noindent \textbf{Step 1.}\quad   If $|  x_{\bot}|\geq 2^{m+l+4\epsilon M_{t^{ }}}$.

\medskip

For this subcase, we have $\forall y\in B(0, 2^{-k-c(m,k,l)+\epsilon M_{t^{  }}/10 }),|  x_{\bot} -   y_{\bot} |\sim |  x_{\bot}|$. From   the estimate of kernels in  \eqref{sep22eqn61},  the volume of support of $\omega, v$, the estimate \eqref{nov24eqn41}  if $|  v_{\bot}|\geq 2^{(\alpha_t + \epsilon)M_t}$, and the Jacobian of changing coordinates $(\theta, \phi)\longrightarrow (z,w)$  in  \eqref{march18eqn66},  we have 
\be\label{sep23eqn100}
\begin{split}
 &\big|\mathbf{P}\big(  H_{k,j;n,l,r; p,q;ess}^{\mu,m,i;lin}(t, x,\zeta) \big)\big|\\
 & \lesssim  \| \mathfrak{m}(\cdot, \zeta)\|_{\mathcal{S}^\infty}    2^{ m+k+ n+c(m,k,l )}2^{\max\{n, (\gamma_1-\gamma_2)M_{t^{\star}}\}+6\epsilon M_{t^{\star}}}   \\
  &\quad \times \min\{   2^{m+3k+2c(m,k,l )+2l-j}, |  x_{\bot} |^{-1} {2^{-m-p-q-j}} ,  2^{m+2p+q}\min\{2^{3j+2l},2^{j+2\alpha_t M_t}\}  \}
\\
 &\lesssim \| \mathfrak{m}(\cdot, \zeta)\|_{\mathcal{S}^\infty}     2^{n+\max\{n, (\gamma_1-\gamma_2)M_{t^{\star}}\}+8\epsilon M_{t^{\star}} } \min\big\{ 2^{2k+n-j}, 2^{-l}\big(|  x_{\bot} |^{-1} {2^{-m-p-q-j}} \big)^{1/2}\\
 &\quad \times \big(2^{m+2p+q+2j+l+\alpha_t M_t} \mathbf{1}_{p\leq l+\epsilon M_t} + 2^{m+2p+q+3j+2l} \mathbf{1}_{p\geq l+\epsilon M_t}\big)^{1/2} \big\}\\
  &\lesssim \| \mathfrak{m}(\cdot, \zeta)\|_{\mathcal{S}^\infty}  2^{n+\max\{n, (\gamma_1-\gamma_2)M_{t^{\star}}\}+10\epsilon M_{t^{\star}}} \min\{  |  x_{\bot}|^{-1/2} 2^{j/2+  {\alpha}_t M_t/2 }, 2^{2k +n -j} \}  \\
 & \lesssim  \| \mathfrak{m}(\cdot, \zeta)\|_{\mathcal{S}^\infty} \big[|  x_{\bot}|^{-1/2} 2^{(\gamma_1-\gamma_2)M_{t^{\star}}/2} 2^{ ( \alpha^{\star} -11 \epsilon ) M_{t^{\star}}  } + | x_{\bot}|^{-3/8}2^{3 (\gamma_1-\gamma_2-\epsilon)M_{t^{\star}}/8} \\
 &\quad\times \mathbf{1}_{n\geq (\gamma_1-\gamma_2-40\epsilon)M_{t^{\star}}  }
 \min\{2^{(k+2n)/2 + 2\alpha^{\star} M_{t^{\star}}/ 3 - (\gamma_1-\gamma_2)M_t/24}, 2^{(k+4n)/2 + \alpha^{\star} M_{t^{\star}}}\}\big].
 \end{split}
\ee
 We remark that we used the fact that $| v_{\bot}|/|v|\sim 2^p$ if $p\geq l+\epsilon M_t$ in the above estimate.

Now, we move on to the estimate of the nonlinear parts in   \eqref{2022feb15eqn2}.  From the Cauchy-Schwarz inequality,  the estimate  \eqref{march18eqn31}  in Lemma \ref{conservationlawlemma}, the estimate of kernels in \eqref{sep22eqn61},  the volume of support of $\omega, v$, and the Jacobian of changing coordinates $(\theta, \phi)\longrightarrow (z,w)$  in  \eqref{march18eqn66}, we have 
\be\label{sep23eqn101}
\begin{split}
  &\big| \mathbf{P}\big(H_{k,j;n,l,r; p,q;ess}^{\mu,m,i;non,1}(t, x,\zeta)\big)\big| \\
  & \lesssim \| \mathfrak{m}(\cdot, \zeta)\|_{\mathcal{S}^\infty}  2^{m-j-l}  2^{\max\{n, (\gamma_1-\gamma_2)M_{t^{\star}}\}+6\epsilon M_{t^{\star}}}     \big( 2^{  -2m } 2^{n-c(m,k,l )} 2^{3j+2l} \big)^{1/2}\\
   & \quad \times \big(  \min\{|  x_{\bot} |^{-1} {2^{-m-p-q-j}} ,   2^{m+2p+q}\min\{2^{3j+2l},2^{j+2\alpha_t M_t}\},  2^{m+3k+2c(m,k,l )+2l-j}\}\big)^{1/2} \\
&\lesssim   \| \mathfrak{m}(\cdot, \zeta)\|_{\mathcal{S}^\infty} 2^{ (n-c(m,k,l ))/2+j/2}  2^{\max\{n, (\gamma_1-\gamma_2)M_{t^{\star}}\}  +10\epsilon M_{t^{\star}}} \\
& \quad \times  \big(   \min\{  |  x_{\bot}|^{-1/2} \min\{2^{  {\alpha}_t M_t  +l/2  },  2^{j/2+ {\alpha}_t M_t/2 +l  } \},  2^{2k+l+ c(m,k,l ) -j} \}\big)^{1/2} \\
&   \lesssim \| \mathfrak{m}(\cdot, \zeta)\|_{\mathcal{S}^\infty} \big[|  x_{\bot}|^{-1/4} 2^{(\gamma_1-\gamma_2)M_{t^{\star}}/4} 2^{ ( \alpha^{\star} -11\epsilon ) M_{t^{\star}}  } + |  x_{\bot}|^{-1/8} 2^{(\gamma_1-\gamma_2) M_{t^{\star}}/8-\epsilon M_{t^{\star}}} \\ 
&  \quad \times \min\{  2^{(k+2n)/2+2\alpha^{\star}  M_{t^{\star}}/3- (\gamma_1-\gamma_2) M_{t^{\star}}/24  },2^{(k+4n)/2 + \alpha^{\star}  M_{t^{\star}}  } \}  \mathbf{1}_{n\geq (\gamma_1-\gamma_2-50\epsilon)M_{t^{\star}}  } \big]. \\
\end{split}
\ee

Lastly, we   estimate  $\mathbf{P}\big(H_{k,j;n,l,r; p,q;ess}^{\mu,m,i;non,2}(t, x,\zeta)\big)$. From the Cauchy-Schwarz inequality,    the estimate of kernels in  \eqref{sep22eqn61}  and  \eqref{sep22eqn93}, and the volume of support of $\omega, v$, we have 
\be\label{sep23eqn102}
\begin{split}
 &\big| \mathbf{P}\big(H_{k,j;n,l,r; p,q;ess}^{\mu,m,i;non,2}(t, x,\zeta)\big)\big|\\
  &\lesssim   \| \mathfrak{m}(\cdot, \zeta)\|_{\mathcal{S}^\infty}  2^{m-j } 2^{ \max\{l,p\}+\max\{n, (\gamma_1-\gamma_2)M_{t^{\star}}\}+6\epsilon M_{t^{\star}} }
 \\
 &\quad \times   \big(   2^{m+3k+ 2c(m,k,l  )} 2^{2\min\{n,p\}} 2^{3j+2l} \big)^{1/2}    \big(    2^{m+3k+2c(m,k,l )  +2l-j} \big)^{1/2} \\
&   \lesssim\| \mathfrak{m}(\cdot, \zeta)\|_{\mathcal{S}^\infty}  2^{\max\{n, (\gamma_1-\gamma_2)M_{t^{\star}}\}} 2^{k +\min\{n,p\}+ \max\{l,p\}+ 10\epsilon M_{t^{\star}}} .
\end{split}
\ee

Moreover,    from the Cauchy-Schwarz inequality,  the estimate of kernels in 
 \eqref{sep22eqn61}  and  \eqref{sep22eqn93}, the volume of support of $\omega, v$,  the estimate  \eqref{nov24eqn41}  if $|  v_{\bot}|\geq 2^{(\alpha_t + \epsilon)M_t}$, and the Jacobian of changing coordinates $(\theta, \phi)\longrightarrow (z,w)$  in \eqref{march18eqn66}, we have
\be\label{sep24eqn16}
\begin{split}
&\big| \mathbf{P}\big(H_{k,j;n,l,r; p,q;ess}^{\mu,m,i;non,2}(t, x,\zeta)\big)\big| \\
& \lesssim  \| \mathfrak{m}(\cdot, \zeta)\|_{\mathcal{S}^\infty} 2^{m-j +  \max\{l,p\}+ \max\{n, (\gamma_1-\gamma_2)M_{t^{\star}}\}+10\epsilon M_{t^{\star}}} \\
&\quad \times  \min\{2^{3j+2l},2^{j+2\alpha_t M_t}\}  \big(      2^{m }  2^{p+q+n}  \big)^{1/2} \big( |  x_{\bot} |^{-1} {2^{-m-p-q }}   \big)^{1/2} \\
 &    \lesssim   \| \mathfrak{m}(\cdot, \zeta)\|_{\mathcal{S}^\infty}  |  x_{\bot}|^{-1/2} 2^{\max\{n, (\gamma_1-\gamma_2)M_{t^{\star}}\}} 2^{-k+ \max\{l,p\} +n/2  + 12\epsilon M_{t^{\star}}}\min\{2^{2j}, 2^{2 \alpha^{\star} M_{t^{\star}} -2l }\}  . 
  \end{split}
\ee
After combining the obtained estimates    \eqref{sep23eqn102}  and  \eqref{sep24eqn16},  and the fact that $|  v_{\bot}|/|v|\sim 2^p$ if $p\geq l+\epsilon M_t$ we have 
\be\label{sep24eqn40}
\begin{split}
 &\big| \mathbf{P}\big(H_{k,j;n,l,r; p,q;ess}^{\mu,m,i;non,2}(t, x,\zeta)\big)\big|\\
 &\lesssim  \| \mathfrak{m}(\cdot, \zeta)\|_{\mathcal{S}^\infty} 2^{\max\{n, (\gamma_1-\gamma_2)M_{t^{\star}}\}   + 12\epsilon M_{t^{\star}}} \\
 &\quad\times   \min\{  2^{k +n+ \max\{l,p\} },  | x_{\bot}|^{-1/2} 2^{-k-\max\{l,p\} +2 \alpha^{\star} M_{t^{\star}} +n/2}\}  \\
& \lesssim \| \mathfrak{m}(\cdot, \zeta )\|_{\mathcal{S}^\infty}\big[| x_{\bot}|^{-1/4}2^{(\gamma_1-\gamma_2)M_{t^{\star}}/4} 2^{( \alpha^{\star}-11\epsilon)M_{t^{\star}}}+  |  x_{\bot}|^{-1/8} 2^{(\gamma_1-\gamma_2) M_{t^{\star}}/8-\epsilon M_{t^{\star}}} 
\\
&\quad  \times \min\{  2^{(k+2n)/2 + 2\alpha^{\star} M_{t^{\star}}/3 }, 2^{(k+4n)/2 +  \alpha^{\star} M_{t^{\star}}}\}  \mathbf{1}_{n\geq (\gamma_1-\gamma_2-50\epsilon)M_{t^{\star}}  } \big]. 
  \end{split}
\ee

 \medskip

\noindent \textbf{Step 2.}\quad   If $2^{-k- n+\epsilon M_{t^{\star}}}\leq |  x_{\bot}|\leq  2^{m+l+4\epsilon M_{t^{\star}}} $, $p\geq  l +10\epsilon M_{t^{\star}}$. 

 \medskip

Note that, for this case, we have,    $|  v_{\bot}|/|v|\sim 2^p$ and 
$\forall y\in B(0, 2^{-k- c(m,k,l)+ \epsilon M_{t^{\star}}} ),    | x_{\bot} -   y_{\bot} + (t-s)  \omega_{\bot}| \sim 2^{m+p}. $ From the cylindrical symmetry of the distribution function, the volume of support of $\omega, v$,   and  the estimate of kernels in  \eqref{sep22eqn61},    we have 
\be\label{sep24eqn57}
\begin{split}
 &\big|\mathbf{P}\big(  H_{k,j;n,l,r; p,q;ess}^{\mu,m,i;lin}(t, x,\zeta) \big)\big| \\
 &\lesssim \| \mathfrak{m}(\cdot, \zeta)\|_{\mathcal{S}^\infty} 2^{m+k+n+c(m,k,l)+2\epsilon M_{t^{\star}}  + \max\{n, (\gamma_1-\gamma_2)M_{t^{\star}}\}}\\
 &\quad \times     \min\{   2^{m+3j+2l+2p },  2^{ 2k+ c(m,k,l)+2l-p-j} \}\\
&\lesssim \| \mathfrak{m}(\cdot, \zeta)\|_{\mathcal{S}^\infty} 2^{n+ 5\epsilon M_{t^{\star}}  + \max\{n, (\gamma_1-\gamma_2)M_{t^{\star}}\}}   \\
&\quad \times \min\{2^{m+l+3j+2p}, 2^{-2m-2l-p-j}, 2^{-m-p+k-j}\}\\
&\lesssim \| \mathfrak{m}(\cdot, \zeta)\|_{\mathcal{S}^\infty} 2^{n+ 5\epsilon M_{t^{\star}}  + \max\{n, (\gamma_1-\gamma_2)M_{t^{\star}}\}}\\
&\quad \times  \min\big\{ \big(2^{m+l+3j+2p} \big)^{1/4}\big( 2^{-m-l+k-j}\big)^{3/4}, \big(2^{m+l+3j+2p} \big)^{1/2}\big(2^{-2m-2l-p-j}\big)^{1/2}  \big\}\\
&\lesssim   2^{n+ 5\epsilon M_{t^{\star}}  + \max\{n, (\gamma_1-\gamma_2)M_{t^{\star}}\}}  2^{-(m+l)/2+5\epsilon M_{t^{\star}}} \min\{2^{3k/4+p/2}, 2^{j+p/2}\}
   \| \mathfrak{m}(\cdot, \zeta)\|_{\mathcal{S}^\infty}\\
&\lesssim\| \mathfrak{m}(\cdot, \zeta)\|_{\mathcal{S}^\infty}  |  x_{\bot}|^{-1/2} 2^{(\gamma_1-\gamma_2)M_{t^{\star}}/2-\epsilon M_{t^{\star}}} \\
&\quad \times \big[  2^{(  \alpha^{\star}-10\epsilon)M_{t^{\star}}} +      2^{(k+4n)/2 + 2 \alpha^{\star} M_{t^{\star}}/3  }   \mathbf{1}_{n\geq (\gamma_1-\gamma_2-50\epsilon)M_{t^{\star}}  }   \big].  
\end{split}
\ee

Now, we move on to the estimate of the nonlinear parts in  \eqref{2022feb15eqn2}. From the Cauchy-Schwarz inequality,   the estimate of kernels in  \eqref{sep22eqn61},  the estimate  \eqref{march18eqn31} in Lemma \ref{conservationlawlemma},  the volume of support of $\omega, v$,  and    the cylindrical symmetry of the distribution function, we have
\be\label{sep24eqn58}
\begin{split}
 &\big| \mathbf{P}\big(H_{k,j;n,l,r; p,q;ess}^{\mu,m,i;non,1}(t, x,\zeta)\big)\big|\\
 &\lesssim  \| \mathfrak{m}(\cdot, \zeta)\|_{\mathcal{S}^\infty}   2^{m-j-l} 2^{  \max\{n, (\gamma_1-\gamma_2)M_{t^{\star}}\}+6\epsilon M_{t^{\star}}}  \big( 2^{  -2m } 2^{n-c(m,k,l,n)} 2^{3j+2l} \big)^{1/2} \\
&\quad \times     \big(   \min\{  2^{  2k+c(m,k,l )+2l-p-j}, 2^{m+2l+3j+2p}\} \big)^{1/2} \\
&
\lesssim  \| \mathfrak{m}(\cdot, \zeta)\|_{\mathcal{S}^\infty}  2^{(n-c(m,k,l ))/2} 2^{\max\{n, (\gamma_1-\gamma_2)M_{t^{\star}}\} +j/2+10\epsilon M_{t^{\star}}}  \\ 
& \quad \times \big( \min\{ 2^{-2m-l-p-j}, 2^{m+2l+3j+2p}\} \big)^{1/2}\\
& \lesssim   \| \mathfrak{m}(\cdot, \zeta)\|_{\mathcal{S}^\infty}  2^{(n-l)/2} 2^{\max\{n, (\gamma_1-\gamma_2)M_{t^{\star}}\} +j/2+12\epsilon M_{t^{\star}}} \\ 
& \quad \times \big( \big( 2^{-2m-l-p-j}\big)^{2/3} \big(2^{m+2l+3j+2p}\big)^{1/3} \big)^{1/2} \\
& \lesssim \| \mathfrak{m}(\cdot, \zeta)\|_{\mathcal{S}^\infty} 2^{20\epsilon M_{t^{\star}}}  2^{-(m+l)/2 +n/2+ \max\{n, (\gamma_1-\gamma_2)M_{t^{\star}}\}  + 2j/3 }   \\
&\lesssim    \| \mathfrak{m}(\cdot, \zeta)\|_{\mathcal{S}^\infty} |  x_{\bot}|^{-1/2} 2^{(\gamma_1-\gamma_2)M_{t^{\star}}/2}   2^{(  \alpha^{\star}-10\epsilon)M_{t^{\star}}} .
\end{split}
\ee

Lastly, we   estimate  $\mathbf{P}\big(H_{k,j;n,l,r; p,q;ess}^{\mu,m,i;non,2}(t, x,\zeta)\big)$. From the Cauchy-Schwarz inequality,   the estimate of kernels in  \eqref{sep22eqn61}  and  \eqref{sep22eqn93},  the volume of support of $\omega, v$, and    the cylindrical symmetry of the distribution function and the electromagnetic field, we have
\be\label{sep24eqn59}
\begin{split}
& \big| \mathbf{P}\big(H_{k,j;n,l,r; p,q;ess}^{\mu,m,i;non,2}(t, x,\zeta)\big)\big| \\
& \lesssim   \| \mathfrak{m}(\cdot, \zeta)\|_{\mathcal{S}^\infty}  2^{m-j } 2^{ \max\{l,p\}} 2^{  \max\{n, (\gamma_1-\gamma_2)M_{t^{\star}}\}+6\epsilon M_{t^{\star}}}  \\
&\quad \times  \big(  \min\{   2^{ 2k+c(m,k,l)-p+2l-j}, 2^{m+2l+3j+2p}\}  \big)^{1/2}  \\
& \quad \times  \big( 2^{ 2k+c(m,k,l)-p + 2\min\{n,p\} +3j+2l} \big)^{1/2}  \\
&  \lesssim  \| \mathfrak{m}(\cdot, \zeta)\|_{\mathcal{S}^\infty} 2^{-l/2+j/2+3p/2 } 2^{  \max\{n, (\gamma_1-\gamma_2)M_{t^{\star}}\}+10\epsilon M_{t^{\star}}} \\
&\quad \times  \big( \min\{2^{-2m-l-p-j}, 2^{m+2l+ 3j+2p} \} \big)^{1/2}\\
 &  \lesssim \| \mathfrak{m}(\cdot, \zeta)\|_{\mathcal{S}^\infty} 2^{-(m+l)/2 + 2\alpha^{\star}M_{t^{\star}}/3+3\epsilon M_{t^{\star}}} \\
 &  \lesssim \| \mathfrak{m}(\cdot, \zeta)\|_{\mathcal{S}^\infty}   |  x_{\bot}|^{-1/2} 2^{(\gamma_1-\gamma_2)M_{t^{\star}}/2}   2^{(  \alpha^{\star}-10\epsilon)M_{t^{\star}}}  .
\end{split}
\ee

 \medskip

\noindent \textbf{Step 3.}\quad   If $2^{-k- n+\epsilon M_{t^{\star}}}\leq |  x_{\bot}|\leq  2^{m+l+4\epsilon M_{t^{\star}}}  $, $p\leq  l +10\epsilon M_{t^{\star}}$. 

 \medskip

 From     the estimate of kernels in \eqref{sep22eqn61},      the volume of support of $\omega, v$,  and the estimate  \eqref{nov24eqn41}  if $|  v_{\bot}|\geq 2^{(\alpha_t + \epsilon)M_t}$,   we have 
\be\label{2021dec14eqn23}
\begin{split}
&\big|\mathbf{P}\big(  H_{k,j;n,l,r; p,q;ess}^{\mu,m,i;lin}(t, x,\zeta) \big)\big| \\
& \lesssim \| \mathfrak{m}(\cdot, \zeta)\|_{\mathcal{S}^\infty}  2^{m+k+n+c(m,k,l) + \max\{n, (\gamma_1-\gamma_2)M_{t^{\star}}\}+ 6\epsilon M_{t^{\star}}  }
\\
&\quad  \times \min\big\{     2^{ m+ 3k+ 2c(m,k,l)+2l-j}, 2^{m+ 2l  }\min\{ 2^{ j+2  {\alpha}_t M_{t^{  }}}, 2^{3j+2p}\} \big\}   \\
&\lesssim  \| \mathfrak{m}(\cdot, \zeta)\|_{\mathcal{S}^\infty}  2^{n+ 10\epsilon M_{t^{\star}}  + \max\{n, (\gamma_1-\gamma_2)M_{t^{\star}}\}}  \\
&\quad \times   \min\{ 2^{m+  l  }\min\{ 2^{ j+2  {\alpha}_t M_{t^{  }}}, 2^{3j+2p}\}, 2^{-2m-3l-j} \} \\
&\lesssim \| \mathfrak{m}(\cdot, \zeta)\|_{\mathcal{S}^\infty}  2^{n+ 10\epsilon M_{t^{\star}}  + \max\{n, (\gamma_1-\gamma_2)M_{t^{\star}}\}}  \min\{ 2^{-m-l+j/3}, 2^{-m/2-l/2+j/2+ \alpha^{\star} M_{t^{\star}}/2}\} \\ 
&\lesssim  \| \mathfrak{m}(\cdot, \zeta)\|_{\mathcal{S}^\infty} 2^{n+ 12\epsilon M_{t^{\star}}  + \max\{n, (\gamma_1-\gamma_2)M_{t^{\star}}\}} |  x_{\bot} |^{-1/2} \min\{ 2^{(k+n)/2+j/3}, 2^{ j/2+ \alpha^{\star} M_{t^{\star}}/2}\} \\
 & \lesssim   \| \mathfrak{m}(\cdot, \zeta)\|_{\mathcal{S}^\infty}   |  x_{\bot}|^{-1/2} 2^{(\gamma_1-\gamma_2-\epsilon)M_{t^{\star}}/2}   \big[2^{(\alpha^{\star}  -10\epsilon)M_{t^{\star}}}  +   2^{(k+4n)/2 + 2 \alpha^{\star}  M_{t^{\star}}/3}    \mathbf{1}_{n\geq (\gamma_1-\gamma_2-50\epsilon)M_{t^{\star}}  }    \big]. \\
 \end{split}
\ee

From the Cauchy-Schwarz inequality,   the estimate of kernels in   \eqref{sep22eqn61},  the estimate  \eqref{march18eqn31}  in Lemma \ref{conservationlawlemma},  the volume of support of $\omega, v$,  and    the cylindrical symmetry of the distribution function, we have
\be\label{2021dec14eqn25}
\begin{split}
  &\big| \mathbf{P}\big(H_{k,j;n,l,r; p,q;ess}^{\mu,m,i;non,1}(t, x,\zeta)\big)\big| \\
  &\lesssim   \| \mathfrak{m}(\cdot, \zeta)\|_{\mathcal{S}^\infty} 2^{m-j-l}2^{ \max\{n, (\gamma_1-\gamma_2)M_{t^{\star}}\}+ 6\epsilon M_{t^{\star}} }  \big( 2^{  -2m } 2^{n-c(m,k,l)} 2^{3j+2l} \big)^{1/2} \\
&\quad \times   \big( \min \{ 2^{ m+ 3k+ 2c(m,k,l)+2l-j}  , 2^{m+2p+3j+2l}\}   \big)^{1/2}   \\
&\lesssim \| \mathfrak{m}(\cdot, \zeta)\|_{\mathcal{S}^\infty}  2^{(n- l)/2} 2^{\max\{n, (\gamma_1-\gamma_2)M_{t^{\star}}\}+j/2 +10\epsilon M_{t^{\star}}}   \\
&\quad \times \big( \min\{ 2^{-2m-2l-j}, 2^{m+2l+3j+2p}\} \big)^{1/2} \\ 
&\lesssim \| \mathfrak{m}(\cdot, \zeta)\|_{\mathcal{S}^\infty}  2^{(n- l)/2} 2^{\max\{n, (\gamma_1-\gamma_2)M_{t^{\star}}\}+j/2 +10\epsilon M_{t^{\star}}} \\
&\quad \times  \big( \big(2^{-2m-2l-j}\big)^{2/3} \big(2^{m+2l+3j+2p}\big)^{1/3}\big)^{1/2}  \\ 
& \lesssim \| \mathfrak{m}(\cdot, \zeta)\|_{\mathcal{S}^\infty} 2^{n/2+\max\{n, (\gamma_1-\gamma_2)M_{t^{\star}}\} +30\epsilon M_{t^{\star}}}|  x_{\bot}|^{-1/2}  2^{ 2j/3} \\
& \lesssim \| \mathfrak{m}(\cdot, \zeta)\|_{\mathcal{S}^\infty} |  x_{\bot} |^{-1/2}2^{(\gamma_1-\gamma_2)M_{t^{\star}}/2 } 2^{ (\alpha^{\star}-10 \epsilon) M_{t^{\star}}  }  . 
 \end{split}
\ee 

Lastly, we   estimate  $\mathbf{P}\big(H_{k,j;n,l,r; p,q;ess}^{\mu,m,i;non,2}(t, x,\zeta)\big) $. From the Cauchy-Schwarz inequality,   the estimate of kernels in  \eqref{sep22eqn61}  and \eqref{sep22eqn93}, and   the volume of support of $\omega, v$,   we have
\be\label{sep27eqn61}
\begin{split}
 &\big| \mathbf{P}\big(H_{k,j;n,l,r; p,q;ess}^{\mu,m,i;non,2}(t, x,\zeta)\big)\big| \\
 & \lesssim  \| \mathfrak{m}(\cdot, \zeta)\|_{\mathcal{S}^\infty}  2^{m-j + \max\{n, (\gamma_1-\gamma_2)M_{t^{\star}}\}+\max\{l,p\}+6\epsilon M_{t^{\star}}}     \\
&\quad \times \big( \min \{ 2^{ m+ 3k+ 2c(m,k,l)+2l-j}  , 2^{m+2p+3j+2l}\}   \big)^{1/2}  \big(  2^{ m+ 3k+ 2c(m,k,l)+2p + 3j+2l} \big)^{1/2}  \\ 
&\lesssim    2^{l+\max\{n, (\gamma_1-\gamma_2)M_{t^{\star}}\}+j/2+10\epsilon M_{t^{\star}}} \| \mathfrak{m}(\cdot, \zeta)\|_{\mathcal{S}^\infty} \big( \min\{ 2^{-2m-2l-j}, 2^{m+2l+3j+2p}\} \big)^{1/2} \\ 
 &\lesssim 2^{ \max\{n, (\gamma_1-\gamma_2)M_{t^{\star}}\} +30\epsilon M_{t^{\star}}}|  x_{\bot}|^{-1/2}  2^{ 2j/3} \| \mathfrak{m}(\cdot, \zeta)\|_{\mathcal{S}^\infty}\\
 &\lesssim |  x_{\bot} |^{-1/2}2^{(\gamma_1-\gamma_2)M_{t^{\star}}/2 } 2^{ (\alpha^{\star}-10 \epsilon) M_{t^{\star}}  }  \| \mathfrak{m}(\cdot, \zeta)\|_{\mathcal{S}^\infty}.
 \end{split}
\ee

 \medskip

\noindent \textbf{Step 4.}\quad    If $|  x_{\bot}|\leq 2^{-k- n+2\epsilon M_{t^{\star}}}$. 
 
  \medskip

  From  the estimate of kernels in  \eqref{sep22eqn61},   the volume of support of $v, \omega$, and the estimate  \eqref{nov24eqn41}  if $|  v_{\bot}|\geq 2^{(\alpha_t + \epsilon)M_t}$,  we have
\be\label{2021dec14eqn45}
\begin{split}
 &\big|\mathbf{P}\big(  H_{k,j;n,l,r; p,q;ess}^{\mu,m,i;lin}(t, x,\zeta) \big)\big|  \\
 &\lesssim  \| \mathfrak{m}(\cdot, \zeta)\|_{\mathcal{S}^\infty} 2^{m+k+n+c(m,k,l) +6\epsilon M_{t^{\star}}   }\\
 &\quad \times \min\big\{2^{m+2l} \min\{ 2^{j+2  {\alpha}_t M_{t^{ }}},2^{3j+2p} \}, 2^{m+3k+2c(m,k,l)+2l-j} \}\\ 
&  \lesssim \| \mathfrak{m}(\cdot, \zeta)\|_{\mathcal{S}^\infty} 2^{n+ 10\epsilon M_{t^{\star}}   }   \min\{ 2^{m+  l  } 2^{ 2j+   {\alpha}_t M_{t^{  }}} , 2^{-m+k-l-j}, 2^{2k+n-j}\}  \\
 & \lesssim \| \mathfrak{m}(\cdot, \zeta)\|_{\mathcal{S}^\infty}  2^{n+ 10\epsilon M_{t^{\star}}   }   \big( 2^{m+  l  } 2^{ 2j+   {\alpha}_t M_{t^{  }}}\big)^{1/3}  \big(2^{-m+k-l-j}\big)^{1/3}  \big( 2^{2k+n-j} \big)^{1/3} \\ 
& \lesssim  \| \mathfrak{m}(\cdot, \zeta)\|_{\mathcal{S}^\infty} 2^{10\epsilon M_{t^{\star}}} |  x_{\bot}|^{-1}2^{ \alpha^{\star} M_{t^{\star}}/3}\\
& \lesssim \| \mathfrak{m}(\cdot, \zeta)\|_{\mathcal{S}^\infty} |  x_{\bot}|^{-1}2^{(\gamma_1-\gamma_2)M_{t^{\star}}}  2^{ ( 5\alpha^{\star}/6-\epsilon) M_{t^{\star}} } .
\end{split} 
\ee

 From the Cauchy-Schwarz inequality,   the estimate of kernels in  \eqref{sep22eqn61},  the estimate \eqref{march18eqn31}  in Lemma \ref{conservationlawlemma},  the volume of support of $\omega, v$,  and    the cylindrical symmetry of the distribution function, we have
\be
\begin{split}
  \big| \mathbf{P}\big(H_{k,j;n,l,r; p,q;ess}^{\mu,m,i;non,1}(t, x,\zeta)\big)\big| 
 & \lesssim    \| \mathfrak{m}(\cdot, \zeta)\|_{\mathcal{S}^\infty} 2^{m-j-l}2^{\max\{n, (\gamma_1-\gamma_2)M_{t^{\star}}\}+6\epsilon M_{t^{\star}}}     \\
& \quad \times\big(  2^{ m+ 3k+ 2c(m,k,l)+2l-j}      \big)^{1/2} \big( 2^{  -2m } 2^{n-c(m,k,l,n)} 2^{3j+2l} \big)^{1/2} \\
& \lesssim \| \mathfrak{m}(\cdot, \zeta)\|_{\mathcal{S}^\infty} 2^{10\epsilon M_{t^{\star}}}   2^{k+(l+n)/2} \\
& \lesssim\| \mathfrak{m}(\cdot, \zeta)\|_{\mathcal{S}^\infty}  |  x_{\bot}|^{-1}2^{(\gamma_1-\gamma_2)M_{t^{\star}}}  2^{  \alpha^{\star}  M_{t^{\star}}/2 } . 
\end{split}
\ee

Lastly, we   estimate  $\mathbf{P}\big(H_{k,j;n,l,r; p,q;ess}^{\mu,m,i;non,2}(t, x,\zeta)\big) $. From the Cauchy-Schwarz inequality,   the estimate of kernels in   \eqref{sep22eqn61}  and  \eqref{sep22eqn93}, and   the volume of support of $\omega, v$,  the following estimate holds if $p\leq l +5\epsilon M_{t^{\star}},$ 
\be
\begin{split}
 \big| \mathbf{P}\big(H_{k,j;n,l,r; p,q;ess}^{\mu,m,i;non,2}(t, x,\zeta)\big)\big| 
& \lesssim  \| \mathfrak{m}(\cdot, \zeta)\|_{\mathcal{S}^\infty}  2^{m-j }2^{\max\{n, (\gamma_1-\gamma_2)M_{t^{\star}}\}+ \max\{l,p\}+5\epsilon M_{t^{\star}}}    \\
&\quad \times  \big( 2^{ m+ 3k+ 2c(m,k,l)+2l-j}      \big)^{1/2}   \big(  2^{  m+3k+ 2 c(m,k,l)+ 2p+   3j+2l} \big)^{1/2} \\
& \lesssim     2^{15\epsilon M_{t^{\star}}}   2^{k+l}\| \mathfrak{m}(\cdot, \zeta)\|_{\mathcal{S}^\infty}\\
&  \lesssim  |  x_{\bot}|^{-1}2^{(\gamma_1-\gamma_2)M_{t^{\star}}}  2^{  \alpha^{\star}  M_{t^{\star}}/2 } \| \mathfrak{m}(\cdot, \zeta)\|_{\mathcal{S}^\infty}.
\end{split} 
\ee
Note that, if  $p\geq l+5\epsilon M_{t^{\star}}$, then we have
\[
\forall y \in B(0, 2^{-k - c(m,k,l,n)+\epsilon M_{t^{\star}}}), \quad |  x_{\bot} -  y_{\bot} + (t-s)  \omega_{\bot}| \sim 2^{m+p}. 
\]
From the Cauchy-Schwarz inequality,   the estimate of kernels in  \eqref{sep22eqn61}  and  \eqref{sep22eqn93},  the volume of support of $\omega, v$, the cylindrical symmetry of the distribution function, we have
\be\label{2021dec14eqn46}
\begin{split}
\big| \mathbf{P}\big(H_{k,j;n,l,r; p,q;ess}^{\mu,m,i;non,2}(t, x,\zeta)\big)\big|& \lesssim    \| \mathfrak{m}(\cdot, \zeta)\|_{\mathcal{S}^\infty}  2^{m-j }2^{ 5\epsilon M_{t^{\star}}} \\
&\quad \times  \big(  2^{  2k+ c(m,k,l)+2l-j-p}     \big)^{1/2} \big(  2^{   2k+  c(m,k,l)+  p+   3j+2l} \big)^{1/2}\\
   &\lesssim 2^{10\epsilon M_{t^{\star}}}  2^{k+l } \| \mathfrak{m}(\cdot, \zeta)\|_{\mathcal{S}^\infty}\\
   & \lesssim  |  x_{\bot}|^{-1}2^{(\gamma_1-\gamma_2)M_{t^{\star}}}  2^{  \alpha^{\star}  M_{t^{\star}}/2 } \| \mathfrak{m}(\cdot, \zeta)\|_{\mathcal{S}^\infty}. 
   \end{split}
\ee
 Recall  \eqref{sep24eqn46}. 
To sum up, our desired estimate  \eqref{2022feb16eqn59}  holds after combining the obtained estimates  \eqref{sep23eqn100}--\eqref{2021dec14eqn46}. 
\end{proof}

Now, we proceed to  consider the case   $m+k$ is relatively big, i.e., the case $m+k> -2l+\epsilon M_{t^{\star}}/5 $.  In this case, we mainly use the second formulation in \eqref{sep18eqn50}, i.e., the G-S type decomposition, in Lemma \ref{locdeclemm}. In the following Lemma, we first estimate the $T$-part. 

\begin{lemma}\label{pointestPartIT}
For any $i\in \{2,3,4\}, ( l,r)\in \mathcal{B}_i $, see \eqref{sep18eqn50}.      Under the assumption of  the part (ii) in Theorem \ref{mainresultsfirstpart},     the following estimate holds if   $m+k\geq -2l+ \epsilon M_{t^{\star}}/5, $
  \be\label{2022feb18eqn88}
  \begin{split}
  &  \big| \mathbf{P}\big(   \widetilde{T}_{k,j;n,l,r}^{T; \mu ,m, i}(\mathfrak{m}, E)(t,x,  \zeta ) + \hat{\zeta}\times   \widetilde{T}_{k,j;n,l,r}^{T; \mu,m,i }(\mathfrak{m},B)(t,x, \zeta )\big)  \big| \\
  & \lesssim  \min_{b\in\{1,2\}} \| \mathfrak{m}(\cdot, \zeta)\|_{\mathcal{S}^\infty} \big[ \sum_{a \in \{0, 1/2\} } |  x_{\bot}|^{-a} 2^{a(\gamma_1-\gamma_2)M_{t^\star} }2^{(\alpha^{\star}-10\epsilon)M_{t^\star}}+  \mathbf{1}_{|  x_{\bot}|\leq 2^{-k-n+ b\epsilon M_t}}  \\
  &\quad \times | x_{\bot}|^{-1} 2^{(\gamma_1-\gamma_2)M_{t^\star} } 2^{5\alpha^{\star}M_{t^{\star}} /6}  +   \mathbf{1}_{|  x_{\bot}|\geq 2^{-k-n+ b\epsilon M_t}}  \mathbf{1}_{n\geq  -\alpha^{\star}M_{t^{\star}} /2 } |  x_{\bot}|^{-a} 2^{a(\gamma_1-\gamma_2)M_{t^\star} } \\
  &\quad  \times     \min\{    2^{(k+2n)/2 +2\alpha^{\star} M_{t^{\star} }/3 -(\gamma_1-\gamma_2)M_{t^{\star}}/6 },    2^{(k+4n)/2   + \alpha^{\star} M_{t^{\star} }-(\gamma_1-\gamma_2)M_{t^{\star}}/3 } \} \big]. \\
  \end{split}
  \ee
 
\end{lemma}
\begin{proof}
We still use the decomposition obtained in  \eqref{sep7eqn51}. From the obtained estimate  \eqref{2022feb9eqn22}, it suffices to estimate the essential part $ G_{k,j;n,l,r}^{ess; m, i;p,q}(t,x,\zeta)$ for the case $p,q\in(-10M_t, \epsilon M_t]\cap \Z.$  

Based on the size of $i$, we proceed in two main steps as follows. 

\medskip
\noindent \textbf{Step 1. }\quad  The case $i=2$.

\medskip

Due to fact that $m+k\geq-2l+\epsilon M_{t^{\star}}/5$, recall  \eqref{sep17eqn63},  \eqref{sep4eqn6}, and  \eqref{sep9eqn21}, we know that $|\tilde{v}-\tilde{\zeta}|\sim 2^r \sim 2^l, r\geq  n + 3\epsilon M_t/4.$ 
  Moreover, if $l\notin [(\gamma_1-\gamma_2-\epsilon )M_{t^\star}, (\gamma_1-\gamma_2+\epsilon)M_{t^\star}]$, then we have $|  v_{\bot} |/|v|\sim 2^{\max\{r, (\gamma_1-\gamma_2  )M_{t^\star} \}+j}. $ From the estimate  \eqref{nov24eqn41}, we can rule out the case $\max\{l, (\gamma_1-\gamma_2  )M_{t^\star} \}+j\geq  (\alpha_t+\epsilon) M_t$ if $l\notin [(\gamma_1-\gamma_2-\epsilon )M_{t^\star}, (\gamma_1-\gamma_2+\epsilon)M_{t^\star}]$.

 From the detailed formulas of coefficients $\omega^{m;U}_{j,l}(t-s,v,\omega), c^{q;m,U}_{j,l}(t-s,v,\omega),  c^{err;m,U}_{j,l}(t-s,v,\omega),  U\in \{E, B\} $, in \eqref{july9eqn11}  and  \eqref{sep20eqn44}, the following improved estimate holds for the projection of $\omega^{m;E}_{j,l}(t-s,v,\omega) +\hat{\zeta}\times \omega^{m;B}_{j,l}(t-s,v,\omega)$,
\be\label{sep25eqn15}
\begin{split}
& \big|\mathbf{P}\big( \omega^{m;E}_{j,l,r}(t-s,v,\omega) +\hat{\zeta}\times \omega^{m;B}_{j,l}(t-s,v,\omega)\big)\big| \\
&  +\big|\mathbf{P}\big(c^{a;m,E}_{j,l,r}(t-s,v,\omega)  
+\hat{\zeta}\times c^{a;m,B}_{j,l,r}(t-s,v,\omega)\big)\big|\\
& +\big|\mathbf{P}\big(c^{err;m,E}_{j,l,r}(t-s,v,\omega)  +\hat{\zeta}\times c^{err;m,B}_{j,l,r}(t-s,v,\omega)\big)\big|\\
&\lesssim 2^{  \max\{ l, (\gamma_1-\gamma_2)M_{t^{\star}}\} }  . 
\end{split}
\ee

Based on the possible size of $|  x_{\bot}|$, we proceed in three sub-steps as follows.

\medskip

 \textbf{Step 1A. } \quad If $|  x_{\bot} | \leq 2^{-k-n+2\epsilon M_t}$. 
\medskip

Recall  \eqref{sep9eqn21}. Note that, the volume of support of $\omega$ is bounded by $2^{2\epsilon M_t}(2^{2n}+ 2^{-m-k})$. From  the estimate of coefficients in  \eqref{sep25eqn15},    the estimate of kernels in \eqref{sep6eqn31},      the estimate  \eqref{march18eqn31}  in Lemma \ref{conservationlawlemma}, the volume of support of $\omega, v$, and the estimate  \eqref{nov24eqn41}  if $| v_{\bot}|\geq 2^{(\alpha_t +\epsilon) M_t}$,   we have
 
\be\label{2022feb17eqn21}
\begin{split}
&\big| \mathbf{P}\big( G_{k,j; n,l,r }^{ess; m, i;p,q}(t,x,\zeta)\big)  \big| \\ &\lesssim  \| \mathfrak{m}(\cdot, \zeta)\|_{\mathcal{S}^\infty} 2^{ 5\epsilon M_{t^{\star}}}  
\min\big\{ 2^{-2m-j-2l}, (2^{2n}+2^{-m-k}) \\
 &  \quad \times  \min\{2^{m+j+2\alpha_t M_t},   2^{m  + 3j+2l}\},    (2^{2n}+2^{-m-k}) 2^{m+3k+2n -j}\big\}   \\
   & \lesssim  2^{  10\epsilon M_{t^{\star}}}  2^{k+ n +  \alpha^{\star} M_{t^{\star}}/3  } \| \mathfrak{m}(\cdot, \zeta)\|_{\mathcal{S}^\infty} \\
   &\lesssim\| \mathfrak{m}(\cdot, \zeta)\|_{\mathcal{S}^\infty}  |  x_{\bot}|^{-1} 2^{ (\gamma_1-\gamma_2)M_{t^{\star}}}  2^{5  \alpha^{\star} M_{t^{\star}}/6}.
\end{split}
\ee

\medskip

 \textbf{Step 1B. } \quad  If $    2^{-k-n+ \epsilon M_t}\leq |  x_{\bot} | \leq 2^{\epsilon M_t} \max\{2^m , 2^{-k-n}\}$. 

 \medskip

 From the estimate of coefficients in  \eqref{sep25eqn15},  the estimate of kernels in  \eqref{sep6eqn31}, the estimate   \eqref{march18eqn31}  in Lemma \ref{conservationlawlemma},   the volume of support of $\omega$, and the estimate  \eqref{nov24eqn41}  if $|  v_{\bot}|\geq 2^{(\alpha_t +\epsilon) M_t}$,     the following estimate holds 
\[
\begin{split}
&\big| \mathbf{P}\big( G_{k,j; n,l,r  }^{ess; m, i;p,q}(t,x,\zeta)\big)  \big|\\
&  \lesssim   \| \mathfrak{m}(\cdot, \zeta) \|_{\mathcal{S}^\infty} 2^{ \max\{l, (\gamma_1-\gamma_2)M_{t^{\star}}\}+6\epsilon M_{t^{\star}}}\min\big\{   2^{-2m-j-2l},    (2^{2n} +2^{-m-k})     \\
& \quad \times  2^{ m+ 3k+2n -j} ,  (2^{2n} +2^{-m-k}) 2^{m} \min\{2^{ j+2\alpha_t M_t},  2^{3j+2l}\}  \big\} \\
&   \lesssim  \| \mathfrak{m}(\cdot, \zeta)\|_{\mathcal{S}^\infty} 2^{ \max\{l, (\gamma_1-\gamma_2)M_{t^{\star}}\}+6\epsilon M_{t^{\star}}}    \min\big\{   2^{-2m-j-2l},    2^{ 2k+2n -j},\\
&\quad  2^{ 2n/3}\min\{2^{4\alpha_t M_t/3 +j/3}, 2^{4l/3+5j/3} \} \big\}. \\
  \end{split}
\] 
From the above estimate, we conclude that 
\be\label{2022feb17eqn11}
\begin{split}
&\big| \mathbf{P}\big( G_{k,j; n,l,r  }^{ess; m, i;p,q}(t,x,\zeta)\big)  \big|\\
 &\lesssim  \| \mathfrak{m}(\cdot, \zeta)\|_{\mathcal{S}^\infty} 2^{ \max\{l, (\gamma_1-\gamma_2)M_{t^{\star}}\}+6\epsilon M_{t^{\star}}}    \big(2^{-2m-j-2l}\big)^{1/4} \\
 &\quad \times \min\big\{  \big(2^{ 2n/3}(2^{4\alpha_t M_t/3 +j/3})^{1/4} (2^{4l/3+5j/3})^{3/4}\big)^{1/2}\big( 2^{ 2k+2n -j}\big)^{1/4} , \\
 &\qquad   \big(2^{ 2n/3} \min\{2^{4\alpha_t M_t/3 +j/3}, 2^{4l/3+5j/3} \} \big)^{3/4}   \big\}\\
&\lesssim \| \mathfrak{m}(\cdot, \zeta)\|_{\mathcal{S}^\infty}   2^{ \max\{l, (\gamma_1-\gamma_2)M_{t^{\star}}\}+8\epsilon M_{t^{\star}}} |  x_{\bot}|^{-1/2}\\
&\quad\times  \min\big\{  2^{n/2   }\min\{2^{j+l/2},2^{\alpha_t M_t-l/2}\} , 2^{(k+n)/2+n/3+ (j+\alpha_t M_t)/6}\big\} \\
&\lesssim  \| \mathfrak{m}(\cdot, \zeta)\|_{\mathcal{S}^\infty} | x_{\bot}|^{-1/2} 2^{ (\gamma_1-\gamma_2-\epsilon )M_{t^{\star}}/2} \\
&\quad \times   \big[ 2^{  (\alpha^{\star}-10\epsilon) M_{t^{\star}} }  +  2^{(k+2n)/2+ 2\alpha^{\star}M_{t^{\star}}/3  }  \mathbf{1}_{n\geq  (\gamma_1-\gamma_2-30\epsilon )M_{t^{\star}}} \big].
\end{split}
\ee
 
\medskip

 \textbf{Step 1C. } \quad  If $    |  x_{\bot} | \geq   2^{\epsilon M_t} \max\{2^m , 2^{-k-n}\}$. 

 \medskip

Note that, for this case, we have
\[
\forall y \in B(0, 2^{-k-n+\epsilon M_{t^{  }}/2}), \quad |  x_{\bot} -   y_{\bot} +(t-s)  \omega_{\bot}|\sim | x_{\bot}|, \quad |  x_{\bot}-   y_{\bot}| \sim |   x_{\bot}|. 
\]
From the above estimate, the estimate of kernels in  \eqref{sep6eqn31}, the cylindrical symmetry of the distribution function,   the volume of support of $\omega$, and the estimate  \eqref{nov24eqn41}   if $| v_{\bot}|\geq 2^{(\alpha_t +\epsilon) M_t}$,   the following estimate holds if $m+k\leq -2n,$
\be\label{2022feb17eqn12}
\begin{split}
&\big| \mathbf{P}\big( G_{k,j; n,l,r  }^{ess; m, i;p,q}(t,x,\zeta)\big)  \big| \\
& \lesssim   \| \mathfrak{m}(\cdot, \zeta)\|_{\mathcal{S}^\infty} 2^{6\epsilon M_{t^{\star}}}    \min\big\{   |  x_{\bot}|^{-1} 2^{k+n -j}   ,  2^{-k}\min\{2^{3j+2l}, 2^{j+2  {\alpha}_t M_{t^{  }}}\}   \big\} \\
&\lesssim   \| \mathfrak{m}(\cdot, \zeta)\|_{\mathcal{S}^\infty}   2^{  6\epsilon M_{t^{\star}} }  \min\{   |  x_{\bot} |^{-1}2^{k+n-j},   |  x_{\bot}|^{-1/2} 2^{ n/2}\min\{2^{j+l}, 2^{  {\alpha}_t M_{t^{  }}  }\}  \}\\
&\lesssim  \| \mathfrak{m}(\cdot, \zeta)\|_{\mathcal{S}^\infty}   2^{   6\epsilon M_{t^{\star}} } \min\big\{   |  x_{\bot}|^{-1/2} 2^{ n/2} 2^{  {\alpha}_t M_{t^{  }}  },  \big( |  x_{\bot}|^{-1/2}2^{3(k+n)/2-j}\big)^{1/3}\\
&\quad \times  \big(|  x_{\bot}|^{-1/2} 2^{ n/2}   (2^{j+l})^{1/2}(2^{  {\alpha}_t M_{t^{  }}  })^{1/2}\big)^{2/3} \big\}  \\
&\lesssim   \| \mathfrak{m}(\cdot, \zeta)\|_{\mathcal{S}^\infty} 2^{  8\epsilon M_{t^{\star}}} |  x_{\bot}|^{-1/2} \min\big\{  2^{n/2 +   \alpha^{\star} M_{t^{\star}}  }, 2^{(k+n)/2+n/3+ \alpha_t M_t /3}\big\}\\
&\lesssim  \| \mathfrak{m}(\cdot, \zeta)\|_{\mathcal{S}^\infty} |  x_{\bot}|^{-1/2} 2^{ (\gamma_1-\gamma_2-\epsilon )M_{t^{\star}}/2} \\
&\quad \times   \big[ 2^{  (\alpha^{\star}-10\epsilon) M_{t^{\star}} }  +  2^{(k+2n)/2+ 2\alpha^{\star}M_{t^{\star}}/3  }  \mathbf{1}_{n\geq  (\gamma_1-\gamma_2-30\epsilon )M_{t^{\star}}} \big].
\end{split}
\ee

 If $m+k\geq -2n,$  in addition to   the  strategy used above, we also use   the Jacobian of changing coordinates $(\theta, \phi)\longrightarrow (z(y), w(y))$in \eqref{march18eqn66}  and the fact that  the volume of support of $\omega$ is bounded by $2^{2\epsilon M_t} \min\{2^{2p+q}, 2^{p+q+n}, 2^{2n}\}$ for the case $m+k\geq -2n$. As a result, we have 
\be\label{2022feb17eqn13}
\begin{split}
&\big| \mathbf{P}\big( G_{k,j; n,l,r  }^{ess; m, i;p,q}(t,x,\zeta)\big)  \big| \\
&  \lesssim  \| \mathfrak{m}(\cdot, \zeta) \|_{\mathcal{S}^\infty}  2^{ 6\epsilon M_{t^{\star}}} 
 \min\big\{   |  x_{\bot}|^{-1} 2^{m+2k+n+p+q+n -j}, | x_{\bot}|^{-1} 2^{-m-p-q -j} , \\
& \qquad   2^{m+n+p+q}\min\{2^{3j+2l}, 2^{j+2  {\alpha}_t M_{t^{  }}}\}  \big\}   \\
& \lesssim    \| \mathfrak{m}(\cdot, \zeta) \|_{\mathcal{S}^\infty} 2^{  8\epsilon M_{t^{\star}}} \min\{  |  x_{\bot}|^{-1/2} 2^{n/2 } \min\{2^{j+l}, 2^{  {\alpha}_t M_{t^{  }}  }\}, |  x_{\bot} |^{-1} 2^{k+n-j} \big\}\\
&\lesssim  \| \mathfrak{m}(\cdot, \zeta)\|_{\mathcal{S}^\infty}   2^{   8\epsilon M_{t^{\star}} } \min\big\{   | x_{\bot}|^{-1/2} 2^{ n/2} 2^{  {\alpha}_t M_{t^{  }}  },  \big( | x_{\bot}|^{-1/2}2^{3(k+n)/2-j}\big)^{1/3}\\
&\quad \times  \big(|  x_{\bot}|^{-1/2} 2^{ n/2}   (2^{j+l})^{1/2}(2^{  {\alpha}_t M_{t^{  }}  })^{1/2}\big)^{2/3} \big\} \\
&\lesssim  \| \mathfrak{m}(\cdot, \zeta)\|_{\mathcal{S}^\infty} |  x_{\bot}|^{-1/2} 2^{ (\gamma_1-\gamma_2-\epsilon )M_{t^{\star}}/2} \\
&\quad \times   \big[ 2^{  (\alpha^{\star}-10\epsilon) M_{t^{\star}} }  +  2^{(k+2n)/2+ 2\alpha^{\star}M_{t^{\star}}/3  }  \mathbf{1}_{n\geq  (\gamma_1-\gamma_2-30\epsilon )M_{t^{\star}}} \big].
\end{split}
\ee

\medskip
\noindent \textbf{Step 2. }\quad The case $i=3,4$.

\medskip

 Due to the  cutoff functions $\varphi^i_{j,n}(v, \zeta), i\in\{3,4\}$, see  \eqref{sep4eqn6}, and the fact that $m+k\geq-2l+\epsilon M_{t^{\star}}/5 $,  the support of the $v, \omega, \xi$,    implies that $|\tilde{v}\times \omega|\lesssim  2^{n+\epsilon M_t} +2^{l-\epsilon M_t/20}, |\tilde{v}\times\omega |\sim 2^l$. Since   $ |\omega \times\tilde{\xi}|\sim 2^l$ if $l> -j$, we have $l\leq n +2\epsilon M_t$ if $l> -j$. As a result, we conclude that $l\in [-j, \max\{-j, n+2\epsilon M_t\}]\cap \Z. $

 Moreover,  from the detailed formulas of coefficients $\omega^{m;U}_{j,l}(t-s,v,\omega), c^{q;m,U}_{j,l}(t-s,v,\omega),  c^{err;m,U}_{j,l}(t-s,v,\omega),  U\in \{E, B\} $, in \eqref{july9eqn11}   and  \eqref{sep20eqn44}, the following   estimate holds for the projection of $\omega^{m;E}_{j,l}(t-s,v,\omega) +\hat{\zeta}\times \omega^{m;B}_{j,l}(t-s,v,\omega)$ for the case we are considering,
\be\label{2022feb20eqn31}
\begin{split}
&\big|\mathbf{P}\big( \omega^{m;E}_{j,l,r}(t-s,v,\omega) +\hat{\zeta}\times \omega^{m;B}_{j,l}(t-s,v,\omega)\big)\big| \\
& +\big|\mathbf{P}\big(c^{a;m,E}_{j,l,r}(t-s,v,\omega)  
+\hat{\zeta}\times c^{a;m,B}_{j,l,r}(t-s,v,\omega)\big)\big|\\
&+\big|\mathbf{P}\big(c^{err;m,E}_{j,l,r}(t-s,v,\omega)+\hat{\zeta}\times c^{err;m,B}_{j,l,r}(t-s,v,\omega)\big)\big|\\
& \lesssim 2^{ -\min\{l,n\}+n + \max\{ n, (\gamma_1-\gamma_2)M_{t^{\star}}\} +\epsilon M_t}   . 
\end{split}
\ee

We first rule out the case $ \max\{-j, n+2\epsilon M_t\}=-j$, which implies that $l=-j$. From    the estimate of coefficients in  \eqref{2022feb20eqn31},    the estimate of kernels in  \eqref{sep6eqn31}, we have 
\be\label{2022feb21eqn41}
\begin{split}
 \big| \mathbf{P}\big( G_{k,j; n,l,r  }^{ess; m, i;p,q}(t,x,\zeta)\big)  \big| 
&\lesssim  \| \mathfrak{m}(\cdot, \zeta)\|_{\mathcal{S}^\infty}  2^{ -\min\{l,n\}+n + \max\{ n, (\gamma_1-\gamma_2)M_{t^{\star}}\} +\epsilon M_t}   
  2^{m+2l}2^{3j+2n+4\epsilon M_t}\\
  & \lesssim \| \mathfrak{m}(\cdot, \zeta)\|_{\mathcal{S}^\infty} 2^{30\epsilon M_{t^\star}}. 
  \end{split}
\ee
Therefore, for the rest of the proof, it suffices to consider   the case $ \max\{-j, n+2\epsilon M_t\}>-j$, which implies that $l\leq n + 2\epsilon M_t. $

 Based on the possible size of $p$ and $|  x_{\bot}|$, we proceed in six sub-steps  as follows.

 \medskip

\textbf{Step 2A.} \quad If $p\leq l +5\epsilon M_{t^{  }}$, $  | x_{\bot} |\leq  2^{-k-n+2\epsilon M_{t^{  }}}$.

\medskip
  
From  the estimate of coefficients in \eqref{2022feb20eqn31},    the estimate of kernels in  \eqref{sep6eqn31},      the estimate  \eqref{march18eqn31}  in Lemma \ref{conservationlawlemma}, the volume of support of $\omega, v$,     we have
\be\label{2021dec15eqn23}
\begin{split}
&\big| \mathbf{P}\big( G_{k,j; n,l,r  }^{ess; m, i;p,q}(t,x,\zeta)\big)  \big|\\
&\lesssim 2^{- l +n  +5\epsilon M_{t^{\star}}}  
\min\{   2^{m+2p+3j+2 l}, 2^{m+2p +j+2   {\alpha}_t M_{t^{  }}},  2^{-2m-j-2l},  2^{m+3k+2n+2 l -j}\} \\
 &\lesssim \| \mathfrak{m}(\cdot, \zeta)\|_{\mathcal{S}^\infty}  2^{- l +n  +5\epsilon M_{t^{\star}}}  
  \big( 2^{m+2p+3j+2 l}\big)^{1/6} \big(2^{m+2p +j+2   {\alpha}_t M_{t^{  }}}\big)^{1/6}  \\
  &\quad \times   \big(2^{-2m-j-2l}\big)^{1/3} \big( 2^{m+3k+2n+2 l -j}\big)^{1/3}\\
   & \lesssim \| \mathfrak{m}(\cdot, \zeta)\|_{\mathcal{S}^\infty}   2^{ n  +10\epsilon M_{t^{\star}}}  2^{k+2n/3+  \alpha^{\star} M_{t^{\star}}/3  } \\
   & \lesssim \| \mathfrak{m}(\cdot, \zeta)\|_{\mathcal{S}^\infty} 2^{ 20\epsilon M_{t^{\star}}} |  x_{\bot}|^{-1}  2^{  \alpha^{\star} M_{t^{\star}}   /3 }.
\end{split}
\ee

  \medskip

\textbf{Step 2B.} \quad  If $p\geq l + 5\epsilon M_{t^{  }}, |  x_{\bot}|\leq  2^{ -k-n+2\epsilon M_{t^{  }}}$. 

  \medskip

Note that, for this case, we have $|  v_{\bot} |/|v|\sim 2^p.$ Moreover,  $\forall y \in B(0, 2^{-k-n+\epsilon M_{t^{  }}}), |  x_{\bot} -  y_{\bot} + (t-s)  \omega_{\bot}|\sim 2^{m+p}. $ From  the estimate of coefficients in  \eqref{2022feb20eqn31},    the estimate of kernels in \eqref{sep6eqn31}, the cylindrical symmetry of the distribution function,   and the estimate  \eqref{march18eqn31}  in Lemma \ref{conservationlawlemma},   we have
\be\label{sep25eqn13}
\begin{split}
& \big| \mathbf{P}\big( G_{k,j; n,l,r  }^{ess; m, i;p,q}(t,x,\zeta)\big)  \big| \\
& \lesssim \| \mathfrak{m}(\cdot, \zeta)\|_{\mathcal{S}^\infty} 2^{-l+n  +5\epsilon M_{t^{\star}}}  
\min\big\{   2^{m+3j+2l+2 \min\{n,p\} }, 2^{-2m-j-2l},  2^{m+3k+2n+2l-k-n-m-p-j}\big\} \\ 
& \lesssim \| \mathfrak{m}(\cdot, \zeta)\|_{\mathcal{S}^\infty}  2^{-l+n  +5\epsilon M_{t^{\star}}}  
\big(   2^{m+3j+2l+2 \min\{n,p\} }\big)^{1/3} \big(2^{-2m-j-2l}\big)^{1/6} \big( 2^{2k+n+2l -p-j}\big)^{1/2} \\
  &  \lesssim \| \mathfrak{m}(\cdot, \zeta)\|_{\mathcal{S}^\infty} 2^{ n+    10\epsilon M_{t^{\star}}}    2^{k+n/2+l/3+p/6+j/3} \\
   & \lesssim \| \mathfrak{m}(\cdot, \zeta)\|_{\mathcal{S}^\infty}  | x_{\bot}|^{-1} 2^{ \alpha^{\star}  M_{t^{\star}}/3+20\epsilon M_{t^{\star}}}.
   \end{split}
\ee

  \medskip

\textbf{Step 2C.} \quad If $p\leq l + 5\epsilon M_{t^{  }}$, $2^{-k-n+\epsilon M_{t^{  }}}\leq |  x_{\bot} |< 2^{\epsilon M_{t^{  }}}\max\{  2^{m+p }, 2^{-k-n}\}$.

 \medskip 

From the estimate of coefficients in  \eqref{2022feb20eqn31},  the estimate of kernels in  \eqref{sep6eqn31}, the estimate  \eqref{march18eqn31}  in Lemma \ref{conservationlawlemma},   the volume of support of $\omega$,   we have
\be\label{sep25eqn31}
\begin{split}
&\big| \mathbf{P}\big( G_{k,j;n,l,r }^{ess; m, i;p,q}(t,x,\zeta)\big)  \big| \\
& \lesssim \| \mathfrak{m}(\cdot, \zeta)\|_{\mathcal{S}^\infty} 2^{- l+n + \max\{n, (\gamma_1-\gamma_2)M_{t^{\star}}\}+4\epsilon M_{t^{\star}}} \\
&\quad \times   \min\{2^{m+2p+3j+2l}, 2^{m+2l+j+2  {\alpha}_t M_{t^{  }}}, 2^{-2m-j-2l}, 2^{m+3k+2n+2l-j}\} \\
& \lesssim\| \mathfrak{m}(\cdot, \zeta)\|_{\mathcal{S}^\infty}  2^{- l+n  + \max\{n, (\gamma_1-\gamma_2)M_{t^{\star}}\} +4\epsilon M_{t^{\star}}} (2^{-2m-j-2l})^{1/2} \\
&\quad \times \min\{(2^{m+ p+2j+  {\alpha}_t M_{t^{  }}+2l})^{1/2} ,(2^{m+ 2p+3j+2l})^{1/3}  (2^{m+3k+2n+2l-j})^{1/6}\}   \\ 
&\lesssim \| \mathfrak{m}(\cdot, \zeta)\|_{\mathcal{S}^\infty} 2^{n+10\epsilon M_{t^{\star}} +  \max\{n, (\gamma_1-\gamma_2)M_{t^{\star}}\}} \\
&\quad \times \min\big\{2^{-(m+l)/2+(j+  \alpha^{\star}  M_{t^{\star}})/2}, 2^{-m/2+k/2+j/3+n/3-l/3}\big\}\\
&\lesssim   \| \mathfrak{m}(\cdot, \zeta)\|_{\mathcal{S}^\infty} |  x_{\bot} |^{-1/2}  2^{(\gamma_1-\gamma_2-\epsilon )M_{t^{\star}}/2}\\
&\quad \times  \big[   2^{(  \alpha^{\star} -10\epsilon)M_{t^{\star}}}+  2^{(k+4n)/2 + 2\alpha^{\star}  M_{t^{\star}}  /3}   \mathbf{1}_{n\geq  (\gamma_1-\gamma_2-30\epsilon )M_{t^{\star}}}  \big]. 
\end{split}
\ee

 \medskip

\textbf{Step 2D.} \quad If $p\leq l + 5\epsilon M_{t^{ }}$, $|  x_{\bot} |\geq  2^{\epsilon M_{t^{  }}}\max\{  2^{m+p }, 2^{-k-n}\}$.

 \medskip

Note that, for this case, we have
\[
\forall y \in B(0, 2^{-k-n+\epsilon M_{t^{\star}}/2}), \quad |  x_{\bot} -   y_{\bot} +(t-s)  \omega_{\bot}|\sim | x_{\bot}|, \quad | x_{\bot}-   y_{\bot}| \sim |   x_{\bot}|. 
\]
From the above estimate, the estimate of kernels in  \eqref{sep6eqn31}, the cylindrical symmetry of the distribution function, the estimate  \eqref{march18eqn31}  in Lemma \ref{conservationlawlemma}, the Jacobian of changing coordinates $(\theta, \phi)\longrightarrow (z(y), w(y))$in  \eqref{march18eqn66},    the volume of support of $\omega$,  and the estimate  \eqref{nov24eqn41} if $|  v_{\bot}|\geq 2^{(\alpha_t +\epsilon) M_t}$,   we have
 \be\label{sep25eqn32}
\begin{split}
&\big| \mathbf{P}\big( G_{k,j; n,l,r }^{ess; m, i;p,q}(t,x,\zeta)\big)  \big| \\
&\lesssim   \| \mathfrak{m}(\cdot, \zeta)\|_{\mathcal{S}^\infty} 2^{- l+n + \max\{n, (\gamma_1-\gamma_2)M_{t^{\star}}\}+4\epsilon M_{t^{\star}}}   \min\big\{  |  x_{\bot}|^{-1} {2^{-m-p-q -j}}{},\\
 &\quad  |  x_{\bot}|^{-1} 2^{m+2k+n+2p+q-j} , 2^{m+2p+q }\min\{2^{3j+2l}, 2^{j+2  {\alpha}_t M_{t^{  }}}\} ,  2^{-2m-j-2l}  \big\}\\
 & \lesssim   \| \mathfrak{m}(\cdot, \zeta)\|_{\mathcal{S}^\infty}  2^{- l+n } 2^{  \max\{n, (\gamma_1-\gamma_2)M_{t^{\star}}\}+4\epsilon M_{t^{\star}} }\\ 
&\quad\times   \min\{|  x_{\bot}|^{-1/2} 2^{(j+  {\alpha}_t M_{t^{  }})/2+l/2+p/2} , |  x_{\bot}|^{-1}  2^{k+n/2+p/2-j}, 2^{4p/3+ 7j/6 +\alpha_t M_t/2+l/6 } \} \\
& \lesssim     \| \mathfrak{m}(\cdot, \zeta)\|_{\mathcal{S}^\infty} |  x_{\bot} |^{-1/2}  2^{(\gamma_1-\gamma_2-\epsilon )M_{t^{\star}}/2} \big[   2^{(  \alpha^{\star} -10\epsilon)M_{t^{\star}}}\\ 
& \quad  + \min\{  2^{(k+2n)/2 +  2\alpha^{\star}  M_{t^{\star}}/3   },  2^{(k+4n)/2 +  \alpha^{\star}  M_{t^{\star}} } \}  \mathbf{1}_{n\geq  (\gamma_1-\gamma_2-30\epsilon )M_{t^{\star}}}  \big].
 \end{split}
\ee

  \medskip

\textbf{Step 2E.} \quad   If $p\geq l + 5\epsilon M_{t^{  }}$, $2^{-k-n+\epsilon M_{t^{  }}}\leq |  x_{\bot} |< 2^{\epsilon M_{t^{  }}}\max\{  2^{m+p }, 2^{-k-n}\}$.

 \medskip

 Note that, for this case, we have $| v_{\bot} |/|v|\sim 2^p$, i.e., $p+j\leq (\alpha_t + \epsilon)M_t+2$,  and  $\forall y\in B(0, 2^{-k-n+\epsilon M_{t^{  }}/2})$, we have $|  x_{\bot} -   y_{\bot}|\sim |  x_{\bot}|.$ Moreover,  for any fixed $v\in \R^3$, we have
 $Vol(supp(\phi))\lesssim 2^{l-p}.$ From   the estimate of coefficients in  \eqref{2022feb20eqn31}, the estimate of kernels in  \eqref{sep6eqn31},  the volume of support of $\omega, v$,    and the Jacobian of changing coordinates $(\theta, \phi)\longrightarrow (z,w)$  in \eqref{march18eqn66},   we have 
 \be\label{sep25eqn10}
 \begin{split}
  \big| \mathbf{P}\big( G_{k,j; n,l,r }^{ess; m, i;p,q}(t,x,\zeta)\big)  \big|  
 & \lesssim   \| \mathfrak{m}(\cdot, \zeta)\|_{\mathcal{S}^\infty}  2^{- l+n + \max\{n, (\gamma_1-\gamma_2)M_{t^{\star}}\}+4\epsilon M_{t^{\star}}}   \\ 
 &\quad \times \min\{2^{-2m-j-2l}, 2^{m+3j+2l+2p+q},  |  x_{\bot}|^{-1} 2^{-m-p-q -j}   \}   \\
 & \lesssim \| \mathfrak{m}(\cdot, \zeta)\|_{\mathcal{S}^\infty}2^{n+\max\{n, (\gamma_1-\gamma_2)M_{t^{\star}}\}+10\epsilon M_{t^{\star}}}  \\
 &\quad \times 
\min\{  |  x_{\bot}|^{-1/2} 2^{ j+p/2}, |  x_{\bot}|^{-1/2} 2^{3k/2-j+p/2} \}\\
 &\lesssim \| \mathfrak{m}(\cdot, \zeta)\|_{\mathcal{S}^\infty}2^{n+\max\{n, (\gamma_1-\gamma_2)M_{t^{\star}}\}+10\epsilon M_{t^{\star}}} |  x_{\bot}|^{-1/2} \\
 &\quad \times 
\min\{   2^{ j+p/2},  \big(2^{3k/2-j+p/2}\big)^{1/3}\big(  2^{ j+p/2}\big)^{2/3} \}\\
  &\lesssim \| \mathfrak{m}(\cdot, \zeta)\|_{\mathcal{S}^\infty} |  x_{\bot}|^{-1/2} 2^{(\gamma_1-\gamma_2-\epsilon )M_{t^{\star}}/2}\\
  &\quad \times   \big[ 2^{  (\alpha^{\star} -10\epsilon) M_{t^{\star}}    }    + 2^{(k+4n)/2 + 2 \alpha^{\star}  M_{t^{\star}}/3 }  \mathbf{1}_{n\geq  (\gamma_1-\gamma_2-30\epsilon )M_{t^{\star}}}  \big].\\
  \end{split}
\ee
In the above estimate, we used the fact that $-m-2l\leq k $.
 
\medskip

\textbf{Step 2F.} \quad  If $p\geq l + 5\epsilon M_{t^{\star}}$,  $|  x_{\bot} |\geq  2^{\epsilon M_{t^{  }}}\max\{  2^{m+p }, 2^{-k-n}\}$.

\medskip

 As in the previous case,   we have $|  v_{\bot} |/|v|\sim 2^p$, i.e., $p+j\leq (\alpha_t + \epsilon)M_t+2$ and   for any fixed $v\in \R^3$, we have
 $Vol(supp(\phi))\lesssim 2^{l-p}.$ Moreover,   $\forall y\in B(0, 2^{-k-n+\epsilon M_{t^{\star}}/2})$, we have $|  x_{\bot} -   y_{\bot}|\sim | x_{\bot}|, |  x_{\bot} -   y_{\bot} +(t-s)  \omega_{\bot}|\sim |   x_{\bot} |$.  

From   the estimate of coefficients in  \eqref{2022feb20eqn31}, the estimate of kernels in  \eqref{sep6eqn31},  the volume of support of $\omega, v$,   the Jacobian of changing coordinates $(\theta, \phi)\longrightarrow (z,w)$  in  \eqref{march18eqn66}, the estimate  \eqref{march18eqn31}  in Lemma \ref{conservationlawlemma},   the cylindrical symmetry of the distribution function,  and the Jacobian of changing coordinates $(y_1, y_2, \theta)\longrightarrow x-y+(t-s)\omega, $  we have 
\[
\begin{split}
 &\big| \mathbf{P}\big( G_{k,j;l,n}^{ess; m, i;p,q}(t,x,\zeta)\big)  \big|\\
 & \lesssim  \| \mathfrak{m}(\cdot, \zeta)\|_{\mathcal{S}^\infty}   2^{- l+n + \max\{n, (\gamma_1-\gamma_2)M_{t^{\star}}\}+4\epsilon M_{t^{\star}}}  \min\big\{2^{m+3j+2l+2p+q},    \\
 &\quad    | x_{\bot} |^{-1}  {2^{-m-p-q-j}} ,|  x_{\bot} |^{-1} 2^{  3k+2n-k-n-k-n+l-p-j   }  ,  2^{-2m-j-2l} \big\} \\
 &\lesssim   \| \mathfrak{m}(\cdot, \zeta)\|_{\mathcal{S}^\infty} 2^{ n + \max\{n, (\gamma_1-\gamma_2)M_{t^{\star}}\}+4\epsilon M_{t^{\star}}}  \min\big\{ \big(2^{m+3j+l+2p+q}\big)^{1/2} \big( |  x_{\bot} |^{-1}  {2^{-m-l-p-q-j}} \big)^{1/2},\\
& \quad  \big(2^{m+3j+l+2p+q}\big)^{3/4} \big(2^{-2m-j-3l}\big)^{1/4} , |  x_{\bot} |^{-1} 2^{  k -p-j   } \big\}\\
 \end{split} 
 \]
 
From the above estimate, we conclude that
\be\label{2021dec20eqn51}
\begin{split}
 &\big| \mathbf{P}\big( G_{k,j;l,n}^{ess; m, i;p,q}(t,x,\zeta)\big)  \big|\\
 &  \lesssim \| \mathfrak{m}(\cdot, \zeta)\|_{\mathcal{S}^\infty} 2^{n+\max\{n, (\gamma_1-\gamma_2)M_{t^{\star}}\}+10\epsilon M_{t^{\star}}} 
\min\{  |  x_{\bot}|^{-1/2} 2^{j+p/2 }, 2^{2j+3p/2},  |  x_{\bot} |^{-1} 2^{  k-p -j     } \} \\
 &\lesssim \| \mathfrak{m}(\cdot, \zeta)\|_{\mathcal{S}^\infty} | x_{\bot}|^{-1/2} 2^{(\gamma_1-\gamma_2-\epsilon )M_{t^{\star}}/2}\big[ 2^{  (\alpha^{\star} -10\epsilon) M_{t^{\star}}    }+  \mathbf{1}_{n\geq  (\gamma_1-\gamma_2-30\epsilon )M_{t^{\star}}} \\
 &\quad \times  \min\{2^{(k+2n)/2 +2\alpha^{\star}  M_{t^{\star}}/ 3 -(\gamma_1-\gamma_2)M_{t^{\star}}/6}, 2^{(k+4 n)/2 + \alpha^{\star}  M_{t^{\star}} }\}   \big].
   \end{split} 
\ee
To sum up, the desired estimate holds  \eqref{2022feb18eqn88}   after combining the obtained estimates  \eqref{2022feb17eqn21},  \eqref{2022feb17eqn11}, and  \eqref{2022feb17eqn13}   for the case $i=2$ and the obtained estimates  \eqref{2022feb21eqn41} ,  \eqref{2021dec15eqn23} ,  \eqref{sep25eqn31}, \eqref{sep25eqn32}, \eqref{sep25eqn10} , and  \eqref{2021dec20eqn51} for the case $i=3,4$.
\end{proof}

\begin{lemma}\label{pointestPartIIS1}
For any $i\in \{2,3,4\}, ( l,r)\in \mathcal{B}_i $, see  \eqref{sep18eqn50}.      Under the assumption of  the part (ii) in Theorem \ref{mainresultsfirstpart},      the following estimate holds if   $m+k\geq -2l+ \epsilon M_{t^{\star}}/5, $
\be\label{2022feb22eqn91}
\begin{split}
\big|  \mathbf{P} \big(\widetilde{K}^{m;p,q,2;1}_{k,j ; n,l,r}(t, x, \zeta )\big)\big|&\lesssim    \| \mathfrak{m}(\cdot, \zeta)\|_{\mathcal{S}^\infty} \big[ \sum_{a \in \{0, 1/4, 1/2\} } |  x_{\bot}|^{-a} 2^{a(\gamma_1-\gamma_2)M_{t^\star} }2^{(\alpha^{\star}-10\epsilon)M_{t^\star}}\\
&\quad +  |  x_{\bot}|^{-a}2^{a(\gamma_1-\gamma_2 )M_{t^\star}  +30\epsilon  M_{t^\star}}   \mathbf{1}_{|  x_{\bot}|\geq 2^{-k-n+ b\epsilon M_t}}  \mathbf{1}_{n\geq -(\alpha^\star/2 + 30\epsilon)M_{t^\star}}\\
&\quad \times    \min\{ 2^{(k+2n)/2  + \alpha^\star M_{t^\star}/2 },2^{(k+4n)/2  +  \alpha^\star M_{t^\star}  } \}  \big], \\
\sum_{i=3,4}\big|  \mathbf{P} \big(\widetilde{K}^{m;p,q,i;1}_{k,j ; n,l,r}(t, x, \zeta )\big)\big|& \lesssim  \min_{b\in\{1,2\}} \| \mathfrak{m}(\cdot, \zeta)\|_{\mathcal{S}^\infty} \big[ \sum_{a \in \{0, 1/2\} } |  x_{\bot}|^{-a} 2^{a(\gamma_1-\gamma_2)M_{t^\star} }2^{(\alpha^{\star}-10\epsilon)M_{t^\star}}\\
& \quad  +  \mathbf{1}_{|  x_{\bot}|\leq 2^{-k-n+ b\epsilon M_t}} |  x_{\bot}|^{-1} 2^{(\gamma_1-\gamma_2)M_{t^\star} } 2^{ \alpha^{\star}M_{t^{\star}} /2}  +  2^{35\epsilon M_{t^\star}}    \\
 &  \quad \times   \big(   \mathbf{1}_{n\in \mathcal{N}_{t^{\star} }^1} \min\big\{  \sum_{a\in \mathcal{T}}  |  x_{\bot}|^{-a} 2^{a(\gamma_1-\gamma_2)M_{t^\star} } 2^{(k+4n)/2   + \alpha^{\star} M_{t^{\star} }},\\
 &\quad  \sum_{a\in \mathcal{T}}  |  x_{\bot}|^{-a} 2^{a(\gamma_1-\gamma_2)M_{t^\star} }   2^{(k+2n)/2 +2\alpha^{\star} M_{t^{\star} }/3 -(\gamma_1-\gamma_2)M_{t^{\star}}/6 } \big\}   \\
 &\quad +  \mathbf{1}_{n\in \mathcal{N}_{t^{\star} }^2} \min\big\{ \sum_{a\in \mathcal{T}}  |  x_{\bot}|^{-a} 2^{a(\gamma_1-\gamma_2)M_{t^\star} }   2^{(k+3n)/2 +3\alpha^{\star} M_{t^{\star} }/4  }, \\
 &\quad   \sum_{a\in \mathcal{T}}  |   x_{\bot}|^{-a} 2^{a(\gamma_1-\gamma_2)M_{t^\star} } 2^{(k+4n)/2 + \alpha^{\star} M_{t^{\star} } }  \big\}  \big) \mathbf{1}_{| x_{\bot}|\geq 2^{-k-n+ b\epsilon M_t}}   \big],\\
\end{split}
\ee
 where the sets  $\mathcal{T} $ and $\mathcal{N}_{t^{\star} }^i, i\in\{1, 2\}$ are defined in \eqref{2024oct8eqn8}. 
\end{lemma}
\begin{proof}
We   use the decomposition obtained in  \eqref{aug8eqn31}. From the obtained estimate  \eqref{2022feb9eqn22}, it suffices to estimate the essential part $\widetilde{K}^{m;p,q,i;1}_{k,j ; n,l,r}(t, x, \zeta )$ for the case $p,q\in(-10M_t, \epsilon M_t]\cap \Z.$  

Based on the size of $i$, we proceed in two main steps as follows. .

\medskip

\noindent \textbf{Step 1.}\quad If $i=2$. 

\medskip 

Recall  \eqref{sep17eqn63},  \eqref{sep4eqn6}, and  \eqref{sep9eqn21}. Due to fact that $m+k\geq-2l+\epsilon M_{t^{\star}}/5$,  we know that $|\tilde{v}-\tilde{\zeta}|\sim 2^r \sim 2^l, r\geq  n + 3\epsilon M_t/4.$ 
  Moreover, if $l\notin [(\gamma_1-\gamma_2-\epsilon )M_{t^\star}, (\gamma_1-\gamma_2+\epsilon)M_{t^\star}]$, then we have $| v_{\bot} |/|v|\sim 2^{\max\{l, (\gamma_1-\gamma_2  )M_{t^\star} \} }. $ From the estimate  \eqref{nov24eqn41}, we can rule out the case $\max\{l, (\gamma_1-\gamma_2  )M_{t^\star} \}+j\geq  (\alpha_t+\epsilon) M_t$ if $l\notin [(\gamma_1-\gamma_2-\epsilon )M_{t^\star}, (\gamma_1-\gamma_2+\epsilon)M_{t^\star}]$. 

  Note that, comparing with the obtained estimate  \eqref{2022feb10eqn1},  the following improved estimate holds for the projection of the coefficients, 
  \be\label{2022feb17eqn1}
\big|   \nabla_v \big(\frac{ \mathbf{P}\big(  m_{E}(v, \omega )  + \hat{ \zeta}\times  m_{B}(v, \omega )\big) }{1+\hat{v}\cdot \omega}       \varphi_{j,n}^{i; r}(v, \zeta)  \varphi_{l; r}(\tilde{v}+\omega )    \big)   \big| \lesssim 2^{-j -l + \max\{l, (\gamma_1-\gamma_2)M_{t^\star}\} +\epsilon M_t}. 
  \ee

  Based on the possible size of $|  x_{\bot} |$, we  proceed in  three sub-steps as follows. 

  \medskip

\textbf{Step 1A.}\quad   If $|  x_{\bot} |\leq 2^{-k-n+2\epsilon M_t}.$

  \medskip

Recall   \eqref{aug8eqn31}. From the estimate \eqref{2022feb17eqn1}, the estimate of kernel in  \eqref{sep6eqn31},    the estimate  \eqref{march18eqn31}  in Lemma \ref{conservationlawlemma}, the volume of support of $\omega, v$, 
    the Cauchy-Schwartz inequality,    we have
 \be\label{2022feb18eqn78}
\begin{split}
\big|  \mathbf{P} \big(\widetilde{K}^{m;p,q,i;1}_{k,j ; n,l,r}(t, x, \zeta )\big)\big|&\lesssim \| \mathfrak{m}(\cdot, \zeta)\|_{\mathcal{S}^\infty} 2^{m-j -l + \max\{l, (\gamma_1-\gamma_2)M_{t^\star}\} +\epsilon M_t} \big(2^{-2m+3j+2l} \big)^{1/2} \\
&\quad \times   \big( \min\{(2^{2n}+2^{-m-k})2^{m+3k+2n-j}, 2^{-2m-j-2l}\} \big)^{1/2} \\
&\lesssim \| \mathfrak{m}(\cdot, \zeta)\|_{\mathcal{S}^\infty} 2^{k+n+10\epsilon M_t}\\
& \lesssim  \| \mathfrak{m}(\cdot, \zeta)\|_{\mathcal{S}^\infty}  |  x_{\bot}|^{-1} 2^{(\gamma_1-\gamma_2)M_{t^\star} } 2^{ ( \alpha^{\star} /2-\epsilon) M_{t^{\star}}}. 
\end{split}
\ee

  \medskip

\textbf{Step 1B.}\quad  If $ 2^{-k-n+\epsilon M_t}\leq |  x_{\bot} |\leq 2^{\epsilon M_t}\max\{2^{m}, 2^{k+n}\}.$
  \medskip

Recall  \eqref{aug8eqn31}. From the estimate  \eqref{2022feb17eqn1}, the estimate of kernel in  \eqref{sep6eqn31},    the estimate  \eqref{march18eqn31}  in Lemma \ref{conservationlawlemma}, the volume of support of $\omega, v$, 
    the Cauchy-Schwartz inequality,    we have
\be\label{2022feb18eqn79}
\begin{split}
&\big|  \mathbf{P} \big(\widetilde{K}^{m;p,q,i;1}_{k,j ; n,l,r}(t, x, \zeta )\big)\big|\\
&\lesssim    \| \mathfrak{m}(\cdot, \zeta)\|_{\mathcal{S}^\infty}2^{m-j -l + \max\{l, (\gamma_1-\gamma_2)M_{t^\star}\} +\epsilon M_t} \big(2^{-2m+3j+2l} \big)^{1/2}    \\
&\quad \times\big( \min\{(2^{2n}+2^{-m-k})2^{m+3k+2n-j}, 2^{-2m-j-2l}, (2^{2n}+2^{-m-k})  2^{m+3j+2l}\} \big)^{1/2} \\
&\lesssim \| \mathfrak{m}(\cdot, \zeta)\|_{\mathcal{S}^\infty} 2^{   \max\{l, (\gamma_1-\gamma_2)M_{t^\star}\} +\epsilon M_t +j/2} \big[  \mathbf{1}_{m+k\geq -2n}  \min\big\{ 2^{-2m-j-2l},  \\
&\quad \big( 2^{m+3k+4n-j}\big)^{2/3} \big(2^{-2m-j-2l}\big)^{1/3}, \big(2^{-2m-j-2l}\big)^{1/3} \big( 2^{2n+m+3j+2l}\big)^{2/3} \big\}\\
&\quad   + \mathbf{1}_{m+k\leq -2n} \min\big\{ 2^{-2m-j-2l},  2^{2k+2n-j} ,  \big(2^{2k+2n-j}\big)^{1/3} \big( 2^{-k+3j+2l}\big)^{2/3} \big\} \big]^{1/2}  \\
&\lesssim  \| \mathfrak{m}(\cdot, \zeta)\|_{\mathcal{S}^\infty}  2^{ \max\{l, (\gamma_1-\gamma_2)M_{t^\star}\} +10\epsilon M_t} \min\{ 2^{-m-l},  2^{k+n }, 2^{4j/3 + n/3+2l/3} \}\\
 & \lesssim\| \mathfrak{m}(\cdot, \zeta)\|_{\mathcal{S}^\infty}  |  x_{\bot}|^{-1/2} 2^{(\gamma_1-\gamma_2-\epsilon)M_{t^\star}/2}\big[ 2^{ (\alpha^{\star}-10\epsilon) M_{t^{\star}} }+ 2^{(k+2n)/2 + \alpha^{\star}M_{t^\star}/3} \mathbf{1}_{n\geq -\iota M_{t^\star}} \big] .
  \end{split}
\ee

 \medskip

\textbf{Step 1C.}\quad  If $  |  x_{\bot} |\geq 2^{\epsilon M_t}\max\{2^{m}, 2^{k+n}\}.$

\medskip

Note that, for any $y\in B(0, 2^{-k-n+\epsilon M_t/2})$, we have $|  x_{\bot} -   y_{\bot} + (t-s)  \omega_{\bot}|\sim | x_{\bot}|$. From the estimate  \eqref{2022feb17eqn1}, the Cauchy-Schwartz inequality, the estimate of kernel in  \eqref{sep6eqn31}, the cylindrical symmetry of the distribution function,     the estimate  \eqref{march18eqn31}  in Lemma \ref{conservationlawlemma}, the volume of support of $\omega, v$, and  the estimate  \eqref{nov24eqn41}  if $| v_{\bot}|\geq 2^{(\alpha_t +\epsilon) M_t}$, if $m+k\leq -2n,$ we have 
\be\label{2024oct20eqn1}
\begin{split}
&\big|  \mathbf{P} \big(\widetilde{K}^{m;p,q,i;1}_{k,j ; n,l,r}(t, x, \zeta )\big)\big|\\
& \lesssim \| \mathfrak{m}(\cdot, \zeta)\|_{\mathcal{S}^\infty} 2^{m-j -l + \max\{l, (\gamma_1-\gamma_2)M_{t^\star}\} +4\epsilon M_t}  \big(2^{-2m } \min\{2^{3j+2l}, 2^{j+2\alpha_t M_t}\}\big)^{1/2} \\
&\quad \times\big( \min\{| x_{\bot}|^{-1} 2^{k+n-j},  2^{-k}\min\{2^{3j+2l}, 2^{j+2\alpha_t M_t}\}\} \big)^{1/2}  \\
& \lesssim \| \mathfrak{m}(\cdot, \zeta)\|_{\mathcal{S}^\infty} 2^{m-j -l + \max\{l, (\gamma_1-\gamma_2)M_{t^\star}\} +4\epsilon M_t}  \min\big\{ \big( |  x_{\bot}|^{-1} 2^{k+n-j}\big)^{1/2}\big(2^{-2m+3j+2l}\big)^{1/2}, \\
&\quad   \big( |  x_{\bot}|^{-1} 2^{k+n-j}\big)^{1/4} \big( 2^{-k}\min\{2^{3j+2l}, 2^{j+2\alpha_t M_t}\} \big)^{1/4} \big(2^{-2m }\min\{2^{3j+2l}, 2^{2j+\alpha_t M_t+l}\}\big)^{1/2} \big\}\\ 
&\lesssim     \| \mathfrak{m}(\cdot, \zeta)\|_{\mathcal{S}^\infty}  2^{ \max\{l, (\gamma_1-\gamma_2)M_{t^\star}\} +10\epsilon M_t}\\
&\quad \times  \min\{ |  x_{\bot}|^{-1/2} 2^{(k+n)/2}, |  x_{\bot}|^{-1/4} 2^{n/4}\min\{2^{j+l/2}, 2^{\alpha_t M_t-l/2}\} \}   \\
&\lesssim  \| \mathfrak{m}(\cdot, \zeta)\|_{\mathcal{S}^\infty} \big[  | x_{\bot}|^{-1/4} 2^{(\gamma_1-\gamma_2)M_{t^\star}/4} 2^{(\alpha^\star-11\epsilon)M_{t^\star}} \\
&\quad  +  |  x_{\bot}|^{-1/2} 2^{(\gamma_1-\gamma_2-\epsilon)M_{t^\star}/2} 2^{(k+2n)/2+ \alpha^\star M_{t^\star}/2 } \mathbf{1}_{n\geq (\gamma_1-\gamma_2-90\epsilon)M_{t^\star}} \big].\\
\end{split}
\ee
Similarly, if $m+k\geq -2n,$ we have 
\be\label{2024Dec8eqn1}
\begin{split}
 &\big|  \mathbf{P} \big(\widetilde{K}^{m;p,q,i;1}_{k,j ; n,l,r}(t, x, \zeta )\big)\big| \\
& \lesssim \| \mathfrak{m}(\cdot, \zeta)\|_{\mathcal{S}^\infty} 2^{m-j -l + \max\{l, (\gamma_1-\gamma_2)M_{t^\star}\} +4\epsilon M_t}   \big(2^{-2m } \min\{2^{3j+2l}, 2^{j+2\alpha_t M_t}\}\big)^{1/2} \\
&\quad \times  \big( \min\{|  x_{\bot}|^{-1} 2^{m+2k+3n-j},  2^{m+2n}\min\{2^{3j+2l}, 2^{j+2\alpha_t M_t}\}, 2^{-2m-j-2l}\} \big)^{1/2} .
\end{split}
\ee
From the above estimate, we have
\be\label{2022feb18eqn80}
\begin{split}
\eqref{2024Dec8eqn1}& \lesssim \| \mathfrak{m}(\cdot, \zeta)\|_{\mathcal{S}^\infty}  2^{m-j -l + \max\{l, (\gamma_1-\gamma_2)M_{t^\star}\} +4\epsilon M_t} \\
&\quad \times  \min\big\{ \big( (|  x_{\bot}|^{-1} 2^{m+2k+3n-j} )^{2/3}(2^{-2m-j-2l})^{1/3} \big)^{1/2}  \big( 2^{-2m  +3j+2l}\big)^{1/2} , \\
&\qquad   \big( (2^{m+3j+2l+2n} )^{2/3}(2^{-2m-j-2l})^{1/3} \big)^{1/2}\big( 2^{-2m  +3j+2l}\big)^{1/2} , \\
&\qquad   \big( 2^{-2m  +5j/3+4\alpha_t M_t/3+2l/3}\big)^{1/2} \big( (2^{m+j+2\alpha_t M_t+2n} )^{2/3}(2^{-2m-j-2l})^{1/3} \big)^{1/2} \big\}\\ 
& \lesssim \| \mathfrak{m}(\cdot, \zeta)\|_{\mathcal{S}^\infty}    2^{ \max\{l, (\gamma_1-\gamma_2)M_{t^\star}\} +10\epsilon M_t} \\
&\quad \times  \min\big\{ |  x_{\bot}|^{-1/3} 2^{2k/3+n-l/3}, 2^{4j/3+2n/3+l/3}, 2^{4\alpha_t M_t/3+2n/3-l}\} \big\}\\
&\lesssim  \| \mathfrak{m}(\cdot, \zeta)\|_{\mathcal{S}^\infty} \big[ 2^{(\alpha^{\star}-11\epsilon)M_{t^\star}} + \mathbf{1}_{n\geq - (\alpha^\star/2 +30\epsilon) M_{t^\star}}\\
&\quad \times |  x_{\bot}|^{-1/4}2^{(\gamma_1-\gamma_2 )M_{t^\star}/4} 2^{(k+2n)/2  + \alpha^\star M_{t^\star}/2 -30\epsilon M_{t^\star} }   \big]. 
\end{split}
\ee
After combining the obtained   estimates  \eqref{2022feb18eqn78}--\eqref{2022feb18eqn80}, our desired first estimate  in \eqref{2022feb22eqn91}  holds. 

\medskip

\noindent \textbf{Step 2.}\quad If $i=3,4$. 

\medskip

Due to the  cutoff functions $\varphi^i_{j,n}(v, \zeta), i\in\{3,4\}$, see  \eqref{sep4eqn6}, and the fact that $m+k\geq-2l+\epsilon M_{t^{\star}}/5$,  the support of the $v, \omega, \xi$,    implies that $|\tilde{v}\times \omega|\lesssim  2^{n+\epsilon M_t} +2^{l-\epsilon M_t/20}, |\tilde{v}\times\omega |\sim 2^l$. Since   $ |\omega \times\tilde{\xi}|\sim 2^l$ if $l> -j$, we have $l\leq n +2\epsilon M_t$ if $l> -j$. As a result, we know that $l\in [-j, \max\{-j, n+2\epsilon M_t\}]\cap \Z. $

Moreover, if $l\leq n-2\epsilon M_t$ and $n> -M_t$, then we have $|\tilde{\zeta}-\tilde{v}|\sim |\tilde{\zeta}\times \tilde{\xi} |\sim 2^n$.
Therefore,   if $n\leq   (\gamma_1-\gamma_2)M_{t^{\star}}-2\epsilon M_{t^{\star}} $ or $ l\leq n-2\epsilon M_t, n \geq (\gamma_1-\gamma_2)M_{t^{\star}}+2\epsilon M_{t^{\star}}$,   we have $|  v_{\bot} |/|v|\sim 2^{\max\{n, (\gamma_1-\gamma_2)M_{t^\star}\}}$. From the estimate  \eqref{nov24eqn41}, it suffices to consider the case $\max\{l, (\gamma_1-\gamma_2  )M_{t^\star} \}+j\leq   (\alpha_t+\epsilon) M_t$ if $n\leq   (\gamma_1-\gamma_2)M_{t^{\star}}-2\epsilon M_{t^{\star}} $ or $ l\leq n-2\epsilon M_t, n \geq (\gamma_1-\gamma_2)M_{t^{\star}}+2\epsilon M_{t^{\star}}$.

 Recall  \eqref{2022feb9eqn71}. As a result of direct computations,  the following improved estimate holds for the projection of the coefficients, 
  \be\label{2022feb17eqn51}
  \begin{split}
& \big|   \nabla_v \big(\frac{ \mathbf{P}\big(  m_{E}(v, \omega )  + \hat{ \zeta}\times  m_{B}(v, \omega )\big) }{1+\hat{v}\cdot \omega}       \varphi_{j,n}^{i; r}(v, \zeta)  \varphi_{l; r}(\tilde{v}+\omega )    \big)   \big|\\
& \lesssim 2^{-j - l-\min\{l,n\} +n + \max\{n, (\gamma_1-\gamma_2)M_{t^\star}\} +2\epsilon M_t}. 
\end{split}
  \ee

Based on the possible size of $p$ and $|  x_{\bot}|$, we proceed in four sub-steps   as follows. 

\medskip

 \textbf{Step 2A.}\quad If $|  x_{\bot}|\leq  2^{-k-n+2\epsilon M_{t^{  }}}$.

\medskip

 From the estimate of kernels in  \eqref{sep6eqn31}, the Cauchy-Schwarz inequality,   the estimate  \eqref{march18eqn31}  in Lemma \ref{conservationlawlemma},   the volume of support of $v$,  we have
\be\label{sep26eqn55}
\begin{split}
\big|\mathbf{P}(\widetilde{ {K} }^{m;p,q,1}_{k,j ;l,n})(t, x, \zeta)\big|&\lesssim \| \mathfrak{m}(\cdot, \zeta)\|_{\mathcal{S}^\infty}     2^{m-j-l-\min\{l,n\}+n+ \max\{ n , (\gamma_1-\gamma_2) M_{t^{\star}}\}+6\epsilon M_{t^{\star}}} \\
&\quad \times  \big(   2^{-2m-j-2l}  \big)^{1/2}\big(2^{-2m+3j+2  \min\{l,n\} } \big)^{1/2} \\
 & \lesssim \| \mathfrak{m}(\cdot, \zeta)\|_{\mathcal{S}^\infty} 2^{ n   +6\epsilon M_{t^{\star}}}  2^{-m-2l}   \\
 &\lesssim \| \mathfrak{m}(\cdot, \zeta)\|_{\mathcal{S}^\infty}  2^{k+n+8\epsilon M_{t^\star}} \\
 &\lesssim \| \mathfrak{m}(\cdot, \zeta)\|_{\mathcal{S}^\infty} 2^{10\epsilon M_{t^{\star}}} |  x_{\bot}|^{-1} .
  \end{split}
 \ee

 \medskip

 \textbf{Step 2B.}\quad If $2^{-k-n+\epsilon M_{t^{  }}}\leq |  x_{\bot}|< 2^{ \epsilon M_{t^{  }}} \max\{2^{m+p},2^{-k-n} \}$, $p\leq l + \epsilon M_{t^{  }}$. 

\medskip

  From the estimate of kernels in  \eqref{sep6eqn31}, the Cauchy-Schwarz inequality,   the estimate  \eqref{march18eqn31}  in Lemma \ref{conservationlawlemma},   the volume of support of $\omega, v$,    we have
\be\label{2021dec16eqn31}
\begin{split}
&\big|  \mathbf{P} \big(\widetilde{K}^{m;p,q,i;1}_{k,j ; n,l,r}(t, x, \zeta )\big)\big| \\& \lesssim\| \mathfrak{m}(\cdot, \zeta)\|_{\mathcal{S}^\infty}   2^{m-j -l-\min\{l,n\}+n+ \max\{ n , (\gamma_1-\gamma_2) M_{t^{\star}}\}+3\epsilon M_{t^{\star}}}    \\
& \quad \times \big(  \min\{ 2^{-2m-j-2l}, 2^{m+ 2p+ 3j+2l} \} \big)^{1/2}  \big(2^{-2m+3j+2\min\{l,n\}} \big)^{1/2} \\
&\lesssim   \| \mathfrak{m}(\cdot, \zeta)\|_{\mathcal{S}^\infty} 2^{j/2-l+n+ \max\{ n , (\gamma_1-\gamma_2) M_{t^{\star}}\}+3\epsilon M_{t^{\star}}} \\
&\quad \times  \big(\min\big\{2^{-2m-j-2l}, (2^{-2m-j-2l})^{1/3}(2^{m+ 2p+ 3j+2l})^{2/3} \big\}\big)^{1/2}\\
&    \lesssim \| \mathfrak{m}(\cdot, \zeta)\|_{\mathcal{S}^\infty} 2^{ n+ \max\{ n , (\gamma_1-\gamma_2) M_{t^{\star}}\}+6\epsilon M_{t^{\star}}}   \min\{2^{-m-2l},       2^{4j/3} \}  \\
&\lesssim \| \mathfrak{m}(\cdot, \zeta)\|_{\mathcal{S}^\infty}  2^{ n+ \max\{ n , (\gamma_1-\gamma_2) M_{t^{\star}}\}+10\epsilon M_{t^{\star}}}   \min\{ 2^{k}, |  x_{\bot}|^{-2},        2^{4j/3} \}  \\
& \lesssim \| \mathfrak{m}(\cdot, \zeta)\|_{\mathcal{S}^\infty} |  x_{\bot}|^{-1/2} 2^{(\gamma_1-\gamma_2-\epsilon )M_{t^{\star}}/2 } \big[ 2^{  ({\alpha}^{\star}-10\epsilon)M_{t^{\star}}}+ 2^{10\epsilon M_{t^{\star}} } \\
&\quad \times \min\{ 2^{(k+2n)/2  + 7 \alpha^{\star}M_{t^{\star}}/12}, 2^{(k+4n)/2  +  \alpha^{\star}M_{t^{\star}}   } \}   \mathbf{1}_{n\geq (\gamma_1-\gamma_2-40\epsilon)M_{t^{\star}}  }\big]. 
\end{split}
\ee

 \medskip
 
 \textbf{Step 2C.}\quad   If  $2^{-k-n+\epsilon M_{t^{  }}}\leq |  x_{\bot}|< 2^{ \epsilon M_{t^{  }}} \max\{2^{m+p},2^{-k-n} \}$, $p\geq  l + \epsilon M_{t^{  }}$.  

  \medskip

Note that, for any $y\in B(0, 2^{-k-n+\epsilon M_{t^{  }}/2 }),$ we have $|  x_{\bot} -   y_{\bot}|\sim |  x_{\bot}|$ and $|  v_{\bot}|/|v|\sim 2^p$.  From the estimate  \eqref{2022feb17eqn51},  the Jacobian of changing coordinates $(\theta, \phi)\longrightarrow (z,w)$  in  \eqref{march18eqn66},   the estimate of kernels in \eqref{sep6eqn31}, the Cauchy-Schwarz inequality,   the estimate  \eqref{march18eqn31}  in Lemma \ref{conservationlawlemma}, and the volume of support of $\omega,v$,    we have 
\be\label{2021dec16eqn26}
\begin{split}
&\big|  \mathbf{P} \big(\widetilde{K}^{m;p,q,i;1}_{k,j ; n,l,r}(t, x, \zeta )\big)\big| \\ & \lesssim \| \mathfrak{m}(\cdot, \zeta )\|_{\mathcal{S}^\infty} 2^{m-j-l-\min\{l,n\}+n+ \max\{ n , (\gamma_1-\gamma_2) M_{t^{\star}}\}+4\epsilon M_{t^{\star}}} \big(2^{-2m+3j+2\min\{l,n\}} \big)^{1/2}  \\
& \quad \times    \big(\min\{2^{m+2p+q+  3j+2l}, |  x_{\bot}|^{-1} 2^{-m-p-q-j} ,    2^{-2m-j-2l}\}\big)^{1/2}\\
& \lesssim \| \mathfrak{m}(\cdot, \zeta )\|_{\mathcal{S}^\infty} 2^{j/2-l +n+ \max\{ n , (\gamma_1-\gamma_2) M_{t^{\star}}\}+4\epsilon M_{t^{\star}}}\\ 
&\quad  \times \big( \min\big\{(2^{m+2p+q+  3j+2l})^{1/2} (|  x_{\bot}|^{-1} 2^{-m-p-q-j})^{1/2}, 2^{-2m-j-2l} \big\}\big)^{1/2}\\
 &  \lesssim \| \mathfrak{m}(\cdot, \zeta )\|_{\mathcal{S}^\infty} 2^{ n+ \max\{ n , (\gamma_1-\gamma_2) M_{t^{\star}}\}+6\epsilon M_{t^{\star}} }     \min\big\{ |  x_{\bot}|^{-1/4} 2^{j-l/2+p/4}, 2^{k} \big\} . \\
\end{split}
\ee

Based on the possible size of $n$, we proceed in three sub-steps as follows.

 \medskip
 
 \quad \textbf{Step 2C-a.}\quad   If $n\leq    (\gamma_1-\gamma_2-2\epsilon)M_{t^{\star}} $.

  \medskip

From the obtained estimate  \eqref{2021dec16eqn26}, we have 
\be\label{2021dec16eqn38}
\begin{split}
&\big|  \mathbf{P} \big(\widetilde{K}^{m;p,q,i;1}_{k,j ; n,l,r}(t, x, \zeta )\big)\big| \\
& \lesssim\| \mathfrak{m}(\cdot, \zeta )\|_{\mathcal{S}^\infty}  2^{ n+  \max\{ n , (\gamma_1-\gamma_2) M_{t^{\star}}\}+10\epsilon M_{t^{\star}} }   \\
&\quad \times   \min\big\{ | x_{\bot}|^{-1/2} 2^{3j/2+p/2}, |  x_{\bot}|^{-1/4} 2^{k/4+ 3j/4 + \alpha^{\star} M_{t^{\star} }/4 } , 2^k\big\} \\
&  \lesssim \| \mathfrak{m}(\cdot, \zeta )\|_{\mathcal{S}^\infty}2^{n +12\epsilon M_{t^{\star}}}  \min\big\{ \big( | x_{\bot}|^{-1/4} 2^{k/4+ (\gamma_1-\gamma_2)M_{t^{\star}}/4 + \alpha^{\star} M_{t^{\star} } }\big)^{2/3} (2^k)^{1/3}  , \\
&\quad | x_{\bot}|^{-1/2}2^{3\alpha M_{t^{\star}}/2} ,  |  x_{\bot}|^{-1/4} 2^{k/4+ (\gamma_1-\gamma_2)M_{t^{\star}}/4 + \alpha^{\star} M_{t^{\star} } } \big\}\\
&    \lesssim \| \mathfrak{m}(\cdot, \zeta )\|_{\mathcal{S}^\infty} \big[ |  x_{\bot}|^{-1/2}2^{(\gamma_1-\gamma_2)M_{t^{\star}}/2} 2^{( \alpha^{\star}   -11\epsilon)M_{t^{\star}}} + 2^{12 \epsilon M_{t^{\star}}} \mathbf{1}_{n\geq  (-  \alpha^{\star} +\gamma_1-\gamma_2) M_{t^{\star} }/2 -30\epsilon M_{t^{\star}}} \\
& \quad \times  \min\{ |  x_{\bot}|^{-1/6}   2^{(k+2n)/2 +2\alpha^{\star} M_{t^{\star} }/3 },  |  x_{\bot}|^{-1/4} 2^{(\gamma_1-\gamma_2)M_{t^{\star}}/4} 2^{(k+4n)/4  + \alpha^{\star} M_{t^{\star} }} \}\big]. \\
\end{split}
\ee

 \medskip
 
 \quad \textbf{Step 2C-b.}\quad  If $n\in [    (\gamma_1-\gamma_2-2\epsilon)M_{t^{\star}},  (\gamma_1-\gamma_2+2\epsilon)M_{t^{\star}}] $. 

  \medskip

From the obtained estimate  \eqref{2021dec16eqn26}  and the fact that $|  v_{\bot}|\sim 2^{p+j},$  we have 
 \be\label{2021dec21eqn21}
\begin{split}
 \big|  \mathbf{P} \big(\widetilde{K}^{m;p,q,i;1}_{k,j ; n,l,r}(t, x, \zeta )\big)\big|  
& \lesssim \| \mathfrak{m}(\cdot, \zeta )\|_{\mathcal{S}^\infty} 2^{2 n +12\epsilon M_{t^{\star}} }     \min\big\{   |  x_{\bot}|^{-3/4} 2^{k/4+  j/4 + 3\alpha^{\star} M_{t^{\star} }/4 } , 2^k\big\} \\
& \lesssim \| \mathfrak{m}(\cdot, \zeta )\|_{\mathcal{S}^\infty} 2^{2 n +12\epsilon M_{t^{\star}} }    \big(  |  x_{\bot}|^{-3/4} 2^{k/4+  j/4 + 3\alpha^{\star} M_{t^{\star} }/4 }\big)^{2/3} \big( 2^k \big)^{1/3}  \\
&\lesssim  \| \mathfrak{m}(\cdot, \zeta )\|_{\mathcal{S}^\infty} |  x_{\bot}|^{-1/2} 2^{(\gamma_1-\gamma_2)M_{t^{\star}}/2 +15\epsilon M_{t^{\star}} }  \\
&\quad \times  \min\{   2^{(k+3n)/2 + \alpha^{\star} M_{t^{\star} }/2 + M_{t^{\star}}/6},   2^{(k+4n)/2 + \alpha^{\star} M_{t^{\star} } }  \}  .
\end{split}
\ee

 \medskip
 
 \quad \textbf{Step 2C-c.}\quad  If $n\geq    (\gamma_1-\gamma_2+2 \epsilon)M_{t^{\star}} $. 

 \medskip

Again, from the obtained estimate  \eqref{2021dec16eqn26}, we have 
\be\label{2021dec16eqn35}
\begin{split}
&\big|  \mathbf{P} \big(\widetilde{K}^{m;p,q,i;1}_{k,j ; n,l,r}(t, x, \zeta )\big)\big|\\ &\lesssim \| \mathfrak{m}(\cdot, \zeta )\|_{\mathcal{S}^\infty} 2^{  6\epsilon M_{t^{\star}}}\big[2^{ 2n }   \min\big\{  |  x_{\bot}|^{-1/4} 2^{k/4+n/4+j},  2^{k} \big\} \mathbf{1}_{l\leq n -\epsilon M_{t^\star}}  \\
&\quad +       \min\big\{ |  x_{\bot}|^{-1/4} 2^{\alpha_t M_t+ 3n/4}, 2^{k+2n} \big\}  \big] \\
& \lesssim \| \mathfrak{m}(\cdot, \zeta )\|_{\mathcal{S}^\infty}  2^{ 10\epsilon M_{t^{\star}}}\big[|  x_{\bot}|^{-1/6} 2^{k/2+7n/6 }\min\{2^{(2\alpha^\star M_{t^\star}+n)/3}, 2^{n+2M_t/3}\} \\
&\quad  + |  x_{\bot}|^{-1/8}2^{(k+2n)/2+ \alpha^\star M_{t^\star}/2+3n/8}\big]\\
  & \lesssim \sum_{a\in\{1/8,1/6\}}   \| \mathfrak{m}(\cdot, \zeta )\|_{\mathcal{S}^\infty}  |  x_{\bot}|^{-a} 2^{ a (\gamma_1-\gamma_2)M_{t^\star}+  20\epsilon  M_{t^{\star} } }\\
  &\quad \times  \min\{ 2^{(k+3n)/2+2 \alpha^{\star} M_{t^{\star} }/3 - (\gamma_1-\gamma_2)M_{t^\star}/6 }, 2^{(k+4n)/2 + \alpha^{\star} M_{t^{\star} } }\} . 
  \end{split}
\ee
Therefore, after combining the above obtained estimates  \eqref{2021dec16eqn38},  \eqref{2021dec21eqn21},  and  \eqref{2021dec16eqn35},  the final estimate in the case we are considering is stated as follows, 
\be\label{2021dec16eqn50}
\begin{split}
&\big|  \mathbf{P} \big(\widetilde{K}^{m;p,q,i;1}_{k,j ; n,l,r}(t, x, \zeta )\big)\big| \\
& \lesssim \| \mathfrak{m}(\cdot, \zeta )\|_{\mathcal{S}^\infty} \big[\sum_{a\in\{1/4, 1/2\}} |  x_{\bot}|^{-a}2^{a(\gamma_1-\gamma_2)M_{t^{\star}} } 2^{( \alpha^{\star}   -11\epsilon)M_{t^{\star}}} + 2^{30 \epsilon M_{t^{\star}}}   \\
& \quad \times \mathbf{1}_{n\in \mathcal{N}_{t^{\star} }^1}   \min\big\{ \sum_{a\in \mathcal{A} }  |  x_{\bot}|^{-a} 2^{a(\gamma_1-\gamma_2)M_{t^{\star}} }2^{(k+2n)/2 +2\alpha^{\star} M_{t^{\star} }/3 -(\gamma_1-\gamma_2)M_{t^{\star}}/6 } ,\\
& \quad    \sum_{a\in \mathcal{A}}   |  x_{\bot}|^{-a} 2^{a(\gamma_1-\gamma_2)M_{t^{\star}} } 2^{(k+4n)/2   + \alpha^{\star} M_{t^{\star} }} \big\}+ |  x_{\bot}|^{-1/2} 2^{(\gamma_1-\gamma_2)M_{t^{\star}}/2 +15\epsilon M_{t^{\star}} }   \\
& \quad \times    \mathbf{1}_{n\in \mathcal{N}_{t^{\star} }^2}  \min\{   2^{(k+3n)/2 +3\alpha^{\star} M_{t^{\star} }/4  },   2^{(k+4n)/2 + \alpha^{\star} M_{t^{\star} } }  \}  \big]\\
\end{split}
\ee 
where $\mathcal{A}:=\{1/8, 1/6, 1/4,1/2\}$ and the index sets $\mathcal{N}_{t^{\star} }^i, i\in\{1,2\},$ are defined in   \eqref{2024oct8eqn8}. 

 \medskip
 
 \textbf{Step 2D.}\quad If $   |  x_{\bot}|\geq  2^{ \epsilon M_{t^{  }}} \max\{2^{m+p},2^{-k-n} \}$ . 

  \medskip

 Note that, for any $y\in B(0, 2^{-k-n+\epsilon M_{t^{  }}/2}),$ we have $|  x_{\bot} -   y_{\bot}+(t-s)  \omega_{\bot} |\sim |  x_{\bot}|$.    From the Cauchy-Schwarz inequality,  the estimate of kernels in  \eqref{sep6eqn31},   the estimate  \eqref{march18eqn31}  in Lemma \ref{conservationlawlemma}, and the volume of support of $v$, the estimate  \eqref{nov24eqn41}  if $|  v_{\bot}|\geq 2^{(\alpha_t+\epsilon)M_t},$ we have
\be\label{sep26eqn28}
   \big|  \mathbf{P} \big(\widetilde{K}^{m;p,q,i;1}_{k,j ; n,l,r}(t, x, \zeta )\big)\big| \lesssim \| \mathfrak{m}(\cdot, \zeta )\|_{\mathcal{S}^\infty} 2^{j/2-l +n+  \max\{ n , (\gamma_1-\gamma_2) M_{t^{\star}}\}+4\epsilon M_{t^{\star}}}  
 \big( J^{m;p }_{k,j ;n,l}(t, x, \zeta )  \big)^{1/2} . 
\ee
   where $J^{m;p }_{k,j ;n,l}(t, x, \zeta )$ is defined in  \eqref{sep26eqn21}. From the obtained estimates  \eqref{sep21eqn71}, \eqref{sep21eqn51},  \eqref{sep21eqn88},  \eqref{sep26eqn23}, the cylindrical symmetry of the distribution function,   the Jacobian of changing coordinates $(y_1,y_2,\theta)\rightarrow x-y+(t-s)\omega$, and the volume of support of $\omega$,  we have 
\[
\begin{split}
|J^{m;p }_{k,j ;n,l}(t, x, \zeta )|&\lesssim 2^{-10M_{t^{\star}}-10j} + |  x_{\bot}|^{-1} 2^{m+ 2k + l -m-k -l +\max\{l,p\} +(l-p)_{-}-j+\epsilon M_{t^{\star}}}\\
& \lesssim  2^{-10M_{t^{\star}} } +|  x_{\bot}|^{-1} 2^{ k+l-j+\epsilon M_{t^{\star}}}.
\end{split}
\]
After combining the above estimate and the obtained estimate  \eqref{sep26eqn28}, we have
\be\label{2021dec16eqn32}
 \big|  \mathbf{P} \big(\widetilde{K}^{m;p,q,i;1}_{k,j ; n,l,r}(t, x, \zeta )\big)\big|\lesssim   2^{ n+  \max\{ n , (\gamma_1-\gamma_2) M_{t^{\star}}\} +5\epsilon M_{t^{\star}}}   | x_{\bot}|^{-1/2} 2^{ k/2  -l/2 }     \| \mathfrak{m}(\cdot, \zeta)\|_{\mathcal{S}^\infty}. 
\ee

Moreover, from the estimate  \eqref{2022feb17eqn51},  the Jacobian of changing coordinates $(\theta, \phi)\longrightarrow (z,w)$  in \eqref{march18eqn66},   the estimate of kernels in  \eqref{sep6eqn31}, the Cauchy-Schwarz inequality,   the estimate  \eqref{march18eqn31} in Lemma \ref{conservationlawlemma}, and the volume of support of $\omega,v$, we have
\be\label{2021dec16eqn33}
\begin{split}
&\big|  \mathbf{P} \big(\widetilde{K}^{m;p,q,i;1}_{k,j ; n,l,r}(t, x, \zeta )\big)\big| \\
& \lesssim  \| \mathfrak{m}(\cdot, \zeta)\|_{\mathcal{S}^\infty}  2^{m-j- l-\min\{l,n\} +n +  \max\{ n , (\gamma_1-\gamma_2) M_{t^{\star}}\} + 4\epsilon M_{t^{\star}} } \\ 
 & \quad  \times\big( \min\{ 2^{m+2p+q+  3j+2l}, |  x_{\bot}|^{-1} 2^{-m-p-q-j}  \} \big)^{1/2}   \big(  2^{-2m+3j+2\min\{l,n\}} \big)^{1/2}  \\
 & \lesssim \| \mathfrak{m}(\cdot, \zeta)\|_{\mathcal{S}^\infty} 2^{ n+  \max\{ n , (\gamma_1-\gamma_2) M_{t^{\star}}\}+ 7\epsilon M_{t^{\star}}}    |  x_{\bot}|^{-1/4} 2^{j-l/2+p/4} .
 \end{split}
\ee
After combining the obtained estimates  \eqref{2021dec16eqn32}  and \eqref{2021dec16eqn33}, we have
\be\label{2021dec16eqn36}
\begin{split}
&\big|  \mathbf{P} \big(\widetilde{K}^{m;p,q,i;1}_{k,j ; n,l,r}(t, x, \zeta )\big)\big| \\
& \lesssim  \| \mathfrak{m}(\cdot, \zeta )\|_{\mathcal{S}^\infty} \big[ |  x_{\bot}|^{-1/4}2^{ (\gamma_1-\gamma_2)M_{t^{\star}}/4} 2^{( \alpha^{\star}   -10 \epsilon)M_{t^{\star}}}  + |  x_{\bot}|^{-1/2} 2^{(\gamma_1-\gamma_2)M_{t^{\star}}/2 +30\epsilon M_{t^{\star}} } \\
 &\quad  \times \big(  \mathbf{1}_{n\in \mathcal{N}_{t^{\star} }^1} \min\{   2^{(k+2n)/2 +2\alpha^{\star} M_{t^{\star} }/3 -(\gamma_1-\gamma_2)M_{t^{\star}}/6 }, 2^{(k+4n)/2   + \alpha^{\star} M_{t^{\star} }} \}\\
 &\qquad  +    \mathbf{1}_{n\in \mathcal{N}_{t^{\star} }^2}  \min\{   2^{(k+3n)/2 +3\alpha^{\star} M_{t^{\star} }/4  },   2^{(k+4n)/2 + \alpha^{\star} M_{t^{\star} } }  \} \big) \big].
\end{split}
\ee 

To sum up,  our desired second estimate  in \eqref{2022feb22eqn91}  holds after   combining the estimates \eqref{sep26eqn55}, \eqref{2021dec16eqn31}, \eqref{2021dec16eqn50}  and  \eqref{2021dec16eqn36}. 
\end{proof}

\begin{lemma}\label{pointestPartIIS2}
For   $( l,r)\in \mathcal{B}_2 $, see  \eqref{sep18eqn50}.      Under the assumption of  the part (ii) in Theorem \ref{mainresultsfirstpart},     the following estimate holds if $m+k\geq -2l+ \epsilon M_{t^{\star}}/5, $  
 \be\label{2022feb18eqn21}
 \begin{split}
 &\big| \mathbf{P}\big(  \widetilde{K}^{m;p,q,2;2}_{k,j ; n,l,r}(t, x, \zeta ) \big) \big|\\
 &\lesssim \min_{b\in\{1,2\}}    \| \mathfrak{m}(\cdot, \zeta)\|_{\mathcal{S}^\infty} \Big[ \sum_{a\in\{ 0, 1/4,1/2\}}  |  x_{\bot}|^{-a }   2^{a(\gamma_1-\gamma_2)M_{t^{\star}} }   2^{  ({\alpha}^{\star}-10\epsilon) M_{t^{\star}}   }  \\
 &\quad + \mathbf{1}_{|  x_{\bot}|\leq 2^{-k-n+b\epsilon M_t}} |  x_{\bot}|^{-1} 2^{(\gamma_1-\gamma_2)M_{t^\star} } 2^{  \alpha^{\star} M_{t^\star}/2}
 + \mathbf{1}_{|  x_{\bot}|\geq 2^{-k-n+b\epsilon M_t}} \big[  2^{    40\epsilon M_{t^{\star}}} \\
  &\quad \times   \min\{ |  x_{\bot}|^{-1/2} 2^{7 \alpha^{\star}  M_{t^{\star}} /8}, |  x_{\bot}|^{-1} 2^{3 \alpha^{\star}  M_{t^{\star}} /8  } \}+ 2^{20\epsilon M_{t^{\star}} }  |  x_{\bot}|^{-1/2} \mathbf{1}_{n \geq -(  {\alpha}^{\star}/2 + 30\epsilon) M_{t^{\star}}  }\\
  &\quad \times   2^{(\gamma_1-\gamma_2)M_{t^{\star}}/2}   \min\{  2^{(k+2n)/2  +  5 {\alpha}^{\star} M_{t^{\star}}/8  },2^{(k+4n)/2 + {\alpha}^{\star} M_{t^{\star}}-(\gamma_1-\gamma_2) M_{t^{\star}} /3 } \}    \big]\Big].\\
  \end{split}
\ee
\end{lemma}
\begin{proof}
We still use the decomposition obtained in  \eqref{sep7eqn61}. From the obtained estimate  \eqref{2022feb9eqn22}, it suffices to estimate the essential part $\widetilde{K}^{m;p,q,2;1}_{k,j ; n,l,r}(t, x, \zeta )$ for the case $p,q\in(-10M_t, \epsilon M_t]\cap \Z.$  

Based on the possible size of $m+2j+2l$ and $n$, we proceed in two main steps   as follows. 

\medskip

\noindent \textbf{Step 1.}\quad  If  $  m+2j+2l\leq  \alpha^{\star}  M_{t^{\star}}  +30\epsilon M_{t^{\star}}$  \textup{or} $ n > -( \alpha^{\star}    /2+30\epsilon) M_{t^{\star}}  $. 

 \medskip

 Based on the possible size of $|  x_{\bot} |$, we proceed in three sub-steps  as follows. 

\medskip

  \textbf{Step 1A.}\quad If $|  x_{\bot} |\leq 2^{-k-n+2\epsilon M_t}.$

  \medskip

 Recall  \eqref{aug8eqn31}.  Similar to the obtained estimates  \eqref{oct8eqn31}--\eqref{2022feb14eqn3}, after splitting into three sub-cases based on the size of $m+k$ and $m+p$ and using the same strategy used in each sub-cases, we have 
\be\label{2022feb18eqn51}
\begin{split}
\big|  \mathbf{P}\big( \widetilde{K}^{m;p,q,i;2}_{k,j ; n,l,r}(t, x, \zeta ) \big) \big| & \lesssim  \| \mathfrak{m}(\cdot, \zeta)\|_{\mathcal{S}^\infty}  2^{k+n +10\epsilon M_t}\\
& \lesssim  \| \mathfrak{m}(\cdot, \zeta)\|_{\mathcal{S}^\infty} |  x_{\bot}|^{-1} 2^{(\gamma_1-\gamma_2)M_{t^\star} } 2^{  (\alpha^{\star} /2-\epsilon) M_{t^\star}}. 
\end{split}
\ee

\medskip

  \textbf{Step 1B.}\quad  If $ 2^{-k-n+ \epsilon M_t}\leq |  x_{\bot} | < 2^{\epsilon M_t}\max\{2^{m},2^{-k-n}\}$.

\medskip

    Similar to the obtained estimates  \eqref{2024oct17eqn81}--\eqref{2022feb14eqn11}, after splitting into three sub-cases based on the size of $m+k$ and $m+p$ and using the same strategy used in each sub-cases, we have 
\be\label{2022feb18eqn52}
\begin{split}
 \big|  \mathbf{P}\big( \widetilde{K}^{m;p,q,i;2}_{k,j ; n,l,r}(t, x, \zeta )\big)   \big| & \lesssim   \| \mathfrak{m}(\cdot, \zeta)\|_{\mathcal{S}^\infty}   \big[  2^{    40\epsilon M_{t^{\star}}}\min\{ |  x_{\bot}|^{-1/2} 2^{7 \alpha^{\star}  M_{t^{\star}} /8}, |  x_{\bot}|^{-1} 2^{3 \alpha^{\star}  M_{t^{\star}} /8  } \}  \\
 &\quad  +    |  x_{\bot}|^{-1/2}   2^{(\gamma_1-\gamma_2)M_{t^{\star}}/2} \big( 2^{  ({\alpha}^{\star}-10\epsilon) M_{t^{\star}}   }  
  + 2^{20\epsilon M_{t^{\star}} } \mathbf{1}_{n \geq -(  {\alpha}^{\star}/2 + 30\epsilon) M_{t^{\star}}  }  \\
&\quad \times \min\{  2^{(k+2n)/2  +  5 {\alpha}^{\star} M_{t^{\star}}/8  },2^{(k+4n)/2 + {\alpha}^{\star} M_{t^{\star}}-(\gamma_1-\gamma_2) M_{t^{\star}} /3 } \} \big) \big].
\end{split}
 \ee

\medskip

  \textbf{Step 1C.}\quad  If $   |  x_{\bot} | \geq 2^{\epsilon M_t}\max\{2^{m},2^{-k-n}\}$.

\medskip

Similar to the obtained estimates  \eqref{aug22eqn10}--\eqref{aug22eqn55}, after splitting into two sub-cases based on the size of $m+k$  and using the same strategy used in each sub-cases, we have 
\be\label{2022feb18eqn53}
\begin{split}
\big|  \mathbf{P}\big( \widetilde{K}^{m;p,q,i;2}_{k,j ; n,l,r}(t, x, \zeta )\big)   \big|& \lesssim   \| \mathfrak{m}(\cdot, \zeta)\|_{\mathcal{S}^\infty}     |   x_{\bot}|^{-1/2}   2^{ (\gamma_1-\gamma_2)M_t/2}\big[ 2^{  ({\alpha}^{\star}-10\epsilon) M_{t^{\star}}   }   + \mathbf{1}_{n \geq -({\alpha}^{\star}/2 + 30) M_{t^{\star}}   } \\
& \quad \times   \min\{  2^{(k+2n)/2  +    {\alpha}^{\star} M_{t^{\star}}/2  },2^{(k+4n)/2 + {\alpha}^{\star} M_{t^{\star}}  } \}  \big].\\
\end{split}
\ee

\medskip

\noindent \textbf{Step 2.}\quad  If  $  m+2j+2l\geq  \alpha^{\star}  M_{t^{\star}}  +30\epsilon M_{t^{\star}}$  \textup{and } $ n \leq  -( \alpha^{\star}    /2+30\epsilon) M_{t^{\star}}  $. 

\medskip

Same as what we did in the proof of Lemma \ref{projlargeregime4}, we localize the magnetic  field further by  using the decomposition in  \eqref{july5eqn1}.  As a result, we have
\be\label{2022feb18eqn22}
 \big|    \widetilde{K}^{m;p,q,i;2}_{k,j ; n,l,r}(t, x, \zeta )    \big| \lesssim    \sum_{ (\tilde{m}, \tilde{k}, \tilde{j}, \tilde{l}) \in \mathcal{S}_1(t)\cup  \mathcal{S}_2(t) }  \big|  \widetilde{K}^{m;p,q; \tilde{m}}_{k,j ; n,l  ; \tilde{k};\tilde{j}, \tilde{l}}(t, x, \zeta )   \big|,
\ee
where $ \widetilde{K}^{m;p,q; \tilde{m}}_{k,j ; n,l  ; \tilde{k};\tilde{j}, \tilde{l}}(t, x, \zeta ) $ is defined in  \eqref{2022feb11eqn1}. Similar to the obtained estimate  \eqref{2022feb14eqn36},  after putting the localized magnetic field in $L^\infty_x$ and using    the first estimate  in  \eqref{aug4eqn10}  in Proposition \ref{meanLinfest}, the following estimate holds if if $(\tilde{m}, \tilde{k}, \tilde{j}, \tilde{l}) \in \mathcal{S}_1(t)$,
\be\label{2022feb18eqn54}
\begin{split}
\sum_{(\tilde{m}, \tilde{k}, \tilde{j}, \tilde{l}) \in \mathcal{S}_1(t) } |\widetilde{K}^{m;p,q; \tilde{m}}_{k,j ; n,l  ; \tilde{k};\tilde{j}, \tilde{l}}(t, x, \zeta ) |  &     \lesssim 2^{ 15\epsilon M_t } 2^{-(m+2j+2l)} 2^{ 2  {\alpha}^{\star}  M_t }  \| \mathfrak{m}(\cdot, \zeta)\|_{\mathcal{S}^\infty} \\
&  \lesssim \| \mathfrak{m}(\cdot, \zeta)\|_{\mathcal{S}^\infty} 2^{(\alpha^{\star}-10\epsilon)M_t}.
\end{split}
\ee
 Moreover,  after splitting into three cases based on  the possible size of $|  x_{\bot}|$ and rerunning the arguments used in obtaining  the  estimates  \eqref{2022feb18eqn31}--\eqref{2022feb14eqn37},  the following estimate holds $(\tilde{m}, \tilde{k}, \tilde{j}, \tilde{l}) \in \mathcal{S}_2(t)$.  
\be\label{2022feb18eqn55}
  \sum_{(\tilde{m}, \tilde{k}, \tilde{j}, \tilde{l}) \in \mathcal{S}_2(t) } |\widetilde{K}^{m;p,q; \tilde{m}}_{k,j ; n,l  ; \tilde{k};\tilde{j}, \tilde{l}}(t, x, \zeta ) |     \lesssim  \sum_{a\in\{ 1/4,1/2\}}  |  x_{\bot}|^{-a}2^{ a(\gamma_1-\gamma_2) M_{t^{\star}}     } 2^{   ( \alpha^{\star} -10\epsilon) M_{t^{\star}}  }.
\ee
To sum up, our desired estimate  \eqref{2022feb18eqn21}  holds from the obtained estimates  \eqref{2022feb18eqn51}--\eqref{2022feb18eqn55}.
\end{proof}

\begin{lemma}\label{pointestPartIIS4}
Let $ i\in\{3,4\}, ( l,r)\in \mathcal{B}_i $, see  \eqref{sep18eqn50}.      Under the assumption of  the part (ii) in Theorem \ref{mainresultsfirstpart},      the following estimate holds if $m+k\geq -2l+ \epsilon M_{t^{\star}}/5, $ 
\be\label{sep26eqn57}
\begin{split}
\big|  \mathbf{P}\big(  \widetilde{K}^{m;p,q,i;2}_{k,j ; n,l,r}(t, x, \zeta )\big) \big|& \lesssim  \min_{b\in\{1,2\}} \| \mathfrak{m}(\cdot, \zeta)\|_{\mathcal{S}^\infty}\big[\sum_{a\in \{0,1/4, 1/2\}} |  x_{\bot} |^{-a}2^{a (\gamma_1-\gamma_2)M_{t^{\star}} + ( \alpha^{\star}-10\epsilon )M_{t^{\star}}} \\
&\quad  +   |  x_{\bot}|^{-1 } 2^{ (\gamma_1-\gamma_2)M_{t^{\star}} +  \alpha^{\star} M_{t^{\star}}/2} \mathbf{1}_{|  x_{\bot} |\leq 2^{-k-n+b\epsilon M_{t^{\star}}}}  \\
&\quad   + 2^{30\epsilon   M_{t^{\star}} } \mathbf{1}_{n\geq -  \alpha^{\star}    M_{t^{\star}}/2 }\mathbf{1}_{|  x_{\bot} |\geq 2^{-k-n+b\epsilon M_{t^{\star}}}} |  x_{\bot} |^{-a}2^{a (\gamma_1-\gamma_2)M_{t^{\star}} }\\
&\quad\times  \min\{ 2^{(k+2n)/2+2\alpha^{\star}M_{t^{\star}}/3},2^{(k+4n)/2+\alpha^{\star}M_{t^{\star}}}\}   \big]. \\
\end{split}
\ee 
\end{lemma}
\begin{proof}
 
 Recall  \eqref{aug8eqn31}. Due to the  cutoff functions $\varphi^i_{j,n}(v, \zeta), i\in\{3,4\}$, see  \eqref{sep4eqn6}, and the fact that $m+k\geq-2l+\epsilon M_{t^{\star}}/5$,  the support of the $v, \omega, \xi$,    implies that $|\tilde{v}\times \omega|\lesssim  2^{n+\epsilon M_t} +2^{l-\epsilon M_t/20}, |\tilde{v}\times\omega |\sim 2^l$. Since   $ |\omega \times\tilde{\xi}|\sim 2^l$ if $l> -j$, we have $l\leq n +2\epsilon M_t$ if $l> -j$, i.e., $l\in [-j, \max\{-j, n+2\epsilon M_t\}]\cap \Z. $ Moreover,   besides the rough estimate  obtained in  \eqref{2022feb17eqn51}, the following improved estimate holds for the $\p_{v_3}$ derivative of the coefficient,
    \be\label{2022feb18eqn70}
    \begin{split}
  &\big|   \p_{v_3} \big(\frac{ \mathbf{P}\big(  m_{E}(v, \omega )  + \hat{ \zeta}\times  m_{B}(v, \omega )\big) }{1+\hat{v}\cdot \omega}       \varphi_{j,n}^{i; r}(v, \zeta)  \varphi_{l; r}(\tilde{v}+\omega )    \big)   \big| \\
  &  \lesssim 2^{-j - l-\min\{l,n\}+n + \max\{n, (\gamma_1-\gamma_2)M_{t^\star}\}+\max\{l,p\} +2\epsilon M_t}. \\
  \end{split}
    \ee

 Based on the sizes of $ |  x_{\bot}|$ and $p$, we proceed in three steps   as follows. 

\medskip
\noindent \textbf{Step 1.} \quad If  $|  x_{\bot} |\leq 2^{-k-n+2\epsilon M_{t^{  }}}$ and $p\leq l +5\epsilon M_{t^{ }}$.
\medskip

From the estimates  \eqref{2022feb17eqn51}  and  \eqref{2022feb18eqn70},     the estimate of kernels in  \eqref{sep6eqn31}, the Cauchy-Schwarz inequality,   the estimate  \eqref{march18eqn31}  in Lemma \ref{conservationlawlemma},      the volume of support of $  \omega,v$, and    the Jacobian of changing coordinates $(y_1,y_2,\theta)\rightarrow x-y+(t-s)\omega$, we have
\be\label{sep26eqn91}
\begin{split}
\big|  \mathbf{P}\big(  \widetilde{K}^{m;p,q,i;2}_{k,j ; n,l,r}(t, x, \zeta )\big) \big| & \lesssim \| \mathfrak{m}(\cdot, \zeta)\|_{\mathcal{S}^\infty} 2^{m-j-  \min\{l,n\} +n+\max\{l,p \}+ \max\{ n , (\gamma_1-\gamma_2) M_{t^{\star}}\} +6\epsilon M_{t^{\star}}}    \\
 &\quad \times   \big(  2^{-2m-j-2l}  \big)^{1/2} \big(2^{2k+ n}  2^{3j+2\min\{l,n\}} \big)^{1/2} \\
 & \lesssim \| \mathfrak{m}(\cdot,  \zeta)\|_{\mathcal{S}^\infty} 2^{n+ \max\{ n , (\gamma_1-\gamma_2) M_{t^{\star}}\} +10\epsilon M_{t^{\star}}}  2^{ k +n/2}   \\
 & \lesssim    \| \mathfrak{m}(\cdot,   \zeta)\|_{\mathcal{S}^\infty} |  x_{\bot}|^{-1 } 2^{ 20 \epsilon M_{t^{\star}} }   . \\
\end{split}
\ee

\medskip
\noindent \textbf{Step 2.} \quad  If  $|  x_{\bot} |\leq 2^{-k-n+2\epsilon M_{t^{ }}}$   and $p\geq l +5\epsilon M_{t^{  }}$.
\medskip

Note that, for this case, we have $\forall y \in B(0, 2^{-k-n+\epsilon M_{t^{  }}/2}), \quad |  x_{\bot} -   y_{\bot} +(t-s)  \omega_{\bot}|\sim 2^{m+p} $.   From  the estimates   \eqref{2022feb17eqn51}  and  \eqref{2022feb18eqn70}, the estimate of kernels in  \eqref{sep6eqn31}, the Cauchy-Schwarz inequality,       the Jacobian of changing coordinates $(y_1,y_2,\theta)\rightarrow x-y+(t-s)\omega$, the cylindrical symmetry of the solution, and the volume of support of $ v$, we have
\be\label{oct9eqn61}
\begin{split}
&\big|  \mathbf{P}\big(  \widetilde{K}^{m;p,q,i;2}_{k,j ; n,l,r}(t, x, \zeta )\big) \big| \\
& \lesssim \| \mathfrak{m}(\cdot, \zeta)\|_{\mathcal{S}^\infty} 2^{m-j-   \min\{l,n\}+n+\max\{l,p\}+  \max\{ n , (\gamma_1-\gamma_2) M_{t^{\star}}\} +6\epsilon M_{t^{\star}}}\\
&\quad \times  \big( 2^{m+3k+2n}\frac{2^{-k-n+\epsilon M_{t^{\star}}}}{2^{m+p} } 2^{-m-k-n  } 2^{-j} \big)^{1/2} \\ 
&\quad  \times       \big(2^{m+3k+2n}\frac{2^{-k-n+\epsilon M_{t^{\star}}}}{2^{m+p} } 2^{-m-k-n  } 2^{3j+2 \min\{l,n\} } \big)^{1/2} \\
& \lesssim  \| \mathfrak{m}(\cdot,   \zeta)\|_{\mathcal{S}^\infty}  |  x_{\bot}|^{-1 } 2^{ 20 \epsilon M_{t^{\star}} }    . \\
\end{split}
\ee

\medskip
\noindent \textbf{Step 3.} \quad If  $  |  x_{\bot} |\geq  2^{-k-n+\epsilon M_{t^{ }}} $. 
\medskip

    For this case, we have $ 
\forall y \in B(0, 2^{-k-n+\epsilon M_{t^{\star}}/2}), \quad | x_{\bot} -   y_{\bot} |\sim |  x_{\bot}|. 
$    From the estimates  \eqref{2022feb17eqn51}  and  \eqref{2022feb18eqn70}, the estimate of kernels in  \eqref{sep6eqn31}, the Cauchy-Schwarz inequality,     the Jacobian of changing coordinates $(\theta, \phi)\longrightarrow (z,w)$  in  \eqref{march18eqn66}, and the volume of support of $\omega,v$, we have
\be\label{sep26eqn81}
\begin{split}
  &\big|  \mathbf{P}\big(  \widetilde{K}^{m;p,q,i;2}_{k,j ; n,l,r}(t, x, \zeta )\big) \big|\\
  & \lesssim  \| \mathfrak{m}(\cdot, \zeta)\|_{\mathcal{S}^\infty}  2^{m-j-  \min\{l,n\} +n+\max\{l,p\}+  \max\{ n , (\gamma_1-\gamma_2) M_{t^{\star}}\} +6\epsilon M_{t^{\star}}}  \\
 &\quad  \times   \big(  2^{m+2p+q+3j+2 \min\{l,n\} } \big)^{1/2} \big(|  x_{\bot} |^{-1}2^{-m-p-q}2^{3j+2l} \big)^{1/2} \\
 & \lesssim   \| \mathfrak{m}(\cdot, \zeta)\|_{\mathcal{S}^\infty} 2^{ m+n +\max\{l,p\}  + \max\{ n , (\gamma_1-\gamma_2) M_{t^{\star}}\} +8\epsilon M_{t^{\star}} }    |  x_{\bot} |^{-1/2} 2^{2j+l+p/2}. \\
 \end{split}
\ee

From the above estimate, we rule out  the case $\max\{-j, n+2\epsilon M_t\}=-j$ easily. Hence, it remains to consider the case $\max\{-j, n+2\epsilon M_t\}>-j$, which implies that $l\leq n + 2\epsilon M_t.$ Moreover, as $|\tilde{v}-\tilde{\zeta}|\leq 2^{n+\epsilon M_t}$, we know that $|  v_{\bot}|\sim 2^{j+(\gamma_1-\gamma_2)M_{t^\star}}$ if $n\leq (\gamma_1-\gamma_2-2\epsilon)M_{t^\star}$. From the estimate   \eqref{nov24eqn41}, it suffices to consider the case $j+(\gamma_1-\gamma_2)M_{t^\star}\leq (\alpha_t+\epsilon)M_t $ if $n\leq (\gamma_1-\gamma_2-2\epsilon)M_{t^\star}$.

Alternatively,  from the estimates  \eqref{2022feb17eqn51}  and \eqref{2022feb18eqn70},  the estimate of kernels in \eqref{sep6eqn31}, the Cauchy-Schwarz inequality,   the estimate  \eqref{march18eqn31}  in Lemma \ref{conservationlawlemma},    the Jacobian of changing coordinates $(y_1,y_2,\theta)\rightarrow x-y+(t-s)\omega$, and the volume of support of $\omega,v$, we have 
\be\label{sep26eqn82}
\begin{split}
& \big|  \mathbf{P}\big(  \widetilde{K}^{m;p,q,i;2}_{k,j ; n,l,r}(t, x, \zeta )\big) \big|\\
&\lesssim  \| \mathfrak{m}(\cdot, \zeta)\|_{\mathcal{S}^\infty} 2^{m-j-  \min\{l,n\}  +n+\max\{l,p\}+  \max\{ n , (\gamma_1-\gamma_2) M_{t^{\star}}\} +6\epsilon M_{t^{\star}}} \big( 2^{-2m-j-2l} \big)^{1/2} \\
& \quad \times \big(\min\{  |  x_{\bot} |^{-1} {2^{-m-p-q}}  , 2^{2k+n +  q }, 2^{m+3k+2n+2p+q}\}  2^{3j+2 \min\{l,n\}  } \big)^{1/2}   \\
& \lesssim \| \mathfrak{m}(\cdot, \zeta)\|_{\mathcal{S}^\infty}  2^{ - l+n+\max\{l,p\}   +  \max\{ n , (\gamma_1-\gamma_2) M_{t^{\star}}\}  +8\epsilon M_{t^{\star}} }    |  x_{\bot} |^{-1/4}  2^{-(m+p)/4}   2^{ (2k+n)/4 } .\\
\end{split}
\ee

  Based on the size of $p$ and $|  x_{\bot}|$, we proceed in  four sub-steps as follows.

 \medskip
  \textbf{Step 3A.} \quad  If $p\leq l+2\epsilon M_{t^{  }}$,  $|  x_{\bot} |\leq 2^{m+p+\epsilon M_t}$. 
 \medskip

  After combining the obtained estimates  \eqref{sep26eqn81}  and   \eqref{sep26eqn82}, we have 
\be\label{2022jan9eqn80}
\begin{split}
 \big|  \mathbf{P}\big(  \widetilde{K}^{m;p,q,i;2}_{k,j ; n,l,r}(t, x, \zeta )\big) \big|& \lesssim  \| \mathfrak{m}(\cdot, \zeta)\|_{\mathcal{S}^\infty} |  x_{\bot}|^{-1/2} 2^{(\gamma_1-\gamma_2-\epsilon)M_{t^{\star}}/2}\\
 &\quad \times \big[ 2^{(\alpha^{\star}-10\epsilon)M_{t^{\star}} }  + 2^{(k+4n)/2+\alpha^{\star}M_{t^{\star}}/2} \mathbf{1}_{n\geq - \alpha^{\star}M_{t^{\star}}/2    }\big]. 
 \end{split}
\ee

  \medskip
  \textbf{Step 3B.} \quad  If  $p\leq l+2\epsilon M_{t^{  }}$,  $|  x_{\bot} |\geq 2^{m+p+\epsilon M_t}$. 
  \medskip

Note that,  $\forall y \in B(0, 2^{-k-n+\epsilon M_{t^{  }}/10}), |  x_{\bot} -   y_{\bot} + (t-s)  \omega_{\bot}|\sim |  x_{\bot} |  $. From   the estimates  \eqref{2022feb17eqn51}  and  \eqref{2022feb18eqn70},  the estimate of kernels in  \eqref{sep6eqn31}, the cylindrical symmetry of the distribution function,  the Cauchy-Schwarz inequality,   the estimate  \eqref{march18eqn31}   in Lemma \ref{conservationlawlemma}, we have
\[
\begin{split}
 &\big|  \mathbf{P}\big(  \widetilde{K}^{m;p,q,i;2}_{k,j ; n,l,r}(t, x, \zeta )\big) \big|\\&\lesssim   \| \mathfrak{m}(\cdot, \zeta)\|_{\mathcal{S}^\infty}   2^{m-j-  \min\{l,n\}  +n+\max\{l,p\}+  \max\{ n , (\gamma_1-\gamma_2) M_{t^{\star}}\} +6\epsilon M_{t^{\star}}}   \big( 2^{-2m-j-2l} \big)^{1/2}  \\
 &\quad \times  \big(\min\{  |  x_{\bot} |^{-1} {2^{-m-p-q}}  ,  |  x_{\bot} |^{-1} 2^{m+2k+ n+2p+q}\} 2^{3j+2 \min\{l,n\}  } \big)^{1/2}    \\
 &\lesssim  |  x_{\bot} |^{-1/2} 2^{(2k+n+p)/4 +n+  \max\{ n , (\gamma_1-\gamma_2) M_{t^{\star}}\} +16\epsilon M_{t^{\star}}}\| \mathfrak{m}(\cdot, \zeta)\|_{\mathcal{S}^\infty}  
 \end{split}
\]
  After combining the above estimate and  the obtained estimate  \eqref{sep26eqn81}, we have 
\be\label{2022feb20eqn61}
\begin{split}
 \big|  \mathbf{P}\big(  \widetilde{K}^{m;p,q,i;2}_{k,j ; n,l,r}(t, x, \zeta )\big) \big| &\lesssim \| \mathfrak{m}(\cdot, \zeta)\|_{\mathcal{S}^\infty} |  x_{\bot}|^{-1/2} 2^{(\gamma_1-\gamma_2)M_{t^{\star}}/2} \big[ 2^{(\alpha^{\star}-10\epsilon)M_{t^{\star}} } \\
 &\quad  + 2^{(k+4n)/2+\alpha^{\star}M_{t^{\star}}/2} \mathbf{1}_{n\geq - \alpha^{\star}M_{t^{\star}}/2 - 20\epsilon M_{t^{\star}}  }\big]. 
 \end{split}
\ee

 \medskip
  \textbf{Step 3C.} \quad  If $p\geq l+2\epsilon M_{t^{  }}$,   $| x_{\bot} |\in [2^{m+p-\epsilon M_t}, 2^{m+p+\epsilon M_t}]$. 
 \medskip

 Note that, for this sub-case, we have $|  v_{\bot}|/|v|\sim 2^p$.   After combining the obtained estimates   \eqref{sep26eqn81}  and   \eqref{sep26eqn82}, we have 
 \be\label{2022jan9eqn81}
 \begin{split}
 & \big|  \mathbf{P}\big(  \widetilde{K}^{m;p,q,i;2}_{k,j ; n,l,r}(t, x, \zeta )\big) \big| \\
  & \lesssim   \| \mathfrak{m}(\cdot, \zeta)\|_{\mathcal{S}^\infty}\big[  2^{(\alpha^{\star}-11\epsilon)M_{t^{\star}} }+   \mathbf{1}_{n\geq - \alpha^{\star}M_{t^{\star}}/2  }  2^{  n +p    + \max\{ n , (\gamma_1-\gamma_2) M_{t^{\star}}\}   + 10\epsilon M_{t^{\star}}  } \\
 &\quad  \times  | x_{\bot}|^{-1/2}  \min\big\{ \big(2^{m+ 2j+l  +p/2}\big)^{1/3} \big(   2^{-l/2} 2^{(m+k)/4} 2^{ (2k+n)/4  } \big)^{2/3},  2^{m/2 + j+ p/4 +(2k+n)/8} \big\}\big]\\
  &\lesssim \| \mathfrak{m}(\cdot, \zeta)\|_{\mathcal{S}^\infty}\big[  2^{(\alpha^{\star}-11\epsilon)M_{t^{\star}} }  + 2^{30\epsilon M_{t^{\star}}}  \mathbf{1}_{n\geq - \alpha^{\star}M_{t^{\star}}/2   }\\
  &\quad \times  \min\{ 2^{(k+2n)/2+2\alpha^{\star}M_{t^{\star}}/3},  2^{(k+4n)/2+\alpha^{\star}M_{t^{\star}}}\} \big]. \\
\end{split}
\ee

 \medskip
  \textbf{Step 3D.} \quad If  If $p\geq l+2\epsilon M_{t^{  }}$, $|  x_{\bot} |\notin [2^{m+p-\epsilon M_t}, 2^{m+p+\epsilon M_t}]$. 
 \medskip

Note that, for this sub-case, we have $|  v_{\bot}|/|v|\sim 2^p$ and  $\forall y \in B(0, 2^{-k-n+\epsilon M_{t^{}}/10}), | x_{\bot} -   y_{\bot} + (t-s)  \omega_{\bot}|\sim \max\{2^{m+p}, |  x_{\bot} | \}.$   From the estimates  \eqref{2022feb17eqn51}  and  \eqref{2022feb18eqn70},  the estimate of kernels in  \eqref{sep6eqn31}, the Cauchy-Schwarz inequality, the cylindrical symmetry of the distribution function,   the Jacobian of changing coordinates $(y_1,y_2,\theta)\rightarrow x-y+(t-s)\omega$,      and the volume of support of $\omega,v$, we have 
\be\label{2022jan9eqn51}
\begin{split}
 & \big|  \mathbf{P}\big(  \widetilde{K}^{m;p,q,i;2}_{k,j ; n,l,r}(t, x, \zeta )\big) \big| \\
 & \lesssim   \| \mathfrak{m}(\cdot, \zeta)\|_{\mathcal{S}^\infty} 2^{m-j-   \min\{l,n\}  +n+\max\{l,p\}+ \max\{ n , (\gamma_1-\gamma_2) M_{t^{\star}}\} +6\epsilon M_{t^{\star}}}   \\
& \quad  \times \big( |  x_{\bot}|^{-1} 2^{-m-p-q-j} \big)^{1/2} \big(   (\max\{2^{m+p}, | x_{\bot} | \})^{-1}  2^{k+q } 2^{3j+2 \min\{l,n\}  } \big)^{1/2}  \\
& \lesssim \| \mathfrak{m}(\cdot, \zeta)\|_{\mathcal{S}^\infty} |  x_{\bot} |^{-1/2} 2^{k/2+n + \max\{ n , (\gamma_1-\gamma_2) M_{t^{\star}}\} +10\epsilon M_{t^{\star}}}.\\
\end{split}
\ee

 After combining the obtained estimates  \eqref{sep26eqn81}     and  \eqref{2022jan9eqn51}, we have 
\be\label{2022jan9eqn82}
\begin{split}
 \big|  \mathbf{P}\big(  \widetilde{K}^{m;p,q,i;2}_{k,j ; n,l,r}(t, x, \zeta )\big) \big| &\lesssim  \| \mathfrak{m}(\cdot, \zeta)\|_{\mathcal{S}^\infty} |  x_{\bot}|^{-1/2} 2^{ (\gamma_1-\gamma_2-\epsilon )M_{t^{\star}}/2 } \\
 &\quad \times   \big[  2^{(\alpha^{\star}-10\epsilon)M_{t^{\star}} } 
+    2^{(k+2n)/2+ \alpha^{{\star}} M_{t^{\star}}/4 }  \mathbf{1}_{n\geq - \alpha^{\star}M_{t^{\star}}/2  }\big]. 
\end{split}
\ee
To sum up, our desired estimate  \eqref{sep26eqn57}  holds from the obtained estimates 
 \eqref{sep26eqn91},  \eqref{oct9eqn61},   \eqref{2022jan9eqn80},  \eqref{2022jan9eqn81}, and  \eqref{2022jan9eqn82}.

\end{proof}
To sum up, the main result of this subsection can be summarized in the following Proposition. 
\begin{proposition}\label{pointpartIIfinal}
        Under the assumption of  the part (ii) in Theorem \ref{mainresultsfirstpart},  for any    $ j\in [0, (1+2\epsilon)M_{t_{ }}]\cap\Z,     n\in [-M_{t_{ }}, 2]\cap \Z,$  we have 
\be\label{2022feb23eqn11}
\begin{split}
&\big|\mathbf{P}\big( T_{k,j;n}^{\mu,2}( \mathfrak{m}, E)(t,x,\zeta) +  \hat{\zeta}\times T_{k,j;n}^{\mu,2}( \mathfrak{m}, B)(t,x, \zeta)\big) \big|\\
&+ 2^{(\gamma_1-\gamma_2)M_{t^{\star}}} \big|  T_{k,j;n}^{\mu,2}( \mathfrak{m}, E)(t,x,\zeta) +  \hat{\zeta}\times T_{k,j;n}^{\mu,2}( \mathfrak{m}, B)(t,x, \zeta)  \big|\\
&\lesssim \min_{b\in\{1,2\}}  \| \mathfrak{m}(\cdot, \zeta )\|_{\mathcal{S}^\infty} \Big[\sum_{a\in \mathcal{T}} |  x_{\bot}|^{-a}2^{a(\gamma_1-\gamma_2)M_{t^{\star}} } 2^{( \alpha^{\star}   -10\epsilon)M_{t^{\star}}} \\
& \quad +    \mathbf{1}_{|  x_{\bot} |\leq 2^{-k-n+b\epsilon M_{t^{\star}}}}    |  x_{\bot}|^{-1} 2^{ (\gamma_1-\gamma_2)M_{t^{\star}} } 2^{  5\alpha^{\star} M_{t^{\star}}   /6 }+  2^{    40\epsilon M_{t^{\star}}}  \mathbf{1}_{|  x_{\bot}|\geq 2^{-k-n+b\epsilon M_t}}  \\
&\quad\times  \big( \min\{ |  x_{\bot}|^{-1/2} 2^{7 \alpha^{\star}  M_{t^{\star}} /8}, |  x_{\bot}|^{-1} 2^{3 \alpha^{\star}  M_{t^{\star}} /8  } \} +  \mathbf{1}_{n\geq -(\alpha^\star/2 + 30\epsilon)M_{t^\star}}  |  x_{\bot}|^{-a}2^{a(\gamma_1-\gamma_2 )M_{t^\star}   } \\
&\quad \times     \min\{2^{(k+2n)/2+ 2\alpha^{\star}M_{t^\star}/3- (\gamma_1-\gamma_2)M_{t^{\star} }/6}, 2^{(k+2n)/2+  \alpha^{\star}M_{t^\star} - (\gamma_1-\gamma_2)M_{t^{\star} }/3 }\}   \big) \Big],  
\end{split}
\ee
  
 \be\label{sep27eqn59}
 \begin{split}
 &\sum_{i=3,4}\big|\mathbf{P}\big( T_{k,j;n}^{\mu,i}( \mathfrak{m}, E)(t,x,\zeta) +  \hat{\zeta}\times T_{k,j;n}^{\mu,i}( \mathfrak{m}, B)(t,x, \zeta)\big) \big| \\
 &+ 2^{(\gamma_1-\gamma_2)M_{t^{\star}}} \big|  T_{k,j;n}^{\mu,i}( \mathfrak{m}, E)(t,x,\zeta) +  \hat{\zeta}\times T_{k,j;n}^{\mu,i}( \mathfrak{m}, B)(t,x, \zeta)  \big|\\
 &\lesssim \min_{b\in\{1,2\}}  \| \mathfrak{m}(\cdot, \zeta )\|_{\mathcal{S}^\infty} \Big[\sum_{a\in \mathcal{T}} |  x_{\bot}|^{-a}2^{a(\gamma_1-\gamma_2)M_{t^{\star}} } 2^{( \alpha^{\star}   -10\epsilon)M_{t^{\star}}}+ \mathbf{1}_{| x_{\bot}|\leq 2^{-k-n+b\epsilon M_{t^{\star}}}} \\
 & \quad  \times        |  x_{\bot}|^{-1} 2^{ (\gamma_1-\gamma_2)M_{t^{\star}} } 2^{  5\alpha^{\star} M_{t^{\star}}   /6 } +   2^{    40\epsilon M_{t^{\star}}} \mathbf{1}_{|  x_{\bot}|\geq 2^{-k-n+b\epsilon M_{t^{\star}}}}  \\
 &\quad \times   \big[ \min\{ |  x_{\bot}|^{-1/2} 2^{7 \alpha^{\star}  M_{t^{\star}} /8}, | x_{\bot}|^{-1} 2^{3 \alpha^{\star}  M_{t^{\star}} /8  } \}  +    \mathbf{1}_{n\in \mathcal{N}_{t^{\star}}^1  } \min\big\{ \sum_{a\in\mathcal{T}}  |  x_{\bot}|^{-a}   2^{a(\gamma_1-\gamma_2)M_{t^{\star}} } \\
 &\quad \times 2^{(k+2n)/2 +2\alpha^{\star} M_{t^{\star} }/3 -(\gamma_1-\gamma_2)M_{t^{\star}}/6 },   
\sum_{a\in\mathcal{T}}  |  x_{\bot}|^{-a} 2^{a(\gamma_1-\gamma_2)M_{t^{\star}} } 2^{(k+4n)/2   + \alpha^{\star} M_{t^{\star} }  - (\gamma_1-\gamma_2)M_{t^{\star} }/3 }  \big\}  \\
& \quad  + | x_{\bot}|^{-1/2} 2^{(\gamma_1-\gamma_2)M_{t^{\star}}/2   } \mathbf{1}_{n\in \mathcal{N}_{t^{\star}}^2  } \min\{   2^{(k+3n)/2 +3\alpha^{\star} M_{t^{\star} }/4  },   2^{(k+4n)/2 + \alpha^{\star} M_{t^{\star} }   - (\gamma_1-\gamma_2)M_{t^{\star} }/3}  \}      \big]  \Big].  \\
\end{split}
\ee
where $ \mathcal{T}:=   \{0,1/8,1/6,1/4,3/8,1/2\}$. 
 
\end{proposition}

\begin{proof}
Recall the estimate of symbols in \eqref{sep22eqn44},   \eqref{2022feb19eqn51},    \eqref{sep25eqn15}, and  \eqref{2022feb20eqn31}. The main difference between the estimate of the projection onto the horizontal plane of the localized acceleration force  $\mathbf{P}\big(T_{k;n}^{\mu,2}( \mathfrak{m} , E)(t,x, \zeta  ) +  \hat{\zeta}\times T_{k;n}^{\mu,2}( \mathfrak{m}, B)(t,x, \zeta ) \big)$ and the localized acceleration force   $ T_{k;n}^{\mu,2}( \mathfrak{m} , E)(t,x, \zeta  ) +  \hat{\zeta}\times T_{k;n}^{\mu,2}( \mathfrak{m}, B)(t,x, \zeta )  $ lies in the extra gain of factor of size $(  2^{\max\{l,n\}} + 2^{(\gamma_1 -\gamma_2)M_{t^{\star}}})$. Therefore,  the estimate of $2^{(\gamma_1-\gamma_2)M_{t^{\star}}}\big| T_{k;n}^{\mu,2}( \mathfrak{m} , E)(t,x, \zeta  ) +  \hat{\zeta}\times T_{k;n}^{\mu,2}( \mathfrak{m}, B)(t,x, \zeta )  \big|$ follows the same argument of the estimate of $\big| \mathbf{P}\big(T_{k;n}^{\mu,2}( \mathfrak{m} , E)(t,x, \zeta  ) +  \hat{\zeta}\times T_{k;n}^{\mu,2}( \mathfrak{m}, B)(t,x, \zeta ) \big)\big|$.

Recall the decomposition in  \eqref{sep18eqn50}   and the estimate of the S-part in   \eqref{sep7eqn61}. The desired estimate  \eqref{2022feb23eqn11}  holds from  the estimate  \eqref{2022feb16eqn53}  in Lemma \ref{smallkfull1}, the estimate  \eqref{2022feb18eqn88}  in Lemma \ref{pointestPartIT}, the estimate  \eqref{2022feb22eqn91}  in Lemma \ref{pointestPartIIS1} and  the estimate  \eqref{2022feb18eqn21}  in Lemma \ref{pointestPartIIS2}. The desired estimate  \eqref{sep27eqn59}  holds from  the estimate  \eqref{2022feb16eqn58}  in Lemma \ref{smallkfull2}, the estimate  \eqref{2022feb16eqn59}  in Lemma \ref{smallprojec1},   the estimate  \eqref{2022feb18eqn88}  in Lemma \ref{pointestPartIT}, the estimate  \eqref{2022feb22eqn91}  in Lemma \ref{pointestPartIIS1},   and the estimate  \eqref{sep26eqn57}  in Lemma \ref{pointestPartIIS4}.
 
\end{proof}

 Since the notations in previous section are carried out in this section.  For notational consistency with Section \ref{linfacceloc}, we adhere to its established convention. Supplementary notations employed herein—excluding those listed in Table \ref{tablesection4}—are systematically delineated in Table    \ref{tablesection5}.
\begin{table}[H]
\centering
 
\begin{tabular}{ |c|c|c|c| } 
 \hline
 Notation & Definition    & Remarks \\ 
\hline
    $\gamma_1$ & Theorem \ref{maintheoremellipitic}[(ii)]   & Measuring $|\zeta_{\bot}|$ \\
    $\gamma_2$ &   & Measuring $|\zeta |$ \\
       \hline
 $\mathbf{P}(v)$   & Definition \ref{varioureldef} &  $\mathbf{P}$ is mainly used for long notations\\ 
  \hline
\end{tabular} 
\caption{Essential notations in section \ref{horizonestpotw}.}\label{tablesection5}
\end{table}

\section{Proof of Theorem \ref{maintheoremellipitic} and Theorem \ref{mainresultsfirstpart} }\label{proofpartofmaintheorems}

In this section, we collect the estimates we have obtained in section \ref{linfacceloc} and section \ref{horizonestpotw} to prove our main results in Theorem \ref{maintheoremellipitic} and Theorem \ref{mainresultsfirstpart}.

The desired estimate \eqref{2022feb25eqn1} in Theorem \ref{maintheoremellipitic} follows directly from the obtained estimates \eqref{sep5eqn100} and \eqref{sep8eqn2} in Lemma \ref{ellplinf}. The desired estimates \eqref{2024oct27eqn1} and \eqref{2022feb24eqn1} follow directly from the estimates \eqref{sep6eqn1} and \eqref{2024nov14eqn41} in Lemma \ref{ellpointestpartI}. Hence concluding the proof of Theorem \ref{maintheoremellipitic}. 

Recall the decomposition \eqref{oct7eqn1}. The desired estimate \eqref{2024oct8eqn1} in Theorem \ref{mainresultsfirstpart} follows directly from the obtained estimate \eqref{sep27eqn1} in Proposition 
\ref{finalestfirst}, the obtained estimates \eqref{2022feb12eqn21} and \eqref{2022feb12eqn22} in Proposition \ref{finalestsmallprop} and the estimate \eqref{sep5eqn100} in Lemma \ref{ellplinf}. 

Moreover, the desired estimates \eqref{2024oct8eqn2} and \eqref{2024oct8eqn5} in Theorem 
 \ref{mainresultsfirstpart} follow  directly from the obtained estimate \eqref{sep27eqn21} in Proposition 
\ref{pointestfinalpartI}, the obtained estimates \eqref{2022feb23eqn11} and \eqref{sep27eqn59} in Proposition \ref{pointpartIIfinal} and the estimate \eqref{sep6eqn1} in Lemma \ref{ellplinf}.  Hence concluding the proof of Theorem \ref{mainresultsfirstpart}.

\end{document}